\documentclass[11pt,a4paper]{article}
\usepackage[numbers,square]{natbib}
\setcitestyle{citesep={,}}
\usepackage[T1]{fontenc}
\usepackage[utf8]{inputenc}
\usepackage[english]{babel}
\usepackage{lmodern,microtype}
\usepackage[a4paper,margin=26mm,headheight=15pt,footskip=12mm]{geometry}
\usepackage{amsmath,amssymb,amsthm,mathtools,mathrsfs}
\usepackage{array,booktabs,longtable}
\usepackage{enumerate}
\usepackage[hypertexnames=false,colorlinks=true,linkcolor=blue!45!black,citecolor=blue!45!black,urlcolor=blue!45!black]{hyperref}
\usepackage{xcolor}
\hypersetup{pdftitle={Abstract Integration in Net Convergence Structures},pdfauthor={Alexandre Reggiolli Teixeira}}
\numberwithin{equation}{section}
\allowdisplaybreaks[2]
\newtheorem{theorem}{Theorem}[section]
\newtheorem{lemma}[theorem]{Lemma}
\newtheorem{proposition}[theorem]{Proposition}
\newtheorem{corollary}[theorem]{Corollary}
\theoremstyle{definition}
\newtheorem{definition}[theorem]{Definition}

\theoremstyle{remark}
\newtheorem{remark}[theorem]{Remark}

\newtheorem{example}[theorem]{Example}

\usepackage{pgf}
\definecolor{dedicationink}{HTML}{42627E}
\definecolor{dedicationline}{HTML}{B5CBDD}
\definecolor{dedicationpale}{HTML}{E0EAF2}
\definecolor{dedicationwash}{HTML}{FBFCFE}
\makeatletter
\renewcommand{\@maketitle}{%
  \begingroup\centering
    \setlength{\parskip}{0pt}%
    {\LARGE\@title\par}%
    \vskip.6em
    {\normalsize\lineskip.35em
      \begin{tabular}[t]{c}\@author\end{tabular}\par}%
    \vskip.4em
    {\small\@date\par}%
  \endgroup
  \par\vskip.2em
}
\makeatother
\newcommand{\dedicationbud}[3]{%
  \pgfsetfillcolor{dedicationwash}
  \pgfsetstrokecolor{dedicationline}
  \pgfsetlinewidth{.35pt}
  \pgfpathcircle{\pgfpoint{#1mm}{#2mm}}{#3mm}
  \pgfusepath{fill,stroke}%
}
\newcommand{\dedicationleaf}[3]{%
  \begin{pgfscope}
    \pgftransformshift{\pgfpoint{#1mm}{#2mm}}
    \pgftransformrotate{#3}
    \pgfsetfillcolor{dedicationpale}
    \pgfpathmoveto{\pgfpoint{0mm}{0mm}}
    \pgfpathcurveto{\pgfpoint{1.1mm}{1.05mm}}{\pgfpoint{3mm}{1.1mm}}{\pgfpoint{4mm}{0mm}}
    \pgfpathcurveto{\pgfpoint{2.8mm}{-.8mm}}{\pgfpoint{1.1mm}{-.9mm}}{\pgfpoint{0mm}{0mm}}
    \pgfusepath{fill}
  \end{pgfscope}%
}

\definecolor{dedicationrose}{HTML}{C8102E}
\definecolor{dedicationrosedeep}{HTML}{8F102B}
\definecolor{dedicationrosebright}{HTML}{E3263D}
\definecolor{dedicationroselight}{HTML}{F26A73}
\newcommand{\dedicationrosepetal}[3]{%
  \begin{pgfscope}
    \pgftransformrotate{#1}
    \pgftransformscale{#2}
    \pgfsetfillcolor{#3}
    \pgfsetstrokecolor{dedicationrosedeep}
    \pgfsetlinewidth{.10pt}
    \pgfpathmoveto{\pgfpoint{-.45mm}{-.15mm}}
    \pgfpathcurveto{\pgfpoint{-1.10mm}{.05mm}}{\pgfpoint{-1.30mm}{.78mm}}{\pgfpoint{-.75mm}{1.12mm}}
    \pgfpathcurveto{\pgfpoint{-.32mm}{1.45mm}}{\pgfpoint{.48mm}{1.2mm}}{\pgfpoint{.80mm}{.63mm}}
    \pgfpathcurveto{\pgfpoint{.34mm}{.83mm}}{\pgfpoint{-.32mm}{.58mm}}{\pgfpoint{-.45mm}{-.15mm}}
    \pgfusepath{fill,stroke}
    \pgfsetstrokecolor{dedicationroselight}
    \pgfsetlinewidth{.13pt}
    \pgfpathmoveto{\pgfpoint{-.88mm}{.88mm}}
    \pgfpathcurveto{\pgfpoint{-.50mm}{1.20mm}}{\pgfpoint{.12mm}{1.13mm}}{\pgfpoint{.44mm}{.87mm}}
    \pgfusepath{stroke}
  \end{pgfscope}%
}
\newcommand{\dedicationroseflower}[4]{%
  \begin{pgfscope}
    \pgftransformshift{\pgfpoint{#1mm}{#2mm}}
    \pgftransformrotate{#3}
    \pgftransformscale{#4}
    \pgfsetfillcolor{dedicationrose}
    \pgfsetstrokecolor{dedicationrosedeep}
    \pgfsetlinewidth{.12pt}
    \pgfpathmoveto{\pgfpoint{0mm}{1.75mm}}
    \pgfpathcurveto{\pgfpoint{.55mm}{1.95mm}}{\pgfpoint{.95mm}{1.50mm}}{\pgfpoint{.94mm}{1.25mm}}
    \pgfpathcurveto{\pgfpoint{1.62mm}{1.36mm}}{\pgfpoint{1.92mm}{.75mm}}{\pgfpoint{1.54mm}{.3mm}}
    \pgfpathcurveto{\pgfpoint{1.98mm}{-.34mm}}{\pgfpoint{1.45mm}{-1mm}}{\pgfpoint{1mm}{-1.03mm}}
    \pgfpathcurveto{\pgfpoint{.72mm}{-1.79mm}}{\pgfpoint{.03mm}{-1.91mm}}{\pgfpoint{-.43mm}{-1.4mm}}
    \pgfpathcurveto{\pgfpoint{-1.11mm}{-1.66mm}}{\pgfpoint{-1.82mm}{-1.1mm}}{\pgfpoint{-1.43mm}{-.57mm}}
    \pgfpathcurveto{\pgfpoint{-1.91mm}{.05mm}}{\pgfpoint{-1.53mm}{.86mm}}{\pgfpoint{-.98mm}{.99mm}}
    \pgfpathcurveto{\pgfpoint{-1.12mm}{1.52mm}}{\pgfpoint{-.48mm}{1.92mm}}{\pgfpoint{0mm}{1.75mm}}
    \pgfusepath{fill,stroke}
    \dedicationrosepetal{8}{1.12}{dedicationrosebright}
    \dedicationrosepetal{82}{1.12}{dedicationrose}
    \dedicationrosepetal{154}{1.12}{dedicationrosebright}
    \dedicationrosepetal{226}{1.12}{dedicationrose}
    \dedicationrosepetal{298}{1.12}{dedicationrosebright}
    \dedicationrosepetal{43}{.61}{dedicationrose}
    \dedicationrosepetal{163}{.61}{dedicationrosebright}
    \dedicationrosepetal{283}{.61}{dedicationrose}
    \pgfsetstrokecolor{dedicationrosedeep}
    \pgfsetlinewidth{.19pt}
    \pgfpathmoveto{\pgfpoint{.23mm}{.22mm}}
    \pgfpathcurveto{\pgfpoint{-.08mm}{.48mm}}{\pgfpoint{-.45mm}{.20mm}}{\pgfpoint{-.27mm}{-.10mm}}
    \pgfpathcurveto{\pgfpoint{-.06mm}{-.35mm}}{\pgfpoint{.29mm}{-.12mm}}{\pgfpoint{.10mm}{.08mm}}
    \pgfusepath{stroke}
  \end{pgfscope}%
}

\newcommand{\dedicationmedallion}[4]{%
  \begin{pgfscope}
    \pgftransformshift{\pgfpoint{#1mm}{#2mm}}
    \pgfsetfillcolor{dedicationwash}
    \pgfpathcircle{\pgfpoint{0mm}{0mm}}{#3mm}
    \pgfusepath{fill}
    \pgfsetstrokecolor{dedicationline}
    \pgfsetlinewidth{.38pt}
    \pgfpathmoveto{\pgfpointpolar{22}{#3mm}}
    \pgfpatharc{22}{158}{#3mm}
    \pgfpathmoveto{\pgfpointpolar{202}{#3mm}}
    \pgfpatharc{202}{338}{#3mm}
    \pgfusepath{stroke}
    \pgfsetstrokecolor{dedicationpale}
    \pgfsetlinewidth{.25pt}
    \pgfpathmoveto{\pgfpointpolar{158}{#3mm}}
    \pgfpatharc{158}{202}{#3mm}
    \pgfpathmoveto{\pgfpointpolar{338}{#3mm}}
    \pgfpatharc{338}{382}{#3mm}
    \pgfusepath{stroke}
    \pgftext[at={\pgfpoint{0mm}{0mm}}]{%
      \parbox{37mm}{\centering\color{dedicationink}%
        \fontsize{8.5}{10.3}\selectfont\itshape
        \setlength{\parindent}{0pt}\setlength{\parskip}{0pt}%
        #4\par}}%
  \end{pgfscope}%
}
\newcommand{\integrationdedication}{%
  \par\noindent\makebox[\textwidth][c]{%
  \begin{pgfpicture}
    \pgfpathrectangle{\pgfpoint{0mm}{0mm}}{\pgfpoint{158mm}{80mm}}
    \pgfusepath{use as bounding box}
    \pgfsetstrokecolor{dedicationline}
    \pgfsetlinewidth{.4pt}
    \pgfsetroundcap
    \pgfpathmoveto{\pgfpoint{79mm}{77.8mm}}
    \pgfpathcurveto{\pgfpoint{65mm}{79mm}}{\pgfpoint{44mm}{75.5mm}}{\pgfpoint{27mm}{74mm}}
    \pgfpathcurveto{\pgfpoint{26.5mm}{66mm}}{\pgfpoint{26mm}{59mm}}{\pgfpoint{27mm}{48mm}}
    \pgfpathcurveto{\pgfpoint{26mm}{32mm}}{\pgfpoint{26mm}{23mm}}{\pgfpoint{27mm}{15mm}}
    \pgfpathcurveto{\pgfpoint{40mm}{14mm}}{\pgfpoint{63mm}{6mm}}{\pgfpoint{79mm}{6mm}}
    \pgfusepath{stroke}
    \pgfpathmoveto{\pgfpoint{79mm}{77.8mm}}
    \pgfpathcurveto{\pgfpoint{93mm}{79mm}}{\pgfpoint{114mm}{75.5mm}}{\pgfpoint{131mm}{74mm}}
    \pgfpathcurveto{\pgfpoint{131.5mm}{66mm}}{\pgfpoint{132mm}{59mm}}{\pgfpoint{131mm}{48mm}}
    \pgfpathcurveto{\pgfpoint{132mm}{32mm}}{\pgfpoint{132mm}{23mm}}{\pgfpoint{131mm}{15mm}}
    \pgfpathcurveto{\pgfpoint{118mm}{14mm}}{\pgfpoint{95mm}{6mm}}{\pgfpoint{79mm}{6mm}}
    \pgfusepath{stroke}
    \pgfpathmoveto{\pgfpoint{79mm}{77.8mm}}
    \pgfpathcurveto{\pgfpoint{76.8mm}{62mm}}{\pgfpoint{81mm}{49mm}}{\pgfpoint{79mm}{35mm}}
    \pgfpathcurveto{\pgfpoint{78mm}{20mm}}{\pgfpoint{79mm}{11mm}}{\pgfpoint{79mm}{1.5mm}}
    \pgfusepath{stroke}
    \pgfsetstrokecolor{dedicationpale}
    \pgfsetlinewidth{.32pt}
    \pgfpathmoveto{\pgfpoint{27mm}{74mm}}
    \pgfpathcurveto{\pgfpoint{44mm}{61mm}}{\pgfpoint{55mm}{41mm}}{\pgfpoint{79mm}{35mm}}
    \pgfpathcurveto{\pgfpoint{101mm}{30mm}}{\pgfpoint{117mm}{19mm}}{\pgfpoint{131mm}{15mm}}
    \pgfusepath{stroke}
    \pgfpathmoveto{\pgfpoint{131mm}{74mm}}
    \pgfpathcurveto{\pgfpoint{114mm}{61mm}}{\pgfpoint{103mm}{41mm}}{\pgfpoint{79mm}{35mm}}
    \pgfpathcurveto{\pgfpoint{57mm}{30mm}}{\pgfpoint{41mm}{19mm}}{\pgfpoint{27mm}{15mm}}
    \pgfusepath{stroke}
    \dedicationleaf{53}{77}{153}
    \dedicationleaf{105}{77}{27}
    \dedicationleaf{78.3}{68.5}{138}
    \dedicationleaf{79}{66}{-35}
    \dedicationleaf{43}{12.5}{151}
    \dedicationleaf{53}{10.5}{213}
    \dedicationleaf{115}{12.5}{29}
    \dedicationleaf{105}{10.5}{-33}
    \dedicationroseflower{53.5}{76.7}{-20}{1.00}
    \dedicationroseflower{104.5}{76.7}{25}{1.00}
    \dedicationroseflower{78.3}{68.5}{-25}{.78}
    \dedicationroseflower{79}{64.7}{35}{.78}
    \dedicationroseflower{43}{12.5}{-15}{.96}
    \dedicationroseflower{53}{10.5}{30}{.84}
    \dedicationroseflower{115}{12.5}{20}{.96}
    \dedicationroseflower{105}{10.5}{-30}{.84}
    \dedicationbud{79}{77.8}{1.5}
    \dedicationbud{27}{74}{1.05}
    \dedicationbud{131}{74}{1.05}
    \dedicationbud{27}{15}{1.5}
    \dedicationbud{131}{15}{1.5}
    \dedicationbud{79}{6}{1.5}
    \dedicationbud{79}{1.5}{.7}
    \dedicationmedallion{27}{48}{24.5}{%
      To the memory of\\[1.5pt]
      {\normalfont\bfseries Edward J. McShane},\\[3pt]
      whose vision opened the way to a general theory of abstract
      integration and continues to illuminate the paths beyond it.}
    \dedicationmedallion{79}{35}{26}{%
      To the memory of\\[1.5pt]
      {\normalfont\bfseries C. S. H\"onig},\\[3pt]
      whose work helped give integration theory deep roots in Brazilian
      soil and a voice of its own. May his memory be a blessing, and may
      the light of his teaching endure from generation to generation.}
    \dedicationmedallion{131}{48}{24.5}{%
      To {\normalfont\bfseries M\'arcia Federson},\\[3pt]
      for the inspiration and generous support that have sustained my
      own mathematical journey, and for being the foremost voice of
      integration theory in our homeland.}
  \end{pgfpicture}}\par
  \vspace{.5em}%
}

\begin{document}

\title{Measure Theory in Convergence Spaces \\ \Large Part II: Abstract Integration in Net Convergence Structures}
\author{Alexandre Reggiolli Teixeira\\[3pt]
\small Departamento de Estat\'istica, Universidade Estadual de Campinas\\
\small Rua S\'ergio Buarque de Holanda, 651, Campinas--SP, Brazil\\
\small \texttt{riccimath123@gmail.com}, \texttt{a163407@dac.unicamp.br}}
\date{October 2026}
\maketitle
\integrationdedication
\begin{abstract}
The construction of a coherent and general integration procedure is the central question
of this second part of our study of measure theory in convergence spaces. We define, and
establish convergence theorems for such process of integration, called the net Riemann
integration, in convergence groups and commutative monoids, preserving the prescribed
convergence throughout. The underlying interchange theorem admits noncommutative
groups and requires right completeness only at the cardinality of the error family. Under
intrinsic stability assumptions, countable error control therefore needs only sequential
completeness, even for arbitrary integration nets, from which we also stablish a general
version of Henstock’s lemma. The resulting uniform theorems reaches incomplete ranges
in which no complete subspace contains a tail of the approximating integrals. In the
ordered locally solid convergence theory side, countability of the inner direction replaces
countable accessibility of the approximants, and local existence arguments yield monotone
and dominated convergence and, under additional controller assumptions, a further extension
through increasing integrable barriers with pointwise coverage and uniformly vanishing excess
primitives. This constructs the limiting integral without a common integrable envelope and
gives uniform control of primitive variation. We finish with two appendices on new and
general results of abstract metrization of net convergences and related topics.
\end{abstract}

\clearpage
\tableofcontents

\section{Introduction}
\label{sec:introduction}

Since the constructions of \citet{bochneroriginalarticle},
\citet{birkhofforiginalarticle} and \citet{pettisoriginalarticle}, and
the subsequent development of bilinear integration by
\citet{bartle1956general}, integration theory has acquired a considerable
variety of procedures, whose differences concern the functions and set
functions admitted, the algebraic structure of their values, and the
approximations through which an integral is obtained. The extent of this
diversity can be appreciated in the treatment of gauge integration by
\citet{henstockmagnum} and in the ordered theories collected by
\citet{integralmeasureandordering} and
\citet{boccutovrabelovariecan}, while the semigroup-valued theory of
\citet{sionsemigroup} provides a broad setting in which the algebra of
these constructions can itself be examined. Across these developments,
additivity, restriction to subsets, and continuity under approximation
occur in forms whose resemblance invites us to examine the operations
responsible for them, both to determine the scope of their generalization
and to understand how the different integral constructions are related.

In this article, the second in our series on measure theory in
convergence spaces, we turn to integration after developing the
extension of measures in the first part \citep{TeixeiraExtension},
following the Riemannian construction in which finite sums associated
with sets, tags and coefficients are arranged in a net and their
convergence determines the integral. The net integral of
\citet{szaznetintegralconvergence}, together with the refinement methods
of \citet{sionsemigroup} and \citet{fleischerchange}, provides a
starting point for a formulation in which the approximation procedure
and the convergence of its values can be varied within a common
framework. We seek assumptions on these data from which an integral
can be constructed and its limit theorems obtained, so that a
specialization identifies the original finite sums and verifies the
conditions permitting their passage to the limit.

To follow these finite sums through their different limiting
procedures, we work in convergence spaces, where the algebraic
operations can be related directly to the prescribed convergent nets
while topological convergence and the various order convergences
remain within the same account. We use the net formulation
developed by \citet{convergencestructuresvanderwalt} and the locally
solid structures studied by \citet{vanderwaldlocallysolid}, together
with the uniform convergence framework of \citet{CookFischer} and
\citet{LechickiZieminska} for the comparisons over a family of
integration domains. A substantial part of the elementary theory
admits commutative monoid values, while the group and lattice
operations needed in the subsequent arguments are introduced where
their particular consequences enter the construction.

Within this setting, uniform approximation brings together two
limiting operations: the integrand varies through one family of
approximations, while each of its integrals is reached through a
further family of finite records. An estimate on those records relates
the two operations before either limit has been taken, leaving us to
determine what allows the estimate to remain effective at their
endpoints. Although completeness enters when an endpoint must be
constructed, a prior question remains even when its existence is
known, namely whether passage to that endpoint preserves the error
family on which the comparison rests. We begin the convergence theory
with a single sequence of continuous functions and a finite positive
scalar measure, through which this distinction becomes explicit.

Even with a complete, Hausdorff and locally convex range admitting
absorbing convex null bases, these functions converge uniformly to
zero in the original convergence while their integrals fail to
converge to zero, although the limiting function is integrable.
Because the construction applies to every finite-sampling procedure
satisfying the stated consistency requirement on the continuous
curves involved, changing the records cannot remove the failure while
retaining those curves and their norm-consistent values. Following
the finite approximation errors to their limits reveals the
obstruction: they form a null family whose adherence no longer tends
to zero, suggesting that stability should first be sought in the
actual limits of nets contained in the error sets. Moreover, the
example satisfies the absorption property used later to combine
independently convergent rows, so its obstruction is located precisely
in the passage to known endpoints.

We develop this observation into two universal characterizations for
separated abelian convergence groups. When the limiting function is
assumed integrable, stability of null error families under one-step
adherence is necessary and sufficient for the parameter-uniform
conclusion over the prescribed biadditive systems with continuous
finite responses. When its integral is to be constructed, the further
requirements are completeness and absorption of independently null
inner errors by a null outer error family, and we recover each of
these requirements from unit-coefficient sampling systems admitted by
the same integration framework. The hypotheses concern nets and their
tail sets before any candidate integral is known, while the necessity
proofs identify the role of each hypothesis within the universal
assertion, whose quantification over the processes and the parameter
families is retained explicitly.

Following this construction more closely also identifies how much
of each property a particular approximation uses. An error base of
cardinality at most an infinite cardinal $\kappa$ requires adherence
and absorption only at that cardinality and permits the endpoint to
be constructed from a Cauchy net with at most $\kappa$ indices, while
the original approximation index and all integration indices remain
unchanged, even when the selected endpoint net is not cofinal in the
approximation index. Thus countable error control requires only
sequential completeness and tests the two stability properties on
countable outer error bases, although the nets within each integration
row may be arbitrary. Two passages through one-step adherence yield
the endpoint estimate, after which absorption constructs the candidate
integral; neither idempotence of adherence nor a countable base for
the original convergence is assumed in this argument.

To connect this characterization with the abstract convergence
theorems developed below, we derive its hypotheses from the algebra
of their original bounds. For ordered controls, finite intersections
and removal of a null added bound prove that the controlling sets
are convergence-closed and that every attained family of candidate
bounds has infimum zero, giving adherence stability and unrestricted
absorption. For Dedekind $\sigma$-complete weakly
$\sigma$-distributive regulator controls, the same closedness
argument preserves adherence, while bounded fusion of the separately
chosen row regulators establishes countable absorption. This second
verification does not identify order and regulator convergence for
arbitrary nets, and it retains the independent regulators that the
different rows provide. Both verifications refer to the original
tolerance relations, so the intrinsic characterization is established
without local convexity or an abstract metrization theorem.

Once these properties of the range have been identified, the finite
coefficients bring us back to the particular integration procedure.
Here repeated tags are first grouped and the actual finite words are
bounded on the original input control sets, so that their response
to a null perturbation follows from finite intersections of input
controls and continuity of the resulting bound, before an integral
of the limiting function is introduced. These estimates give the
bilinear uniform theorem and its order and lattice-normed
specializations, while comparisons with row-dependent residuals
provide another route in which normalization supplies the required
asymptotic estimate. We retain this latter route with its own
hypotheses, since the normalization used in the ordered uniform
theorem need not supply a global finite-response bound, and the
broader comparison theorems also apply to uniform spaces and monoids
for which subtraction is unavailable.

We begin in Section~\ref{sec:foundations} with the prescribed net
convergence and the algebraic structures used by the integral,
specifying tail equivalence and saturation so that every later
comparison is expressed by actual nets and their tail images.
Section~\ref{sectionnetriemannintegral} defines the net Riemann
integral and establishes its elementary properties, including
composition with continuous homomorphisms, additivity, subset
integration and the integration of simple functions, with the
required conditions on records stated independently of the integral
being constructed. Section~\ref{sec:uniform-limits} begins with the
obstruction and develops the intrinsic interchange and universal
characterization, followed by verification from ordered and
regulator controls and by the finite coefficient certificates.
The noncommutative sufficient theorem uses the right convergence
structure, while the support examples show that the cardinal
completeness requirement does not imply full completeness or
containment of an endpoint tail in a complete subspace. The section
also includes the ordered theorem of \citet{monteirouniform} and
the disjoint quasi-uniform specialization associated with
\citet{millingtonproduct}, together with the asymptotic and conjugate
variants for their respective coefficient hypotheses.

From these global comparisons we pass in Section~\ref{lc:section}
to division systems and a Henstock lemma through which estimates on
full records reach partial records, providing the additional
structure needed for monotone and dominated convergence. Here the construction
of integrals for lattice operations is part of the argument, and
the local estimates must be preserved across compatible cells if
the resulting primitives are to retain uniform control. We develop
both the integrable-envelope and bounded-primitives approaches in
connection with the gauge and ordered theories of
\citet{henstockmagnum}, \citet{integralmeasureandordering} and
\citet{boccutovrabelovariecan}, keeping their additional lattice
hypotheses distinct from those of the uniform characterization.
Increasing integrable barriers then permit approximation beyond a
common integrable envelope, because their pointwise coverage and the
uniformly vanishing excess primitives construct the limiting
modulus before it is used as an envelope. The controller selection
assumptions enter at that construction step, and
Theorem~\ref{vt:vitali} and Example~\ref{vt:strictness} record the
resulting extension and its strictness within the stated controller
class.

To support these verifications and clarify the relations among their
hypotheses, we develop the required order and representation theory
in the appendices. Appendix~\ref{app:order}
contains the regulator calculus, bounded fusion and positive
transport, with the distinctions between sequential and arbitrary-net
arguments retained in each statement, while
Appendix~\ref{app:metrization} studies exact representations and
their intrinsic conditions in the original convergence. Its results
include the characterization of invariant partially ordered group
representations, closed-bound realizations, preservation of
continuity under uniform limits, and the effects of changing the
space of values, in relation to the work of
\citet{Reichel1977distance}, \citet{rr1977total},
\citet{FeldmanPorter1973} and \citet{BBW2025countability}.
The examples there distinguish representation, adherence,
completeness and absorption, complementing the integration
obstruction in the main text without replacing the intrinsic proofs.

These distinctions also prepare a broader comparison of integral
constructions, being developed in \emph{A Galaxy of Integrals}
\citep{TeixeiraGalaxy}, where the
approximation principles studied here will be used to examine which
different prescriptions of sums lead to the same integral and which
convergence properties remain attached to the prescription. In the
present article we establish the elementary and convergence theory
that supports such a comparison, with the universal characterization
making explicit why the same finite estimates can yield a uniform
theorem in some convergences and fail after passage to their limits
in others.

We write $\operatorname{Fin}(T)$ for the finite subsets of $T$.
For an element of a lattice-ordered group,
\[
 |x|=x\vee(-x),\qquad x^+=x\vee0,\qquad x^-=(-x)\vee0.
\]
Multiplication by a positive integer denotes repeated addition.
Scalar actions are specified only when present in the data, and
every integral is defined using the prescribed convergence of its
range.

\section{Net convergence and faithful representations}\label{sec:foundations}

\subsection{Nets and the original convergence}
\label{fd:net-foundations}

To follow the successive choices of cells, tags and coefficients through
the integration procedure, we describe its approximations by nets and
take their limits in the prescribed convergence of the range. Our first
task is therefore to fix conventions under which these choices and their
limits can be treated together, using the net formulation of
\cite{convergencestructuresvanderwalt} and, in the ordered setting, the
locally solid convergence framework of \cite{vanderwaldlocallysolid}.

A directed set is a nonempty directed preorder. For a net
$x=(x_i)_{i\in I}$ in a nonempty set $X$, put
\[
 T_i(x)=\{x_j:j\ge i\},\qquad
 \operatorname{Tail}(x)
 =\{B\subseteq X:(\exists i\in I)\ T_i(x)\subseteq B\}.
\]
Thus $B\in\operatorname{Tail}(x)$ abbreviates the assertion that
$x_i\in B$ eventually; allowing a preorder permits repeated indices
without changing that assertion, so that $\operatorname{Tail}(x)$
records only statements about the given net.

\begin{definition}[Admissible nets and saturation]
\label{fd:admissible}
A net $y=(y_k)_{k\in K}$ is a quasi-subnet of $x$ if
\[
 (\forall i\in I)(\exists k\in K)\quad T_k(y)\subseteq T_i(x).
\]
Two nets are tail equivalent if each is a quasi-subnet of the other.
Let $\mathcal B$ be a nonempty family of nonempty subsets of $X$
satisfying
\[
 (\forall B,C\in\mathcal B)(\exists D\in\mathcal B)
 \quad D\subseteq B\cap C.
\]
Its saturation is the net
\[
 \begin{gathered}
 D_{\mathcal B}=\{(B,b):B\in\mathcal B,\ b\in B\},
 \qquad n_{\mathcal B}(B,b)=b,\\
 (B,b)\preceq(C,c)\iff C\subseteq B.
 \end{gathered}
\]
We write $\langle B:B\in\mathcal B\rangle$ for this net. For a
family whose finite intersections are nonempty, the same brackets mean
the saturation of its nonempty finite intersections. A family of sets
is null when its saturation converges to zero.
The admissible nets are these saturations and the constant nets on
singleton directed sets. Any other net is interpreted through a
tail-equivalent admissible representative.
\end{definition}

For the saturation just defined, directedness of $D_{\mathcal B}$
follows because a member $D\subseteq B\cap C$ and a point $d\in D$
give $(D,d)\succeq(B,b),(C,c)$. Moreover, its tail at $(B,b)$ has
image exactly $B$, since every later value belongs to $B$, while
$(B,u)\succeq(B,b)$ for every $u\in B$. These observations give
\begin{equation}\label{fd:tail-representation}
 \begin{gathered}
 \operatorname{Tail}(n_{\mathcal B})
 =\{A\subseteq X:(\exists B\in\mathcal B)\ B\subseteq A\},\\
 \operatorname{Tail}(n_{\operatorname{Tail}(x)})
 =\operatorname{Tail}(x).
 \end{gathered}
\end{equation}
Every saturation index $D_{\mathcal B}$ is a subset of
$\mathcal P(X)\times X$, and a constant net at $x$ has the canonical
saturation representative indexed by $\{(\{x\},x)\}$.
Consequently the possible saturation representatives form a set,
while replacing a net by such a representative changes no tail
assertion and imposes no restriction on the original nets of
integration records.

\begin{definition}[Net convergence]
\label{netconvergencestructure}
A convergence $\lambda$ on $X$ assigns limits to admissible nets and
satisfies:
\begin{enumerate}[(i)]
\item Every constant net with value $x$ converges to $x$.
\item Every quasi-subnet of a net converging to $x$ converges to $x$.
\item If $x_i\to_\lambda x$ and $y_i\to_\lambda x$ on the same
      directed set, then
      $z_i\in\{x_i,y_i\}$ for every $i$ implies $z_i\to_\lambda x$.
\end{enumerate}
Tail-equivalent nets have the same limits. The convergence is separated
if each net has at most one limit. When $\mathcal B$ is a directed
family as in Definition~\ref{fd:admissible}, the notation
$\mathcal B\to_\lambda x$ abbreviates $n_{\mathcal B}\to_\lambda x$.
\end{definition}

\begin{proposition}[Tail calculus]\label{fd:tail-calculus}
Let $x=(x_i)_{i\in I}$ and $y=(y_k)_{k\in K}$ be nets in $X$, and
let $T:X\to Y$. Then
\[
 \operatorname{Tail}(y)\supseteq\operatorname{Tail}(x),\quad
 x\to_\lambda a
 \quad\Longrightarrow\quad y\to_\lambda a.
\]
If both nets converge to $a$, their mixture
\[
 m_{i,k,0}=x_i,\qquad m_{i,k,1}=y_k,
 \qquad (i,k,\varepsilon)\preceq(i',k',\varepsilon')
 \iff i\preceq i',\ k\preceq k'
\]
also converges to $a$, and
\[
 \operatorname{Tail}(m)
 =\operatorname{Tail}(x)\cap\operatorname{Tail}(y).
\]
The image net satisfies
\[
 \operatorname{Tail}(Tx)
 =\{B\subseteq Y:T^{-1}(B)\in\operatorname{Tail}(x)\}.
\]
\end{proposition}
\begin{proof}
For every $i\in I$, the set $T_i(x)$ belongs to
$\operatorname{Tail}(x)\subseteq\operatorname{Tail}(y)$, so there
is $k\in K$ with $T_k(y)\subseteq T_i(x)$. Thus $y$ is a
quasi-subnet of $x$, and axiom~(ii) proves the first assertion.
On the product index
$I\times K\times\{0,1\}$, each coordinate net is a quasi-subnet of
its original net. Axiom~(iii) therefore gives $m\to a$. Its tail at
$(i_0,k_0,\varepsilon)$ is $T_{i_0}(x)\cup T_{k_0}(y)$, so
\[
 \begin{aligned}
 B\in\operatorname{Tail}(m)
 &\iff(\exists i_0,k_0)\quad
 T_{i_0}(x)\cup T_{k_0}(y)\subseteq B\\
 &\iff B\in\operatorname{Tail}(x)\cap\operatorname{Tail}(y).
 \end{aligned}
\]
Finally,
\[
 (\exists i_0)(\forall i\ge i_0)\ Tx_i\in B
 \iff
 (\exists i_0)(\forall i\ge i_0)\ x_i\in T^{-1}(B),
\]
which proves the image assertion.
\end{proof}

In algebraic constructions we sometimes denote the derived eventual
sets of a specified net by a capital letter, such as
$\mathcal F=\operatorname{Tail}(x)$. In that notation
$\mathcal F\to a$ means $x\to a$. We put
\[
 [a]=\operatorname{Tail}((a)),\qquad
 [A]=\operatorname{Tail}(\langle A\rangle)
       \quad(\varnothing\ne A\subseteq X),
\]
where $(a)$ is the constant net and $\langle A\rangle$ is the
saturation of $\{A\}$. Thus square brackets in these expressions
record eventual subsets, while angle brackets always denote a net.
If a saturation is used in an operation on eventual-set records, we
write its $\operatorname{Tail}$ explicitly. An image $T_*\mathcal F$
means the eventual sets of
$(Tx_i)$, a product $\mathcal F\times\mathcal G$ those of
$(x_i,y_k)_{(i,k)}$, and, in an additive group,
$\mathcal F+\mathcal G$ those of $(x_i+y_k)_{(i,k)}$.
An intersection $\mathcal F\cap\mathcal G$ is represented by the
mixture in Proposition~\ref{fd:tail-calculus}. Whenever a common
refinement is needed, its representative is the saturation of the
finite intersections of the prescribed tail sets, and its construction
requires those intersections to be nonempty. Every later statement
using these abbreviations is therefore a statement about the displayed
nets or their saturations.

A map $T:(X,\lambda_X)\to(Y,\lambda_Y)$ is continuous if
\[
 x_i\to_{\lambda_X}x\quad\Longrightarrow\quad
 Tx_i\to_{\lambda_Y}Tx.
\]
Products carry coordinatewise convergence.

\begin{definition}[Convergence groups]\label{productconvergence}
An abelian convergence group $(G,\lambda)$ is an abelian group equipped
with a convergence for which addition and inversion are continuous.
Equivalently,
\[
 x_i\to x,\quad y_k\to y
 \quad\Longrightarrow\quad
 (x_i+y_k)_{(i,k)}\to x+y,
 \qquad -x_i\to-x.
\]
A net $x_i$ is Cauchy if
\[
 (x_i-x_k)_{(i,k)\in I\times I}\longrightarrow_\lambda0.
\]
The group is complete if every Cauchy net converges. Sequential
completeness requires this assertion only for Cauchy sequences.
\end{definition}

Completeness will enter only where a limit must be constructed;
regardless of that property, every convergence group satisfies
\[
 x_i\to x
 \iff x_i-x\to0,
 \qquad
 x_i\to x
 \Longrightarrow
 x_i-x_k=(x_i-x)-(x_k-x)\to0.
\]
If a vector-space structure is present, continuity of its scalar action
is specified separately wherever needed. The group arguments do not use
such an action.

\begin{definition}[Commutative convergence monoids]
\label{fd:convergence-monoid}
A commutative convergence monoid is a commutative monoid $(Z,+,0)$
equipped with a convergence such that
\[
 x_i\to x,\quad y_k\to y
 \quad\Longrightarrow\quad
 (x_i+y_k)_{(i,k)}\to x+y.
\]
A monoid homomorphism preserves addition and $0$. Separation refers to
uniqueness of original limits. Whenever a Cauchy or uniform assertion
is made for a monoid, its uniform convergence structure is specified
separately and is required to induce the stated point convergence.
No difference map is assumed on a monoid.
\end{definition}

\begin{definition}[Locally solid convergence]\label{fd:locally-solid}
Let $G$ be an abelian lattice-ordered group. A group convergence
$\lambda$ on $G$ is locally solid if
\[
 y_i\to_\lambda0,\qquad |x_i|\le|y_i|\quad(i\in I)
 \quad\Longrightarrow\quad x_i\to_\lambda0.
\]
For $B\subseteq G$, put
\[
 \operatorname{sol}(B)=\{x\in G:(\exists b\in B)\ |x|\le|b|\}.
\]
\end{definition}

\begin{lemma}[Domination on arbitrary tails]\label{fd:solid-tail}
Let $\lambda$ be locally solid and let $y_k\to_\lambda0$. If
\[
 (\forall k_0)(\exists i_0)(\forall i\ge i_0)(\exists k\ge k_0)
 \quad |x_i|\le|y_k|,
\]
then $x_i\to_\lambda0$.
\end{lemma}
\begin{proof}
On
\[
 D=\{(k,u):k\in K,\ u\in G,\ |u|\le|y_k|\},
 \qquad (k,u)\preceq(k',u')\iff k\preceq k',
\]
consider the nets $v_{(k,u)}=u$ and $w_{(k,u)}=y_k$. The preorder
is directed because $(k',0)$ belongs to it for every $k'$, and the
net $w$ is tail equivalent to $y$. Thus $w\to0$, and local solidity
gives $v\to0$. Its tail above $(k_0,u_0)$ is
\[
 \{u:(\exists k\ge k_0)\ |u|\le|y_k|\}
 =\operatorname{sol}(T_{k_0}(y)).
\]
The assumed domination says that every such tail contains a tail of
$x$. Hence $x$ is a quasi-subnet of $v$ and converges to zero.
\end{proof}

With domination available on arbitrary tails, the familiar lattice
inequalities also show that, in every separated locally solid convergence
group, lattice operations are continuous and order inequalities pass to
existing limits. We shall use the estimates
\[
 \bigl||x|-|y|\bigr|\le|x-y|,
 \qquad |x^+-y^+|\le|x-y|,
 \qquad |(x\vee z)-(y\vee w)|\le|x-y|+|z-w|.
\]
For example, if $x_i\ge0$ and $x_i\to x$, then
\[
 0\le x^-\le|x-x_i|\to0,
 \qquad x^-=0,
 \qquad x\ge0.
\]
Here the middle equality follows by applying local solidity to the
constant net $x^-$ and then using separation; subtraction of two
convergent nets yields the corresponding assertion for a varying
inequality.

\subsection{Original convergence and faithful order representations}
\label{lt:representation-section}

Although the intrinsic convergence theorems are formulated directly in
the original convergence groups, their applications often supply ordered
bounds through which the hypotheses can be verified. To relate these
bounds to the original convergence, we now introduce the corresponding
notation, including the pair convergence needed for Cauchy nets and for
estimates uniform over parameter families. Throughout this comparison,
assertions quantify over nets, with images and saturations understood
as in Definition~\ref{fd:admissible} and Proposition~\ref{fd:tail-calculus}.

Order convergence and its regulator form are recalled in the order
appendix. For a positive map $r$ and an admissible net $(y_\alpha)$,
\[
 r(y_\alpha)\longrightarrow_o0
 \iff
 (\exists u_\theta\downarrow0)(\forall\theta)(\exists\alpha_0)
 (\forall\alpha\ge\alpha_0)\quad r(y_\alpha)\le u_\theta.
\]
Reading this condition from the tail sets of $(y_\alpha)$ gives the
form recorded in Definition~\ref{lt:order-null}; we write $\chi_j=o$ or $\chi_j=D$
for the specified controller convergence.

\begin{definition}[Lattice-valued gauges]\label{fd:lattice-gauges}
Let $G$ be an abelian group and $F$ an abelian lattice-ordered group.
A symmetric lattice pseudonorm is a map $\rho:G\to F^+$ with
\[
 \rho(0)=0,\qquad \rho(-x)=\rho(x),\qquad
 \rho(x+y)\le\rho(x)+\rho(y).
\]
It is a lattice norm if $\rho(x)=0$ implies $x=0$. A family
$(\rho_j)_{j\in J}$ is separating if
\[
 [\rho_j(x)=0\ (j\in J)]\quad\Longrightarrow\quad x=0.
\]
On a lattice-ordered group, a gauge is solid if
$|x|\le|y|$ implies $\rho(x)\le\rho(y)$.
No scalar homogeneity is required. Asymmetric gauges, when used, are
identified explicitly and retain both orientations in their estimates.
\end{definition}

For uniform convergence structures we use the relational formulation
of \cite{CookFischer,LechickiZieminska}, expressed through pair nets.
If $\Phi=(z_\alpha,w_\alpha)$ and $\Psi=(u_\beta,v_\beta)$ are pair
nets, their transposes are obtained by reversing the coordinates.
Their relational composition is the saturation of
\[
 T_\alpha(\Phi)\circ T_\beta(\Psi)
 =\{(z,v):(\exists w)\ (z,w)\in T_\alpha(\Phi),
                              \ (w,v)\in T_\beta(\Psi)\},
\]
when all these sets are nonempty. They form a decreasing family on the
product of the two index sets, so this saturation is defined.

A uniform convergence structure $\mathcal U$ is a class of pair nets,
identified up to tail equivalence, which contains the constant net
$(z,z)$ for every $z\in Z$ and is stable under quasi-subnets, finite
mixtures, transposition, and every relational composition just defined.
The notation $\Phi\in\mathcal U$ may equally use the eventual-set
record of the pair net, under the explicit convention preceding this
subsection. The structure is globally diagonal when the saturation of
$\{\Delta_Z\}$ belongs to $\mathcal U$, and it is separated when
membership of a constant pair $(z,w)$ implies $z=w$.
Its point convergence and Cauchy nets are specified by
\[
 z_\alpha\to z
 \iff (z_\alpha,z)_\alpha\in\mathcal U,
 \qquad
 (z_\alpha)\text{ is Cauchy}
 \iff (z_\alpha,z_\beta)_{(\alpha,\beta)}\in\mathcal U.
\]
A member of $\mathcal U$ is called an originally uniformly null pair
net. The canonical group structure is globally diagonal, as verified
in Proposition~\ref{lt:canonical-group} from the group convergence
axioms.

\begin{definition}[Faithful lattice representation]\label{lt:faithful}
Let $(Z,\mathcal U)$ be a separated uniform convergence space. For each
$j\in J$, let $F_j$ be an abelian lattice-ordered group and let
$d_j:Z^2\to F_j^+$ satisfy
\begin{equation}\label{lt:pseudometric}
 d_j(z,z)=0,\qquad d_j(z,w)=d_j(w,z),\qquad
 d_j(z,v)\le d_j(z,w)+d_j(w,v).
\end{equation}
Let $\chi_j$ be a specified null convergence on $F_j$. The family
$(d_j,\chi_j)_{j\in J}$ is faithful for $\mathcal U$ if, for every
admissible pair net,
\begin{equation}\label{fd:faithful-pair-net}
 (z_\alpha,w_\alpha)_\alpha\in\mathcal U
 \iff d_j(z_\alpha,w_\alpha)\longrightarrow_{\chi_j}0
 \quad(j\in J).
\end{equation}
For an arbitrary pair net $\Phi=(z_\alpha,w_\alpha)_\alpha$, this
can also be written
\begin{equation}\label{lt:pair-faithfulness}
 \Phi\in\mathcal U
 \quad\Longleftrightarrow\quad
 d_j\Phi:=(d_j(z_\alpha,w_\alpha))_\alpha
 \xrightarrow{\chi_j}0\quad(j\in J).
\end{equation}
An \emph{order representation} has $\chi_j=o$ for every $j$.
The regulator representations used later have $\chi_j=D$.
\end{definition}

By quantifying over all pair nets,
\eqref{lt:pair-faithfulness} includes point convergence, Cauchyness and
parameter-uniform convergence within the same representation. A
representation formulated only for $d_j(z_\alpha,z)$ accounts for point
convergence, and its extension to \eqref{lt:pair-faithfulness} therefore
requires a separate assertion concerning arbitrary pair nets.

For a nonempty parameter set $\mathcal A$, uniform convergence of
$z_\alpha^A$ to $z^A$ means original uniform nullity of the pair net
\[
 (z_\alpha^A,z^A)_{(\alpha,A)\in I\times\mathcal A},
 \qquad
 (\alpha,A)\preceq(\beta,B)\iff\alpha\preceq\beta.
\]
Because every tail contains every parameter, the order controls in
Proposition~\ref{lt:three-translations} bound all parameters on the same
tail, without requiring a supremum of their values. That proposition
establishes the exact translations for point, Cauchy and parameter-uniform
convergence.

For an abelian convergence group $(Z,\lambda)$, we use its original
canonical pair convergence
\[
 \delta(z,w)=z-w,\qquad
 \mathcal U_\lambda
 =\{(z_\alpha,w_\alpha)_\alpha:(z_\alpha-w_\alpha)_\alpha\to_\lambda0\}.
\]
Proposition~\ref{lt:canonical-group} proves the uniform convergence axioms,
global diagonal stability and recovery of $\lambda$. It also establishes
that a faithful null representation $\rho_j$ yields the faithful invariant
pair metrics $d_j(z,w)=\rho_j(z-w)$.

For these gauges, symmetry and subadditivity give the useful inequality
\begin{equation}\label{lt:reverse-triangle}
 |\rho_j(z)-\rho_j(w)|\le\rho_j(z-w).
\end{equation}
Indeed,
\[
 \rho_j(z)-\rho_j(w)\le\rho_j(z-w),\qquad
 \rho_j(w)-\rho_j(z)\le\rho_j(w-z)=\rho_j(z-w).
\]
Thus order convergence of $\rho_j(z_\alpha-z)$ gives order convergence
of $\rho_j(z_\alpha)-\rho_j(z)$. The corresponding assertion for
regulator convergence uses its solid null convergence.

\subsubsection{Exact values and lattice realizations}

To place these representations in their wider context, the appendix
constructs the required value objects, beginning with every $T_1$ net
convergence in Theorem~\ref{amp:point-metrization} and passing to
separated globally diagonal uniform convergence structures in
Corollary~\ref{amu:uniform-metrization}. For the canonical group
structure, Theorem~\ref{amg:group-metrization} gives an invariant
metric, and in each construction the tail values of the admissible
net determine its vanishing radius.

To obtain faithful lattice representations, additional algebraic
properties enter through Proposition~\ref{lt:algebraic-resolution},
which derives such a representation from preservation of joins,
downward-directed infima and subadditivity in a compatible exact value model.
Since these properties are not asserted for the universal incidence
monoids, passage to lattice carriers and stability under limits are
established separately in Propositions~\ref{lt:resolved-carrier} and
\ref{lt:carrier-completeness-sandwich}, while
Proposition~\ref{fd:intrinsic-order-realizations} derives the order and
unbounded-order realizations directly from lattice-group operations.
This development also contains the bounded product argument and the
complete-carrier criterion, with weak distributivity, completeness and
solidity subject throughout to their explicitly stated hypotheses.

For a faithful invariant representation in a partially ordered abelian
group, Theorem~\ref{nx:canonical-po} gives an exact existence criterion
through a closure condition computed from the original convergence.
Here the representing group need not have binary meets, so the monoid
operations used in order interchange are identified separately in
Lemma~\ref{mg:calculus}; the intrinsic properties \ref{nx:CB} and
\ref{nx:P}, together with the functional equivalence in
Theorem~\ref{nx:uniform-limit}, then provide a second route to the
limit argument. Throughout the partially ordered group constructions,
the infimum of every order regulator is zero in the whole group, with
its lower-bound condition also tested against elements incomparable
with zero.

\section{The abstract net Riemann integral}\label{sectionnetriemannintegral}

Having fixed the language in which limits are to be taken, we turn to
the finite expressions whose limits define the integral. A coefficient
product determines their algebra, while an indexing net describes the
successive approximations; together they extend the construction of
\cite{szaznetintegralconvergence,fleischerchange}, with the integral and
its elementary properties expressed throughout in the original
convergence structures of the preceding section.

\begin{definition}[Coefficient product]\label{bilinearproduct}
Let $(X,\lambda_X)$, $(Y,\lambda_Y)$ and $(Z,\lambda_Z)$ be commutative
convergence monoids in Definition~\ref{fd:convergence-monoid}. A coefficient product is a map
\[
 \bullet:X\times Y\longrightarrow Z
\]
whose partial maps are continuous zero-preserving monoid homomorphisms.
Thus, for all admissible variables,
\[
\begin{gathered}
(x+x')\bullet y=x\bullet y+x'\bullet y,\qquad
x\bullet(y+y')=x\bullet y+x\bullet y',\\
0\bullet y=0,\qquad x\bullet0=0.
\end{gathered}
\]
The zero identities are explicit. Additivity alone does not imply them
in a noncancellative monoid. Abelian convergence groups are included.
If the three monoids are positive cones, or the three groups are lattice
ordered, positivity means
\[
 x\in X_+,\ y\in Y_+\quad\Longrightarrow\quad x\bullet y\in Z_+.
\]
A coefficient product between modules over a common ring is bilinear when it
also respects that scalar multiplication. Scalar compatibility is used only
in assertions about the corresponding module operations.
\end{definition}

\begin{definition}[Nets for integration and integration structures]
\label{integrationstructure}
Let $\mathcal H\subseteq\mathcal P(\Omega)$ be a paving. For each
$A\in\mathcal H$, an integration net consists of a directed set $\Gamma_A$
and finite records
\[
 R^A_\gamma=((C^A_{\gamma r},t^A_{\gamma r}))_{r\in I^A_\gamma},
 \qquad C^A_{\gamma r}\in\mathcal H,\quad t^A_{\gamma r}\in\Omega,
 \qquad\gamma\in\Gamma_A.
\]
The index family $I^A_\gamma$ is finite, and repetitions are retained. Additional
conditions on cells, tags, disjointness and refinement are properties of the
specified integration nets. An integration structure consists of this family
of nets and a coefficient product with separated target $Z$.

The symbol $R\sqcup Q$ denotes concatenation of finite records, with disjoint
index copies. For an empty record every sum below is $0$. If
$\varnothing\in\mathcal H$, its integration net is the constant empty record.
\end{definition}

\begin{definition}[Net Riemann integral]\label{netriemannpartiallyordered}
For $f:\Omega\to X$ and $\mu:\mathcal H\to Y$, put
\[
 S^A_\gamma(f,\mu)
   =\sum_{r\in I^A_\gamma}
       f(t^A_{\gamma r})\bullet\mu(C^A_{\gamma r}).
\]
The function $f$ is net Riemann integrable on $A$ if there is $J_f(A)\in Z$
with
\[
 S^A_\gamma(f,\mu)\xrightarrow[\gamma\in\Gamma_A]{\lambda_Z}J_f(A).
\]
Separation makes this element unique, and we write
\[
 J_f(A)=\int_A f\,d\mu.
\]
Integrability on every $A\in\mathcal H$ is called hereditary integrability.
For a fixed integration structure and a fixed $\mu$, its class is denoted by
$L_{NR}(\mu)$.
\end{definition}

\begin{definition}[Abstract kernel integral]
\label{abstractnetriemannpartiallyordered}
Let $Z$ be any separated commutative convergence monoid. For
$U:\Omega\times\mathcal H\to Z$, define
\[
 S^{U,A}_\gamma=\sum_{r\in I^A_\gamma}
 U(t^A_{\gamma r},C^A_{\gamma r}).
\]
If $S^{U,A}_\gamma\to_{\lambda_Z}J_U(A)$, its abstract net Riemann integral
is $\int_AU=J_U(A)$. In particular,
\[
 U_f(t,C)=f(t)\bullet\mu(C)
 \quad\Longrightarrow\quad
 S^{U_f,A}_\gamma=S^A_\gamma(f,\mu).
\]
\end{definition}

At this stage every record contributes only a finite sum, so its
behaviour under products and continuous homomorphisms follows directly
from the corresponding algebra. Countable decompositions will enter
later, through the approximation structure and the order properties of
the range.

\begin{proposition}[Products and continuous homomorphisms]
\label{nr:functorial}
Let $(Z_l)_{l\in L}$ be separated commutative convergence monoids and let
$Z=\prod_{l\in L}Z_l$ carry product convergence. For
$U:\Omega\times\mathcal H\to Z$ and $z\in Z$,
\[
 S^{U,A}_\gamma\longrightarrow z
 \quad\Longleftrightarrow\quad
 S^{\pi_lU,A}_\gamma\longrightarrow\pi_lz\quad(l\in L).
\]
Consequently $U$ is integrable on $A$ exactly when every coordinate is
integrable there.

If $T:Z\to Z'$ is a continuous zero-preserving homomorphism between
separated commutative convergence monoids, then
\[
 U\text{ integrable on }A
 \quad\Longrightarrow\quad
 \int_A(T\circ U)=T\!\left(\int_AU\right).
\]
\end{proposition}
\begin{proof}
For every $l$ and every $\gamma$,
\[
 \pi_lS^{U,A}_\gamma
 =\sum_{r\in I^A_\gamma}\pi_lU(t^A_{\gamma r},C^A_{\gamma r})
 =S^{\pi_lU,A}_\gamma.
\]
Product convergence gives the first equivalence. If the coordinate integrals
are $z_l$, the element $(z_l)_{l\in L}$ belongs to $\prod_lZ_l$ and is the
integral of $U$.
For the second assertion,
\[
 S^{T\circ U,A}_\gamma
 =\sum_{r\in I^A_\gamma}T(U(t^A_{\gamma r},C^A_{\gamma r}))
 =T(S^{U,A}_\gamma)
 \longrightarrow T(J_U(A)).
\]
\end{proof}

\begin{corollary}[Additivity in the integrand]
\label{linearityofnetriemannintegral}
The class $L_{NR}(\mu)$ is a submonoid of $X^\Omega$, and
\[
 J_{f+g}(A)=J_f(A)+J_g(A),\qquad J_0(A)=0.
\]
If $X,Z$ are abelian convergence groups, this class is a subgroup and
$J_{-f}(A)=-J_f(A)$. In the module setting of
Definition~\ref{bilinearproduct}, if multiplication by every fixed scalar
is continuous on $Z$, it is a submodule and
\[
 J_{af+bg}(A)=aJ_f(A)+bJ_g(A).
\]
\end{corollary}
\begin{proof}
On the common directed set $\Gamma_A$, the coefficient identities give
\[
 S^A_\gamma(f+g,\mu)=S^A_\gamma(f,\mu)+S^A_\gamma(g,\mu),
 \qquad S^A_\gamma(0,\mu)=0.
\]
Joint continuity of addition and separation yield
\[
 S^A_\gamma(f+g,\mu)\longrightarrow J_f(A)+J_g(A),
 \qquad J_{f+g}(A)=J_f(A)+J_g(A).
\]
For groups, the identities and continuity required for inverses are
\[
 S^A_\gamma(-f,\mu)=-S^A_\gamma(f,\mu)
 \longrightarrow-J_f(A).
\]
For modules they are
\[
 S^A_\gamma(af+bg,\mu)=aS^A_\gamma(f,\mu)+bS^A_\gamma(g,\mu)
 \longrightarrow aJ_f(A)+bJ_g(A).
\]
These statements hold on every $A\in\mathcal H$.
\end{proof}

\begin{proposition}[Order properties]\label{propertiesnetriemannintegral}
Suppose that $X,Y,Z$ are abelian lattice-ordered groups, that the coefficient
product is positive, that $\mu(C)\ge0$ for every $C$, and that $\lambda_Z$
is separated and locally solid. For integrable functions on $A$,
\[
 f\ge0\Longrightarrow J_f(A)\ge0,
 \qquad f\le g\Longrightarrow J_f(A)\le J_g(A).
\]
If $f$ and $|f|$ are both integrable on $A$, then
\[
 |J_f(A)|\le J_{|f|}(A).
\]
\end{proposition}
\begin{proof}
The positive cone is stable under original convergence. Indeed,
$z_\alpha\ge0$ and $z_\alpha\to z$ give
\[
 0\le z^-\le|z_\alpha-z|.
\]
Local solidity makes the constant net $z^-$ converge to $0$. Separation
implies $z^-=0$. Thus $z\ge0$.
If $f\ge0$, positivity gives $S^A_\gamma(f,\mu)\ge0$, hence $J_f(A)\ge0$.
If $f\le g$, Corollary~\ref{linearityofnetriemannintegral} and the first
assertion give
\[
 0\le J_{g-f}(A)=J_g(A)-J_f(A).
\]
Finally, $-|f|\le f\le|f|$ gives
\[
 -J_{|f|}(A)\le J_f(A)\le J_{|f|}(A),
\]
which is equivalent to the asserted modulus inequality.
\end{proof}

To pass from integration on a set to integration on its subsets, we
must also examine how the corresponding records can be related, since
the existence of the first integral alone does not determine the
others. The conditions below express the cancellation of a common
complement directly in terms of those records, before an integrand
is fixed.

\begin{definition}[Tail pasting]\label{nr:tail-pasting}
Let $A,D,B\in\mathcal H$, where $A\cap D=\varnothing$ and $B=A\cup D$.
The integration nets admit tail pasting for $(A,D)$ if both conditions hold.
\begin{enumerate}
\item For $C=A$ and for $C=D$, with $C'=B\setminus C$, and every
$\beta_0,\beta'_0\in\Gamma_B$, there are
$\alpha_0,\alpha'_0\in\Gamma_C$ such that
\[
 \begin{gathered}
 \alpha\ge\alpha_0,\quad\alpha'\ge\alpha'_0
 \\
 \Longrightarrow\quad
 \exists\delta\in\Gamma_{C'}\ \exists\beta\ge\beta_0\
 \exists\beta'\ge\beta'_0:\qquad
 R^B_\beta=R^C_\alpha\sqcup R^{C'}_\delta,\quad
 R^B_{\beta'}=R^C_{\alpha'}\sqcup R^{C'}_\delta.
 \end{gathered}
\]
Record equality here and below permits a permutation of its finite indices.
\item For every $\beta_0\in\Gamma_B$, there are
$\alpha_0\in\Gamma_A$ and $\delta_0\in\Gamma_D$ such that
\[
 \alpha\ge\alpha_0,\quad\delta\ge\delta_0
 \quad\Longrightarrow\quad
 \exists\beta\ge\beta_0:\quad
 R^B_\beta=R^A_\alpha\sqcup R^D_\delta.
\]
\end{enumerate}
\end{definition}

\begin{theorem}[Subset integration]\label{subsetintegration}
Suppose $Z$ is complete for its original group convergence and the integration
nets admit tail pasting for $(A,D)$, where $B=A\sqcup D$. If a kernel $U$ is
integrable on $B$, then it is integrable on $A$ and on $D$, and
\[
 J_U(B)=J_U(A)+J_U(D).
\]
\end{theorem}
\begin{proof}
Write $x^C_\alpha=S^{U,C}_\alpha$. Since
$x^B_\beta\to J_U(B)$, product reindexing and subtraction give
\[
 x^B_\beta-x^B_{\beta'}\longrightarrow0
 \qquad((\beta,\beta')\in\Gamma_B\times\Gamma_B).
\]
For $C\in\{A,D\}$ and each pair of $B$-tail indices, the first tail-pasting
condition gives a pair of $C$-tail indices. Its record identities imply
\[
 \begin{aligned}
 x^C_\alpha-x^C_{\alpha'}
 &=(x^C_\alpha+x^{C'}_\delta)
     -(x^C_{\alpha'}+x^{C'}_\delta)\\
 &=x^B_\beta-x^B_{\beta'}.
 \end{aligned}
\]
Therefore
\[
 \begin{split}
 \{x^C_\alpha-x^C_{\alpha'}:
       \alpha\ge\alpha_0,\ \alpha'\ge\alpha'_0\}
 \subseteq
 \{x^B_\beta-x^B_{\beta'}:
       \beta\ge\beta_0,\ \beta'\ge\beta'_0\}.
 \end{split}
\]
This inclusion makes the $C$-difference net a quasi-subnet of the
$B$-difference net, hence null by
Definition~\ref{netconvergencestructure}. Consequently $x^A$ and $x^D$
are Cauchy, and completeness supplies their limits $J_U(A)$ and $J_U(D)$.

On $\Gamma_A\times\Gamma_D$, joint addition gives
\[
 y_{\alpha,\delta}=x^A_\alpha+x^D_\delta
 \longrightarrow J_U(A)+J_U(D).
\]
The second tail-pasting condition gives
\[
 \{y_{\alpha,\delta}:\alpha\ge\alpha_0,\delta\ge\delta_0\}
 \subseteq\{x^B_\beta:\beta\ge\beta_0\}
\]
for every $\beta_0$. Thus $y$ is a quasi-subnet of $x^B$ and also converges to
$J_U(B)$. Separation proves the stated identity.
\end{proof}

\begin{proposition}[Subset integration with invariant monoid distances]
\label{nr:monoid-subset}
Let $Z$ be a separated commutative convergence monoid whose original
uniform convergence is complete. Suppose symmetric distances
$d_j:Z^2\to M_j$ faithfully represent point and pair convergence by
order convergence, and satisfy
\begin{equation}\label{nr:invariant-laws}
\begin{gathered}
 d_j(z,z)=0,\qquad d_j(z,w)\le d_j(z,v)+d_j(v,w),\\
 d_j(z+c,w+c)=d_j(z,w).
\end{gathered}
\end{equation}
The controllers are partially ordered commutative monoids. Order nullity
means eventual domination by every member of a decreasing control net
with infimum zero. For a positive cone in a partially ordered group,
these infima are taken in the entire group. No completeness, finite-meet
property, or continuity of controller addition is assumed here.
If the records satisfy
Definition~\ref{nr:tail-pasting} and a kernel $U$ is integrable on
$B=A\sqcup D$, then it is integrable on $A,D$ and
\[
 J_U(B)=J_U(A)+J_U(D).
\]
\end{proposition}
\begin{proof}
Put $x^C_\gamma=S^{U,C}_\gamma$ and $b=J_U(B)$. Original
convergence and symmetry give the uniformly null pair nets
\[
 (x^B_\beta,b)_\beta\in\mathcal U_Z,
 \qquad (b,x^B_{\beta'})_{\beta'}\in\mathcal U_Z.
\]
Their composed tail relations satisfy
\[
 (T_{\beta_0}(x^B)\times\{b\})
 \circ(\{b\}\times T_{\beta'_0}(x^B))
 =T_{\beta_0}(x^B)\times T_{\beta'_0}(x^B).
\]
Thus the composed saturation is tail equivalent to
$(x^B_\beta,x^B_{\beta'})_{(\beta,\beta')}$, which is uniformly null
by the composition axiom. Pair faithfulness gives
\[
 d_j(x^B_\beta,x^B_{\beta'})\longrightarrow_o0.
\]
Fix $j$, $C\in\{A,D\}$, and an attained bound $u$ from a decreasing
control for this pair net. There are $\beta_0,\beta'_0$ such that
\[
 \beta\ge\beta_0,\ \beta'\ge\beta'_0
 \quad\Longrightarrow\quad d_j(x^B_\beta,x^B_{\beta'})\le u.
\]
The first tail-pasting clause supplies $\alpha_0,\alpha'_0$. For every
$\alpha\ge\alpha_0$, $\alpha'\ge\alpha'_0$, its completion identities give
indices $\delta,\beta,\beta'$ with
\[
 x^B_\beta=x^C_\alpha+x^{C'}_\delta,\qquad
 x^B_{\beta'}=x^C_{\alpha'}+x^{C'}_\delta,
 \qquad\beta\ge\beta_0,\ \beta'\ge\beta'_0.
\]
Consequently
\[
 d_j(x^C_\alpha,x^C_{\alpha'})
 =d_j(x^B_\beta,x^B_{\beta'})\le u.
\]
Every control member is therefore attained by the $C$-pair net.
Inverse pair faithfulness makes that net Cauchy in the original uniform
convergence. Completeness supplies $J_U(A)$ and $J_U(D)$.

Joint continuity gives
\[
 x^A_\alpha+x^D_\delta\longrightarrow J_U(A)+J_U(D).
\]
For each $B$-tail, the second tail-pasting clause provides a product tail
whose sum image is contained in that $B$-tail image. The sum net is a
quasi-subnet of $x^B$ and hence also converges to $J_U(B)$. Separation
proves the identity.
\end{proof}

\begin{remark}[Cancellation and the scope of the range extension]
\label{nr:monoid-cancellation}
The invariant faithful family forces cancellation in $Z$:
\[
 x+c=y+c\Longrightarrow
 d_j(x,y)=d_j(x+c,y+c)=0\ (j)\Longrightarrow x=y.
\]
For the final implication the constant pair has every represented
distance zero, hence converges to the diagonal. Separation applies.
No inverse element is constructed or used. In contrast,
Theorem~\ref{subsetintegration} needs no ordered representation of the
original group convergence. The finite-sum assertions of
Proposition~\ref{nr:functorial} and
Corollary~\ref{linearityofnetriemannintegral} also apply to
noncancellative commutative monoids.
\end{remark}

\begin{corollary}[Finite additivity]\label{additivitynetriemannintegral}
Let $Z$ be a separated commutative convergence monoid. Let $\mathcal H$
be a ring of sets and suppose the second clause of
Definition~\ref{nr:tail-pasting} holds for every disjoint pair of its members. If $U$ is hereditarily integrable, then
$J_U:\mathcal H\to Z$ is finitely additive. Completeness is unnecessary for
this assertion. Neither inverse elements nor an invariant representation
are needed.
\end{corollary}
\begin{proof}
All three integrals already exist. On $\Gamma_A\times\Gamma_D$,
\[
 S^{U,A}_\alpha+S^{U,D}_\delta\longrightarrow J_U(A)+J_U(D).
\]
The second tail-pasting condition makes this sum net a quasi-subnet of
$(S^{U,A\sqcup D}_\beta)_\beta$, which converges to $J_U(A\sqcup D)$.
Separation gives
\[
 J_U(A\sqcup D)=J_U(A)+J_U(D).
\]
The empty-record convention gives $J_U(\varnothing)=0$. Induction on the
number of disjoint members gives finite additivity.
\end{proof}

\begin{proposition}[A geometric sufficient condition for tail pasting]
\label{nr:product-pasting}
Suppose records over $B=A\sqcup D$ are indexed by
$\Gamma_B=\Gamma_A\times\Gamma_D$ with the product order and satisfy
\[
 R^B_{(\alpha,\delta)}=R^A_\alpha\sqcup R^D_\delta.
\]
Then tail pasting holds. The assertion remains valid if each of the three
record nets is replaced by a tail-equivalent record net.
\end{proposition}
\begin{proof}
Let $\beta_0=(\alpha_0,\delta_0)$ and
$\beta'_0=(\alpha'_0,\delta'_0)$. Directedness supplies
$\delta\ge\delta_0,\delta'_0$. For
$\alpha\ge\alpha_0$ and $\alpha'\ge\alpha'_0$,
\[
 (\alpha,\delta)\ge\beta_0,\qquad
 (\alpha',\delta)\ge\beta'_0.
\]
These indices give the common-complement identities for $C=A$.
For $C=D$, let $\alpha_*\ge\alpha_0,\alpha'_0$. For every
$\delta\ge\delta_0$ and $\delta'\ge\delta'_0$,
\[
 (\alpha_*,\delta)\ge\beta_0,\qquad
 (\alpha_*,\delta')\ge\beta'_0,
\]
and the two records have the common complement $R^A_{\alpha_*}$.
For the second condition, $\alpha\ge\alpha_0$ and
$\delta\ge\delta_0$ give $(\alpha,\delta)\ge\beta_0$ and
$R^B_{(\alpha,\delta)}=R^A_\alpha\sqcup R^D_\delta$.
For a net of records write
\[
 \mathcal T_C(\gamma_0)=\{R^C_\gamma:\gamma\ge\gamma_0\},
 \qquad
 \mathcal U\sqcup\mathcal V
 =\{R\sqcup Q:R\in\mathcal U,\ Q\in\mathcal V\}.
\]
Let $\widetilde{\mathcal T}_B$ and
$\widetilde{\mathcal T}'_B$ be tails of a replacement $B$-net.
Tail equivalence supplies indices with
\[
 \begin{aligned}
 \mathcal T_A(\alpha_0)\sqcup\mathcal T_D(\delta_0)
    &\subseteq\widetilde{\mathcal T}_B,\\
 \mathcal T_A(\alpha'_0)\sqcup\mathcal T_D(\delta'_0)
    &\subseteq\widetilde{\mathcal T}'_B.
 \end{aligned}
\]
For some $\delta_*\ge\delta_0,\delta'_0$, tail equivalence also gives
replacement tails satisfying
\[
 \begin{gathered}
 \widetilde{\mathcal T}_A\subseteq\mathcal T_A(\alpha_0),\qquad
 \widetilde{\mathcal T}'_A\subseteq\mathcal T_A(\alpha'_0),\\
 \varnothing\ne\widetilde{\mathcal T}_D
              \subseteq\mathcal T_D(\delta_*).
 \end{gathered}
\]
Consequently, for every
$(R,R',Q)\in\widetilde{\mathcal T}_A\times
 \widetilde{\mathcal T}'_A\times\widetilde{\mathcal T}_D$,
\[
 R\sqcup Q\in\widetilde{\mathcal T}_B,\qquad
 R'\sqcup Q\in\widetilde{\mathcal T}'_B.
\]
This proves the first condition for $C=A$. For $C=D$, let
$\alpha_*\ge\alpha_0,\alpha'_0$ and take replacement tails with
\[
 \begin{gathered}
 \varnothing\ne\widetilde{\mathcal T}_A
                   \subseteq\mathcal T_A(\alpha_*),\\
 \widetilde{\mathcal T}_D\subseteq\mathcal T_D(\delta_0),\qquad
 \widetilde{\mathcal T}'_D\subseteq\mathcal T_D(\delta'_0).
 \end{gathered}
\]
For $R\in\widetilde{\mathcal T}_A$,
$Q\in\widetilde{\mathcal T}_D$ and
$Q'\in\widetilde{\mathcal T}'_D$, the original product-tail
inclusions give
\[
 R\sqcup Q\in\widetilde{\mathcal T}_B,\qquad
 R\sqcup Q'\in\widetilde{\mathcal T}'_B.
\]
The record $R$ is the common complement, proving the other case.
For the second condition, a prescribed replacement $B$-tail contains
an original product tail, and each coordinate tail contains a
replacement tail. Hence
\[
 \widetilde{\mathcal T}_A\sqcup\widetilde{\mathcal T}_D
 \subseteq\mathcal T_A\sqcup\mathcal T_D
 \subseteq\widetilde{\mathcal T}_B.
\]
\end{proof}

Having related integration on disjoint sets, we next ask when the
records recover the prescribed masses of those sets. For indicators,
this amounts to a relation between tags and coefficient masses that can
be specified directly in the given measure data in $Y$, before the
corresponding integral in $Z$ has been established.

\begin{definition}[Regular coefficient masses]\label{nr:regular-integrator}
Suppose $A,B,A\cap B\in\mathcal H$. The coefficient mass net is
\[
 m^A_\gamma(B)
   =\sum_{\substack{r\in I^A_\gamma\\t^A_{\gamma r}\in B}}
       \mu(C^A_{\gamma r})\in Y.
\]
The integrator is regular at $(A,B)$ if
\[
 m^A_\gamma(B)\longrightarrow_{\lambda_Y}\mu(A\cap B).
\]
If $\mathcal H$ is closed under finite intersections, regularity means
regularity at every $(A,B)\in\mathcal H^2$.
\end{definition}

\begin{definition}[Exact indicator resolution]\label{nr:indicator-resolution}
The integration structure resolves the indicator of $B$ on $A$ for $\mu$
if, for some $\gamma_0\in\Gamma_A$,
\[
 m^A_\gamma(B)=\mu(A\cap B)\qquad(\gamma\ge\gamma_0).
\]
Exact indicator resolution means this condition for every admissible pair
$(A,B)$. It implies regularity: every one of the indicated mass nets is
eventually constant at its required value.
\end{definition}

\begin{proposition}[Finite partitions verify indicator resolution]
\label{nr:indicator-partitions}
Let $\mathcal H$ be a ring containing $\Omega$ and let $\mu$ be finitely
additive. Suppose that, for every $A,B\in\mathcal H$, all sufficiently late
records on $A$ consist of a finite partition of $A$ whose cells lie either
in $B$ or in $\Omega\setminus B$, with each tag belonging to its cell.
Then exact indicator resolution holds.
\end{proposition}
\begin{proof}
For every such record,
\[
 t^A_{\gamma r}\in B\quad\Longleftrightarrow\quad
 C^A_{\gamma r}\subseteq B,
 \qquad
 A\cap B=\bigsqcup_{t^A_{\gamma r}\in B}C^A_{\gamma r}.
\]
Finite additivity gives the defining identity.
\end{proof}

\begin{proposition}[The finite-partition construction]
\label{nr:canonical-partitions}
Let $\mathcal H$ be a ring of subsets of $\Omega$ containing $\Omega$.
For $A\in\mathcal H$, let $\Gamma_A$ consist of all tagged finite
partitions of $A$ into nonempty members of $\mathcal H$. Order records
by refinement of their partitions, independently of their tags. The
integration record at $\gamma$ is $R^A_\gamma=\gamma$.
For $A=\varnothing$, use the singleton empty record. Then:
\begin{enumerate}[(i)]
\item each $\Gamma_A$ is directed,
\item tail pasting holds for every $B=A\sqcup D$ in $\mathcal H$, and
\item every finitely additive $\mu:\mathcal H\to Y$ has exact indicator
      resolution on these nets.
\end{enumerate}
These assertions depend only on finite partitions and the coefficient
masses.
\end{proposition}
\begin{proof}
For a finite partition $\pi$ into nonempty cells, its set of possible
tag assignments is
\[
 \operatorname{Tag}(\pi)=\prod_{C\in\pi}C\ne\varnothing.
\]
Thus every such partition occurs as the underlying partition of a
record. Given partitions $\pi,\sigma$ of $A$, put
\[
 \pi\wedge\sigma
 =\{C\cap E:C\in\pi,\ E\in\sigma,\ C\cap E\ne\varnothing\}.
\]
This is a finite partition of $A$ in $\mathcal H$, refining both
$\pi$ and $\sigma$. Every element of
$\operatorname{Tag}(\pi\wedge\sigma)$ produces a common upper bound
of the two original records. This proves~(i).

For~(ii), let $B=A\sqcup D$ and let $\beta_0,\beta'_0\in\Gamma_B$
have partitions $\pi,\pi'$. For $C\in\{A,D\}$ and $C'=B\setminus C$,
put
\[
 \begin{aligned}
 \sigma_C&=\{C\cap E\cap E':E\in\pi,\ E'\in\pi',\
                         C\cap E\cap E'\ne\varnothing\},\\
 \sigma_{C'}&=\{C'\cap E\cap E':E\in\pi,\ E'\in\pi',\
                         C'\cap E\cap E'\ne\varnothing\}.
 \end{aligned}
\]
For every record $\alpha_0$ on $\sigma_C$, every record $\delta$ on
$\sigma_{C'}$, and every $\alpha,\alpha'\ge\alpha_0$, the records
\[
 \beta=R^C_\alpha\sqcup R^{C'}_\delta,
 \qquad
 \beta'=R^C_{\alpha'}\sqcup R^{C'}_\delta
\]
belong to $\Gamma_B$. Their partitions refine $\pi$ and $\pi'$, so
\[
 \beta\ge\beta_0,\qquad\beta'\ge\beta'_0.
\]
The same complement record $R^{C'}_\delta$ occurs in both identities.
This verifies the first tail-pasting condition for each $C\in\{A,D\}$.

For the second condition, let $\beta_0$ have partition $\pi$. Define
\[
 \pi|_A=\{A\cap E:E\in\pi,\ A\cap E\ne\varnothing\},\qquad
 \pi|_D=\{D\cap E:E\in\pi,\ D\cap E\ne\varnothing\}.
\]
For records $\alpha_0,\delta_0$ on these two partitions,
\[
 \alpha\ge\alpha_0,\quad\delta\ge\delta_0
 \quad\Longrightarrow\quad
 R^A_\alpha\sqcup R^D_\delta\in\Gamma_B,
 \qquad R^A_\alpha\sqcup R^D_\delta\ge\beta_0.
\]
This proves~(ii).

For~(iii), fix $A,B\in\mathcal H$ and form the partition
\[
 \pi(A,B)=\{A\cap B,A\setminus B\}\setminus\{\varnothing\}.
\]
Every record beyond a tagged version of $\pi(A,B)$ has each cell
contained in $B$ or in its complement. Since its tags belong to their
cells,
\[
 A\cap B=\bigsqcup_{t^A_{\gamma r}\in B}C^A_{\gamma r},
 \qquad
 m^A_\gamma(B)
 =\sum_{t^A_{\gamma r}\in B}\mu(C^A_{\gamma r})
 =\mu(A\cap B).
\]
These are finite-additivity identities, valid on an entire record tail.
They prove exact indicator resolution.
\end{proof}

\begin{proposition}[Integration of simple functions]
\label{integrationofsimplefunctions}
Let $X,Y,Z$ be commutative monoids with separately zero-preserving
additive coefficient product, and give $Y,Z$ convergence structures, with
$Z$ separated and addition jointly continuous. Only continuity of
$y\mapsto x_i\bullet y$ is required below. For $x\in X$ and
$B\subseteq\Omega$, put
\[
 x\mathbf1_B(t)=
 \begin{cases}x,&t\in B,\\0,&t\notin B.\end{cases}
\]
Let $f=\sum_{i=1}^m x_i\mathbf1_{B_i}$, where $x_i\in X$ and
$B_i\in\mathcal H$. If $\mu$ is regular at $(A,B_i)$ for each $i$, then
$f$ is integrable on $A$ and
\[
 J_f(A)=\sum_{i=1}^m x_i\bullet\mu(A\cap B_i).
\]
In particular, regularity on $\mathcal H^2$ makes all these simple functions
hereditarily integrable. Disjointness of the $B_i$ is unnecessary.
\end{proposition}
\begin{proof}
Biadditivity and finite distributivity give, for every $\gamma\in\Gamma_A$,
\[
 \begin{aligned}
 S^A_\gamma(f,\mu)
 &=\sum_{r\in I^A_\gamma}\sum_{i=1}^m
       (x_i\mathbf1_{B_i}(t^A_{\gamma r}))\bullet\mu(C^A_{\gamma r})\\
 &=\sum_{i=1}^m x_i\bullet
       \sum_{t^A_{\gamma r}\in B_i}\mu(C^A_{\gamma r})\\
 &=\sum_{i=1}^m x_i\bullet m^A_\gamma(B_i).
 \end{aligned}
\]
For every $i$, regularity and continuity of the partial coefficient map give
\[
 x_i\bullet m^A_\gamma(B_i)
 \longrightarrow_{\lambda_Z}x_i\bullet\mu(A\cap B_i).
\]
Joint continuity of finite addition gives the stated limit.
Under exact indicator resolution, directedness supplies a common upper
bound for the finitely many corresponding tail indices. The displayed sum
is then eventually equal to the asserted integral.
\end{proof}

\begin{corollary}[Simple functions on partition nets]
\label{nr:partition-simple-functions}
Under Proposition~\ref{nr:canonical-partitions}, every finitely additive
$\mu:\mathcal H\to Y$ and every
$f=\sum_{i=1}^m x_i\mathbf1_{B_i}$ satisfy
\[
 \int_A f\,d\mu=\sum_{i=1}^m x_i\bullet\mu(A\cap B_i)
 \qquad(A\in\mathcal H).
\]
The finite-sum net is eventually equal to this value. Neither
completeness of $Z$, continuity of coefficient maps, nor an assumed
convergence of coefficient masses is needed.
\end{corollary}
\begin{proof}
Let $\pi$ be the finite partition consisting of the nonempty sets
\[
 A\cap D_1\cap\cdots\cap D_m,
 \qquad D_i\in\{B_i,\Omega\setminus B_i\}.
\]
For every record beyond a tagged version of $\pi$, the finite-mass
identities in Proposition~\ref{nr:canonical-partitions} hold
simultaneously for $B_1,\ldots,B_m$. Hence
\[
 S^A_\gamma(f,\mu)
 =\sum_{i=1}^m x_i\bullet m^A_\gamma(B_i)
 =\sum_{i=1}^m x_i\bullet\mu(A\cap B_i).
\]
An eventually constant net is a quasi-subnet of the corresponding
constant net. The original convergence axioms give the displayed
integral.
\end{proof}

Before passing from finite sums to limits of integrands, it is useful
to separate the two coefficient properties contained in preservation of
positive countable suprema. Their roles in the convergence theorems are
different, and the following equivalence allows each later argument to
invoke precisely the property it needs.

\begin{proposition}[Order properties of coefficient maps]
\label{nr:coefficient-order}\label{nr:coefficient-split}
Let $L:X\to Z$ be positive and additive between abelian lattice-ordered
groups. The identity
\[
 L\!\left(\bigvee_nx_n\right)=\bigvee_nLx_n
 \qquad(x_n\ge0, \bigvee_nx_n\text{ exists in }X)
\]
is equivalent to the conjunction
\[
 L|x|=|Lx|\quad(x\in X),\qquad
 x_n\downarrow0\Longrightarrow Lx_n\downarrow0.
\]
No completeness of either group is assumed.
\end{proposition}
\begin{proof}
Suppose the displayed countable suprema are preserved. The alternating
sequence $x^+,x^-,x^+,x^-,\ldots$ has supremum $|x|$. Writing
$a=Lx^+$ and $b=Lx^-$ gives
\[
 L|x|=a\vee b=a+b.
\]
The identity $a+b=(a\vee b)+(a\wedge b)$ gives $a\wedge b=0$, and
\[
 |Lx|=|a-b|=a+b=L|x|.
\]
If $x_n\downarrow0$, then $x_1-x_n\uparrow x_1$, whence
\[
 Lx_1=\bigvee_nL(x_1-x_n)=\bigvee_n(Lx_1-Lx_n),
 \qquad\bigwedge_nLx_n=0.
\]
Positivity makes $(Lx_n)$ decreasing.

Conversely, suppose both coefficient properties hold. Modulus
preservation gives preservation of positive parts: the identity
$2x^+=x+|x|$ implies
\[
 2Lx^+=Lx+L|x|=Lx+|Lx|=2(Lx)^+,
 \qquad Lx^+=(Lx)^+.
\]
Here doubling is injective in an abelian lattice-ordered group. Therefore
\[
 L(x\vee y)=L\bigl(y+(x-y)^+\bigr)
 =Ly+(Lx-Ly)^+=Lx\vee Ly.
\]
Let $x_n\ge0$ and $x=\bigvee_nx_n$ exist, and put
$q_N=\bigvee_{n\le N}x_n$. Then
\[
 x-q_N\downarrow0,\qquad Lx-Lq_N\downarrow0,
 \qquad Lq_N=\bigvee_{n\le N}Lx_n.
\]
The inequalities $x_n\le x$ give $Lx_n\le Lx$ for every $n$.
If $z\ge Lx_n$
for every $n$, then $z\ge Lq_N$ for every $N$, and
$z\ge\sup_NLq_N=Lx$. Thus $Lx=\sup_nLx_n$.
\end{proof}

\begin{remark}[Sequential normality does not imply modulus preservation]
\label{nr:coefficient-strictness}
Let $H\ne\{0\}$ be an abelian lattice-ordered group and define
$L:H\times H\to H$ by $L(a,b)=a+b$. The map is positive and additive.
If $(a_n,b_n)\downarrow(0,0)$, then $a_n,b_n\ge0$ and both coordinates
decrease to zero. For a lower bound $c$ of the product sum net,
\[
 c-b_m\le a_n\ (n,m)\Longrightarrow c\le b_m\ (m)
 \Longrightarrow c\le0.
\]
Thus $(a_n+b_m)_{(n,m)}\downarrow0$. For each $(n,m)$,
$N=\max\{n,m\}$ satisfies $a_N+b_N\le a_n+b_m$. Consequently
$a_n+b_n\downarrow0$. For $u>0$ and
$x=(u,-u)$,
\[
 |Lx|=0,\qquad L|x|=2u>0.
\]
Thus the sequentially order-continuous coefficient class strictly contains
the class satisfying both properties of
Proposition~\ref{nr:coefficient-order}.
\end{remark}

These elementary properties bring us to the question that finite
additivity alone does not settle: under which conditions does uniform
approximation of integrands survive the passage from their record sums
to their integrals? We shall express that approximation through the
original pair convergence, first identifying the obstruction and then
the conditions that permit construction of the candidate integral and,
in the respective limit theorems, the required bounds on its primitive.

\section{Limit interchange and uniform convergence of the integrals}
\label{uc:section}\label{sec:uniform-limits}

With the elementary integral identities established, we turn to the
comparison of two limiting operations, beginning where their relation
can already be expressed by finite sums. Biadditivity writes the
difference of two record sums through the pointwise differences of
their integrands, so uniform approximation supplies a finite comparison
before an integral of the limiting function has been constructed.
Whether this comparison reaches the known integrals depends on how
the original convergence treats its error sets, and the example below
shows that completeness and local convexity do not by themselves
preserve that comparison.
Following the errors through the two limits will lead us to intrinsic
stability operations, whose sufficiency is proved at the cardinality
of the error family and whose necessity is recovered from sampling
systems within the same finite-record framework. We then return to
the ordered and regulator hypotheses of the subsequent comparisons
to verify those operations directly from their bounds.

\subsection{A uniform approximation obstruction}
\label{ob:section}

Uniform approximation leaves a specific passage to be justified: a
small set of pointwise errors may control the finite sums, but those
estimates must still survive the limits that define the integrals.
To examine this passage, we retain completeness and local convexity in
the range while arranging for the individual integrals to exist, so
that the obstruction can be located in an operation on the error sets.
Throughout the construction the measure on $[0,1]$ is ordinary Lebesgue
measure, and every convergence assertion concerning the vector values
refers to the convergence defined below.

\begin{definition}[Consistency of finite sampling]
\label{ob:consistency}
A finite-sampling procedure assigns to a curve $g:[0,1]\to B$ a directed
family of sums
\[
 S_\gamma g=\sum_{r=1}^{m_\gamma}a_{\gamma r}g(t_{\gamma r}),
 \qquad m_\gamma<\infty,\quad
 a_{\gamma r}\in\mathbb R,\quad t_{\gamma r}\in[0,1].
\]
It is norm-consistent on continuous curves if, for every real Banach
space $B$ and every norm-continuous $g:[0,1]\to B$, this net converges in
norm to the ordinary Riemann integral of $g$.  Integrability in a
prescribed convergence means convergence of the same sums in that
convergence.  The record family may depend on the curve.
\end{definition}

By stating consistency independently of the curves to be constructed,
Definition~\ref{ob:consistency} leaves the choice of partitions, tags and
other record parameters with the procedure. In the proof, only its
application to those particular continuous curves will be needed.

\begin{theorem}[A common obstruction for finite sampling]
\label{ob:counterexample}
There are a real vector space $E$ with a complete separated vector net
convergence $\lambda$ and continuous maps $f_n:[0,1]\to E$ such that
every null net in $E$ is a quasi-subnet of a null net whose tail-value
sets are absorbing, balanced and convex, and
\[
 f_n\longrightarrow0\quad\hbox{uniformly on }[0,1],
 \qquad
 \bigcup_{n\ge1}f_n([0,1])\quad\hbox{is bounded in }E.
\]
For every procedure in Definition~\ref{ob:consistency}, each $f_n$ is
integrable and the zero function is integrable, whereas
\[
 I f_n\not\longrightarrow_\lambda I0=0.
\]
Uniformity and completeness are understood in the original canonical
difference convergence, in which a pair net $(u_i,v_i)$ is uniformly
null exactly when $u_i-v_i\to_\lambda0$.
\end{theorem}

\begin{proof}
Let
\[
 E=\ell^1(\mathbb N\times\mathbb N),\qquad
 M=c_{00}(\mathbb N\times\mathbb N),\qquad
 y_n=\frac1n\sum_{k=1}^{\infty}2^{-2k-1}e_{n,k},
\]
where $e_{n,k}$ is the coordinate vector at $(n,k)$.  The series
converges in $E$ and
\begin{equation}\label{ob:y-norm}
 \|y_n\|_1=\frac1n\sum_{k=1}^{\infty}2^{-2k-1}
 =\frac1{6n}.
\end{equation}
For any distinct $n_1,\ldots,n_m$, a linear combination
$\sum_{j=1}^m c_jy_{n_j}$ with $c_j\ne0$ has infinitely many nonzero
coordinates in its $n_j$-th row.  Consequently its membership in $M$
forces every $c_j$ to vanish, and the cosets $y_n+M$ are linearly
independent in $E/M$.  We can therefore define a linear map on
$M+\operatorname{span}\{y_n:n\ge1\}$ by
\[
 q|_M=0,\qquad q(y_n)=n,
\]
and extend it algebraically to a linear map $q:E\to\mathbb R$.

For every net $(x_i)_{i\in I}$, prescribe
\begin{equation}\label{ob:net-convergence}
 x_i\longrightarrow_\lambda x
 \quad\Longleftrightarrow\quad
 \|x_i-x\|_1\longrightarrow0
 \quad\hbox{and}\quad
 (\exists i_0\in I)(\exists R<\infty)
 (\forall i\ge i_0)\ |q(x_i-x)|\le R.
\end{equation}
Constant nets satisfy this condition, and passage to a quasi-subnet
preserves both of its clauses.  If $x_i\to_\lambda x$ and
$z_i\to_\lambda x$ on the same directed set, and
$w_i\in\{x_i,z_i\}$ for every $i$, then
\[
 \|w_i-x\|_1\le\max\{\|x_i-x\|_1,\|z_i-x\|_1\},\qquad
 |q(w_i-x)|\le\max\{|q(x_i-x)|,|q(z_i-x)|\}.
\]
On a common tail the second bound is finite, while the first tends to
zero.  These estimates verify the finite-mixing axiom of the original
net convergence.  Uniqueness of norm limits proves separation.

For $x_i\to_\lambda x$ and $z_j\to_\lambda z$, the inequalities
\[
 \|(x_i+z_j)-(x+z)\|_1
 \le\|x_i-x\|_1+\|z_j-z\|_1,
\]
\[
 |q((x_i+z_j)-(x+z))|
 \le |q(x_i-x)|+|q(z_j-z)|
\]
prove continuity of addition on the product directed set, and replacing
$x_i-x$ by its negative proves continuity of inversion.  If
$r_j\to r$ in $\mathbb R$, then $(r_j)$ is bounded on a tail and
\[
 \|r_jx_i-rx\|_1
 \le |r_j|\,\|x_i-x\|_1+|r_j-r|\,\|x\|_1,
\]
\[
 q(r_jx_i-rx)
 =r_jq(x_i-x)+(r_j-r)q(x).
\]
On the product directed set the first right-hand side tends to zero,
while the second is bounded on a tail, proving joint continuity of
scalar multiplication and making $E$ a convergence vector space.

For $\varepsilon,R>0$, put
\begin{equation}\label{ob:disks}
 D_{\varepsilon,R}
 =\{u\in E:\|u\|_1\le\varepsilon, |q(u)|\le R\}.
\end{equation}
Each set is balanced and convex by the triangle inequalities and the
linearity of $q$.  It is absorbing because, for $u\in E$ and
\[
 c>\max\{\|u\|_1/\varepsilon,|q(u)|/R\},
\]
one has $u\in cD_{\varepsilon,R}$.  Fix $R>0$ and index the net
$v_{\varepsilon,u}=u$ by all pairs
$(\varepsilon,u)$ with $\varepsilon>0$ and $u\in D_{\varepsilon,R}$,
ordered by
\[
 (\varepsilon,u)\preceq(\eta,v)
 \quad\Longleftrightarrow\quad \eta\le\varepsilon.
\]
Its tail-value set at $(\varepsilon,u)$ is exactly
$D_{\varepsilon,R}$, so it converges to zero by
\eqref{ob:net-convergence}.  Given any $x_i\to_\lambda0$, there are
$i_0$ and $R>0$ such that $|q(x_i)|\le R$ for $i\ge i_0$.
For each $\varepsilon>0$, a common upper index for $i_0$ and a norm
threshold places the corresponding tail of $(x_i)$ inside
$D_{\varepsilon,R}$.  Hence $(x_i)$ is a quasi-subnet of this null
net of absorbing convex disks, as asserted.

To prove completeness, let $(x_i)_{i\in I}$ be Cauchy for the
canonical difference convergence, so that
$(x_i-x_j)_{(i,j)\in I\times I}\to_\lambda0$.  It is norm Cauchy.
For each integer $n\ge1$ there is $i_n'\in I$ with
\[
 i,j\ge i_n'\quad\Longrightarrow\quad\|x_i-x_j\|_1\le2^{-n}.
\]
Select indices $i_n$ recursively with $i_1\ge i_1'$ and
$i_{n+1}\ge i_n,i_{n+1}'$.  The norm Cauchy sequence $(x_{i_n})$
has a norm limit $x\in E$.  For $i\ge i_n$ and $m\ge n$,
\[
 \|x_i-x_{i_m}\|_1\le2^{-n},
\]
and passage to the norm limit in $m$ gives
$\|x_i-x\|_1\le2^{-n}$.  Thus the entire original net converges in
norm to $x$.  Its difference convergence also supplies $i_*,j_*\in I$
and $R<\infty$ such that
\[
 i\ge i_*,\quad j\ge j_*
 \quad\Longrightarrow\quad |q(x_i-x_j)|\le R.
\]
An upper bound $k$ of $i_*,j_*$ satisfies, for every $i\ge k$,
\[
 |q(x_i-x)|
 \le |q(x_i-x_k)|+|q(x_k-x)|
 \le R+|q(x_k-x)|.
\]
Equation~\eqref{ob:net-convergence} therefore gives
$x_i\to_\lambda x$, proving completeness without a restriction on
the original directed set.

With completeness and local convexity established, we construct curves
whose disjoint coordinate directions keep every finite sample in $M$.
For $k\ge1$, set
\[
 I_k=[1-2^{-(k-1)},1-2^{-k}],
\]
and let $\varphi_k$ be the triangular function on $I_k$ that vanishes
at both endpoints and has height $2^{-k}$ at its midpoint.  Its area is
\[
 \int_{I_k}\varphi_k(t)\,dt
 =\tfrac12\,2^{-k}\,2^{-k}=2^{-2k-1}.
\]
Define
\begin{equation}\label{ob:curves}
 f_n(t)=\frac1n\varphi_k(t)e_{n,k}\quad(t\in I_k),
 \qquad f_n(1)=0.
\end{equation}
Since the values agree at adjacent endpoints, continuity can be checked
through an estimate valid across the intervals. Write $I_k=[a_k,b_k]$;
each $\varphi_k$ has slopes $2$ and $-2$, and hence
\[
 |\varphi_k(s)-\varphi_k(t)|\le2|s-t|,
 \qquad
 0\le\varphi_k(t)\le2\min\{t-a_k,b_k-t\}.
\]
For $s\in I_k$, $t\in I_l$ with $k<l$, one has $b_k\le a_l$ and
\[
 \|f_n(t)-f_n(s)\|_1
 =\frac{\varphi_k(s)+\varphi_l(t)}n
 \le\frac2n\bigl((b_k-s)+(t-a_l)\bigr)
 \le\frac2n(t-s).
\]
For $s,t$ in one interval, the slope bound gives the same estimate.
For $t=1$, the bound
$\varphi_k(s)\le2(b_k-s)\le2(1-s)$ gives it as well. Therefore
\begin{equation}\label{ob:lipschitz-curves}
 \|f_n(t)-f_n(s)\|_1\le\frac2n|t-s|
 \qquad(s,t\in[0,1]),
\end{equation}
and $f_n$ is norm-continuous on the entire interval. All its values belong to $M$, and
\begin{equation}\label{ob:curve-bounds}
 q(f_n(t))=0,
 \qquad \|f_n(t)\|_1\le\frac1{2n}
 \quad(n\ge1, t\in[0,1]).
\end{equation}
If $t_i\to t$, then $f_n(t_i)-f_n(t)$ converges to zero in norm and
its $q$-value is zero at every index.  This proves continuity into the
prescribed convergence.

Order $\mathbb N\times[0,1]$ by its first coordinate.  The net
$(f_n(t))_{(n,t)}$ tends to zero by \eqref{ob:curve-bounds}, which is
the required uniform convergence for the canonical difference
structure.  Write $B=\bigcup_{n\ge1}f_n([0,1])$.  For every scalar
net $r_i\to0$, the net $(r_i b)_{(i,b)\in I\times B}$, ordered by
the $i$ coordinate, satisfies
\[
 \|r_i b\|_1\le\tfrac12|r_i|,\qquad q(r_i b)=0.
\]
It is null in $\lambda$, proving the stated boundedness of $B$.

To identify the ordinary norm Riemann integral of $f_n$ as $y_n$,
truncate \eqref{ob:curves} after the first $K$ hats. The resulting
continuous function $f_n^{(K)}$ has integral
\[
 y_n^{(K)}=\frac1n\sum_{k=1}^{K}2^{-2k-1}e_{n,k}
\]
and uniform error at most $2^{-K-1}/n$.  Every ordinary positive
partition sum of $f_n-f_n^{(K)}$ has norm at most this same number,
because its coefficients have sum one.  Therefore, for every such
record $\gamma$,
\[
 \|S_\gamma f_n-y_n\|_1
 \le \frac{2^{-K-1}}n
     +\|S_\gamma f_n^{(K)}-y_n^{(K)}\|_1
     +\|y_n^{(K)}-y_n\|_1.
\]
In this estimate the last term is explicitly
\[
 \|y_n-y_n^{(K)}\|_1
 =\frac1n\sum_{k>K}2^{-2k-1}=\frac{2^{-2K}}{6n}.
\]
For completeness, the truncated function is also $2/n$-Lipschitz:
the preceding within-hat and between-hat bounds remain valid when
the later hats are zero. Its integral is the finite sum $y_n^{(K)}$,
computed coordinatewise from the displayed triangular areas. If
$\gamma$ is a tagged partition of mesh at most $\delta$, integration
of this Lipschitz bound on each cell gives
\[
 \|S_\gamma f_n^{(K)}-y_n^{(K)}\|_1
 \le\sum_{[a,b]\in\gamma}\frac2n\delta(b-a)
 =\frac{2\delta}n.
\]
Consequently, for every $K\ge1$ and every such partition,
\[
 \|S_\gamma f_n-y_n\|_1
 \le\frac{2^{-K-1}}n+\frac{2\delta}n
                         +\frac{2^{-2K}}{6n}.
\]
Given $\varepsilon>0$, an integer $K$ with
$2^{-K-1}/n+2^{-2K}/(6n)<\varepsilon/2$ and a positive
$\delta<n\varepsilon/4$ make this bound less than $\varepsilon$.
A finite partition of mesh below $\delta$ exists, and every refinement
has mesh below $\delta$. This verifies norm convergence of the
ordinary partition sums to $y_n$ on their original refinement net.

For a procedure satisfying Definition~\ref{ob:consistency}, every
finite sum $S_\gamma f_n$ belongs to $M$.  Hence norm consistency
and the defining values of $q$ give
\[
 \|S_\gamma f_n-y_n\|_1\longrightarrow0,
 \qquad q(S_\gamma f_n-y_n)=-n.
\]
Because the second expression is bounded for each fixed $n$,
\eqref{ob:net-convergence} gives $S_\gamma f_n\to_\lambda y_n$;
yet the integral values satisfy \eqref{ob:y-norm} and $q(y_n)=n$,
so they do not converge to zero in $\lambda$. Every finite sum of the
zero function is zero, proving its integrability and completing the
construction.
\end{proof}

What permits each row to converge is precisely what fails to control
the collection of endpoints: the algebraic map $q$ allows a bounded
error for each fixed row, although the endpoint errors are unbounded
as the row varies. Here $q$ specifies which norm-convergent nets the
target admits, while the integral values come directly from the
vector-valued sums; no continuous dual or scalar representation of
an integral is involved.

\begin{proposition}[Persistence for partial Saks defects]
\label{ob:saks}
For the curves in Theorem~\ref{ob:counterexample}, let
$J_n(t)=\int_0^t f_n(s)\,ds$ be the ordinary norm primitive.  For a
finite family $Q$ of nondegenerate tagged intervals
$([a,b],t)$ with pairwise disjoint interiors and $t\in[a,b]$, and
numbers $\theta_{a,b,t}$ with $|\theta_{a,b,t}|\le1$, put
\[
 R_{n,Q,\theta}
 =\sum_{([a,b],t)\in Q}\theta_{a,b,t}
       \bigl((b-a)f_n(t)-J_n(b)+J_n(a)\bigr).
\]
For every fixed $n$, these defects converge to zero in $\lambda$
uniformly over $Q$ and $\theta$ as the maximum interval length tends
to zero, whereas $J_n(1)\not\to_\lambda0$ as $n\to\infty$.
\end{proposition}
\begin{proof}
For $t<1$, the interval $[0,t]$ meets only finitely many nonzero
hats, so $J_n(t)\in M$ and $q(J_n(t))=0$.  At $t=1$,
$J_n(1)=y_n$ and $q(J_n(1))=n$.  Define
\[
 \omega_n(\delta)
 =\sup\{\|f_n(s)-f_n(t)\|_1:s,t\in[0,1],\ |s-t|\le\delta\}.
\]
Equation~\eqref{ob:lipschitz-curves} gives
$0\le\omega_n(\delta)\le2\delta/n$, and hence
$\omega_n(\delta)\to0$ as $\delta\downarrow0$.  For $b-a\le\delta$
and $t\in[a,b]$, positive partition sums on $[a,b]$, followed by
their norm limit, give
\[
 \|(b-a)f_n(t)-J_n(b)+J_n(a)\|_1
 \le(b-a)\omega_n(\delta).
\]
Since the intervals of $Q$ have total length at most one, summing
these bounds gives
\[
 \|R_{n,Q,\theta}\|_1\le\omega_n(\delta)
 \qquad\bigl(\max_{([a,b],t)\in Q}(b-a)\le\delta\bigr).
\]
At most one interval in $Q$ has right endpoint $1$, and no
nondegenerate interval has left endpoint $1$.  Since $q(f_n(t))=0$ and $q(J_n(t))=0$ for $t<1$, applying $q$
to the finite expression gives
\[
 q(R_{n,Q,\theta})
 =-n\sum_{([a,1],t)\in Q}\theta_{a,1,t},
 \qquad |q(R_{n,Q,\theta})|\le n,
\]
where the sum has at most one term.
For $\delta>0$, let $\mathcal D_{n,\delta}$ consist of zero and
all these defects with maximum interval length at most $\delta$.
The net indexed by $(\delta,r)$, $r\in\mathcal D_{n,\delta}$,
ordered by decreasing $\delta$, has exactly
$\mathcal D_{n,\delta}$ as its tail-value sets.  The two displayed
bounds prove that this net is null by
\eqref{ob:net-convergence}.  The endpoint assertion is
$J_n(1)=y_n$ and \eqref{ob:y-norm} together with $q(y_n)=n$.
\end{proof}

Having seen that the obstruction persists for partial Saks defects,
we can locate it more precisely by following the null error families
through the limiting operation. Two properties emerge from this
calculation, which we state in net form before giving their intrinsic
definitions. For $A\subseteq E$, let $a_E(A)$ be the set of limits of
nets all of whose values belong to $A$.

\begin{proposition}[Adherence failure and absorption of independent rows]
\label{ob:diagnosis}
In the convergence of Theorem~\ref{ob:counterexample}, the decreasing
sets
\[
 A_N=\{u\in M:\|u\|_1\le N^{-1}\}
\]
have a null saturation net, but the saturation of their one-step
adherences is not null.  In contrast, let $(A_\eta)_{\eta\in\Lambda}$
be any decreasing family of nonempty sets whose saturation tends to
zero, and let $(z_i)_{i\in I}$ be a net.  If, for every $\eta$,
there is a net $(h_j^\eta)_{j\in J_\eta}\to_\lambda0$ such that
\begin{equation}\label{ob:absorption-premise}
 (\forall j\in J_\eta)(\exists i_{\eta,j}\in I)
 (\forall i\ge i_{\eta,j})\quad
 z_i\in A_\eta+\{h_k^\eta:k\ge j\},
\end{equation}
then $z_i\to_\lambda0$.
\end{proposition}
\begin{proof}
For a decreasing nonempty family $(B_\eta)$, its saturation is the
net $b_{\eta,u}=u$, indexed by $u\in B_\eta$ and ordered by the
$\eta$ coordinate.  Its tail-value set at $(\eta,u)$ is exactly
$B_\eta$, since every later set is contained in $B_\eta$ and every
point of $B_\eta$ occurs at the same $\eta$ level.  In particular,
the saturation of $(A_N)$ is null because its values have norm at
most $N^{-1}$ and have $q$-value zero.

Every limit in $\lambda$ is a norm limit, and therefore
$a_E(A_N)\subseteq\{x:\|x\|_1\le N^{-1}\}$.  Conversely, if
$\|x\|_1\le N^{-1}$, let $x^{(m)}$ be its restriction to
$\{1,\ldots,m\}\times\{1,\ldots,m\}$.  Then
\[
 x^{(m)}\in A_N,\qquad
 \|x^{(m)}-x\|_1\longrightarrow0,
 \qquad q(x^{(m)}-x)=-q(x).
\]
The sequence converges to $x$ in $\lambda$, giving
\begin{equation}\label{ob:adherence-balls}
 a_E(A_N)=\{x\in E:\|x\|_1\le N^{-1}\}.
\end{equation}
For every $N$ and every $L>0$, an integer
$n>\max\{N/6,L\}$ satisfies $y_n\in a_E(A_N)$ and
$q(y_n)=n>L$.  Each tail-value set of the saturation of
$(a_E(A_N))$ thus has unbounded $q$-image, so this net fails the
second clause of \eqref{ob:net-convergence}.

Now suppose \eqref{ob:absorption-premise} holds.  Given
$\varepsilon>0$, nullity of the outer saturation supplies an index
$\eta_\varepsilon$ with
\[
 \|a\|_1\le\varepsilon/2\qquad(a\in A_{\eta_\varepsilon}).
\]
The inner null net at that index has a threshold $j_\varepsilon$
such that
$\|h_k^{\eta_\varepsilon}\|_1\le\varepsilon/2$ for
$k\ge j_\varepsilon$.  Equation~\eqref{ob:absorption-premise}
then places a tail of $(z_i)$ inside the norm ball of radius
$\varepsilon$.  This proves $\|z_i\|_1\to0$.

To verify the remaining clause of the convergence, outer nullity
supplies one index $\eta_*$ and a finite $R$ such that $|q(a)|\le R$
for all $a\in A_{\eta_*}$, while the null net at that fixed index
supplies $j_*$ and a finite $S$ with $|q(h_k^{\eta_*})|\le S$ for
$k\ge j_*$. On the tail given by \eqref{ob:absorption-premise},
\[
 |q(z_i)|\le R+S.
\]
Together with norm convergence to zero, this proves
$z_i\to_\lambda0$.
\end{proof}

\begin{proposition}[Finite responses within the obstruction]
\label{ob:finite-response}
For the same range, ordinary positive partition coefficients have
continuous finite responses.  More generally, let $\Sigma(B)$ consist
of zero and all finite words
\[
 \sum_{r=1}^m a_r x_r,
 \qquad a_r\ge0,\quad \sum_{r=1}^m a_r\le1,\quad x_r\in B.
\]
If $(B_\eta)$ is decreasing and its saturation is null in $\lambda$,
then the saturation of $(\Sigma(B_\eta))$ is also null in $\lambda$.
\end{proposition}
\begin{proof}
For every $\varepsilon>0$, a threshold $\eta_\varepsilon$ gives
$\|x\|_1\le\varepsilon$ for all $x\in B_{\eta_\varepsilon}$.
Each response from that set satisfies
\[
 \left\|\sum_{r=1}^m a_r x_r\right\|_1
 \le\sum_{r=1}^m a_r\|x_r\|_1\le\varepsilon.
\]
There are also one index $\eta_*$ and a finite $R$ such that
$|q(x)|\le R$ for all $x\in B_{\eta_*}$.  Every corresponding
response satisfies
\[
 \left|q\left(\sum_{r=1}^m a_r x_r\right)\right|
 \le\sum_{r=1}^m a_r|q(x_r)|\le R.
\]
The zero response satisfies both inequalities.  Since the response
sets decrease, these are the two tail conditions in
\eqref{ob:net-convergence} for their saturation.
\end{proof}

Taken together, these calculations locate the failure in the
adherence passage of \eqref{ob:adherence-balls}, even though the finite
coefficient response is continuous and independent inner errors can be
absorbed. Its scope extends beyond positive partitions, because all
finite samples of the constructed curves lie in $M$ and norm
consistency fixes their limits; any additional requirement restoring
uniform interchange in every such range must therefore exclude at
least one of these curves or change its prescribed value. We are thus
led to distinguish the passages that preserve estimates from those
that construct a limiting integral, and to formulate the corresponding
intrinsic conditions before verifying them from the abstract
assumptions already available for the integration theorems.

\subsection{Adherence and absorption in the original convergence}
\label{uc:intrinsic-mechanism}

To understand the preceding obstruction, we return to the estimate on finite
records and ask what remains of its error set after the two sum nets have
passed to their known integrals. Their limiting difference belongs to the
adherence of that set, so this first passage already requires a stability
property of the original convergence; constructing an integral that is not
yet known to exist will require a further passage through independently
convergent rows. We shall keep these two operations separate throughout
the argument, writing $\langle B_\lambda:\lambda\in\Lambda\rangle$
for the saturated net of a decreasing family of nonempty sets, with the
meaning fixed in Section~\ref{fd:net-foundations}. Nullity refers to this
net, while domination by the family means eventual membership in each of
its sets, and neither convention imposes a countability condition on the
saturated index, whose cardinality may exceed that of $\Lambda$.

\begin{definition}[Right pair convergence and cardinal completeness]
\label{ci:right-structure}
Let $G$ be a separated convergence group with identity $e$. Multiplication
and inversion are continuous for its prescribed convergence, and the group
need not be abelian. A pair net $(x_i,y_i)$ is right uniformly null when
\[
 \delta_r(x_i,y_i):=x_i y_i^{-1}\longrightarrow e.
\]
A net $(g_i)_{i\in I}$ is right Cauchy when
\[
 (g_i g_j^{-1})_{(i,j)\in I\times I}\longrightarrow e.
\]
For an infinite cardinal $\kappa$, $G$ is right $\kappa$-complete if
every right Cauchy net with a displayed directed index of cardinality at
most $\kappa$ converges. Right sequential completeness has the
corresponding meaning for sequences. For maps $f_p,f:S\to G$, right
uniform convergence means nullity of
\[
 \left\langle
 \{f_p(s)f(s)^{-1}:p\ge p_0,\ s\in S\}:p_0\in P
 \right\rangle.
\]
Uniform convergence on the empty set is automatic. Cardinality refers to
the displayed index before a net is replaced by a tail-equivalent
representative.
\end{definition}

\begin{lemma}[The original convergence and countable completeness]
\label{ci:countable-completeness}
The right pair convergence is a globally diagonal uniform convergence
structure, expressed through pair nets, and induces the original point
convergence.
Right $\aleph_0$-completeness is equivalent to right sequential completeness.
\end{lemma}
\begin{proof}
For every $x,y,z\in G$,
\[
 \delta_r(x,x)=e,\qquad
 \delta_r(y,x)=\delta_r(x,y)^{-1},\qquad
 \delta_r(x,z)=\delta_r(x,y)\delta_r(y,z).
\]
A net in $\Delta_G$ therefore has constant difference $e$. If a
pair net $(x_i,y_i)$ has null difference, continuity of inversion
gives $\delta_r(y_i,x_i)\to e$. For the composition axiom, let
$(R_\alpha)_{\alpha\in H}$ and $(T_\beta)_{\beta\in K}$ be
decreasing nonempty families of pair sets with null difference
saturations, and suppose $R_\alpha\circ T_\beta\ne\varnothing$
for all $(\alpha,\beta)$. The product of their two difference
saturations has tail-value sets
\[
 \delta_r(R_\alpha)\delta_r(T_\beta),
 \qquad(\alpha,\beta)\in H\times K,
\]
and converges to $e$ by continuity of multiplication. If
$(x,z)\in R_\alpha\circ T_\beta$, a witnessing $y$ satisfies
$(x,y)\in R_\alpha$ and $(y,z)\in T_\beta$, whence
\[
 \delta_r(x,z)=\delta_r(x,y)\delta_r(y,z)
 \in\delta_r(R_\alpha)\delta_r(T_\beta).
\]
Thus every tail of that product contains the corresponding tail of
$\langle\delta_r(R_\alpha\circ T_\beta)\rangle$, proving its
nullity by the quasi-subnet axiom.

If a pair net $v$ is a quasi-subnet of a pair net $u$, then for
each tail $T_i(u)$ there is a tail $T_k(v)\subseteq T_i(u)$,
and consequently
\[
 T_k(\delta_r v)=\delta_r(T_k(v))
 \subseteq\delta_r(T_i(u))=T_i(\delta_r u).
\]
This proves stability under quasi-subnets. For a finite mixing
$w_i\in\{u_i,v_i\}$ one has
$\delta_r(w_i)\in\{\delta_r(u_i),\delta_r(v_i)\}$, so the
finite-mixing axiom for the point convergence proves the corresponding
pair axiom. If $R_i=T_i(u)$ are the tails of a right uniformly null
pair net, the image of $R_i\cup\Delta_G$ under $\delta_r$ is
$\delta_r(R_i)\cup\{e\}$. Proposition~\ref{fd:tail-calculus}
represents these sets as the tails of the mixture of $\delta_r u$
with the constant net $e$. That mixture is null, which proves global
diagonality. Finally, continuity of $v\mapsto vx^{-1}$ and
$v\mapsto vx$ gives
\[
 (x_i,x)\longrightarrow\Delta_G
 \iff x_i x^{-1}\longrightarrow e
 \iff x_i\longrightarrow x.
\]

Right $\aleph_0$-completeness applies to the directed index
$\mathbb N$, and hence implies right sequential completeness.
Conversely, let $I$ be countable, let $(g_i)_{i\in I}$ be right
Cauchy, and enumerate $I$ as $(a_n)_{n\ge1}$, repeating indices if
$I$ is finite. Directedness permits the recursion
\[
 i_1\ge a_1,\qquad i_{n+1}\ge i_n,a_{n+1}.
\]
If $i=a_m$, then $i_n\ge i$ for every $n\ge m$. Thus
$n\mapsto i_n$ is increasing and cofinal. For any $(i,j)\in I^2$,
indices $m,n$ with $i_m\ge i$ and $i_n\ge j$ show that
$(m,n)\mapsto(i_m,i_n)$ is cofinal into $I^2$. Hence
$(g_{i_m}g_{i_n}^{-1})_{m,n}\to e$, and sequential completeness
supplies $g\in G$ with $g_{i_n}\to g$.

On $I\times\mathbb N$, the map $(i,n)\mapsto(i,i_n)$ is
increasing and cofinal into $I^2$: above a prescribed $(i_0,j_0)$,
it suffices to take $i\ge i_0$ and $n\ge n_0$ with
$i_{n_0}\ge j_0$. Therefore
\[
 (g_i g_{i_n}^{-1})_{(i,n)}\longrightarrow e,\qquad
 (g_{i_n})_{(i,n)}\longrightarrow g,
\]
and continuity of multiplication gives
\[
 (g_i g_{i_n}^{-1})g_{i_n}=g_i\longrightarrow g
 \quad\text{on }I\times\mathbb N.
\]
Since the tail-value set above $(i_0,n_0)$ is exactly
$\{g_i:i\ge i_0\}$, this product-index net is tail equivalent to
the original net, and hence $g_i\to g$ on $I$.
\end{proof}

\begin{definition}[One-step adherence stability]
\label{ci:adherence}
For a subset $B$ of a convergence space $X$, put
\[
 a_X(B)=\{x\in X:\text{some net with every value in }B
                         \text{ converges to }x\},\qquad
 a_X^2(B)=a_X(a_X(B)).
\]
For an infinite cardinal $\kappa$, a convergence group $G$ satisfies
\textnormal{$(R_\kappa)$} when every decreasing nonempty family
$(B_\lambda)_{\lambda\in\Lambda}$ satisfies
\begin{equation}\label{ci:R}
 |\Lambda|\le\kappa,\qquad
 \langle B_\lambda\rangle\longrightarrow e
 \quad\Longrightarrow\quad
 \langle a_G(B_\lambda)\rangle\longrightarrow e.
\end{equation}
The unrestricted condition is denoted by \textnormal{(R)}. We write
\textnormal{$(R_\omega)$} for $\kappa=\aleph_0$. One-step adherence
is not assumed idempotent.
\end{definition}

In this form, adherence stability belongs to the convergence-theoretic
regularity calculus studied by Wilde and Kent \cite{WildeKent1999}, with
the cardinal restriction recording how many error levels the approximation
provides. It leaves the adhering nets arbitrary, as they must be when
their limits are integrals on prescribed integration nets, and can be
compared with the stronger requirement that the entire convergence-closed
hull of each tail remain small.

\begin{definition}[Closed-base stability]\label{nx:CB}
A subset $C$ of a convergence space $X$ is closed when every net in
$C$ that converges in $X$ has its limit in $C$. Let $K_X(A)$ be the
intersection of all closed supersets of $A$. The space satisfies
\textnormal{(CB)} when every net $u=(u_i)\to x$ satisfies
\begin{equation}\label{nx:CB-equation}
 \langle K_X(T_i(u)):i\in I\rangle\longrightarrow x.
\end{equation}
The displayed family is downward directed and consists of nonempty
sets, so its saturation is defined.
\end{definition}

\begin{lemma}[Adherence calculus]\label{ci:adherence-calculus}
The operation $a_X$ is extensive and monotone, and
$a_X(\varnothing)=\varnothing$. For a continuous map $T:X\to Y$,
\begin{equation}\label{ci:adherence-image}
 T(a_X(B))\subseteq a_Y(T(B)).
\end{equation}
In a convergence group,
\begin{align}
 a_G(B)^{-1}&=a_G(B^{-1}),\label{ci:adherence-inverse}\\
 u\,a_G(B)\,v&=a_G(uBv),\label{ci:adherence-translate}\\
 a_G(B)a_G(C)&\subseteq a_G(BC).\label{ci:adherence-product}
\end{align}
If $K_X(B)$ is the intersection of the convergence-closed supersets of
$B$, then
\[
 B\subseteq a_X(B)\subseteq a_X^2(B)\subseteq K_X(B).
\]
Consequently closed-base stability \textnormal{(CB)} implies
\textnormal{(R)}. Under \textnormal{$(R_\kappa)$}, both
$\langle a_G(B_\lambda)\rangle$ and
$\langle a_G^2(B_\lambda)\rangle$ are null whenever
$\langle B_\lambda\rangle$ is null and $|\Lambda|\le\kappa$.
\end{lemma}
\begin{proof}
Constant nets prove extensivity. A net with all its values in $B$ also
has all its values in $C$ when $B\subseteq C$, proving monotonicity.
A net has a nonempty index and cannot have every value in the empty set.
For $x\in a_X(B)$, let $b_i\in B$ converge to $x$. Continuity gives
$T(b_i)\to T(x)$ with $T(b_i)\in T(B)$, proving
\eqref{ci:adherence-image}. Apply \eqref{ci:adherence-image} to inversion to obtain
$a_G(B)^{-1}\subseteq a_G(B^{-1})$, and apply it again with
$B^{-1}$ in place of $B$ to obtain
$a_G(B^{-1})^{-1}\subseteq a_G(B)$. Inverting the latter inclusion
gives the reverse inclusion. For $T(x)=uxv$, both $T$ and
$T^{-1}(x)=u^{-1}xv^{-1}$ are continuous, so
\[
 T(a_G(B))\subseteq a_G(T(B)),\qquad
 T^{-1}(a_G(T(B)))\subseteq a_G(B).
\]
Applying $T$ to the second inclusion proves equality. If $b_i\in B$ and
$c_j\in C$ converge to $b$ and $c$, respectively, then
$b_i c_j\in BC$ on the product index and $b_i c_j\to bc$, proving
\eqref{ci:adherence-product}.

If $C$ is closed and $B\subseteq C$, then every net in $B$
converging to $x$ is a net in $C$, and $x\in C$. Intersecting over
these sets gives $a_X(B)\subseteq K_X(B)$. Arbitrary intersections
of closed sets are closed: a net in their intersection is a net in
each set, so its limit belongs to every set. In particular $K_X(B)$
is closed, and therefore
\[
 a_X^2(B)\subseteq a_X(K_X(B))=K_X(B).
\]
Under \textnormal{(CB)}, apply \eqref{nx:CB-equation} to the null
saturation of $(B_\lambda)$. Its tail at an index $(B_\lambda,b)$
has image $B_\lambda$, so the resulting closed-tail saturation has
tail sets $K_G(B_\lambda)$. It is tail equivalent to
$\langle K_G(B_\lambda)\rangle$, which is consequently null.
For every $\lambda$,
\[
 a_G(B_\lambda)\subseteq K_G(B_\lambda),
\]
so the saturation of the adherence sets is a quasi-subnet of this
null net. This proves \textnormal{(R)}. Finally, monotonicity makes the adherence families
decreasing on the same index $\Lambda$, so two successive applications
of \eqref{ci:R} give the last assertion.
\end{proof}

\begin{definition}[Absorption of independently null rows]
\label{ci:two-scale}
For an infinite cardinal $\kappa$, a convergence group $G$ satisfies
\textnormal{$(P^\times_\kappa)$} if the following implication holds.
Let $(A_\lambda)_{\lambda\in\Lambda}$ be a decreasing family of
nonempty sets with $|\Lambda|\le\kappa$ and
$\langle A_\lambda\rangle\to e$, and let $(z_i)_{i\in I}$ be a net.
For each $\lambda$, suppose that a net
$(h^\lambda_j)_{j\in J_\lambda}\to e$ satisfies
\begin{equation}\label{ci:P}
 \forall j\in J_\lambda\ \exists i_0\in I\ \forall i\ge i_0:
 \quad z_i\in
 A_\lambda\{h^\lambda_k:k\ge j\}.
\end{equation}
Then $z_i\to e$. The unrestricted condition is
\textnormal{$(P_\times)$}. In additive abelian groups the corresponding
notations are \textnormal{$(P_\kappa)$} and \textnormal{(P)}, with set
addition replacing set multiplication.
\end{definition}

Each row in \eqref{ci:P} retains its own directed set and convergence
thresholds, so absorption concerns the compatibility of errors whose
inner convergence has been established separately. Its formulation
therefore precedes any choice of integration process or proposed limiting
integral: restricting the outer family to a base of at most $\kappa$
sets limits how many rows are combined, while leaving both the size of
each $J_\lambda$ and the candidate index $I$ unrestricted.

\subsection{Interchange with a prescribed error base}

With these operations available, we can follow a finite comparison
through its successive limits, first obtaining a Cauchy family of
endpoints by adherence and then applying adherence once more to recover
the entire approximation net from a smaller auxiliary index. This
reduction places the use of completeness at the cardinality of the error
base; only after its limit has been constructed do we use absorption to
establish convergence of the candidate integration row.

\begin{theorem}[Cardinal intrinsic interchange in convergence groups]
\label{ci:interchange}
Let $\kappa$ be an infinite cardinal, and let $G$ be separated, right
$\kappa$-complete, and satisfy \textnormal{$(R_\kappa)$} and
\textnormal{$(P^\times_\kappa)$}. Let $P$ and $\Gamma_A$, for
$A\in\mathcal A\ne\varnothing$, be directed sets, and suppose
\begin{equation}\label{ci:rows}
 x^A_{p,\gamma}\longrightarrow y_p^A
 \qquad(p\in P,\ A\in\mathcal A).
\end{equation}
Let $(Q_\lambda)_{\lambda\in\Lambda}$ decrease through nonempty
subsets of $G$, with
\begin{equation}\label{ci:null-errors}
 |\Lambda|\le\kappa,\qquad
 \langle Q_\lambda:\lambda\in\Lambda\rangle\longrightarrow e.
\end{equation}
For each $\lambda$, let $P_\lambda\subseteq P$ contain a tail of
$P$, independently of $A$, and assume
\begin{equation}\label{ci:comparison}
 p\in P_\lambda
 \quad\Longrightarrow\quad
 x^A_\gamma(x^A_{p,\gamma})^{-1}\in Q_\lambda
 \quad\text{eventually in }\gamma.
\end{equation}
Then unique $z^A\in G$ satisfy
\[
 x^A_\gamma\longrightarrow z^A,\qquad
 y_p^A\longrightarrow z^A\quad\text{right uniformly in }A,
\]
and
\begin{equation}\label{ci:endpoint}
 p\in P_\lambda
 \quad\Longrightarrow\quad
 y_p^A(z^A)^{-1}\in
 a_G^2(Q_\lambda^{-1}Q_\lambda).
\end{equation}
Without \textnormal{$(P^\times_\kappa)$}, the endpoint family, its
uniform convergence, and \eqref{ci:endpoint} still exist. A countable
error base requires only right sequential completeness.
\end{theorem}
\begin{proof}
Put
\[
 S_\lambda=Q_\lambda^{-1}Q_\lambda,\qquad
 C_\lambda=a_G(S_\lambda),\qquad
 D_\lambda=a_G(C_\lambda).
\]
To see that the saturation of $(S_\lambda)$ is null, take two error
levels $\lambda,\eta$ and observe that any $\nu\ge\lambda,\eta$ satisfies
$S_\nu\subseteq Q_\lambda^{-1}Q_\eta$, so this saturation is a
quasi-subnet of the product of the inverse and the original saturated
error nets, which converges to $e$. Because the error family remains indexed by $\Lambda$,
two applications of \textnormal{$(R_\kappa)$} give
\begin{equation}\label{ci:null-adherences}
 \langle C_\lambda\rangle\longrightarrow e,
 \qquad \langle D_\lambda\rangle\longrightarrow e.
\end{equation}
These families are decreasing and contain $e$.

Fix $A$, $\lambda$, and $p,q\in P_\lambda$. On a common inner tail,
\[
 u_\gamma=x^A_\gamma(x^A_{p,\gamma})^{-1}\in Q_\lambda,
 \qquad
 v_\gamma=x^A_\gamma(x^A_{q,\gamma})^{-1}\in Q_\lambda,
\]
and therefore
\[
 x^A_{p,\gamma}(x^A_{q,\gamma})^{-1}
 =u_\gamma^{-1}v_\gamma\in S_\lambda.
\]
Restricting to that tail and using
\eqref{ci:rows} together with the group operations gives the limit
$y_p^A(y_q^A)^{-1}$, and hence
\begin{equation}\label{ci:row-difference}
 p,q\in P_\lambda
 \quad\Longrightarrow\quad
 y_p^A(y_q^A)^{-1}\in C_\lambda.
\end{equation}

To construct the smaller index from these comparisons, for each
$\lambda$ fix $r_\lambda\in P$ such that
$p\ge r_\lambda$ implies $p\in P_\lambda$. Fix an upper-bound
map $b$ on the nonempty finite subsets of $P$ and an index $r_*\in P$.
Define $p_F$, for $F\in\operatorname{Fin}(\Lambda)$, by recursion
on the finite cardinality of $F$:
\[
 p_\varnothing=r_*,\qquad
 p_F=b\bigl(\{r_\lambda:\lambda\in F\}
             \cup\{p_E:E\subsetneq F\}\bigr)
 \quad(F\ne\varnothing).
\]
For $F\ne\varnothing$, the displayed union contains
$p_\varnothing$ and has at most $|F|+2^{|F|}-1$ members.
Thus $b$ is applicable at every step. If $E\subsetneq F$, the
point $p_E$ is one of its arguments, and if $\lambda\in F$, the
point $r_\lambda$ is another argument. Consequently
\begin{equation}\label{ci:finite-indices}
 E\subseteq F\Longrightarrow p_E\le p_F,
 \qquad
 \lambda\in F\Longrightarrow p_F\ge r_\lambda
                         \Longrightarrow p_F\in P_\lambda.
\end{equation}
For an infinite $\Lambda$,
$|\operatorname{Fin}(\Lambda)|=|\Lambda|$, since the singleton map
embeds $\Lambda$ into $\operatorname{Fin}(\Lambda)$ and
\[
 |\operatorname{Fin}(\Lambda)|
 \le\sum_{n<\omega}|\Lambda|^n=|\Lambda|.
\]
For finite $\Lambda$, $\operatorname{Fin}(\Lambda)$ is finite.
In both cases its cardinality is at most $\kappa$. Set $u_F^A=y_{p_F}^A$. For $F,H\supseteq\{\lambda\}$,
\eqref{ci:row-difference} gives
\[
 u_F^A(u_H^A)^{-1}\in C_\lambda.
\]
More explicitly, the tail of the right-difference net above
$(\{\lambda\},\{\lambda\})$ satisfies
\[
 \{u_F^A(u_H^A)^{-1}:F,H\supseteq\{\lambda\}\}
 \subseteq C_\lambda.
\]
Every tail of $\langle C_\lambda\rangle$ is one of these controlling
sets, so the right-difference net is a quasi-subnet of that null
saturation. We have therefore obtained a right Cauchy net $(u_F^A)_F$
on an index of cardinality at most $\kappa$, for which right
$\kappa$-completeness supplies
\begin{equation}\label{ci:reduced-limit}
 u_F^A\longrightarrow z^A.
\end{equation}

To return to the original approximation index, let $p\in P_\lambda$.
Then $y_p^A(u_F^A)^{-1}\in C_\lambda$ for
every $F\supseteq\{\lambda\}$. The net on that tail converges to
$y_p^A(z^A)^{-1}$ by \eqref{ci:reduced-limit}, proving
\eqref{ci:endpoint}. In particular,
\[
 \{y_p^A(z^A)^{-1}:p\ge r_\lambda,\ A\in\mathcal A\}
 \subseteq D_\lambda.
\]
Set
\[
 E_{p_0}=\{y_p^A(z^A)^{-1}:p\ge p_0,\ A\in\mathcal A\}.
\]
These sets $E_{p_0}$ are decreasing and nonempty, and the preceding
inclusion gives $E_{r_\lambda}\subseteq D_\lambda$ for every
$\lambda$. Thus $\langle E_{p_0}\rangle$ is a quasi-subnet of
$\langle D_\lambda\rangle$ and converges to $e$. This is the stated
right uniform convergence. For a fixed $A$, the tail of
$(y_p^A(z^A)^{-1})_p$ above $p_0$ is contained in $E_{p_0}$,
so $y_p^A(z^A)^{-1}\to e$, and right translation gives
$y_p^A\to z^A$. If another family $w^A$ has the asserted uniform
endpoint convergence, the same restriction gives $y_p^A\to w^A$.
Separation implies $w^A=z^A$ for every $A$. This argument uses only
\eqref{ci:finite-indices}, not cofinality of $F\mapsto p_F$ in $P$.

For the candidate row, put $R_\lambda=Q_\lambda D_\lambda$. The
family decreases and its saturation is null, since for
$\nu\ge\lambda,\eta$,
$R_\nu\subseteq Q_\lambda D_\eta$, and multiplication is continuous.
Thus
\begin{equation}\label{ci:outer-null}
 \langle R_\lambda\rangle\longrightarrow e.
\end{equation}
Fix $A$ and $\lambda$, put $p=r_\lambda$, and set
\[
 c_\lambda=y_p^A(z^A)^{-1}\in D_\lambda,
 \qquad h^\lambda_\gamma
 =c_\lambda^{-1}
   [x^A_{p,\gamma}(y_p^A)^{-1}]c_\lambda
 \longrightarrow e.
\]
On a tail supplied by \eqref{ci:comparison}, the factorization
\[
 \begin{aligned}
 x^A_\gamma(z^A)^{-1}
 &=\bigl[x^A_\gamma(x^A_{p,\gamma})^{-1}c_\lambda\bigr]
      h^\lambda_\gamma\in R_\lambda h^\lambda_\gamma
 \end{aligned}
\]
holds. For every $\gamma_0\in\Gamma_A$, a common upper bound of
$\gamma_0$ and this comparison threshold gives
\[
 x^A_\gamma(z^A)^{-1}
 \in R_\lambda
       \{h^\lambda_\eta:\eta\ge\gamma_0\}
 \quad\text{eventually in }\gamma.
\]
We may now apply absorption: the outer family in \eqref{ci:outer-null}
is decreasing, nonempty, and indexed by $\Lambda$, with
$|\Lambda|\le\kappa$, while for each $\lambda$ the displayed net
$h^\lambda$ is null. Together with these facts, the last containment
verifies every quantifier in \eqref{ci:P}
for the candidate net $(x^A_\gamma(z^A)^{-1})_\gamma$.
Property \textnormal{$(P^\times_\kappa)$} therefore gives
$x^A_\gamma(z^A)^{-1}\to e$, and right translation gives
$x^A_\gamma\to z^A$. The construction of $z^A$, the endpoint
inclusion, and uniform endpoint convergence were obtained before this
application, which proves the assertion without absorption. If
$|\Lambda|\le\aleph_0$, the only use of completeness was for
$(u_F^A)_{F\in\operatorname{Fin}(\Lambda)}$, whose index is
countable. Lemma~\ref{ci:countable-completeness} supplies precisely
that completeness from right sequential completeness.
\end{proof}

\begin{corollary}[The unrestricted abelian specialization]
\label{ci:abelian-recovery}\label{mg:intrinsic-interchange}
Let $Z$ be a complete separated abelian convergence group satisfying
\textnormal{(R)} and \textnormal{(P)}. Retain the rows and null error
family of Theorem~\ref{ci:interchange}, without a cardinal restriction,
and assume the additive comparisons
\begin{equation}\label{mg:intrinsic-comparison}
 p\in P_\lambda\quad\Longrightarrow\quad
 x^A_\gamma-x^A_{p,\gamma}\in Q_\lambda
 \quad\text{eventually in }\gamma.
\end{equation}
Then its conclusions hold, with
\begin{equation}\label{mg:intrinsic-endpoint}
 p\in P_\lambda\quad\Longrightarrow\quad
 y_p^A-z^A\in a_Z^2(Q_\lambda-Q_\lambda)
                 \subseteq K_Z(Q_\lambda-Q_\lambda).
\end{equation}
In particular, \textnormal{(CB)} is a sufficient replacement for
\textnormal{(R)}.
\end{corollary}
\begin{proof}
Take $\kappa\ge\max(\aleph_0,|\Lambda|)$. Completeness gives
$\kappa$-completeness, and \textnormal{(R)} gives
\textnormal{$(R_\kappa)$}. Unrestricted absorption gives
\textnormal{$(P_\kappa)$}. Theorem~\ref{ci:interchange}, written
additively, yields the limits and
\[
 y_p^A-z^A\in a_Z^2(Q_\lambda-Q_\lambda)
             \subseteq K_Z(Q_\lambda-Q_\lambda).
\]
If \textnormal{(CB)} is assumed in place of \textnormal{(R)},
Lemma~\ref{ci:adherence-calculus} first gives \textnormal{(R)}
for every null family, after which the same application of
Theorem~\ref{ci:interchange} has all the required hypotheses.
\end{proof}

\subsection{Finite bilinear responses and uniform integration}
\label{uc:finite-intrinsic}

Having isolated the passages through limits, we turn to the finite
comparison that supplies their common error family. It must be obtained
from the prescribed records before a limit for the candidate function
is available, with equal tags combined because a function assigns only
one value to each tag; the resulting response family then records exactly
the perturbations which the actual finite sums can produce.

Let $X,Z$ be separated abelian convergence groups. Write the finite sums
on the original integration nets as
\[
 S_\gamma^Af=
 \sum_{r=1}^{m(A,\gamma)}L_{A,\gamma,r}f(t_{A,\gamma,r}),
 \qquad L_{A,\gamma,r}\in\operatorname{Hom}(X,Z).
\]
For the coefficient products of Definition~\ref{bilinearproduct}, we have
$L_{A,\gamma,r}x=b(x,\mu(C_{A,\gamma,r}))$, with $b$ biadditive
and with continuous partial maps, although the definitions below apply
directly to any prescribed family of continuous additive coefficient maps.

\begin{definition}[Finite response family]\label{uc:response}\label{lt:csv}
For each record put
\[
 T_{A,\gamma}=\{t_{A,\gamma,r}:1\le r\le m(A,\gamma)\},
 \qquad
 L_t^{A,\gamma}=\sum_{r:t_{A,\gamma,r}=t}L_{A,\gamma,r}.
\]
For $0\in B\subseteq X$, define
\[
 \Sigma(B)=\{0\}\cup\bigcup_{A,\gamma}
 \left\{\sum_{t\in T_{A,\gamma}}L_t^{A,\gamma}x_t:
                                        x_t\in B\right\}.
\]
The coefficients have continuous finite semivariance at cardinal
$\kappa$ if
\[
 \begin{gathered}
 B_\lambda\downarrow,\quad 0\in B_\lambda,\quad |\Lambda|\le\kappa,
 \quad\langle B_\lambda\rangle\longrightarrow0\text{ in }X\quad\Longrightarrow\quad
 \langle\Sigma(B_\lambda)\rangle\longrightarrow0\text{ in }Z.
 \end{gathered}
\]
Continuous finite semivariance without qualification means the same
assertion for every directed error family.
\end{definition}

\begin{lemma}[Actual responses and finite subwords]
\label{lt:reduced-words}
For every $0\in B\subseteq X$,
\begin{equation}\label{uc:response-identity}
 \Sigma(B)=\{0\}\cup
 \{S_\gamma^Ah:h[\Omega]\subseteq B,\ A\in\mathcal A,
                                      \gamma\in\Gamma_A\}.
\end{equation}
Adjoining subwords of reduced coefficient words leaves $\Sigma(B)$
unchanged. In particular, partial responses of disjoint internally tagged
records belong to the same response family.
\end{lemma}
\begin{proof}
Each function $h$ in the right-hand side assigns values
$x_t=h(t)\in B$ to the distinct tags, so its response belongs to the
left-hand side. Conversely, a prescribed assignment $x_t\in B$ on the
finite set $T_{A,\gamma}$ extends to $\Omega$ by the value $0\in B$,
and the extended function has the prescribed response. For a subword on
$T'\subseteq T_{A,\gamma}$, put $x_t=0$ on its complement. Then
\[
 \sum_{t\in T'}L_t^{A,\gamma}x_t
 =\sum_{t\in T_{A,\gamma}}L_t^{A,\gamma}x_t\in\Sigma(B).
\]
Disjoint cells with internal tags have distinct tags, so each partial
record is such a subword.
\end{proof}

\begin{theorem}[Uniform approximation with a prescribed error base]
\label{cu:uniform-cardinal}
Let $\kappa$ be infinite, and let $X,Z$ be separated abelian convergence
groups. Suppose that $Z$ is $\kappa$-complete and satisfies
\textnormal{$(R_\kappa)$} and \textnormal{$(P_\kappa)$}, and that the
coefficients have continuous finite semivariance at cardinal $\kappa$.
Let $(f_p)_{p\in P}$ be integrable on every $A\in\mathcal A\ne\varnothing$
on the prescribed net $\Gamma_A$, and let $f:\Omega\to X$.
Suppose that a decreasing family $(B_\lambda)_{\lambda\in\Lambda}$
satisfies
\begin{equation}\label{cu:input-null-base}
 |\Lambda|\le\kappa,\qquad 0\in B_\lambda,\qquad
 \langle B_\lambda\rangle\longrightarrow0,
\end{equation}
and
\begin{equation}\label{cu:input-thresholds}
 \forall\lambda\ \exists p_\lambda\ \forall p\ge p_\lambda\quad
 \forall t\in\Omega:\quad f_p(t)-f(t)\in B_\lambda.
\end{equation}
Then $f$ is integrable on every $A$, and
\begin{equation}\label{cu:uniform-conclusion}
 I_Af_p\longrightarrow I_Af\quad\text{uniformly in }A.
\end{equation}
For each $A$, both iterated limits exist and satisfy
\begin{equation}\label{cu:iterated-limits}
 \lim_p\lim_{\gamma\in\Gamma_A}S_\gamma^Af_p
 =I_Af
 =\lim_{\gamma\in\Gamma_A}\lim_pS_\gamma^Af_p.
\end{equation}
More precisely,
\begin{equation}\label{cu:endpoint-inclusion}
 p\ge p_\lambda\quad\Longrightarrow\quad
 I_Af_p-I_Af\in a_Z^2\bigl(\Sigma(B_\lambda)-\Sigma(B_\lambda)\bigr)
 \qquad(A\in\mathcal A).
\end{equation}
\end{theorem}
\begin{proof}
For every actual record, additivity and grouping of equal tags give
\[
 S_\gamma^Af_p-S_\gamma^Af
 =\sum_{t\in T_{A,\gamma}}L_t^{A,\gamma}(f_p(t)-f(t)).
\]
Thus
\begin{equation}\label{cu:finite-comparison}
 p\ge p_\lambda\quad\Longrightarrow\quad
 S_\gamma^Af_p-S_\gamma^Af\in\Sigma(B_\lambda)
 \qquad(A,\gamma).
\end{equation}
If $\eta\ge\lambda$, then every assignment with values in
$B_\eta$ also has values in $B_\lambda$, and hence
$\Sigma(B_\eta)\subseteq\Sigma(B_\lambda)$. Each response set
contains zero. Applying finite semivariance at cardinal $\kappa$
to \eqref{cu:input-null-base} gives
$\langle\Sigma(B_\lambda)\rangle\to0$; continuity of inversion
also gives $\langle-\Sigma(B_\lambda)\rangle\to0$.
Apply Theorem~\ref{ci:interchange} in additive notation with
\[
 x^A_{p,\gamma}=S_\gamma^Af_p,\qquad
 x^A_\gamma=S_\gamma^Af,\qquad
 y_p^A=I_Af_p,\qquad Q_\lambda=-\Sigma(B_\lambda),
\]
and $P_\lambda=\{p:p\ge p_\lambda\}$. The known rows are exactly
the assumed integrability statements, and \eqref{cu:finite-comparison}
is the required comparison on every record. The theorem constructs
$z^A$ with $S_\gamma^Af\to z^A$ and $I_Af_p\to z^A$ uniformly in
$A$. By the definition of the original integral, $z^A=I_Af$.
Its endpoint inclusion is \eqref{cu:endpoint-inclusion}, since
$-Q_\lambda+Q_\lambda=\Sigma(B_\lambda)-\Sigma(B_\lambda)$.
For fixed $A,\gamma$, equation \eqref{cu:finite-comparison} gives
\[
 \{S_\gamma^Af_p-S_\gamma^Af:p\ge p_\lambda\}
 \subseteq\Sigma(B_\lambda)\qquad(\lambda\in\Lambda).
\]
Thus the difference net is a quasi-subnet of
$\langle\Sigma(B_\lambda)\rangle$, and
\[
 \lim_p S_\gamma^Af_p=S_\gamma^Af.
\]
For these fixed-record limits, the convergence already constructed yields
\[
 \lim_{\gamma\in\Gamma_A}\bigl(\lim_pS_\gamma^Af_p\bigr)
 =\lim_{\gamma\in\Gamma_A}S_\gamma^Af=I_Af.
\]
For the other order, the assumed integrability gives
$\lim_{\gamma\in\Gamma_A}S_\gamma^Af_p=I_Af_p$ for each $p$,
and the endpoint convergence gives
\[
 \lim_p\bigl(\lim_{\gamma\in\Gamma_A}S_\gamma^Af_p\bigr)
 =\lim_p I_Af_p=I_Af.
\]
Together, these calculations give the two equalities in
\eqref{cu:iterated-limits}.
\end{proof}

\begin{corollary}[Sequential completeness for arbitrary integration nets]
\label{cu:uniform-sequential}
Suppose $Z$ is separated and sequentially complete and satisfies
\textnormal{$(R_\omega)$} and \textnormal{$(P_\omega)$}. Assume
continuous finite semivariance at cardinal $\aleph_0$. If the known
integrable approximants admit decreasing sets $B_n\ni0$ such that
\[
 \langle B_n\rangle\to0,\qquad
 \forall n\ \exists p_n\ \forall p\ge p_n\ \forall t:
                  f_p(t)-f(t)\in B_n,
\]
then all conclusions of Theorem~\ref{cu:uniform-cardinal} hold.
In particular, they hold for a uniformly convergent sequence of
integrable functions. The original approximation index, integration
nets, and parameter family remain unrestricted.
\end{corollary}
\begin{proof}
Lemma~\ref{ci:countable-completeness} gives $\aleph_0$-completeness,
and the displayed family verifies
\eqref{cu:input-null-base}--\eqref{cu:input-thresholds} with
$\Lambda=\mathbb N$. For a uniformly convergent sequence put
\[
 B_n=\{0\}\cup\{f_m(t)-f(t):m\ge n,\ t\in\Omega\}.
\]
If $\Omega$ is empty, these sets equal $\{0\}$. Otherwise the
error-tail saturation is null by uniform convergence. Mixing that net
with the constant zero net gives a null net whose tail values are the
error-tail sets with zero adjoined, so $\langle B_n\rangle\to0$.
The implication $m\ge n\Rightarrow f_m(t)-f(t)\in B_n$ verifies
the threshold requirement.
\end{proof}

When the limiting function is already integrable, the construction just
completed is unnecessary: the difference of the two sum nets can pass to
its limit from a single error set. The resulting estimate uses adherence only once,
which isolates the part of the uniform theorem that survives without
completeness or absorption.

\begin{theorem}[Preservation when the limit is integrable]
\label{uc:known-limit}
Retain the coefficient and approximation assumptions of
Theorem~\ref{cu:uniform-cardinal}, but assume only
\textnormal{$(R_\kappa)$} on the separated target $Z$. If $f$ is
already integrable on every $A$, then
\[
 I_Af_p\longrightarrow I_Af\quad\text{uniformly in }A,
\]
with the sharper one-step inclusion
\[
 p\ge p_\lambda\quad\Longrightarrow\quad
 I_Af_p-I_Af\in a_Z(\Sigma(B_\lambda)).
\]
No completeness or absorption assumption is required.
\end{theorem}
\begin{proof}
Equation \eqref{cu:finite-comparison} holds before either limit is
taken. For fixed $p\ge p_\lambda$ and $A$, the two integrability
assumptions give
\[
 S_\gamma^Af_p-S_\gamma^Af\longrightarrow I_Af_p-I_Af.
\]
Every value lies in $\Sigma(B_\lambda)$, proving the displayed
adherence inclusion. The response family is null and indexed by at
most $\kappa$ levels, so \textnormal{$(R_\kappa)$} makes its
adherence saturation null. For
\[
 E_{p_0}=\{I_Af_p-I_Af:p\ge p_0,\ A\in\mathcal A\},
\]
the adherence inclusion gives
$E_{p_\lambda}\subseteq a_Z(\Sigma(B_\lambda))$ for every
$\lambda$. Hence each tail of the null adherence saturation contains
a tail of $\langle E_{p_0}\rangle$. The latter saturation is null,
which is uniform convergence over the entire parameter family.
\end{proof}

\subsection{Universal characterization by elementary sampling}
\label{uc:universal-section}

We have obtained the uniform theorem by following errors through finite
responses, adherence and absorption, and can now ask whether these
passages are forced by the theorem itself. Requiring its assertion
throughout the prescribed finite-record class gives an affirmative
answer, provided that the tested family of domains remains among the
data: a single approximation level must be able to test every adherence
limit of its error set while the nets realizing those limits keep their
own directed indices. After specifying this quantification, we construct
three sampling systems that recover adherence stability, completeness
and absorption from the corresponding uniform assertions.

\begin{definition}[The universal assertions at cardinal $\kappa$]
\label{uc:universal-class}
Fix a separated abelian convergence group $Z$ and an infinite cardinal
$\kappa$. Consider all data consisting of a separated abelian convergence
input group $X$, a coefficient product, a paving $\mathcal H$ of an
underlying set $\Omega$, a coefficient function $\mu$, and prescribed
integration nets as in
Definitions~\ref{bilinearproduct}--\ref{integrationstructure}, with
target $Z$ and continuous finite semivariance at cardinal $\kappa$.
For each such system, consider every nonempty tested family
$\mathcal A\subseteq\mathcal H$, every directed approximation index
$P$, and all maps $f_p,f:\Omega\to X$ satisfying
\eqref{cu:input-null-base}--\eqref{cu:input-thresholds} for some
error index $\Lambda$ of cardinality at most $\kappa$.

The \emph{universal preservation assertion} requires, whenever every
$f_p$ and $f$ is integrable on every $A\in\mathcal A$, that
\begin{equation}\label{uc:universal-preservation-formula}
 \left\langle
 \{I_Af_r-I_Af:r\ge p,\ A\in\mathcal A\}:p\in P
 \right\rangle\longrightarrow0.
\end{equation}
The \emph{universal construction assertion} assumes only that every
$f_p$ is integrable on every $A\in\mathcal A$ and requires the
existence of $z_A\in Z$ satisfying
\begin{equation}\label{uc:universal-construction-formula}
 \begin{gathered}
 S^A_\gamma f\longrightarrow z_A\qquad(A\in\mathcal A),\\
 \left\langle
 \{I_Af_r-z_A:r\ge p,\ A\in\mathcal A\}:p\in P
 \right\rangle\longrightarrow0.
 \end{gathered}
\end{equation}
For each $A$, its prescribed net $\Gamma_A$ is the same for all
$f_p$ and for $f$. Separation identifies $z_A$ with $I_Af$ when
\eqref{uc:universal-construction-formula} holds.
\end{definition}

To recover these properties, we use sampling systems whose unit
coefficients leave the errors unchanged. Let $Q$ be
a nonempty set, let $J_q$ be a nonempty directed set for each $q\in Q$,
and put
\[
 \Omega=\coprod_{q\in Q}(\{q\}\times J_q),\qquad
 \Omega_q=\{q\}\times J_q,\qquad
 \mathcal H=\{\varnothing\}\cup\{\Omega_q:q\in Q\}.
\]
Take $X=Z$, give the coefficient group $\mathbb Z$ its discrete
convergence, and put $b(x,m)=mx$, where integer multiplication is
formed from the group operations of $Z$. For each fixed $m$, the map
$x\mapsto mx$ is continuous by finite addition and inversion. For
fixed $x$, a convergent net of integers is eventually constant, so
$m\mapsto mx$ is continuous. Both maps are additive and preserve zero.
Set
\[
 \mu(\varnothing)=0,\qquad \mu(\Omega_q)=1,\qquad
 R_j^{\Omega_q}=((\Omega_q,(q,j)))\quad(j\in J_q),
\]
and prescribe the empty record on $\varnothing$. The tested family is
$\mathcal A=\{\Omega_q:q\in Q\}$. Its sums and response sets are
\begin{equation}\label{uc:unit-sampling}
 S_j^{\Omega_q}h=h(q,j),\qquad
 \Sigma(B)=\{0\}\cup\bigcup_{q,j}\{x:x\in B\}=B
 \quad(0\in B\subseteq Z).
\end{equation}
In particular, for every decreasing null family $B_\lambda\ni0$,
\[
 \langle\Sigma(B_\lambda)\rangle
 =\langle B_\lambda\rangle\longrightarrow0,
\]
without a restriction on the cardinality of its index. A single
sampling net $t_\gamma\in\Omega$ gives the same identities by taking
$\mathcal H=\{\varnothing,\Omega\}$, the single record
$((\Omega,t_\gamma))$, and $\mu(\Omega)=1$.

In these systems, inactive sampling domains will contribute zero, so
we record why adjoining zero to a null error family preserves its
convergence. If $a_d$ is its saturation, the net on
$D\times\{0,1\}$ defined by
\[
 \widetilde a_{d,0}=a_d,\qquad \widetilde a_{d,1}=0,\qquad
 (d,\varepsilon)\preceq(d',\varepsilon')\iff d\preceq d'
\]
converges to zero by finite mixing. Its tail-value sets are the
original tail-value sets with zero adjoined. Therefore
\begin{equation}\label{uc:null-zero-adjunction}
 \langle A_\lambda\rangle\longrightarrow0
 \quad\Longrightarrow\quad
 \langle A_\lambda\cup\{0\}\rangle\longrightarrow0.
\end{equation}
Inversion gives the corresponding implication for
$(-A_\lambda)\cup\{0\}$.

\begin{proposition}[The sampling test for adherence]
\label{uc:necessity-adherence}
Suppose that the preservation conclusion in
Definition~\ref{uc:universal-class} holds for every unit-coefficient
sampling system with an input error base of cardinality at most
$\kappa$. Then $Z$ satisfies \textnormal{$(R_\kappa)$}. In
particular, the universal preservation assertion at cardinal $\kappa$
implies \textnormal{$(R_\kappa)$}.
\end{proposition}
\begin{proof}
Let $(A_\lambda)_{\lambda\in\Lambda}$ be decreasing and nonempty,
with $|\Lambda|\le\kappa$ and $\langle A_\lambda\rangle\to0$.
Give $P=\Lambda\times\mathbb N$ the product order and define
\[
 Q=\{((\lambda,n),y):\lambda\in\Lambda,\ n\in\mathbb N,
                              \ y\in a_Z(A_\lambda)\}.
\]
For $q=((\lambda,n),y)$ there is a net
$(x_j^q)_{j\in J_q}$ with every value in $A_\lambda$ and limit $y$.
This collection can be selected from a set of possible representatives:
replace each witnessing net by the saturation of its nonempty
actual tail-value sets, all contained in $A_\lambda$. These families
belong to $\mathcal P(\mathcal P(Z))$, their saturation indices are
subsets of $\mathcal P(Z)\times Z$, and the replacement has the
same limits by \eqref{fd:tail-representation}. Thus a set-indexed
choice supplies
\begin{equation}\label{uc:adherence-witnesses}
 x_j^q\in A_\lambda\quad(j\in J_q),\qquad
 x_j^q\longrightarrow y
 \quad\bigl(q=((\lambda,n),y)\bigr).
\end{equation}
Use these $J_q$ in the system \eqref{uc:unit-sampling}. For $r\in P$,
define
\[
 f_r(q,j)=
 \begin{cases}
 x_j^q,&q=(r,y)\text{ for some }y,\\
 0,&q=(p,y)\text{ with }p\ne r,
 \end{cases}
 \qquad f(q,j)=0.
\]
For fixed $q=(p,y)$ and $r$, the sampling net is either the net
in \eqref{uc:adherence-witnesses} or the constant zero net. Hence
\begin{equation}\label{uc:adherence-integrals}
 I_{\Omega_q}f_r=
 \begin{cases}y,&p=r,\\0,&p\ne r,
 \end{cases}
 \qquad I_{\Omega_q}f=0.
\end{equation}
Put $B_\lambda=A_\lambda\cup\{0\}$ and
$p_\lambda=(\lambda,1)$. By
\eqref{uc:null-zero-adjunction}, $\langle B_\lambda\rangle\to0$.
If $r=(\eta,n)\ge p_\lambda$, then $\eta\ge\lambda$, and every
$(q,j)\in\Omega$ satisfies
\[
 f_r(q,j)-f(q,j)\in\{0\}\cup A_\eta
 \subseteq\{0\}\cup A_\lambda=B_\lambda.
\]
Thus \eqref{cu:input-null-base}--\eqref{cu:input-thresholds} hold
with the index $\Lambda$ and with the displayed thresholds.

For $p_0=(\lambda_0,n_0)$, define the endpoint-tail set
\[
 E_{p_0}=\{I_{\Omega_q}f_r:r\ge p_0,\ q\in Q\}.
\]
If an element of $E_{p_0}$ is nonzero, then
\eqref{uc:adherence-integrals} places it in $a_Z(A_\eta)$ for some
$\eta\ge\lambda_0$, and therefore in $a_Z(A_{\lambda_0})$.
For the reverse inclusion, each $y\in a_Z(A_{\lambda_0})$ is
obtained by taking $r=p_0$ and $q=(p_0,y)$. To obtain zero, fix
$a\in A_{\lambda_0}$, take $r=p_0$ and
$q=((\lambda_0,n_0+1),a)$, and use the inactive case of
\eqref{uc:adherence-integrals}. We have proved
\begin{equation}\label{uc:adherence-endpoint-tails}
 E_{(\lambda_0,n_0)}=a_Z(A_{\lambda_0})\cup\{0\}.
\end{equation}
By the assumed preservation conclusion,
$\langle E_p:p\in P\rangle\to0$, while the family on the right of
\eqref{uc:adherence-endpoint-tails} is independent of $n_0$ and
contains $a_Z(A_{\lambda_0})$. For every endpoint tail $E_{p_0}$,
the saturation of $a_Z(A_\lambda)$ consequently has a tail contained
in $E_{p_0}$, namely $a_Z(A_{\lambda_0})$. It is a quasi-subnet
of the null endpoint saturation and converges to zero, as required
by \textnormal{$(R_\kappa)$}.
\end{proof}

\begin{proposition}[The sampling test for completeness]
\label{uc:necessity-completeness}
Suppose that the construction conclusion in
Definition~\ref{uc:universal-class} holds for every unit-coefficient
sampling system with an input error base of cardinality at most
$\kappa$. Then $Z$ is $\kappa$-complete. In particular, the universal
construction assertion at cardinal $\kappa$ implies
$\kappa$-completeness.
\end{proposition}
\begin{proof}
Let $(x_i)_{i\in I}$ be a Cauchy net with $|I|\le\kappa$.
On $\Omega=I$, take the identity sampling net and the single-domain
system of \eqref{uc:unit-sampling}. For $p\in I$, put
\[
 f_p(i)=\begin{cases}x_p,&i\ge p,\\x_i,&i\not\ge p,\end{cases}
 \qquad f(i)=x_i.
\]
For $i\ge p$, its sampling sum is $S_i f_p=x_p$, so the net of
sums is eventually constant and
\begin{equation}\label{uc:cauchy-test-integrals}
 If_p=x_p.
\end{equation}
Define
\[
 B_k=\{x_p-x_i:p,i\ge k\}\qquad(k\in I).
\]
These sets are nonempty, decrease with $k$, and contain zero by taking
$p=i=k$. The difference net $d_{p,i}=x_p-x_i$ on $I\times I$
converges to zero by the Cauchy assumption. Its tail above $(k,k)$
is $B_k$, while for each $(p_0,i_0)$ an upper bound $k\ge p_0,i_0$
satisfies
\[
 B_k\subseteq\{x_p-x_i:p\ge p_0,\ i\ge i_0\}.
\]
Thus the difference net and $\langle B_k:k\in I\rangle$ are
tail equivalent, and the latter is null. For $p\ge k$ and $i\in I$,
\[
 f_p(i)-f(i)=
 \begin{cases}x_p-x_i\in B_k,&i\ge p,\\0\in B_k,&i\not\ge p.
 \end{cases}
\]
With index $I$, cardinality at most
$\kappa$, and threshold $p_k=k$, this input error base permits the
construction conclusion to supply $z\in Z$ with $S_i f\to z$.
Since $S_i f=f(i)=x_i$ for every $i$, the original Cauchy net
converges to $z$.
\end{proof}

Recovering absorption requires a further construction, because its rows
need not share a directed index. We place their finite tail constraints
on a common sampling space, retaining each coordinate convergence while
making the candidate coordinate tail equivalent to the original net.

\begin{lemma}[Sampling realization of independently null rows]
\label{uc:absorption-realization}
Let $\Lambda$ be a nonempty set and let $A_\lambda\subseteq Z$ be
nonempty for every $\lambda\in\Lambda$. Let $(z_i)_{i\in I}$ be a
net, and for each $\lambda$ let
$(h_j^\lambda)_{j\in J_\lambda}\to0$. Suppose
\begin{equation}\label{uc:realization-hypothesis}
 \forall\lambda\in\Lambda\ \forall j\in J_\lambda\
 \exists i(\lambda,j)\in I\ \forall i\ge i(\lambda,j):
 \quad z_i\in A_\lambda+\{h_k^\lambda:k\ge j\}.
\end{equation}
There is a nonempty set $T$, a sampling net
$(t_\gamma)_{\gamma\in\Gamma}$ in $T$, and maps
$\pi,\pi_\lambda:T\to Z$ such that
\begin{equation}\label{uc:realization-conclusions}
 \begin{gathered}
 \pi(t)-\pi_\lambda(t)\in A_\lambda
       \quad(t\in T,\ \lambda\in\Lambda),\\
 \pi_\lambda(t_\gamma)\longrightarrow0
       \quad(\lambda\in\Lambda),\\
 (\pi(t_\gamma))_{\gamma\in\Gamma}
       \text{ is tail equivalent to }(z_i)_{i\in I}.
 \end{gathered}
\end{equation}
Neither directedness nor a convergence assumption on the family
$(A_\lambda)$ is needed for this realization.
\end{lemma}
\begin{proof}
Write
\[
 Z_{i_0}=\{z_i:i\ge i_0\},\qquad
 H_j^\lambda=\{h_k^\lambda:k\ge j\},
\]
and fix $a_\lambda\in A_\lambda$ for every $\lambda$. Define
\begin{equation}\label{uc:realization-space}
 T=\{(x,(y_\lambda)_{\lambda\in\Lambda})\in Z\times Z^\Lambda:
                     x-y_\lambda\in A_\lambda\ (\lambda\in\Lambda)\},
 \qquad \pi(x,y)=x,\quad \pi_\lambda(x,y)=y_\lambda.
\end{equation}
To impose finitely many row conditions at a time, take as constraints all triples
$c=(i_0,F,(j_\lambda)_{\lambda\in F})$, where $i_0\in I$,
$F\subseteq\Lambda$ is finite and $j_\lambda\in J_\lambda$.
For such a constraint, put
\begin{equation}\label{uc:realization-constraint}
 C_c=\{(x,y)\in T:x\in Z_{i_0},\quad
                        y_\lambda\in H_{j_\lambda}^\lambda
                                      \ (\lambda\in F)\}.
\end{equation}
Each $C_c$ is nonempty. In fact, an upper bound $i_1$ of
$\{i_0\}\cup\{i(\lambda,j_\lambda):\lambda\in F\}$ satisfies
\[
 z_{i_1}\in Z_{i_0}\cap
                 \bigcap_{\lambda\in F}(A_\lambda+H_{j_\lambda}^\lambda).
\]
Set $x=z_{i_1}$. For each $\lambda\in F$ there is
$y_\lambda\in H_{j_\lambda}^\lambda$ with
$x-y_\lambda\in A_\lambda$, and for $\lambda\notin F$ set
$y_\lambda=x-a_\lambda$. These coordinates define a point of $C_c$.

For $c=(i_0,F,j)$ and $c'=(i'_0,F',j')$, set
\[
 c\preceq c'\quad\Longleftrightarrow\quad
 i_0\le i'_0,\quad F\subseteq F',\quad
 j_\lambda\le j'_\lambda\ (\lambda\in F).
\]
This order is directed: for two constraints, take an upper bound of
their $I$-indices, take the union of their finite coordinate sets,
and in each shared coordinate take an upper bound in $J_\lambda$.
On an unshared coordinate retain its given index. The resulting
constraint is above both, and $c\preceq c'$ implies $C_{c'}\subseteq C_c$.
Define
\[
 \Gamma=\{(c,t):c\text{ is a constraint},\ t\in C_c\},\qquad
 (c,t)\preceq(c',t')\iff c\preceq c',\qquad
 t_{(c,t)}=t.
\]
Nonemptiness of every $C_c$ and directedness of the constraint order
make $\Gamma$ directed. Its tail-value set at $(c,t)$ is exactly
$C_c$: the decreasing property gives containment in $C_c$, and every
$t'\in C_c$ occurs at the later index $(c,t')$. Thus
\begin{equation}\label{uc:realization-tail-identity}
 T_{(c,t)}((t_\gamma))=C_c.
\end{equation}

Fix $\lambda$ and $j\in J_\lambda$. A constraint with
$F=\{\lambda\}$ and $j_\lambda=j$ has
$\pi_\lambda(C_c)\subseteq H_j^\lambda$. By
\eqref{uc:realization-tail-identity}, the coordinate net
$\pi_\lambda(t_\gamma)$ therefore has a tail in every tail-value
set of $h^\lambda$. It is a quasi-subnet of the null net
$h^\lambda$, proving the second assertion in
\eqref{uc:realization-conclusions}. The first assertion is part of
the definition of $T$.

For the last assertion, we calculate the projection of every sampling
tail. Its exact value is
\begin{equation}\label{uc:realization-projection}
 \pi(C_c)=Z_{i_0}\cap
                 \bigcap_{\lambda\in F}(A_\lambda+H_{j_\lambda}^\lambda).
\end{equation}
Membership in $C_c$ gives the inclusion from left to right. For
$x$ in the right-hand side, each $\lambda\in F$ admits a witness
$y_\lambda\in H_{j_\lambda}^\lambda$ with
$x-y_\lambda\in A_\lambda$. For $\lambda\notin F$, the assignment
$y_\lambda=x-a_\lambda$ supplies the remaining coordinates. The
resulting $(x,y)$ belongs to $C_c$, proving the reverse inclusion.
With the threshold $i_1$ chosen above, this identity gives
\begin{equation}\label{uc:realization-tail-bounds}
 Z_{i_1}\subseteq\pi(C_c)\subseteq Z_{i_0}.
\end{equation}
For every original tail $Z_{i_0}$, the constraint
$c=(i_0,\varnothing,\varnothing)$ has $\pi(C_c)=Z_{i_0}$.
Consequently $\pi(t_\gamma)$ is a quasi-subnet of $z_i$. Conversely,
for every sampling index $(c,t)$, the first inclusion in
\eqref{uc:realization-tail-bounds} places an original tail inside
$\pi(C_c)=T_{(c,t)}((\pi(t_\gamma)))$, so $z_i$ is a quasi-subnet
of $\pi(t_\gamma)$. This proves tail equivalence.
\end{proof}

\begin{proposition}[The sampling test for absorption]
\label{uc:necessity-absorption}
Suppose that the construction conclusion in
Definition~\ref{uc:universal-class} holds for every unit-coefficient
sampling system with an input error base of cardinality at most
$\kappa$. Then $Z$ satisfies \textnormal{$(P_\kappa)$}. In
particular, the universal construction assertion at cardinal $\kappa$
implies \textnormal{$(P_\kappa)$}.
\end{proposition}
\begin{proof}
Let $A_\lambda$ decrease on a directed set $\Lambda$ with
$|\Lambda|\le\kappa$ and $\langle A_\lambda\rangle\to0$.
Let $z_i$ and null nets $h_j^\lambda$ satisfy the additive form of
\eqref{ci:P}. For each $\lambda,j$, that hypothesis supplies the
threshold $i(\lambda,j)$ in
\eqref{uc:realization-hypothesis}. Apply
Lemma~\ref{uc:absorption-realization} to these sets and nets, and use
its sampling net $t_\gamma\in T$ in the single-domain system
\eqref{uc:unit-sampling}. Put
\[
 f_\lambda=\pi_\lambda,\qquad f=\pi,\qquad
 B_\lambda=(-A_\lambda)\cup\{0\}.
\]
By the second assertion in \eqref{uc:realization-conclusions},
\[
 S_\gamma f_\lambda=\pi_\lambda(t_\gamma)\longrightarrow0,
 \qquad If_\lambda=0.
\]
Here the family $B_\lambda$ is decreasing and null by inversion and
\eqref{uc:null-zero-adjunction}, while for $\eta\ge\lambda$ and $t\in T$,
the first assertion in \eqref{uc:realization-conclusions} gives
\[
 f_\eta(t)-f(t)=\pi_\eta(t)-\pi(t)
       \in-A_\eta\subseteq-A_\lambda\subseteq B_\lambda.
\]
Thus the approximation index is $P=\Lambda$, the error index is
$\Lambda$, and its threshold at $\lambda$ is $p_\lambda=\lambda$.
Applying the construction conclusion gives $z\in Z$ with
$\pi(t_\gamma)\to z$ and $If_\lambda=0\to z$; since the constant
zero net also converges to zero, separation gives $z=0$.
Finally, the last assertion in \eqref{uc:realization-conclusions} transfers
$\pi(t_\gamma)\to0$ to $z_i\to0$. This is
\textnormal{$(P_\kappa)$}.
\end{proof}

These sampling tests recover each intrinsic requirement from the part
of the uniform theorem that uses it, and can now be combined with the
sufficiency results to give the two characterizations.

\begin{theorem}[Universal preservation and adherence]
\label{uc:universal-known}
For a separated abelian convergence group $Z$ and an infinite cardinal
$\kappa$, the universal preservation assertion in
Definition~\ref{uc:universal-class} holds if and only if $Z$ satisfies
\textnormal{$(R_\kappa)$}.
\end{theorem}
\begin{proof}
Suppose $Z$ satisfies \textnormal{$(R_\kappa)$}, and fix any of the
systems and approximation data quantified in
Definition~\ref{uc:universal-class}. The coefficient hypothesis is
continuous finite semivariance at cardinal $\kappa$, the input error
sets satisfy \eqref{cu:input-null-base}--\eqref{cu:input-thresholds},
and both $f_p$ and $f$ are integrable on every tested domain.
Theorem~\ref{uc:known-limit}, with this same $X,Z,\mathcal A,P$ and
error family $B_\lambda$, gives
\[
 p\ge p_\lambda\quad\Longrightarrow\quad
 I_Af_p-I_Af\in a_Z(\Sigma(B_\lambda))\qquad(A\in\mathcal A).
\]
Continuous finite semivariance and \textnormal{$(R_\kappa)$} give
\[
 \langle\Sigma(B_\lambda)\rangle\longrightarrow0,
 \qquad\langle a_Z(\Sigma(B_\lambda))\rangle\longrightarrow0.
\]
For every $\lambda$, the endpoint tail beginning at $p_\lambda$ is
contained in $a_Z(\Sigma(B_\lambda))$. Its saturation is a
quasi-subnet of the displayed null saturation, which proves
\eqref{uc:universal-preservation-formula} for the arbitrary system.

Conversely, \eqref{uc:unit-sampling} verifies the coefficient condition
for every unit-coefficient sampling system, so universal preservation
includes every system in Proposition~\ref{uc:necessity-adherence}.
That proposition applied to each decreasing null family with at most
$\kappa$ indices gives \textnormal{$(R_\kappa)$}.
\end{proof}

\begin{theorem}[Universal construction of the limiting integral]
\label{uc:universal-constructive}
For a separated abelian convergence group $Z$ and an infinite cardinal
$\kappa$, the universal construction assertion in
Definition~\ref{uc:universal-class} holds if and only if $Z$ is
$\kappa$-complete and satisfies \textnormal{$(R_\kappa)$} and
\textnormal{$(P_\kappa)$}.
\end{theorem}
\begin{proof}
Assume the three stated properties of $Z$. For an arbitrary system in
Definition~\ref{uc:universal-class}, its coefficient and approximation
assumptions are those of Theorem~\ref{cu:uniform-cardinal}, and its
known integrals are the row limits required there. That theorem gives
$S^A_\gamma f\to I_Af$ on every prescribed $\Gamma_A$ and uniform
convergence $I_Af_p\to I_Af$. Taking $z_A=I_Af$ yields both assertions
in \eqref{uc:universal-construction-formula}.

For the converse, suppose universal construction holds. If $f$ is
already integrable, its prescribed sum net has both the known limit
$I_Af$ and the constructed limit $z_A$. Separation gives
\[
 S^A_\gamma f\longrightarrow I_Af,\ z_A
 \quad\Longrightarrow\quad I_Af=z_A
 \qquad(A\in\mathcal A).
\]
Substituting these equalities into
\eqref{uc:universal-construction-formula} gives
\eqref{uc:universal-preservation-formula}. Universal preservation and
Theorem~\ref{uc:universal-known} then imply
\textnormal{$(R_\kappa)$}. The unit systems satisfy the coefficient
condition by \eqref{uc:unit-sampling}, so construction also has the
hypothesis of Proposition~\ref{uc:necessity-completeness} and of
Proposition~\ref{uc:necessity-absorption}. Their conclusions are
$\kappa$-completeness and \textnormal{$(P_\kappa)$}, respectively.
\end{proof}

\begin{corollary}[Sequential and unrestricted characterizations]
\label{uc:universal-sequential}
For sequences of integrable functions and arbitrary prescribed
integration nets, universal parameter-uniform preservation is equivalent
to \textnormal{$(R_\omega)$}, and universal parameter-uniform
construction is equivalent to sequential completeness together with
\textnormal{$(R_\omega)$} and \textnormal{$(P_\omega)$}. The
coefficient class in these assertions has continuous finite semivariance
at cardinal $\aleph_0$.

For arbitrary approximation nets and the coefficient class with
unrestricted continuous finite semivariance, the corresponding
characterizations are \textnormal{(R)}, and completeness together with
\textnormal{(R)} and \textnormal{(P)}, respectively.
\end{corollary}
\begin{proof}
We first prove the two sequential assertions. Suppose
$f_n\to f$ uniformly on $\Omega$, and put
\begin{equation}\label{uc:sequential-actual-errors}
 B_N=\{0\}\cup\{f_n(t)-f(t):n\ge N,\ t\in\Omega\}.
\end{equation}
If $\Omega=\varnothing$, then $B_N=\{0\}$. Otherwise uniform
convergence makes the saturation of the error tails in
\eqref{uc:sequential-actual-errors} null, and
\eqref{uc:null-zero-adjunction} gives $\langle B_N\rangle\to0$.
These sets decrease, have index $\mathbb N$, and satisfy
\[
 n\ge N\quad\Longrightarrow\quad
 f_n(t)-f(t)\in B_N\qquad(t\in\Omega).
\]
Under \textnormal{$(R_\omega)$}, Theorem~\ref{uc:known-limit}
applies with $\kappa=\aleph_0$, $\Lambda=P=\mathbb N$ and
$p_N=N$ whenever $f$ is integrable, proving sequential preservation.
For construction, sequential completeness gives
$\aleph_0$-completeness by Lemma~\ref{ci:countable-completeness}.
Together with \textnormal{$(R_\omega)$} and
\textnormal{$(P_\omega)$}, the same error sets and thresholds satisfy
all assumptions of Theorem~\ref{cu:uniform-cardinal}, which constructs
the integrals of $f$ and proves their uniform convergence.

For the converse, we give sampling tests whose approximation index is
$\mathbb N$. We will need an explicit reduction of a nonempty
countable directed set $\Lambda$. Enumerate it as
$(\eta_n)_{n\ge1}$, with repetitions if it is finite. Fix a map
$v:\Lambda\times\Lambda\to\Lambda$ satisfying
$v(\lambda,\eta)\ge\lambda,\eta$, and recursively define
\begin{equation}\label{uc:countable-cofinal-indices}
 \lambda_1=\eta_1,\qquad
 \lambda_{n+1}=v(\lambda_n,\eta_{n+1}).
\end{equation}
Induction gives
\[
 \lambda_{n+1}\ge\lambda_n,\qquad
 \lambda_n\ge\eta_j\quad(1\le j\le n).
\]
Thus for each $\lambda\in\Lambda$ there is $j$ with
$\eta_j=\lambda$, and $\lambda_n\ge\lambda$ for every $n\ge j$.
For a decreasing family $A_\lambda$, put $D_n=A_{\lambda_n}$.
Then $D_{n+1}\subseteq D_n$, every $D_n$ is a member of the original
family, and
\begin{equation}\label{uc:cofinal-error-families}
 \forall\lambda\in\Lambda\ \exists j\ \forall n\ge j:
 D_n\subseteq A_\lambda,\qquad
 a_Z(D_n)\subseteq a_Z(A_\lambda).
\end{equation}
For each original level $\lambda$, choose its enumerating position
$j$ with $\eta_j=\lambda$. The two families then satisfy
\[
 D_j\subseteq A_\lambda,\qquad
 a_Z(D_j)\subseteq a_Z(A_\lambda).
\]
In the other direction, each sequential level $n$ has the original
level $\lambda_n$ with
\[
 A_{\lambda_n}=D_n,\qquad
 a_Z(A_{\lambda_n})=a_Z(D_n).
\]
Thus every tail of $\langle A_\lambda\rangle$ contains a tail of
$\langle D_n\rangle$ and every tail of $\langle D_n\rangle$ is
an original tail, proving their tail equivalence. The second inclusion
and equality give these two containments for the adherence families
as well, so $\langle a_Z(D_n)\rangle$ and
$\langle a_Z(A_\lambda)\rangle$ are tail equivalent.

Assume sequential preservation, and let $A_\lambda$ be an arbitrary
nonempty decreasing null family on a countable directed index.
The sets $D_n$ above have null saturation. Define
\[
 Q=\{(n,y):n\ge1, y\in a_Z(D_n)\}.
\]
For each $q=(n,y)$, select an admissible witnessing net
$x_j^q\in D_n$ converging to $y$ and form the disjoint sampling
domains $\Omega_q=\{q\}\times J_q$ of
\eqref{uc:unit-sampling}. The representatives can be taken from the
set of saturations specified in
\eqref{uc:adherence-witnesses}. On their disjoint union define
\[
 g_m(q,j)=
 \begin{cases}x_j^q,&q=(m,y)\text{ for some }y,\\0,&q=(n,y),\ n\ne m,
 \end{cases}
 \qquad g(q,j)=0.
\]
For every $q=(n,y)$, the sampling sums give
\begin{equation}\label{uc:sequential-adherence-integrals}
 I_{\Omega_q}g_m=
 \begin{cases}y,&m=n,\\0,&m\ne n,
 \end{cases}
 \qquad I_{\Omega_q}g=0.
\end{equation}
For $m\ge N$ and every sampling point, its error is either zero or
an element of $D_m\subseteq D_N$. Hence
\[
 g_m(q,j)-g(q,j)\in D_N\cup\{0\}\qquad(m\ge N),
 \qquad \langle D_N\cup\{0\}\rangle\longrightarrow0.
\]
This proves the required uniform convergence of the sequence $g_m$.
Every endpoint in the tail $m\ge N$ belongs to
$a_Z(D_N)\cup\{0\}$ by \eqref{uc:sequential-adherence-integrals}
and monotonicity. Conversely, $y\in a_Z(D_N)$ is attained at
$m=N$, $q=(N,y)$, while zero is attained at $m=N$ and
$q=(N+1,a)$ for any $a\in D_{N+1}$. Therefore
\[
 \{I_{\Omega_q}g_m:m\ge N, q\in Q\}=a_Z(D_N)\cup\{0\}.
\]
Sequential preservation makes the saturation of these sets null.
Their contained family $a_Z(D_N)$ is null as a quasi-subnet, and
\eqref{uc:cofinal-error-families} transfers this convergence to
$\langle a_Z(A_\lambda)\rangle$. Thus \textnormal{$(R_\omega)$}
is necessary.

Assume now sequential construction. It implies sequential preservation
because the known and constructed integral of the same sum net agree
by separation, so the preceding test gives
\textnormal{$(R_\omega)$}. To obtain sequential completeness, let
$(x_i)_{i\ge1}$ be a Cauchy sequence. Use $\Omega=\mathbb N$ with
identity sampling and define
\[
 g_n(i)=\begin{cases}x_n,&i\ge n,\\x_i,&i<n,\end{cases}
 \qquad g(i)=x_i,\qquad
 B_N=\{x_n-x_i:n,i\ge N\}.
\]
The sum $S_i g_n=x_n$ for $i\ge n$ proves $Ig_n=x_n$.
The sets $B_N$ are the diagonal tails of the Cauchy difference net
$(x_n-x_i)_{(n,i)\in\mathbb N^2}$, and these tails are cofinal:
the tail at $(n_0,i_0)$ contains $B_{\max(n_0,i_0)}$.
Consequently $\langle B_N\rangle\to0$, with $0\in B_N$.
For every $n\ge N$ and $i\in\mathbb N$,
\[
 g_n(i)-g(i)=
 \begin{cases}x_n-x_i\in B_N,&i\ge n,\\0\in B_N,&i<n.
 \end{cases}
\]
Sequential construction therefore makes $g$ integrable. Its sampling
net is $S_i g=x_i$, so the original Cauchy sequence converges.

To obtain \textnormal{$(P_\omega)$}, let $A_\lambda,z_i,h_j^\lambda$
satisfy its defining assumptions, with $\Lambda$ countable.
Use \eqref{uc:countable-cofinal-indices} and set
\[
 D_n=A_{\lambda_n},\qquad
 K_n=J_{\lambda_n},\qquad
 k_j^n=h_j^{\lambda_n}\quad(j\in K_n).
\]
Then $\langle D_n\rangle\to0$, each $k^n$ is null, and the
original absorption antecedent gives
\[
 \forall n\ \forall j\in K_n\ \exists i(n,j)\ \forall i\ge i(n,j):
 \quad z_i\in D_n+\{k_l^n:l\ge j\}.
\]
Apply Lemma~\ref{uc:absorption-realization} with index set
$\mathbb N$, sets $D_n$ and rows $k^n$. Its underlying set and maps are
\[
 T=\{(x,(y_n))\in Z\times Z^{\mathbb N}:x-y_n\in D_n\ (n\ge1)\},
 \qquad \pi(x,y)=x,\quad\pi_n(x,y)=y_n.
\]
For its sampling net $t_\gamma$, the lemma proves
$\pi_n(t_\gamma)\to0$ and tail equivalence of
$\pi(t_\gamma)$ with $z_i$. Set $g_n=\pi_n$ and $g=\pi$ on the
single sampling domain $T$. Their finite sums satisfy
\[
 S_\gamma g_n=\pi_n(t_\gamma)\longrightarrow0,\qquad Ig_n=0,
\]
and, for $n\ge N$ and $t\in T$,
\[
 g_n(t)-g(t)\in-D_n\subseteq-D_N
       \subseteq(-D_N)\cup\{0\}.
\]
Inversion and \eqref{uc:null-zero-adjunction} make the last error
family null. Sequential construction supplies $z\in Z$ such that
$S_\gamma g=\pi(t_\gamma)\to z$ and $Ig_n=0\to z$.
Separation gives $z=0$, and the proved tail equivalence then gives
$z_i\to0$. This establishes \textnormal{$(P_\omega)$} without
restricting $I$, any $J_\lambda$, or the sampling index $\Gamma$.

We finally prove the unrestricted assertions. For a uniformly convergent
net $f_p\to f$, put
\[
 B_{p_0}=\{0\}\cup
       \{f_p(t)-f(t):p\ge p_0, t\in\Omega\},\qquad
 \kappa=\max(\aleph_0,|P|).
\]
The family is decreasing and null by uniform convergence and
\eqref{uc:null-zero-adjunction}, and its threshold at $p_0$ is
$p_0$ itself. Unrestricted finite semivariance implies finite
semivariance at $\kappa$. Under \textnormal{(R)}, its
$\kappa$-restricted form permits Theorem~\ref{uc:known-limit} to
prove preservation. If $Z$ is complete and satisfies
\textnormal{(R)} and \textnormal{(P)}, it is $\kappa$-complete and
satisfies \textnormal{$(R_\kappa)$} and
\textnormal{$(P_\kappa)$}, so
Theorem~\ref{cu:uniform-cardinal} proves construction for these
same data.

Conversely, fix an arbitrary infinite cardinal $\kappa$. Every
unit-coefficient sampling system has unrestricted continuous finite
semivariance by \eqref{uc:unit-sampling}. Hence unrestricted
preservation includes the unit systems with error index of cardinality
at most $\kappa$, and Proposition~\ref{uc:necessity-adherence} gives
\textnormal{$(R_\kappa)$}. If unrestricted construction holds, it
includes the construction conclusion for those unit systems as well.
Propositions~\ref{uc:necessity-completeness} and
\ref{uc:necessity-absorption} give $\kappa$-completeness and
\textnormal{$(P_\kappa)$}; preservation, obtained by separation,
gives \textnormal{$(R_\kappa)$}. For any particular directed error
index $\Lambda$, the cardinal
$\max(\aleph_0,|\Lambda|)$ is among these cardinals, proving
unrestricted \textnormal{(R)} and \textnormal{(P)}. For any Cauchy
net on a directed set $I$, the cardinal $\max(\aleph_0,|I|)$ is
among them, so that net converges. This proves completeness.
\end{proof}

Through these equivalences, the prescribed-record class, including
the parameterized sampling systems used in the converse arguments,
determines precisely which operations universal uniform passage
requires. Unit coefficients suffice in each necessity proof because
the sampling points and tested domains retain the information in the
arbitrary error family. We may therefore return to the concrete
integration theorems and derive these operations from their independent
assumptions on the original convergence and the actual finite
coefficients.

To see how the group hypotheses operate beyond the abelian setting,
we first place them in an example governed by coordinatewise order
convergence, where both the
noncommutative product and the failure of topologicity can be followed
through the coordinates.

\begin{example}[A complete non-topological nonabelian target]
\label{ci:heisenberg}
Let $E=\mathbb R^{\mathbb N}$, with coordinatewise algebraic
operations and order convergence, and let $\mathbf1=(1,1,\ldots)$.
On $G=E^3$ define
\begin{equation}\label{ci:heisenberg-product}
 (a,b,c)(a',b',c')=(a+a',b+b',c+c'+ab').
\end{equation}
The identity is $(0,0,0)$, and
\begin{equation}\label{ci:heisenberg-inverse}
 (a,b,c)^{-1}=(-a,-b,-c+ab).
\end{equation}
Equip $G$ with order convergence in each of its three coordinates.
This is a separated, right complete convergence group satisfying
\textnormal{(CB)}, \textnormal{(R)}, and
\textnormal{($P_\times$)}. Its convergence is not induced by a
topology.

For $g=(a,b,c)$, $h=(a',b',c')$ and
$k=(a'',b'',c'')$, both $(gh)k$ and $g(hk)$ have first two
coordinates $a+a'+a''$ and $b+b'+b''$. Their third coordinates are,
respectively,
\[
 \begin{aligned}
 (c+c'+ab')+c''+(a+a')b''
 &=c+c'+c''+ab'+ab''+a'b'',\\
 c+(c'+c''+a'b'')+a(b'+b'')
 &=c+c'+c''+ab'+ab''+a'b''.
 \end{aligned}
\]
Hence $(gh)k=g(hk)$. Formula \eqref{ci:heisenberg-product} gives
$g(0,0,0)=(0,0,0)g=g$, and
\[
 \begin{aligned}
 (a,b,c)(-a,-b,-c+ab)&=(0,0,c-c+ab-ab)=(0,0,0),\\
 (-a,-b,-c+ab)(a,b,c)&=(0,0,-c+ab+c-ab)=(0,0,0).
 \end{aligned}
\]
With the identity and inverse formulas verified, noncommutativity is
witnessed by
\[
 (\mathbf1,0,0)(0,\mathbf1,0)
       =(\mathbf1,\mathbf1,\mathbf1),\qquad
 (0,\mathbf1,0)(\mathbf1,0,0)
       =(\mathbf1,\mathbf1,0).
\]

If $u_i\to_o u$ and $v_i\to_o v$ in $E$, a tail of $(v_i)$
satisfies $|v_i|\le w$ for some $w\in E_+$. On that tail,
\[
 |u_iv_i-uv|
 \le |u_i-u|w+|u|\,|v_i-v|.
\]
Let $d_\alpha\downarrow0$ and $e_\beta\downarrow0$ control
$u_i-u$ and $v_i-v$, respectively. The family
\[
 r_{\alpha,\beta}=wd_\alpha+|u|e_\beta
\]
is decreasing on the product control index. At every coordinate $n$,
\[
 \inf_{\alpha,\beta}r_{\alpha,\beta}(n)
 =w(n)\inf_\alpha d_\alpha(n)
   +|u(n)|\inf_\beta e_\beta(n)=0.
\]
For a fixed $(\alpha,\beta)$, intersect the tail on which
$|v_i|\le w$ with the two tails giving
$|u_i-u|\le d_\alpha$ and $|v_i-v|\le e_\beta$.
On this common tail, $|u_iv_i-uv|\le r_{\alpha,\beta}$,
which proves $u_iv_i\to_o uv$. For two independently indexed
factors, the same bounds hold on their product directed set.
The inequalities
\[
 |(u_i+v_j)-(u+v)|\le d_\alpha+e_\beta,
 \qquad |-u_i+u|=|u_i-u|
\]
prove continuity of addition and negation. Each coordinate in
\eqref{ci:heisenberg-product} and
\eqref{ci:heisenberg-inverse} is a finite composition of these
continuous operations, so multiplication and inversion on $G$ are
continuous. If a net in $G$ has two limits, every real coordinate
has both corresponding coordinate limits. Uniqueness of limits in
$\mathbb R$ identifies the two points, proving separation.

For $g=(a,b,c)$, define
\[
 N(g)=|a|\vee|b|\vee|c|,\qquad
 B_r=\{g:N(g)\le r\}\quad(r\in E_+).
\]
Writing $h=(a',b',c')$, the three coordinates satisfy
\[
 \begin{aligned}
 |a+a'|&\le N(g)+N(h),\\
 |b+b'|&\le N(g)+N(h),\\
 |c+c'+ab'|&\le N(g)+N(h)+N(g)N(h).
 \end{aligned}
\]
Their supremum therefore satisfies
\begin{equation}\label{ci:heisenberg-radius}
 N(gh)\le N(g)+N(h)+N(g)N(h).
\end{equation}
A net $(g_i)$ converges to the identity exactly when there is a
nonempty downward directed $D\subseteq E_+$ with $\inf D=0$ such that
$g_i\in B_r$ eventually for every $r\in D$. Indeed, if the three coordinate nets are controlled by
$d_\alpha\downarrow0$, $e_\beta\downarrow0$ and
$f_\theta\downarrow0$, use
\[
 D=\{d_\alpha+e_\beta+f_\theta:(\alpha,\beta,\theta)\}.
\]
It is downward directed, and its infimum is zero at every real
coordinate. The intersection of the three witnessing tails gives
$N(g_i)\le d_\alpha+e_\beta+f_\theta$. Conversely,
$N(g_i)\le r$ gives all three coordinate inequalities with the
same bound $r$, so $D$ regulates each coordinate. Each $B_r$ is convergence-closed, since a net in $B_r$
converging to $g$ has, at every real coordinate, a limit in the closed
interval $[-r(n),r(n)]$. Thus $N(g)\le r$.

Let $A_\lambda\downarrow$ be an originally null family. Its controls
$D\downarrow0$ satisfy, for every $r\in D$,
\[
 A_\lambda\subseteq B_r\quad\text{eventually in }\lambda,
 \qquad K_G(A_\lambda)\subseteq K_G(B_r)=B_r
                  \quad\text{on the same tail}.
\]
These inequalities give order-null controls for each coordinate of
$\langle K_G(A_\lambda)\rangle$, proving its nullity. If
$g_i\to g$, apply this assertion to
$A_i=T_i((g_jg^{-1})_j)=T_i((g_j)_j)g^{-1}$.
Right translation and its inverse preserve closed sets, and hence
\[
 K_G(A_i)=K_G(T_i((g_j)_j))g^{-1}.
\]
The nullity just proved and translation by $g$ give
$\langle K_G(T_i((g_j)_j))\rangle\to g$, which is
\textnormal{(CB)}. Finally
$a_G(A_\lambda)\subseteq K_G(A_\lambda)$ makes the adherence
saturation a quasi-subnet of the null closed-hull saturation, proving
\textnormal{(R)}.

For \textnormal{($P_\times$)}, let $A_\lambda\downarrow$ be null,
let $(z_i)$ be a candidate net, and suppose that for every $\lambda$
there is a null net $(h^\lambda_t)$ such that
\[
 \forall t_0:\quad
 z_i\in A_\lambda\{h^\lambda_t:t\ge t_0\}
                       \quad\text{eventually in }i.
\]
Let $D\downarrow0$ control $A_\lambda$. For each $r\in D$, fix an
index $\lambda(r)$ with $A_{\lambda(r)}\subseteq B_r$, and let
$E_r\downarrow0$ control $h^{\lambda(r)}$. For $s\in E_r$ a tail of
this inner net lies in $B_s$, so the candidate net eventually belongs
to
\[
 A_{\lambda(r)}B_s\subseteq B_rB_s\subseteq B_{r+s+rs}.
\]
Consequently
\[
 \mathscr R=\{u\in E_+:z_i\in B_u\text{ eventually}\}
\]
is nonempty. It is closed under finite meets because a common upper
index for two witnessing tails gives membership in
$B_u\cap B_v=B_{u\wedge v}$. Every lower bound $w\in E$ of
$\mathscr R$ satisfies
\[
 w(n)\le r(n)+(1+r(n))s(n)
       \qquad(r\in D,\ s\in E_r,\ n\in\mathbb N).
\]
For fixed $r,n$, the infimum over $s\in E_r$ is $r(n)$, and the
subsequent infimum over $r\in D$ is zero. Thus every lower bound of $\mathscr R$ is at most zero, while zero
is itself a lower bound. Since $\mathscr R$ is nonempty and closed
under finite meets, $\mathscr R\downarrow0$. By its definition,
for every $u\in\mathscr R$ there is an index $i_u$ such that
$N(z_i)\le u$ for all $i\ge i_u$. This is the defining order-null
control for all three coordinates of $z_i$, proving $z_i\to e$
and \textnormal{($P_\times$)}.

To verify right completeness, first let $(u_i)$ be order Cauchy
in $E$. For a coordinate $n$ and $\varepsilon>0$, an order-Cauchy control
$D\downarrow0$ contains $d$ with $d(n)<\varepsilon$.
The tail inequality $|u_i-u_j|\le d$ then gives
$|u_i(n)-u_j(n)|<\varepsilon$ on a square tail. Hence each real
coordinate net is Cauchy and has a limit $u(n)$. If $d\in D\downarrow0$ regulates the differences, then
some $i_d$ satisfies
\[
 |u_i-u_j|\le d\qquad(i,j\ge i_d).
\]
For fixed $i\ge i_d$ and $n$, passage to the real limit in $j$
gives $|u_i(n)-u(n)|\le d(n)$. Thus $|u_i-u|\le d$ on the
same tail, and $u_i\to_o u$.
Now let $g_i=(a_i,b_i,c_i)$ be right Cauchy in $G$. A direct
calculation gives
\begin{equation}\label{ci:heisenberg-cauchy}
 g_ig_j^{-1}
 =\bigl(a_i-a_j,b_i-b_j,
             c_i-c_j+(a_j-a_i)b_j\bigr).
\end{equation}
The first two coordinates are order null on the product index.
The preceding argument gives $a_i\to_o a$ and $b_i\to_o b$.
There are $j_0$ and $w\in E_+$ with $|b_j|\le w$ for
$j\ge j_0$. Hence
\[
 |(a_i-a_j)b_j|\le |a_i-a_j|w\longrightarrow_o0.
\]
The third coordinate of \eqref{ci:heisenberg-cauchy} is order
null, so
\[
 c_i-c_j
 =\bigl[c_i-c_j+(a_j-a_i)b_j\bigr]+(a_i-a_j)b_j
 \longrightarrow_o0.
\]
There is $c\in E$ with $c_i\to_o c$. Therefore
$g_i\to(a,b,c)$ in the original convergence.

To complete the example, we examine whether its convergence can come
from a topology. A sequence in $E$ converges in order exactly when it
converges coordinatewise, with the nontrivial implication witnessed,
for a coordinatewise null sequence $(u_m)$, by
\[
 d_k(n)=\sup_{m\ge k}|u_m(n)|<\infty,
 \qquad d_k\downarrow0,
 \qquad |u_m|\le d_k\quad(m\ge k).
\]
Every order-convergent net converges coordinatewise. Consequently,
every set closed in the usual product topology on
$\mathbb R^{\mathbb N}$ is order-closed. Conversely, if $C$ is
order-closed and $x$ is in its product closure, for each positive
integer $m$ the set
\[
 \{u\in C:|u(n)-x(n)|<1/m\text{ for }1\le n\le m\}
\]
is nonempty. A sequence taking one element from each of these sets
converges coordinatewise, hence in order, to $x$. Thus $x\in C$.
Thus the order-closed sets are exactly the product-closed sets,
which forces any topology inducing order convergence to equal the
product topology.

For finite $F\subseteq\mathbb N$, put
\[
 V_F=\{u\in E:u(n)=0\text{ for every }n\in F\}.
\]
The saturated value net of $V_F$, indexed by pairs $(F,u)$ with
$u\in V_F$ and ordered by inclusion of the finite first coordinates,
converges coordinatewise to zero. Its tail above $(F,u)$ has value set
exactly $V_F$, and this set is not order bounded because
$\{v(n):v\in V_F\}=\mathbb R$ for every $n\notin F$.
Every order-null net has an order-bounded tail, so this saturated net
is not order null. This contradicts equality with product convergence.
Finally, the subgroup
$\{(a,0,0):a\in E\}$ has exactly this induced convergence.
If the convergence of $G$ were topological, its restriction to
that subgroup would be topological. Therefore the convergence
of $G$ is not topological.
\end{example}

Within this group, the dependence on the right uniform structure
also becomes visible at the level of the conclusion: the cancellation
in a right endpoint error need not survive when that error is written
in the opposite order or as a coordinate difference.

\begin{example}[Right uniform limits beyond additive coordinate limits]
\label{ci:strict-uniformity}
In the group of Example~\ref{ci:heisenberg}, for positive integers
$n,p$ put
\[
 z^n=(0,n\mathbf1,0),\qquad
 y_p^n=
 \left(-\frac1p\mathbf1,n\mathbf1,-\frac np\mathbf1\right).
\]
Let $x_\gamma^n=z^n$ and $x_{p,\gamma}^n=y_p^n$ on any
nonempty directed inner index. The rows are constant and converge
to their stated endpoints. Formula
\eqref{ci:heisenberg-inverse} gives
\[
 (y_p^n)^{-1}=\left(\frac1p\mathbf1,-n\mathbf1,0\right),
 \qquad
 x_\gamma^n(x_{p,\gamma}^n)^{-1}
       =\left(\frac1p\mathbf1,0,0\right).
\]
Thus \eqref{ci:comparison} holds with
\[
 Q_k=\{(t\mathbf1,0,0):0\le t\le1/k\},\qquad
 P_k=\{p:p\ge k\}.
\]
The balls $B_{\mathbf1/k}$ contain $Q_k$, so the saturated value
net of $(Q_k)$ converges to $e$. The right endpoint errors are
\[
 y_p^n(z^n)^{-1}
 =\left(-\frac1p\mathbf1,0,0\right),
\]
which converge to $e$ uniformly in $n$.
In contrast, the left endpoint errors and the additive coordinate
errors are both
\[
 (z^n)^{-1}y_p^n
 =y_p^n-z^n
 =\left(-\frac1p\mathbf1,0,-\frac np\mathbf1\right).
\]
For every $p_0$, the third coordinates in the uniform tail
$\{(z^n)^{-1}y_p^n:p\ge p_0,\ n\ge1\}$ are unbounded even at
the first coordinate of $E$. No such tail is order-bounded.
The saturated error net cannot be null, since every null net in $G$
has a tail contained in an order-bounded coordinate box. Consequently left uniform
convergence fails, and additive uniform convergence in $E^3$
also fails. The right uniform conclusion of
Theorem~\ref{ci:interchange} therefore cannot be recovered by
applying the abelian theorem to the additive coordinate group.
\end{example}

% BEGIN CARDINAL EXAMPLES: ce:
\subsection{Cardinal completeness and the absence of a complete carrier}

Returning to the size of the approximation, we recall that the
cardinal restriction in Theorem~\ref{ci:interchange} belongs to the
error base, leaving the number of approximants and parameters
unrestricted. This distinction raises the question whether the weaker
completeness assumption merely conceals a complete subspace containing
the endpoints; the following construction separates completeness at
the stated cardinal from both completeness of the whole range and the
existence of such a subspace.

\begin{example}[Restricted supports and completeness levels]
\label{ce:restricted-support}
Let $\kappa$ be an infinite cardinal, let $|I|>\kappa$, and put
\[
 E=\mathbb R^{\mathbb N},\qquad
 \mathbf e=(1,1,\ldots)\in E,
\]
\[
 X_\kappa=\left\{x\in E^I:
       |\operatorname{supp}x|\le\kappa\right\},
 \qquad
 \operatorname{supp}x=\{i\in I:x_i\ne0\}.
\]
Equip $E$ with order convergence and $X_\kappa$ with coordinatewise
order convergence:
\begin{equation}\label{ce:original-convergence}
 x_\alpha\longrightarrow x
 \quad\Longleftrightarrow\quad
 |x_{\alpha,i}-x_i|\longrightarrow_o0\text{ in }E
 \quad(i\in I).
\end{equation}
Then $X_\kappa$ is a separated abelian convergence group satisfying
\textnormal{(CB)}, \textnormal{(R)} and \textnormal{(P)}. It is
$\kappa$-complete, but it is not $\kappa^+$-complete. Its original
convergence is not topological. There are constant known rows with a
countable null error base for which no complete subspace of $X_\kappa$
contains even one tail of the row endpoints.
\end{example}

\begin{proof}
The set $X_\kappa$ is an additive subgroup of $E^I$: its zero
vector has empty support, inversion preserves support, and
\[
 \operatorname{supp}(x+y)
 \subseteq\operatorname{supp}x\cup\operatorname{supp}y,
 \qquad
 |\operatorname{supp}x\cup\operatorname{supp}y|
 \le\kappa+\kappa=\kappa.
\]
The order on $E$ is coordinatewise. Every nonempty order-bounded subset
$H\subset E$ has supremum
\[
 (\sup H)(m)=\sup_{h\in H}h(m)\qquad(m\in\mathbb N).
\]
For positive decreasing controls, $u_\theta\downarrow0$ in $E$ means
\[
 u_\theta(m)\downarrow0\qquad(m\in\mathbb N).
\]
For a net $(x_\alpha)$ in $X_\kappa$, nullity in
\eqref{ce:original-convergence} means that for every $i\in I$ there
are controls $u^i_\theta\downarrow0$ such that
\begin{equation}\label{ce:ball-regulators}
 x_\alpha\in B_i(u^i_\theta)\text{ eventually for every }\theta,
 \qquad B_i(u)=\{x\in X_\kappa:|x_i|\le u\}.
\end{equation}
The constant zero net has the zero control, and every quasi-subnet
retains each eventual containment. If $u_\theta,v_\eta\downarrow0$,
then coordinatewise $u_\theta+v_\eta\downarrow0$, and
\[
 \begin{gathered}
 B_i(u_\theta)\cup B_i(v_\eta)
       \subseteq B_i(u_\theta+v_\eta),\\
 B_i(u_\theta)+B_i(v_\eta)
       \subseteq B_i(u_\theta+v_\eta),\qquad
 -B_i(u_\theta)=B_i(u_\theta).
 \end{gathered}
\]
On common tails these inclusions verify finite mixing and continuity
of addition and inversion, with independent factors lifted to their
product directed set. Translation defines convergence at every point.
If $x_\alpha$ converges to both $x$ and $y$, its real coordinate nets
converge to both $x_i(m)$ and $y_i(m)$. Uniqueness of real limits gives
$x_i(m)=y_i(m)$ for every $i,m$, so $x=y$ and the convergence is
separated. The coordinate bounds $\rho_i(x)=|x_i|\in E_+$ are
symmetric, subadditive and detect exactly the original null nets by
\eqref{ce:ball-regulators}.

Each map $x\mapsto x_i(m)$ preserves original convergence. If a net
$x_\alpha\in B_i(u)$ converges to $x$, closedness of the real
interval $[-u(m),u(m)]$ gives $|x_i(m)|\le u(m)$ for every $m$.
Thus $x\in B_i(u)$ and every $B_i(u)$ is convergence-closed. For a
decreasing null family $A_\lambda$, its coordinate controls give
\[
 A_\lambda\subseteq B_i(u^i_\theta)\text{ eventually},\qquad
 K(A_\lambda)\subseteq B_i(u^i_\theta)\text{ on the same tail}.
\]
For every $i$, the same decreasing controls $u^i_\theta$ therefore
bound the $i$th coordinate of $\langle K(A_\lambda)\rangle$,
which proves its original nullity. For a convergent net $x_\alpha\to x$,
apply this assertion to $A_\alpha=T_\alpha((x_\beta)_\beta)-x$.
Translation and its inverse preserve closed sets, so
\[
 K(A_\alpha)=K(T_\alpha((x_\beta)_\beta))-x.
\]
Translating the null closed-hull saturation back by $x$ proves
\textnormal{(CB)} at $x$. The inclusions
$a_{X_\kappa}(A_\lambda)\subseteq K(A_\lambda)$ make the
adherence saturation a quasi-subnet of the null closed-hull
saturation, proving \textnormal{(R)}.

For \textnormal{(P)}, let $A_\lambda\downarrow$ be null and let
$(z_d)$ be a net for which each $\lambda$ admits a null inner net
$(h^\lambda_t)$ with
\[
 \forall t_0:\quad z_d\in A_\lambda+
       \{h^\lambda_t:t\ge t_0\}\quad\text{eventually in }d.
\]
Fix $i\in I$ and let $u_\theta\downarrow0$ control the $i$th
coordinate of the outer family. For each $\theta$, there is
$\lambda(\theta)$ such that
$A_{\lambda(\theta)}\subseteq B_i(u_\theta)$. Let
$v_{\theta,\eta}\downarrow_\eta0$ control the $i$th coordinate of
$h^{\lambda(\theta)}$. For every $\eta$, an inner tail lies in
$B_i(v_{\theta,\eta})$, so the candidate net eventually lies in
\[
 A_{\lambda(\theta)}+B_i(v_{\theta,\eta})
       \subseteq B_i(u_\theta+v_{\theta,\eta}).
\]
The attained-bound set
\[
 D_i=\{w\in E_+:z_d\in B_i(w)\text{ eventually}\}
\]
is nonempty and closed under finite meets, because two witnessing
tails have a common upper index and
$B_i(w)\cap B_i(w')=B_i(w\wedge w')$. If $c\in E$ is a lower
bound of $D_i$, then
\[
 c\le u_\theta+v_{\theta,\eta}\quad(\theta,\eta),\qquad
 c\le u_\theta\quad(\theta),\qquad c\le0,
\]
where both passages to the infimum are coordinatewise real
calculations. Thus $D_i\downarrow0$, and its balls control the
$i$th coordinate of $z_d$. This holds for every $i$, proving
$z_d\to0$ and \textnormal{(P)}. In this abelian group the
multiplicative property \textnormal{($P_\times$)} has the same
content.

Let $(x_d)_{d\in D}$ be a Cauchy net in $X_\kappa$ with
$|D|\le\kappa$. The set
\[
 J=\bigcup_{d\in D}\operatorname{supp}x_d
\]
satisfies $|J|\le\kappa\cdot\kappa=\kappa$. For each $i\in I$,
the net $|x_{d,i}-x_{e,i}|$, indexed by $D\times D$, is order
null. Thus there are $u^i_\theta\downarrow0$ and indices
$d(i,\theta)$ such that
\begin{equation}\label{ce:cauchy-control}
 d,e\ge d(i,\theta)
 \quad\Longrightarrow\quad
 |x_{d,i}-x_{e,i}|\le u^i_\theta.
\end{equation}
For each $m\in\mathbb N$, the real net $(x_{d,i}(m))_d$ is
Cauchy. Real completeness supplies $b_i(m)$ with
\[
 x_{d,i}(m)\longrightarrow b_i(m).
\]
Fix $d\ge d(i,\theta)$. Passing to the real limit in the $e$
coordinate of \eqref{ce:cauchy-control} gives
\[
 |x_{d,i}(m)-b_i(m)|\le u^i_\theta(m)
       \qquad(m\in\mathbb N).
\]
Hence $b_i=(b_i(m))_m\in E$ and
\[
 d\ge d(i,\theta)
 \quad\Longrightarrow\quad
 |x_{d,i}-b_i|\le u^i_\theta,
 \qquad x_{d,i}\longrightarrow_o b_i.
\]
If $i\notin J$, all $x_{d,i}$ vanish, and $b_i=0$. Therefore
$b=(b_i)_{i\in I}\in X_\kappa$ and $x_d\to b$ originally.
This proves $\kappa$-completeness.

Let $J_0\subseteq I$ have cardinality $\kappa^+$. For
$F\in\mathcal P_{\mathrm{fin}}(J_0)$, define
\[
 u_F(i)=
 \begin{cases}
 \mathbf e,&i\in F,\\
 0,&i\notin F.
 \end{cases}
\]
The inclusion-directed index set has cardinality $\kappa^+$.
For each $i\in J_0$, the $i$th coordinate equals $\mathbf e$ once
$F\supseteq\{i\}$, and each coordinate outside $J_0$ is zero.
For $i\in J_0$ and $F,H\supseteq\{i\}$, one has
$u_F(i)-u_H(i)=0$; for $i\notin J_0$ the same identity holds
for every $F,H$. Thus every coordinate of the difference net on
$\mathcal P_{\mathrm{fin}}(J_0)^2$ is eventually zero and is
order null, so $(u_F)$ is Cauchy in the original convergence. If $u_F\to u$
in $X_\kappa$, coordinate uniqueness would give
\[
 u_i=\mathbf e\quad(i\in J_0),\qquad
 u_i=0\quad(i\notin J_0).
\]
This contradicts $|\operatorname{supp}u|\le\kappa$.
Consequently $X_\kappa$ is not $\kappa^+$-complete.

Having separated the two completeness levels, we turn to the nature
of the convergence itself. For a sequence $(z_n)$ in $E$,
coordinatewise convergence to $z$ implies
\[
 w_N(m)=\sup_{n\ge N}|z_n(m)-z(m)|<\infty,
 \qquad w_N\downarrow0\text{ in }E,
 \qquad |z_n-z|\le w_N\quad(n\ge N).
\]
Thus coordinatewise convergence and order convergence agree for
sequences in $E$. They do not agree for arbitrary nets.
Indeed, let
\[
 D_0=\{(F,z,n):F\subset\mathbb N\text{ finite},\ z\in E,
                     z|_F=0,\ n\in\mathbb N\},
\]
ordered by equality together with
\[
 (F,z,n)<(F',z',n')
 \quad\Longleftrightarrow\quad F\subseteq F',\ n<n'.
\]
The net $v_{(F,z,n)}=z$ converges coordinatewise to zero. Its tail
above $(F_0,z_0,n_0)$ contains every $z$ vanishing on $F_0$,
because $(F_0,z,n_0+1)$ belongs to that tail. For every $w\in E_+$,
an index $m\notin F_0$ and the vector supported at $m$ with value
$w(m)+1$ give a member of that tail outside $[-w,w]$. No tail is
order bounded. Therefore $(v_d)$ does not converge in order.

Suppose a topology $\tau$ induced order convergence on $E$. The
coordinate maps would be $\tau$-continuous, so $\tau$ would contain
the product topology $\pi$. If a $\tau$-open set $U$ containing
$z$ were not a $\pi$-neighborhood of $z$, then for every $n\ge1$
there would be $z_n\notin U$ with
\[
 |z_n(m)-z(m)|<1/n\qquad(1\le m\le n).
\]
This sequence converges coordinatewise, hence in order, to $z$.
Its $\tau$-convergence contradicts $z_n\notin U$ for all $n$.
Thus $\tau=\pi$, contradicting the net $(v_d)$ above. The
convergence on $E$ is not topological. A fixed coordinate embedding
$E\longrightarrow X_\kappa$ preserves and reflects the specified
convergences. A topology inducing the convergence on $X_\kappa$
would therefore induce one on $E$, which has just been excluded.

Finally, put
\[
 B=\{v\in X_\kappa:0\le v_i\le\mathbf e\ (i\in I)\},
 \qquad
 P=\mathcal P_{\mathrm{fin}}(I)\times(\mathbb N_{\ge1}\times B),
\]
where the first factor has its inclusion order and the second
factor has equality together with $(n,v)<(m,w)$ whenever $n<m$.
This product is a directed partial order: two indices
$(H,n,v)$ and $(K,m,w)$ have the common upper bound
\[
 (H\cup K,\max(n,m)+1,0).
\]
It has no cofinal subset
of cardinality at most $\kappa$. Indeed, if $S\subseteq P$ and
$|S|\le\kappa$, then
\[
 \left|\bigcup_{(H,n,v)\in S}H\right|\le\kappa<|I|.
\]
For an index $i$ outside this union, no element of $S$ dominates
$(\{i\},1,0)$.
Take one parameter $A$, any nonempty inner directed set $\Gamma$,
and define
\[
 y_{(H,n,v)}=n^{-1}v,\qquad
 x_{(H,n,v),\gamma}=n^{-1}v,\qquad x_\gamma=0.
\]
Every known row is constant. The decreasing sets
\[
 Q_n=\{x\in X_\kappa:|x_i|\le n^{-1}\mathbf e\ (i\in I)\}
\]
contain $0$ and have an originally null saturated value net, since
$n^{-1}\mathbf e\downarrow0$ in every coordinate copy of $E$.
For $m\ge n$ and $v\in B$,
\[
 x_\gamma-x_{(H,m,v),\gamma}=-m^{-1}v\in Q_n
       \qquad(H\subseteq I\text{ finite}).
\]
For each $n$, let $P_n$ consist of the indices $(H,m,v)$ with
$m\ge n$. It contains the tail above $(\varnothing,n,0)$, and
the preceding comparison holds for every index in $P_n$ and every
$\gamma$. The target is sequentially complete because it is
$\kappa$-complete with $\kappa\ge\aleph_0$, and it satisfies
\textnormal{(R)} and \textnormal{(P)} by the preceding proofs.
Thus all hypotheses of Theorem~\ref{ci:interchange} hold with
$\Lambda=\mathbb N$ and the displayed constant rows. The resulting
candidate row is identically zero, so its limit is zero by
separation.

Let $p_0=(H_0,n_0,v_0)$ be arbitrary. The endpoint tail above $p_0$
contains $(n_0+1)^{-1}B$. If a complete subspace $C\subseteq
X_\kappa$ contained this tail, it would contain
\[
 (n_0+1)^{-1}u_F
       \qquad(F\in\mathcal P_{\mathrm{fin}}(J_0)).
\]
This net is Cauchy in the restricted uniform convergence of $C$.
A limit in $C$ would, by the inclusion into $X_\kappa$, have
coordinate $(n_0+1)^{-1}\mathbf e$ at every $i\in J_0$. Such a
point does not belong to $X_\kappa$. No complete subspace contains
the tail. Since $p_0$ was arbitrary, no removal of an initial part
of the outer net supplies a complete carrier.
\end{proof}
% END CARDINAL EXAMPLES: ce:

\subsection{Ordered comparisons in the original pair convergence}

Having identified the operations needed by the group argument, we
return to the original bounds through which they can be verified.
Working first with the algebra of those bounds also retains a
comparison theorem for uniform spaces whose points carry no algebraic
operation, because the hypotheses can be checked on finite control
expressions before any of the endpoints sought by the comparison has
been constructed.

\begin{definition}[Ordered monoid controls]\label{mg:controls}
Let $(M,\le,\oplus,0)$ be a commutative monoid with a partial order.
Assume that $0$ is its least element, that $\oplus$ is isotone in both
variables, and that every pair has a meet. Throughout the monoid
interchange results, the following additional law is imposed:
\begin{equation}\label{mg:zero-translation}
 D\subseteq M,\quad D\downarrow0
 \quad\Longrightarrow\quad
 \inf_{d\in D}(r\oplus d)=r\qquad(r\in M).
 \tag{Z}
\end{equation}
The notation $D\downarrow0$ requires $D$ to be nonempty and downward
directed and its infimum in $M$ to exist and equal $0$. A positive
control net $u_\theta\downarrow0$ has the corresponding meaning.
A net $a_\alpha$ in $M$ is order null if
\[
 \exists u_\theta\downarrow0\quad
 \forall\theta\ \exists\alpha_\theta\ \forall\alpha\ge\alpha_\theta:
 a_\alpha\le u_\theta.
\]
For a saturated family of values, this criterion means that every
control bounds all the values on one tail of its saturation net.
\end{definition}

\begin{lemma}[Finite control calculus]\label{mg:calculus}
Let $M$ satisfy Definition~\ref{mg:controls} and
\eqref{mg:zero-translation}. The following assertions hold.
\begin{enumerate}
\item If $D,E\downarrow0$, then
$D\oplus E=\{d\oplus e:d\in D,e\in E\}\downarrow0$.
For every positive integer $n$, $\{nd:d\in D\}\downarrow0$, where
$nd$ denotes the $n$-fold monoid sum. Finite sums of order-null
nets, on a common index set or on a product index set, are order null.
\item If $a_\alpha\to_o0$ and $c\le b\oplus a_\alpha$
eventually, then $c\le b$.
\item Let $e_\alpha(t)\in M$, where $\alpha$ belongs to a directed
set and $t$ belongs to a nonempty parameter set. Put
\[
 R=\{r\in M:\exists\alpha_r\ \forall\alpha\ge\alpha_r\
                  \forall t:\ e_\alpha(t)\le r\}.
\]
If $R$ is nonempty, it is closed under finite meets. If a nonempty
subset $B\subseteq R$ has infimum $0$, then $R\downarrow0$ and
$e_\alpha(t)\to_o0$ uniformly in $t$.
\item Suppose $a_\theta\downarrow0$ and, for every $\theta$,
$b_{\theta,\eta}\downarrow_\eta0$. If
$a_\theta\oplus b_{\theta,\eta}\in R$ for all $\theta,\eta$,
then $R\downarrow0$.
\end{enumerate}
The assertions apply to saturation nets as well as to the original
indexed nets, with every control attained on an actual tail. No completeness of $M$ is required.
\end{lemma}
\begin{proof}
For $d_1,d_2\in D$ and $e_1,e_2\in E$, downward directedness gives
$d\in D$, $e\in E$ with $d\le d_1,d_2$ and $e\le e_1,e_2$. Thus
$d\oplus e\le d_1\oplus e_1,d_2\oplus e_2$. If $c$ is a lower
bound of $D\oplus E$, then
\[
 c\le d\oplus e\quad(d\in D,e\in E)
 \Longrightarrow c\le\inf_{e\in E}(d\oplus e)=d\quad(d\in D)
 \Longrightarrow c\le0.
\]
Since $0$ is itself a lower bound, the infimum is $0$, and this
binary calculation extends to finite lists as follows. For
$D_1,\ldots,D_n\downarrow0$, define
$S_1=D_1$ and $S_{k+1}=S_k\oplus D_{k+1}$. The binary argument
just proved gives $S_k\downarrow0$ successively for $k=1,\ldots,n$,
and
\[
 S_n=\{d_1\oplus\cdots\oplus d_n:d_k\in D_k\}.
\]
For $d_1,\ldots,d_n\in D$, downward directedness supplies $d\in D$
with $d\le d_i$ for every $i$, and
\[
 nd\le d_1\oplus\cdots\oplus d_n.
\]
The diagonal sums are therefore downward cofinal in the independent
sums and have infimum $0$: any lower bound of the diagonal sums
is a lower bound of $S_n$, and $0$ bounds both families below.
For order-null nets $a^{(k)}_\alpha$ on a common directed set, let
$u^{(k)}_{\theta_k}\downarrow0$ be their controls and let
$\alpha(k,\theta_k)$ be the corresponding thresholds. Given
$\theta=(\theta_1,\ldots,\theta_n)$, take
$\alpha_\theta\ge\alpha(k,\theta_k)$ for every $k$. Then
\[
 \alpha\ge\alpha_\theta\quad\Longrightarrow\quad
 \bigoplus_{k=1}^n a^{(k)}_\alpha
 \le\bigoplus_{k=1}^n u^{(k)}_{\theta_k}.
\]
The right side decreases to zero on the product control index.
For nets on independent sets $I_k$, the product threshold is the
tuple $(\alpha(k,\theta_k))_{k=1}^n\in\prod_k I_k$, and the
identical inequality holds for all tuples above it.

For (2), a regulator $u_\theta\downarrow0$ for $a_\alpha$ gives
$c\le b\oplus u_\theta$ for every $\theta$: the controlling tail
and the assumed comparison tail have nonempty intersection. Hence
\[
 c\le\inf_\theta(b\oplus u_\theta)=b.
\]
For $r,s\in R$, a common upper index for $\alpha_r,\alpha_s$ gives
$e_\alpha(t)\le r\wedge s$ for every parameter on that tail.
Consequently $r\wedge s\in R$. If $B\subseteq R$ has infimum $0$,
a lower bound of $R$ is a lower bound of $B$ and is at most $0$.
Thus $\inf R=0$. Indexing $R$ by reverse order supplies the
regulator asserted in (3).

Under the assumptions of (4), $R$ is nonempty. Every lower bound
$c$ of $R$ satisfies
\[
 c\le a_\theta\oplus b_{\theta,\eta}\quad(\theta,\eta),
 \qquad
 c\le\inf_\eta(a_\theta\oplus b_{\theta,\eta})=a_\theta
       \quad(\theta),
 \qquad c\le\inf_\theta a_\theta=0.
\]
This proves (4). For a saturation net, finitely many eventual bounds hold on
a common tail by directedness, so every displayed intersection
used in the argument is attained by that net.
\end{proof}

\begin{definition}[Faithful monoid pair representation]\label{mg:representation}
Let $(Z,\mathcal U_Z)$ be a separated uniform convergence space.
For each $j\in J$, let $M_j$ satisfy
Definition~\ref{mg:controls}, and let $d_j:Z^2\to M_j$ satisfy
\[
 d_j(z,z)=0,\qquad d_j(z,w)=d_j(w,z),\qquad
 d_j(z,v)\le d_j(z,w)\oplus_j d_j(w,v).
\]
The family is a faithful order representation of the original pair
convergence if, for every admissible pair net,
\begin{equation}\label{uc:faithfulness}
 (z_\alpha,w_\alpha)\longrightarrow_{\mathcal U_Z}\Delta_Z
 \quad\Longleftrightarrow\quad
 d_j(z_\alpha,w_\alpha)\longrightarrow_o0\qquad(j\in J).
\end{equation}
A faithful lattice-group order representation is the special case
$M_j=(F_j)_+$, $\oplus_j=+$, with order infima interpreted in $F_j$.
A faithful regulator representation uses lattice-group controllers
and replaces $o$ by $D$ in the pair equivalence.
\end{definition}

For an abelian lattice-ordered group $F$, its positive cone satisfies
\eqref{mg:zero-translation}. Indeed, if $D\subseteq F_+$ and
$\inf_FD=0$, then
\[
 c\le r+d\quad(d\in D)
 \Longrightarrow c-r\le d\quad(d\in D)
 \Longrightarrow c-r\le0\Longrightarrow c\le r.
\]
The two possible infimum conventions agree here: if $\inf_{F_+}D=0$
and $v\in F$ is a lower bound of $D$, then $v\vee0\in F_+$ is
also a lower bound, so $v\vee0=0$ and $v\le0$. This argument uses
the lattice operation. No equivalence of these infimum conventions
is asserted for a merely partially ordered group.

Products, restrictions to tails, and the saturation nets of
Section~\ref{fd:net-foundations} retain the required admissible
representatives. For a decreasing nonempty family $(B_\eta)$ and a
control $r$, the saturation criterion is
\[
 \langle B_\eta\rangle\text{ is eventually bounded by }r
 \quad\Longleftrightarrow\quad
 (\exists\eta_0)(\forall\eta\ge\eta_0)(\forall b\in B_\eta):
 \ b\le r.
\]
For a parameterized endpoint family, put
\[
 Q_{p_0}=\{(y_p^A,z^A):p\ge p_0,\ A\in\mathcal A\}.
\]
Its pair saturation is originally null exactly when every represented
family
\[
 \{d_j(y_p^A,z^A):p\ge p_0,\ A\in\mathcal A\},\qquad p_0\in P,
\]
has order-null saturation, by \eqref{uc:faithfulness}. The quantifier
over $A$ is therefore uniform on each tail. Completeness below always
refers to the original structure.
Proposition~\ref{lt:represented-completeness} expresses it through the
originally null pair nets, including the sequential case.

\begin{theorem}[Parameter-uniform order interchange]\label{lt:order-transfer}
Suppose $Z$ is complete and has a faithful monoid pair representation
whose controllers satisfy \eqref{mg:zero-translation}. Let $P$
and $\Gamma_A$, $A\in\mathcal A\ne\varnothing$, be directed sets.
Suppose
\begin{equation}\label{lt:known-rows}
 x^A_{p,\gamma}\longrightarrow y_p^A
 \qquad(\gamma\in\Gamma_A),\qquad p\in P.
\end{equation}
Let $x^A_\gamma\in Z$. For every $j$, suppose
$\beta^j_\theta\downarrow0$ and sets $P^j_\theta\subseteq P$
containing tails of $P$ satisfy
\begin{equation}\label{lt:comparison}
 \begin{gathered}
 r^{A,j}_{p,\gamma}\in M_j,\qquad
 r^{A,j}_{p,\gamma}\longrightarrow_o0,\\
 p\in P^j_\theta\quad\Longrightarrow\quad
 d_j(x^A_\gamma,x^A_{p,\gamma})
 \le\beta^j_\theta\oplus_j r^{A,j}_{p,\gamma}
 \quad\text{eventually in }\gamma.
 \end{gathered}
\end{equation}
The sets $P^j_\theta$ and bounds $\beta^j_\theta$ are independent
of $A$. Then unique $z^A\in Z$ satisfy
\begin{equation}\label{lt:transfer-conclusions}
 x^A_\gamma\longrightarrow z^A,\qquad
 y_p^A\longrightarrow z^A\quad\text{uniformly in }A,
\end{equation}
and
\begin{equation}\label{lt:sharp-endpoints}
 p\in P^j_\theta\quad\Longrightarrow\quad
 d_j(y_p^A,z^A)\le\beta^j_\theta.
\end{equation}
No algebraic operation on $Z$ is required.
\end{theorem}
\begin{proof}
Fix $A,j$, $p\in P^j_\theta$, and $q\in P^j_\eta$.
In this proof $\oplus$ denotes $\oplus_j$.
On a common inner tail, symmetry and the triangle inequalities give
\begin{align*}
 d_j(y_p^A,y_q^A)
 &\le d_j(y_p^A,x^A_{p,\gamma})
      \oplus d_j(x^A_{p,\gamma},x^A_\gamma)
      \oplus d_j(x^A_\gamma,x^A_{q,\gamma})
      \oplus d_j(x^A_{q,\gamma},y_q^A)\\
 &\le\beta^j_\theta\oplus\beta^j_\eta
      \oplus r^{A,j}_{p,\gamma}\oplus r^{A,j}_{q,\gamma}
      \oplus d_j(y_p^A,x^A_{p,\gamma})
      \oplus d_j(x^A_{q,\gamma},y_q^A).
\end{align*}
Here point faithfulness makes the last two terms order null, while
the hypothesis gives the same conclusion for the residuals. Their
fourfold sum is therefore order null by Lemma~\ref{mg:calculus}(1),
and part (2) removes it from the fixed endpoint comparison to give
\begin{equation}\label{lt:endpoint-pairs}
 d_j(y_p^A,y_q^A)\le\beta^j_\theta\oplus_j\beta^j_\eta.
\end{equation}
For every $j,\theta$, there exists $p_0(j,\theta)$ such that
\[
 p\ge p_0(j,\theta)\Longrightarrow p\in P^j_\theta.
\]
Therefore
\[
 p,q\ge p_0(j,\theta),\quad A\in\mathcal A
 \quad\Longrightarrow\quad
 d_j(y_p^A,y_q^A)\le\beta^j_\theta\oplus\beta^j_\theta.
\]
Lemma~\ref{mg:calculus}(1) gives
$\beta^j_\theta\oplus\beta^j_\theta\downarrow0$.
For every fixed $A$, pair faithfulness gives
\[
 (y_p^A,y_q^A)_{(p,q)}\longrightarrow_{\mathcal U_Z}\Delta_Z.
\]
Thus $(y_p^A)_p$ is Cauchy, so completeness gives $z^A$ with
$y_p^A\to z^A$, and separation makes $z^A$ unique.

Fix $p\in P^j_\theta$. For each $\eta$, restrict $q$ to a tail
contained in $P^j_\eta$. Then
\[
 d_j(y_p^A,z^A)
 \le\beta^j_\theta\oplus\beta^j_\eta\oplus d_j(y_q^A,z^A).
\]
Point faithfulness and Lemma~\ref{mg:calculus}(2) give
$d_j(y_p^A,z^A)\le\beta^j_\theta\oplus\beta^j_\eta$ for every
$\eta$. Equation~\eqref{mg:zero-translation} gives
\[
 d_j(y_p^A,z^A)
 \le\inf_\eta(\beta^j_\theta\oplus\beta^j_\eta)
 =\beta^j_\theta.
\]

For fixed $A,j$, put
\[
 e_\gamma=d_j(x^A_\gamma,z^A),\qquad
 h_{p,\gamma}=r^{A,j}_{p,\gamma}
                 \oplus d_j(x^A_{p,\gamma},y_p^A),
\]
\[
 T_p=\{u\in M_j:h_{p,\gamma}\le u\text{ eventually}\},\qquad
 R=\{v\in M_j:e_\gamma\le v\text{ eventually}\}.
\]
Each $h_p$ is order null. Hence $T_p$ contains a zero regulator,
is nonempty, and has infimum $0$. Its finite meets belong to $T_p$
by Lemma~\ref{mg:calculus}(3), so $T_p\downarrow0$.
For $p\in P^j_\theta$, the comparison tail and the triangle
inequality give
\[
 \begin{aligned}
 e_\gamma
 &\le d_j(x^A_\gamma,x^A_{p,\gamma})
       \oplus d_j(x^A_{p,\gamma},y_p^A)
       \oplus d_j(y_p^A,z^A)\\
 &\le
 \beta^j_\theta\oplus\beta^j_\theta\oplus h_{p,\gamma}
 \quad\text{eventually}.
 \end{aligned}
\]
Thus $R$ contains
$\beta^j_\theta\oplus\beta^j_\theta\oplus u$ for every
$p\in P^j_\theta$ and $u\in T_p$, and is nonempty. If $c$ is a
lower bound of $R$, then for every $\theta$ and every
$p\in P^j_\theta$,
\[
 c\le\inf_{u\in T_p}
  (\beta^j_\theta\oplus\beta^j_\theta\oplus u)
  =\beta^j_\theta\oplus\beta^j_\theta.
\]
The right-hand side decreases to $0$, so $c\le0$.
Lemma~\ref{mg:calculus}(3) gives $e_\gamma\to_o0$.
Point faithfulness yields $x^A_\gamma\to z^A$.
Finally, \eqref{lt:sharp-endpoints} gives
\[
 \{d_j(y_p^A,z^A):p\ge p_0(j,\theta),\ A\in\mathcal A\}
 \subseteq\{u\in M_j:u\le\beta^j_\theta\}.
\]
Since $\beta^j_\theta\downarrow0$, the saturation of these represented
endpoint tails is order null for every $j$. Equation
\eqref{uc:faithfulness} applied to the original pair saturation proves
$y_p^A\to z^A$ uniformly in $A$.
\end{proof}

\begin{example}[A noncancellative controller]\label{mg:truncated}
On $M=[0,1]$, put $a\oplus b=\min\{1,a+b\}$ with the usual order.
Associativity is the identity
\[
 (a\oplus b)\oplus c=\min\{1,a+b+c\}=a\oplus(b\oplus c).
\]
For $D\downarrow0$ and $r<1$, each $\varepsilon>0$ has a
$d\in D$ with $d<\varepsilon$, whence
\[
 r\le\inf_{d\in D}(r\oplus d)\le r+\varepsilon.
\]
For $r=1$ every translated value equals $1$. Thus
\eqref{mg:zero-translation} holds. Cancellation fails since
$1\oplus0=1\oplus1$. In a normed space,
\[
 d(x,y)=\min\{1,\|x-y\|\}
\]
is symmetric and invariant, and its triangle inequality follows from
\[
 \min\{1,\|x-z\|\}
 \le\min\{1,\|x-y\|+\|y-z\|\}
 =d(x,y)\oplus d(y,z).
\]
For $0<r<1$, the condition $d(x,y)\le r$ is equivalent to
$\|x-y\|\le r$. Consequently the two distances have exactly the
same uniformly null pair nets. This example verifies a strict weakening
of controller cancellation, without asserting necessity of
\eqref{mg:zero-translation} for every possible interchange theorem.
\end{example}

\begin{corollary}[A complete represented subspace]\label{lt:carrier-transfer}
The theorem applies on a complete subspace $C\subseteq Z$ when the
faithful representation is defined on $C$ and every candidate value,
row value, and row endpoint belongs to $C$. Its limits and uniform
endpoint convergence hold in the original ambient structure.
\end{corollary}
\begin{proof}
For the restricted structure, the representing maps are
$d_j|_{C^2}:C^2\to M_j$ with the same controller operations and
zero-translation law. The arrays $x^A_{p,\gamma}$, $x^A_\gamma$ and
$y_p^A$ take values in $C$ by hypothesis, and their row convergence
and displayed comparison inequalities hold in that restricted
structure. Applying Theorem~\ref{lt:order-transfer} with target $C$
therefore produces $z^A\in C$ with
\[
 x^A_\gamma\to_C z^A,\qquad
 \left\langle\{(y_p^A,z^A):p\ge p_0,\ A\in\mathcal A\}
       :p_0\in P\right\rangle\longrightarrow_{\mathcal U_C}\Delta_C.
\]
For the inclusion $\iota:C\hookrightarrow Z$, the restricted
structure is defined by the implications
$u_\alpha\to_C u\Rightarrow\iota(u_\alpha)\to_Z\iota(u)$ and
pair-nullity in $C\Rightarrow$ pair-nullity of the image in $Z$.
Applying them to the two displayed nets proves the ambient conclusions.
\end{proof}

\begin{proposition}[Moving positive majorants]\label{lt:moving-majorant}
For the lattice-group specialization $M_j=(F_j)_+$ of
Theorem~\ref{lt:order-transfer}, the comparison hypothesis is verified by
\[
 d_j(x^A_\gamma,x^A_{p,\gamma})\le W^{A,j}_{p,\gamma},\qquad
 W^{A,j}_{p,\gamma}\longrightarrow_o v_p^{A,j},\qquad
 0\le v_p^{A,j}\le b_p^j,\qquad b_p^j\longrightarrow_o0.
\]
After construction of the limits,
$d_j(y_p^A,z^A)\le v_p^{A,j}$.
\end{proposition}
\begin{proof}
Set $r^{A,j}_{p,\gamma}=|W^{A,j}_{p,\gamma}-v_p^{A,j}|$.
If $\beta^j_\theta\downarrow0$ controls $b_p^j$, then
$P^j_\theta=\{p:b_p^j\le\beta^j_\theta\}$ contains a tail and
\[
 d_j(x^A_\gamma,x^A_{p,\gamma})
 \le v_p^{A,j}+r^{A,j}_{p,\gamma}
 \le\beta^j_\theta+r^{A,j}_{p,\gamma}\qquad(p\in P^j_\theta).
\]
Theorem~\ref{lt:order-transfer} gives both limits. For fixed $p$,
\[
 d_j(y_p^A,z^A)\le v_p^{A,j}+r^{A,j}_{p,\gamma}
 +d_j(y_p^A,x^A_{p,\gamma})+d_j(x^A_\gamma,z^A).
\]
Here $r^{A,j}_{p,\gamma}\to_o0$ by its definition, while point
faithfulness for the two constructed limits gives
\[
 d_j(y_p^A,x^A_{p,\gamma})\to_o0,
 \qquad d_j(x^A_\gamma,z^A)\to_o0.
\]
Their sum $h_\gamma$ is order null by Lemma~\ref{mg:calculus}(1).
If $u_\eta\downarrow0$ controls $h_\gamma$, intersection of the
comparison tail and the tail $h_\gamma\le u_\eta$ gives
\[
 d_j(y_p^A,z^A)-v_p^{A,j}\le u_\eta\qquad(\eta).
\]
Taking the infimum in $F_j$ yields
$d_j(y_p^A,z^A)\le v_p^{A,j}$.
\end{proof}

\subsection{Regulator interchange}

In regulator convergence, the finite comparisons still remove the
residual terms from each fixed endpoint inequality, but recovering the
candidate requires the separately chosen regulators to be brought
under one control. We begin with sequential approximation, where
bounded fusion performs this passage while leaving the inner nets
on their original directed sets.

\begin{theorem}[Sequential regulator interchange]\label{lt:d-transfer}
Suppose $Z$ is separated, sequentially complete, and faithfully
$D$-represented by symmetric $d_j$ with Dedekind
$\sigma$-complete WSD controllers $F_j$. Let $\Gamma_A$ be directed,
and suppose
\[
 x^A_{n,\gamma}\longrightarrow y_n^A\qquad(n\ge1).
\]
Let $b_{j,n}\in(F_j)_+$ be independent of $A$, and suppose
\[
 b_{j,n}\longrightarrow_o0,\qquad
 r^{A,j}_{n,\gamma}\longrightarrow_D0,\qquad
 d_j(x^A_\gamma,x^A_{n,\gamma})
 \le b_{j,n}+r^{A,j}_{n,\gamma}
 \quad\text{eventually in }\gamma,
\]
where all residuals are positive. Then unique $z^A$ satisfy
\[
 x^A_\gamma\longrightarrow z^A,\qquad
 y_n^A\longrightarrow z^A\quad\text{uniformly in }A,\qquad
 d_j(y_n^A,z^A)\le b_{j,n}.
\]
\end{theorem}
\begin{proof}
For every $j$, Lemma~\ref{lt:sequence-envelope} gives
\[
 v_{j,N}=\sup_{n\ge N}b_{j,n},\qquad v_{j,N}\downarrow0.
\]
Fix $A,j$ and suppress those indices in the first comparison. For
fixed $n,m$, take an inner index beyond the two comparison thresholds.
Symmetry and the triangle inequality then give
\[
 \begin{aligned}
 d(y_n,y_m)
 &\le d(y_n,x_{n,\gamma})+d(x_{n,\gamma},x_\gamma)
       +d(x_\gamma,x_{m,\gamma})+d(x_{m,\gamma},y_m)\\
 &\le b_n+b_m+g_{n,m,\gamma},\\
 g_{n,m,\gamma}
 &=r_{n,\gamma}+r_{m,\gamma}
       +d(y_n,x_{n,\gamma})+d(x_{m,\gamma},y_m).
 \end{aligned}
\]
In this decomposition, the first two terms defining $g_{n,m,\gamma}$
are $D$-null by hypothesis, and point faithfulness for the two known
row limits makes the last two terms $D$-null as well. Thus
Lemma~\ref{lt:d-calculus} gives
$g_{n,m,\gamma}\to_o0$. If $u_\theta\downarrow0$ controls this
net, a common upper index of its control threshold and the comparison
threshold gives
\[
 d(y_n,y_m)-b_n-b_m\le u_\theta\qquad(\theta).
\]
Taking the infimum yields
\[
 d(y_n,y_m)\le b_n+b_m,
 \qquad n,m\ge N\Longrightarrow d(y_n,y_m)\le2v_N.
\]
For coordinate $j$, the array $c^j_{rs}=2v_{j,s}$ is bounded by
$2v_{j,1}$ and has decreasing zero rows. For every selector $\varphi$,
\[
 n,m\ge\varphi(1)\quad\Longrightarrow\quad
 d_j(y_n^A,y_m^A)\le2v_{j,\varphi(1)}\le H_{c^j}(\varphi).
\]
Thus the endpoint product net is originally null by pair faithfulness.
For each $A$, sequential completeness and separation give one
$z^A\in Z$ with $y_n^A\to z^A$.

Fix $A,j,n$ again. The triangle inequality and the endpoint-pair
bound give
\[
 d_j(y_n^A,z^A)
 \le b_{j,n}+b_{j,m}+d_j(y_m^A,z^A)\qquad(m\ge1).
\]
The last term is $D$-null by point faithfulness, hence order null,
and $b_{j,m}\to_o0$ by hypothesis. If $u_\theta\downarrow0$ controls
their sum, the displayed inequality gives
\[
 d_j(y_n^A,z^A)-b_{j,n}\le u_\theta\qquad(\theta),
 \qquad d_j(y_n^A,z^A)\le b_{j,n}.
\]

For the candidate limit fix $A,j$ and write
\[
 e_\gamma=d_j(x^A_\gamma,z^A),\qquad
 h_{n,\gamma}=r^{A,j}_{n,\gamma}
                +d_j(x^A_{n,\gamma},y_n^A).
\]
Both summands in $h_n$ are $D$-null, so
Lemma~\ref{lt:d-calculus} supplies a regulator $a^{(n+1)}$ for
$h_n$. For every $n$, the comparison and endpoint bounds give
\begin{equation}\label{lt:d-candidate-comparison}
 \begin{aligned}
 e_\gamma
 &\le d_j(x^A_\gamma,x^A_{n,\gamma})
       +d_j(x^A_{n,\gamma},y_n^A)+d_j(y_n^A,z^A)\\
 &\le2b_{j,n}+h_{n,\gamma}
 \le2v_{j,n}+h_{n,\gamma}
 \qquad(\gamma\ge\gamma_n)
 \end{aligned}
\end{equation}
for some comparison threshold $\gamma_n\in\Gamma_A$.
Put $a^{(1)}_{rs}=2v_{j,s}$ and $\iota(r)=1$. A threshold
$\gamma_h$ for the selector $\iota$ in the regulator of $h_1$
and an index $\gamma_w\ge\gamma_1,\gamma_h$ give
\[
 e_\gamma\le w:=2v_{j,1}+H_{a^{(2)}}(\iota)
 \qquad(\gamma\ge\gamma_w).
\]
Define
\[
 B_{1s}=0,\qquad
 B_{rs}=w\wedge2^r\sum_{1\le k<r}a^{(k)}_{r-k,s}
 \quad(r\ge2).
\]
Lemma~\ref{lt:fusion} proves that this is one bounded regulator and
that its selector bound dominates the capped shifted sum of every
finite initial list of the $a^{(k)}$.

Fix $\varphi\in\mathbb N^{\mathbb N}$, and put
\[
 n=\varphi(2),\qquad \eta(r)=\varphi(r+n+1).
\]
The term with $r=1$ in the first shifted selector supremum gives
\[
 2v_{j,n}=a^{(1)}_{1,\varphi(2)}
 \le H_{a^{(1)}}(r\mapsto\varphi(r+1)).
\]
Let $\gamma_h(n,\eta)$ satisfy
$h_{n,\gamma}\le H_{a^{(n+1)}}(\eta)$ whenever
$\gamma\ge\gamma_h(n,\eta)$, and choose an upper index
\[
 \gamma_\varphi\ge\gamma_w,\gamma_n,\gamma_h(n,\eta).
\]
For every $\gamma\ge\gamma_\varphi$, positivity of the intervening
summands and the fusion inequality give
\[
 \begin{aligned}
 e_\gamma
 &\le w\wedge\left[
 H_{a^{(1)}}(r\mapsto\varphi(r+1))
 +H_{a^{(n+1)}}(r\mapsto\varphi(r+n+1))\right]\\
 &\le w\wedge\sum_{k=1}^{n+1}
 H_{a^{(k)}}(r\mapsto\varphi(r+k))
 \le H_B(\varphi).
 \end{aligned}
\]
Thus $d_j(x^A_\gamma,z^A)\to_D0$. Since $A,j$ were arbitrary,
point faithfulness gives $x^A_\gamma\to z^A$ for every $A$.

For uniformity of the endpoints, the regulator
$d^j_{rs}=v_{j,s}$ is independent of $A$, and
\[
 n\ge\varphi(1),\quad A\in\mathcal A
 \quad\Longrightarrow\quad
 d_j(y_n^A,z^A)\le b_{j,n}\le v_{j,\varphi(1)}
 \le H_{d^j}(\varphi).
\]
Every represented endpoint saturation is therefore $D$-null.
Pair faithfulness gives $y_n^A\to z^A$ uniformly in $A$.
\end{proof}

For a general outer net, completeness of the range allows the
comparison method to use a countable selection of rows for each
represented coordinate, without requiring that selection to exhaust
the approximation index. The next theorem makes this passage with
full completeness in place of the preceding sequential assumption.

\begin{theorem}[Countably accessible regulator comparisons]
\label{lt:d-anchors}
Suppose $Z$ is separated, complete, and faithfully $D$-represented,
with the controller hypotheses of Theorem~\ref{lt:d-transfer}.
Let $P,\Gamma_A$ be directed and
$x^A_{p,\gamma}\to y_p^A$. For each $j$, let $c^j$ be a regulator
and let $P^j_\varphi\subseteq P$ contain a tail, independently of $A$.
Suppose
\[
 r^{A,j}_{p,\gamma}\longrightarrow_o0,\qquad
 p\in P^j_\varphi\Longrightarrow
 d_j(x^A_\gamma,x^A_{p,\gamma})
 \le H_{c^j}(\varphi)+r^{A,j}_{p,\gamma}
 \quad\text{eventually}.
\]
For each $j$, suppose a sequence $(p_{j,n})_n$ satisfies
\[
 \{p_{j,n}:n\ge1\}\cap P^j_\varphi\ne\varnothing
 \qquad(\varphi\in\mathbb N^{\mathbb N}),
\]
and each residual belonging to $p_{j,n}$ is $D$-null.
Then $x^A_\gamma\to z^A$, $y_p^A\to z^A$ uniformly in $A$, and
\[
 p\in P^j_\varphi\Longrightarrow
 d_j(y_p^A,z^A)\le H_{c^j}(\varphi).
\]
\end{theorem}
\begin{proof}
Fix $A,j$ and write $d,c,P_\varphi$ for the corresponding data.
For $p\in P_\varphi$ and $q\in P_\psi$, an index beyond both
comparison thresholds gives
\[
 \begin{aligned}
 d(y_p,y_q)
 &\le H_c(\varphi)+H_c(\psi)+g_{p,q,\gamma},\\
 g_{p,q,\gamma}
 &=r_{p,\gamma}+r_{q,\gamma}
   +d(y_p,x_{p,\gamma})+d(x_{q,\gamma},y_q).
 \end{aligned}
\]
The residuals are order null by hypothesis, while point faithfulness
and Lemma~\ref{lt:d-calculus} make the two row errors order null.
Lemma~\ref{lt:fixed-bound}, with this finite sum as its null term,
gives
\[
 d(y_p,y_q)\le H_c(\varphi)+H_c(\psi).
\]
For each $j,\varphi$, fix $p^j_\varphi$ with
$\{p:p\ge p^j_\varphi\}\subseteq P^j_\varphi$. Then
\[
 p,q\ge p^j_\varphi
 \quad\Longrightarrow\quad
 d_j(y_p^A,y_q^A)\le2H_{c^j}(\varphi)=H_{2c^j}(\varphi).
\]
The bounded regulator $2c^j$ controls the endpoint product net in
coordinate $j$. Pair faithfulness and completeness give
$y_p^A\to z^A$ in the original structure, and separation makes the
limit unique.

For fixed $A,j$ and $p\in P_\varphi$, the triangle inequality gives
\[
 d(y_p,z)\le H_c(\varphi)+H_c(\psi)+d(y_q,z)
 \qquad(q\ge p^j_\psi).
\]
Point faithfulness makes $d(y_q,z)$ $D$-null and hence order null.
Lemma~\ref{lt:fixed-bound} therefore gives, for every $\psi$,
\[
 d(y_p,z)-H_c(\varphi)\le H_c(\psi).
\]
Weak $\sigma$-distributivity says that the infimum of the right sides
is $0$, so $d(y_p,z)\le H_c(\varphi)$.

Write $p_n=p_{j,n}$. For each $n$, the stipulated $D$-null residual
at $p_n$ and the $D$-null row error give a regulator $q^{(n)}$ for
\[
 h_{n,\gamma}=r^{A,j}_{p_n,\gamma}
                +d(x_{p_n,\gamma},y_{p_n}).
\]
The candidate error $e_\gamma=d(x_\gamma,z)$ satisfies
\[
 p_n\in P_\psi\quad\Longrightarrow\quad
 e_\gamma\le2H_c(\psi)+h_{n,\gamma}
 \quad(\gamma\ge\delta(n,\psi))
\]
for a comparison threshold $\delta(n,\psi)$. With $\iota(r)=1$,
accessibility makes
\[
 n_0=\min\{n:p_n\in P_\iota\}
\]
well-defined. The regulator $q^{(n_0)}$ supplies a threshold
$\delta_h$ for $\iota$. For an upper index
$\delta_w\ge\delta(n_0,\iota),\delta_h$,
\[
 \gamma\ge\delta_w\quad\Longrightarrow\quad
 e_\gamma\le w:=2H_c(\iota)+H_{q^{(n_0)}}(\iota).
\]
Set $a^{(1)}=2c$, $a^{(n+1)}=q^{(n)}$, and define
\[
 B_{1s}=0,\qquad B_{rs}
 =w\wedge2^r\sum_{1\le k<r}a^{(k)}_{r-k,s}\quad(r\ge2).
\]
Lemma~\ref{lt:fusion} makes $B$ a regulator with its stated shifted
finite-sum bound.

For a selector $\varphi$, define
\[
 \psi(r)=\varphi(r+1),\qquad
 n=\min\{m:p_m\in P_\psi\},\qquad
 \eta(r)=\varphi(r+n+1).
\]
Again $n$ exists by accessibility. Take $\delta_h(n,\eta)$ with
$h_{n,\gamma}\le H_{q^{(n)}}(\eta)$ for
$\gamma\ge\delta_h(n,\eta)$, and let
\[
 \delta_\varphi\ge\delta_w,\delta(n,\psi),\delta_h(n,\eta).
\]
For $\gamma\ge\delta_\varphi$,
\[
 \begin{aligned}
 e_\gamma
 &\le w\wedge[2H_c(\psi)+H_{q^{(n)}}(\eta)]\\
 &=w\wedge\left[
 H_{a^{(1)}}(r\mapsto\varphi(r+1))
 +H_{a^{(n+1)}}(r\mapsto\varphi(r+n+1))\right]\\
 &\le w\wedge\sum_{k=1}^{n+1}
 H_{a^{(k)}}(r\mapsto\varphi(r+k))
 \le H_B(\varphi).
 \end{aligned}
\]
Thus $e_\gamma\to_D0$ for every $A,j$, and point faithfulness gives
$x^A_\gamma\to z^A$.

Finally,
\[
 p\ge p^j_\varphi,\quad A\in\mathcal A
 \quad\Longrightarrow\quad
 d_j(y_p^A,z^A)\le H_{c^j}(\varphi).
\]
The thresholds $p^j_\varphi$ and regulators $c^j$ are independent of
$A$. Pair faithfulness applied to the endpoint saturation therefore
gives uniform convergence in the original structure.
\end{proof}

Countable selection in Theorem~\ref{lt:d-anchors} concerns the
approximating rows, yet the inner direction offers another place where
countability can enter. By placing the restriction there, we can allow
the outer comparisons to remain inaccessible to every countable
family of rows.

\begin{theorem}[Regulator interchange through the inner direction]
\label{ic:inner-regulator-interchange}
Let $Z$ be a separated complete uniform convergence space with a faithful
$D$-representation $d_j:Z^2\to(F_j)_+$, $j\in J$, where the $d_j$ are
symmetric and satisfy the triangle inequality, and the $F_j$ are
Dedekind $\sigma$-complete WSD abelian lattice-ordered groups.
Let $P$ and $\Gamma_A$, $A\in\mathcal A\ne\varnothing$, be directed,
and suppose
\[
 x^A_{p,\gamma}\longrightarrow y_p^A
 \qquad(p\in P,\ A\in\mathcal A).
\]
For each $j$, let $c^j$ be a regulator and let
$P^j_\varphi\subseteq P$ contain a tail of $P$ for every
$\varphi\in\mathbb N^{\mathbb N}$, independently of $A$.
Suppose the positive residuals satisfy
\begin{equation}\label{ic:inner-comparison}
 \begin{gathered}
 r^{A,j}_{p,\gamma}\longrightarrow_o0,\\
 p\in P^j_\varphi\quad\Longrightarrow\quad
 d_j(x^A_\gamma,x^A_{p,\gamma})
 \le H_{c^j}(\varphi)+r^{A,j}_{p,\gamma}
 \quad\text{eventually in }\gamma.
 \end{gathered}
\end{equation}
Assume either that every $\Gamma_A$ is countable, or that every
$\Gamma_A$ has a countable cofinal subset and every $F_j$ is Dedekind
complete. Then unique $z^A\in Z$ satisfy
\begin{equation}\label{ic:inner-conclusions}
 x^A_\gamma\longrightarrow z^A,
 \qquad y_p^A\longrightarrow z^A\quad\text{uniformly in }A,
\end{equation}
and
\begin{equation}\label{ic:inner-endpoint-bound}
 p\in P^j_\varphi\quad\Longrightarrow\quad
 d_j(y_p^A,z^A)\le H_{c^j}(\varphi).
\end{equation}
No countable accessibility of the sets $P^j_\varphi$ is required.
\end{theorem}
\begin{proof}
For a regulator $a$ in a WSD controller, the selector bounds form a
downward-directed family. Indeed, the selector
$\chi(r)=\max\{\varphi(r),\psi(r)\}$ satisfies
\[
 H_a(\chi)\le H_a(\varphi),H_a(\psi),
 \qquad \inf_\varphi H_a(\varphi)=0.
\]
Consequently every positive $D$-null net is order null. For decreasing
zero controls $u_\theta$ and $v_\eta$, translation invariance of the
order gives
\[
 b\le u_\theta+v_\eta\ (\theta,\eta)
 \ \Longrightarrow\ b-u_\theta\le\inf_\eta v_\eta=0\ (\theta)
 \ \Longrightarrow\ b\le\inf_\theta u_\theta=0.
\]
Their product sum decreases to zero. Induction treats every finite
sum of positive order-null nets. Moreover, if $s\le t+h_\gamma$
eventually and $h_\gamma\to_o0$ with controls $u_\theta\downarrow0$,
then intersecting the two tails gives
\begin{equation}\label{ic:fixed-bound-removal}
 s-t\le u_\theta\quad(\theta),\qquad s-t\le0.
\end{equation}

Fix $A,j$ and suppress them until the final paragraph. Let
$p\in P_\varphi$ and $q\in P_\psi$. On a common inner tail,
\begin{align*}
 d(y_p,y_q)
 &\le d(y_p,x_{p,\gamma})+d(x_{p,\gamma},x_\gamma)
       +d(x_\gamma,x_{q,\gamma})+d(x_{q,\gamma},y_q)\\
 &\le H_c(\varphi)+H_c(\psi)+r_{p,\gamma}+r_{q,\gamma}
       +d(y_p,x_{p,\gamma})+d(x_{q,\gamma},y_q).
\end{align*}
The last two terms are $D$-null by point faithfulness and are therefore
order null. The sum of the last four terms is order null. Equation
\eqref{ic:fixed-bound-removal} yields
\begin{equation}\label{ic:inner-two-endpoints}
 d(y_p,y_q)\le H_c(\varphi)+H_c(\psi)
 \qquad(p\in P_\varphi,\ q\in P_\psi).
\end{equation}
For every $\varphi$ there exists $p_\varphi$ with
$\{p:p\ge p_\varphi\}\subseteq P_\varphi$. Hence
\[
 p,q\ge p_\varphi\quad\Longrightarrow\quad
 d(y_p,y_q)\le2H_c(\varphi)=H_{2c}(\varphi).
\]
In every coordinate $j$, the regulator $2c^j$ controls the endpoint
product net. Pair faithfulness makes $(y_p^A)_p$ Cauchy in the original
structure. Completeness and separation give its unique limit $z^A$.

Fix $p\in P_\varphi$. For $q\ge p_\psi$,
\[
 d(y_p,z)\le H_c(\varphi)+H_c(\psi)+d(y_q,z).
\]
Point faithfulness gives $d(y_q,z)\to_D0$, hence
$d(y_q,z)\to_o0$. Equation~\eqref{ic:fixed-bound-removal} gives
\[
 d(y_p,z)\le H_c(\varphi)+H_c(\psi)\qquad(\psi),
 \qquad d(y_p,z)-H_c(\varphi)\le\inf_\psi H_c(\psi)=0.
\]
This proves~\eqref{ic:inner-endpoint-bound}.

Put
\[
 e_\gamma=d(x_\gamma,z),\qquad
 h_{p,\gamma}=r_{p,\gamma}+d(x_{p,\gamma},y_p).
\]
Each $h_p$ is order null. For $p\in P_\varphi$, the comparison and
the endpoint bound imply
\begin{equation}\label{ic:inner-error}
 e_\gamma\le2H_c(\varphi)+h_{p,\gamma}
 \quad\text{eventually in }\gamma.
\end{equation}
Define
\[
 R=\{u\in F_+:e_\gamma\le u\text{ eventually in }\gamma\}.
\]
An attained upper bound for one $h_p$, with $p\in P_\varphi$, and
\eqref{ic:inner-error} give $R\ne\varnothing$. Intersecting two
controlling tails gives
\[
 u,v\in R\quad\Longrightarrow\quad u\wedge v\in R.
\]
Let $b\in F$ be a lower bound of $R$. For every $\varphi$ and every
$p\in P_\varphi$, let $u_{p,\theta}\downarrow0$ regulate $h_p$.
Equation~\eqref{ic:inner-error} gives
\[
 2H_c(\varphi)+u_{p,\theta}\in R\quad(\theta),
 \qquad b-2H_c(\varphi)\le\inf_\theta u_{p,\theta}=0.
\]
Since $P_\varphi$ contains a nonempty tail,
\[
 b\le2H_c(\varphi)=H_{2c}(\varphi)\quad(\varphi),
 \qquad b\le\inf_\varphi H_{2c}(\varphi)=0.
\]
Thus $R\downarrow0$. With its reverse order, $R$ regulates $e_\gamma$,
and
\begin{equation}\label{ic:inner-order-null}
 e_\gamma\longrightarrow_o0.
\end{equation}

There exist $w\in F_+$ and $\gamma_*$ such that
$e_\gamma\le w$ whenever $\gamma\ge\gamma_*$. Both alternatives
in the statement provide a countable cofinal subset of $\Gamma_A$.
Let $(\alpha_N)_{N\ge1}$ enumerate such a subset, with repetitions
when necessary. Directedness provides a sequence $(\gamma_N)$ with
\[
 \gamma_1\ge\gamma_*,\alpha_1,
 \qquad \gamma_{N+1}\ge\gamma_N,\alpha_{N+1}\quad(N\ge1).
\]
This sequence is increasing and cofinal. The extrema
\begin{equation}\label{ic:inner-tail-envelopes}
 v_N=\sup\{e_\gamma:\gamma\ge\gamma_N\}
\end{equation}
exist and satisfy $0\le v_{N+1}\le v_N\le w$. Under the first
alternative, each set in~\eqref{ic:inner-tail-envelopes} is nonempty,
countable, and bounded above, so Dedekind $\sigma$-completeness
supplies its supremum. Under the second alternative, Dedekind
completeness supplies the supremum without a cardinality restriction
on that set.

Let $t_\lambda\downarrow0$ regulate
\eqref{ic:inner-order-null}. For each $\lambda$ there exists
$\delta_\lambda$ such that
$\gamma\ge\delta_\lambda\Longrightarrow e_\gamma\le t_\lambda$.
Cofinality supplies $N(\lambda)$ with
$\gamma_{N(\lambda)}\ge\delta_\lambda$. Therefore
\[
 v_{N(\lambda)}\le t_\lambda.
\]
If $b\le v_N$ for every $N$, then
$b\le t_\lambda$ for every $\lambda$, and $b\le0$. Hence
$v_N\downarrow0$. The double sequence
\[
 a_{rs}=v_s\qquad(r,s\ge1)
\]
is a regulator bounded by $w$. For every selector $\eta$,
\[
 \gamma\ge\gamma_{\eta(1)}\quad\Longrightarrow\quad
 e_\gamma\le v_{\eta(1)}
 \le\sup_{r\ge1}v_{\eta(r)}=H_a(\eta).
\]
Thus $e_\gamma\to_D0$. This construction is valid for every $A,j$.
Point faithfulness gives $x_\gamma^A\to z^A$.

For every $j,\varphi$, the index $p^j_\varphi$ is independent of $A$,
and the already established endpoint estimate gives
\[
 p\ge p^j_\varphi,\ A\in\mathcal A
 \quad\Longrightarrow\quad
 d_j(y_p^A,z^A)\le H_{c^j}(\varphi).
\]
The same regulator $c^j$ controls the saturated endpoint pairs in
coordinate $j$. Pair faithfulness gives the uniform assertion
in~\eqref{ic:inner-conclusions}.
\end{proof}

Since sequences of records satisfy the inner hypothesis, they provide
a setting in which to test whether it differs from the
preceding selection condition. In the following example, every
countable selection of adequate outer rows fails, even after the
comparison regulator is changed within the specified faithful
representation.

\begin{example}[Inner sequences with inaccessible outer comparisons]
\label{ic:inaccessible-example}
Let
\[
 I=\mathbb N_0^{\mathbb N},\qquad F=\mathbb R^I,
 \qquad d(u,v)=|u-v|,
\]
with coordinatewise order, and give $Z=F$ its $D$-convergence and the
canonical group pair convergence. This is a separated complete
faithfully $D$-represented space. The controller $F$ is Dedekind
complete and WSD but does not have the countable sup property.
There are an outer directed set $P$, convergent inner sequences
$(x_{p,k})_k$, and a nonconstant candidate sequence $(x_k)_k$ satisfying
Theorem~\ref{ic:inner-regulator-interchange}, for which no countable
family of the given outer rows satisfies the accessibility requirement
of Theorem~\ref{lt:d-anchors} under any comparison regulator for $d$.
\end{example}
\begin{proof}
Bounded nonempty suprema and infima in $F$ are coordinatewise real
suprema and infima, so $F$ is Dedekind complete. For a regulator $a$,
an index $i\in I$, and $\varepsilon>0$, each set
\[
 \{s\ge1:a_{rs}(i)\le\varepsilon\}
\]
is nonempty. Define $\varphi(r)$ as its least member. Then
\[
 H_a(\varphi)(i)=\sup_r a_{r,\varphi(r)}(i)\le\varepsilon,
 \qquad \inf_\varphi H_a(\varphi)(i)=0.
\]
This proves WSD. For $\mathbf1_B$ denoting the indicator of
$B\subseteq I$,
\[
 \sup_{i\in I}\mathbf1_{\{i\}}=\mathbf1_I,
 \qquad
 \sup_{i\in C}\mathbf1_{\{i\}}=\mathbf1_C\ne\mathbf1_I
 \quad(C\subseteq I\text{ countable}).
\]
Thus the countable sup property fails.

The constant zero net has the zero regulator, and quasi-subnets
preserve every eventual selector bound. If two null nets have
regulators $a$ and $b$, the matrix $c_{rs}=2(a_{rs}+b_{rs})$ is
bounded and each row decreases to zero. For fixed $\varphi$, put
$C=H_{a+b}(\varphi)$. Then
\[
 a_{r,\varphi(r)}\le C,\qquad b_{s,\varphi(s)}\le C,
 \qquad H_a(\varphi)+H_b(\varphi)\le2C=H_c(\varphi).
\]
The two original bounds hold on a common tail, or on a product tail
for independently indexed nets. The displayed inequality controls
their sum, and it also controls every finite mixing of their values.
Inversion leaves moduli unchanged. These facts establish the net
convergence axioms and continuity of the group operations, while the
canonical pair convergence is detected by $d(u,v)=|u-v|$ by its
definition. If $z_\alpha$ converges to both $x$ and $y$, regulators
$a,b$ for its two errors give, on a common tail for each pair of
selectors,
\[
 |x-y|\le |x-z_\alpha|+|z_\alpha-y|
       \le H_a(\varphi)+H_b(\psi).
\]
Taking the two coordinatewise infima yields $|x-y|=0$ by WSD, so
$x=y$.

Let $(u_\alpha)_{\alpha\in D}$ be a Cauchy net. A regulator $a$
for its difference net gives, for each $\varphi$, two indices whose
common upper bound $\alpha_\varphi$ satisfies
\begin{equation}\label{ic:product-cauchy-diameter}
 \alpha,\beta\ge\alpha_\varphi
 \quad\Longrightarrow\quad
 |u_\alpha-u_\beta|\le H_a(\varphi).
\end{equation}
For each coordinate $i$ and each $\varepsilon>0$, WSD supplies a
selector $\varphi$ with $H_a(\varphi)(i)<\varepsilon$. Equation
\eqref{ic:product-cauchy-diameter} therefore makes $(u_\alpha(i))$
a real Cauchy net. Its limit $z(i)$ exists by real completeness,
and these limits define $z\in\mathbb R^I=F$.

Fix $\varphi$, $\alpha\ge\alpha_\varphi$, $i\in I$, and
$\delta>0$. Coordinate convergence supplies $\beta_0$ such that
$|u_\beta(i)-z(i)|<\delta$ for $\beta\ge\beta_0$. For an index
$\beta\ge\beta_0,\alpha_\varphi$, the Cauchy bound gives
\[
 |u_\alpha(i)-z(i)|
 \le |u_\alpha(i)-u_\beta(i)|+|u_\beta(i)-z(i)|
 \le H_a(\varphi)(i)+\delta.
\]
The infimum over $\delta>0$ yields
$|u_\alpha(i)-z(i)|\le H_a(\varphi)(i)$. This holds at every
coordinate, so
\[
 \alpha\ge\alpha_\varphi
 \quad\Longrightarrow\quad |u_\alpha-z|\le H_a(\varphi).
\]
The original regulator $a$ consequently proves $u_\alpha\to_Dz$,
establishing completeness for arbitrary Cauchy nets.

Define a regulator by
\[
 c_{rs}(i)=\mathbf1_{\{s\le i(r)\}},
 \qquad 0\le c_{rs}\le\mathbf1_I,
 \qquad c_{rs}\downarrow_s0.
\]
Let $P=\mathbb N^{\mathbb N}$ with its coordinatewise order, and put
\[
 y_p=H_c(p),\qquad
 x_k=k^{-1}\mathbf1_I,\qquad
 x_{p,k}=y_p+k^{-1}\mathbf1_I,
 \qquad P_\varphi=\{p:p\ge\varphi\}.
\]
The matrix $b_{rs}=s^{-1}\mathbf1_I$ is a regulator, and
\[
 k\ge\eta(1)\quad\Longrightarrow\quad
 d(x_{p,k},y_p)=k^{-1}\mathbf1_I
 \le\eta(1)^{-1}\mathbf1_I\le H_b(\eta).
\]
Thus each inner row converges to $y_p$. The comparison error is
independent of $k$, and
\[
 d(x_k,x_{p,k})=y_p=H_c(p)\le H_c(\varphi)
 \qquad(p\in P_\varphi).
\]
All residuals are zero, every $P_\varphi$ is a tail, and the inner
index is $\mathbb N$. The theorem applies and yields $z=0$.
The same regulator $b$ directly verifies $x_k\to_D0$, and the
displayed comparison verifies $y_p\to_D0$.

Let $(p_n)_{n\ge1}$ be any countable family of outer rows, including
possible repetitions. Define $i\in I$ by $i(n)=p_n(n)$. For every $n$,
\begin{equation}\label{ic:persistent-coordinate}
 y_{p_n}(i)
 =\sup_{r\ge1}\mathbf1_{\{p_n(r)\le i(r)\}}
 \ge\mathbf1_{\{p_n(n)\le i(n)\}}=1.
\end{equation}
Since $y_{p_n}\le\mathbf1_I$, the value in
\eqref{ic:persistent-coordinate} equals $1$. In particular,
\[
 \left(\inf_{n\ge1}y_{p_n}\right)(i)=1.
\]
Suppose a regulator $a$ and valid comparison cuts for these rows
admitted countable accessibility. For every selector $\eta$, there
would be an $n$ and a positive order-null residual $r_{n,k}$ with
\[
 y_{p_n}=d(x_k,x_{p_n,k})\le H_a(\eta)+r_{n,k}
 \quad\text{eventually in }k.
\]
The fixed-bound calculation~\eqref{ic:fixed-bound-removal} gives
$y_{p_n}\le H_a(\eta)$. Hence
\[
 \inf_n y_{p_n}\le H_a(\eta)\quad(\eta),
 \qquad \inf_n y_{p_n}\le\inf_\eta H_a(\eta)=0,
\]
contradicting~\eqref{ic:persistent-coordinate}. No replacement of the
comparison regulator and no selection of countably many of these rows
can therefore restore the accessibility hypothesis for the displayed
faithful representation. The argument concerns this specified
representation and does not impose a restriction on all other possible
representations of the same convergence space.
\end{proof}

If the controllers are super Dedekind complete and WSD, the distinction
between these two sources of countability disappears from the
interchange argument: Proposition~\ref{lt:order-d} identifies order
and regulator convergence, so that Theorem~\ref{lt:order-transfer}
applies without the countable accessibility requirement. Throughout
these passages the convergence remains the original one, and neither
a diagonal axiom nor its idempotence is required.

\subsection{Verification by ordered and regulator controls}
\label{uv:verification}

We can now return from the comparison estimates to the intrinsic
conditions and identify where each is supplied by the original
assumptions. Beginning with the target, we write its given bounds as
relations, so that their finite-intersection and composition laws
reveal the required stability without introducing another convergence;
the coefficient estimates will then be verified independently. For
the group statements in this subsection we use the canonical pair
convergence, in which a pair net is null exactly when its difference
net converges to zero, while the comparison theorems for arbitrary
uniform convergence spaces retain their stated scope independently
of this specialization.

\begin{lemma}[The original bounds and their translation sections]
\label{uv:translation-sections}
Let $Z$ be an abelian convergence group. For each $j\in J$, let $M_j$
be an ordered commutative monoid with least element $0$, isotone
addition $\oplus_j$ and finite meets, and suppose that the original
pair bounds are relations $\mathscr R_j(r)\subseteq Z^2$, $r\in M_j$,
with
\begin{align}
 \Delta_Z&\subseteq\mathscr R_j(r),&
 \mathscr R_j(r)^{-1}&=\mathscr R_j(r),\label{uv:relation-symmetry}\\
 \mathscr R_j(r)\circ\mathscr R_j(s)
 &\subseteq\mathscr R_j(r\oplus_j s),&
 \mathscr R_j(\inf_{\ell}r_\ell)
 &=\bigcap_{\ell}\mathscr R_j(r_\ell)
 \label{uv:relation-laws}
\end{align}
whenever the displayed infimum exists. Suppose that a pair net is
originally null exactly when, for every $j$, it eventually belongs to
each relation $\mathscr R_j(r)$ in some family $D_j\downarrow0$.
Put
\[
 C_j(r)=\{u\in Z:(z+u,z)\in\mathscr R_j(r)\text{ for every }z\in Z\}.
\]
Then
\begin{align}
 0&\in C_j(r),&-C_j(r)&=C_j(r),\label{uv:section-symmetry}\\
 C_j(r)+C_j(s)&\subseteq C_j(r\oplus_j s),&
 C_j(\inf_\ell r_\ell)&=\bigcap_\ell C_j(r_\ell).
 \label{uv:section-laws}
\end{align}
For every net $(x_i)_{i\in I}$ in $Z$,
\begin{equation}\label{uv:section-null-test}
 x_i\longrightarrow0
 \quad\Longleftrightarrow\quad
 \forall j\ \exists D_j\downarrow0\
 \forall r\in D_j\ \exists i_r\
 \forall i\ge i_r:\ x_i\in C_j(r).
\end{equation}
The same equivalence holds with $C_j(r)$ replaced by the zero section
\[
 V_j(r)=\{u\in Z:(u,0)\in\mathscr R_j(r)\}.
\]
Alternatively, replace the order-control detection assumption by the
regulator criterion
\[
 (z_i,w_i)\longrightarrow\Delta_Z
 \quad\Longleftrightarrow\quad
 \forall j\ \exists a^j\ \forall\varphi\ \exists i_\varphi\
 \forall i\ge i_\varphi:\quad
 (z_i,w_i)\in\mathscr R_j(H_{a^j}(\varphi)),
\]
where $a^j$ is a regulator in the stated lattice-group controller.
Then the section criteria hold with $D_j$ replaced by the selector
bounds $\{H_{a^j}(\varphi):\varphi\in\mathbb N^{\mathbb N}\}$ of
one regulator $a^j$ for each $j$.
\end{lemma}
\begin{proof}
The diagonal inclusion gives $0\in C_j(r)$. If $u\in C_j(r)$, then,
for every $z\in Z$, the pair $(z,z-u)$ belongs to $\mathscr R_j(r)$.
Its reverse belongs to the same relation, so
$(z-u,z)\in\mathscr R_j(r)$ and $-u\in C_j(r)$. If
$u\in C_j(r)$ and $v\in C_j(s)$, then
\[
 (z+u+v,z+v)\in\mathscr R_j(r),\qquad
 (z+v,z)\in\mathscr R_j(s)\qquad(z\in Z).
\]
Composition gives $u+v\in C_j(r\oplus_j s)$. Finally, the condition
that $(z+u,z)$ belong to every $\mathscr R_j(r_\ell)$ for every $z$
is equivalent, by \eqref{uv:relation-laws}, to membership in
$\mathscr R_j(\inf_\ell r_\ell)$ for every $z$. This proves the
intersection identity. It also implies monotonicity: if $r\le s$, then
$r=r\wedge s$, and hence
\[
 \mathscr R_j(r)=\mathscr R_j(r)\cap\mathscr R_j(s),\qquad
 C_j(r)=C_j(r)\cap C_j(s).
\]
In particular, $C_j(r)\subseteq C_j(s)$ whenever $r\le s$.

Order $I\times Z$ by
$(i,z)\preceq(i',z')$ if and only if $i\le i'$, and consider
\[
 \bigl(z+x_i,z\bigr)_{(i,z)\in I\times Z}.
\]
This is a net on a directed preorder, and its difference net is
$(x_i)_{(i,z)}$, whose tail-value sets equal those of $(x_i)_i$.
It is therefore originally null exactly when $x_i\to0$. Moreover,
\[
 \begin{aligned}
 &(z+x_i,z)\in\mathscr R_j(r)
       \text{ eventually in }I\times Z\\
 &\qquad\Longleftrightarrow
 (\exists i_r)(\forall i\ge i_r)(\forall z\in Z)\
       (z+x_i,z)\in\mathscr R_j(r)\\
 &\qquad\Longleftrightarrow
 (\exists i_r)(\forall i\ge i_r)\ x_i\in C_j(r).
 \end{aligned}
\]
The assumed pair criterion proves \eqref{uv:section-null-test}.
For the pair net $(x_i,0)$, its difference net is $(x_i)$ and
\[
 (x_i,0)\in\mathscr R_j(r)\quad\Longleftrightarrow\quad x_i\in V_j(r),
\]
which proves the zero-section criterion. Under a regulator description,
the assumed pair criterion reads
\[
 \forall j\ \exists a^j\ \forall\varphi\ \exists i_\varphi\
 \forall i\ge i_\varphi:\quad
 (z+x_i,z)\in\mathscr R_j(H_{a^j}(\varphi))\quad(z\in Z).
\]
The order on $I\times Z$ makes the last universal quantifier part of
every tail, so this is equivalent to
$x_i\in C_j(H_{a^j}(\varphi))$ for all $i\ge i_\varphi$.
Replacing the translated pair by $(x_i,0)$ gives the corresponding
statement for $V_j$.
\end{proof}

For the bounds in Definition~\ref{mg:representation}, the relations
in the preceding lemma are the sets on which the corresponding
original bound is at most $r$. Their symmetry and composition laws
are the stated finite comparison inequalities, while the intersection
identity is the order-theoretic equivalence between being below every
$r_\ell$ and being below their infimum. From these relations the lemma
derives translation-compatible sets through the original canonical
pair convergence, without imposing invariance on the bounds themselves,
and this also explains why finite-response estimates may use the
larger zero sections $V_j(r)$.

\begin{theorem}[Ordered controls verify intrinsic stability]
\label{uv:order-verification}\label{mg:intrinsic-properties}
Under the order hypotheses of Lemma~\ref{uv:translation-sections},
suppose that each $M_j$ satisfies \eqref{mg:zero-translation}.
Then the sets $C_j(r)$ are convergence-closed, the original convergence
has closed-base stability, and $Z$ satisfies \textnormal{(R)} and
\textnormal{(P)}. No completeness of the control monoids or of $Z$
is required for these conclusions.
\end{theorem}
\begin{proof}
For $D\downarrow0$ in $M_j$, the section intersection identity and
\eqref{mg:zero-translation} give
\begin{equation}\label{uv:zero-removal-sections}
 \bigcap_{d\in D}C_j(r\oplus_j d)
 =
 C_j\!\left(\inf_{d\in D}(r\oplus_j d)\right)
 =C_j(r).
\end{equation}
Suppose $x_\alpha\in C_j(r)$ and $x_\alpha\to x$. Continuity of
subtraction gives $x-x_\alpha\to0$, so
\eqref{uv:section-null-test} supplies $D\downarrow0$ such that, for
each $d\in D$, there is $\alpha_d$ satisfying
$x-x_\alpha\in C_j(d)$ for $\alpha\ge\alpha_d$. For any index on
that tail,
\[
 x=x_\alpha+(x-x_\alpha)\in C_j(r\oplus_j d).
\]
This is a fixed membership assertion for every $d\in D$, so
\eqref{uv:zero-removal-sections} gives $x\in C_j(r)$ and proves
closedness.

Let $(A_\lambda)_{\lambda\in\Lambda}$ be a decreasing null family
of nonempty sets, with nullity interpreted through its saturated
value net. Fix $j$ and use \eqref{uv:section-null-test} to obtain
$D_j\downarrow0$. For each $r\in D_j$, there is $\lambda_r$ such
that
\[
 \lambda\ge\lambda_r\quad\Longrightarrow\quad
 A_\lambda\subseteq C_j(r).
\]
For $\lambda\ge\lambda_r$, every net in $A_\lambda$ is a net in the
closed set $C_j(r)$, and its limit therefore belongs to $C_j(r)$.
Consequently
\[
 \lambda\ge\lambda_r\quad\Longrightarrow\quad
 a_Z(A_\lambda)\subseteq a_Z(C_j(r))=C_j(r).
\]
The same controls $D_j$ therefore prove that the saturated value
net of $a_Z(A_\lambda)$ converges to zero. This is
\textnormal{(R)}. If $K_Z(A)$ is the intersection of the
convergence-closed supersets of $A$, closedness gives also
$K_Z(A_\lambda)\subseteq C_j(r)$ on that tail, proving the
closed-base assertion.

For a family converging to $z$, subtracting $z$ produces a null
family. Translation by $z$ and its inverse preserve convergent nets and
therefore preserve closed sets. If $C\supseteq A$ is closed, then
$C-z\supseteq A-z$ is closed, so
$K_Z(A-z)\subseteq C-z$. Intersecting over $C$ gives
$K_Z(A-z)\subseteq K_Z(A)-z$. Applying this inclusion to $A-z$
and the translation $-z$ gives the reverse inclusion. Thus
\[
 K_Z(A-z)=K_Z(A)-z.
\]
The zero-point assertion applied to $A_\lambda-z$ gives
$\langle K_Z(A_\lambda)-z\rangle\to0$, and continuity of translation
then gives $\langle K_Z(A_\lambda)\rangle\to z$.

To prove \textnormal{(P)}, let $A_\lambda\downarrow$ be a null
family, let $z_i$ be a net, and suppose that for every $\lambda$
there is a null net $(h^\lambda_t)_{t\in T_\lambda}$ such that
\begin{equation}\label{uv:absorption-premise}
 \forall t_0\in T_\lambda\quad
 z_i\in A_\lambda+
       \{h^\lambda_t:t\ge t_0\}
 \quad\text{eventually in }i.
\end{equation}
Fix $j$. The outer null family has controls $D_j\downarrow0$.
For each $r\in D_j$, there is an index $\lambda(r)$ with
$A_{\lambda(r)}\subseteq C_j(r)$. The net $h^{\lambda(r)}$ has
controls $E_r\downarrow0$ in coordinate $j$, so for every
$e\in E_r$ there is $t(r,e)$ such that
\[
 t\ge t(r,e)\quad\Longrightarrow\quad
 h^{\lambda(r)}_t\in C_j(e).
\]
Applying \eqref{uv:absorption-premise} at this inner index gives
\[
 z_i\in C_j(r)+C_j(e)
       \subseteq C_j(r\oplus_j e)
 \quad\text{eventually}.
\]
Consequently the attained-bound set
\[
 T_j=\{u\in M_j:z_i\in C_j(u)\text{ eventually}\}
\]
contains every $r\oplus_j e$ with $r\in D_j$, $e\in E_r$.
It is nonempty. If $u,v\in T_j$, a common upper index for their
two witnessing tails gives
$z_i\in C_j(u)\cap C_j(v)=C_j(u\wedge v)$ there, so $u\wedge v$
belongs to $T_j$.

If $c$ is a lower bound of $T_j$, then for each $r\in D_j$,
\[
 c\le r\oplus_j e\quad(e\in E_r)
 \quad\Longrightarrow\quad
 c\le\inf_{e\in E_r}(r\oplus_j e)=r.
\]
Hence $c\le\inf D_j=0$. Since $0$ is a lower bound and $T_j$ is
closed under finite meets, $T_j\downarrow0$. The sets $C_j(u)$,
$u\in T_j$, are attained by the candidate net itself, and
\eqref{uv:section-null-test} now gives $z_i\to0$.
\end{proof}

In the ordered verification, an independently null added control can
be removed after the finite bounds have been combined. The closedness
argument remains available for regulator convergence, whereas recovery
of the candidate now asks us to combine the independently chosen
regulators, a task performed by bounded fusion for a countable family.
Thus countability enters the next theorem through the outer errors,
with every inner net retaining its original directed set.

\begin{theorem}[Regulator controls verify countable absorption]
\label{uv:regulator-verification}
Suppose the canonical pair convergence of an abelian convergence group
$Z$ is faithfully described by the regulator bounds of
Lemma~\ref{uv:translation-sections}, with Dedekind
$\sigma$-complete weakly $\sigma$-distributive abelian lattice-ordered
control groups $F_j$. Then the original convergence has closed-base
stability and satisfies \textnormal{(R)} and
\textnormal{($P_{\aleph_0}$)}.
\end{theorem}
\begin{proof}
Fix $j$ and write
\[
 H_a(\varphi)=\sup_{r\ge1}a_{r,\varphi(r)}
 \qquad(\varphi\in\mathbb N^{\mathbb N}).
\]
For two selectors $\varphi,\psi$, the selector
$\theta(r)=\max\{\varphi(r),\psi(r)\}$ satisfies
$H_a(\theta)\le H_a(\varphi),H_a(\psi)$, and weak
$\sigma$-distributivity gives $\inf_\varphi H_a(\varphi)=0$.
Translation in $F_j$ preserves this infimum, so
\begin{equation}\label{uv:regulator-zero-removal}
 \bigcap_\varphi C_j(r+H_a(\varphi))=C_j(r).
\end{equation}
If $x_\alpha\in C_j(r)$ and $x_\alpha\to x$, a regulator $a$
for $x-x_\alpha$ gives, for every selector $\varphi$,
$x-x_\alpha\in C_j(H_a(\varphi))$ eventually. The addition
inclusion gives $x\in C_j(r+H_a(\varphi))$, and
\eqref{uv:regulator-zero-removal} yields $x\in C_j(r)$.
Thus each section is convergence-closed.

For a decreasing originally null family $A_\lambda$, its original
regulator $a^j$ in coordinate $j$ gives
\[
 \forall\varphi\ \exists\lambda_\varphi\
 \forall\lambda\ge\lambda_\varphi:\quad
 A_\lambda\subseteq C_j(H_{a^j}(\varphi)).
\]
Closedness proves the same inclusions for $a_Z(A_\lambda)$ and for
$K_Z(A_\lambda)$. The regulator criterion of
Lemma~\ref{uv:translation-sections}, applied in every coordinate,
therefore proves \textnormal{(R)} and closed-base stability.
For a family converging to a point $z$, the translated family
$A_\lambda-z$ is null, and the equality
$K_Z(A_\lambda-z)=K_Z(A_\lambda)-z$, obtained by translating closed
supersets in both directions, proves the same assertion at $z$.

For countable absorption, let $A_n\downarrow$ be originally null,
and suppose that a candidate net $z_i$ and null inner nets
$h^n=(h^n_t)_{t\in T_n}$ satisfy
\begin{equation}\label{uv:countable-absorption-premise}
 \forall n\ge1\ \forall t_0\in T_n:\quad
 z_i\in A_n+\{h^n_t:t\ge t_0\}
 \quad\text{eventually in }i.
\end{equation}
Fix $j$ and suppress it from the notation. Let $a^{(1)}$ regulate
the saturated outer net, and let $a^{(n+1)}$ regulate $h^n$.
The choices are independent, and their bounds need not be common.
With $\iota(r)=1$, there is $n_0$ such that
\[
 A_{n_0}\subseteq C(H_{a^{(1)}}(\iota)).
\]
There is also $t_0\in T_{n_0}$ such that
$h^{n_0}_t\in C(H_{a^{(n_0+1)}}(\iota))$ whenever $t\ge t_0$.
Equation~\eqref{uv:countable-absorption-premise} supplies the attained
cap
\begin{equation}\label{uv:attained-cap}
 z_i\in C(w)\text{ eventually},\qquad
 w=H_{a^{(1)}}(\iota)+H_{a^{(n_0+1)}}(\iota).
\end{equation}

Define
\begin{equation}\label{uv:fused-regulator}
 \beta_{1s}=0,\qquad
 \beta_{rs}
 =w\wedge\left(2^r\sum_{1\le k<r}a^{(k)}_{r-k,s}\right)
 \quad(r\ge2).
\end{equation}
For fixed $r$ the sum is finite and decreases to zero as $s$
increases, and every $\beta_{rs}$ lies between $0$ and $w$.
Lemma~\ref{lt:fusion} proves that $\beta$ is a regulator and that
\begin{equation}\label{uv:fused-bound}
 w\wedge
 \sum_{k=1}^{N}
 H_{a^{(k)}}\bigl(r\mapsto\varphi(r+k)\bigr)
 \le H_\beta(\varphi)
 \qquad(N\ge1,\ \varphi\in\mathbb N^{\mathbb N}).
\end{equation}
Fix a selector $\varphi$ and put $\theta(r)=\varphi(r+1)$.
The outer regulator gives an integer $n$ with
$A_n\subseteq C(H_{a^{(1)}}(\theta))$. For
$\psi(r)=\varphi(r+n+1)$, the inner regulator gives
$t_\varphi\in T_n$ with
\[
 t\ge t_\varphi
 \quad\Longrightarrow\quad
 h^n_t\in C(H_{a^{(n+1)}}(\psi)).
\]
Using \eqref{uv:countable-absorption-premise} for this $n$ and
$t_\varphi$, followed by \eqref{uv:section-laws}, gives
\[
 z_i\in
 C\bigl(H_{a^{(1)}}(\theta)+H_{a^{(n+1)}}(\psi)\bigr)
 \quad\text{eventually}.
\]
Intersect this tail with the tail in \eqref{uv:attained-cap}.
The finite-meet identity and positivity of the omitted terms give
\[
 \begin{aligned}
 z_i&\in
 C\bigl(w\wedge
       [H_{a^{(1)}}(\theta)+H_{a^{(n+1)}}(\psi)]\bigr)\\
 &\subseteq
 C\left(w\wedge\sum_{k=1}^{n+1}
       H_{a^{(k)}}(r\mapsto\varphi(r+k))\right)
 \subseteq C(H_\beta(\varphi))
 \quad\text{eventually}.
 \end{aligned}
\]
Thus one regulator $\beta$ controls the candidate net in coordinate
$j$. This construction applies to every $j$, so the original
regulator criterion gives $z_i\to0$.

Now let $(A_\lambda)_{\lambda\in\Lambda}$ have a nonempty countable
directed index, and enumerate $\Lambda$ as $(\lambda_n)_{n\ge1}$,
repeating indices if $\Lambda$ is finite. Directedness gives
\[
 \nu_1=\lambda_1,\qquad
 \nu_{n+1}\ge\nu_n,\lambda_{n+1}\qquad(n\ge1).
\]
For every $\lambda\in\Lambda$, there is $m$ with $\lambda=\lambda_m$,
and then
\[
 n\ge m\quad\Longrightarrow\quad
 \nu_n\ge\lambda,\qquad A_{\nu_n}\subseteq A_\lambda.
\]
Conversely each $A_{\nu_n}$ is one member of the original family.
These inclusions identify the two saturated nets' eventual sets, and
in particular $\langle A_{\nu_n}\rangle\to0$. Put
$\widetilde A_n=A_{\nu_n}$,
$\widetilde T_n=T_{\nu_n}$ and
$\widetilde h^n_t=h^{\nu_n}_t$. Each $\widetilde h^n$ is null, and
the original absorption premise at $\lambda=\nu_n$ gives
\[
 \forall n\ \forall t_0\in\widetilde T_n\ \exists i(n,t_0)\
 \forall i\ge i(n,t_0):\quad
 z_i\in\widetilde A_n+
       \{\widetilde h^n_t:t\ge t_0\}.
\]
Thus the sequential calculation above, with these displayed families
and the unchanged candidate net $(z_i)$, supplies in every coordinate
$j$ a regulator $\beta^j$ and, for each selector $\varphi$, an index
$i(j,\varphi)$ such that
\[
 i\ge i(j,\varphi)\quad\Longrightarrow\quad
 z_i\in C_j(H_{\beta^j}(\varphi)).
\]
The original regulator criterion gives $z_i\to0$, proving
\textnormal{($P_{\aleph_0}$)} for every countable directed outer base.
\end{proof}

By retaining the regulators attached separately to arbitrary inner
nets, the proof of Theorem~\ref{uv:regulator-verification} makes
precise why countable fusion restricts the outer errors alone. It
reaches the original regulator convergence directly, without
identifying it with order convergence for arbitrary nets; when the
stronger hypotheses ensuring that identification are available,
Theorem~\ref{uv:order-verification} also supplies unrestricted
absorption.

\begin{corollary}[Topological, order and lattice-normed ranges]
\label{uv:standard-cases}
The preceding intrinsic properties hold in the following original
convergences.
\begin{enumerate}
\item Every separated topological abelian group satisfies
\textnormal{(R)} and \textnormal{(P)}.
\item Every vector lattice with order convergence satisfies
\textnormal{(R)} and \textnormal{(P)}. More generally, let
$q:Z\to F_+$ be a lattice norm, where $F$ is a vector lattice and
the original convergence is
\[
 z_\alpha\longrightarrow z
 \quad\Longleftrightarrow\quad q(z_\alpha-z)\longrightarrow_o0.
\]
Then $Z$ satisfies \textnormal{(R)} and \textnormal{(P)}.
\item If $F$ is Dedekind $\sigma$-complete and weakly
$\sigma$-distributive, the original convergence
\[
 z_\alpha\longrightarrow z
 \quad\Longleftrightarrow\quad q(z_\alpha-z)\longrightarrow_D0
\]
satisfies \textnormal{(R)} and \textnormal{($P_\omega$)}.
\end{enumerate}
In each case, sequential completeness of the range and countable
finite semivariance give the constructive uniform sequence theorem
of Corollary~\ref{cu:uniform-sequential}, and known-limit preservation
requires only the coefficient hypothesis. In item~(3), sequential
completeness is supplied by Proposition~\ref{lt:d-complete} when
$Z=F$ and $q(z)=|z|$.
\end{corollary}
\begin{proof}
For item~(1), every neighborhood $W$ of zero contains a closed
neighborhood. Indeed, there is a symmetric open neighborhood $V$
with $V+V\subseteq W$, and $\overline V\subseteq V+V$: if
$x\in\overline V$, the neighborhood $x+V$ meets $V$, so
$x\in V-V=V+V$. If $A_\lambda\downarrow$ is null, it is eventually
contained in $\overline V$, and so is $a_Z(A_\lambda)$, since
$\overline V$ contains the limits of all its nets. Varying $W$
proves \textnormal{(R)}. For absorption, fix $W$ and $V$ as above.
An outer level $\lambda$ satisfies $A_\lambda\subseteq V$, and
the associated null inner net has a tail contained in $V$. The
absorption premise then places the candidate net eventually in
$V+V\subseteq W$. This holds for every $W$, proving
\textnormal{(P)}.

For item~(2), use the original bound relations
\[
 \mathscr R(u)=\{(z,w):q(z-w)\le u\}\qquad(u\in F_+).
\]
The lattice-norm identities $q(0)=0$, $q(-z)=q(z)$ and
$q(z+w)\le q(z)+q(w)$ give the diagonal, symmetry and composition
laws. If an infimum of bounds exists, being below every bound is
equivalent to being below that infimum, so the intersection law
holds as well. The original convergence supplies the exact
pair criterion by its definition. Finally, for $D\downarrow0$ and
a lower bound $v$ of $r+D$,
\[
 v-r\le d\quad(d\in D)
 \quad\Longrightarrow\quad v-r\le0,
\]
so $\inf(r+D)=r$. Theorem~\ref{uv:order-verification} applies.
Taking $Z=F$ and $q(z)=|z|$ proves the native order case. Only
infima already supplied by decreasing zero controls were used,
so no Dedekind completeness is required for this stability assertion.

For item~(3), the defining equivalence for the original convergence
and the canonical difference criterion give, for every pair net,
\[
 (z_\alpha,w_\alpha)\longrightarrow\Delta_Z
 \quad\Longleftrightarrow\quad
 \exists a\ \forall\varphi\ \exists\alpha_\varphi\
 \forall\alpha\ge\alpha_\varphi:\quad
 q(z_\alpha-w_\alpha)\le H_a(\varphi).
\]
Equivalently, that pair net is eventually in
$\mathscr R(H_a(\varphi))$ for each selector of one regulator. The
relations already verified in item~(2) therefore meet every hypothesis
of Theorem~\ref{uv:regulator-verification}, which gives
\textnormal{(R)} and \textnormal{($P_\omega$)}.
The constructive assertion in all three cases is
Corollary~\ref{cu:uniform-sequential} with these verified target
properties. The known-limit assertion is
Theorem~\ref{uc:known-limit}. The coefficient condition remains the
independent finite-word condition in both applications.
\end{proof}

\subsection{Finite coefficient estimates}
\label{uv:coefficient-verification}

With the required operations in the range verified, it remains to
follow the errors through the coefficients, before any candidate
integral is constructed. For this purpose, the response family
$\Sigma$ retains the values of the actual reduced finite words, and
the next criterion relates input and output bounds using only the
finitely many input coordinates needed by each output estimate.
In these finite-word formulas, the join of an empty family of positive
bounds is $0$, and an empty record has response $0$.

\begin{proposition}[Finite bounds in the original convergence]
\label{uv:finite-set-certificate}
Let $E_i,F_j$ be the lattice-ordered controller groups of the
original input and output bounds. Let $U_i(u)\subseteq X$ and
$V_j(v)\subseteq Z$ be those
input and output bound sets, with $0\in U_i(u)$ and $0\in V_j(v)$,
isotone in their arguments, and suppose that originally null nets
are exactly those admitting the corresponding order controls.
For each $j$, let $I_j$ be finite and let
$\omega_j:\prod_{i\in I_j}(E_i)_+\to(F_j)_+$ be isotone with
\[
 u_\theta\downarrow0
 \quad\Longrightarrow\quad
 \omega_j(u_\theta)\downarrow0.
\]
If the actual reduced words satisfy
\begin{equation}\label{uv:finite-set-bound}
 \Sigma\left(\bigcap_{i\in I_j}U_i(u_i)\right)
 \subseteq V_j(\omega_j((u_i)_{i\in I_j})),
\end{equation}
then the coefficients have continuous finite semivariance.

Alternatively, suppose that the original input and output null
nets are detected by regulator bounds in these sets. The same
conclusion holds when the controller groups are Dedekind
$\sigma$-complete and each $\omega_j$ is isotone, subadditive,
vanishes at zero and preserves decreasing zero sequences.
\end{proposition}
\begin{proof}
Let $B_\eta\ni0$ be a decreasing originally null family. Fix $j$.
For each $i\in I_j$, the input criterion supplies controls
$u_{i,\theta_i}\downarrow0$ and indices $\eta(i,\theta_i)$ such that
\[
 \eta\ge\eta(i,\theta_i)
 \quad\Longrightarrow\quad
 B_\eta\subseteq U_i(u_{i,\theta_i}).
\]
On the finite product $\Theta=\prod_{i\in I_j}\Theta_i$ of the
control index sets put $u_\theta=(u_{i,\theta_i})_{i\in I_j}$.
This net decreases coordinatewise. If $v$ is a lower bound of all
$u_\theta$, then, for each fixed $i$, varying $\theta_i$ while fixing
the other finitely many coordinates gives
\[
 v_i\le u_{i,\theta_i}\quad(\theta_i\in\Theta_i),\qquad v_i\le0.
\]
Hence $u_\theta\downarrow0$ in the finite product controller. For each
$\theta$, directedness supplies an index $\eta_\theta$ with
$\eta_\theta\ge\eta(i,\theta_i)$ for every $i\in I_j$. Thus
\[
 \eta\ge\eta_\theta
 \quad\Longrightarrow\quad
 \Sigma(B_\eta)\subseteq V_j(\omega_j(u_\theta)).
\]
The controls $\omega_j(u_\theta)$ decrease to zero by assumption.
For every such control, the displayed inclusion bounds all values on
the response saturation tail determined by $\eta_\theta$.
The original output criterion, for every $j$, therefore gives
$\langle\Sigma(B_\eta)\rangle\to0$. If $I_j$ is empty, its product
controller is the one-element group $\{0\}$. Applying the assumed
preservation of decreasing zero nets to its constant zero net gives
$\omega_j(0)=\inf\{\omega_j(0)\}=0$. The empty intersection of
input bound sets is $X$, so \eqref{uv:finite-set-bound} gives
$\Sigma(X)\subseteq V_j(0)$ and covers every $B_\eta$.

For the regulator assertion, fix $j$ and let $a^i$ regulate the
input family in coordinate $i\in I_j$. Their finite product
$a_{rs}=(a^i_{rs})_{i\in I_j}$ is bounded by
$M=(M_i)_{i\in I_j}$, where $M_i$ bounds $a^i$, and
$\inf_s a_{rs}=0$ coordinatewise for each $r$. Hence it is a
regulator, with
\[
 H_a(\varphi)=\bigl(H_{a^i}(\varphi)\bigr)_{i\in I_j}.
\]
For a fixed selector $\varphi$, let $\eta(i,\varphi)$ be the
threshold supplied by input coordinate $i$, and take
$\eta_\varphi\ge\eta(i,\varphi)$ for every $i\in I_j$.
Then
\[
 B_\eta\subseteq
 \bigcap_{i\in I_j}U_i(H_{a^i}(\varphi))
 \quad(\eta\ge\eta_\varphi).
\]
Define
\[
 c_{rs}=\omega_j(M)\wedge2^r\omega_j(a_{rs}).
\]
Lemma~\ref{lt:positive-d-transport} gives a regulator $c$ with
$\omega_j(H_a(\varphi))\le H_c(\varphi)$ for every $\varphi$.
Equation~\eqref{uv:finite-set-bound} therefore yields
\[
 \Sigma(B_\eta)
 \subseteq V_j(\omega_j(H_a(\varphi)))
 \subseteq V_j(H_c(\varphi))
 \quad\text{eventually}.
\]
Thus the one regulator $c$, independent of the selector and of
$\eta$, supplies the output criterion in coordinate $j$, with a tail
threshold for every selector. Doing this for every $j$ proves
$\langle\Sigma(B_\eta)\rangle\to0$. For $I_j=\varnothing$, the
product regulator is zero, $\omega_j(0)=0$ by hypothesis, and
$c_{rs}=0$, so the same assertion holds with every threshold.
\end{proof}

\begin{proposition}[An intrinsic locally solid certificate]
\label{uc:intrinsic-positive-certificate}
Let $X,Z$ be abelian lattice-ordered groups with their original
locally solid group convergences. Suppose the convergence on $X$ is
finite-supremum stable in Definition~\ref{fv:finite-sup}, and let
$A:X\to Z$ be positive, additive and continuous. If every actual
finite coefficient word consists of positive additive maps and
\[
 \sum_rL_rx\le Ax\qquad(x\in X_+),
\]
then the coefficient system has continuous finite semivariance.
\end{proposition}
\begin{proof}
For $0\in B\subseteq X$, put
\[
 V(B)=\left\{\bigvee_{r=1}^m|x_r|:
               m\ge1,\ x_1,\ldots,x_m\in B\right\}.
\]
If $B_\eta$ is decreasing and originally null, finite-supremum
stability gives nullity of the saturated family $V(B_\eta)$.
Continuity of $A$ then gives nullity of $A[V(B_\eta)]$. For
$z=\sum_rL_rx_r$ with $x_r\in B_\eta$, set
$u=\bigvee_r|x_r|\in V(B_\eta)$. Positivity and additivity imply
\[
 |z|\le\sum_r|L_rx_r|
       \le\sum_rL_r|x_r|
       \le\sum_rL_ru\le Au.
\]
The zero response is covered by $u=0$, so
$\Sigma(B_\eta)\subseteq\operatorname{sol}(A[V(B_\eta)])$.
Lemma~\ref{fd:solid-tail} gives nullity of the saturated response
net. This verifies the coefficient property using the original
locally solid convergence and the stated finite-supremum assumption.
\end{proof}

Here finite-supremum stability remains an independent input
hypothesis, verified by Proposition~\ref{fv:gauges} under its
finite-join assumptions and, in the indicated setting, through the
band-covariance criterion of Proposition~\ref{fv:band-covariance}.
To obtain a certificate without that stability assumption, we return
to estimates on the finite words themselves, allowing different
controls on the input and output sides.

\begin{proposition}[Finite sufficient estimates]
\label{uc:finite-certificates}\label{lt:csv-order-test}
Suppose the original input and output order bounds are given by
$(p_i,E_i)$ and $(\rho_j,F_j)$. For each $j$, let $I_j$ be finite,
and let $\omega_j:\prod_{i\in I_j}(E_i)_+\to(F_j)_+$ be isotone
with $\omega_j(u_\theta)\downarrow0$ whenever $u_\theta\downarrow0$.
The inequalities
\begin{equation}\label{uc:finite-estimate}
 \rho_j\left(\sum_tL_t^{A,\gamma}x_t\right)
 \le\omega_j\left((\bigvee_tp_i(x_t))_{i\in I_j}\right)
\end{equation}
for the actual reduced records imply continuous finite semivariance.
\end{proposition}
\begin{proof}
Define
$U_i(u)=\{x:p_i(x)\le u\}$ and
$V_j(v)=\{z:\rho_j(z)\le v\}$. If every $x_t$ belongs to
$\bigcap_{i\in I_j}U_i(u_i)$, then
$\bigvee_tp_i(x_t)\le u_i$ for every $i\in I_j$.
The displayed finite inequality and isotonicity yield
\[
 \rho_j\left(\sum_tL_t^{A,\gamma}x_t\right)
 \le\omega_j((u_i)_{i\in I_j}).
\]
The zero response obeys the same bound since its left side is zero.
Thus \eqref{uv:finite-set-bound} holds, and
Proposition~\ref{uv:finite-set-certificate} proves the assertion.
\end{proof}

\begin{proposition}[Positive coefficient certificates]
\label{uc:positive-certificates}\label{lt:positive-budget}
\label{lt:other-budgets}
The following conditions supply \eqref{uc:finite-estimate}.
The bound maps $\rho_j$ in items~(2)--(3) are the original symmetric
subadditive output bounds.
\begin{enumerate}
\item Let $X,Z$ be abelian lattice-ordered groups with their order
convergences, let each $L_r$ be positive and additive, and suppose
$\sum_rL_r\le A$ for a positive additive map $A:X\to Z$. Then
\[
 \left|\sum_rL_rx_r\right|\le A\left(\bigvee_r|x_r|\right).
\]
\item Let $X=Z$, and suppose positive additive controller maps
$T_{j,r}$ satisfy
\[
 \rho_j(L_rx)\le T_{j,r}(\rho_j(x)),\qquad
 \sum_rT_{j,r}\le T_j.
\]
Then
\[
 \rho_j\left(\sum_rL_rx_r\right)
 \le T_j\left(\bigvee_r\rho_j(x_r)\right).
\]
\item Let $L_Cx=b(x,\mu C)$, where $b$ is biadditive and
$\mu:\mathcal H\to Y_+$. Let $\Omega\in\mathcal H$ be a root and
suppose that the actual coefficient records satisfy the finite
budget
\[
 \sum_r\mu C_r\le\mu\Omega.
\]
This budget is automatic for positive finitely additive $\mu$
whenever the record cells are disjoint and can be completed to a
finite partition of $\Omega$ by cells of $\mathcal H$, in particular
for disjoint cells contained in $\Omega$ when $\mathcal H$ is a ring.
Suppose also that
\[
 \rho_j(b(x,y))\le B_j((p_i(x))_{i\in I_j},y)\qquad(y\ge0),
\]
where $B_j$ is isotone in both variables and additive in its second
variable, with positive values. The second-variable isotonicity is
automatic from positivity and additivity when the coefficient cone is
that of a partially ordered abelian group. Then
\[
 \rho_j\left(\sum_rb(x_r,\mu C_r)\right)
 \le B_j\left((\bigvee_rp_i(x_r))_{i\in I_j},\mu\Omega\right).
\]
The second-variable isotonicity may be omitted when the budget has a
positive residual witness
\[
 \mu\Omega=\sum_r\mu C_r+v,\qquad v\in Y_+.
\]
In particular, finite partition completion provides such a witness.
\end{enumerate}
Order continuity at zero of the displayed right-hand map verifies
continuous finite semivariance. For regulator descriptions with
Dedekind $\sigma$-complete controller groups, isotone subadditive
right-hand maps that vanish at zero and preserve decreasing zero
sequences verify the corresponding assertion.
\end{proposition}
\begin{proof}
For item~(1), put $u=\bigvee_r|x_r|$. Positivity gives
$|L_rx_r|\le L_r|x_r|\le L_ru$, and hence
\[
 \left|\sum_rL_rx_r\right|
 \le\sum_r|L_rx_r|
 \le\sum_rL_ru\le Au.
\]
For item~(2), put $u=\bigvee_r\rho_j(x_r)$. Subadditivity,
positivity of the controller maps and the operator budget give
\[
 \rho_j\left(\sum_rL_rx_r\right)
 \le\sum_rT_{j,r}(\rho_j(x_r))
 \le\sum_rT_{j,r}(u)\le T_j(u).
\]
For item~(3), put $u=(\bigvee_rp_i(x_r))_{i\in I_j}$.
In the ordered-group case, $0\le y\le y'$ implies
$y'-y\in Y_+$ and
\[
 B_j(u,y')=B_j(u,y)+B_j(u,y'-y)\ge B_j(u,y),
\]
which verifies the stated automatic isotonicity. In the monoid case,
it is the explicit order-preserving assumption on the coefficient
variable. The actual finite budget then gives
\[
 \begin{aligned}
 \rho_j\left(\sum_rb(x_r,\mu C_r)\right)
 &\le\sum_r\rho_j(b(x_r,\mu C_r))\\
 &\le\sum_rB_j(u,\mu C_r)
 =B_j\left(u,\sum_r\mu C_r\right)
 \le B_j(u,\mu\Omega).
 \end{aligned}
\]
If the budget instead has the stated residual witness, its last
inequality is verified directly by
\[
 B_j(u,\mu\Omega)
 =B_j\left(u,\sum_r\mu C_r\right)+B_j(u,v)
 \ge B_j\left(u,\sum_r\mu C_r\right).
\]
To verify that budget from the stated partition condition, write
$\Omega=(\bigsqcup_rC_r)\sqcup(\bigsqcup_{\ell=1}^qD_\ell)$ with
all cells in $\mathcal H$. Positive finite additivity gives
\[
 \mu\Omega=\sum_r\mu C_r+\sum_{\ell=1}^q\mu D_\ell
 \ge\sum_r\mu C_r.
\]
When $\mathcal H$ is a ring and the record cells are disjoint,
$D_1=\Omega\setminus\bigcup_rC_r$ belongs to $\mathcal H$ and
supplies such a completion, with the empty remainder omitted.
Each inequality is an instance of \eqref{uc:finite-estimate}.
Proposition~\ref{uc:finite-certificates} proves the order assertion,
and the regulator part of
Proposition~\ref{uv:finite-set-certificate} proves the regulator
assertion using its explicitly constructed output regulator.
\end{proof}

Positive means are included in the operator certificate whenever
the original bound maps have the corresponding homogeneity, for a
finite scalar budget $w_r\ge0$, $\sum_rw_r\le M$ gives
\[
 \rho_j\left(\sum_rw_rx_r\right)
 \le M\bigvee_r\rho_j(x_r).
\]
A finite variation budget provides a further sufficient estimate:
if $B_j$ has the isotonicity and additivity of item~(3), and the
coefficient bound $q$ and the biadditive product satisfy
\[
 \rho_j(b(x,y))\le B_j((p_i(x))_i,q(y)),\qquad
 \sum_rq(\mu C_r)\le v,
\]
then finite subadditivity, isotonicity and additivity give
\[
 \rho_j\left(\sum_rb(x_r,\mu C_r)\right)
 \le B_j((\bigvee_rp_i(x_r))_i,v).
\]
This calculation uses only the displayed finite budget, leaving the
direct whole-word estimate \eqref{uv:finite-set-bound} available
whenever cancellation within the actual records gives a bound without
a variation budget.

\subsection{Finite responses without range subtraction}

The finite estimates just obtained suggest returning to the records
before taking differences: a comparison of two assignments can already
be expressed through their pair bounds. In this form the interchange
argument applies to finite sums and, more generally, to finite response
maps, even when subtraction is unavailable in the input or response
space.

\begin{theorem}[Finite response comparison]\label{mg:finite-response}
Let $X$ and $Z$ be separated uniform convergence spaces. Suppose
maps $e_i:X^2\to E_i$ take values in partially ordered sets with
least element $0$, with order nullity defined by downward zero
regulators, and satisfy
\[
 (x_\alpha,y_\alpha)\longrightarrow_{\mathcal U_X}\Delta_X
 \quad\Longleftrightarrow\quad
 e_i(x_\alpha,y_\alpha)\longrightarrow_o0\quad(i\in I)
\]
for every pair net.
Suppose $Z$ is
complete and has a faithful pair representation $d_j:Z^2\to M_j$
with controllers satisfying Definition~\ref{mg:controls} and
\eqref{mg:zero-translation}.
Let $P$ be the directed approximation index set, and let
$\Gamma_A$, $A\in\mathcal A\ne\varnothing$, be the prescribed
directed record sets. For every record, let
$T_{A,\gamma}\subseteq\Omega$ be its finite tag set, and let
\[
 \mathscr R_{A,\gamma}:X^{T_{A,\gamma}}\longrightarrow Z.
\]
For each $j$, suppose a finite nonempty set $I_j\subseteq I$ and an isotone map
\[
 \omega_j:\prod_{i\in I_j}E_i\longrightarrow M_j
\]
satisfy $\omega_j(0)=0$ and
\[
 u_\theta\downarrow0
 \Longrightarrow \omega_j(u_\theta)\downarrow0,
\]
where the product order is used on its domain. Assume, for every
record and every two assignments $x,y\in X^{T_{A,\gamma}}$, that
\begin{equation}\label{mg:finite-response-estimate}
 \left[
 e_i(x_t,y_t)\le u_i
 \quad(t\in T_{A,\gamma},\ i\in I_j)
 \right]
 \Longrightarrow
 d_j(\mathscr R_{A,\gamma}x,\mathscr R_{A,\gamma}y)
       \le\omega_j((u_i)_{i\in I_j}).
\end{equation}
If $f_p:\Omega\to X$ converges uniformly to $f$ and
\[
 \mathscr R_{A,\gamma}(f_p|_{T_{A,\gamma}})
       \longrightarrow y_p^A\qquad(p,A),
\]
then there are unique $z^A$ such that
\[
 \mathscr R_{A,\gamma}(f|_{T_{A,\gamma}})\longrightarrow z^A,
 \qquad y_p^A\longrightarrow z^A\quad\text{uniformly in }A.
\]
For each fixed record, its response to $f_p$ converges to its response
to $f$, so both iterated limits exist and equal $z^A$.
\end{theorem}
\begin{proof}
Fix $j$. For every $i\in I_j$, input pair faithfulness supplies a
regulator $u_{i,\theta_i}\downarrow0$ and indices $p(i,\theta_i)$
such that
\[
 p\ge p(i,\theta_i),\quad t\in\Omega
 \Longrightarrow e_i(f_p(t),f(t))\le u_{i,\theta_i}.
\]
On the finite product of the regulator index sets, put
$u_\theta=(u_{i,\theta_i})_{i\in I_j}$.
This tuple decreases to $0$: a lower bound has every coordinate
at most every member of the corresponding coordinate regulator,
and therefore every coordinate equals $0$.
Define
\[
 P^j_\theta=\{p:p\ge p(i,\theta_i)\ (i\in I_j)\},
 \qquad \beta^j_\theta=\omega_j(u_\theta).
\]
Directedness of $P$ and finiteness of $I_j$ make $P^j_\theta$
contain a tail. Its definition is independent of $A,\gamma$.
For $p\in P^j_\theta$, \eqref{mg:finite-response-estimate} gives
\[
 d_j\bigl(\mathscr R_{A,\gamma}(f|_{T_{A,\gamma}}),
          \mathscr R_{A,\gamma}(f_p|_{T_{A,\gamma}})\bigr)
 \le\beta^j_\theta
 \qquad(A,\gamma).
\]
The net $\beta^j_\theta$ decreases to zero by the hypothesis on
$\omega_j$. Apply Theorem~\ref{lt:order-transfer} with
\[
 \begin{gathered}
 x^A_{p,\gamma}
   =\mathscr R_{A,\gamma}(f_p|_{T_{A,\gamma}}),\qquad
 x^A_\gamma=\mathscr R_{A,\gamma}(f|_{T_{A,\gamma}}),\\
 r^{A,j}_{p,\gamma}=0,
 \end{gathered}
\]
and the given endpoints $y_p^A$. In this application, the assumed
response convergence supplies the known-row condition, while the
preceding finite estimate supplies the comparison for every inner
index, with the outer sets $P^j_\theta$ independent of $A$.
Since the residual is the constant zero net, it is order null;
the target completeness and controller laws are those assumed here.
All hypotheses of the theorem are therefore verified, giving unique
$z^A$ with
$x^A_\gamma\to z^A$ and $y_p^A\to z^A$ uniformly in $A$.

For each fixed $A,\gamma,j$, every $\theta$ has an outer tail in
$P^j_\theta$, on which
$d_j(x^A_\gamma,x^A_{p,\gamma})\le\beta^j_\theta$.
Thus this distance is order null in $p$. Pair faithfulness and
symmetry give $x^A_{p,\gamma}\to x^A_\gamma$. The two iterated
limits can now be written explicitly as
\[
 \lim_p\lim_\gamma x^A_{p,\gamma}
 =\lim_p y_p^A=z^A,
 \qquad
 \lim_\gamma\lim_p x^A_{p,\gamma}
 =\lim_\gamma x^A_\gamma=z^A.
\]
\end{proof}

For additive finite responses, the preceding estimate can itself be
verified coefficient by coefficient, provided that the separate
bounds fit within one budget independent of the record.

\begin{proposition}[Finite monoid coefficient budgets]\label{mg:finite-budget}
Suppose $Z$ is a commutative monoid and its representing distances
satisfy
\[
 d_j(z+w,z'+w')\le d_j(z,z')\oplus_j d_j(w,w').
\]
Let $e:X^2\to E$ be an input control. For every finite record let
$L_r:X\to Z$ and isotone maps $T_{j,r}:E\to M_j$ satisfy
\[
 d_j(L_rx,L_ry)\le T_{j,r}(e(x,y)),\qquad
 \bigoplus_{r=1}^{m}T_{j,r}(u)\le T_j(u)\quad(u\in E),
\]
where $T_j$ is isotone and independent of the record, $T_j(0)=0$, and
$u_\theta\downarrow0\Longrightarrow T_j(u_\theta)\downarrow0$.
Then the responses $\mathscr R(x)=\sum_{r=1}^{m}L_rx_r$ satisfy
\eqref{mg:finite-response-estimate} with $\omega_j=T_j$.
The same assertion holds for a finite tuple of input controls,
with $E$ their product.
\end{proposition}
\begin{proof}
For an empty record the two responses equal the monoid identity
and their distance is $0$. For a one-term record the asserted
bound is the individual coefficient estimate. If the finite-sum
bound holds for $m-1$ terms, nonexpansiveness of addition gives
\begin{align*}
 d_j\left(\sum_{r=1}^{m}L_rx_r,\sum_{r=1}^{m}L_ry_r\right)
 &\le d_j\left(\sum_{r=1}^{m-1}L_rx_r,
                \sum_{r=1}^{m-1}L_ry_r\right)
           \oplus_j d_j(L_mx_m,L_my_m)\\
 &\le\bigoplus_{r=1}^{m}T_{j,r}(e(x_r,y_r)).
\end{align*}
Induction proves this inequality for every finite record. If
$e(x_r,y_r)\le u$ for every $r$, isotonicity and the budget give
\[
 d_j\left(\sum_rL_rx_r,\sum_rL_ry_r\right)
 \le\bigoplus_rT_{j,r}(e(x_r,y_r))
 \le\bigoplus_rT_{j,r}(u)\le T_j(u).
\]
All inequalities remain valid when $u$ is a tuple and the comparisons
in the arguments are coordinatewise. This proves the assertion.
\end{proof}

\subsection{An algebraic source of known integrable rows}

Before applying the interchange theorems, it is useful to identify
a class of approximants whose integrals are already determined by the
finite coefficient algebra. On partition records, refinement
makes the sums of simple functions exact; by measuring the oscillation left
within each cell, we shall pass from that algebraic core to integrable
functions through estimates on attainable finite partitions.

\begin{definition}[Finite coefficient core]\label{lt:finite-core-system}
Let $\mathcal H$ be an algebra on $\Omega$ and
$L:\mathcal H\to\operatorname{Hom}(X,Z)$ satisfy
\[
 L_\varnothing=0,\qquad L_{C\sqcup D}=L_C+L_D.
\]
For every $A\in\mathcal A\subseteq\mathcal H$, suppose the
original records are finite internally tagged partitions of $A$.
Partitions in this subsection have nonempty cells.
For every finite measurable partition $\Pi$ of $\Omega$, suppose
\begin{equation}\label{lt:refinement-tail}
 T_A(\Pi)=\{\gamma:\text{every cell of the record at }\gamma
                         \text{ is contained in a member of }\Pi\}
\end{equation}
contains a tail of $\Gamma_A$.
\end{definition}

\begin{definition}[Finite packing estimate]\label{lt:packing-control}
Let $p_i:X\to(E_i)_+$ and $\rho_j:Z\to(F_j)_+$ be faithful order
representations. For every $j$, let $I_j$ be finite and let
$\kappa_j:\mathcal H\times\prod_{i\in I_j}(E_i)_+\to(F_j)_+$
satisfy $\kappa_j(C,0)=\kappa_j(\varnothing,u)=0$, isotonicity in $u$,
and
\begin{align}
 \rho_j(L_Cx)&\le\kappa_j(C,(p_i(x))_i),\label{lt:kappa-bound}\\
 \sum_{r=1}^m\kappa_j(D_r,u)&\le\kappa_j(C,u)
 \quad(D_r\subseteq C\text{ disjoint}).\label{lt:kappa-packing}
\end{align}
An oscillation datum for $f$ is a finite partition $\Pi$ and a tuple
$(u_C)_{C\in\Pi}$ satisfying
\[
 p_i(f(s)-f(t))\le(u_C)_i
 \quad(s,t\in C, i\in I_j).
\]
Its cost is $c_j(\Pi,u)=\sum_{C\in\Pi}\kappa_j(C,u_C)$.
Denote all attainable costs by $\mathscr C_j(f)$.
\end{definition}

\begin{lemma}[Finite refinement and reproduction]\label{lt:finite-core}
Every nonempty $\mathscr C_j(f)$ is downward directed.
For $s=\sum_{C\in\Pi}a_C\mathbf1_C$, the element
\[
 J_A(s)=\sum_{C\in\Pi}L_{A\cap C}a_C
\]
is independent of its representation, and
\begin{equation}\label{lt:core-reproduction}
 \gamma\in T_A(\Pi)\quad\Longrightarrow\quad S_\gamma^As=J_A(s).
\end{equation}
\end{lemma}
\begin{proof}
For two data $(\Pi,u),(\Xi,v)$, let
\[
 \Theta=\{C\cap B\ne\varnothing:C\in\Pi,\ B\in\Xi\},\qquad
 w_{C\cap B}=u_C\wedge v_B.
\]
The tuple $w$ bounds the oscillations on $\Theta$, and
\[
 \sum_{C,B}\kappa_j(C\cap B,u_C\wedge v_B)
 \le\sum_C\sum_B\kappa_j(C\cap B,u_C)
 \le\sum_C\kappa_j(C,u_C).
\]
Interchanging $C,u$ and $B,v$ gives the other cost bound.
For a record refining $\Pi$,
\[
 S_\gamma^As
 =\sum_{C\in\Pi}\sum_{D_r\subseteq C}L_{D_r}a_C
 =\sum_{C\in\Pi}L_{A\cap C}a_C.
\]
If also $s=\sum_{B\in\Xi}b_B\mathbf1_B$, then
$a_C=b_B$ whenever $C\cap B\ne\varnothing$, since both equal
$s(t)$ at any $t\in C\cap B$. Finite additivity gives
\[
 \begin{aligned}
 \sum_{C\in\Pi}L_{A\cap C}a_C
 &=\sum_{\substack{C\in\Pi,\,B\in\Xi\\C\cap B\ne\varnothing}}
       L_{A\cap C\cap B}a_C\\
 &=\sum_{\substack{C\in\Pi,\,B\in\Xi\\C\cap B\ne\varnothing}}
       L_{A\cap C\cap B}b_B
 =\sum_{B\in\Xi}L_{A\cap B}b_B.
 \end{aligned}
\]
Empty intersections contribute $L_\varnothing=0$, so this calculation
also covers cells disjoint from $A$.
\end{proof}

Once refinement has made the attainable costs downward directed,
the existence criterion asks whether their infimum is zero. This
zero infimum will provide a Cauchy control
for the original sums and retain an explicit estimate for the resulting
integral.

\begin{theorem}[Existence from attainable oscillations]
\label{lt:oscillation-existence}
Suppose $Z$ is complete in its original uniform convergence and
faithfully order represented. Under Definitions~\ref{lt:finite-core-system}
and \ref{lt:packing-control}, assume
\begin{equation}\label{lt:oscillation-infimum}
 \mathscr C_j(f)\ne\varnothing,\qquad
 \inf\mathscr C_j(f)=0\qquad(j\in J).
\end{equation}
Then $f$ is integrable on every $A\in\mathcal A$. For every
oscillation datum and every tuple
$(a_C)_{C\in\Pi}\in\prod_{C\in\Pi}f[C]$,
\begin{align}
 \rho_j\left(I_Af-\sum_{C\in\Pi}L_{A\cap C}a_C\right)
 &\le\sum_{C\in\Pi}\kappa_j(A\cap C,u_C)
 \le c_j(\Pi,u),\label{lt:oscillation-integral-bound}\\
 \gamma\in T_A(\Pi)\quad\Longrightarrow\quad
 \rho_j(S_\gamma^Af-I_Af)&\le2c_j(\Pi,u).
 \label{lt:oscillation-sum-bound}
\end{align}
The latter convergence is uniform in $A$ when records are restricted
by the same finite partition $\Pi$.
\end{theorem}
\begin{proof}
Put $J_A=\sum_{C\in\Pi}L_{A\cap C}a_C$. For
$\gamma\in T_A(\Pi)$,
\begin{align*}
 S_\gamma^Af-J_A
 &=\sum_{C\in\Pi}\sum_{D_r\subseteq C}L_{D_r}(f(t_r)-a_C),\\
 \rho_j(S_\gamma^Af-J_A)
 &\le\sum_{C\in\Pi}\sum_{D_r\subseteq C}
        \kappa_j(D_r,(p_i(f(t_r)-a_C))_i)\\
 &\le\sum_{C\in\Pi}\kappa_j(A\cap C,u_C)
 \le c_j(\Pi,u).
\end{align*}
Hence
\[
 \gamma,\gamma'\in T_A(\Pi)
 \Longrightarrow
 \rho_j(S_\gamma^Af-S_{\gamma'}^Af)\le2c_j(\Pi,u).
\]
By Lemma~\ref{lt:finite-core}, these costs are downward directed;
to use their doubles as controls, we verify the whole-group infimum
of $2\mathscr C_j(f)$. Let
$b\in F_j$ satisfy $b\le2c$ for every cost $c$. For any costs
$c,d$, a cost $e\le c,d$ gives $b\le2e\le c+d$. Hence
\[
 b-d\le c\quad(c\in\mathscr C_j(f)),\qquad
 b-d\le0\quad(d\in\mathscr C_j(f)),\qquad b\le0.
\]
Since every cost is positive, $0$ is a lower bound and
$\inf 2\mathscr C_j(f)=0$. Indexing these bounds by reverse order
gives a decreasing zero regulator. For each of its values $2c$, a
datum realizing $c$ and a tail contained in $T_A(\Pi)$ give the
corresponding eventual bound for the difference net. Thus
\[
 \rho_j(S_\gamma^Af-S_{\gamma'}^Af)\longrightarrow_o0
 \qquad((\gamma,\gamma')\in\Gamma_A^2)
\]
for every $j$. Pair faithfulness makes the original sum net Cauchy,
and completeness gives its unique original limit $I_Af$.

Let $v_\eta\downarrow0$ control $\rho_j(I_Af-S_\gamma^Af)$.
On the intersection of its controlling tail with $T_A(\Pi)$,
\[
 \rho_j(I_Af-J_A)
 \le v_\eta+\sum_{C\in\Pi}\kappa_j(A\cap C,u_C).
\]
Subtracting the fixed finite coefficient bound and taking the
infimum of $v_\eta$ gives
\[
 \rho_j(I_Af-J_A)
 -\sum_{C\in\Pi}\kappa_j(A\cap C,u_C)\le\inf_\eta v_\eta=0.
\]
This proves \eqref{lt:oscillation-integral-bound}. For
$\gamma\in T_A(\Pi)$, the triangle inequality now gives
\[
 \rho_j(S_\gamma^Af-I_Af)
 \le\rho_j(S_\gamma^Af-J_A)+\rho_j(J_A-I_Af)
 \le2c_j(\Pi,u),
\]
which is \eqref{lt:oscillation-sum-bound}.

For each $\Pi$, let $E_\Pi$ contain all pairs
$(S_\gamma^Af,I_Af)$ with $A\in\mathcal A$ and
$\gamma\in T_A(\Pi)$. Finite common refinement makes these sets
decreasing under partition refinement. For every attainable cost
$c_j(\Pi,u)$, all pairs in $E_\Pi$ satisfy the corresponding bound
$2c_j(\Pi,u)$. The decreasing zero regulator $2\mathscr C_j(f)$ just constructed
bounds this pair saturation: for a cost realized by $(\Pi,u)$,
every partition $\Theta$ refining $\Pi$ has
$E_\Theta\subseteq E_\Pi$, and all pairs in $E_\Pi$ have represented
distance at most $2c_j(\Pi,u)$. Thus each regulator value bounds an
actual partition tail, independently of $A$. Pair faithfulness gives
original uniform nullity of the pair saturation.
\end{proof}

To apply this criterion from the input data, one can either make
the oscillation small on every cell or allow larger oscillation on
an exceptional union whose coefficient cost becomes small. Both
possibilities are expressed by the following finite partition tests.

\begin{proposition}[Input tests for attainable oscillations]
\label{lt:oscillation-tests}
Each of the following conditions implies
\eqref{lt:oscillation-infimum}.
\begin{enumerate}
\item There is a decreasing family $u_\alpha\downarrow0$ such that
$u_\alpha$ bounds the oscillations on every cell of some finite
partition and $\kappa_j(\Omega,u_\alpha)\downarrow0$.
\item For every $\alpha,\beta$, there is a finite partition whose
good cells have oscillations bounded by $u_\alpha$, whose exceptional
cells have oscillations bounded by $v$, and whose exceptional union
$E$ satisfies
\[
 \kappa_j(\Omega,u_\alpha)\le a_\alpha\downarrow0,\qquad
 \kappa_j(E,v)\le b_\beta\downarrow0.
\]
\end{enumerate}
\end{proposition}
\begin{proof}
In the first case, finite packing gives
$c_j(\Pi,u)\le\kappa_j(\Omega,u_\alpha)$.
In the second case,
\[
 c_j(\Pi,u)
 \le\kappa_j(\Omega,u_\alpha)+\kappa_j(E,v)
 \le a_\alpha+b_\beta.
\]
At least one datum exists in each case, so the cost set is nonempty.
Let $c\in F_j$ be a lower bound of all attainable costs. In case~(1),
\[
 c\le\kappa_j(\Omega,u_\alpha)\quad(\alpha)
 \quad\Longrightarrow\quad c\le0.
\]
In case~(2), $c\le a_\alpha+b_\beta$ for every $\alpha,\beta$,
and translation in the lattice group gives
\[
 c-b_\beta\le a_\alpha\quad(\alpha)
 \Longrightarrow c\le b_\beta\quad(\beta)
 \Longrightarrow c\le0.
\]
All costs are positive, so $0$ is itself a lower bound. These
calculations prove the required infimum in both cases.
\end{proof}

For the positive measure certificate,
$\kappa_j(C,u)=B_j(u,\mu C)$ satisfies finite packing. The second
criterion is verified by $\mu E\le\varepsilon_\beta\downarrow0$
and the relevant order continuity, through
\[
 c_j(\Pi,u)
 \le B_j(u_\alpha,\mu\Omega)+B_j(v,\varepsilon_\beta).
\]
Thus the passage from the algebraic core to an integral is governed
by cell oscillations and coefficient packings: the preceding estimates
first make the sums Cauchy and then yield existence in the original
convergence.

\subsection{Uniform approximation and normalization}

With both coefficient estimates and a source of integrable rows
available, we return to uniform approximation on the prescribed nets.
The cardinal theorem matches the strength of the target assumptions
to the available error base, and its unrestricted form follows from
the ordered verification when the integration data are given by such
bounds. After recording this consequence, we turn to normalization,
whose comparisons retain independently null residuals.

\begin{theorem}[Uniform approximation on the prescribed nets]
\label{uc:uniform}\label{lt:uniform}\label{mg:uniform-extension}
Let $X,Z$ be separated abelian convergence groups, and suppose that
$Z$ is complete and satisfies either of the following conditions:
\begin{enumerate}
\item its canonical group pair convergence has a faithful monoid
order representation satisfying \eqref{mg:zero-translation}, or
\item its original convergence satisfies \textnormal{(R)} and
\textnormal{(P)}.
\end{enumerate}
Suppose the coefficients have continuous finite semivariance, every
$f_a$ is integrable on every $A\in\mathcal A$, and
\[
 U_{a_0}=\{0\}\cup
 \{f_a(t)-f(t):a\ge a_0,\ t\in\Omega\}
\]
is a null family in the original convergence of $X$. Then $f$ is
integrable on every $A$, the integrals $I_Af_a$ converge to $I_Af$
uniformly in $A$, and, for each fixed $A$,
\[
 \lim_a\lim_\gamma S_\gamma^Af_a
 =I_Af
 =\lim_\gamma\lim_aS_\gamma^Af_a.
\]
Closed-base stability is sufficient for the adherence assumption in
condition \textnormal{(2)}. By Theorem~\ref{nx:uniform-limit},
uniform-limit stability of continuous maps on every convergence
source is an equivalent replacement for its absorption assumption.
\end{theorem}
\begin{proof}
Let $P$ denote the directed index of the approximants. Under
condition~(1), Theorem~\ref{uv:order-verification} gives
\textnormal{(R)} and \textnormal{(P)} for the original convergence
of $Z$, so it suffices to use these two intrinsic properties in both
cases. The sets $U_{a_0}$ decrease on $P$, contain zero, and are null
by assumption. For $a\ge a_0$ and any original record, biadditivity
gives the finite identity
\[
 S_\gamma^Af_a-S_\gamma^Af
 =\sum_{t\in T_{A,\gamma}}
         L_t^{A,\gamma}(f_a(t)-f(t))\in\Sigma(U_{a_0}).
\]
Here $L_t^{A,\gamma}$ is the sum of the coefficient maps attached
to the occurrences of the tag $t$, so each tag is assigned a single
input value. Continuous finite semivariance gives
\[
 \langle\Sigma(U_{a_0}):a_0\in P\rangle\longrightarrow0.
\]
Set $\kappa=\max\{\aleph_0,|P|\}$. Completeness gives
$\kappa$-completeness, and \textnormal{(R)}, \textnormal{(P)} imply
\textnormal{$(R_\kappa)$}, \textnormal{$(P_\kappa)$}. In
Theorem~\ref{cu:uniform-cardinal}, take the error index to be $P$,
$B_{a_0}=U_{a_0}$, and the approximation threshold at $a_0$ to be
$a_0$ itself. Its input condition reads
\[
 a\ge a_0\quad\Longrightarrow\quad
 f_a(t)-f(t)\in B_{a_0}\qquad(t\in\Omega),
\]
which is the defining inclusion for $U_{a_0}$. Its known-row
condition is the assumed convergence
$S_\gamma^Af_a\to I_Af_a$ on each prescribed $\Gamma_A$.
The theorem therefore constructs $I_Af$ with
\[
 S_\gamma^Af\longrightarrow I_Af,
 \qquad I_Af_a\longrightarrow I_Af
                 \quad\text{uniformly in }A.
\]

For a fixed $A$ and a fixed record $\gamma$, uniform convergence of
$f_a$ implies $f_a(t)\to f(t)$ at each of its finitely many tags.
Continuity of the partial coefficient maps and of finite addition
then gives
\[
 S_\gamma^Af_a
 =\sum_{r\in I^A_\gamma}
       f_a(t^A_{\gamma r})\bullet\mu(C^A_{\gamma r})
 \longrightarrow
 \sum_{r\in I^A_\gamma}
       f(t^A_{\gamma r})\bullet\mu(C^A_{\gamma r})
 =S_\gamma^Af.
\]
Taking the limit in $\gamma$ of this identity gives the right
iterated limit. Taking the limit in $a$ of the known integrals gives
the left iterated limit. Both equal the constructed $I_Af$.
\end{proof}

Uniform control of the finite responses is one way to compare the
sums; another is to follow their coefficient masses to the prescribed
set function. Constant and indicator normalization record precisely
these elementary passages to the limit.

\begin{definition}[Elementary coefficient normalization]
\label{uc:normalization}
Let $b:X\times Y\to F$ be bilinear and $\mu:\mathcal H\to Y$.
Constant normalization on $A$ is
\[
 S_\gamma^A\underline{x}\longrightarrow_D b(x,\mu A)
 \qquad(x\in X).
\]
Indicator normalization on a root $\Omega$ is
\[
 S_\gamma^\Omega(x\mathbf1_A)\longrightarrow_D b(x,\mu A)
 \qquad(x\in X,\ A\in\mathcal A).
\]
\end{definition}

\begin{proposition}[Coefficient verification of normalization]
\label{uc:normalization-certificate}
Assume that $y\mapsto b(x,y)$ is continuous from the original convergence
of $Y$ to $D$-convergence in $F$, for every $x\in X$. Regular coefficient
masses of Definition~\ref{nr:regular-integrator} imply the corresponding
constant and indicator normalizations. Exact indicator resolution verifies
both by finite identities.
\end{proposition}
\begin{proof}
For the original record on $A$, biadditivity gives
\[
 S_\gamma^A\underline{x}
 =b\!\left(x,m_\gamma^A(\Omega)\right),\qquad
 m_\gamma^A(\Omega)\longrightarrow\mu A
 \Longrightarrow S_\gamma^A\underline{x}\longrightarrow_D b(x,\mu A).
\]
For the root record,
\[
 S_\gamma^\Omega(x\mathbf1_A)
 =b\!\left(x,m_\gamma^\Omega(A)\right),\qquad
 m_\gamma^\Omega(A)\longrightarrow\mu A
 \Longrightarrow S_\gamma^\Omega(x\mathbf1_A)\longrightarrow_D b(x,\mu A).
\]
If $\mu$ is finitely additive and the records are full disjoint
partitions, then
\[
 m_\gamma^A(\Omega)=\sum_r\mu C_{\gamma r}=\mu A.
\]
If the root record also refines $\{A,\Omega\setminus A\}$ and its tags
are internal, then
\[
 m_\gamma^\Omega(A)
 =\sum_{t_{\gamma r}\in A}\mu C_{\gamma r}
 =\mu\!\left(\bigsqcup_{t_{\gamma r}\in A}C_{\gamma r}\right)=\mu A.
\]
These equalities prove the exact assertions without a limiting argument.
\end{proof}

For positive bilinear products, normalization turns the finite sum
of an order bound into a moving majorant with a known limit. The
regulator interchange theorem can therefore construct the integral
while preserving the resulting sharp endpoint estimate.

\begin{corollary}[Order-uniform net Riemann convergence]
\label{uc:monteiro-net}\label{lt:monteiro-corollary}\label{uniformconvergencetheoremriemann}
Let $X,Y$ be Riesz spaces and $F$ be Dedekind $\sigma$-complete and WSD, with its original
$D$-convergence. Suppose $b:X\times Y\to F$ is positive bilinear and
\[
 x_k\uparrow x,\ y\ge0\Longrightarrow b(x_k,y)\uparrow b(x,y).
\]
Let $\mu$ be positive and satisfy constant normalization. Let $f_n$
be integrable on $A$, and suppose
\[
 u_n\downarrow0,\qquad |f_n(t)-f(t)|\le u_n\quad(n,t).
\]
Then $f$ is integrable on the original net and
\begin{equation}\label{uc:monteiro-sharp}
 |I_Af_n-I_Af|\le b(u_n,\mu A)\downarrow0.
\end{equation}
If constant normalization and integrability of the approximants hold
on every $A\in\mathcal A$, and indicator normalization holds on
$\Omega$, then
\[
 |I_Af_n-I_Af|\le b(u_n,\mu A)\le b(u_n,\mu\Omega)\downarrow0
 \quad(A\in\mathcal A).
\]
\end{corollary}
\begin{proof}
For fixed $A,n$, put
\[
 W_{n,\gamma}^A=\sum_rb(u_n,\mu C_{\gamma r}),\qquad
 v_n^A=b(u_n,\mu A),\qquad
 r_{n,\gamma}^A=|W_{n,\gamma}^A-v_n^A|.
\]
Positivity gives
\[
 |S_\gamma^Af_n-S_\gamma^Af|
 \le\sum_r|b(f_n(t_{\gamma r})-f(t_{\gamma r}),\mu C_{\gamma r})|
 \le W_{n,\gamma}^A\le v_n^A+r_{n,\gamma}^A.
\]
Constant normalization gives $r_{n,\gamma}^A\to_D0$. Moreover,
\[
 0\le u_1-u_n\uparrow u_1,
 \qquad
 v_n^A=b(u_1,\mu A)-b(u_1-u_n,\mu A)\downarrow0.
\]
Proposition~\ref{lt:d-complete} proves sequential completeness of $F$.
Theorem~\ref{lt:d-transfer} applies with
\[
 d(z,w)=|z-w|,\quad x_{n,\gamma}=S_\gamma^Af_n,\quad
 x_\gamma=S_\gamma^Af,\quad y_n=I_Af_n,\quad b_n=v_n^A.
\]
It proves candidate integrability and \eqref{uc:monteiro-sharp}.

For the family assertion, on root records positivity gives
\[
 0\le S_\gamma^\Omega(u_n\mathbf1_A)
 \le S_\gamma^\Omega\underline{u_n}.
\]
Put $w_n=b(u_n,\mu\Omega)$ and
\[
 e_\gamma=
 |S_\gamma^\Omega(u_n\mathbf1_A)-v_n^A|
 +|S_\gamma^\Omega\underline{u_n}-w_n|.
\]
The two normalization assumptions and Lemma~\ref{lt:d-calculus}
give $e_\gamma\to_o0$. The preceding finite inequalities imply
\[
 -e_\gamma\le v_n^A\le w_n+e_\gamma.
\]
Lemma~\ref{lt:fixed-bound}, applied to each inequality, gives
$0\le v_n^A\le w_n$. The sequential normality of $b$ already used
above gives $w_n\downarrow0$. Thus for every $A,n$ the original
finite comparison strengthens to
\[
 |S_\gamma^Af-S_\gamma^Af_n|
 \le w_n+r_{n,\gamma}^A,
 \qquad r_{n,\gamma}^A\longrightarrow_D0,
\]
with $w_n$ independent of $A$. The parameterized assertion of
Theorem~\ref{lt:d-transfer} now gives uniform endpoint convergence
and the bound $|I_Af_n-I_Af|\le w_n$. The earlier application on
each fixed $A$ retains the sharper bound by $v_n^A$.
\end{proof}

\begin{corollary}[Monteiro--Fernandez uniform convergence]
\label{uc:monteiro-kurzweil}
Let $F$ be a Dedekind $\sigma$-complete WSD Riesz space and
$[a,b]$ a compact interval. If every $f_n:[a,b]\to F$ is
Kurzweil integrable and
\[
 u_n\downarrow0,\qquad |f_n(t)-f(t)|\le u_n\quad(n\ge1,\ t\in[a,b]),
\]
then $f$ is Kurzweil integrable and
\[
 \left|\int_a^bf_n-\int_a^bf\right|\le(b-a)u_n\downarrow0.
\]
This is the uniform convergence assertion of
\cite[Theorem~3.3]{monteirouniform}.
\end{corollary}
\begin{proof}
For $a=b$ all sums and integrals are zero. Henceforth assume $a<b$.
A local constraint at $t$ belongs to $D_t=(0,\infty]$, ordered by
the usual order, with meet $\min$ and largest element $\infty$. Put
\[
 \mathscr V_t(d)=\{[c,e]:t\in[c,e]\subset(t-d,t+d)\},
\]
where $\mathscr V_t(\infty)$ contains every closed subinterval carrying
$t$. Finite-valued gauges are cofinal among these local constraints.
A gauge is a map $\delta:[a,b]\to(0,\infty)$. Let $\mathscr D(\delta)$
be the family of finite tagged divisions
\[
 \begin{gathered}
 P=(([t_{r-1},t_r],\xi_r))_{r=1}^m,\qquad
 a=t_0<\cdots<t_m=b,\\
 \xi_r\in[t_{r-1},t_r],\qquad
 [t_{r-1},t_r]\subset
       (\xi_r-\delta(\xi_r),\xi_r+\delta(\xi_r)).
 \end{gathered}
\]
Bisecting each closed interval gives the branching system of
Proposition~\ref{ls:branch}. Along every infinite branch
$C_n=[a_n,b_n]$,
\[
 a_n\uparrow t,\qquad b_n\downarrow t,\qquad
 b_n-a_n=2^{-n}(b-a),\qquad
 t\in\bigcap_n C_n.
\]
For every $d>0$, some $n$ satisfies $b_n-a_n<d$, and then
$C_n\subset(t-d,t+d)$. Thus the branch condition holds for the local
constraints $C\subset(t-d,t+d)$, and Proposition~\ref{ls:branch} gives
$\mathscr D(\delta)\ne\varnothing$. Set
\[
 \Gamma=\{(\delta,P,k):P\in\mathscr D(\delta),\ k\in\mathbb N\},
 \qquad
 S_{(\delta,P,k)}h=\sum_{r=1}^m(t_r-t_{r-1})h(\xi_r).
\]
Besides equality, order these indices by
\[
 (\delta,P,k)<(\varepsilon,Q,l)
 \iff \varepsilon(t)\le\delta(t)\ (t\in[a,b]),\quad l>k.
\]
The pointwise minimum of two gauges and any larger integer supply a
common upper index, so $\Gamma$ is directed. For every
$\gamma_0=(\delta_0,P_0,k_0)$,
\[
 \{S_\gamma h:\gamma\ge\gamma_0\}
 =\left\{\sum_{r=1}^m(t_r-t_{r-1})h(\xi_r):
                                  P\in\mathscr D(\delta_0)\right\}.
\]
The reverse inclusion uses $(\delta_0,P,k_0+1)$ for each indicated $P$.
Consequently, for every $h$ and $I\in F$,
\[
 \begin{aligned}
 S_\gamma h\longrightarrow_D I
 \iff{}&\exists\text{ regulator }c\ \forall\varphi\in\mathbb N^{\mathbb N}
        \ \exists\delta\ \forall P\in\mathscr D(\delta):\\[-1mm]
 &\left|\sum_{r=1}^m(t_r-t_{r-1})h(\xi_r)-I\right|
      \le H_c(\varphi).
 \end{aligned}
\]
The right side is the Kurzweil definition in
\cite[Section~2, equation~(2)]{monteirouniform}.
Thus each known row $S_\gamma f_n$ has its prescribed integral as limit.
For $x\in F$,
\[
 S_\gamma\underline{x}
 =\left(\sum_{r=1}^m(t_r-t_{r-1})\right)x=(b-a)x.
\]
For every $n,\gamma$,
\[
 |S_\gamma f_n-S_\gamma f|
 \le\sum_{r=1}^m(t_r-t_{r-1})u_n=(b-a)u_n.
\]
Theorem~\ref{lt:d-transfer} applies with zero residual and
$b_n=(b-a)u_n\downarrow0$, with sequential completeness of the target
$F$ supplied by Proposition~\ref{lt:d-complete}. It constructs a limit
$I$ of $S_\gamma f$
and gives $|I_n-I|\le(b-a)u_n$. The displayed equivalence identifies $I$
with the Kurzweil integral of $f$ on the same interval.
\end{proof}

Returning to the original integrator, regularity gives constant and
indicator normalization through
Proposition~\ref{uc:normalization-certificate}, while its
positive-sequence supremum condition gives the displayed sequential
normality by Proposition~\ref{nr:coefficient-order}. Every premise
of the transfer proof is thus supplied by the original algebraic and
approximation data; neither this corollary nor
Theorem~\ref{uc:uniform} requires Henstock estimates or finite lattice
closure.

\subsection{Exceptional finite responses and the disjoint corollary}

Uniform smallness may fail on an exceptional union even when the
finite responses supported there become negligible. To retain this
possibility, we separate each record into its ordinary and exceptional
parts and impose the error condition directly on the latter responses.

Let $\mathcal H$ be a $\sigma$-algebra and $L_C:X\to Z$ additive in
both the argument and disjoint measurable cells. For a countable
disjoint family $\mathcal P_\lambda\subset\mathcal H$, put
$E_\lambda=\bigcup\mathcal P_\lambda$ and
\[
 B_\lambda(g)=\{0\}\cup\left\{
 \sum_{r=1}^mL_{C_r}g(t_r):
 C_r\text{ disjoint},\ t_r\in C_r,
 \ \forall r\ \exists P\in\mathcal P_\lambda:\ C_r\subseteq P
 \right\}.
\]

\begin{definition}[Intrinsic exceptional witnesses]\label{uc:quasi-witness}\label{lt:quasi-definition}
Let $\Lambda$ be directed. Witnesses consist of $a_\lambda$,
$\mathcal P_\lambda$, and $0\in U_\lambda\subset X$, with
\[
 a\ge a_\lambda,\quad t\notin E_\lambda
 \Longrightarrow f(t)-f_a(t)\in U_\lambda.
\]
Set
\[
 B_\lambda^{\rm app}=\bigcup_{a\ge a_\lambda}B_\lambda(f_a),
 \qquad \widehat D_\lambda=\bigcup_{\eta\ge\lambda}D_\eta.
\]
Require the decreasing families $\widehat U_\lambda$,
$\widehat B_\lambda(f)$, and $\widehat B_\lambda^{\rm app}$ to be
null in their original spaces. Each exceptional response is finite.
Each approximating response uses one index $a$.
\end{definition}

\begin{theorem}[Quasi-uniform integral extension]
\label{uc:quasi}\label{lt:quasi}\label{mg:quasi-extension}
Let $Z$ be complete and satisfy either of the two target conditions
of Theorem~\ref{uc:uniform}. Suppose the coefficients have continuous
finite semivariance on disjoint records. Suppose the actual records
are disjoint and internally tagged, and eventually refine, for every
$A,\lambda$, the completed partition
\begin{equation}\label{uc:completed-partition}
 \{A\cap P:P\in\mathcal P_\lambda,\ A\cap P\ne\varnothing\}
 \cup\{A\setminus E_\lambda\},
\end{equation}
with an empty complement omitted. If all $f_a$ are integrable on all
$A$ and Definition~\ref{uc:quasi-witness} holds, then $f$ is integrable
on all $A$ and $I_Af_a\to I_Af$ uniformly in $A$.
\end{theorem}
\begin{proof}
On the localized record tail, partition its finite index set into
$G$, consisting of cells contained in $A\setminus E_\lambda$, and
$B$, consisting of cells contained in members of $\mathcal P_\lambda$.
For $a\ge a_\lambda$, the finite identity is
\begin{align*}
 S_\gamma^Af-S_\gamma^Af_a
 &=\sum_{r\in G}L_{C_r}(f(t_r)-f_a(t_r))
 +\sum_{r\in B}L_{C_r}f(t_r)-\sum_{r\in B}L_{C_r}f_a(t_r)\\
 &\in Q_\lambda:=\Sigma(\widehat U_\lambda)
 +\widehat B_\lambda(f)-\widehat B_\lambda^{\rm app}.
\end{align*}
The good partial response belongs to $\Sigma(\widehat U_\lambda)$
by Lemma~\ref{lt:reduced-words}, since omitted tags can be assigned
zero. The exceptional subwords are disjoint, internally tagged and
subordinate to $\mathcal P_\lambda$, so their values belong to
$B_\lambda(f)$ and $B_\lambda^{\rm app}$, respectively. The three
families defining $Q_\lambda$ decrease, contain zero and have null
saturations, by finite semivariance and
Definition~\ref{uc:quasi-witness}.
To check nullity of $Q_\lambda$, fix three levels
$\lambda_1,\lambda_2,\lambda_3$ and an upper bound $\nu$ of them.
Then
\[
 Q_\nu\subseteq
 \Sigma(\widehat U_{\lambda_1})
 +\widehat B_{\lambda_2}(f)
 -\widehat B_{\lambda_3}^{\rm app}.
\]
Thus its saturation is a quasi-subnet of the sum of the first two
null saturation nets minus the third, indexed independently. Continuity
of addition and inversion gives $\langle Q_\lambda\rangle\to0$.

Under target condition~(1), Theorem~\ref{uv:order-verification}
gives \textnormal{(R)} and \textnormal{(P)}. These properties are
already assumed under condition~(2). Set
\[
 \kappa=\max\{\aleph_0,|\Lambda|\},\qquad
 P_\lambda=\{a:a\ge a_\lambda\}.
\]
For each $A,\lambda$, the localization hypothesis supplies
$\gamma(A,\lambda)$ such that the finite splitting above is valid
for every $\gamma\ge\gamma(A,\lambda)$. Consequently
\[
 a\in P_\lambda,\quad\gamma\ge\gamma(A,\lambda)
 \quad\Longrightarrow\quad
 S_\gamma^Af-S_\gamma^Af_a\in Q_\lambda.
\]
In Theorem~\ref{ci:interchange}, use the additive operation of $Z$,
the error family $(Q_\lambda)$, the candidate row
$x_\gamma^A=S_\gamma^Af$, the known rows
$x_{a,\gamma}^A=S_\gamma^Af_a$, and their limits
$y_a^A=I_Af_a$. The displayed inclusion is its comparison condition,
with outer tail $P_\lambda$ independent of $A$. Its error-index
bound is $|\Lambda|\le\kappa$, its nullity condition was proved
above, and completeness, \textnormal{(R)} and \textnormal{(P)} give
$\kappa$-completeness, \textnormal{$(R_\kappa)$} and
\textnormal{$(P_\kappa)$}. Its conclusions are therefore
\[
 S_\gamma^Af\longrightarrow z^A,
 \qquad I_Af_a\longrightarrow z^A\quad\text{uniformly in }A.
\]
The first limit defines $I_Af=z^A$ on the prescribed net, proving
both assertions.
\end{proof}

\begin{proposition}[Positive exceptional coefficients]
\label{uc:quasi-budget}\label{lt:exceptional-budget}
Suppose faithful input gauges $p_j:X\to(E_j)_+$ and output gauges
$\rho_j:Z\to(F_j)_+$ satisfy
\[
 \rho_j(L_rx)\le T_{j,r}(p_j(x)),
\]
where $T_{j,r}:E_j\to F_j$ are positive additive controller maps. On the localized
original records, let
\[
 \sum_rT_{j,r}\le T_j,\qquad
 \sum_{r\in B}T_{j,r}\le E_{j,\lambda},
\]
where $B$ denotes the exceptional indices. Suppose
\[
 \begin{aligned}
 p_j(f(t)-f_a(t))&\le u_{j,\lambda}
 &&(t\notin E_\lambda, a\ge a_\lambda),\\
 p_j(f(t)),p_j(f_a(t))&\le v_j
 &&(t\in E_\lambda, a\ge a_\lambda).
 \end{aligned}
\]
Then every such record satisfies
\begin{equation}\label{uc:exceptional-estimate}
 \rho_j(S_Rf-S_Rf_a)
 \le T_ju_{j,\lambda}+2E_{j,\lambda}v_j.
\end{equation}
If the two terms on the right are order null, uniformly over their
respective tails in $\lambda$, then they verify the comparison in
Theorem~\ref{lt:order-transfer}. If, in addition, $u_{j,\lambda}\to_o0$ and the exceptional
coefficient bound holds for all finite subordinate words, these
estimates also verify the witness conditions of Definition~\ref{uc:quasi-witness}.
\end{proposition}
\begin{proof}
The finite splitting into good and exceptional indices gives
\begin{align*}
 \rho_j(S_Rf-S_Rf_a)
 &\le\sum_{r\in G}T_{j,r}u_{j,\lambda}
 +\sum_{r\in B}T_{j,r}p_j(f(t_r))
 +\sum_{r\in B}T_{j,r}p_j(f_a(t_r))\\
 &\le T_ju_{j,\lambda}+2E_{j,\lambda}v_j.
\end{align*}
For one exceptional finite word and one integrand
$g\in\{f\}\cup\{f_a:a\ge a_\lambda\}$,
\[
 \rho_j\left(\sum_rL_rg(t_r)\right)
 \le\sum_rT_{j,r}p_j(g(t_r))
 \le E_{j,\lambda}v_j.
\]
Order controls for the right-hand families bound their saturated
tail unions. The minimal good-error sets
\[
 U_\lambda=\{0\}\cup
 \{f(t)-f_a(t):a\ge a_\lambda, t\notin E_\lambda\}
\]
have the input controls $u_{j,\lambda}$. Faithfulness proves the
three required original convergences.
\end{proof}

For positive scalar means, the specialization
$T_{j,r}u=w_ru$, $T_ju=Mu$, and
$E_{j,\lambda}u=\varepsilon_\lambda u$ yields
\[
 \rho_j(S_Rf-S_Rf_a)
 \le Mu_{j,\lambda}+2\varepsilon_\lambda v_j.
\]
In an Archimedean controller vector lattice,
$\varepsilon_\lambda\to0$ implies
$\varepsilon_\lambda v_j\to_o0$, since
$\inf_{n\ge1}n^{-1}v_j=0$.
\begin{proposition}[Positive measure verification of exceptional witnesses]
\label{uc:measure-exceptional-certificate}
Let $(p_i:X\to(E_i)_+)_i$ and $(\rho_j:Z\to(F_j)_+)_j$ be faithful
order representations of the original input and output convergences.
Let $Y$ be an abelian lattice-ordered group, and let
$\mu:\mathcal H\to Y_+$ be finitely additive on a $\sigma$-algebra,
with $\Omega\in\mathcal H$. Suppose
\[
 \rho_j(b(x,y))\le B_j((p_i(x))_{i\in I_j},y),\qquad y\ge0,
\]
where $I_j$ is finite, $B_j$ is isotone in its first variable and
positive additive in its second. Assume the partial maps
\[
 u\longmapsto B_j(u,\mu\Omega),\qquad
 y\longmapsto B_j(v_j,y)
\]
preserve decreasing zero infima for each specified tuple $v_j$.
Suppose the witness data satisfy
\[
 \begin{gathered}
 u_{j,\lambda}\longrightarrow_o0,\qquad
 e_\lambda\longrightarrow_o0,\qquad \mu E_\lambda\le e_\lambda,\\
 (p_i(f(t)-f_a(t)))_{i\in I_j}\le u_{j,\lambda}
       \quad(t\notin E_\lambda,\ a\ge a_\lambda),\\
 (p_i(g(t)))_{i\in I_j}\le v_j
       \quad(t\in E_\lambda,\
             g\in\{f\}\cup\{f_a:a\ge a_\lambda\}).
 \end{gathered}
\]
Then the coefficient family has continuous finite semivariance. If every
input coordinate occurs in at least one $I_j$, the minimal good-error sets
together with the prescribed exceptional families also satisfy
Definition~\ref{uc:quasi-witness}. On every localized record,
\[
 \rho_j(S_Rf-S_Rf_a)
 \le B_j(u_{j,\lambda},\mu\Omega)+2B_j(v_j,e_\lambda)
 \qquad(a\ge a_\lambda).
\]
\end{proposition}
\begin{proof}
For finite disjoint coefficient cells, put
$w=(\bigvee_rp_i(x_r))_{i\in I_j}$. Finite additivity gives
\[
 \begin{aligned}
 \rho_j\!\left(\sum_rb(x_r,\mu C_r)\right)
 &\le\sum_rB_j(w,\mu C_r)\\
 &=B_j\!\left(w,\sum_r\mu C_r\right)
 \le B_j(w,\mu\Omega).
 \end{aligned}
\]
Proposition~\ref{uc:finite-certificates} proves continuous finite
semivariance.
For a finite exceptional word and any indicated $g$,
\[
 \sum_r\mu C_r\le\mu E_\lambda\le e_\lambda,
 \qquad
 \rho_j\!\left(\sum_rb(g(t_r),\mu C_r)\right)
 \le B_j(v_j,e_\lambda).
\]
Let $v^j_\theta\downarrow0$ and $w_\zeta\downarrow0$ control
$u_{j,\lambda}$ and $e_\lambda$. Their attained tails give
\[
 \begin{gathered}
 \lambda\ge\lambda(j,\theta,\zeta)
 \Longrightarrow
 u_{j,\lambda}\le v^j_\theta,\quad e_\lambda\le w_\zeta,\\
 B_j(u_{j,\lambda},\mu\Omega)\le B_j(v^j_\theta,\mu\Omega),
 \qquad B_j(v_j,e_\lambda)\le B_j(v_j,w_\zeta),\\
 B_j(v^j_\theta,\mu\Omega)\downarrow0,\qquad
 B_j(v_j,w_\zeta)\downarrow0.
 \end{gathered}
\]
For $\lambda\ge\lambda(j,\theta,\zeta)$ and every
$\eta\ge\lambda$, the same controlling thresholds apply to
$u_{j,\eta}$ and $e_\eta$. Therefore
\[
 \begin{gathered}
 x\in\widehat U_\lambda\quad\Longrightarrow\quad
 (p_i(x))_{i\in I_j}\le v^j_\theta,\\
 z\in\widehat B_\lambda(f)\cup\widehat B_\lambda^{\rm app}
 \quad\Longrightarrow\quad
 \rho_j(z)\le B_j(v_j,w_\zeta).
 \end{gathered}
\]
For an input coordinate $i$, take a $j$ with $i\in I_j$. The
$i$th coordinate of $v^j_\theta\downarrow0$ supplies its attained
zero regulator. Thus every input coordinate is null on the good-error
saturation, and every output coordinate is null on both exceptional
saturations. Faithfulness gives the three original null convergences. Finally the good/bad splitting gives
\[
 \begin{aligned}
 \rho_j(S_Rf-S_Rf_a)
 &\le\sum_{r\in G}B_j(u_{j,\lambda},\mu C_r)
       +2\sum_{r\in B}B_j(v_j,\mu C_r)\\
 &\le B_j(u_{j,\lambda},\mu\Omega)+2B_j(v_j,e_\lambda).
 \end{aligned}
\]
\end{proof}

For decreasing measurable $E_n$ with $\bigcap_nE_n=\varnothing$, positive
order-countable additivity independently verifies the coefficient decay.
Indeed, with $D_n=E_n\setminus E_{n+1}$,
\[
 E_1=\bigsqcup_{n\ge1}D_n,\qquad
 \sum_{n=1}^{N}\mu D_n=\mu E_1-\mu E_{N+1}\uparrow\mu E_1,\qquad
 \mu E_n\downarrow0.
\]
Combined with the pointwise bounds on the integrands and the displayed
positive coefficient maps, order-countable additivity of the given
measure therefore supplies this instance directly, without assuming
convergence of integrals or of infinite sums of integration errors.

\begin{corollary}[The disjoint quasi-uniform specialization]
\label{uc:millington}\label{lt:millington-manuscript}
\label{quasiuniformconvergencetheorem}
Let $X,Y,Z$ be separated topological vector spaces with $Z$ complete,
let $b:X\times Y\to Z$ be bilinear, and let
$\mu:\mathcal H\to Y$ be finitely additive on a $\sigma$-algebra.
Put $L_Cx=b(x,\mu C)$. Suppose finite semivariation holds:
\begin{equation}\label{finitesemivariationmillington}
 W\in\mathcal N_Z\quad\Longrightarrow\quad
 \exists U\in\mathcal N_X:\quad\Sigma(U)\subseteq W,
\end{equation}
where $\Sigma$ ranges over all finite disjoint internally tagged
coefficient words and $\mathcal N_X,\mathcal N_Z$ denote the
families of neighborhoods of zero in the respective original topologies.
For every $U\in\mathcal N_X,W\in\mathcal N_Z$, suppose there are
$a_0$ and a countable disjoint $\mathcal P\subseteq\mathcal H$ such that
\[
 \begin{gathered}
 a\ge a_0,\quad t\notin\bigcup\mathcal P
 \Longrightarrow f(t)-f_a(t)\in U,\\
 B_{\mathcal P}(f)\subseteq W,\qquad
 \bigcup_{a\ge a_0}B_{\mathcal P}(f_a)\subseteq W.
 \end{gathered}
\]
Assume the original records are disjoint, internally tagged, and
ultimately refine every completed witness partition
\eqref{uc:completed-partition}. If every $f_a$ is integrable on every
$A\in\mathcal H$, then $f$ is integrable on every $A$, and
$I_Af_a\to I_Af$ uniformly in $A$.
\end{corollary}
\begin{proof}
If a decreasing family $B_\eta\ni0$ is originally null, then for
each $U\in\mathcal N_X$ some $\eta_U$ satisfies
$B_\eta\subseteq U$ for $\eta\ge\eta_U$.
Equation~\eqref{finitesemivariationmillington} gives
\[
 W\in\mathcal N_Z\Longrightarrow
 \Sigma(B_\eta)\subseteq\Sigma(U)\subseteq W
 \quad(\eta\ge\eta_U).
\]
Thus the coefficient system has continuous finite semivariance.
For each nonempty $C\in\mathcal H$, a singleton coefficient word
with cell $C$ gives $L_C(U)\subseteq\Sigma(U)$, so
\eqref{finitesemivariationmillington} proves continuity of $L_C$ at
zero. Its additivity gives continuity everywhere, and
$L_\varnothing=0$. To realize these fixed coefficient maps within
Definition~\ref{bilinearproduct} without assuming continuity of
$y\mapsto b(x,y)$, give
\[
 Y_0=\bigoplus_{C\in\mathcal H\setminus\{\varnothing\}}\mathbb Z e_C
\]
the discrete convergence and put
\[
 b_0\!\left(x,\sum_Cn_Ce_C\right)=\sum_Cn_CL_Cx,
 \qquad
 \mu_0(C)=
 \begin{cases}e_C,&C\ne\varnothing,\\0,&C=\varnothing.\end{cases}
\]
The sums are finite. For fixed $y\in Y_0$, the map
$x\mapsto b_0(x,y)$ is a finite sum of continuous homomorphisms.
For fixed $x$, the map on the discrete group $Y_0$ is continuous.
Both partial maps preserve addition and zero, and
$b_0(x,\mu_0(C))=L_Cx$. Thus this coefficient product gives exactly
the original finite sums, response sets and witness conditions.

Let $\Lambda$ consist of all witnesses
$(U,W,a_0,\mathcal P)$, ordered by shrinking $U,W$ and increasing
$a_0$. Given two witnesses, the hypotheses supply a witness for the
intersections of their neighborhoods. Increasing its outer index to
a common upper bound preserves both witness inequalities.
Consequently $\Lambda$ is directed. Its tail unions satisfy
\[
 \widehat U_\lambda\subseteq U_\lambda,\qquad
 \widehat B_\lambda(f)\subseteq W_\lambda,\qquad
 \widehat B_\lambda^{\rm app}\subseteq W_\lambda.
\]
For any prescribed zero neighborhoods $U_*$ and $W_*$, the witness
hypothesis gives $\lambda_*$ with $U_{\lambda_*}=U_*$ and
$W_{\lambda_*}=W_*$. On every later witness, the displayed
inclusions place the three tail unions in $U_*$, $W_*$ and $W_*$,
respectively. Hence their saturations converge to zero in the
original topologies. The prescribed completed-partition condition
is the localization in Theorem~\ref{uc:quasi}.
The original additive convergence of $Z$ satisfies
\textnormal{(R)} and \textnormal{(P)} by
Corollary~\ref{uv:standard-cases}(1), and $Z$ is complete by
assumption. Together with the continuous finite semivariance and
witness nullity just established, these are target condition~(2)
and the coefficient and witness assumptions of Theorem~\ref{uc:quasi}.
The known rows are $S_\gamma^Af_a\to I_Af_a$ on the prescribed
$\Gamma_A$, and the stated completed-partition hypothesis is its
localization requirement. Applying that theorem gives integrability
of $f$ and convergence $I_Af_a\to I_Af$ uniformly in $A$.
\end{proof}

This corollary returns, in a single-valued prescribed-net form, to
the quasi-uniform definition and integral convergence theorem of
\cite[Definition~1.6.1 and Theorem~1.8.1]{millingtonproduct} for
set-valued partitionable functions, with the finite-response witnesses
now verified explicitly. Its conclusion includes integrability of the
candidate for the prescribed net Riemann integral and is limited to
that setting, rather than the full set-valued assertion. In making
the comparison, the complement of the exceptional union completes
the partition, and the refinement condition is imposed only after
this complement has been included.

\begin{corollary}[The uniform specialization on a ring]
\label{uc:uniform-classical}\label{lt:uniform-manuscript}
\label{uniformconvergencetheoremtopological}
Retain the spaces, finite semivariation, disjoint internally tagged
records, and known integrable approximants of
Corollary~\ref{uc:millington}, with $\mathcal H$ only a ring.
If $f_a\to f$ uniformly in the original input convergence, then
$f$ is integrable on every $A\in\mathcal H$, and
$I_Af_a\to I_Af$ uniformly in $A$.
\end{corollary}
\begin{proof}
For every $U\in\mathcal N_X$, uniform convergence supplies $a_U$
with $f_a(t)-f(t)\in U$ for every $a\ge a_U$ and $t\in\Omega$.
Since $0\in U$, this gives
\[
 \{0\}\cup\{f_a(t)-f(t):a\ge a_U,\ t\in\Omega\}\subseteq U.
\]
These inclusions for all zero neighborhoods say that the saturated
input-error family in Theorem~\ref{uc:uniform} converges to zero.
To check its coefficient hypothesis, let $(B_\eta)$ be any decreasing
null input family containing zero. Given $W\in\mathcal N_Z$,
finite semivariation supplies $U\in\mathcal N_X$ with
$\Sigma(U)\subseteq W$. Nullity supplies $\eta_U$ such that
$B_\eta\subseteq U$ whenever $\eta\ge\eta_U$, and hence
\[
 \Sigma(B_\eta)\subseteq\Sigma(U)\subseteq W
                         \qquad(\eta\ge\eta_U).
\]
The response saturation is therefore null. The fixed coefficient
maps are continuous because
$L_C(U)\subseteq\Sigma(U)\subseteq W$ for each nonempty $C$,
while $L_\varnothing=0$. The discrete coefficient group
$Y_0=\bigoplus_{C\in\mathcal H\setminus\{\varnothing\}}\mathbb Z e_C$, with
$b_0(x,\sum_Cn_Ce_C)=\sum_Cn_CL_Cx$,
$\mu_0(C)=e_C$ for $C\ne\varnothing$, and $\mu_0(\varnothing)=0$,
therefore has continuous additive partial maps and reproduces every
original record sum. The target has
\textnormal{(R)} and \textnormal{(P)} by
Corollary~\ref{uv:standard-cases}(1), and is complete by assumption.
Every row is integrable on its prescribed net. These verifications
supply all the hypotheses of Theorem~\ref{uc:uniform}, giving
integrability, uniform endpoint convergence and the two iterated
limits. No partition beyond the finite disjoint records enters this
argument, so a ring suffices.
\end{proof}

When the approximation is uniform, the countable witness partitions
are no longer needed. Both classical specializations concern the
convergence already given on the spaces, and the abstract theorems
continue to apply in their stated convergence structures without
assigning a topology to a non-topological space.

\begin{proposition}[Response-preserving disjointization]
\label{uc:disjointization}
Suppose finite tag sets are measurable. Let
$R=((C_r,t_r))_{r=1}^m$, $t_r\in C_r$, and assume
\[
 B\in\mathcal H,\quad B\subseteq C_r\cap C_s,\quad r\ne s
 \Longrightarrow L_B=0.
\]
Then there is a finite disjoint internally tagged coefficient word
$R'$, not necessarily an admissible integration record, with
$S_{R'}g=S_Rg$ for every $g$. If $R$ is localized at a completed
witness partition, $R'$ is localized at that same partition.
\end{proposition}
\begin{proof}
Let $t_1,\ldots,t_k$ be the distinct tags, $T=\{t_1,\ldots,t_k\}$,
and $C_i^*=\bigcup_{r:t_r=t_i}C_r$. Finite additivity and hereditary
nullity of overlaps give
\[
 L_{C_i^*}=\sum_{r:t_r=t_i}L_{C_r}.
\]
Finite unions of hereditarily coefficient-null sets are hereditarily
coefficient-null: a measurable subset of their union has a finite
disjoint decomposition into subsets of the individual sets.
Hence every $C_i^*\cap C_j^*$, $i\ne j$, is hereditarily null. Define
\[
 D_i=\left((C_i^*\setminus T)\setminus
                 \bigcup_{j<i}(C_j^*\setminus T)\right)\cup\{t_i\}.
\]
Then $t_i\in D_i\subseteq C_i^*$ and $D_i\cap D_j=\varnothing$ for
$i\ne j$. Every point of $C_i^*\setminus D_i$ either lies in an
earlier $C_j^*$, or is a different tag $t_j\in C_j^*$. Therefore
\[
 C_i^*\setminus D_i\subseteq\bigcup_{j\ne i}(C_i^*\cap C_j^*),
 \qquad L_{D_i}=L_{C_i^*}.
\]
For every $g$,
\[
 \sum_iL_{D_i}g(t_i)
 =\sum_i\sum_{r:t_r=t_i}L_{C_r}g(t_r)=S_Rg.
\]
If the original record is localized, every cell carrying $t_i$ lies
in the unique member of the completed partition containing $t_i$.
Its union $C_i^*$ and its subset $D_i$ lie in that member.
\end{proof}

Since the source finite-semivariance and exceptional-response
conditions quantify over all finite disjoint coefficient words, the
equality of sums transfers their estimates to the original records
without any cofinality requirement on the auxiliary words in the
integration net.

For positive measures, null overlaps are hereditary by monotonicity:
$B\subseteq N$ and $\mu N=0$ imply $0\le\mu B\le\mu N=0$.
For a general signed or group-valued set function, $\mu(N)=0$ alone does
not establish hereditary coefficient-nullity. The null-set convention in
\cite[Definition~0.1]{millingtonproduct} is hereditary.
Proposition~\ref{uc:disjointization} identifies a sufficient mechanism
for the overlapping-record version without altering its sums or tags.

\begin{corollary}[Regulator forms of uniform extension]\label{lt:quasi-d}
The uniform and quasi-uniform theorems have the following versions
for a faithful $D$-representation with Dedekind $\sigma$-complete
WSD controllers.
\begin{enumerate}
\item For a sequence of approximants, sequential completeness of $Z$
suffices when the finite comparison is eventually bounded by
$b_{j,n}\to_o0$, independently of the domain.
\item For an arbitrary approximation net, completeness of $Z$ suffices
when its finite comparisons have regulator cuts $P^j_\varphi$ admitting
the countable accessibility in Theorem~\ref{lt:d-anchors}.
\item For super Dedekind complete WSD controllers, the order and
regulator versions coincide.
\end{enumerate}
\end{corollary}
\begin{proof}
Use the original response nets
\[
 x^A_\gamma=S_\gamma^Af,\qquad
 x^A_{a,\gamma}=S_\gamma^Af_a,\qquad
 y_a^A=I_Af_a.
\]
Integrability of each approximant is exactly
$x^A_{a,\gamma}\to y_a^A$ on $\Gamma_A$.
In item~(1), $a=n\in\mathbb N$ and the asserted finite comparison
reads
\[
 d_j(x^A_\gamma,x^A_{n,\gamma})\le b_{j,n}+0
                      \quad\text{eventually in }\gamma.
\]
The bound is order null in $n$ and independent of $A$, and the
constant zero residual is $D$-null. Together with the specified
sequential completeness and faithful $D$-representation, these are
the assumptions of Theorem~\ref{lt:d-transfer}. Its conclusion gives
$x^A_\gamma\to z^A$ and $y_n^A\to z^A$ uniformly in $A$.

In item~(2), the finite comparisons instead have the form
\[
 a\in P^j_\varphi\quad\Longrightarrow\quad
 d_j(x^A_\gamma,x^A_{a,\gamma})
       \le H_{c^j}(\varphi)+0
                 \quad\text{eventually in }\gamma.
\]
The sets $P^j_\varphi$ contain outer tails independently of $A$,
and the assumed countable accessibility is the remaining index
condition of Theorem~\ref{lt:d-anchors}. Its residual requirement
holds on every accessible index because the residual is zero.
That theorem, with the given completeness of $Z$, gives the same
two limits on the prescribed indices.

In item~(3), Proposition~\ref{lt:order-d} gives, for every represented
pair net $(z_\alpha,w_\alpha)$,
\[
 d_j(z_\alpha,w_\alpha)\longrightarrow_D0
 \iff d_j(z_\alpha,w_\alpha)\longrightarrow_o0.
\]
Thus the original point and pair convergences, their Cauchy nets,
and their completeness assertions are unchanged when either
controller description is used. The order bounds in the uniform
and quasi-uniform comparisons therefore satisfy
Theorem~\ref{lt:order-transfer} in this original convergence.
In every item, $x^A_\gamma\to z^A$ defines $I_Af=z^A$, and the
endpoint limit is the asserted uniform convergence of the integrals.
\end{proof}

\subsection{Asymptotic finite responses with monoid values}

To pass from estimates on individual finite responses to conjugacy
estimates, we must follow the bounds through a tail infimum. In a
monoid this passage is supplied by preservation of directed infima
under translation, extending the zero-infimum law used above and
performing the step that, in a lattice group, can be written by
subtracting a limiting tail bound.

\begin{definition}[Translation of directed infima]\label{mg:directed-infima}
Let $M$ be a complete lattice and an ordered commutative monoid with
least element and identity $0$. Translations preserve directed
infima if
\begin{equation}\label{mg:directed-translation}
 r\oplus\inf D=\inf_{d\in D}(r\oplus d)
 \qquad(r\in M, \varnothing\ne D\subseteq M
                  \text{ downward directed}).
 \tag{DI}
\end{equation}
This identity implies \eqref{mg:zero-translation}. All suprema below
are taken in $M$, including its largest element when required.
\end{definition}

\begin{theorem}[Interchange from asymptotic response bounds]
\label{mg:asymptotic-interchange}
Let $Z$ be a complete separated uniform convergence space with a
faithful monoid pair representation $d_j:Z^2\to M_j$. Suppose each
$M_j$ satisfies Definition~\ref{mg:directed-infima}. Let
\[
 x^A_{p,\gamma}\longrightarrow y_p^A
 \qquad(p\in P,\ A\in\mathcal A\ne\varnothing),
\]
where $P$ and every $\Gamma_A$ are directed. For candidate values
$x^A_\gamma$, define
\[
 M^{A,j}_{p,\delta}
   =\sup_{\gamma\ge\delta}d_j(x^A_{p,\gamma},x^A_\gamma),
 \qquad
 R^{A,j}_p=\inf_{\delta\in\Gamma_A}M^{A,j}_{p,\delta}.
\]
Suppose $b_{j,p}\in M_j$ is independent of $A$ and satisfies
\[
 R^{A,j}_p\le b_{j,p},\qquad b_{j,p}\longrightarrow_o0.
\]
Then unique $z^A\in Z$ satisfy
\[
 x^A_\gamma\longrightarrow z^A,\qquad
 y_p^A\longrightarrow z^A\quad\text{uniformly in }A,
\]
and, for every $A,j,p$,
\begin{equation}\label{mg:asymptotic-sharp}
 d_j(y_p^A,z^A)=R^{A,j}_p\le b_{j,p}.
\end{equation}
\end{theorem}
\begin{proof}
Fix $A,j$ and omit these superscripts and subscripts until uniformity
is considered. For $p,q\in P$ and $\delta,\varepsilon\in\Gamma_A$,
all $\gamma$ beyond a common upper index of $\delta,\varepsilon$
satisfy
\begin{align*}
 d(y_p,y_q)
 &\le d(y_p,x_{p,\gamma})
       \oplus d(x_{p,\gamma},x_\gamma)
       \oplus d(x_\gamma,x_{q,\gamma})
       \oplus d(x_{q,\gamma},y_q)\\
 &\le M_{p,\delta}\oplus M_{q,\varepsilon}
       \oplus d(y_p,x_{p,\gamma})
       \oplus d(x_{q,\gamma},y_q).
\end{align*}
The last two terms have an order-null sum by row convergence,
point faithfulness, and Lemma~\ref{mg:calculus}(1).
Part (2) gives
$d(y_p,y_q)\le M_{p,\delta}\oplus M_{q,\varepsilon}$.
The nets $M_{p,\delta}$ and $M_{q,\varepsilon}$ decrease, so
\eqref{mg:directed-translation} gives
\[
 d(y_p,y_q)
 \le\inf_\varepsilon\inf_\delta
             (M_{p,\delta}\oplus M_{q,\varepsilon})
 =R_p\oplus R_q\le b_p\oplus b_q.
\]
Let $v_\theta\downarrow0$ regulate $b_p$, and let $p_\theta$
satisfy $p\ge p_\theta\Longrightarrow b_p\le v_\theta$.
Then
\[
 p,q\ge p_\theta\Longrightarrow
 d(y_p,y_q)\le v_\theta\oplus v_\theta\downarrow0.
\]
This argument holds for every $j$. Pair faithfulness makes
$(y_p^A)_p$ originally Cauchy, and completeness supplies its unique
limit $z^A$. For fixed $p$, the triangle inequality gives
\[
 d(y_p,z)\le b_p\oplus b_q\oplus d(y_q,z).
\]
The last two terms have an order-null sum. Hence
Lemma~\ref{mg:calculus}(2) gives $d(y_p,z)\le b_p$.
The regulators for $b_{j,p}$ and their outer thresholds are
independent of $A$. Therefore this inequality and pair faithfulness
already prove uniform endpoint convergence.

Put
\[
 e_\gamma=d(x_\gamma,z),\qquad
 D_p\downarrow0\text{ a regulator of }d(x_{p,\gamma},y_p),
 \qquad
 \mathscr R=\{r:e_\gamma\le r\text{ eventually}\}.
\]
For each $p,\delta$ and $u\in D_p$, a common inner tail satisfies
\[
 e_\gamma\le d(x_\gamma,x_{p,\gamma})
          \oplus d(x_{p,\gamma},y_p)\oplus d(y_p,z)
       \le M_{p,\delta}\oplus b_p\oplus u.
\]
Thus $\mathscr R$ contains every displayed bound and is nonempty.
If $c$ is a lower bound of $\mathscr R$, then
\begin{align*}
 c&\le\inf_{u\in D_p}(M_{p,\delta}\oplus b_p\oplus u)
       =M_{p,\delta}\oplus b_p\qquad(p,\delta),\\
 c&\le\inf_\delta(M_{p,\delta}\oplus b_p)
       =R_p\oplus b_p\le b_p\oplus b_p\qquad(p).
\end{align*}
For every $\theta$, there is $p\ge p_\theta$, and the last bound
is at most $v_\theta\oplus v_\theta$. Hence $c\le0$.
Lemma~\ref{mg:calculus}(3) gives $e_\gamma\to_o0$.
For every $j$ this is point convergence to zero of the candidate
error, so faithfulness gives $x^A_\gamma\to z^A$.

It remains to identify the asymptotic radius. Fix $p$ and write
\[
 a_\gamma=d(x_{p,\gamma},x_\gamma),\qquad
 \ell=d(y_p,z),\qquad
 h_\gamma=d(x_{p,\gamma},y_p)\oplus d(x_\gamma,z).
\]
The net $h_\gamma$ is order null. Let $U\downarrow0$ regulate it.
For each $u\in U$, there exists $\delta_u$ such that
$h_\gamma\le u$ for $\gamma\ge\delta_u$. The two triangle
inequalities give
\[
 a_\gamma\le\ell\oplus h_\gamma,\qquad
 \ell\le a_\gamma\oplus h_\gamma.
\]
Consequently
\[
 R_p\le M_{p,\delta_u}\le\ell\oplus u\quad(u\in U),
 \qquad R_p\le\inf_{u\in U}(\ell\oplus u)=\ell.
\]
For every $\delta\ge\delta_u$, the nonempty tail at $\delta$
also gives $\ell\le M_{p,\delta}\oplus u$.
The tail $\{\delta:\delta\ge\delta_u\}$ is cofinal in
$\Gamma_A$, and $M_{p,\delta}$ decreases. Therefore
\[
 \inf_{\delta\ge\delta_u}M_{p,\delta}=R_p,
 \qquad
 \ell\le\inf_{\delta\ge\delta_u}(M_{p,\delta}\oplus u)
       =R_p\oplus u.
\]
Taking the infimum over $u\in U$ gives $\ell\le R_p$.
This proves \eqref{mg:asymptotic-sharp}.
\end{proof}

\begin{corollary}[Conjugate extension with monoid values]
\label{mg:monoid-conjugacy}
Let $Z$ be a complete separated convergence abelian group with a
faithful invariant symmetric gauge family $\rho_j:Z\to M_j$,
where every $M_j$ satisfies Definition~\ref{mg:directed-infima}.
For the original finite sums, put
\[
 \mathcal R_j^A(h)
   =\inf_{\delta\in\Gamma_A}
       \sup_{\gamma\ge\delta}\rho_j(S_\gamma^Ah).
\]
Let $E_j$ be a partially ordered set with least element $0$, let
$\widehat E_j=E_j\cup\{\infty\}$, and let an input map $P_j$
take values in $\widehat E_j$. Suppose $T_j:E_j\to M_j$ is
isotone, $T_j(0)=0$, and preserves decreasing zero infima. Extend
it by $\widehat T_j(\infty)=\max M_j$ and
$\widehat T_j|_{E_j}=T_j$. Assume
\[
 \mathcal R_j^A(h)\le\widehat T_j(P_j(h))\qquad(A,j,h).
\]
If every $f_p$ is integrable on every prescribed $A$ and, for every
$j$, $P_j(f_p-f)$ is eventually finite and its finite tail is
order null, then $f$ is integrable on every such $A$,
$I_Af_p\to I_Af$ uniformly in $A$, and
\[
 \rho_j(I_Af_p-I_Af)
   =\mathcal R_j^A(f_p-f)
   \le T_j(P_j(f_p-f))
\]
at every index where the right-hand side is finite-valued.
\end{corollary}
\begin{proof}
Put $h_p=f_p-f$ and $b_{j,p}=\widehat T_j(P_j(h_p))$.
For each $j$, eventual finiteness and order nullity provide a
regulator $u_{j,\theta}\downarrow0$ and outer thresholds with
\[
 p\ge p(j,\theta)
 \Longrightarrow P_j(h_p)\in E_j,\quad
 P_j(h_p)\le u_{j,\theta},\quad
 b_{j,p}\le T_j(u_{j,\theta}).
\]
The right-hand family decreases to $0$, so $b_{j,p}\to_o0$.
Additivity of the finite sums and invariance of the representation
give
\[
 d_j(S_\gamma^Af_p,S_\gamma^Af)
   =\rho_j(S_\gamma^Af_p-S_\gamma^Af)
   =\rho_j(S_\gamma^Ah_p).
\]
Thus the asymptotic radii of
Theorem~\ref{mg:asymptotic-interchange} are
$\mathcal R_j^A(h_p)\le b_{j,p}$. Its known row limits are
$I_Af_p$. The theorem constructs $I_Af$, proves uniform endpoint
convergence, and identifies each radius with
$\rho_j(I_Af_p-I_Af)$. At an index with $P_j(h_p)\in E_j$,
$b_{j,p}=T_j(P_j(h_p))$, giving the asserted bound.
\end{proof}

Returning to the truncated monoid of Example~\ref{mg:truncated},
we can also verify \eqref{mg:directed-translation}. If $a=\inf D$,
isotonicity gives
$r\oplus a\le\inf_{d\in D}(r\oplus d)$. For $\varepsilon>0$
there is $d\in D$ with $d<a+\varepsilon$, unless $a\in D$, in
which case equality is attained directly. Hence
\[
 \inf_{d\in D}(r\oplus d)
 \le\min\{1,r+a+\varepsilon\},\qquad
 \inf_{d\in D}(r\oplus d)\le\min\{1,r+a\}=r\oplus a.
\]
Thus cancellation is unnecessary for this conjugacy theorem, as the
truncated controller already demonstrates. In the lattice formulation
to which we now turn, conditional completeness and the extended-finiteness
conventions are retained, while the group structure justifies
subtraction of a tail infimum directly, independently of the monoid
argument.

\subsection{Conjugate lattice controls}

Within a lattice group, a finite tail envelope allows the asymptotic
bound to be separated from a decreasing zero remainder. We retain
extended values until that envelope is available, so that finiteness
is required only where the comparison actually uses it.

Let the faithful representation be by symmetric $\rho_j$, and let
$F_j$ be Dedekind complete. In $\widehat F_j^+=F_j^+\cup\{\infty\}$,
put
\[
 R_j^A(h)=\inf_{\delta\in\Gamma_A}
              \sup_{\gamma\ge\delta}\rho_j(S_\gamma^Ah).
\]
Suppose the input functional $P_j$ and the fixed integrator satisfy
\begin{equation}\label{uc:conjugacy}
 R_j^A(h)\le T_j(P_j(h)),\qquad
 T_j(u)=u\star_jq_j(\mu),\qquad
 u_\theta\downarrow0\Longrightarrow T_j(u_\theta)\downarrow0,
\end{equation}
on the finite domain of $P_j$, where $T_j$ is isotone and finite valued.

\begin{proposition}[A finite certificate for conjugacy]
\label{uc:conjugacy-certificate}
Suppose every input controller $E_i$ is Dedekind complete and the
finite coefficient estimate \eqref{uc:finite-estimate} holds. Define, in the extended product
cone,
\[
 P_j(h)=\left(\sup_{t\in\Omega}p_i(h(t))\right)_{i\in I_j}.
\]
Whenever this tuple is finite,
\[
 R_j^A(h)\le\omega_j(P_j(h))\qquad(A\in\mathcal A).
\]
Consequently an order-continuous map
$T_j=\omega_j$ verifies the conjugate theorem. The positive measure
and variation budgets give
$T_j(u)=B_j(u,v)$.
\end{proposition}
\begin{proof}
For every finite record and every $i\in I_j$,
\[
 \bigvee_{t\in T_{A,\gamma}}p_i(h(t))\le(P_j(h))_i.
\]
Isotonicity and \eqref{uc:finite-estimate} give
\[
 \rho_j(S_\gamma^Ah)\le\omega_j(P_j(h))
 \quad(\gamma\in\Gamma_A).
\]
Every tail supremum has the same bound. Its infimum gives the stated
conjugacy inequality. The finite measure and variation estimates
already proved give the displayed formula for $T_j$.
\end{proof}

\begin{theorem}[Conjugate integral extension]\label{uc:conjugate}\label{lt:conjugate}
Let $Z$ be complete in its original uniform convergence. Suppose
every $f_a$ is integrable on every $A$, and
\[
 P_j(f_a-f)\text{ is eventually finite},\qquad
 P_j(f_a-f)\longrightarrow_o0\quad(j\in J).
\]
Then $f$ is integrable on every $A$, the primitives converge uniformly
in $A$, and
\begin{equation}\label{uc:conjugate-endpoint}
 \rho_j(I_Af_a-I_Af)=R_j^A(f_a-f)
 \le P_j(f_a-f)\star_jq_j(\mu)
\end{equation}
whenever the right side is finite.
\end{theorem}
\begin{proof}
Put $h_a=f_a-f$ and define
\[
 b_{j,a}=
 \begin{cases}
 T_j(P_j(h_a)),&P_j(h_a)\text{ is finite},\\
 0,&\text{otherwise}.
 \end{cases}
\]
Input order controls
$u_{j,\theta}\downarrow0$ give
\[
 a\ge a(j,\theta)\Longrightarrow
 b_{j,a}\le T_j(u_{j,\theta})=:\beta^j_\theta\downarrow0.
\]
Fix $A,j,a$ on the finite-domain tail and put
\[
 M_\delta=\sup_{\gamma\ge\delta}\rho_j(S_\gamma^Ah_a),
 \qquad R=\inf_\delta M_\delta\le b_{j,a}<\infty.
\]
Lemma~\ref{app:tail-envelope} gives $\delta_0$ with
$M_{\delta_0}<\infty$ and $M_\delta\downarrow R$ on
$\delta\ge\delta_0$. Define
\[
 r_{a,\gamma}^{A,j}=M_\gamma-R\quad(\gamma\ge\delta_0),
 \qquad r_{a,\gamma}^{A,j}=0\quad\text{otherwise}.
\]
Then $r_{a,\gamma}^{A,j}\to_o0$, and on that original tail,
\[
 d_j(S_\gamma^Af,S_\gamma^Af_a)
 =\rho_j(S_\gamma^Ah_a)
 \le R+r_{a,\gamma}^{A,j}
 \le b_{j,a}+r_{a,\gamma}^{A,j}.
\]
For each $j,\theta$, let
\[
 P^j_\theta=\{a:P_j(h_a)\text{ is finite and }b_{j,a}\le\beta^j_\theta\}.
\]
The input controls give a tail contained in $P^j_\theta$, independently
of $A$. Residuals at indices outside the finite domain are assigned zero.
Comparisons are used only for $a\in P^j_\theta$. Hence
\[
 x^A_{a,\gamma}=S_\gamma^Af_a\longrightarrow y_a^A=I_Af_a,
 \qquad
 d_j(S_\gamma^Af,S_\gamma^Af_a)
 \le\beta^j_\theta+r^{A,j}_{a,\gamma}\quad(a\in P^j_\theta)
\]
verify Theorem~\ref{lt:order-transfer}. It constructs $I_Af$ and proves
uniform endpoint convergence. The controlling inner tails may depend on
$A,j,a$. No common inner tail or common finite-domain tail over all $j$
has been used.

For fixed $a$, subtraction gives
$S_\gamma^Ah_a\to I_Af_a-I_Af$. Set
$\ell=\rho_j(I_Af_a-I_Af)$. Then
\[
 |\rho_j(S_\gamma^Ah_a)-\ell|
 \le\rho_j(S_\gamma^Ah_a-(I_Af_a-I_Af))\longrightarrow_o0.
\]
If $v_\eta\downarrow0$ is a control, a cofinal tail satisfies
$\ell-v_\eta\le\rho_j(S_\gamma^Ah_a)\le\ell+v_\eta$.
For every $\eta$, there is $\delta_\eta$ such that
\[
 \delta\ge\delta_\eta\Longrightarrow
 \ell-v_\eta\le
 \sup_{\gamma\ge\delta}\rho_j(S_\gamma^Ah_a)
 \le\ell+v_\eta.
\]
This tail of $\Gamma_A$ is cofinal. Infimizing over it gives
\[
 \ell-v_\eta\le R_j^A(h_a)\le\ell+v_\eta\quad(\eta),
 \qquad |R_j^A(h_a)-\ell|\le v_\eta\downarrow0.
\]
Thus $R_j^A(h_a)=\ell$. Equation
\eqref{uc:conjugacy} proves \eqref{uc:conjugate-endpoint}.
\end{proof}

\begin{definition}[Extended null convergence]\label{extendedconvergence}
For a positive cone $E_+$, write $\widehat E_+=E_+\cup\{\infty\}$.
A net in $\widehat E_+$ is extended null if it is eventually finite
and its eventual finite restriction is order null. Equivalently, it is a quasi-subnet of a finite-valued order-null net.
Indeed, if $x_i$ is a quasi-subnet of $y_k\in E_+$ and
$y_k\to_o0$, then one tail of $y$ lies in $E_+$ and each of its
controlled tails lies below the corresponding zero-regulator value.
The quasi-subnet condition places a tail of $x$ in each of these
sets, proving eventual finiteness and order nullity. Conversely, if
$x_i\in E_+$ for $i\ge i_0$, define $y$ to be the restriction of
$x$ to $I_{i_0}=\{i:i\ge i_0\}$. For every $j\in I_{i_0}$,
$T_j(y)=T_j(x)$, while every original tail contains $T_k(y)$ for
an upper bound $k$ of its initial index and $i_0$. Thus $x$ and $y$
are tail equivalent and $y$ is finite valued and order null.
\end{definition}

\begin{proposition}[Solid target gauges verify symmetry]
\label{uc:solid-symmetry}
Let $Z$ be an abelian lattice-ordered group. If a target gauge satisfies
\[
 |z|\le|w|\Longrightarrow\rho_j(z)\le\rho_j(w),
\]
then it is symmetric. Consequently a faithful solid subadditive target
family meets the symmetry hypothesis of Theorem~\ref{uc:conjugate}
without an extra condition on the input functional $P_j$.
\end{proposition}
\begin{proof}
For every $z\in Z$,
\[
 |z|=|-z|,\qquad
 \rho_j(z)\le\rho_j(-z)\le\rho_j(z),\qquad
 \rho_j(-z)=\rho_j(z).
\]
The other target hypotheses of Theorem~\ref{uc:conjugate} remain the
stated subadditivity, pair faithfulness and completeness assumptions.
\end{proof}

\begin{corollary}[The conjugated-seminorm theorem]
\label{uc:szabo-szaz}\label{lt:conjugate-manuscript}\label{latticenormconvergencetheorem}
Let $Z$ be a complete proper symmetric lattice-normed group with
Dedekind complete controller. Suppose $P,q$ are conjugated in the
sense of \eqref{uc:conjugacy}, $q(\mu)$ is finite, and the positive
partial product $u\mapsto u\star q(\mu)$ is order continuous.
If each $f_a$ is integrable on every prescribed $A$ and
$P(f_a-f)\to0$ in the extended convergence of Definition~\ref{extendedconvergence}, then $f$ is
integrable on those domains and
\[
 \rho(I_Af_a-I_Af)\le P(f_a-f)\star q(\mu)
 \qquad\text{whenever }P(f_a-f)\text{ is finite}.
\]
The endpoint convergence is uniform over every family of domains
with that same conjugate bound. Finiteness of $P(f)$ and an additional
continuity assumption on $\rho$ are unnecessary.
\end{corollary}
\begin{proof}
Extended convergence is defined through a quasi-subnet of a
finite-valued order-null net. Every eventual order bound on that net
therefore bounds an eventual tail of $P(f_a-f)$. Thus the latter is
eventually finite and order null. The continuity of the partial
product $u\mapsto u\star q(\mu)$ gives
$u_\theta\downarrow0\Longrightarrow u_\theta\star q(\mu)\downarrow0$.
Conjugacy gives \eqref{uc:conjugacy} on each original directed set
$\Gamma_A$.
Theorem~\ref{uc:conjugate} applies. Finally,
\[
 |\rho(z)-\rho(w)|\le\rho(z-w)
\]
proves continuity of $\rho$ in its represented convergence.
\end{proof}

Through the preceding corollary, the method of conjugate seminorms
used for Banach-valued net-integral convergence in
\cite{szaznetintegralconvergence} acquires a symmetric lattice-valued
form under Theorem~\ref{uc:conjugate}. To reach that form from a
merely subadditive target gauge, continuous inversion permits the
faithful symmetrization
$\rho^s(z)=\rho(z)\vee\rho(-z)$. Conjugacy for that symmetrization is
verified whenever both oriented input bounds are available:
\[
 R_{\rho^s}^A(h)
 \le T(P(h))\vee T(P(-h)).
\]
Indeed, on a common finite-envelope tail, put
\[
 M_\delta^\pm=\sup_{\gamma\ge\delta}\rho(\pm S_\gamma^Ah),
 \qquad R^\pm=\inf_\delta M_\delta^\pm.
\]
Then
\[
 M_\delta^s=M_\delta^+\vee M_\delta^-,\qquad
 0\le M_\delta^s-(R^+\vee R^-)
 \le(M_\delta^+-R^+)+(M_\delta^--R^-)\downarrow0.
\]
Thus $R_{\rho^s}^A(h)=R^+\vee R^-$ whenever both conjugate
bounds are finite. If either bound is infinite the asserted inequality
holds in the extended positive cone.
Symmetry of $P$, or order nullity of both $P(f_a-f)$ and $P(f-f_a)$,
is sufficient. For asymmetric target gauges, the following propositions
give two independently verified mechanisms: uniform finite estimates
and positive domination of inversion. Sublinearity of $P$ alone does
not supply either oriented null bound.

\begin{proposition}[Finite estimates and continuous inversion]
\label{uc:asymmetric-finite}
Let $Z$ be a complete separated convergence group with its original
additive uniform convergence. Suppose subadditive gauges
$\rho_j:Z\to(F_j)_+$, $\rho_j(0)=0$, faithfully represent its point
and pair convergences by order convergence. Symmetry of $\rho_j$ is
not required. Let $f_a$ be integrable on every prescribed $A$, and
suppose positive elements $b_{j,a}$ and outer indices $a_j$ satisfy
\[
 b_{j,a}\longrightarrow_o0,\qquad
 a\ge a_j\Longrightarrow
 \rho_j\bigl(S_\gamma^A(f_a-f)\bigr)\le b_{j,a}
 \quad(A\in\mathcal A,\ \gamma\in\Gamma_A).
\]
The thresholds $a_j$ may depend on the representing coordinate.
Then $f$ is integrable on every $A$, the convergence
$I_Af_a\to I_Af$ is uniform in $A$, and
\[
 \rho_j(I_Af_a-I_Af)\le b_{j,a}\qquad(a\ge a_j).
\]
The same sharp bound holds at any additional index where its finite-record
comparison is valid.
\end{proposition}
\begin{proof}
Put
\[
 e_{a,\gamma}^A=S_\gamma^A(f_a-f),\qquad
 U_{a_0}=\{0\}\cup
 \{e_{a,\gamma}^A:a\ge a_0,\ A\in\mathcal A,
                                      \ \gamma\in\Gamma_A\}.
\]
For every $j$, order convergence supplies a decreasing control
$u_{j,\theta}\downarrow0$ and indices $a(j,\theta)\ge a_j$ such that
\[
 a\ge a(j,\theta)\Longrightarrow b_{j,a}\le u_{j,\theta}.
\]
Consequently,
\[
 a_0\ge a(j,\theta),\ z\in U_{a_0}
 \Longrightarrow\rho_j(z)\le u_{j,\theta}.
\]
The saturated value net of $U_{a_0}$ is therefore null in every
representing coordinate. Faithfulness gives its original null
convergence. Continuity of inversion gives original null convergence
of the saturation of $-U_{a_0}$. Applying faithfulness to this net
gives, for every $j$, a decreasing control $v_{j,\eta}\downarrow0$
and indices $c(j,\eta)$ satisfying
\[
 a_0\ge c(j,\eta),\ z\in U_{a_0}
 \Longrightarrow\rho_j(-z)\le v_{j,\eta}.
\]
The controlling tail of the saturation has value set $-U_c$ for
some $c$. For $a_0\ge c$, one has $-U_{a_0}\subseteq-U_c$,
so the resulting bound holds for every $z\in U_{a_0}$.

Define
\[
 \rho_j^s(z)=\rho_j(z)\vee\rho_j(-z).
\]
Then $\rho_j^s(-z)=\rho_j^s(z)$ and
\[
 \begin{aligned}
 \rho_j^s(z+w)
 &\le[\rho_j(z)+\rho_j(w)]
       \vee[\rho_j(-z)+\rho_j(-w)]\\
 &\le\rho_j^s(z)+\rho_j^s(w).
 \end{aligned}
\]
For any net $z_\alpha$, original continuity of inversion and
faithfulness imply
\[
 z_\alpha\longrightarrow0
 \Longrightarrow
 \rho_j(z_\alpha),\rho_j(-z_\alpha)\longrightarrow_o0
 \Longrightarrow\rho_j^s(z_\alpha)\longrightarrow_o0.
\]
Conversely, $\rho_j\le\rho_j^s$ and solidity of order nullity
give $\rho_j(z_\alpha)\to_o0$ from
$\rho_j^s(z_\alpha)\to_o0$, and the original faithfulness gives
$z_\alpha\to0$. For an arbitrary pair net this yields
\[
 (z_\alpha,w_\alpha)\longrightarrow\Delta_Z
 \iff z_\alpha-w_\alpha\longrightarrow0
 \iff\rho_j^s(z_\alpha-w_\alpha)\longrightarrow_o0\quad(j\in J),
\]
which proves pair faithfulness of the symmetrized family.

On the product control index $(\theta,\eta)$, put
\[
 \beta_{j,\theta,\eta}=u_{j,\theta}+v_{j,\eta}\downarrow0,
 \qquad
 P_{j,\theta,\eta}
 =\{a:a\ge a(j,\theta),\ a\ge c(j,\eta)\}.
\]
Directedness makes $P_{j,\theta,\eta}$ a set containing an outer
tail. For all its indices,
\[
 \rho_j^s(S_\gamma^Af-S_\gamma^Af_a)
 =\rho_j^s(e_{a,\gamma}^A)
 \le\beta_{j,\theta,\eta}
 \quad(A,\gamma).
\]
Theorem~\ref{lt:order-transfer} applies with
\[
 x_{a,\gamma}^A=S_\gamma^Af_a,\quad
 y_a^A=I_Af_a,\quad x_\gamma^A=S_\gamma^Af,\quad
 d_j(z,w)=\rho_j^s(z-w),\quad r_{a,\gamma}^{A,j}=0.
\]
It constructs $I_Af$ and proves original uniform convergence of the
primitives. For fixed $j,A,a$ at which the finite comparison is valid,
set $e=I_Af_a-I_Af$. Subtraction gives
$e_{a,\gamma}^A\to e$, and
\[
 \rho_j(e)
 \le\rho_j(e-e_{a,\gamma}^A)+\rho_j(e_{a,\gamma}^A)
 \le\rho_j(e-e_{a,\gamma}^A)+b_{j,a}.
\]
The first term is order null by faithfulness. Lemma~\ref{lt:fixed-bound}
gives the asserted endpoint estimate.
\end{proof}

For the preceding proposition, a direct certificate is the finite
coefficient inequality
\[
 \rho_j(S_\gamma^Ah)\le T_j(P_j(h))\qquad(A,\gamma)
\]
whenever $P_j(f_a-f)$ is eventually finite
and order null, and the isotone map $T_j$ preserves decreasing null
controls. For each $j$, its finite-domain tail has some initial index
$a_j$. On that tail put $b_{j,a}=T_j(P_j(f_a-f))$, and assign zero outside
it; the resulting net remains order null, and no common finite-domain
tail across all representing coordinates is required. Because this
finite estimate is stronger than its tail-infimum consequence
$R_j^A(h)\le T_j(P_j(h))$, continuous inversion suffices in the former
case without a quantitative domination of inversion.

\begin{proposition}[Positive domination of inversion]
\label{uc:reciprocal-domination}
Let the original complete separated convergence group be faithfully
represented by possibly asymmetric subadditive gauges
$\rho_j:Z\to(F_j)_+$, with order convergence and the corresponding
pair equivalence. Suppose, for every $j$, a finite set $K_j$ and
positive order-continuous additive maps $A_{jk}:F_k\to F_j$ satisfy
\begin{equation}\label{uc:inversion-bound}
 \rho_j(-z)\le\sum_{k\in K_j}A_{jk}\rho_k(z)\qquad(z\in Z).
\end{equation}
Retain the one-sided conjugacy inequality
$R_j^A(h)\le T_j(P_j(h))$ and the input hypotheses of
Theorem~\ref{uc:conjugate}. Then its integrability and uniform
convergence conclusions hold, and
\[
 \rho_j(I_Af_a-I_Af)=R_j^A(f_a-f)\le b_{j,a},\qquad
 \rho_j(I_Af-I_Af_a)\le\sum_{k\in K_j}A_{jk}b_{k,a},
 \quad b_{j,a}=T_j(P_j(f_a-f)),
\]
on the respective finite-domain tails. The first bound holds at every
index where $b_{j,a}$ is finite.
\end{proposition}
\begin{proof}
Define $\rho_j^s(z)=\rho_j(z)\vee\rho_j(-z)$. It is symmetric,
vanishes at zero, and
\[
 \rho_j^s(z+w)
 \le[\rho_j(z)+\rho_j(w)]\vee[\rho_j(-z)+\rho_j(-w)]
 \le\rho_j^s(z)+\rho_j^s(w).
\]
Moreover,
\[
 \rho_j(z)\le\rho_j^s(z)
 \le\rho_j(z)+\sum_{k\in K_j}A_{jk}\rho_k(z).
\]
Finite sums and positive order-continuous images preserve order
nullity. Thus $\rho_j^s$ faithfully represents the same original
point and pair convergences.

Put $e_{a,\gamma}^A=S_\gamma^A(f_a-f)$. For each finite $b_{k,a}$, set
\[
 M^{A,k}_{a,\delta}=\sup_{\gamma\ge\delta}\rho_k(e^A_{a,\gamma}),
 \qquad R^{A,k}_a=\inf_\delta M^{A,k}_{a,\delta}\le b_{k,a}.
\]
Lemma~\ref{app:tail-envelope} supplies $\delta_0$ with
$M^{A,k}_{a,\delta_0}<\infty$. For $\gamma\ge\delta_0$, put
\[
 r^{A,k}_{a,\gamma}=M^{A,k}_{a,\gamma}-R^{A,k}_a\downarrow0.
\]
Values before this tail may be set to zero. On the displayed tail,
\[
 \rho_k(e_{a,\gamma}^A)
 \le M^{A,k}_{a,\gamma}
 =R^{A,k}_a+r^{A,k}_{a,\gamma}
 \le b_{k,a}+r_{a,\gamma}^{A,k}.
\]
For fixed $j,A$, with $a$ on the common finite-domain tail for
$\{j\}\cup K_j$, intersect the finitely many inner tails indexed by
$\{j\}\cup K_j$. There,
\[
 \rho_j^s(e_{a,\gamma}^A)
 \le b_{j,a}+\sum_{k\in K_j}A_{jk}b_{k,a}
 +r_{a,\gamma}^{A,j}
 +\sum_{k\in K_j}A_{jk}r_{a,\gamma}^{A,k}.
\]
The outer bound is order null and independent of $A$. The last two
terms form an order-null residual. Theorem~\ref{lt:order-transfer},
with the symmetric faithful distances $d_j(z,w)=\rho_j^s(z-w)$,
constructs $I_Af$ and gives uniform endpoint convergence.

Subadditivity in both directions gives
\[
 |\rho_j(z)-\rho_j(w)|
 \le\rho_j(z-w)\vee\rho_j(w-z)=\rho_j^s(z-w).
\]
Put $\ell=\rho_j(I_Af_a-I_Af)$. Since
$e_{a,\gamma}^A\to I_Af_a-I_Af$, the displayed inequality and
faithfulness give a decreasing zero regulator $v_\eta$ and indices
$\delta_\eta$ such that
\[
 \gamma\ge\delta_\eta\quad\Longrightarrow\quad
 \ell-v_\eta\le\rho_j(e_{a,\gamma}^A)\le\ell+v_\eta.
\]
Every nonempty later tail has finite supremum between the same
bounds. Since the indices above $\delta_\eta$ are cofinal,
\[
 \ell-v_\eta\le
 \inf_{\delta\ge\delta_\eta}\sup_{\gamma\ge\delta}
                  \rho_j(e_{a,\gamma}^A)
 =R_j^A(f_a-f)\le\ell+v_\eta.
\]
Taking the infimum of $v_\eta$ gives
$R_j^A(f_a-f)=\ell$. One-sided conjugacy then gives
$\rho_j(I_Af_a-I_Af)\le b_{j,a}$. Applying
\eqref{uc:inversion-bound} to $z=I_Af_a-I_Af$ gives
\[
 \rho_j(I_Af-I_Af_a)
 \le\sum_{k\in K_j}A_{jk}\rho_k(I_Af_a-I_Af)
 \le\sum_{k\in K_j}A_{jk}b_{k,a}
\]
on the common finite-domain tail of the indicated coordinates.
\end{proof}

Symmetric gauges satisfy \eqref{uc:inversion-bound} with
$K_j=\{j\}$ and $A_{jj}=\mathrm{id}$. On a Dedekind-complete Riesz space,
\[
 \rho(z)=(z^+,z^-),\qquad A(u,v)=(v,u)
\]
gives a nonsymmetric example, since $\rho(-z)=A\rho(z)$.
In each instance, the maps in \eqref{uc:inversion-bound} concern
the original gauges alone, so this sufficient condition can be
verified before an integral or candidate limit is known and without
a Cauchy assertion. Its force lies in the displayed domination,
which is not replaced here by continuity of inversion alone.

\section{Local Henstock control and monotone limits}
\label{lc:section}\label{ls:section}

Uniform approximation compares complete integration records, whereas
local approximation requires us to retain some cells and replace others;
to make that passage, we must understand how the error of a full sum
controls the errors carried by its partial divisions. We develop this
question through the finite subdivision calculus of \cite{henstockmagnum},
following the local estimates until they construct limits in the original
convergence. Along the way, the ordered monotone and dominated principles
of \cite[Chapter~5]{integralmeasureandordering} emerge as concrete
specializations, with their approximation conditions distinguished from
the structural properties that allow the local estimates to be combined.

\subsection{Local divisions and the original response nets}

\begin{definition}[Cell divisions and local constraints]\label{ls:calculus}
Let $T$ be a set of tags and $\mathscr C$ a set of cells. Each cell $C$
has a nonempty tag carrier $\tau(C)\subseteq T$ and a family
$\mathfrak D(C)$ of nonempty finite lists of cells. Require
\[
 (C)\in\mathfrak D(C),\qquad
 (C_r)_{r=1}^{q}\in\mathfrak D(C),\quad
 (D_{rs})_{s=1}^{q_r}\in\mathfrak D(C_r)
 \Longrightarrow (D_{rs})_{r,s}\in\mathfrak D(C).
\]
Lists are invariant under permutation. Repeated entries retain their
occurrence indices. Each $t\in T$ carries a meet semilattice $D_t$ with
largest element $1_t$, and sets $\mathscr V_t(d)$ satisfying
\[
 \mathscr V_t(d)\subseteq\{C:t\in\tau(C)\},\quad
 d\le e\Longrightarrow\mathscr V_t(d)\subseteq\mathscr V_t(e),
 \quad \mathscr V_t(1_t)=\{C:t\in\tau(C)\}.
\]
A gauge is an arbitrary element of $\mathscr G=\prod_{t\in T}D_t$.
For $\delta\in\mathscr G$, put
\[
 \mathscr R_A(\delta)=
 \left\{((C_r,t_r))_{r=1}^{q}:
 (C_r)_{r=1}^{q}\in\mathfrak D(A),\
 C_r\in\mathscr V_{t_r}(\delta(t_r))\ (1\le r\le q)\right\}.
\]
A partial record on $A$ is a subword of a tagged division whose cell
list belongs to $\mathfrak D(A)$.
\end{definition}

Nothing in this division calculus requires disjoint carriers: cells and
their coefficients retain their occurrence indices under finite
substitution, so that the argument follows the recorded occurrences.

\begin{proposition}[Finite divisions from branching]\label{ls:branch}
For each $A\in\mathscr C$, suppose there is a rooted tree of cell
occurrences, with root $A$, such that every node has a nonempty finite
list of children belonging to the division family of its cell.
Suppose every infinite branch $(C_n)_{n\ge0}$ satisfies
\begin{equation}
\exists t\in\bigcap_{n\ge0}\tau(C_n)\quad
 \forall d\in D_t\quad\exists n\ge0:\ C_n\in\mathscr V_t(d).
 \label{ls:eq-1}
\end{equation}
Then $\mathscr R_A(\delta)\ne\varnothing$ for every cell $A$ and every
gauge $\delta$.
\end{proposition}
\begin{proof}
Fix $A,\delta$. Call a node $v$ good when
\[
 G_v:=\{t:C_v\in\mathscr V_t(\delta(t))\}\ne\varnothing.
\]
Stop each good node and retain every child of each bad node. If the
resulting tree were infinite, its root would have infinitely many
remaining descendants. Number the children of every node. Recursively,
the first child with infinitely many remaining descendants exists,
since a finite union of finite sets is finite. This constructs an
infinite branch consisting of bad nodes. By \eqref{ls:eq-1},
\[
 t\in\bigcap_n\tau(C_n),\qquad
 C_n\in\mathscr V_t(\delta(t))\quad\text{for some }n,
\]
contradicting $G_{C_n}=\varnothing$. The stopped tree is finite.
Every leaf is good. Its leaves form a division of $A$ by finite
substitution, and
\[
 \prod_{v\text{ leaf}}G_v\ne\varnothing.
\]
Every tuple in this finite product gives a member of
$\mathscr R_A(\delta)$.
\end{proof}

With the branching data of Proposition~\ref{ls:branch} assumed on
every cell, fine records are now available wherever a division leaves
cells to be filled. This is precisely the existence statement needed to
complete a compatible partial record while retaining its chosen entries.

\begin{proposition}[Completion of compatible partial records]
\label{ls:completion}
Let $Q$ be a partial record on $A$, with a witnessing full cell division
$D=(C_1,\ldots,C_q)$. Suppose every entry of $Q$ is $\delta$-fine.
For any gauges $\eta_r$ on the omitted cells, there exist records
\[
 T_r\in\mathscr R_{C_r}(\delta\wedge\eta_r)
 \quad(r\text{ omitted})
\]
such that $Q\sqcup\bigsqcup_rT_r\in\mathscr R_A(\delta)$.
\end{proposition}
\begin{proof}
Let $K\subseteq\{1,\ldots,q\}$ be the retained occurrence set and
$M=\{1,\ldots,q\}\setminus K$. Proposition~\ref{ls:branch} gives
\[
 \mathscr R_{C_r}(\delta\wedge\eta_r)\ne\varnothing\quad(r\in M),
 \qquad
 \prod_{r\in M}\mathscr R_{C_r}(\delta\wedge\eta_r)\ne\varnothing.
\]
For every tuple $(T_r)_{r\in M}$ in this finite product, put
$R=Q\sqcup\bigsqcup_{r\in M}T_r$. Finite substitution and gauge
monotonicity yield
\[
 \operatorname{div}R\in\mathfrak D(A),\qquad
 (B,t)\in T_r\Longrightarrow
 B\in\mathscr V_t(\delta(t)\wedge\eta_r(t))
 \subseteq\mathscr V_t(\delta(t)).
\]
Every occurrence of $Q$ is $\delta$-fine by hypothesis. Therefore
$R\in\mathscr R_A(\delta)$. For $M=\varnothing$, the product is the
singleton containing the empty tuple and $R=Q$.
\end{proof}

In this completion only the residual cells are refined, while the given
partial record is retained; its fineness under subsequent gauges is not
asserted. To follow such completions through convergence, we next record
the fine divisions in a net whose tails preserve their exact meaning.

\begin{definition}[The full record net]\label{ls:record-net}
Put
\[
 \Gamma_A=\{(\delta,R,k):\delta\in\mathscr G,
                 R\in\mathscr R_A(\delta),\ k\in\mathbb N\}.
\]
For distinct indices, set
\[
 (\delta,R,k)\le(\eta,S,l)
 \iff \eta\le\delta\text{ pointwise and }k<l.
\]
Equality is also allowed. The record net is
$\mathcal R_A(\delta,R,k)=R$.
\end{definition}

\begin{proposition}[Exact record tails]\label{ls:tails}
The order on $\Gamma_A$ is directed, and
\begin{equation}
\{\mathcal R_A(\xi):\xi\ge(\delta,R,k)\}
 =\mathscr R_A(\delta).
 \label{ls:eq-2}
\end{equation}
For an originally prescribed record net $(R_\alpha)$, write
$T_{\alpha_0}=\{R_\alpha:\alpha\ge\alpha_0\}$. It is tail equivalent
to $\mathcal R_A$ if and only if
\begin{equation}
\forall\delta\ \exists\alpha_0:\ T_{\alpha_0}\subseteq\mathscr R_A(\delta),
 \qquad
 \forall\alpha_0\ \exists\delta:\ \mathscr R_A(\delta)\subseteq T_{\alpha_0}.
 \label{ls:eq-3}
\end{equation}
Every response map preserves this tail equivalence.
\end{proposition}
\begin{proof}
For $\xi_r=(\delta_r,R_r,k_r)$, $r=1,2$, put
\[
 \eta=\delta_1\wedge\delta_2,\qquad l=1+\max\{k_1,k_2\}.
\]
Proposition~\ref{ls:branch} gives
\[
 S\in\mathscr R_A(\eta)\Longrightarrow
 (\eta,S,l)\ge\xi_1,\xi_2,
 \qquad \mathscr R_A(\eta)\ne\varnothing.
\]
For distinct indices, both $\xi\le\zeta$ and $\zeta\le\xi$ would
require $k_\xi<k_\zeta<k_\xi$. Transitivity is the implication
\[
 \epsilon\le\eta\le\delta,\quad m>l>k
 \Longrightarrow\epsilon\le\delta,\quad m>k.
\]
For $\xi_0=(\delta,R,k)$,
\[
 \xi=(\eta,S,l)\ge\xi_0
 \Longrightarrow S\in\mathscr R_A(\eta)\subseteq\mathscr R_A(\delta),
\]
\[
 S\in\mathscr R_A(\delta)
 \Longrightarrow(\delta,S,k+1)\ge\xi_0.
\]
These two implications prove \eqref{ls:eq-2}. Therefore the original
net is a quasi-subnet of the full record net exactly when
\[
 \forall\delta\ \exists\alpha_0:\quad
 T_{\alpha_0}\subseteq\mathscr R_A(\delta),
\]
and the reverse quasi-subnet relation is the second clause of
\eqref{ls:eq-3}. For any map $V$,
\[
 T_{\alpha_0}\subseteq\mathscr R_A(\delta)
 \Longrightarrow V[T_{\alpha_0}]\subseteq V[\mathscr R_A(\delta)],
\]
\[
 \mathscr R_A(\delta)\subseteq T_{\alpha_0}
 \Longrightarrow V[\mathscr R_A(\delta)]\subseteq V[T_{\alpha_0}].
\]
Thus the two response nets are tail equivalent.
\end{proof}

We henceforth use these full record nets, or original integration nets
satisfying \eqref{ls:eq-3}, taking admissible representatives at the record
level and retaining every decoration carried by a record. With the tails
thus fixed, local constraints can be combined at each tag before any
response is formed.

\begin{theorem}[Local finite synchronization]\label{ls:synchronization}
Let $(\delta_i)_{i\in I}$ be any family of gauges and let
$F:T\to\operatorname{Fin}(I)$. For $c\in\mathscr G$, define
\begin{equation}
\delta_{c,F}(t)=c(t)\wedge\bigwedge_{i\in F(t)}\delta_i(t),
 \qquad \bigwedge\varnothing=1_t.
 \label{ls:eq-4}
\end{equation}
Then $\mathscr R_A(\delta_{c,F})\ne\varnothing$, and every subword
of a record in this family whose tags satisfy $i\in F(t)$ is
$\delta_i$-fine.
In particular, for $p:T\to\mathbb N_{\ge1}$,
\begin{equation}
\delta^{p}(t)=\bigwedge_{m=1}^{p(t)}\delta_m(t)
 \label{ls:eq-5}
\end{equation}
is a gauge, and every subword with $p(t)\ge m$ is $\delta_m$-fine.
\end{theorem}
\begin{proof}
For each $t$, the meet in \eqref{ls:eq-4} is finite and belongs to $D_t$.
Hence $\delta_{c,F}\in\mathscr G$, and
Proposition~\ref{ls:branch} supplies a fine division. If $i\in F(t)$,
\[
 \delta_{c,F}(t)\le\delta_i(t),\qquad
 C\in\mathscr V_t(\delta_{c,F}(t))
 \Longrightarrow C\in\mathscr V_t(\delta_i(t)).
\]
Taking $F(t)=\{1,\ldots,p(t)\}$ proves the last assertion.
\end{proof}

Because the meet is local, the index set $I$ may be uncountable: only
the finite set $F(t)$ enters at each tag, and the full product of the local
fibres makes \eqref{ls:eq-4} legitimate even when
$\sup_t|F(t)|=\infty$.

Once compatible records can be assembled, finite additivity of their
known limits can be derived directly. For this purpose, let $E$ be a set,
let $Z$ be a separated commutative convergence monoid, and let
$L_C:E\to Z$ be arbitrary maps; equivalently, the responses may be formed
from an arbitrary kernel $U(t,C)\in Z$. The argument requires neither
additive inverses nor additivity of $L_C$. Write
\[
 S_Rh=\sum_{r=1}^{q}L_{C_r}h(t_r),\qquad
 J_h(C)=\lim_{\Gamma_C}S_Rh
\]
when the latter limit exists in the original convergence.

\begin{proposition}[Derived finite additivity]\label{ls:additivity}
If $h$ is integrable on $A$ and on every cell of
$(C_1,\ldots,C_q)\in\mathfrak D(A)$, then
\begin{equation}
J_h(A)=\sum_{r=1}^{q}J_h(C_r).
 \label{ls:eq-6}
\end{equation}
\end{proposition}
\begin{proof}
Define
\[
 \Xi=\left\{(\delta,R_1,\ldots,R_q,k):
 R_r\in\mathscr R_{C_r}(\delta),\ k\in\mathbb N\right\},
\]
with equality and the strict order given by gauge refinement and
strict increase of $k$. For $\xi_1,\xi_2\in\Xi$, put
\[
 \eta=\delta_1\wedge\delta_2,\qquad
 l=1+\max\{k_1,k_2\},\qquad
 \prod_{r=1}^{q}\mathscr R_{C_r}(\eta)\ne\varnothing.
\]
Every tuple $(T_1,\ldots,T_q)$ in the displayed product gives the
common upper bound $(\eta,T_1,\ldots,T_q,l)$. At every
$\xi=(\delta,R_1,\ldots,R_q,k)$, finite substitution gives
\[
 R(\xi):=\bigsqcup_{r=1}^{q}R_r\in\mathscr R_A(\delta).
\]
For any gauge $\epsilon$, every index with gauge $\delta\le\epsilon$
satisfies
\[
 R_r\in\mathscr R_{C_r}(\epsilon)\quad(1\le r\le q),\qquad
 R(\xi)\in\mathscr R_A(\epsilon).
\]
The response nets on $\Xi$ are consequently quasi-subnets of the
original full record nets. Their known limits and exact sums are
\[
 S_{R_r}h\longrightarrow J_h(C_r),\qquad
 S_{R(\xi)}h\longrightarrow J_h(A),\qquad
 S_{R(\xi)}h=\sum_{r=1}^{q}S_{R_r}h.
\]
Finite addition gives $S_{R(\xi)}h\to\sum_{r=1}^{q}J_h(C_r)$.
Separation yields $J_h(A)=\sum_{r=1}^{q}J_h(C_r)$.
\end{proof}

\subsection{The original convergence and its lattice representation}

Having identified the response tails, we turn to the convergence in
which their errors are to disappear. For a decreasing family of nonempty
sets $(B_\alpha)_{\alpha\in D}$, let $\mathsf N(B_\alpha)$ denote the
saturated net
\[
 (\alpha,b,k)\longmapsto b,\qquad b\in B_\alpha,\ k\in\mathbb N,
\]
ordered by equality and
\[
 (\alpha,b,k)<(\beta,c,l)\iff\alpha\le\beta,\ k<l.
\]
Its tail image at $(\alpha,b,k)$ is exactly $B_\alpha$:
\[
 c\in B_\alpha\Longrightarrow(\alpha,c,k+1)>(\alpha,b,k),
 \qquad \beta\ge\alpha\Longrightarrow B_\beta\subseteq B_\alpha.
\]
We interpret every product or saturation below through its
tail-equivalent admissible representative from
Definition~\ref{fd:admissible}, with
\eqref{fd:tail-representation} preserving every displayed eventual
bound. For an admissible net $x=(x_\alpha)$, write
$T_\alpha(x)=\{x_\beta:\beta\ge\alpha\}$.

\begin{definition}[Locally solid convergence and faithful representation]
\label{lc:faithful}
Let $Z$ be an abelian lattice-ordered group with a separated group
convergence $\lambda$. For $B\subseteq Z$, put
\[
 \operatorname{sol}B=\{z:\exists b\in B,\ |z|\le|b|\}.
\]
Local solidity means that, for every admissible null net,
\[
 x_\alpha\longrightarrow_\lambda0
 \Longrightarrow
 \mathsf N\bigl(\operatorname{sol}T_\alpha(x)\bigr)
 \longrightarrow_\lambda0.
\]
The convergence is sequentially order compatible when
\[
 z_n\downarrow0\Longrightarrow z_n\to_\lambda0.
\]
A faithful solid lattice representation consists of maps
$\rho_j:Z\to(F_j)_+$ satisfying
\[
 \rho_j(0)=0,\quad \rho_j(-z)=\rho_j(z),\quad
 \rho_j(z+w)\le\rho_j(z)+\rho_j(w),
\]
\[
 |z|\le|w|\Longrightarrow\rho_j(z)\le\rho_j(w),
\]
and, for every admissible pair net $((z_\alpha,w_\alpha))$,
\begin{equation}\label{lc:pair-faithfulness}
 (z_\alpha,w_\alpha)\longrightarrow_{\mathcal U_\lambda}\Delta_Z
 \iff z_\alpha-w_\alpha\longrightarrow_\lambda0
 \iff \rho_j(z_\alpha-w_\alpha)\longrightarrow_D0\quad(j).
\end{equation}
Here $\mathcal U_\lambda$ is the original difference uniform convergence,
and $D$ denotes the regulator convergence in the order appendix.
\end{definition}

In applying the metrizability theorem through
\eqref{lc:pair-faithfulness}, we use both implications and retain the
actual representation throughout, so that the solidness of its norms
and the controller properties introduced below must be verified for
that representation. Any replacement of a controller by a completion
would therefore require a separate proof of faithfulness.

\begin{proposition}[Order compatibility of the representing norms]
\label{lc:norm-normality}
If $\lambda$ is sequentially order compatible and every $F_j$ is weakly
$\sigma$-distributive, then
\[
 z_n\downarrow0\Longrightarrow\rho_j(z_n)\downarrow0.
\]
\end{proposition}
\begin{proof}
Solidness gives $\rho_j(z_{n+1})\le\rho_j(z_n)$, while sequential
order compatibility and faithfulness give a regulator $a$ satisfying
\[
 \forall\varphi\ \exists N(\varphi)\ \forall n\ge N(\varphi):
 \quad \rho_j(z_n)\le H_a(\varphi).
\]
If $v\in F_j$ is a lower bound for the whole sequence, then
$\rho_j(z_n)\ge0$ implies $v^+=v\vee0\le\rho_j(z_n)$ for every $n$.
In particular,
\[
 0\le v^+\le\rho_j(z_{N(\varphi)})\le H_a(\varphi)
 \quad(\varphi\in\mathbb N^{\mathbb N}),
 \qquad
 v^+\le\inf_\varphi H_a(\varphi)=0.
\]
Thus every lower bound $v$ in the whole controller group satisfies
$v\le0$, and $0$ itself is a lower bound. This proves
$\inf_n\rho_j(z_n)=0$.
\end{proof}

The role of order compatibility can already be seen on constants,
which provide an independent necessity test. In a normalized
constant-function integration structure,
\[
 S_R\mathbf z=z\quad(R),\qquad I\mathbf z=z.
\]
For $z_n\downarrow0$, the increasing constants $h_n=-\mathbf z_n$
have primitives bounded above by zero, so a monotone theorem valid for
these data gives $-z_n\to_\lambda0$, hence $z_n\to_\lambda0$.

\subsection{Positive responses and partial Henstock estimates}

To carry the local estimates into the ordered theory, let $E$ be a
Dedekind $\sigma$-complete abelian lattice-ordered group and let $Z$
also be Dedekind $\sigma$-complete, with a sequentially order compatible
locally solid convergence and the faithful representation of
Definition~\ref{lc:faithful}. Each controller $F_j$ is Dedekind
$\sigma$-complete and weakly $\sigma$-distributive; stronger controller
conditions will enter only in the bounded-primitives branch below.
Fix a root cell $A$ and positive additive coefficient maps $L_C:E\to Z$
satisfying
\[
 L_C=\sum_{r=1}^mL_{C_r}
 \quad((C_r)_{r=1}^m\in\mathfrak D(C)),
 \qquad e_k\downarrow0\Longrightarrow L_Ae_k\downarrow0.
\]
Write $\mathscr C_A$ for the cells occurring in divisions of $A$,
including finite substitutions. For compatible partial records,
\begin{equation}\label{lc:coefficient-bound}
 0\le\sum_{(C,t)\in Q}L_C\le L_A.
\end{equation}
For a tagged record $R$, put $S_Rh=\sum_{(C,t)\in R}L_Ch(t)$.
Known hereditary original integrals are denoted by $J_h(C)$.

Finite additivity has already been established in
Proposition~\ref{ls:additivity}, and positivity follows from the same
original convergence: if $h\ge0$ and $S_Rh\to_\lambda J_h(C)$, then
\[
 0\le J_h(C)^-\le|S_Rh-J_h(C)|.
\]
Local solidity therefore makes the constant value $J_h(C)^-$ converge
to zero, and separation gives $J_h(C)^-=0$.

Before using that order, however, we isolate the partial Henstock
estimates that follow from completion of finite responses and invariant
distances alone, so that this part of the argument remains available
without a lattice order on $E$ or $Z$.

\begin{lemma}[Removal of a residual response in a monoid]
\label{lc:monoid-completion}
Let $Z$ be a commutative monoid and let $d:Z^2\to M$ be symmetric,
translation invariant, and satisfy the triangle inequality. The
controller $M$ is either the positive cone of a partially ordered
abelian group, with decreasing control infima taken in the entire group,
or a partially ordered commutative monoid with least identity $0$,
isotone addition, and
\begin{equation}\label{lc:controller-translation-inf}
 v_\theta\downarrow0\quad\Longrightarrow\quad
 \inf_\theta(u+v_\theta)=u\qquad(u\in M).
\end{equation}
Let $u\in M$. Suppose $x,a,b\in Z$ and
$d(y_\beta,b)\to_o0$. If
\[
 \forall\beta_0\ \exists\beta\ge\beta_0:
 d(x+y_\beta,a+b)\le u,
\]
then $d(x,a)\le u$.
\end{lemma}
\begin{proof}
Let $v_\theta\downarrow0$ control $d(y_\beta,b)$. For every $\theta$,
there is $\beta_\theta$ such that
\[
 \beta\ge\beta_\theta\Longrightarrow d(y_\beta,b)\le v_\theta.
\]
The hypothesis applied at $\beta_\theta$ gives a value $\beta$ for which
both this inequality and the completion bound hold. Therefore
\[
\begin{aligned}
 d(x,a)&=d(x+b,a+b)\\
 &\le d(x+b,x+y_\beta)+d(x+y_\beta,a+b)\\
 &=d(b,y_\beta)+d(x+y_\beta,a+b)\le v_\theta+u.
\end{aligned}
\]
For a partially ordered group controller, $d(x,a)-u$ is a lower bound for every
$v_\theta$ in the entire partially ordered group, so $d(x,a)-u\le0$.
For a monoid controller, \eqref{lc:controller-translation-inf} gives
$d(x,a)\le\inf_\theta(u+v_\theta)=u$.
\end{proof}

\begin{lemma}[Partial Henstock estimate for a monoid range]
\label{lc:monoid-partial}
Let $Z$ be a separated commutative convergence monoid carrying invariant
symmetric pseudometrics $d_j:Z^2\to M_j$ with the controller assumptions
of Lemma~\ref{lc:monoid-completion}. Assume original convergence
$z_\alpha\to z$ implies $d_j(z_\alpha,z)\to_o0$.
Let $U$ have known hereditary original integrals and fix $j$ and
$u\in M_j$. If
\[
 d_j(S_RU,J_U(A))\le u\qquad(R\in\mathscr R_A(\eta)),
\]
then every compatible $\eta$-fine partial record $Q$ satisfies
\[
 d_j\!\left(\sum_{(C,t)\in Q}U(t,C),
             \sum_{(C,t)\in Q}J_U(C)\right)\le u.
\]
\end{lemma}
\begin{proof}
Let $(B_1,\ldots,B_m)$ be the omitted occurrences of a witnessing full
cell division. For $m=0$, the stated bound is the assumed full-record
bound. Suppose $m\ge1$, and set
\[
 x=\sum_{(C,t)\in Q}U(t,C),\qquad
 a=\sum_{(C,t)\in Q}J_U(C),\qquad b=\sum_{r=1}^mJ_U(B_r).
\]
Proposition~\ref{ls:additivity} gives $J_U(A)=a+b$. Use the directed set
\[
 \mathcal I=\{(\epsilon,T_1,\ldots,T_m,k):
 \epsilon\le\eta,\ T_r\in\mathscr R_{B_r}(\epsilon),\ k\in\mathbb N\},
\]
with the strict refinement order of Definition~\ref{ls:record-net}.
For two indices, the common gauge is their pointwise meet and the new
integer is one plus their maximum. Proposition~\ref{ls:branch} makes
all residual record families nonempty. Thus $\mathcal I$ is directed.
For $\beta=(\epsilon,T_1,\ldots,T_m,k)$, put
\[
 y_\beta=\sum_{r=1}^mS_{T_r}U,\qquad
 R_\beta=Q\sqcup\bigsqcup_{r=1}^mT_r.
\]
Every residual coordinate is a quasi-subnet of its full record net.
Finite addition and Proposition~\ref{ls:completion} give
\[
 y_\beta\longrightarrow b,\qquad R_\beta\in\mathscr R_A(\eta),
 \qquad S_{R_\beta}U=x+y_\beta.
\]
Hence $d_j(y_\beta,b)\to_o0$ and
$d_j(x+y_\beta,a+b)\le u$ for every $\beta$. Lemma~\ref{lc:monoid-completion}
proves the assertion.
\end{proof}

\begin{lemma}[Cofinal completion of a fixed partial response]
\label{lc:cofinal-completion}
Let $Z$ be an abelian convergence group with a faithful symmetric
subadditive lattice representation in WSD controllers. Let the known
integrals of $h$ on $A$ and a residual domain $B$ be $J_A,J_B$.
Fix a finite partial response $S_Qh$, a value $J_Q$, and a root tail
$\gamma\ge\gamma_0$ such that
\[
 J_A=J_Q+J_B,
 \qquad \rho_j(S_\gamma^Ah-J_A)\le u\quad(\gamma\ge\gamma_0).
\]
Suppose the response records admit cofinal completion:
\[
 \forall\beta_0\in\Gamma_B\quad
 \exists\beta\ge\beta_0\quad\exists\gamma\ge\gamma_0:
 \qquad S_\gamma^Av=S_Qv+S_\beta^Bv
 \quad\text{for every }v.
\]
Then $\rho_j(S_Qh-J_Q)\le u$.
\end{lemma}
\begin{proof}
Put $d_j(z,w)=\rho_j(z-w)$, $x=S_Qh$, $a=J_Q$, $b=J_B$, and
$y_\beta=S^B_\beta h$. Faithfulness and
Lemma~\ref{lt:d-calculus} give
\[
 d_j(y_\beta,b)=\rho_j(S^B_\beta h-J_B)\longrightarrow_D0
 \quad\Longrightarrow\quad d_j(y_\beta,b)\longrightarrow_o0.
\]
Every prescribed residual tail has a completion satisfying
\[
 d_j(x+y_\beta,a+b)
 =\rho_j(S^A_\gamma h-J_A)\le u.
\]
Lemma~\ref{lc:monoid-completion} yields
$d_j(x,a)=\rho_j(S_Qh-J_Q)\le u$.
\end{proof}

\begin{corollary}[Henstock estimate on a prescribed net]
\label{henstocklemma}
Let $S_\gamma^Ah$ be the original net Riemann sums. Suppose $h$ is
integrable on the domains below. For a record
$R_{\gamma_0}=((C_r,t_r))_{r=1}^m$ and $K\subseteq\{1,\ldots,m\}$,
write $Q=((C_r,t_r))_{r\in K}$ and let $B$ be the residual domain.
Assume the primitive decomposition
\[
 J_h(A)=\sum_{r\in K}J_h(C_r)+J_h(B)
\]
and the following record property:
\[
 \forall\beta_0\in\Gamma_B\quad\exists\beta\ge\beta_0\quad
 \exists\gamma\ge\gamma_0:
 \qquad R_\gamma=Q\sqcup R_\beta^B.
\]
If the original convergence has a faithful symmetric lattice norm
$\rho$ with WSD controller and
\[
 \rho(S_\gamma^Ah-J_h(A))\le H_a(\varphi)
 \qquad(\gamma\ge\gamma_0),
\]
then
\[
 \rho\!\left(\sum_{r\in K}[L_{C_r}h(t_r)-J_h(C_r)]\right)
 \le H_a(\varphi).
\]
The same conclusion holds simultaneously for every member of a faithful
family of symmetric lattice pseudonorms.
\end{corollary}
\begin{proof}
The record concatenation gives, for every function $v$,
\[
 S_\gamma^Av=\sum_{r\in K}L_{C_r}v(t_r)+S_\beta^Bv.
\]
In Lemma~\ref{lc:cofinal-completion}, put
\[
 S_Qv=\sum_{r\in K}L_{C_r}v(t_r),\qquad
 J_Q=\sum_{r\in K}J_h(C_r),\qquad u=H_a(\varphi).
\]
Hereditary integrability supplies the residual regulator. The assumed
primitive decomposition and record completion verify its other premises.
The lemma gives the asserted estimate for each represented coordinate.
\end{proof}

\begin{lemma}[Partial Henstock estimate in a lattice norm]
\label{lc:partial}
Suppose $h$ has known hereditary original integrals. Fix $j$ and suppose
\[
 \rho_j(S_Rh-J_h(A))\le u
 \quad(R\in\mathscr R_A(\eta)).
\]
For every compatible $\eta$-fine partial record $Q$,
\[
 \rho_j\left(\sum_{(C,t)\in Q}[L_Ch(t)-J_h(C)]\right)\le u.
\]
\end{lemma}
\begin{proof}
For $d_j(z,w)=\rho_j(z-w)$, symmetry, invariance, and the triangle
inequality are the defining properties of the lattice pseudonorm.
Original convergence implies regulator convergence of these distances.
Lemma~\ref{lt:d-calculus} converts it to order convergence. Apply
Lemma~\ref{lc:monoid-partial} to $U(t,C)=L_Ch(t)$. Its bound is
\[
 d_j\!\left(\sum_{(C,t)\in Q}L_Ch(t),
             \sum_{(C,t)\in Q}J_h(C)\right)
 =\rho_j\!\left(\sum_{(C,t)\in Q}[L_Ch(t)-J_h(C)]\right)\le u.
\]
\end{proof}

\begin{definition}[Two input approximation conditions]\label{lc:profiles}
Countable order profiles mean
\[
 e_k^r\downarrow_k0\quad(r\ge1),\qquad\jmath:T\to\mathbb N,
\]
\[
 n\ge\theta(t,k)\Longrightarrow
 |h_n(t)-h(t)|\le e_k^{\jmath(t)}.
\]
Common regulator approximation means that a regulator $a$ in $E$
satisfies
\[
 n\ge\theta(t,\varphi)\Longrightarrow
 |h_n(t)-h(t)|\le H_a(\varphi)
 \quad(t,\varphi).
\]
\end{definition}

\begin{proposition}[Countably indexed sources of input profiles]
\label{lc:profile-certificate}
If the tag set is countable and $h_n(t)\to_o h(t)$ at every tag in a
Dedekind $\sigma$-complete group $E$, the countable-profile condition
holds. More generally, it holds when the tag set has a countable cover
$T=\bigcup_qT_q$ and each $T_q$ admits an order-null sequence $e_k^q$
that eventually bounds $|h_n(t)-h(t)|$ separately at each $t\in T_q$.
\end{proposition}
\begin{proof}
For a fixed tag $t$, an eventual order bound $v$ and its finite initial
segment give a bound for the whole sequence:
\[
 |h_n(t)-h(t)|\le
 v\vee\bigvee_{1\le r<N}|h_r(t)-h(t)|\qquad(n\ge1).
\]
Thus
\[
 e_k(t)=\sup_{n\ge k}|h_n(t)-h(t)|
\]
exists. If $v_\alpha\downarrow0$ controls pointwise order convergence,
then for every $\alpha$ there is $k$ with $e_k(t)\le v_\alpha$.
Consequently $e_k(t)\downarrow0$. Enumerating the tags as $(t_q)$,
with first occurrence used for repetitions, gives
\[
 e_k^q=e_k(t_q),\qquad\jmath(t_q)=q,\qquad\theta(t_q,k)=k.
\]
For a countable cover, set $\jmath(t)=\min\{q:t\in T_q\}$. The least
integer after which the stated bound holds is $\theta(t,k)$. These
assignments verify every quantifier in Definition~\ref{lc:profiles}.
\end{proof}

\subsection{Envelope monotone convergence in the original convergence}
\label{fv:envelope-section}

With partial errors controlled, we can ask how an integrable envelope
allows the estimates for different approximants to be combined. Let $E$
and $Z$ be Dedekind $\sigma$-complete abelian lattice-ordered groups,
and retain the original separated locally solid group convergence
$\lambda$, its original difference uniform convergence, and the faithful
solid lattice representation of Definition~\ref{lc:faithful}. Assume
that $\lambda$ is sequentially order compatible and that every controller
$F_j$ is a Dedekind $\sigma$-complete, weakly $\sigma$-distributive
abelian lattice-ordered group, so that
\begin{equation}\label{fv:rho-faithful}
 z_\alpha\longrightarrow_\lambda0
 \quad\Longleftrightarrow\quad
 \rho_j(z_\alpha)\longrightarrow_D0\quad(j\in\mathcal J).
\end{equation}
Throughout this argument, integer multiples denote repeated addition;
neither additivity nor finite-supremum preservation of the gauges is
required.

On the integration side, we retain the local division calculus, finite
substitution, compatible partial records, gauge completion, and the
original response-tail correspondence, with positive additive
coefficients satisfying
\begin{equation}\label{fv:envelope-coefficients}
 L_C=\sum_{D\in\mathcal D}L_D
 \quad(\mathcal D\in\mathfrak D(C)),\qquad
 x_k\downarrow0\ \Longrightarrow\ L_Ax_k\downarrow0.
\end{equation}
Here $A$ is the root cell, and $\mathscr C_A$ consists of the cells
occurring in a division of $A$. On every compatible partial record,
\[
 x\in E_+\quad\Longrightarrow\quad
 0\le\sum_{(C,t)\in Q}L_Cx\le L_Ax.
\]
All saturation and product nets use the admissible representatives
of Definition~\ref{fd:admissible}.

Let $g:T\to E_+$ be hereditarily integrable on $\mathscr C_A$, and
let the hereditarily integrable functions $h_m$ satisfy
\begin{equation}\label{fv:known-envelope}
 |h_m(t)|\le g(t)\qquad(m\ge1,\ t\in T).
\end{equation}
Write $G(C)=J_g(C)$ and $J_m(C)=J_{h_m}(C)$. If
$z_\alpha\ge0$ and $z_\alpha\to_\lambda z$, then
\[
 0\le z^-\le|z_\alpha-z|.
\]
Local solidity makes the constant net at $z^-$ converge to zero, and
separation gives $z^-=0$. Applying this implication to the known limits
of the positive finite sums of $g$, $g-h_m$, and $g+h_m$ therefore gives
\begin{equation}\label{fv:primitive-envelope}
 G(C)\ge0,\qquad |J_m(C)|\le G(C).
\end{equation}
All these known primitives are finitely additive.

\begin{lemma}[Bounds obtained from the known envelope]
\label{fv:derived-caps}
For each $j$ there are a gauge $\eta_j$ and elements
$w_j,v_j\in(F_j)_+$ such that every $\eta_j$-fine compatible partial
root record $Q=((C_r,t_r))$ and every assignment of positive integer
labels $\nu_r$ satisfy
\begin{equation}\label{fv:mixed-cap}
 \rho_j\!\left(
  \sum_r[L_{C_r}h_{\nu_r}(t_r)-J_{\nu_r}(C_r)]
 \right)\le w_j.
\end{equation}
If $|h|\le g$, then
\begin{equation}\label{fv:data-cap}
 \rho_j\!\left(\sum_rL_{C_r}|h(t_r)-h_{\nu_r}(t_r)|\right)
 \le v_j.
\end{equation}
\end{lemma}
\begin{proof}
An attained regulator bound for the known root integral of $g$ gives
a gauge $\eta_j$ and $u_j\in(F_j)_+$ with
\[
 \rho_j(S_Rg-G(A))\le u_j
 \qquad(R\in\mathscr R_A(\eta_j)).
\]
Gauge completion extends $Q$ to such a record $R$. Positivity and
finite additivity give
\[
 0\le S_Qg\le S_Rg,\qquad
 \sum_{(C,t)\in Q}G(C)\le G(A).
\]
For $E_Q=\sum_r[L_{C_r}h_{\nu_r}(t_r)-J_{\nu_r}(C_r)]$,
\begin{align*}
 |E_Q|
 &\le\sum_r\bigl[L_{C_r}g(t_r)+G(C_r)\bigr]
 \le S_Rg+G(A),\\
 \rho_j(E_Q)
 &\le\rho_j(S_Rg+G(A))\\
 &\le\rho_j(S_Rg-G(A))+\rho_j(2G(A))
 \le u_j+\rho_j(2G(A)).
\end{align*}
Thus $w_j=u_j+\rho_j(2G(A))$ satisfies
\eqref{fv:mixed-cap}.

For $D_Q=\sum_rL_{C_r}|h(t_r)-h_{\nu_r}(t_r)|$,
\[
 0\le D_Q\le2S_Qg\le2S_Rg,
\]
and hence
\begin{align*}
 \rho_j(D_Q)
 &\le\rho_j(2S_Rg)\\
 &\le2\rho_j(S_Rg-G(A))+\rho_j(2G(A))\\
 &\le2u_j+\rho_j(2G(A)).
\end{align*}
This proves \eqref{fv:data-cap} with
$v_j=2u_j+\rho_j(2G(A))$.
\end{proof}

\begin{lemma}[Mixed signed Henstock estimates]
\label{fv:envelope-mixed}
For each $j$ there is a regulator $B^j$ in $F_j$ such that
\[
 \forall\varphi\in\mathbb N^{\mathbb N}\ \forall p:T\to\mathbb N
 \ \exists\delta_{j,\varphi,p}\le\eta_j
\]
and, for every $C\in\mathscr C_A$, every
$Q=((C_r,t_r))\in\mathscr R_C(\delta_{j,\varphi,p})$, and every
$1\le\nu_r\le p(t_r)$,
\begin{equation}\label{fv:envelope-mixed-bound}
 \rho_j\!\left(
  \sum_r[L_{C_r}h_{\nu_r}(t_r)-J_{\nu_r}(C_r)]
 \right)\le H_{B^j}(\varphi).
\end{equation}
\end{lemma}
\begin{proof}
For every known row $m$, the signed partial Henstock lemma supplies
a regulator $a^{j,m}$ and gauges $\delta_{j,m,\psi}$ such that
\[
 \rho_j\!\left(
 \sum_{(C,t)\in Q}[L_Ch_m(t)-J_m(C)]\right)
 \le H_{a^{j,m}}(\psi)
\]
for every compatible $\delta_{j,m,\psi}$-fine partial root record.
No finite-supremum operation on defects is used here.

For fixed $\varphi,p$, put
\[
 \psi_{m,\varphi}(i)=\varphi(i+m),\qquad
 \delta_{j,\varphi,p}(t)
 =\eta_j(t)\wedge
   \bigwedge_{1\le m\le p(t)}\delta_{j,m,\psi_{m,\varphi}}(t).
\]
The meet is finite at every tag. The subword with label $m$ is fine
for $\delta_{j,m,\psi_{m,\varphi}}$. With
$M=\max_r\nu_r$, subadditivity and
Lemma~\ref{fv:derived-caps} give
\[
 \rho_j(E_Q)
 \le w_j\wedge
 \sum_{m=1}^{M}H_{a^{j,m}}(\psi_{m,\varphi}).
\]
The bounded-fusion lemma gives \eqref{fv:envelope-mixed-bound}. An
explicit controller is
\[
 B^j_{1s}=0,\qquad
 B^j_{rs}
 =w_j\wedge2^r\sum_{1\le m<r}a^{j,m}_{r-m,s}
 \quad(r\ge2).
\]
Each finite sum decreases to zero in $s$, and
$0\le B^j_{rs}\le w_j$. The controller is independent of $p,C,Q$.
\end{proof}

\begin{lemma}[Transfer of countable order profiles]
\label{fv:envelope-data}
Suppose $|h|\le g$, and let
\[
 e_k^q\downarrow0\text{ in }E_+,\qquad
 \jmath:T\to\mathbb N,\qquad
 \theta:T\times\mathbb N\to\mathbb N
\]
satisfy
\[
 n\ge\theta(t,k)\quad\Longrightarrow\quad
 |h_n(t)-h(t)|\le e_k^{\jmath(t)}.
\]
For every $j$ there is a regulator $D^j$ such that, with
\[
 \ell_\varphi(t)=\theta(t,\varphi(\jmath(t)+1)),
\]
every $\eta_j$-fine compatible partial root record $Q=((C_r,t_r))$
and every $\nu_r\ge\ell_\varphi(t_r)$ satisfy
\begin{equation}\label{fv:envelope-data-bound}
 \rho_j\!\left(\sum_rL_{C_r}|h(t_r)-h_{\nu_r}(t_r)|\right)
 \le H_{D^j}(\varphi).
\end{equation}
\end{lemma}
\begin{proof}
Define
\[
 r^{j,q}_k=\rho_j(L_Ae_k^q).
\]
Root normality and Proposition~\ref{lc:norm-normality} give
$r^{j,q}_k\downarrow0$. For the subword with profile label $q$,
\[
 0\le D_{Q,q}:=
 \sum_{r:\jmath(t_r)=q}L_{C_r}|h(t_r)-h_{\nu_r}(t_r)|
 \le L_Ae_{\varphi(q+1)}^q.
\]
Let $K$ be the largest active profile label. Solidity, subadditivity,
and \eqref{fv:data-cap} give
\[
 \rho_j\!\left(\sum_qD_{Q,q}\right)
 \le v_j\wedge\sum_{q=1}^{K}r^{j,q}_{\varphi(q+1)}.
\]
Apply bounded fusion to the regulators
$c^{j,q}_{ik}=r^{j,q}_k$. They satisfy
\[
 r^{j,q}_{\varphi(q+1)}
 \le H_{c^{j,q}}(i\mapsto\varphi(i+q)).
\]
One resulting regulator is
\[
 D^j_{1s}=0,\qquad
 D^j_{rs}=v_j\wedge2^r\sum_{1\le q<r}r^{j,q}_s
 \quad(r\ge2).
\]
It is bounded by $v_j$, and each row decreases to zero.
\end{proof}

\begin{lemma}[Transfer of a common input regulator]
\label{fv:envelope-common-data}
The conclusion of Lemma~\ref{fv:envelope-data} also holds if its
profile hypothesis is replaced by
\begin{equation}\label{fv:source-common-input}
 \exists a\ \forall\psi\ \forall t\ \exists n_0(t,\psi)
 \ \forall n\ge n_0(t,\psi):
 \quad |h_n(t)-h(t)|\le H_a(\psi),
\end{equation}
where $a$ is a regulator in $E$. The lower selector may then depend
on $j,\varphi$.
\end{lemma}
\begin{proof}
For $x\in E_+$, put $A_j(x)=\rho_j(L_Ax)$. Positivity, solidity,
and subadditivity give
\[
 A_j(0)=0,\quad x\le y\Longrightarrow A_j(x)\le A_j(y),
 \quad A_j(x+y)\le A_j(x)+A_j(y).
\]
Root normality and Proposition~\ref{lc:norm-normality} give
\[
 x_s\downarrow0\Longrightarrow L_Ax_s\downarrow0
 \Longrightarrow A_j(x_s)\downarrow0.
\]
Let $u\in E_+$ bound $a$. Lemma~\ref{lt:positive-d-transport}
therefore applies to $A_j$ and gives the explicit regulator
\[
 D^j_{rs}=\rho_j(L_Au)\wedge2^r\rho_j(L_Aa_{rs}),
 \qquad
 \rho_j(L_AH_a(\varphi))\le H_{D^j}(\varphi).
\]
Set
\[
 \ell_{j,\varphi}(t)=
 \min\{N\ge1:\ |h_n(t)-h(t)|\le H_a(\varphi)
                       \text{ for every }n\ge N\}.
\]
The set defining this minimum is nonempty by
\eqref{fv:source-common-input}. For a compatible partial record and
labels $\nu_r\ge\ell_{j,\varphi}(t_r)$,
\[
 0\le\sum_rL_{C_r}|h(t_r)-h_{\nu_r}(t_r)|
 \le\sum_rL_{C_r}H_a(\varphi)\le L_AH_a(\varphi).
\]
Consequently,
\[
 \rho_j\!\left(\sum_rL_{C_r}|h(t_r)-h_{\nu_r}(t_r)|\right)
 \le\rho_j(L_AH_a(\varphi))\le H_{D^j}(\varphi).
\]
This estimate holds on every compatible partial record. No envelope
bound and no gauge restriction are needed for this input alternative.
\end{proof}

\begin{theorem}[Envelope monotone convergence]
\label{fv:envelope-monotone}
Assume \eqref{fv:known-envelope}, and let $h_n\le h_{n+1}$.
Define $h(t)=\bigvee_nh_n(t)$. Suppose either the countable order
profiles of Lemma~\ref{fv:envelope-data} or the common input
regulator of Lemma~\ref{fv:envelope-common-data} is available.
Then $h$ is hereditarily integrable in the original convergence,
\begin{equation}\label{fv:envelope-integral-sup}
 J_h(C)=\bigvee_nJ_n(C)\qquad(C\in\mathscr C_A),
\end{equation}
and
\begin{equation}\label{fv:envelope-endpoints}
 0\le J_h(C)-J_n(C)
 \le J_h(A)-J_n(A)\downarrow0.
\end{equation}
For each $j$ there is a regulator $R^j$ such that
\[
 \forall\varphi\ \exists\Delta_{j,\varphi}\quad
 \forall C\in\mathscr C_A\quad
 \forall Q\in\mathscr R_C(\Delta_{j,\varphi}):
 \quad
 \rho_j(S_Qh-J_h(C))\le H_{R^j}(\varphi).
\]
Both the candidate-sum estimate and convergence of the primitives are
uniform over the compatible cells in the original uniform convergence.
\end{theorem}
\begin{proof}
The input envelope and Dedekind $\sigma$-completeness give
$-g(t)\le h_n(t)\le h(t)\le g(t)$. Positivity of known integrals
gives $J_n(C)\le J_{n+1}(C)$, and
\eqref{fv:primitive-envelope} defines
\[
 J(C)=\bigvee_nJ_n(C)\in[-G(C),G(C)].
\]
For a finite division $(C_r)_{r=1}^{q}$ of $C$,
\[
 J_n(C)=\sum_{r=1}^{q}J_n(C_r)
 \le\sum_{r=1}^{q}J(C_r).
\]
Conversely, for every finite tuple $(n_r)$, its maximum $N$ satisfies
\[
 \sum_{r=1}^{q}J_{n_r}(C_r)
 \le\sum_{r=1}^{q}J_N(C_r)=J_N(C)\le J(C).
\]
Successive countable suprema in finitely many variables give
$J(C)=\sum_rJ(C_r)$. Define
\[
 q_N=J(A)-J_N(A)\downarrow0.
\]
The additive charge $J-J_N$ is positive. A witnessing root division
for $C$ yields
\begin{equation}\label{fv:envelope-primitive-gap}
 0\le J(C)-J_N(C)\le q_N.
\end{equation}

Fix $j,\varphi$ and put
\[
 N=\varphi(1),\qquad
 p(t)=\max\{N,\ell_{j,\varphi}(t)\},\qquad
 \Delta_{j,\varphi}=\delta_{j,\varphi,p}.
\]
For the profile alternative, $\ell_{j,\varphi}=\ell_\varphi$.
Let $Q=((C_r,t_r))$ be a $\Delta_{j,\varphi}$-fine record on $C$.
Set
\[
 \nu_r=p(t_r),\qquad
 M_Q=\sum_rL_{C_r}h_{\nu_r}(t_r),\qquad
 P_Q=\sum_rJ_{\nu_r}(C_r).
\]
The data and mixed estimates give
\[
 \rho_j(S_Qh-M_Q)\le H_{D^j}(\varphi),\qquad
 \rho_j(M_Q-P_Q)\le H_{B^j}(\varphi).
\]
Monotonicity and finite additivity give
\[
 J_N(C)=\sum_rJ_N(C_r)
 \le P_Q\le\sum_rJ(C_r)=J(C).
\]
Equation~\eqref{fv:envelope-primitive-gap} and solidity imply
\[
 \rho_j(P_Q-J(C))\le\rho_j(q_{\varphi(1)}).
\]
Therefore
\begin{equation}\label{fv:envelope-three-errors}
 \rho_j(S_Qh-J(C))
 \le H_{D^j}(\varphi)+H_{B^j}(\varphi)
      +\rho_j(q_{\varphi(1)}).
\end{equation}
Proposition~\ref{lc:norm-normality} gives $\rho_j(q_s)\downarrow0$.
The matrix
\[
 R^j_{is}=3\bigl(D^j_{is}+B^j_{is}+\rho_j(q_s)\bigr)
\]
is an order-bounded regulator. Its diagonal supremum dominates the
right-hand side of \eqref{fv:envelope-three-errors}.
Define the decreasing pair families
\[
 \mathcal P_\delta=
 \{(S_Qh,J(C)):C\in\mathscr C_A,\ Q\in\mathscr R_C(\delta)\}.
\]
For every $j,\varphi$,
\[
 (x,y)\in\mathcal P_{\Delta_{j,\varphi}}
 \Longrightarrow\rho_j(x-y)\le H_{R^j}(\varphi).
\]
Thus every coordinate of the difference of
$\mathsf N(\mathcal P_\delta)$ is $D$-null. Pair faithfulness gives
\[
 \mathsf N(\mathcal P_\delta)
 \longrightarrow_{\mathcal U_\lambda}\Delta_Z.
\]
For fixed $C$, the pair net $(S_Qh,J(C))$ is a quasi-subnet of this
saturated net. Consequently $S_Qh\to_\lambda J(C)$. The original
response-tail correspondence and separation give
$J_h(C)=J(C)$, proving \eqref{fv:envelope-integral-sup}.

For endpoint convergence, define
\[
 \mathcal E_N=
 \{(J_n(C),J_h(C)):n\ge N,\ C\in\mathscr C_A\}.
\]
Equation~\eqref{fv:envelope-primitive-gap} gives
\[
 (x,y)\in\mathcal E_N\Longrightarrow
 \rho_j(x-y)\le\rho_j(q_N)\downarrow_N0.
\]
The regulator $c^j_{is}=\rho_j(q_s)$ satisfies
\[
 (x,y)\in\mathcal E_{\varphi(1)}\Longrightarrow
 \rho_j(x-y)\le H_{c^j}(\varphi).
\]
Pair faithfulness now yields
$\mathsf N(\mathcal E_N)\to_{\mathcal U_\lambda}\Delta_Z$.
This is uniform primitive convergence over all compatible cells.
Its positive bound is \eqref{fv:envelope-endpoints}.
\end{proof}

What makes the construction possible is the combination of one known
envelope tail with the signed partial Henstock estimates: together they
supply the mixed controller, while the original input profiles or the
original common regulator supply the data controller. No order-bounded
sum of a countable family of row errors or common row regulator has
been assumed, and finite-supremum stability and lattice homomorphism
properties of the coefficients have not entered the argument.

\subsection{Finite lattice suprema and the variation of local defects}
\label{fv:section}

To pass from signed partial errors to their variation, we must also
control the finite lattice suprema that arise when the signs are
separated. Retain the original separated locally solid group convergence
$\lambda$ on $Z$, its faithful solid lattice-valued gauges $\rho_j$,
and the local division calculus, with weakly $\sigma$-distributive
controller groups. The additional requirement is a property of the given
convergence on $Z$.

\begin{definition}[Finite-supremum stability]\label{fv:finite-sup}
For a nonempty $B\subseteq Z$, define
\[
 \mathfrak J(B)=
 \left\{z\in Z:
 |z|\le\bigvee_{r=1}^{m}|b_r|
 \text{ for some }m\ge1\text{ and }b_1,\ldots,b_m\in B\right\}.
\]
The original convergence is finite-supremum stable if
\[
 x_\alpha\longrightarrow_\lambda0
 \Longrightarrow
 \mathsf N\bigl(\mathfrak J(T_\alpha(x))\bigr)
 \longrightarrow_\lambda0
\]
for every admissible net.
\end{definition}

In this formulation the lattice supremum is continuous uniformly over
finite arities, a requirement placed on the fixed original convergence
and verified below either through a representation or directly in the
order models.

\begin{proposition}[A representation criterion]\label{fv:gauges}
Suppose the faithful solid gauges satisfy
\begin{equation}\label{fv:join-gauge}
 \rho_j(x\vee y)=\rho_j(x)\vee\rho_j(y)
 \qquad(x,y\in Z_+).
\end{equation}
Then $\lambda$ is finite-supremum stable.
\end{proposition}
\begin{proof}
Let $x_\alpha\to_\lambda0$. For each $j$, faithfulness gives a
regulator $a^j$ such that, for every $\varphi$, some $\alpha(j,\varphi)$
satisfies
\[
 \rho_j(x_\beta)\le H_{a^j}(\varphi)
 \qquad(\beta\ge\alpha(j,\varphi)).
\]
For $z\in\mathfrak J(T_{\alpha(j,\varphi)}(x))$, a finite family
$b_r\in T_{\alpha(j,\varphi)}(x)$ satisfies
\[
 \rho_j(z)
 \le\rho_j\!\left(\bigvee_{r=1}^{m}|b_r|\right)
 =\bigvee_{r=1}^{m}\rho_j(b_r)
 \le H_{a^j}(\varphi).
\]
Thus the same regulator controls the saturated finite-supremum net.
Inverse faithfulness proves its original null convergence.
\end{proof}

\begin{proposition}[Intrinsic order models]\label{fv:order-models}
Original order convergence and original regulator convergence are
finite-supremum stable. Original unbounded order convergence has the
same property.
\end{proposition}
\begin{proof}
For $u\ge0$,
\[
 B\subseteq[-u,u]\quad\Longrightarrow\quad
 \mathfrak J(B)\subseteq[-u,u].
\]
Every attained order bound for a null net remains an attained bound
for its saturated finite-supremum net. Every attained regulator bound
$H_a(\varphi)$ remains unchanged. These observations prove the first
two assertions directly from their definitions.

For unbounded order convergence, fix $u\in Z_+$. If a tail set $B$
satisfies $|b|\wedge u\le v$ for every $b\in B$, then
\[
 z\in\mathfrak J(B)
 \Longrightarrow
 |z|\wedge u
 \le\left(\bigvee_{r=1}^{m}|b_r|\right)\wedge u
 =\bigvee_{r=1}^{m}(|b_r|\wedge u)\le v.
\]
The same decreasing order bounds apply for every $u$.
\end{proof}

\begin{lemma}[The family of signed partial defects]\label{fv:signed-defects}
Let $h$ have known hereditary original integrals $J_h(C)$. For every gauge $\delta$, let $\mathscr Q_A(\delta)$ denote its
compatible fine partial records. Define
\[
 \mathscr P_\delta(h)=\{0\}\cup
 \left\{\sum_{(C,t)\in Q}[L_Ch(t)-J_h(C)]:
 Q\in\mathscr Q_A(\delta)\right\}.
\]
The saturated net $\mathsf N(\mathscr P_\delta(h))$ converges to zero in $\lambda$.
\end{lemma}
\begin{proof}
Fix $j$. The known root integral supplies a regulator $a^j$ and,
for every $\varphi$, a gauge $\eta_{j,\varphi}$ satisfying
\[
 \rho_j(S_Rh-J_h(A))\le H_{a^j}(\varphi)
 \qquad(R\in\mathscr R_A(\eta_{j,\varphi})).
\]
Lemma~\ref{lc:partial} gives
\[
 \rho_j(z)\le H_{a^j}(\varphi)
 \qquad(z\in\mathscr P_{\eta_{j,\varphi}}(h)).
\]
Hence $\rho_j(\mathsf N(\mathscr P_\delta(h)))$ is $D$-null.
Inverse faithfulness gives the original null convergence.
\end{proof}

\begin{theorem}[Variation obtained from signed Henstock control]
\label{fv:variation-henstock}
Suppose the original convergence is finite-supremum stable. For every
known hereditarily integrable $h$, the sets
\[
 \mathscr V_\delta(h)=\{0\}\cup
 \left\{\sum_{(C,t)\in Q}|L_Ch(t)-J_h(C)|:
 Q\in\mathscr Q_A(\delta)\right\}
\]
satisfy $\mathsf N(\mathscr V_\delta(h))\to_\lambda0$.
\end{theorem}
\begin{proof}
Put $\mathscr N_\delta=\mathscr P_\delta(h)-\mathscr P_\delta(h)$.
Finite gauge meets make the diagonal family $\mathscr N_\delta$
cofinal in the differences of pairs of tail sets of the saturated nets. Group continuity
and Lemma~\ref{fv:signed-defects} give
\[
 \mathsf N(\mathscr N_\delta)\to_\lambda0,
 \qquad
 \mathsf N(\mathfrak J(\mathscr N_\delta))\to_\lambda0.
\]
For a compatible partial record $Q=((C_r,t_r))_{r=1}^{m}$, put
$e_r=L_{C_r}h(t_r)-J_h(C_r)$. Finite lattice distributivity gives
\[
 \sum_{r=1}^{m}|e_r|
 =\sum_{r=1}^{m}(e_r\vee(-e_r))
 =\bigvee_{\varepsilon\in\{-1,1\}^{m}}
   \sum_{r=1}^{m}\varepsilon_re_r.
\]
For each sign vector, both subwords $Q_+$ and $Q_-$ remain compatible
and $\delta$-fine. Thus
\[
 \sum_r\varepsilon_re_r
 =\sum_{r:\varepsilon_r=1}e_r-
  \sum_{r:\varepsilon_r=-1}e_r
 \in\mathscr N_\delta.
\]
Consequently $\mathscr V_\delta(h)\subseteq
\mathfrak J(\mathscr N_\delta)$. The saturated net of variation defects is a quasi-subnet of the
latter null net. Quasi-subnet stability gives its original null convergence.
\end{proof}

\subsection{Band covariance and finite lattice operations}
\label{fv:band-section}

To recognize when a representation has the finite-lattice properties
used above, we examine its compatibility with band projections. Let $Z$
and $F$ be abelian lattice-ordered groups; a solid symmetric lattice
pseudonorm is a map $\rho:Z\to F_+$ satisfying
\[
 \rho(0)=0,\quad \rho(-z)=\rho(z),\quad
 \rho(z+w)\le\rho(z)+\rho(w),\quad
 |z|\le|w|\Longrightarrow\rho(z)\le\rho(w).
\]
These conditions imply $\rho(z)=\rho(|z|)$, without requiring positive
homogeneity.

\begin{definition}[Band covariance]\label{fv:covariance}
Let $\mathscr B_Z$ and $\mathscr B_F$ be Boolean algebras of band
projections, including the identity and zero projections. Suppose
$\mathscr B_Z$ contains the principal projections needed below.
A solid symmetric lattice pseudonorm is band covariant if there is
a unital Boolean homomorphism
$\Theta:\mathscr B_Z\to\mathscr B_F$ such that
\[
 \rho(Pz)=\Theta(P)\rho(z)
 \qquad(P\in\mathscr B_Z,\ z\in Z).
\]
\end{definition}

\begin{proposition}[Domination implies covariance]
\label{fv:domination}
In Definition~\ref{fv:covariance}, the weaker inequalities
\[
 \rho(Pz)\le\Theta(P)\rho(z)
 \qquad(P\in\mathscr B_Z,\ z\in Z)
\]
already imply all the required equalities.
\end{proposition}
\begin{proof}
Write $Q=\Theta(P)$. Then
\[
 \rho(z)\le\rho(Pz)+\rho((I-P)z)
 \le Q\rho(z)+(I-Q)\rho(z)=\rho(z).
\]
The two positive differences
$Q\rho(z)-\rho(Pz)$ and $(I-Q)\rho(z)-\rho((I-P)z)$ sum to zero.
Each difference is therefore zero.
\end{proof}

\begin{proposition}[Finite-join preservation]
\label{fv:joins}\label{fv:band-covariance}
Every band-covariant solid symmetric lattice pseudonorm satisfies
\[
 \rho(x\vee y)=\rho(x)\vee\rho(y)\qquad(x,y\in Z_+).
\]
\end{proposition}
\begin{proof}
Let $P=P_{(x-y)^+}$ and $Q=\Theta(P)$. The band decomposition gives
\[
 Px\ge Py,\qquad (I-P)x\le(I-P)y,
 \qquad x\vee y=Px+(I-P)y.
\]
Solidity and covariance imply
\[
 Q\rho(x)\ge Q\rho(y),\qquad
 (I-Q)\rho(x)\le(I-Q)\rho(y).
\]
Covariance applied to $x\vee y$ gives
\[
 Q\rho(x\vee y)=Q\rho(x),\qquad
 (I-Q)\rho(x\vee y)=(I-Q)\rho(y).
\]
Consequently,
\[
 \rho(x\vee y)
 =Q\rho(x)+(I-Q)\rho(y)
 =\rho(x)\vee\rho(y).
\]
\end{proof}

\begin{proposition}[Variation from signed partial sums]
\label{fv:variation}
Suppose $\rho$ preserves finite joins on $Z_+$. For
$e_1,\ldots,e_N\in Z$ and $a\in F_+$, assume
\[
 \rho\!\left(\sum_{r\in S}e_r\right)\le a
 \qquad(S\subseteq\{1,\ldots,N\}).
\]
Then
\[
 \rho\!\left(\sum_{r=1}^{N}|e_r|\right)\le2a.
\]
\end{proposition}
\begin{proof}
Finite distributivity of addition over joins gives
\[
 \sum_{r=1}^{N}e_r^+
 =\bigvee_{S\subseteq\{1,\ldots,N\}}\sum_{r\in S}e_r
 =\bigvee_{S\subseteq\{1,\ldots,N\}}
 \left(\sum_{r\in S}e_r\right)^+.
\]
Hence
\[
 \rho\!\left(\sum_re_r^+\right)
 =\bigvee_S\rho\!\left(\left(\sum_{r\in S}e_r\right)^+\right)
 \le\bigvee_S\rho\!\left(\sum_{r\in S}e_r\right)\le a.
\]
Replacing each $e_r$ by $-e_r$ proves
$\rho(\sum_re_r^-)\le a$. Subadditivity now gives
\[
 \rho\!\left(\sum_r|e_r|\right)
 =\rho\!\left(\sum_re_r^++\sum_re_r^-\right)\le2a.
\]
\end{proof}

\begin{proposition}[Purely ordered examples of compatible pseudonorms]
\label{fv:examples}
Let $T:Z\to F$ be a positive additive lattice homomorphism and $u\in F_+$.
The maps
\[
 \rho_T(z)=T|z|,\qquad \rho_{T,u}(z)=T|z|\wedge u
\]
are solid symmetric lattice pseudonorms preserving finite joins on
$Z_+$. For $T=I$, truncated modulus is band covariant with $\Theta=I$.
\end{proposition}
\begin{proof}
The inequalities
\[
 T|z+w|\le T|z|+T|w|,\qquad
 (a+b)\wedge u\le(a\wedge u)+(b\wedge u)
 \quad(a,b,u\ge0)
\]
prove subadditivity. Symmetry and solidity are preserved by $T$ and by
the truncation map. For $x,y\ge0$,
\[
 T(x\vee y)=Tx\vee Ty,\qquad
 (Tx\vee Ty)\wedge u=(Tx\wedge u)\vee(Ty\wedge u).
\]
For a band projection $P$ on $Z$,
\[
 |Pz|\wedge u=P|z|\wedge Pu=P(|z|\wedge u).
\]
This proves the last assertion. These are algebraic certificates.
Faithfulness of a selected family is verified separately.
\end{proof}

\begin{proposition}[A faithful unbounded-order realization]
\label{fv:uo}
Let $Z$ be a Dedekind $\sigma$-complete WSD abelian lattice-ordered
group with the countable sup property on principal order intervals.
Define its original convergence directly on nets by
\[
 z_\alpha\longrightarrow_{\mathrm{uo}}z
 \quad\Longleftrightarrow\quad
 |z_\alpha-z|\wedge u\longrightarrow_o0
 \quad(u\in Z_+).
\]
Pair nets converge to the diagonal when their differences converge to
zero in this sense. Then $\rho_u(z)=|z|\wedge u$, $u\in Z_+$, is a faithful regulator
representation on point and pair nets. This convergence is separated,
locally solid, and sequentially order compatible. It has intrinsic finite-supremum
stability. If $Z$ also satisfies lateral $\sigma$-completeness, its controllers satisfy
Definition~\ref{lc:calibration}.
\end{proposition}
\begin{proof}
Proposition~\ref{fd:intrinsic-order-realizations} gives the separated
locally solid group convergence and its faithful order-valued point
and pair realization. Proposition~\ref{lt:order-d} applies because
$Z$ is Dedekind $\sigma$-complete, WSD, and has the countable sup
property on each principal interval. For every admissible net and
pair net, respectively,
\[
 z_\alpha\to_{\mathrm{uo}}z
 \iff\rho_u(z_\alpha-z)\to_o0\ (u)
 \iff\rho_u(z_\alpha-z)\to_D0\ (u),
\]
\[
 (z_\alpha,w_\alpha)\to_{\mathcal U_{\mathrm{uo}}}\Delta_Z
 \iff\rho_u(z_\alpha-w_\alpha)\to_D0\quad(u\in Z_+).
\]
For sequential order compatibility and finite-supremum stability,
\[
 z_s\downarrow0\Longrightarrow z_s\wedge u\downarrow0\quad(u),
\]
\[
 x,y\ge0\Longrightarrow
 \rho_u(x\vee y)=(x\vee y)\wedge u
 =(x\wedge u)\vee(y\wedge u)
 =\rho_u(x)\vee\rho_u(y).
\]
Proposition~\ref{fv:gauges} gives finite-supremum stability. Saturated
solid-tail convergence is Lemma~\ref{fd:solid-tail}.
Finally, Proposition~\ref{lc:component-covers} gives the component-cover
condition from the principal-interval countable sup property. The
additional lateral $\sigma$-completeness supplies the remaining
condition of Definition~\ref{lc:calibration}.
\end{proof}

Through these verifications, algebraic properties of the representing
pseudonorms yield the corresponding intrinsic convergence property.
Finite-join preservation and band covariance must nevertheless be checked
for the representation in use, since an arbitrary metrizing family is
not asserted to possess them.

\subsection{Resolution of finite divisions}

Alongside the convergence of local errors, we need a way to resolve
a prescribed division without changing the finite response. A division
$(E_s)_{s\in S}$ refines $(D_r)_{r\in R}$ when there is an occurrence
partition $S=\coprod_{r\in R}S_r$ such that
$(E_s)_{s\in S_r}\in\mathfrak D(D_r)$ for every $r$.

\begin{definition}[Resolution of a finite division]\label{dt:resolution}
For every cell $C$ and $\mathcal D\in\mathfrak D(C)$, require a
gauge $\eta_{\mathcal D}$ such that, for every gauge $\delta$ and
every $R\in\mathscr R_C(\delta\wedge\eta_{\mathcal D})$, there is
$R'\in\mathscr R_C(\delta)$ satisfying:
\begin{enumerate}
\item its cell division refines $\mathcal D$ and the cell division of
$R$, and
\item every cell of $R$ is replaced by a division of that cell, with
its original tag retained on each replacement cell.
\end{enumerate}
\end{definition}
Coefficient additivity turns this condition on divisions and tags into
the response identity
\begin{equation}\label{dt:sum-preserved}
 S_{R'}u=S_Ru\qquad\text{for every function }u:T\to E.
\end{equation}

For a finite division $\mathcal D$ of $C$, write
\[
 \mathscr R_C(\delta,\mathcal D)
 =\{R\in\mathscr R_C(\delta):
       \operatorname{div}R\succeq\mathcal D\}.
\]

\begin{proposition}[Equality of sum tails]\label{cr:response-tails}
Suppose finite-division resolution, Definition~\ref{dt:resolution},
holds and the coefficients are additive on divisions. For every
function $h$, every division $\mathcal D$, and every gauge $\delta$,
\begin{equation}\label{cr:tail-sandwich}
 \{S_Rh:R\in\mathscr R_C(\delta\wedge\eta_{\mathcal D})\}
 \subseteq
 \{S_Rh:R\in\mathscr R_C(\delta,\mathcal D)\}
 \subseteq
 \{S_Rh:R\in\mathscr R_C(\delta)\}.
\end{equation}
Consequently, the full gauge net and the net imposing both gauge and
division refinement have identical sum-tail families up to cofinal
inclusion. They define the same integral in every faithful order or
regulator representation considered here.
\end{proposition}
\begin{proof}
For $R\in\mathscr R_C(\delta\wedge\eta_{\mathcal D})$, finite
resolution gives $R'\in\mathscr R_C(\delta,\mathcal D)$. Every
parent occurrence $(B,t)$ is replaced by a division $(E_s,t)$, and
\[
 \sum_sL_{E_s}h(t)=L_Bh(t).
\]
Summation over the parent occurrences gives $S_{R'}h=S_Rh$.
This proves the first inclusion in \eqref{cr:tail-sandwich}.
The second inclusion holds by the definition of
$\mathscr R_C(\delta,\mathcal D)$.

Write
\[
 V_\delta=\{S_Rh:R\in\mathscr R_C(\delta)\},\qquad
 V_{\delta,\mathcal D}
 =\{S_Rh:R\in\mathscr R_C(\delta,\mathcal D)\}.
\]
The first family gives the value tails of the full gauge net. The
refinement net described below has value tails
$V_{\delta,\operatorname{div}R}$. Each of its tails contains the
full-gauge tail $V_{\delta\wedge\eta_{\operatorname{div}R}}$ by
\eqref{cr:tail-sandwich}. Conversely, for each gauge $\delta$,
fine-record existence supplies an $R\in\mathscr R_C(\delta)$,
and then
\[
 V_{\delta,\operatorname{div}R}\subseteq V_\delta.
\]
Every tail of either response net therefore contains a tail of the
other response net. They are tail equivalent and have exactly the
same limits for the prescribed convergence. This proves the assertion
for each of its faithful order or regulator representations as well.
The equality $S_{R'}h=S_Rh$ holds for every $h$ under the same
refinement, since the coefficient identities precede evaluation at $h$.
\end{proof}

To make the refinement order explicit, index its net by
$(\delta,R,k)$, where $R\in\mathscr R_C(\delta)$ and
$k\in\mathbb N$. Equality and the strict relation
\[
 (\delta,R,k)<(\epsilon,S,l)
 \iff \epsilon\le\delta,\quad
 \operatorname{div}S\succeq\operatorname{div}R,\quad l>k
\]
define a directed order. For two indices with divisions
$\mathcal D_1,\mathcal D_2$, let $\delta$ be the meet of their gauges.
Fine-record existence supplies a record fine for
$\delta\wedge\eta_{\mathcal D_1}\wedge\eta_{\mathcal D_2}$.
Resolution at $\mathcal D_1$ preserves
$\delta\wedge\eta_{\mathcal D_2}$-fineness. Subsequent resolution
at $\mathcal D_2$ preserves $\delta$-fineness. The final division
refines both prescribed divisions. Increasing the clock gives an
upper bound. The tail at $(\delta,R,k)$ has record image exactly
$\mathscr R_C(\delta,\operatorname{div}R)$.

With this ordering convention, later indices in the gauge/refinement
example have smaller gauges and finer cell divisions, and
\eqref{cr:tail-sandwich} supplies the required equality of sum tails
without requiring equality of the record tails with those of the pure
gauge net.

\begin{definition}[Countable division resolution]\label{cr:geometry}
For each cell $C$, let $(\mathcal P_n(C))_{n\ge1}$ be a refining
sequence of finite divisions. Require the following two properties.
\begin{enumerate}
\item For every $n$, the division $\mathcal P_n(C)$ admits the
tag-preserving resolution of Definition~\ref{dt:resolution}.
\item For every finite division $\mathcal D$ of $C$ and every gauge
$\delta$, there exist an index $n$, a subword
$\mathcal B\subseteq\mathcal P_n(C)$, divisions
$\mathcal E_B\in\mathfrak D(B)$ for $B\in\mathcal B$, and tags
$t_B$, such that
\[
 \mathcal E=
 (\mathcal P_n(C)\setminus\mathcal B)
 \mathbin{\sqcup}\mathop{\bigsqcup}_{B\in\mathcal B}\mathcal E_B
\]
refines $\mathcal D$, and both partial records
\[
 Q=\bigl((B,t_B)\bigr)_{B\in\mathcal B},\qquad
 Q'=\bigl((E,t_B)\bigr)_{B\in\mathcal B,\ E\in\mathcal E_B}
\]
are $\delta$-fine. The tag $t_B$ belongs to the carrier of $B$ and
to the carrier of every child in $\mathcal E_B$.
\end{enumerate}
The empty subword is permitted. All operations act on occurrence
indices. Neither property involves functions, coefficients, or
integrals.
\end{definition}

\begin{proposition}[Dyadic resolution of interval divisions]
\label{cr:intervals}
Let cells be nondegenerate interval pieces, with closed interval tag
carriers and an additive endpoint convention. Divisions are finite
subdivisions at interior cuts. A record is $\delta$-fine when
\[
 \overline I\subset(t-\delta(t),t+\delta(t)),\qquad t\in\overline I,
\]
where $\delta(t)>0$ is arbitrary. These divisions satisfy finite
resolution and countable division resolution.
\end{proposition}
\begin{proof}
Fix $C$ with carrier $[c,d]$. Let $\mathcal P_n(C)$ be its dyadic
division with mesh $m_n=(d-c)2^{-n}$. Finite common refinements are
obtained by taking unions of finite cut sets.

Fine-record existence is verified by the dyadic branching criterion:
every infinite branch has nested compact carriers of diameters tending
to zero, with common point $t$. Eventually its carrier is contained in
$(t-\delta(t),t+\delta(t))$. Proposition~\ref{ls:branch} supplies a
finite fine division.

For a fixed finite cut set $B$, let
\[
 \eta_B(t)=\tfrac12\min\{|t-b|:b\in B,\ b\ne t\},
\]
with a fixed positive value when the set is empty. For a tagged cell
$(I,t)$ fine for $\delta\wedge\eta_B$ and
$b\in B\cap\operatorname{int}I$ with $b\ne t$,
\[
 |b-t|<\eta_B(t)\le\tfrac12|b-t|,
\]
a contradiction. Thus
\[
 B\cap\operatorname{int}I\subseteq\{t\}.
\]
When $t\in B\cap\operatorname{int}I$, let $I_-,I_+$ be the two
children obtained by subdivision at $t$. Then
\[
 t\in\overline{I_-}\cap\overline{I_+},\qquad
 \overline{I_\pm}\subseteq\overline I
 \subset(t-\delta(t),t+\delta(t)).
\]
When this intersection is empty, retain $I$. Applying this operation
to every occurrence gives a $\delta$-fine record whose division
contains every cut in $B$ and whose children retain the parent tags.
This verifies Definition~\ref{dt:resolution}.

For Definition~\ref{cr:geometry}(2), fix a division $\mathcal D$,
let $B$ be its finite set of interior cuts, and fix $\delta$.
If $B=\varnothing$, the unchanged dyadic division suffices. Otherwise,
for all sufficiently large $n$,
\[
 m_n<\min_{t\in B}\delta(t),
\]
and, when $B$ contains at least two points,
\[
 m_n<\min\{|s-t|:s,t\in B,\ s\ne t\}.
\]
Thus each dyadic cell contains at most one cut in its interior.
Let $\mathcal B$ consist of the cells containing such a cut. For
$Q\in\mathcal B$, denote its unique interior cut by $t_Q$.
Subdivide $Q$ there into two nondegenerate children and tag the parent
and both children at $t_Q$.
For every $x\in\overline Q$,
\[
 |x-t_Q|\le m_n<\delta(t_Q).
\]
The parent and both children are $\delta$-fine. The parent word is a
subword of $\mathcal P_n(C)$, while the child word is a subword of the
division with the union of the two cut sets. That division refines
$\mathcal D$. Cuts already at dyadic endpoints cause no subdivision.
Both required words are therefore compatible and fine.
\end{proof}

Only the two boundary words are subject to the fineness requirement
used here, so the argument imposes no uniform positive lower bound on
an arbitrary gauge over the entire interval. Moreover, cancellation at a
common tag remains exact in the presence of atoms, provided the cell
coefficients obey their additive endpoint convention.

\subsection{Primitive variation in the original convergence}

We can now bring division resolution and partial Henstock control
together to construct the variation of a primitive. Retain the rooted
division calculus, its original record tails, and the signed partial
Henstock theorem, and let $E$ and $Z$ be abelian lattice-ordered groups,
with $Z$ Dedekind $\sigma$-complete. Its separated additive convergence
$\lambda$ is assumed locally solid, sequentially order compatible, and
finite-supremum stable in the sense of Definition~\ref{fv:finite-sup},
while the given convergence and uniform convergence retain their faithful
solid lattice representation $\rho_j:Z\to(F_j)_+$ from the preceding
section. Neither additivity nor finite-supremum preservation of the
individual maps $\rho_j$ is required here.

For every cell $C$, let $L_C:E\to Z$ be positive and additive, with
\begin{equation}\label{fv:geometry-coefficients}
 L_C|x|=|L_Cx|,
 \qquad
 L_C=\sum_{D\in\mathcal D}L_D
 \quad(\mathcal D\in\mathfrak D(C)).
\end{equation}
Assume one of the following structural conditions:
\begin{enumerate}
\item countable division resolution, Definition~\ref{cr:geometry}, holds
on every cell, or
\item finite-division resolution, Definition~\ref{dt:resolution}, holds
on every cell, and $Z$ is Dedekind complete with the countable sup property
on each principal order interval.
\end{enumerate}
Under the first alternative, the division sequence is geometric and
fixed in advance; under the second, we shall construct it from the known
primitive.

For $\diamond\in\{+,-,\mathrm{abs}\}$, write
\[
 x^+=x\vee0,\qquad x^-=(-x)\vee0,
 \qquad x^{\mathrm{abs}}=|x|.
\]
The lattice identities needed below are
\begin{equation}\label{fv:geometry-lattice-identities}
 0\le x^\diamond\le|x|,
 \quad
 (x+y)^\diamond\le x^\diamond+y^\diamond,
 \quad
 |x^\diamond-y^\diamond|\le|x-y|.
\end{equation}
For positive parts, the last inequality follows from the two bounds
\[
 x^+\le y^++(x-y)^+\le y^++|x-y|,
 \qquad
 y^+\le x^++|x-y|.
\]
Negative parts give the same bounds after replacing $x,y$ by $-x,-y$.
For moduli, the two lattice triangle inequalities give
$|x|\le|y|+|x-y|$ and $|y|\le|x|+|x-y|$.
The positive-part subadditivity follows from
$x+y\le x^++y^+$ and $0\le x^++y^+$, and negative parts and moduli have
the corresponding inequalities.

Equation~\eqref{fv:geometry-coefficients} gives all three identities
\begin{equation}\label{fv:geometry-polarity-transport}
 L_C(x^\diamond)=(L_Cx)^\diamond.
\end{equation}
Indeed, the group identity $2x^+=x+|x|$ yields
\[
 2L_C(x^+)=L_Cx+|L_Cx|=2(L_Cx)^+.
\]
Doubling is injective in an abelian lattice-ordered group. Hence
$L_C(x^+)=(L_Cx)^+$. Applying this equality to $-x$ proves the
negative-part identity. No division in either group is used.

\begin{theorem}[Countable realization and lattice operations]
\label{fv:geometry-variation}
\label{fv:finite-lattice-closure}\label{cr:lattice-closure}\label{dt:lattice-closure}
Let $h$ and $g\ge0$ be hereditarily integrable on $\mathscr C_A$,
with $|h|\le g$. Under the preceding hypotheses, for every
$C\in\mathscr C_A$ and $\diamond\in\{+,-,\mathrm{abs}\}$,
\begin{equation}\label{fv:geometry-variation-formula}
 V_h^\diamond(C)
 :=\sup_{n\ge1}
 \sum_{B\in\mathcal P_n(C)}J_h(B)^\diamond
 =\sup_{\mathcal D\in\mathfrak D(C)}
 \sum_{D\in\mathcal D}J_h(D)^\diamond
 \in Z.
\end{equation}
In the second structural alternative, the refining sequence
$\mathcal P_n(C)$ in \eqref{fv:geometry-variation-formula} may depend on
$h$ and $\diamond$. The functions $h^+$, $h^-$, and $|h|$ are hereditarily integrable in
the original convergence, and
\begin{equation}\label{fv:geometry-polarity-integrals}
 J_{h^\diamond}(C)=V_h^\diamond(C).
\end{equation}
Consequently, the hereditarily integrable functions dominated by
integer multiples of $g$ are closed under finite lattice operations.
\end{theorem}

\begin{proof}
Write $G=J_g$. Positivity of the prescribed finite sums and local
solidity preserve inequalities between existing integrals. In
particular,
\[
 -G(C)\le J_h(C)\le G(C),
 \qquad |J_h(C)|\le G(C).
\]
For clarity, if two integrable functions $r\le s$ are given, every
record satisfies $S_Rr\le S_Rs$. Thus
\[
 (J_r(C)-J_s(C))^+
 \le |J_r(C)-S_Rr|+|S_Rs-J_s(C)|\longrightarrow_\lambda0.
\]
Solidity makes the constant on the left $\lambda$-null. Separation
makes it zero. This proves the order passage just used.

Fix $C$ and $\diamond$. For a finite division $\mathcal D$, set
\[
 v(\mathcal D)=\sum_{D\in\mathcal D}J_h(D)^\diamond.
\]
Finite additivity of the known primitives gives
\begin{equation}\label{fv:geometry-variation-bound}
 0\le v(\mathcal D)
 \le\sum_{D\in\mathcal D}|J_h(D)|
 \le\sum_{D\in\mathcal D}G(D)=G(C).
\end{equation}
If $\mathcal E$ refines $\mathcal D$, then
\[
 \begin{split}
 v(\mathcal D)
 &=\sum_{D\in\mathcal D}
 \left(\sum_{\substack{E\in\mathcal E\\E\text{ in }D}}
 J_h(E)\right)^\diamond\\
 &\le\sum_{E\in\mathcal E}J_h(E)^\diamond
 =v(\mathcal E).
 \end{split}
\]
In the second structural alternative, the nonempty set
$\{v(\mathcal D):\mathcal D\in\mathfrak D(C)\}\subseteq[0,G(C)]$
has a supremum $V$. The countable sup property supplies divisions
$\mathcal D_k$ satisfying
\[
 V=\sup_kv(\mathcal D_k).
\]
Finite-division resolution gives common finite refinements: a record
fine for $\eta_{\mathcal D}\wedge\eta_{\mathcal E}$ resolves first
against $\mathcal D$ and then against $\mathcal E$, retaining fineness
for the second resolution. Finite recursion therefore gives refining
divisions $\mathcal P_n(C)$ with
\[
 \mathcal P_n(C)\succeq\mathcal D_k\quad(k\le n).
\]
Consequently,
\[
 \sup_{k\le n}v(\mathcal D_k)
 \le v(\mathcal P_n(C))\le V,
 \qquad \sup_nv(\mathcal P_n(C))=V.
\]
Each $\mathcal P_n(C)$ admits tag-preserving resolution by the same
finite-division hypothesis.

In the first alternative, use the prescribed geometric sequence.
In either alternative put $v_n=v(\mathcal P_n(C))$.
Therefore $v_n\uparrow$ and $v_n\le G(C)$. Dedekind
$\sigma$-completeness constructs
\begin{equation}\label{fv:geometry-countable-candidate}
 V=\sup_nv_n\in Z,
 \qquad q_n=V-v_n\downarrow0.
\end{equation}
For the first alternative, this construction has not used a supremum
over all finite divisions.

For a compatible partial record $Q$, define its positive defect
variation by
\[
 \mathcal V_h(Q)
 =\sum_{(B,t)\in Q}|J_h(B)-L_Bh(t)|.
\]
Let
\[
 \mathcal W_\delta(C)
 =\{0\}\cup
 \{\mathcal V_h(Q):Q\text{ is a compatible }
       \delta\text{-fine partial record in }C\}.
\]
Theorem~\ref{fv:variation-henstock} proves, in the original
convergence,
\begin{equation}\label{fv:geometry-variation-null}
 \mathsf N(\mathcal W_\delta(C))\longrightarrow_\lambda0
 \qquad(\delta\text{ decreases}).
\end{equation}
This family is obtained from the known integral of $h$ by the signed
Henstock estimate and finite-supremum stability of $\lambda$.

In the first structural alternative, let $\mathcal D$ be an arbitrary
finite division of $C$. For every gauge $\delta$,
Definition~\ref{cr:geometry} supplies an index $n$,
a subword $\mathcal B\subseteq\mathcal P_n(C)$, child divisions
$\mathcal E_B$, and common tags $t_B$, such that
\[
 \mathcal E=
 (\mathcal P_n(C)\setminus\mathcal B)
 \sqcup\bigsqcup_{B\in\mathcal B}\mathcal E_B
 \succeq\mathcal D,
\]
and both compatible partial records
\[
 Q=((B,t_B))_{B\in\mathcal B},
 \qquad
 Q'=((E,t_B))_{B\in\mathcal B,\ E\in\mathcal E_B}
\]
are $\delta$-fine. For each replaced parent, coefficient additivity
and \eqref{fv:geometry-polarity-transport} give the exact identity
\begin{equation}\label{fv:geometry-same-tag-cancellation}
 \sum_{E\in\mathcal E_B}(L_Eh(t_B))^\diamond
 =\sum_{E\in\mathcal E_B}L_E(h(t_B)^\diamond)
 =L_B(h(t_B)^\diamond)
 =(L_Bh(t_B))^\diamond.
\end{equation}
Finite additivity and
\eqref{fv:geometry-lattice-identities} give
\[
 0\le
 \sum_{E\in\mathcal E_B}J_h(E)^\diamond-J_h(B)^\diamond.
\]
Subtracting and adding the two equal quantities in
\eqref{fv:geometry-same-tag-cancellation} gives
\[
 \begin{split}
 &\sum_{E\in\mathcal E_B}J_h(E)^\diamond-J_h(B)^\diamond\\
 &\quad\le
 \sum_{E\in\mathcal E_B}|J_h(E)-L_Eh(t_B)|
 +|L_Bh(t_B)-J_h(B)|.
 \end{split}
\]
The unchanged occurrences cancel in $v(\mathcal E)-v_n$. Hence
\begin{equation}\label{fv:geometry-variation-gain}
 0\le v(\mathcal E)-v_n
 \le\mathcal V_h(Q')+\mathcal V_h(Q).
\end{equation}
Since $v(\mathcal D)\le v(\mathcal E)$ and $v_n\le V$,
\begin{equation}\label{fv:geometry-fixed-positive-bound}
 0\le (v(\mathcal D)-V)^+
 \le\mathcal V_h(Q')+\mathcal V_h(Q).
\end{equation}

For any two gauges,
\[
 \mathcal W_{\delta_1\wedge\delta_2}(C)
 +\mathcal W_{\delta_1\wedge\delta_2}(C)
 \subseteq\mathcal W_{\delta_1}(C)+\mathcal W_{\delta_2}(C).
\]
Equation~\eqref{fv:geometry-variation-null}, product convergence,
addition, and local solidity give
\[
 \mathsf N\!\left(\operatorname{sol}
 [\mathcal W_\delta(C)+\mathcal W_\delta(C)]\right)
 \longrightarrow_\lambda0.
\]
For the fixed $d=(v(\mathcal D)-V)^+$,
\eqref{fv:geometry-fixed-positive-bound} gives
\[
 d\in\operatorname{sol}
 [\mathcal W_\delta(C)+\mathcal W_\delta(C)]\quad(\delta).
\]
The constant net at $d$ is a quasi-subnet of this saturated null net.
It also converges to $d$. Separation yields $d=0$. Thus
\[
 v(\mathcal D)\le V
 \qquad(\mathcal D\in\mathfrak D(C)).
\]
Conversely, every upper bound for all $v(\mathcal D)$ bounds every
$v_n$, and therefore bounds $V$. This proves
\eqref{fv:geometry-variation-formula}.

In both alternatives, let $\eta_n$ resolve $\mathcal P_n(C)$. If $R$ is fine for
$\delta\wedge\eta_n$, its tag-preserving resolution $R'$ satisfies
\[
 R'\text{ is }\delta\text{-fine},
 \qquad \operatorname{div}R'\succeq\mathcal P_n(C),
 \qquad S_R(h^\diamond)=S_{R'}(h^\diamond).
\]
The last identity uses coefficient additivity on each replaced cell,
with its tag retained. The proved supremum identity gives
\[
 v_n\le v(\operatorname{div}R')\le V.
\]
Equation~\eqref{fv:geometry-polarity-transport} and the lattice
Lipschitz inequality give
\[
 \begin{split}
 |S_{R'}(h^\diamond)-v(\operatorname{div}R')|
 &=\left|\sum_{(B,t)\in R'}
       [(L_Bh(t))^\diamond-J_h(B)^\diamond]\right|\\
 &\le\sum_{(B,t)\in R'}|L_Bh(t)-J_h(B)|
 =\mathcal V_h(R').
 \end{split}
\]
Consequently,
\begin{equation}\label{fv:geometry-candidate-error}
 |S_R(h^\diamond)-V|\le\mathcal V_h(R')+q_n.
\end{equation}
For each $j$, \eqref{fv:geometry-variation-null} supplies a regulator
$a^j$ and gauges $\delta(j,\varphi)$ with
\[
 w\in\mathcal W_{\delta(j,\varphi)}(C)
 \Longrightarrow\rho_j(w)\le H_{a^j}(\varphi).
\]
The sequence $\rho_j(q_s)$ decreases to zero by
Proposition~\ref{lc:norm-normality}. Set
\[
 b^j_{is}=2\bigl(a^j_{is}+\rho_j(q_s)\bigr).
\]
This is a bounded regulator. In
\eqref{fv:geometry-candidate-error}, take $n=\varphi(1)$ and
$\delta=\delta(j,\varphi)$. Then
\[
 \begin{split}
 \rho_j(S_R(h^\diamond)-V)
 &\le H_{a^j}(\varphi)+\rho_j(q_{\varphi(1)})\\
 &\le H_{b^j}(\varphi)
 \end{split}
\]
for every $R\in\mathscr R_C(\delta(j,\varphi)\wedge
\eta_{\varphi(1)})$. Hence
\[
 \rho_j(S_R(h^\diamond)-V)\longrightarrow_D0\quad(j),
 \qquad S_R(h^\diamond)\longrightarrow_\lambda V.
\]
The last implication is inverse faithfulness. The equality of response
tails in Proposition~\ref{cr:response-tails} identifies the same limit
for the refinement net. Every prescribed net satisfying the original
tail correspondence has this limit. Thus
$J_{h^\diamond}(C)=V$ in the original convergence.

The argument applies to every $C\in\mathscr C_A$ and all three
polarities. Primitive additivity now applies to each hereditarily
integrable function $h^\diamond$ and gives
\[
 V_h^\diamond(C)=\sum_{D\in\mathcal D}V_h^\diamond(D)
 \qquad(\mathcal D\in\mathfrak D(C)).
\]
For integrable $r,s$ with $|r|,|s|\le mg$, $m\in\mathbb N$,
\[
 |r-s|\le2mg,
 \qquad
 r\vee s=s+(r-s)^+,
 \qquad
 r\wedge s=r-(r-s)^+.
\]
The function $2mg$ is positive and hereditarily integrable.
The proved positive-part assertion applies to $r-s$. Additivity of
the integral gives hereditary integrability of $r\vee s$ and
$r\wedge s$. Finite iteration proves finite lattice closure.
\end{proof}

\subsection{Dominated limits and uniform primitive bounds}

Finite lattice closure allows the preceding monotone construction to be
applied to the upper and lower approximations of a dominated family. In
this passage, one integrable envelope controls both constructions, and
their primitive bounds will control the intervening functions uniformly
over compatible cells.

\begin{theorem}[Dominated convergence in the original locally solid convergence]
\label{fv:dominated}\label{lc:dominated}
Retain the local division, coefficient, order, and faithful representation
hypotheses of Theorem~\ref{fv:envelope-monotone}. Assume finite-supremum stability of the original convergence,
modulus preservation of the coefficients, and either structural alternative
of Theorem~\ref{fv:finite-lattice-closure}. The controllers need only
Dedekind $\sigma$-completeness and WSD.
Let $g\ge0$ and $f_n$,
$n\ge1$, be hereditarily integrable,
with $|f_n|\le g$. Suppose either:
\begin{enumerate}
\item $f_n$ approximates $f$ through the countable order profiles of
Definition~\ref{lc:profiles}, or
\item there is an input regulator $a$ with
\[
 \inf_{\varphi\in\mathbb N^{\mathbb N}}H_a(\varphi)=0,
 \qquad
 n\ge\theta(t,\varphi)\Longrightarrow
 |f_n(t)-f(t)|\le H_a(\varphi).
\]
\end{enumerate}
The displayed infimum condition is automatic when $E$ is WSD.
Then $f$ is hereditarily integrable. There is a sequence $d_N\downarrow0$
in $Z$ such that
\begin{equation}\label{lc:dominated-uniform}
 |J_{f_n}(C)-J_f(C)|\le d_N
 \qquad(n\ge N,\ C\in\mathscr C_A).
\end{equation}
The primitive convergence is uniform in the original uniform
convergence represented by the lattice distances.
\end{theorem}

\begin{proof}
In the profile case, for every $t,k$ and all sufficiently large $n$,
\[
 (|f(t)|-g(t))^+
 \le |f(t)-f_n(t)|\le e_k^{\jmath(t)}.
\]
The infimum over $k$ gives $|f|\le g$. In the regulator case the same
inequality has right-hand side $H_a(\varphi)$. Its infimum over
$\varphi$ is zero. Both alternatives also give pointwise order
convergence. For the regulator alternative, the controls are
\[
 c_K=\bigwedge_{\varphi\in K}H_a(\varphi),
 \qquad
 K\subset\mathbb N^{\mathbb N}\text{ finite and nonempty},
 \qquad c_K\downarrow0.
\]

For $m\ge N$, define
\[
 u_{N,m}=\bigvee_{N\le r\le m}f_r,
 \qquad l_{N,m}=\bigwedge_{N\le r\le m}f_r,
 \qquad
 u_N=\bigvee_{r\ge N}f_r,
 \qquad l_N=\bigwedge_{r\ge N}f_r.
\]
The countable suprema and infima exist pointwise in $[-g,g]$.
Theorem~\ref{fv:finite-lattice-closure} gives hereditary integrability
of $u_{N,m}$ and $l_{N,m}$. At each finite iteration its positive-part
argument is dominated by $2g$, whose hereditary integrability follows
from additivity.
Write $G(C)=J_g(C)$. Positivity gives
\begin{equation}\label{lc:finite-envelope-bounds}
 \begin{aligned}
 -G(C)&\le J_{u_{N,m}}(C)\le G(C),\\
 -G(C)&\le J_{l_{N,m}}(C)\le G(C).
 \end{aligned}
\end{equation}

Fix $t,k$ in the profile alternative and put
$e=e_k^{\jmath(t)}$. If
$m\ge\max\{N,\theta(t,k)\}$ and $r>m$, then
\[
 f_r(t)\le f(t)+e\le f_m(t)+2e\le u_{N,m}(t)+2e,
\]
\[
 f_r(t)\ge f(t)-e\ge f_m(t)-2e\ge l_{N,m}(t)-2e.
\]
Taking the countable supremum and infimum gives
\begin{equation}\label{lc:inner-envelope-profiles}
 0\le u_N(t)-u_{N,m}(t)\le2e,
 \qquad
 0\le l_{N,m}(t)-l_N(t)\le2e.
\end{equation}
For fixed $N$, the sequences $u_{N,m}$ and $-l_{N,m}$ are increasing.
They and their pointwise limits are dominated in modulus by the known
integrable function $g$.
Theorem~\ref{fv:envelope-monotone}, applied separately to these two sequences with
profiles $2e_k^j$ and thresholds $\max\{N,\theta(t,k)\}$, gives
hereditary integrability of $u_N$ and $l_N$, and
\begin{equation}\label{lc:inner-envelope-primitives}
 U_N(C):=J_{u_N}(C)=\sup_{m\ge N}J_{u_{N,m}}(C),
 \qquad
 L_N(C):=J_{l_N}(C)=\inf_{m\ge N}J_{l_{N,m}}(C).
\end{equation}

For the common-regulator alternative, fix $t,\varphi$ and put
$e=H_a(\varphi)$. For
$m\ge\max\{N,\theta(t,\varphi)\}$ the two inequalities preceding
\eqref{lc:inner-envelope-profiles} hold with this $e$.
Thus both increasing sequences have the common input regulator $2a$,
and $g$ is a common positive integrable envelope for their terms and
pointwise limits.
The common-regulator form of Theorem~\ref{fv:envelope-monotone} proves
\eqref{lc:inner-envelope-primitives} in this alternative as well.

For every $n\ge N$,
\[
 (f(t)-u_N(t))^+\le |f(t)-f_n(t)|,
 \qquad
 (l_N(t)-f(t))^+\le |f_n(t)-f(t)|.
\]
Pointwise order convergence makes both fixed positive quantities zero.
In particular,
\[
 -g\le l_N\le l_{N+1}\le f\le u_{N+1}\le u_N\le g.
\]
If $N\ge\theta(t,k)$ in the profile case, every $r\ge N$ satisfies
$f(t)-e_k^{\jmath(t)}\le f_r(t)\le f(t)+e_k^{\jmath(t)}$.
Therefore
\begin{equation}\label{lc:outer-envelope-profiles}
 0\le f(t)-l_N(t)\le e_k^{\jmath(t)},
 \qquad
 0\le u_N(t)-f(t)\le e_k^{\jmath(t)}.
\end{equation}
For the regulator case, $N\ge\theta(t,\varphi)$ gives the two
corresponding inequalities with right-hand side $H_a(\varphi)$.
Consequently,
\[
 \begin{gathered}
 l_N\uparrow f,\qquad -u_N\uparrow-f,\\
 -G(A)\le L_N(A)\le G(A),\qquad
 -G(A)\le -U_N(A)\le G(A).
 \end{gathered}
\]
The functions $l_N$, $-u_N$, and their pointwise limits are dominated
in modulus by $g$. Theorem~\ref{fv:envelope-monotone}, applied to
$l_N$ and then to $-u_N$, proves
hereditary integrability of $f$ and $-f$. It yields
\[
 J_f(C)=\sup_NL_N(C),
 \qquad
 J_{-f}(C)=\sup_N[-U_N(C)].
\]
Additivity of the prescribed sums, continuity of inversion, and
uniqueness of their limits give $J_{-f}(C)=-J_f(C)$. Hence
\begin{equation}\label{lc:outer-envelope-primitives}
 L_N(C)\uparrow J_f(C),\qquad U_N(C)\downarrow J_f(C).
\end{equation}
The four applications use the input profiles $2e,2e,e,e$, or the
input regulators $2a,2a,a,a$, respectively.

For $n\ge N$, positivity gives
\[
 \begin{aligned}
 L_N(C)&\le J_{f_n}(C)\le U_N(C),\\
 L_N(C)&\le J_f(C)\le U_N(C).
 \end{aligned}
\]
Thus
\[
 |J_{f_n}(C)-J_f(C)|\le U_N(C)-L_N(C).
\]
The primitive $U_N-L_N$ is positive and additive on every division.
If $C$ occurs in a division $\mathcal D$ of $A$, then
\[
 0\le U_N(C)-L_N(C)
 \le\sum_{D\in\mathcal D}[U_N(D)-L_N(D)]
 =U_N(A)-L_N(A).
\]
Set $d_N=U_N(A)-L_N(A)$. Equation~\eqref{lc:outer-envelope-primitives}
gives $d_N\downarrow0$, and proves \eqref{lc:dominated-uniform}.
Proposition~\ref{lc:norm-normality} gives
$\rho_j(d_N)\downarrow0$ for every $j$. Put
\[
 \mathcal E_N=
 \{(J_{f_n}(C),J_f(C)):n\ge N,\ C\in\mathscr C_A\},
 \qquad b^j_{is}=\rho_j(d_s).
\]
Each $b^j$ is an order-bounded regulator, and
\[
 (x,y)\in\mathcal E_{\varphi(1)}
 \Longrightarrow\rho_j(x-y)
 \le\rho_j(d_{\varphi(1)})\le H_{b^j}(\varphi).
\]
Therefore all represented differences of
$\mathsf N(\mathcal E_N)$ converge to zero. The inverse pair
implication in \eqref{lc:pair-faithfulness} gives
\[
 \mathsf N(\mathcal E_N)
 \longrightarrow_{\mathcal U_\lambda}\Delta_Z.
\]
This is the asserted uniform convergence in the original structure.
\end{proof}

\subsection{Mixed control without an integrable envelope}

An integrable envelope supplied the common bounds in the preceding
construction; without one, we turn to the controllers to construct
summable selections from the individual estimates. We therefore impose
the additional order conditions of Definition~\ref{lc:calibration} on
each controller, both here and in Theorem~\ref{lc:levi}. These
conditions belong to this extension: Theorems~\ref{fv:envelope-monotone}
and~\ref{fv:dominated} retain their envelope hypotheses and do not
require them.

\begin{proposition}[Mixed control derived from individual integrals]
\label{lc:mixed}\label{lc:mixed-monoid}
Let $Z$ be a separated commutative convergence monoid with symmetric
invariant pseudometrics $d_j:Z^2\to(F_j)_+$. Suppose original point convergence implies regulator
convergence of every $d_j$, and every $F_j$ satisfies
Definition~\ref{lc:calibration}. Let kernels $U_n$ have known hereditary
original integrals $J_n$. For every $j$ there exist $B_{j,s}\downarrow0$
and gauges $\eta_{j,n,s}$ such that
\[
 \delta_{j,s,p}(t)=\bigwedge_{n\le p(t)}\eta_{j,n,s}(t)
 \qquad(p:T\to\mathbb N_{\ge1})
\]
satisfies
\begin{equation}\label{lc:mixed-monoid-bound}
 d_j\!\left(\sum_rU_{\nu_r}(t_r,C_r),
             \sum_rJ_{\nu_r}(C_r)\right)\le B_{j,s}
\end{equation}
for all compatible $\delta_{j,s,p}$-fine partial records with
$1\le\nu_r\le p(t_r)$. In particular, for an abelian group,
$U_n(t,C)=L_Ch_n(t)$ and $d_j(z,w)=\rho_j(z-w)$ give
\begin{equation}\label{lc:mixed-bound}
 \rho_j\left(\sum_r[L_{C_r}h_{\nu_r}(t_r)-J_{\nu_r}(C_r)]\right)
 \le B_{j,s}.
\end{equation}
\end{proposition}
\begin{proof}
For the known root integral of $U_n$, let $a^{j,n}$ regulate the metric
error. Corollary~\ref{lc:regulator-calibration} supplies selectors
$\psi_{j,n,s}$ and $B_{j,s}\downarrow0$ such that
\[
 \sum_{n\in M}H_{a^{j,n}}(\psi_{j,n,s})\le B_{j,s}
 \qquad(M\subset\mathbb N\text{ finite}).
\]
Let $\eta_{j,n,s}$ attain the corresponding full-record bound.
Lemma~\ref{lt:d-calculus} converts the known residual regulator limits
into order limits, so Lemma~\ref{lc:monoid-partial} gives
\[
 d_j\!\left(\sum_{(C,t)\in Q_n}U_n(t,C),
             \sum_{(C,t)\in Q_n}J_n(C)\right)
 \le H_{a^{j,n}}(\psi_{j,n,s})
\]
for every compatible $\eta_{j,n,s}$-fine partial record $Q_n$.
For a record in the assertion, put
\[
 M=\{\nu_r:r\},\qquad
 x_n=\sum_{r:\nu_r=n}U_n(t_r,C_r),\qquad
 y_n=\sum_{r:\nu_r=n}J_n(C_r).
\]
Its $n$th subword is $\eta_{j,n,s}$-fine, since
$\delta_{j,s,p}(t_r)\le\eta_{j,n,s}(t_r)$ whenever $\nu_r=n$.
For any two pairs in $Z$, invariance and the triangle inequality give
\[
\begin{aligned}
 d_j(x_1+x_2,y_1+y_2)
 &\le d_j(x_1+x_2,y_1+x_2)+d_j(y_1+x_2,y_1+y_2)\\
 &=d_j(x_1,y_1)+d_j(x_2,y_2).
\end{aligned}
\]
Finite induction, including the empty sum, yields
\[
 d_j\!\left(\sum_{n\in M}x_n,\sum_{n\in M}y_n\right)
 \le\sum_{n\in M}d_j(x_n,y_n)
 \le\sum_{n\in M}H_{a^{j,n}}(\psi_{j,n,s})\le B_{j,s}.
\]
This proves \eqref{lc:mixed-monoid-bound}. Substituting the displayed
kernel and invariant group distances gives \eqref{lc:mixed-bound}.
\end{proof}

\begin{proposition}[Data transfer through the original lattice norms]
\label{lc:data}
Under either input condition, every $j$ admits $V_{j,s}\downarrow0$ and
integers $\ell_{j,s}(t)$ such that
\[
 \nu_r\ge\ell_{j,s}(t_r)\Longrightarrow
 \rho_j\left(\sum_rL_{C_r}[h(t_r)-h_{\nu_r}(t_r)]\right)\le V_{j,s}
\]
for every compatible partial record.
\end{proposition}
\begin{proof}
For profiles, Proposition~\ref{lc:norm-normality} gives
\[
 v_{r,k}=\rho_j(L_Ae_k^r)\downarrow_k0.
\]
Theorem~\ref{lc:series-selection} supplies $k(r,s)$ with
$V_{j,s}=\sum_rv_{r,k(r,s)}\downarrow0$. Put
$\ell_{j,s}(t)=\theta(t,k(\jmath(t),s))$. Grouping by the finite set $M$
of active profiles, solidness gives
\[
 \begin{aligned}
 &\rho_j\left(\sum_rL_{C_r}[h(t_r)-h_{\nu_r}(t_r)]\right)\\
 &\quad\le\sum_{q\in M}
 \rho_j\left(\sum_{\jmath(t_r)=q}L_{C_r}|h(t_r)-h_{\nu_r}(t_r)|\right)\\
 &\quad\le\sum_{q\in M}\rho_j(L_Ae_{k(q,s)}^q)\le V_{j,s}.
 \end{aligned}
\]

For a common regulator $a$, apply Theorem~\ref{lc:series-selection} to
$v_{i,k}=\rho_j(L_Aa_{ik})\downarrow0$. It supplies selectors $\psi_s$
with
\[
 \sum_i\rho_j(L_Aa_{i,\psi_s(i)})=V_{j,s}\downarrow0.
\]
For $y_N=\bigvee_{i\le N}a_{i,\psi_s(i)}$ and
$y=H_a(\psi_s)$,
\[
 y_N\uparrow y,\qquad
 L_A(y-y_N)\downarrow0,
\]
\[
 \rho_j(L_Ay_N)
 \le\sum_{i\le N}\rho_j(L_Aa_{i,\psi_s(i)})\le V_{j,s},
\]
\[
 \rho_j(L_Ay)\le V_{j,s}+\rho_j(L_A(y-y_N)).
\]
The final summand decreases to zero, so $\rho_j(L_Ay)\le V_{j,s}$.
With $\ell_{j,s}(t)=\theta(t,\psi_s)$,
\[
 \left|\sum_rL_{C_r}[h(t_r)-h_{\nu_r}(t_r)]\right|\le L_Ay.
\]
Solidness proves the required estimate.
\end{proof}

\subsection{The monotone limit and uniformity of the primitives}

The mixed and data estimates can now be combined with the order bound on
the known primitives. Their supremum provides the candidate primitive,
and the remaining task is to show that the original response nets
converge to it, with a bound that remains valid as the cell varies.

\begin{theorem}[Bounded-primitives monotone theorem]\label{lc:levi}
Assume the standing local division, coefficient, order, and faithful
representation conditions, together with Definition~\ref{lc:calibration}
for every controller $F_j$. Let $\lambda$ be locally solid and
sequentially order compatible. Suppose $h_n$ have hereditary original
integrals, and
\[
 h_n(t)\uparrow h(t)\quad(t),\qquad
 \{J_n(A):n\ge1\}\text{ is order bounded above}.
\]
Assume either condition of Definition~\ref{lc:profiles}.
Then $h$ has hereditary original integrals,
\[
 J_h(C)=\bigvee_nJ_n(C)\quad(C\in\mathscr C_A),
\]
and
\begin{equation}\label{lc:levi-uniform}
 q_N:=J_h(A)-J_N(A)\downarrow0,\qquad
 0\le J_h(C)-J_n(C)\le q_N
 \quad(n\ge N,\ C\in\mathscr C_A).
\end{equation}
Consequently the primitives converge uniformly in the original
uniform convergence structure.
\end{theorem}
\begin{proof}
Positivity and finite additivity give
\[
 0\le J_n(C)-J_1(C)\le J_n(A)-J_1(A)
 \le K:=\bigvee_m[J_m(A)-J_1(A)].
\]
Dedekind $\sigma$-completeness defines $J(C)=\bigvee_nJ_n(C)$.
For $(C_r)_{r=1}^m\in\mathfrak D(C)$,
\[
 J(C)=\bigvee_n\sum_rJ_n(C_r)
 \le\sum_rJ(C_r).
\]
For any $(n_1,\ldots,n_m)$, the integer $N=\max_rn_r$ gives
\[
 \sum_rJ_{n_r}(C_r)\le\sum_rJ_N(C_r)\le J(C).
\]
Taking the finite iterated supremum yields
$\sum_rJ(C_r)\le J(C)$. Therefore $J$ is finitely additive.
The positive additive difference $J-J_N$ satisfies
\[
 0\le J(C)-J_N(C)\le q_N:=J(A)-J_N(A),\qquad q_N\downarrow0.
\]

Fix $j,s$. Propositions~\ref{lc:mixed} and~\ref{lc:data} give
$B_{j,s},V_{j,s},\ell_{j,s}$ and gauges $\eta_{j,n,s}$. Put
\[
 p_{j,s}(t)=\max\{s,\ell_{j,s}(t)\},\qquad
 \delta_{j,s}(t)=\bigwedge_{n\le p_{j,s}(t)}\eta_{j,n,s}(t).
\]
For $R=((C_r,t_r))\in\mathscr R_C(\delta_{j,s})$, set
\[
 M_R=\sum_rL_{C_r}h_{p_{j,s}(t_r)}(t_r),\qquad
 P_R=\sum_rJ_{p_{j,s}(t_r)}(C_r).
\]
Every record on a root-compatible cell is a compatible partial record
on $A$. Thus
\[
 \rho_j(S_Rh-M_R)\le V_{j,s},\qquad
 \rho_j(M_R-P_R)\le B_{j,s}.
\]
Monotonicity and finite additivity give
\[
 J_s(C)\le P_R\le J(C),\qquad
 |P_R-J(C)|\le q_s.
\]
Hence
\begin{equation}\label{lc:final-estimate}
 \rho_j(S_Rh-J(C))
 \le V_{j,s}+B_{j,s}+\rho_j(q_s)=:E_{j,s},
 \qquad E_{j,s}\downarrow0.
\end{equation}
The last decrease uses Proposition~\ref{lc:norm-normality}.
Define the regulator $c^j_{ik}=E_{j,k}$. For each selector $\varphi$,
\[
 R\in\mathscr R_C(\delta_{j,\varphi(1)})
 \Longrightarrow
 \rho_j(S_Rh-J(C))\le E_{j,\varphi(1)}\le H_{c^j}(\varphi).
\]
For fixed $C$, these inequalities and the exact response tails give
\[
 \rho_j(S_Rh-J(C))\longrightarrow_D0\quad(j),
 \qquad S_Rh\longrightarrow_\lambda J(C).
\]
The second statement uses the inverse implication in
\eqref{lc:pair-faithfulness}. Separation gives $J_h(C)=J(C)$.
The same regulators also control the varying-cell pair families
\[
 \mathcal P_\delta=
 \{(S_Rh,J_h(C)):C\in\mathscr C_A,\ R\in\mathscr R_C(\delta)\},
\]
\[
 (x,y)\in\mathcal P_{\delta_{j,\varphi(1)}}
 \Longrightarrow\rho_j(x-y)\le H_{c^j}(\varphi).
\]
Thus $\mathsf N(\mathcal P_\delta)\to_{\mathcal U_\lambda}\Delta_Z$.
For the endpoint families
\[
 \mathcal E_N=
 \{(J_n(C),J_h(C)):n\ge N,\ C\in\mathscr C_A\},
\]
\[
 (x,y)\in\mathcal E_N\Longrightarrow
 \rho_j(x-y)\le\rho_j(q_N)\downarrow_N0.
\]
The rank-one regulators $b^j_{is}=\rho_j(q_s)$ and pair faithfulness
give $\mathsf N(\mathcal E_N)\to_{\mathcal U_\lambda}\Delta_Z$.
Thus both conclusions use the original uniform convergence.
\end{proof}

Here Dedekind $\sigma$-completeness of $Z$ constructs the candidate
$J(C)=\bigvee_nJ_n(C)$, while the order hypotheses on the controllers
construct the decreasing bounds $B_{j,s}$ and $V_{j,s}$ that establish
convergence to it in the original structure. In this way the proof
requires neither a separate completeness assumption on $\lambda$ nor
a prescribed common output regulator or integrable envelope.

\subsection{Increasing barriers and convergence of primitive variation}
\label{vt:section}

Once bounded primitives suffice for the monotone construction, a single
integrable envelope can be replaced by an increasing sequence of
integrable barriers. We shall control the discarded parts through
primitives of the known functions, using their decay first to construct
the integral of the limit and then to control the variation of its
primitive. This extension retains the controller assumptions of
Definition~\ref{lc:calibration} throughout the subsection, and therefore
belongs to the controller class of Theorem~\ref{lc:levi}; the weaker
controller assumptions of Theorem~\ref{fv:dominated} remain sufficient
for that theorem itself.

\begin{lemma}[Truncation in an abelian lattice-ordered group]
\label{vt:truncation}
Let $E$ be an abelian lattice-ordered group. For $a\in E_+$ and $x\in E$, put
\[
 T_a x=(x\vee(-a))\wedge a,
 \qquad M_a x=|x|\wedge a,
 \qquad R_a x=(|x|-a)^+.
\]
Then
\begin{equation}\label{vt:truncation-identities}
 |T_a x|=M_a x,
 \qquad |x-T_a x|=R_a x,
 \qquad M_a x+R_a x=|x|,
\end{equation}
and
\begin{equation}\label{vt:truncation-lipschitz}
 |T_a x-T_a y|\le |x-y|,
 \qquad |M_a x-M_a y|\le |x-y|.
\end{equation}
For $0\le a\le b$,
\begin{equation}\label{vt:barrier-increment}
 0\le M_bx-M_ax\le R_ax.
\end{equation}
If $|x|\le a$, then $T_ax=x$ and $M_ax=|x|$.
\end{lemma}

\begin{proof}
Translations preserve finite meets and joins, inversion interchanges them,
and the lattice of an abelian lattice-ordered group is distributive. The
identities
\[
 (u\vee v)+(u\wedge v)=u+v,
 \qquad u\vee v=u+(v-u)^+,
 \qquad u\wedge v=u-(u-v)^+
\]
will be used with parentheses wherever sums and lattice operations occur
together. In particular, for $u,v\ge0$ with $u\wedge v=0$,
\[
 u\vee v=u+v,
 \qquad |u-v|=u+v.
\]
Distributivity and $-a\le0\le a$ give
\[
 \begin{aligned}
 (T_ax)^+
 &=\bigl((x\vee(-a))\wedge a\bigr)\vee0\\
 &=\bigl(x\vee(-a)\vee0\bigr)\wedge(a\vee0)
 =x^+\wedge a.
 \end{aligned}
\]
Moreover,
\[
 \begin{aligned}
 -T_ax
 &=\bigl((-x)\wedge a\bigr)\vee(-a)\\
 &=\bigl((-x)\vee(-a)\bigr)\wedge\bigl(a\vee(-a)\bigr)
 =T_a(-x),
 \end{aligned}
\]
so $(T_ax)^-=x^-\wedge a$. The two positive elements on the right are
disjoint, since their meet is at most $x^+\wedge x^-=0$. Consequently,
\[
 \begin{aligned}
 |T_ax|
 &=(x^+\wedge a)+(x^-\wedge a)\\
 &=(x^+\wedge a)\vee(x^-\wedge a)
 =(x^+\vee x^-)\wedge a
 =|x|\wedge a.
 \end{aligned}
\]
Put $p=(x-a)^+$ and $q=(-x-a)^+$. Positivity of $a$ gives
\[
 0\le p\le x^+,
 \qquad 0\le q\le x^-,
 \qquad p\wedge q=0.
\]
The join and meet identities give
\[
 x\vee(-a)=x+q,
\]
\[
 \begin{aligned}
 \bigl((x\vee(-a))-a\bigr)^+
 &=\bigl((x-a)\vee(-2a)\bigr)\vee0\\
 &=(x-a)\vee0=p,
 \end{aligned}
\]
and hence $T_ax=x+q-p$. Therefore
\[
 \begin{aligned}
 |x-T_ax|
 &=|p-q|=p\vee q\\
 &=\bigl((x-a)\vee0\bigr)\vee\bigl((-x-a)\vee0\bigr)\\
 &=\bigl((x\vee(-x))-a\bigr)\vee0
 =(|x|-a)^+.
 \end{aligned}
\]
The formula $|x|\wedge a=|x|-(|x|-a)^+$ proves the third identity
in \eqref{vt:truncation-identities}.

For the Lipschitz estimates, set $u=|x-y|$. Then
$x\le y+u$ and $y\le x+u$. For every $c\in E$,
\[
 x\vee c\le(y+u)\vee(c+u)=(y\vee c)+u,
\]
\[
 x\wedge c\le(y+u)\wedge(c+u)=(y\wedge c)+u.
\]
Interchanging $x$ and $y$ bounds the negative of each difference by $u$.
Thus
\[
 |(x\vee c)-(y\vee c)|\le u,
 \qquad |(x\wedge c)-(y\wedge c)|\le u.
\]
The first inequality followed by the second, with $c=-a$ and $c=a$,
respectively, proves the assertion for $T_a$. The inequalities
\[
 |x|\le|y|+|x-y|,
 \qquad |y|\le|x|+|x-y|
\]
give $\bigl||x|-|y|\bigr|\le|x-y|$. The meet estimate then proves
the assertion for $M_a$.
Finally, $a\le b$ gives $M_ax\le M_bx\le|x|$, and
\[
 0\le M_bx-M_ax\le|x|-M_ax=R_ax.
\]
If $|x|\le a$, then $-a\le x\le a$, so both the join and the meet
in $T_ax$ leave $x$ unchanged. The formula for $M_ax$ is immediate
from its defining meet.
\end{proof}

\begin{theorem}[Increasing-barrier convergence and primitive variation]
\label{vt:vitali}
Retain the local division calculus, finite substitution and completion,
and original response-tail correspondence of
Theorem~\ref{fv:dominated}. Let $E$ and $Z$ be Dedekind
$\sigma$-complete abelian lattice-ordered groups. Let the original
separated group convergence $\lambda$ on $Z$ be locally solid,
sequentially order compatible, and finite-supremum stable. Retain its
faithful solid lattice representation from Definition~\ref{lc:faithful}.
Every controller $F_j$ is assumed to satisfy
Definition~\ref{lc:calibration}. Thus it is Dedekind $\sigma$-complete,
laterally $\sigma$-complete and weakly $\sigma$-distributive, and every
family of components of a positive element with that element as its
supremum has a countable subfamily with the same supremum.

The coefficient maps $L_C:E\to Z$ are positive, additive and
modulus preserving, and satisfy
\[
 L_C=\sum_{D\in\mathcal D}L_D
 \quad(\mathcal D\in\mathfrak D(C)),
 \qquad e_s\downarrow0\Longrightarrow L_Ae_s\downarrow0.
\]
Assume either countable division resolution on every cell, or
finite-division resolution on every cell together with Dedekind
completeness of $Z$ and the countable sup property on each of its
principal order intervals, as specified in
Theorem~\ref{fv:finite-lattice-closure}.

Let $f_n,f:T\to E$ and $g_k:T\to E_+$ satisfy the following conditions.
Every $f_n$, $|f_n|$ and $g_k$ is hereditarily integrable in the original
convergence, and
\begin{equation}\label{vt:local-barriers}
 g_k\le g_{k+1},
 \qquad
 \forall t\in T\ \exists k(t)\ge1:
 \quad |f(t)|\le g_{k(t)}(t).
\end{equation}
Assume either countable order profiles
for $f_n\to f$ in Definition~\ref{lc:profiles}, or a regulator $a$ in
$E$ satisfying
\begin{equation}\label{vt:common-input}
 \inf_{\varphi\in\mathbb N^{\mathbb N}}H_a(\varphi)=0,
 \qquad
 n\ge\theta(t,\varphi)
 \Longrightarrow |f_n(t)-f(t)|\le H_a(\varphi).
\end{equation}
Suppose that $v_k\downarrow0$ in $Z$ and
\begin{equation}\label{vt:known-tail-bound}
 J_{(|f_n|-g_k)^+}(A)\le v_k
 \qquad(n,k\ge1).
\end{equation}
The functions in \eqref{vt:known-tail-bound} are integrable by the
preceding hypotheses and finite lattice closure.

Then $f$ and $|f|$ are hereditarily integrable. The countable suprema
\begin{equation}\label{vt:variation-control}
 d_N:=\bigvee_{n\ge N}J_{|f_n-f|}(A)
\end{equation}
exist in $Z$ and satisfy $d_N\downarrow0$. For every finite family
$(C_r)_{r=1}^q$ forming a subfamily of a division of $A$, including
finite substitutions,
\begin{equation}\label{vt:finite-variation}
 \sum_{r=1}^q|J_{f_n}(C_r)-J_f(C_r)|\le d_N
 \qquad(n\ge N).
\end{equation}
The primitives therefore converge uniformly in the original uniform
convergence, and the sums in \eqref{vt:finite-variation} converge to
zero in the original convergence uniformly over all such finite
families.
The local bound in \eqref{vt:local-barriers} is guaranteed, in
particular, by the condition involving only the given functions
\begin{equation}\label{vt:data-local-barriers}
 \forall t\in T\ \exists\ell(t)\ge1\ \forall n\ge1:
 \quad |f_n(t)|\le g_{\ell(t)}(t).
\end{equation}
\end{theorem}

\begin{proof}
We first record how inequalities are preserved by the original
convergence. If $z_\alpha\ge0$ and $z_\alpha\to_\lambda z$, then
\[
 0\le z^-\le|z_\alpha-z|.
\]
Local solidity makes the constant net at $z^-$ null, and separation
gives $z^-=0$. Hence $z\ge0$. For existing integrals of functions
$u\le w$, the positive response net $S_R(w-u)$ has limit
$J_w(C)-J_u(C)$, so
\begin{equation}\label{vt:positive-integrals}
 u\le w\Longrightarrow J_u(C)\le J_w(C).
\end{equation}
Additivity of existing integrals is obtained by adding their response
nets on their common directed record set. In particular,
\[
 J_{u+w}(C)=J_u(C)+J_w(C),
 \qquad J_{-u}(C)=-J_u(C).
\]
Their finite additivity on cells is
Proposition~\ref{ls:additivity}.

To verify the stated sufficient condition, suppose
\eqref{vt:data-local-barriers} holds. For a tag $t$, let $\ell(t)$
be the least positive integer satisfying that condition. The inequality
\[
 |f(t)|\le|f_n(t)|+|f(t)-f_n(t)|
 \le g_{\ell(t)}(t)+|f(t)-f_n(t)|
\]
gives
\[
 0\le\bigl(|f(t)|-g_{\ell(t)}(t)\bigr)^+
 \le|f(t)-f_n(t)|.
\]
In the profile case, for every $s\ge1$ and every
$n\ge\theta(t,s)$ the last member is bounded by
$e_s^{\jmath(t)}$, whose infimum is zero. In the regulator case, for
every $\varphi$ and every $n\ge\theta(t,\varphi)$ it is bounded by
$H_a(\varphi)$, whose infimum is zero by
\eqref{vt:common-input}. Therefore
\[
 |f(t)|\le g_{\ell(t)}(t).
\]
This proves that \eqref{vt:data-local-barriers} implies the local
bound in \eqref{vt:local-barriers}. Under either formulation, the
integer
\[
 k(t)=\min\{k\ge1:|f(t)|\le g_k(t)\}
\]
is defined and satisfies
\begin{equation}\label{vt:limit-local-bound}
 |f(t)|\le g_{k(t)}(t)
 \qquad(t\in T).
\end{equation}

For $n,k\ge1$, define pointwise
\[
 h_{n,k}=M_{g_k}f_n=|f_n|\wedge g_k,
 \qquad p_{n,k}=T_{g_k}f_n,
 \qquad r_{n,k}=R_{g_k}f_n=(|f_n|-g_k)^+.
\]
Here and below the subscript $g_k$ means that the truncating element
at tag $t$ is $g_k(t)$. The positive function
$w_{n,k}=|f_n|+g_k$ is hereditarily integrable. The functions
$f_n$, $|f_n|$, $g_k$ and $-g_k$ are bounded in modulus by
$w_{n,k}$. Finite lattice closure in
Theorem~\ref{fv:finite-lattice-closure} proves the hereditary
integrability of $h_{n,k}$ and $p_{n,k}$. The function
$|f_n|-g_k$ is integrable and bounded in modulus by $w_{n,k}$, so
the same theorem proves the hereditary integrability of $r_{n,k}$.
Thus every primitive in \eqref{vt:known-tail-bound} is already
defined at this stage.

Put
\[
 h_k=|f|\wedge g_k,
 \qquad p_k=T_{g_k}f.
\]
For each fixed $k$, Lemma~\ref{vt:truncation} gives
\[
 0\le h_{n,k}\le g_k,
 \qquad |p_{n,k}|\le g_k,
\]
\[
 |h_{n,k}(t)-h_k(t)|\le|f_n(t)-f(t)|,
 \qquad |p_{n,k}(t)-p_k(t)|\le|f_n(t)-f(t)|.
\]
In the profile alternative, both differences are bounded by
$e_s^{\jmath(t)}$ for $n\ge\theta(t,s)$. In the regulator
alternative, both are bounded by $H_a(\varphi)$ for
$n\ge\theta(t,\varphi)$. Thus the original input controls, their
thresholds, and their zero-infimum conditions remain unchanged.
The known integrable envelope for each of these two sequences is
$g_k$. Theorem~\ref{fv:dominated} gives hereditary original
integrals of $h_k$ and $p_k$ and, for every cell $C$,
\begin{equation}\label{vt:fixed-barrier-limits}
 J_{h_{n,k}}(C)\longrightarrow_\lambda J_{h_k}(C),
 \qquad
 J_{p_{n,k}}(C)\longrightarrow_\lambda J_{p_k}(C).
\end{equation}

Write $G_k=J_{g_k}$. The decomposition in
\eqref{vt:truncation-identities} and
\eqref{vt:known-tail-bound} give
\[
 \begin{aligned}
 0\le J_{|f_n|}(A)
 &=J_{h_{n,1}}(A)+J_{r_{n,1}}(A)\\
 &\le G_1(A)+v_1=:M.
 \end{aligned}
\]
Since $0\le h_{n,k}\le|f_n|$,
\[
 0\le J_{h_{n,k}}(A)\le M.
\]
For fixed $k$, the positive net
$M-J_{h_{n,k}}(A)$ converges to $M-J_{h_k}(A)$ by
\eqref{vt:fixed-barrier-limits}. Preservation of positivity yields
\begin{equation}\label{vt:bounded-truncation-primitives}
 0\le J_{h_k}(A)\le M
 \qquad(k\ge1).
\end{equation}
The monotonicity of the barriers and
\eqref{vt:limit-local-bound} give
\[
 0\le h_k\le h_{k+1}\le|f|,
 \qquad
 k\ge k(t)\Longrightarrow h_k(t)=|f(t)|.
\]
In particular, $h_k(t)\uparrow|f(t)|$ for every tag. This sequence
has the countable order profiles
\[
 e_s^r=0\quad(r,s\ge1),
 \qquad \jmath(t)=1,
 \qquad \theta(t,s)=k(t).
\]
Every row $e_s^r$ decreases to zero. The hereditary integrability
of $h_k$, the bound
\eqref{vt:bounded-truncation-primitives}, and these explicit profiles
verify the hypotheses of Theorem~\ref{lc:levi}. Its conclusion gives
\begin{equation}\label{vt:modulus-constructed}
 |f|\text{ is hereditarily integrable},
 \qquad
 J_{h_k}(C)\uparrow J_{|f|}(C)
 \quad(C\in\mathscr C_A).
\end{equation}
This is the use of the additional controller conditions in
Definition~\ref{lc:calibration}.

Lemma~\ref{vt:truncation} also gives
\[
 |p_k|=h_k\le|f|,
 \qquad
 k\ge k(t)\Longrightarrow p_k(t)=f(t).
\]
The functions $p_k$ already have their hereditary integrals from
\eqref{vt:fixed-barrier-limits}. Their positive integrable envelope
is now the function $|f|$ constructed in
\eqref{vt:modulus-constructed}. Their profiles are again the zero
profiles just displayed. Theorem~\ref{fv:dominated} therefore
proves the hereditary original integrability of $f$.

We next transfer the known tail bound to the limit. If $m\ge k$,
\eqref{vt:barrier-increment} gives
\[
 0\le h_{n,m}-h_{n,k}\le r_{n,k}.
\]
All three functions have known integrals. Consequently,
\[
 0\le J_{h_{n,m}}(A)-J_{h_{n,k}}(A)\le v_k.
\]
For fixed $m,k$, both endpoint sequences converge by
\eqref{vt:fixed-barrier-limits}. Group continuity and preservation
of positivity yield
\[
 0\le J_{h_m}(A)-J_{h_k}(A)\le v_k.
\]
The increasing supremum in \eqref{vt:modulus-constructed} gives
\[
 0\le J_{|f|}(A)-J_{h_k}(A)\le v_k.
\]
The function
\[
 r_k:=|f|-h_k=(|f|-g_k)^+
\]
is integrable by additivity of its two known integrals. Hence
\begin{equation}\label{vt:limit-tail-bound}
 0\le J_{r_k}(A)\le v_k
 \qquad(k\ge1).
\end{equation}

For fixed $k$, define
\[
 B_{n,k}=|p_{n,k}-p_k|.
\]
The difference $p_{n,k}-p_k$ is hereditarily integrable and has
positive integrable envelope $2g_k$. Finite lattice closure gives
the hereditary integrability of $B_{n,k}$, and
\[
 0\le B_{n,k}\le2g_k,
 \qquad B_{n,k}(t)\le|f_n(t)-f(t)|.
\]
Thus $B_{n,k}\to0$ has the same input profiles or the same
regulator as $f_n\to f$. Theorem~\ref{fv:dominated} gives a
sequence $q_{N,k}\downarrow_N0$ in $Z$ such that
\[
 0\le J_{B_{n,k}}(A)\le q_{N,k}
 \qquad(n\ge N).
\]
The same primitives satisfy $J_{B_{n,k}}(A)\le2G_k(A)$ for
every $n$. Dedekind $\sigma$-completeness defines
\begin{equation}\label{vt:bounded-part-null}
 b_{N,k}:=\bigvee_{n\ge N}J_{B_{n,k}}(A),
 \qquad
 0\le b_{N,k}\le q_{N,k},
 \qquad b_{N,k}\downarrow_N0.
\end{equation}

The difference $f_n-f$ is hereditarily integrable and satisfies
$|f_n-f|\le|f_n|+|f|$. The envelope on the right is positive
and hereditarily integrable. Finite lattice closure therefore
proves that $|f_n-f|$ is hereditarily integrable. The two residual
identities in Lemma~\ref{vt:truncation} give
\[
 \begin{aligned}
 |f_n-f|
 &\le|f_n-p_{n,k}|+|p_{n,k}-p_k|+|p_k-f|\\
 &=r_{n,k}+B_{n,k}+r_k.
 \end{aligned}
\]
Positivity, additivity, \eqref{vt:known-tail-bound}, and
\eqref{vt:limit-tail-bound} imply
\begin{equation}\label{vt:three-part-error}
 a_n:=J_{|f_n-f|}(A)
 \le2v_k+J_{B_{n,k}}(A)
 \qquad(n,k\ge1).
\end{equation}
In particular,
\[
 0\le a_n\le2v_1+2G_1(A)
 \qquad(n\ge1).
\]
Thus all the countable suprema $d_N$ in
\eqref{vt:variation-control} exist. They form a decreasing positive
sequence, so its infimum $d=\bigwedge_Nd_N$ also exists. For each
$k$, \eqref{vt:bounded-part-null} and
\eqref{vt:three-part-error} give
\[
 0\le d\le d_N\le2v_k+b_{N,k}
 \qquad(N\ge1).
\]
Translation preserves existing infima. Taking the infimum in $N$
therefore gives $0\le d\le2v_k$. To verify that the infimum of
these last bounds is zero without scalar division, a positive
element $w\le2v_k$ for every $k$ also satisfies
\[
 w\le2v_{\max\{k,\ell\}}\le v_k+v_\ell
 \qquad(k,\ell\ge1).
\]
Taking the infimum in $\ell$ and then in $k$ gives $w\le0$.
Applying this argument to $w=d$ proves $d=0$, and hence
$d_N\downarrow0$.

For every cell $C$, positivity applied to
$-|f_n-f|\le f_n-f\le|f_n-f|$ gives
\[
 |J_{f_n}(C)-J_f(C)|\le J_{|f_n-f|}(C).
\]
Let $(C_r)_{r=1}^q$ be a subfamily of a root division
$\mathcal D$, with finite substitutions allowed. Positivity and
finite additivity give, for $n\ge N$,
\[
 \begin{aligned}
 \sum_{r=1}^q|J_{f_n}(C_r)-J_f(C_r)|
 &\le\sum_{r=1}^qJ_{|f_n-f|}(C_r)\\
 &\le\sum_{C\in\mathcal D}J_{|f_n-f|}(C)\\
 &=J_{|f_n-f|}(A)=a_n\le d_N.
 \end{aligned}
\]
This proves \eqref{vt:finite-variation}. No supremum over cells or
over divisions has been taken.

For the original uniform convergence, set
\[
 \mathcal E_N=
 \{(J_{f_n}(C),J_f(C)):n\ge N,\ C\in\mathscr C_A\}.
\]
Proposition~\ref{lc:norm-normality} gives
$\rho_j(d_N)\downarrow_N0$ for every $j$, and
\[
 (x,y)\in\mathcal E_N
 \Longrightarrow\rho_j(x-y)\le\rho_j(d_N).
\]
The matrix $c^j_{is}=\rho_j(d_s)$ is an order-bounded regulator.
For every selector $\varphi$,
\[
 (x,y)\in\mathcal E_{\varphi(1)}
 \Longrightarrow\rho_j(x-y)
 \le\rho_j(d_{\varphi(1)})\le H_{c^j}(\varphi).
\]
The saturated pair net therefore satisfies
\[
 \mathsf N(\mathcal E_N)
 \longrightarrow_{\mathcal U_\lambda}\Delta_Z
\]
by the inverse implication in \eqref{lc:pair-faithfulness}.
For the finite variation sums, let $\mathcal V_N$ consist of all
sums on the left of \eqref{vt:finite-variation} with $n\ge N$,
including the empty sum. These sets are nonempty and decreasing,
and
\[
 z\in\mathcal V_N\Longrightarrow
 0\le z\le d_N,
 \qquad \rho_j(z)\le\rho_j(d_N).
\]
The same regulators and point faithfulness give
$\mathsf N(\mathcal V_N)\to_\lambda0$. This proves the asserted
uniformity over all compatible finite families in the original
convergence.
\end{proof}

\begin{corollary}[Recovery of integrable domination]
\label{vt:dominated-recovery}
Within the controller class of Definition~\ref{lc:calibration},
Theorem~\ref{vt:vitali} contains the dominated convergence theorem
and supplies the finite-variation estimate
\eqref{vt:finite-variation}.
\end{corollary}

\begin{proof}
Let the hypotheses of Theorem~\ref{fv:dominated} hold with positive
integrable envelope $g$, and impose
Definition~\ref{lc:calibration} on every controller. Finite lattice
closure, with $h=f_n$ and envelope $g$, gives the hereditary
integrability of $|f_n|$. Set
\[
 g_k=g\quad(k\ge1),
 \qquad k(t)=1\quad(t\in T),
 \qquad v_k=0\quad(k\ge1).
\]
Then
\[
 |f_n|\le g_{k(t)},
 \qquad (|f_n|-g_k)^+=0,
 \qquad J_{(|f_n|-g_k)^+}(A)=0\le v_k.
\]
The input profiles or common regulator are exactly those already
assumed. Every hypothesis of Theorem~\ref{vt:vitali} has therefore
been verified, and its conclusions give both hereditary
integrability and \eqref{vt:finite-variation}.
\end{proof}

\begin{example}[Vanishing tails without an integrable common envelope]
\label{vt:strictness}
Consider real-valued interval integration on $A=[0,1]$, using the
full gauge record nets of Proposition~\ref{cr:intervals}, with
$L_Cx=|C|x$, the ordinary separated real convergence, and
$\rho(x)=|x|$. The interval division calculus satisfies the
resolution condition in Proposition~\ref{cr:intervals}. The
controller $\mathbb R$ satisfies
Definition~\ref{lc:calibration}. In particular, a disjoint family
of positive real numbers has at most one nonzero member, and the
components of $u>0$ are $0$ and $u$.

Set
\[
 c=\frac6{\pi^2},
 \qquad s_0=0,
 \qquad s_n=c\sum_{r=1}^{n}\frac1{r^2},
 \qquad I_n=[s_{n-1},s_n),
\]
\[
 f_n=n\mathbf1_{I_n},
 \qquad f=0,
 \qquad g_k=k\mathbf1_A.
\]
The intervals $I_n$ are pairwise disjoint, have length $c/n^2$,
and have union $[0,1)$. Each $f_n$ and $|f_n|=f_n$ is
hereditarily integrable, with
\[
 J_{f_n}(C)=n|C\cap I_n|
 \qquad(C\subseteq A\text{ an interval}).
\]
For completeness, only cells meeting an endpoint of $I_n$ can
contribute to the discrepancy between a tagged sum and this
formula. In a partition of mesh less than $\delta$, the total
length of those cells is at most $4\delta$, and their discrepancy
in absolute value is at most $4n\delta$. A constant local radius
$\delta/2$ makes every fine cell have length less than $\delta$.
Consequently, for every $\varepsilon>0$, the constant local radius
$\varepsilon/(16n)$ makes every full fine record have error at
most $\varepsilon/2$. The same estimate holds after restriction
to any interval $C$, and proves convergence along the full gauge
record nets.

At a tag $t\in I_m$, the value $f_n(t)$ is zero for $n\ne m$.
At $t=1$, every value is zero. Consequently,
\[
 k(t)=m\quad(t\in I_m),
 \qquad k(1)=1
\]
verifies \eqref{vt:local-barriers}. The input profiles can all be
zero, with threshold $m+1$ on $I_m$ and threshold $1$ at the
remaining tag. The tails have explicit primitives
\[
 J_{(f_n-g_k)^+}(A)
 =c\frac{(n-k)^+}{n^2}.
\]
If $n\le k$, this quantity is zero. If $n>k$, then
\[
 n^2-4k(n-k)=(n-2k)^2\ge0,
 \qquad
 c\frac{n-k}{n^2}\le\frac{c}{4k}.
\]
Equality holds at $n=2k$. Thus
\[
 \sup_{n\ge1}J_{(f_n-g_k)^+}(A)
 =\frac{c}{4k}=:v_k\downarrow0.
\]
All the hypotheses of Theorem~\ref{vt:vitali} hold, and in this
case its variation controls are
\[
 d_N=\sup_{n\ge N}J_{|f_n|}(A)
 =\sup_{n\ge N}\frac cn
 =\frac cN\downarrow0.
\]

Suppose that a positive integrable function $q$ dominated every
$f_n$. Disjointness of the intervals would give, for every $N$,
\[
 0\le\sum_{n=1}^{N}f_n\le q.
\]
Positivity and additivity of existing integrals would then imply
\[
 J_q(A)\ge\sum_{n=1}^{N}J_{f_n}(A)
 =c\sum_{n=1}^{N}\frac1n.
\]
For $N=2^m$, grouping the summands from $2^{r-1}+1$ through
$2^r$ gives
\[
 \sum_{n=1}^{2^m}\frac1n
 \ge1+\sum_{r=1}^{m}2^{r-1}\frac1{2^r}
 =1+\frac m2.
\]
No finite real number $J_q(A)$ satisfies these inequalities for
every $m$. Hence no positive integrable common envelope exists.
The increasing-barrier theorem therefore covers data excluded
by integrable domination, even in the real interval model.
\end{example}

Throughout the barrier construction, the assumptions concern pointwise
bounds and positive primitives of discarded parts of the given
functions, with the local bound verifiable through
\eqref{vt:data-local-barriers} and no integrability condition imposed on
the candidate limit. From those data we first obtain the integral of
$|f|$ by bounded primitives, then the integral of $f$, and finally the
uniform variation bound. This route leaves open a complementary question:
what follows when a common approximation estimate for the known
integrals is supplied directly?

\subsection{The equiintegrability implication in the original range}

\begin{theorem}[Sequential equiintegrability transfer]
\label{lc:equiintegrability}
Let $Z$ be separated and sequentially complete in its original uniform
convergence. Let symmetric lattice-valued pseudometrics
$d_j:Z^2\to(F_j)_+$ faithfully represent point and pair convergence
through regulator convergence, where every $F_j$ is Dedekind
$\sigma$-complete and WSD. Let $\Gamma$ be directed. Suppose
\[
 x_{n,\gamma}\longrightarrow y_n\quad(\gamma\in\Gamma),
 \qquad
 x_{n,\gamma}\longrightarrow x_\gamma\quad(n\to\infty).
\]
For each $j$, suppose a regulator $b^j$ satisfies
\begin{equation}\label{lc:equi-row-bound}
 \forall\varphi\ \exists\gamma(j,\varphi)\quad
 \forall\gamma\ge\gamma(j,\varphi)\ \forall n:
 \quad d_j(x_{n,\gamma},y_n)\le H_{b^j}(\varphi).
\end{equation}
Then a unique $y\in Z$ satisfies $y_n\to y$ and $x_\gamma\to y$.
Moreover,
\begin{equation}\label{lc:equi-candidate-bound}
 \gamma\ge\gamma(j,\varphi)
 \Longrightarrow d_j(x_\gamma,y)\le H_{b^j}(\varphi).
\end{equation}
\end{theorem}

\begin{proof}
Fix $j$, a selector $\varphi_0$, and
$\gamma_0\ge\gamma(j,\varphi_0)$. Faithfulness and WSD give
\[
 d_j(x_{n,\gamma_0},x_{\gamma_0})\longrightarrow_o0.
\]
An order-null sequence is order bounded. Indeed, one eventual order
bound, joined with the finitely many preceding values, bounds every
term. The triangle inequality gives
\[
 d_j(y_n,y_m)
 \le2H_{b^j}(\varphi_0)
 +d_j(x_{n,\gamma_0},x_{\gamma_0})
 +d_j(x_{\gamma_0},x_{m,\gamma_0}).
\]
Thus all these pair distances are order bounded. Dedekind
$\sigma$-completeness supplies
\[
 q_{j,N}=\sup_{n,m\ge N}d_j(y_n,y_m),
 \qquad q_j=\inf_Nq_{j,N}.
\]
For an arbitrary $\varphi$, fix
$\gamma_\varphi\ge\gamma(j,\varphi)$. The same column convergence
and the sequential upper-envelope lemma give
\[
 r_{j,N}^{\gamma_\varphi}
 =\sup_{n\ge N}d_j(x_{n,\gamma_\varphi},x_{\gamma_\varphi})
 \downarrow0.
\]
Equation~\eqref{lc:equi-row-bound} yields
\[
 q_{j,N}\le2H_{b^j}(\varphi)+2r_{j,N}^{\gamma_\varphi},
 \qquad
 0\le q_j\le2H_{b^j}(\varphi).
\]
The last inequality holds for every selector. WSD gives $q_j=0$.
For the regulator $a^j_{is}=q_{j,s}$ and each selector $\varphi$,
\[
 n,m\ge\varphi(1)\quad\Longrightarrow\quad
 d_j(y_n,y_m)\le q_{j,\varphi(1)}\le H_{a^j}(\varphi).
\]
Thus the pair net $(y_n,y_m)_{(n,m)\in\mathbb N^2}$ is null in
all represented coordinates. Pair faithfulness makes $(y_n)$
Cauchy in the original uniform convergence. Original sequential
completeness supplies $y\in Z$ with $y_n\to y$, and separation
makes this limit unique.

Fix $j,\varphi$ and $\gamma\ge\gamma(j,\varphi)$. The triangle
inequality gives, for every $n$,
\[
 d_j(x_\gamma,y)
 \le H_{b^j}(\varphi)+e_n,\qquad
 e_n=d_j(x_\gamma,x_{n,\gamma})+d_j(y_n,y).
\]
Faithfulness, symmetry, and regulator-to-order convergence give
$e_n\to_o0$. Let $v_\theta\downarrow0$ be order bounds for this
sequence, with $e_n\le v_\theta$ for $n\ge N(\theta)$. Then
\[
 d_j(x_\gamma,y)-H_{b^j}(\varphi)\le v_\theta
 \quad(\theta),
 \qquad
 d_j(x_\gamma,y)-H_{b^j}(\varphi)\le\inf_\theta v_\theta=0.
\]
This proves \eqref{lc:equi-candidate-bound}. For every $j$, the
same regulator $b^j$ now bounds $d_j(x_\gamma,y)$ on the tail
$\gamma\ge\gamma(j,\varphi)$ for every selector $\varphi$.
Point faithfulness gives $x_\gamma\to y$.
\end{proof}

For prescribed integral sums, the substitutions are
\[
 x_{n,\gamma}=S_\gamma f_n,
 \qquad y_n=J_{f_n},
 \qquad x_\gamma=S_\gamma f.
\]
Finiteness of each record allows pointwise original convergence and
continuity of each fixed coefficient map to give the required
fixed-column convergence. The estimate \eqref{lc:equi-row-bound}
concerns only known integrals, from which the preceding theorem derives
endpoint Cauchyness and constructs the candidate integral under its
stated hypotheses. In Theorem~\ref{lc:levi}, by comparison, the local
mixed approximation estimates themselves were obtained from the data and
bounded primitives.

\subsection{Verification of known-row conditions}

When neither an envelope nor the calibrated-controller construction is
available, the mixed Henstock estimate can still be tested through an
order-bounded series of known-row regulator bounds, or through an order
bound on the candidate finite sums relative to known primitives. We
verify these two sufficient tests next, keeping their additional
quantitative premises separate from the preceding theorems; neither
test involves an integral of the limit function. Throughout this
subsection, $Z=F$ is a Dedekind $\sigma$-complete weakly
$\sigma$-distributive abelian lattice-ordered group whose original point
and pair convergences are faithfully represented by modulus
$D$-convergence. Retain the established local division, completion, and
finite synchronization hypotheses, with positive additive coefficients
$L_C:E\to F$ that are additive on divisions and with $L_A$ sequentially
order continuous. All known rows are integrable on the root and its
compatible cells.

\begin{proposition}[Summable known-row bounds]\label{lc:source-summable}
Let $a^{(n)}=(a_{nis})$ regulate the known integral of $h_n$ on $A$. Put
\[
 u_n=\bigvee_{i\ge1}a_{ni1}.
\]
Suppose the partial sums $\sum_{n=1}^Nu_n$ are order bounded.
Then the local mixed
Henstock errors admit a single regulator, without an integrable envelope.
\end{proposition}
\begin{proof}
For the constant selector $\mathbf1$,
\[
 \sum_{n=1}^N\bigvee_i a_{ni,\mathbf1(i+n+1)}
 =\sum_{n=1}^Nu_n.
\]
Conversely, for every selector $\varphi$,
\[
 0\le\sum_{n=1}^N\bigvee_i a_{ni,\varphi(i+n+1)}
 \le\sum_{n=1}^Nu_n.
\]
Dedekind $\sigma$-completeness therefore defines
\[
 W=\bigvee_N\sum_{n=1}^Nu_n\in F_+.
\]
Fix $\varphi$. For each $n$, use the known-row gauge at the selector
$i\mapsto\varphi(i+n+1)$. Local finite synchronization, Lemma~\ref{lc:partial}, and
Proposition~\ref{fv:variation} give, for every fine record, compatible
subword $Q$, and labels
$1\le\nu_r\le p(t_r)$,
\[
 \sum_{r\in Q}|L_{C_r}h_{\nu_r}(t_r)-J_{\nu_r}(C_r)|
 \le 2\sum_{n\in M(Q)}H_{a^{(n)}}(i\mapsto\varphi(i+n+1))
 \le2W,
\]
where $M(Q)=\{\nu_r:r\in Q\}$ is finite.
Set $c^{(1)}=0$ and $c^{(n+1)}=2a^{(n)}$. Lemma~\ref{lt:fusion},
with cap $2W$, gives the regulator
\[
 B_{1s}=0,\qquad
 B_{rs}=2W\wedge
 2^r\sum_{1\le m<r}c^{(m)}_{r-m,s}\quad(r\ge2).
\]
The finite active-label sum is bounded by the corresponding initial
sum. Its minimum with $2W$ dominates the mixed variation. Consequently,
\[
 \sum_{r\in Q}|L_{C_r}h_{\nu_r}(t_r)-J_{\nu_r}(C_r)|
 \le H_B(\varphi).
\]
The cap, the regulator, and the estimate are independent of $p,Q$, and
the number of active labels.
\end{proof}

\begin{proposition}[A bounded-response criterion]
\label{lc:source-bounded-response}
Suppose $h_n$ is increasing, the known primitives are finitely additive,
and $h_n$ tends to $h$ with the prescribed common input regulator.
Suppose there is $z\in F_+$ such that
\begin{equation}\label{lc:literal-ii}
 |J_n(A)-S_Rh|\le z
 \qquad(n\ge1,\ R\text{ an original full record on }A).
\end{equation}
Then the primitive increments and the local mixed Henstock variations
are order bounded. Bounded fusion of the individual known-row
regulators consequently gives a regulator for the mixed errors.
\end{proposition}
\begin{proof}
Fix one original full record $R_0$. Positivity gives
\[
 0\le J_n(A)-J_1(A)
 \le |J_n(A)-S_{R_0}h|+|S_{R_0}h-J_1(A)|\le2z.
\]
Thus
\[
 K=\bigvee_n[J_n(A)-J_1(A)]\le2z.
\]
For the known integral of $h_1$, a fixed gauge $\eta$ and $u\in F_+$
satisfy
\[
 |S_Rh_1-J_1(A)|\le u\qquad(R\text{ $\eta$-fine}).
\]
Lemma~\ref{lc:partial} and Proposition~\ref{fv:variation} give
\[
 \sum_{r\in Q}|L_{C_r}h_1(t_r)-J_1(C_r)|\le2u.
\]
Let $Q$ be a subword of an $\eta$-fine full record $R$, with arbitrary
finite labels $\nu_r\ge1$, and put $M=\max_{r\in Q}\nu_r$ when
$Q\ne\varnothing$. Empty sums are zero. Positivity and additivity give
\[
 0\le\sum_{r\in Q}[J_{\nu_r}(C_r)-J_1(C_r)]\le K,
\]
\[
 0\le\sum_{r\in Q}L_{C_r}(h_{\nu_r}-h_1)(t_r)
 \le S_R(h_M-h_1).
\]
For the common input regulator $a$, Lemma~\ref{dt:transport} gives a
regulator $d$ with $L_AH_a(\varphi)\le H_d(\varphi)$. Put
\[
 N(R,\varphi)=\max_{(C,t)\in R}\theta(t,\varphi).
\]
Then
\[
 n\ge N(R,\varphi)\Longrightarrow
 |S_Rh_n-S_Rh|\le\sum_{(C,t)\in R}L_CH_a(\varphi)
 =L_AH_a(\varphi)\le H_d(\varphi).
\]
Thus $S_Rh_n\to_D S_Rh$, and WSD gives $S_Rh_n\to_o S_Rh$.
For fixed $M$ and $n\ge M$,
\[
 0\le(S_Rh_M-S_Rh)^+\le|S_Rh_n-S_Rh|.
\]
The fixed positive left side is zero. Therefore
\[
 S_R(h_M-h_1)
 \le S_Rh-S_Rh_1
 \le z+u.
\]
Consequently,
\begin{align*}
 \sum_{r\in Q}|L_{C_r}h_{\nu_r}(t_r)-J_{\nu_r}(C_r)|
 &\le\sum_{r\in Q}L_{C_r}(h_{\nu_r}-h_1)(t_r)
       +\sum_{r\in Q}[J_{\nu_r}(C_r)-J_1(C_r)]\\
 &\quad+\sum_{r\in Q}|L_{C_r}h_1(t_r)-J_1(C_r)|\\
 &\le K+z+3u.
\end{align*}
Set $W=K+z+3u$. For every $m\ge1$, let $a^{(m)}$ be a regulator
for the known root integral of $h_m$. For each selector $\psi$, its
root bound and the partial Henstock estimate give a gauge
$\eta_{m,\psi}$ for which
\[
 \sum_{(C,t)\in Q_m}|L_Ch_m(t)-J_m(C)|
 \le2H_{a^{(m)}}(\psi)
\]
whenever $Q_m$ is a compatible $\eta_{m,\psi}$-fine partial
record. Here the modulus pseudonorm preserves finite joins, so
Proposition~\ref{fv:variation} applies to all subwords of $Q_m$.
Put $c^{(m)}=2a^{(m)}$ and define
\[
 B_{1s}=0,\qquad
 B_{rs}=W\wedge2^r\sum_{1\le m<r}c^{(m)}_{r-m,s}
 \quad(r\ge2).
\]
Every entry belongs to $[0,W]$, and each row decreases to zero
because the sum has finitely many order-null rows. Thus $B$ is a
regulator. For $\varphi\in\mathbb N^{\mathbb N}$ and
$p:T\to\mathbb N_{\ge1}$, set
\[
 \psi_{m,\varphi}(i)=\varphi(i+m),\qquad
 \delta_{\varphi,p}(t)
 =\eta(t)\wedge\bigwedge_{1\le m\le p(t)}
                         \eta_{m,\psi_{m,\varphi}}(t).
\]
Let $Q=((C_r,t_r))$ be a compatible $\delta_{\varphi,p}$-fine
partial record, with labels $1\le\nu_r\le p(t_r)$. Since
$\delta_{\varphi,p}\le\eta$, Proposition~\ref{ls:completion}
extends $Q$ to an $\eta$-fine full record. The cap $W$ just proved
therefore applies to this $Q$. For nonempty $Q$, put
$M=\max_r\nu_r$. Its subword with label $m$ satisfies
$\delta_{\varphi,p}(t_r)\le\eta_{m,\psi_{m,\varphi}}(t_r)$,
so the preceding variation bound and the independently proved cap give
\[
 \begin{aligned}
 \sum_r|L_{C_r}h_{\nu_r}(t_r)-J_{\nu_r}(C_r)|
 &\le W\wedge\sum_{m=1}^{M}
                  H_{c^{(m)}}(\psi_{m,\varphi})\\
 &\le H_B(\varphi),
 \end{aligned}
\]
where the last inequality is Lemma~\ref{lt:fusion} with the
specified cap $W$ and the specified family $(c^{(m)})$. For empty
$Q$ the sum is zero and has the same bound. The regulator $B$ is
independent of $p$, the record, and its labels, as required.
\end{proof}

For comparison, the distinct known-row condition
\[
 |J_n(A)-S_Rh_n|\le z\qquad(n,R)
\]
also supplies a cap. If $K=\sup_n[J_n(A)-J_1(A)]<\infty$, then
\[
 S_R(h_M-h_1)\le K+2z,
 \qquad
 \sum_{r\in Q}|L_{C_r}h_1(t_r)-J_1(C_r)|\le2z,
\]
so the mixed variation is at most $2K+4z$. This condition is not
identified with \eqref{lc:literal-ii}.

\begin{corollary}[A monotone criterion from known-row bounds]\label{lc:source-monotone}
Let $h_n$ be increasing and tend to $h$ with a common input regulator.
Suppose its root primitives are order bounded above and the known-row
regulators satisfy Proposition~\ref{lc:source-summable}.
Then $h$ is hereditarily integrable and
\[
 J_h(C)=\bigvee_nJ_n(C),\qquad
 0\le J_h(C)-J_n(C)\le J_h(A)-J_n(A)\downarrow0.
\]
The convergence is uniform over compatible cells in the original
pair convergence. The alternative
\eqref{lc:literal-ii} gives the same conclusions.
These hypotheses give a local-division extension of the monotone
criterion in \cite[Theorem~5.4.2, p.~87]{integralmeasureandordering}.
\end{corollary}
\begin{proof}
Positivity and root additivity give
\[
 0\le J_n(C)-J_1(C)\le J_n(A)-J_1(A)\le K,
 \qquad K=\bigvee_n[J_n(A)-J_1(A)].
\]
Dedekind $\sigma$-completeness defines $J(C)=\bigvee_nJ_n(C)$.
For a finite division $(C_r)_{r=1}^m$ of $C$ and a tuple
$(n_1,\ldots,n_m)$, put $N=\max_rn_r$. Then
\[
 \sum_{r=1}^mJ_{n_r}(C_r)\le\sum_{r=1}^mJ_N(C_r)=J_N(C)\le J(C),
 \qquad J_n(C)=\sum_{r=1}^mJ_n(C_r)\le\sum_{r=1}^mJ(C_r).
\]
Taking the bounded countable suprema yields
\[
 J(C)=\sum_{r=1}^mJ(C_r),\qquad
 0\le J(C)-J_N(C)\le q_N:=J(A)-J_N(A)\downarrow0.
\]
Let $a$ be the given common input regulator. Lemma~\ref{dt:transport}
supplies a regulator $D$ with
\[
 L_AH_a(\varphi)\le H_D(\varphi),
 \qquad\psi_\varphi=\varphi.
\]
Let $\ell_\varphi(t)$ be the least eventual threshold for the
common-regulator inequality at $(t,\psi_\varphi)$, and set
\[
 p_\varphi(t)=\max\{\varphi(1),\ell_\varphi(t)\}.
\]
Proposition~\ref{lc:source-summable}, or
Proposition~\ref{lc:source-bounded-response}, supplies the mixed
regulator $B$. Apply its local finite synchronization with
$p=p_\varphi$. On a resulting fine record $R=((C_r,t_r))$ of $C$,
\begin{align*}
 S_Rh-J(C)
 &=\sum_rL_{C_r}(h-h_{p_\varphi(t_r)})(t_r)\\
 &\quad+\sum_r[L_{C_r}h_{p_\varphi(t_r)}(t_r)
                    -J_{p_\varphi(t_r)}(C_r)]\\
 &\quad-\sum_r[J(C_r)-J_{p_\varphi(t_r)}(C_r)].
\end{align*}
The three terms satisfy
\[
 \sum_rL_{C_r}|h-h_{p_\varphi(t_r)}|(t_r)
 \le L_CH_a(\psi_\varphi)\le L_AH_a(\psi_\varphi)
 \le H_D(\varphi),
\]
\[
 \sum_r|L_{C_r}h_{p_\varphi(t_r)}(t_r)-J_{p_\varphi(t_r)}(C_r)|
 \le H_B(\varphi),
\]
\[
 0\le\sum_r[J(C_r)-J_{p_\varphi(t_r)}(C_r)]
 \le J(C)-J_{\varphi(1)}(C)\le q_{\varphi(1)}.
\]
Thus $c_{is}=3(D_{is}+B_{is}+q_s)$ is a regulator and
\[
 |S_Rh-J(C)|\le H_D(\varphi)+H_B(\varphi)+q_{\varphi(1)}
 \le H_c(\varphi).
\]
Faithfulness proves $J_h(C)=J(C)$. The already established bound
$0\le J(C)-J_n(C)\le q_n$ and the rank-one regulator $b_{is}=q_s$
prove uniform convergence over compatible cells in the original pair
convergence.

For coefficients $L_Cx=x\bullet\mu(C)$, positivity and division
additivity hold for a positive additive set function and an additive
product. Preservation of bounded countable positive suprema gives
sequential order continuity by applying it to $x_1-x_k\uparrow x_1$. Under \eqref{lc:literal-ii},
Proposition~\ref{lc:source-bounded-response} also supplies $K\le2z$.
\end{proof}

\subsection{Regulator convergence and the ordered corollaries}
\label{hc:section}

We conclude by making these constructions explicit for the ordered
corollaries of \cite[Chapter~5]{integralmeasureandordering}. Let $E$ be
Dedekind $\sigma$-complete and let $F$ be Dedekind $\sigma$-complete
and WSD, with the original point and pair convergences of $F$ given by
modulus $D$-convergence. Each prescribed record is finite, with
coefficients
\[
 L_Cx=x\bullet\mu(C)
\]
that are positive, additive, and sequentially order continuous:
\[
 x_s\downarrow0\quad\Longrightarrow\quad L_Cx_s\downarrow0.
\]
Where finite lattice operations on integrands enter, we also impose
preservation of modulus; Proposition~\ref{nr:coefficient-order}
identifies its conjunction with sequential order continuity as
preservation of existing bounded countable positive suprema. All these
requirements concern the coefficient maps and the original convergence.

\begin{proposition}[Fixed-record transfer of a common input regulator]
\label{hc:fixed-record}
Under these coefficient conditions, if $f_n$ approximates $f$ with a
common input regulator, then, for every fixed original record,
\[
 S_Rf_n\longrightarrow_D S_Rf.
\]
No common approximation threshold over different records is required.
\end{proposition}
\begin{proof}
For an empty record both sums are zero. Otherwise write
$R=((C_r,t_r))_{r=1}^m$, $m\ge1$, and $L_rx=x\bullet\mu(C_r)$.
Let $M\in E_+$ bound the input regulator $a$. For each $r$,
Lemma~\ref{lt:positive-d-transport} gives the regulator
\[
 b^r_{is}=L_rM\wedge2^iL_ra_{is},\qquad
 L_rH_a(\varphi)\le H_{b^r}(\varphi).
\]
Its rows decrease to zero by sequential order continuity, and its entries
are bounded by $L_rM$. Define
\[
 c_{is}=m\sum_{r=1}^mb^r_{is},\qquad
 N(R,\varphi)=\max_{1\le r\le m}\theta(t_r,\varphi).
\]
Then $0\le c_{is}\le m\sum_rL_rM$ and $c_{is}\downarrow_s0$ for each
$i$, so $c$ is a regulator. For $n\ge N(R,\varphi)$,
\[
\begin{aligned}
 |S_Rf_n-S_Rf|
 &\le\sum_{r=1}^m|L_r(f_n(t_r)-f(t_r))|\\
 &\le\sum_{r=1}^mL_r|f_n(t_r)-f(t_r)|\\
 &\le\sum_{r=1}^mL_rH_a(\varphi)
 \le\sum_{r=1}^mH_{b^r}(\varphi).
\end{aligned}
\]
Put $C_\varphi=\sup_i\sum_{r=1}^mb^r_{i,\varphi(i)}$.
For every $r$, $H_{b^r}(\varphi)\le C_\varphi$. Consequently
\[
 |S_Rf_n-S_Rf|\le mC_\varphi=H_c(\varphi).
\]
This is the required regulator convergence. The proof uses the positive
inequality $|L_rx|\le L_r|x|$, not modulus preservation.
\end{proof}

\begin{corollary}[Known-row uniform approximation]
\label{hc:old-bounded}\label{secondmastertheorem}\label{hc:ordered-specializations}
Let $f_n$ have original net Riemann integrals on a domain $A$, and let
$f_n$ approximate $f$ with a common input regulator. Suppose a regulator
$b$ in $F$ satisfies
\begin{equation}\label{hc:known-row}
 \forall\varphi\quad\exists\gamma_0(A,\varphi)\quad
 \forall\gamma\ge\gamma_0(A,\varphi)\quad\forall n:\qquad
 |J_{f_n}(A)-S_\gamma^Af_n|\le H_b(\varphi).
\end{equation}
Then the original net Riemann integral of $f$ on $A$ exists, and
\[
 J_{f_n}(A)\longrightarrow_DJ_f(A),\qquad
 J_{f_n}(A)\longrightarrow_oJ_f(A).
\]
If $f_n\le f_{n+1}$, then $J_f(A)=\sup_nJ_{f_n}(A)$.
\end{corollary}
\begin{proof}
Proposition~\ref{lt:d-complete} gives sequential completeness of $F$
in its given regulator convergence. Proposition~\ref{hc:fixed-record}
gives convergence of every fixed column. In
Theorem~\ref{lc:equiintegrability}, put
\[
 x_{n,\gamma}=S_\gamma^Af_n,\qquad y_n=J_{f_n}(A),\qquad
 x_\gamma=S_\gamma^Af,\qquad d(z,w)=|z-w|.
\]
Known-row integrability gives the row limits. Equation~\eqref{hc:known-row}
gives the uniform row bound. The theorem constructs a unique $J_f(A)$
and gives convergence of both the candidate sums and the known endpoints.
WSD converts endpoint regulator convergence into order convergence.

For $f_n\le f_{n+1}$, positivity gives
$J_{f_n}(A)\le J_{f_m}(A)$ for $m\ge n$. Passage to the order limit in
$m$ yields $J_{f_n}(A)\le J_f(A)$. If $v\ge J_{f_n}(A)$ for every $n$,
then $v-J_{f_n}(A)\to_o v-J_f(A)$ and $v-J_f(A)\ge0$. Therefore
$J_f(A)=\sup_nJ_{f_n}(A)$.
\end{proof}

So far, the uniform approximation estimate for the known endpoints
$J_{f_n}$ in Corollary~\ref{hc:old-bounded} has supplied one route to
the ordered limit theorem. To recover the monotone and dominated
principles associated with \cite[Chapter~5]{integralmeasureandordering}
from bounded data and local divisions, we now verify the structural
conditions that allow their mixed estimates to be derived.

\begin{proposition}[Coefficient verification for the positive product]
\label{hc:source-coefficients}
For positive additive coefficients, additivity on full divisions gives
exact integration of constants. If a coefficient additionally preserves
existing bounded countable positive suprema, it preserves the modulus. In the full-division
coefficient calculus,
\[
 J_{\underline u}(C)=L_Cu=u\bullet\mu(C)
 \qquad(u\in E_+).
\]
\end{proposition}
\begin{proof}
Under the additional supremum hypothesis,
Proposition~\ref{nr:coefficient-order} gives $L_C|x|=|L_Cx|$.
Independently of that hypothesis, every full record of $C$ satisfies
\[
 S_R\underline u=\sum_{D\in\operatorname{div}R}L_Du=L_Cu.
\]
The response net is constant, so separation identifies its integral
with $L_Cu$.
\end{proof}

\begin{proposition}[Finite lattice stability of regulator convergence]
\label{hc:lattice-continuity}
Modulus $D$-convergence is preserved by the modulus and by finite
lattice suprema and infima.
\end{proposition}
\begin{proof}
The lattice identities give
\[
 \bigl||x|-|y|\bigr|\le|x-y|,
 \qquad
 |(x\vee y)-(x'\vee y')|\le|x-x'|+|y-y'|.
\]
If regulators $a,b$ control the two positive errors, then
\[
 H_a(\varphi)+H_b(\varphi)\le H_{2(a+b)}(\varphi).
\]
More explicitly, suppose $x^r_\alpha\to_D x^r$ for
$r=1,\ldots,m$, with regulators $a^r$ on a common directed index.
For each $\varphi$, let $\alpha_r(\varphi)$ attain its $r$th
bound and let $\alpha(\varphi)$ be an upper bound of those
finitely many indices. Repeated application of the displayed
lattice inequality gives
\[
 \left|\bigvee_{r=1}^{m}x^r_\alpha-\bigvee_{r=1}^{m}x^r\right|
 \le\sum_{r=1}^{m}|x^r_\alpha-x^r|
 \le\sum_{r=1}^{m}H_{a^r}(\varphi)
 \quad(\alpha\ge\alpha(\varphi)).
\]
The bounded matrix $c_{is}=m\sum_{r=1}^{m}a^r_{is}$ has rows
decreasing to zero, and
\[
 \sum_{r=1}^{m}H_{a^r}(\varphi)
 \le m\sup_i\sum_{r=1}^{m}a^r_{i,\varphi(i)}
 =H_c(\varphi).
\]
This proves the finite-supremum assertion with an explicit regulator.
For independently indexed nets, use their finite product index
and the tuple $(\alpha_r(\varphi))_{r=1}^{m}$ as the threshold.
Symmetry preserves the same bounds under negation, and
\[
 \bigwedge_{r=1}^{m}x^r_\alpha
 =-\bigvee_{r=1}^{m}(-x^r_\alpha)
\]
proves the infimum assertion. The modulus inequality retains the
regulator of its input difference without alteration.
\end{proof}

These continuity estimates describe how lattice operations behave
under convergence; to pass from them to integrability of pointwise
finite lattice combinations, we use the finite-lattice-closure result
proved above, including its additional argument for the integral.

\begin{corollary}[Bounded monotone and dominated data without a row estimate]
\label{hc:bounded-geometric}
Let $E$ and $F$ be Dedekind $\sigma$-complete WSD lattice-ordered groups.
Let the local division structure satisfy the branching condition.
The integration nets are the full record nets, or prescribed nets
satisfying \eqref{ls:eq-3}. Let the positive additive coefficients be
additive on divisions and let the root coefficient $L_A$ be sequentially
order continuous. Equip $F$ with its original
modulus regulator convergence. Suppose $f_n$ are hereditarily integrable
and approximate $f$ with a common input regulator. Then:
\begin{enumerate}
\item if $f_n\le f_{n+1}$ and $|f_1|\le u_1$, $|f|\le u_f$, then
$f$ is hereditarily integrable and
\[
 J_f(C)=\sup_nJ_{f_n}(C),\qquad
 0\le J_f(C)-J_{f_n}(C)\le J_f(A)-J_{f_n}(A)\downarrow0.
\]
\item if, in addition, $L_C|x|=|L_Cx|$ for every compatible cell,
$|f_n|\le h\le u$ for some $u\in E_+$, and either geometric alternative
of Theorem~\ref{fv:finite-lattice-closure} holds, then $f$
is hereditarily integrable and there is $d_N\downarrow0$ such that
\[
 |J_{f_n}(C)-J_f(C)|\le d_N
 \qquad(n\ge N,\ C\in\mathscr C_A).
\]
\end{enumerate}
Both endpoint limits are uniform in the original pair convergence.
\end{corollary}
\begin{proof}
For $\rho(z)=|z|$, point and pair faithfulness are the defining regulator
convergences. Every solid or finite-supremum saturation preserves an
attained bound $|z|\le H_b(\varphi)$. Local solidity and finite-supremum
stability therefore hold. For $z_s\downarrow0$, the regulator
$b_{is}=z_s$ gives sequential order compatibility.
Root normality is assumed. For a compatible cell $C$, positivity and
coefficient additivity give
\[
 0\le L_Cx_s\le L_Ax_s\downarrow0
 \quad(x_s\downarrow0),
\]
so each compatible coefficient is sequentially order continuous.
For positive constants, every full response equals $L_Cu$, by
Proposition~\ref{hc:source-coefficients}. No modulus identity is used
in the monotone assertion.

Let $a$ be the common input regulator. WSD in $E$ gives
$\inf_\varphi H_a(\varphi)=0$. For the increasing sequence, for each
$n,t,\varphi$ and every $m\ge\max\{n,\theta(t,\varphi)\}$,
\[
 (f_n(t)-f(t))^+\le|f_m(t)-f(t)|\le H_a(\varphi).
\]
Thus $f_n(t)\le f(t)$. If $v\ge f_n(t)$ for all $n$, then for
$n\ge\theta(t,\varphi)$,
\[
 (f(t)-v)^+\le|f(t)-f_n(t)|\le H_a(\varphi).
\]
Hence $f(t)\le v$, proving $f_n(t)\uparrow f(t)$. The bound
\[
 |f_n(t)|\le u_1\vee u_f
\]
gives the hereditarily integrable constant envelope
\[
 g=\underline{u_1\vee u_f},\qquad
 J_g(C)=L_C(u_1\vee u_f),\qquad |f_n|\le g.
\]
The substitution into Theorem~\ref{fv:envelope-monotone} is
\[
 h_n=f_n,\quad h=f,\quad Z=F,\quad \rho=|\cdot|,
 \quad G(C)=L_C(u_1\vee u_f).
\]
The common input regulator is the given $a$. Root normality and the
faithful locally solid representation were verified above. The theorem
therefore gives
\[
 J_f(C)=\sup_nJ_{f_n}(C),\qquad
 0\le J_f(C)-J_{f_n}(C)
 \le J_f(A)-J_{f_n}(A)\downarrow0.
\]

For the dominated case, the additional modulus assumption gives
$L_C|x|=|L_Cx|$ for every compatible coefficient, and
\[
 g=\underline u,\qquad S_Rg=L_Cu\quad(R\in\mathscr R_C),
 \qquad J_g(C)=L_Cu,
\]
\[
 |f_n|\le h\le g,\qquad
 \inf_\varphi H_a(\varphi)=0,
 \qquad L_C|x|=|L_Cx|.
\]
Finite-supremum stability of the original $D$-convergence was verified
at the start of this proof. The stated division hypothesis supplies
one of the two alternatives of Theorem~\ref{fv:finite-lattice-closure}.
Thus all premises of Theorem~\ref{fv:dominated} hold for
$(E,Z,\lambda,\rho,g,a)=(E,F,D,|\cdot|,\underline u,a)$.
Its construction gives
\[
 d_N\downarrow0,\qquad
 |J_{f_n}(C)-J_f(C)|\le d_N
 \quad(n\ge N,\ C\in\mathscr C_A).
\]
Pair faithfulness identifies both endpoint conclusions with the
original uniform convergence. For interval divisions, the independent
geometric verification is Proposition~\ref{cr:intervals}.
\end{proof}

\section{Conclusion}
\label{sec:conclusion}

By following finite records through the prescribed convergence of
their values, we have examined the integration procedures considered
here within a common net Riemann construction, in which the
coefficient product determines the elementary identities and the
convergence determines their passage to limits. This formulation has
led from additivity and subset integration to a series of convergence
theorems, retaining the original approximation nets and making
explicit the additional hypotheses used when group or lattice
operations become available in the range. Through the ordered and
gauge specializations, these arguments connect with the theories of \citet{sionsemigroup},
\citet{henstockmagnum} and \citet{integralmeasureandordering},
and the monoid formulations identify those parts of the
construction that require no subtraction.

Returning to the question raised by uniform approximation, we recover
in the abelian group setting, from elementary sampling systems,
precisely the target conditions that suffice for passage
from finite comparisons to their limits. Adherence stability governs
the passage to known endpoints, as the complete locally convex
obstruction demonstrates, while construction of the limiting integral
additionally uses completeness and absorption of the independently
convergent row errors. Because each property can be restricted to
the cardinality of the available outer error family, countable error
control requires only sequential completeness even when the original
approximation and integration nets have arbitrary size. This cardinal
distinction persists in the noncommutative sufficient theorem through
the right convergence structure, and the support examples show that
the resulting range need not contain a complete subspace carrying an
endpoint tail.

To understand how these intrinsic conditions enter the abstract
theorems, we have returned to their original bounds: finite
intersections and removal of a null added bound supply ordered
adherence and absorption, while bounded fusion combines the unrelated
regulators of the countable regulator argument. Alongside these
verifications, the finite coefficient estimates are checked
independently, and the normalization route retains its row-dependent
residuals where a global response bound is unavailable. The familiar
topological, ordered and lattice-normed consequences thereby retain
their given convergence, with the proofs identifying the operations
that make the passage to limits possible.

Passing from uniform approximation to the local results requires
further information about divisions and lattice operations, which
the Henstock estimates supply before monotone and dominated limits
are constructed. Once these local estimates are available, increasing
integrable barriers extend the approximation beyond a common
integrable envelope by controlling the omitted parts through their
excess primitives, with the controller selection hypotheses retained
at the construction of the limiting modulus. Through the resulting
primitive-variation estimates and the disjoint spike example, this
extension is distinguished from the envelope theorem while the
weaker assumptions of that earlier result are preserved.

Within the representation theory of the appendices, we obtain further
verifications and comparisons of these convergence properties, while
its examples distinguish the existence of an exact representation
from the additional stability needed for uniform passage. These
distinctions will also guide the broader comparison developed in
\emph{A Galaxy of Integrals} \citep{TeixeiraGalaxy}, where the
same prescribed finite sums can be followed through different
integral constructions and their convergence theories. The
characterization established here supplies an intrinsic reference
for that comparison, determining which operations on their error
families account for a universal uniform theorem and which
additional data belong to the individual procedure.

\appendix

\section{Order convergence, regulators, and countable selection}
\label{app:order}
\label{lt:order-section}

To follow the order estimates used in the integration theorems, we
collect the control and selection arguments in abelian lattice-ordered
groups, where integer multiplication is understood as repeated addition
and no scalar field is needed. For such a group $F$, write
\[
 F_+=\{x\in F:x\ge0\},\qquad x^+=x\vee0,\qquad
 x^-=(-x)\vee0,\qquad |x|=x\vee(-x).
\]
We use the group form of the standard lattice identities treated in
\cite{vulikh,zaanenluxemburg1,zaanen2}, under which translation and
inversion preserve the corresponding existing order suprema and infima,
and finite lattice operations preserve existing bounded increasing
suprema and decreasing infima.

\subsection{Order controls and the distinction between nets and sequences}

\begin{definition}[Order completeness]
\label{dedekindcomplete}
\label{app:completeness}
The group $F$ is Dedekind $\sigma$-complete if every nonempty countable
subset bounded above has a supremum. It is Dedekind complete if the
same assertion holds for arbitrary nonempty bounded-above subsets.
It has the countable sup property if each nonempty bounded-above subset
whose supremum exists has a nonempty countable subset with the same
supremum. It is super Dedekind complete if it is Dedekind complete and
has the countable sup property. Finite subsets count as countable.
It is laterally $\sigma$-complete if every countable disjoint family in
$F_+$ has a supremum, without an a priori order bound.
\end{definition}

A net is order bounded when all of its values belong to some interval
$[-u,u]$, $u\in F_+$. Eventual order boundedness requires this only on a
tail. For a function $h:T\to F$, $|h|$ denotes the pointwise modulus.

\begin{definition}[Order convergence]
\label{usualorderconvergence}
\label{pavlakosorderconvergence}
A net $(x_\alpha)_{\alpha\in I}$ order converges to $x\in F$, written
$x_\alpha\xrightarrow{o}x$, if there is a decreasing net
$(u_\theta)_{\theta\in\Theta}$ in $F_+$ such that
\[
 \inf_{\theta\in\Theta}u_\theta=0,\qquad
 \forall\theta\ \exists\alpha_\theta\ \forall\alpha\ge\alpha_\theta:
 \quad |x_\alpha-x|\le u_\theta.
\]
A decreasing sequence of controls gives \emph{sequential order control},
even when the original net has an arbitrary directed index set.
\end{definition}

\begin{definition}[Order controls on tail values]\label{lt:order-null}
For a net $(x_i)_{i\in I}$ in an abelian lattice-ordered group $F$,
order nullity is equivalent to
\begin{equation}\label{lt:order-tail}
 \exists u_\theta\downarrow0\quad
 \forall\theta\ \exists i_\theta\quad
 T_{i_\theta}(x)\subseteq[-u_\theta,u_\theta].
\end{equation}
For a positive-valued map $q:Y\to F_+$ and a net $(y_i)$ in $Y$,
this reads
\[
 q(y_i)\xrightarrow{o}0
 \iff
 \exists u_\theta\downarrow0\quad
 \forall\theta\ \exists i_\theta\ \forall i\ge i_\theta:
 \quad q(y_i)\le u_\theta.
\]
The notation $u_\theta\downarrow0$ requires an existing infimum equal
to zero. Dedekind completeness is not part of this definition.
For a directed family $(B_\lambda)$ of nonempty subsets, its saturation
is order null exactly when, for every $\theta$, some $B_\lambda$ is
contained in $[-u_\theta,u_\theta]$.
\end{definition}

If $u_\theta\downarrow0$ and $v_\eta\downarrow0$, then
\begin{equation}\label{lt:sum-controls}
 (u_\theta+v_\eta)_{(\theta,\eta)}\downarrow0.
\end{equation}
For an arbitrary common lower bound $h$, the fixed-$\eta$ inequalities give
\[
 h-v_\eta\le u_\theta\quad(\theta,\eta)
 \ \Longrightarrow\
 h-v_\eta\le0\quad(\eta)
 \ \Longrightarrow\ h\le0.
\]
Since the sum family is downward directed on the product of its
control index sets and zero is a lower bound, this proves
\eqref{lt:sum-controls}; the calculation requires only a partially
ordered abelian group, provided both infima are taken in that whole
group. It follows that finite sums of positive order-null nets on their
product index are order null, while a common upper index for finitely
many thresholds in \eqref{lt:order-tail} supplies the simultaneous
finite estimates used below.

This use of a product index also reconciles the control formulation
with the two-sided definition, in which the lower and upper bounding
nets may initially have independent index sets. Placing them on their
product gives, in a lattice-ordered group, the condition
\[
 a_\gamma\uparrow x,\quad b_\eta\downarrow x,\quad
 a_\gamma\le x_\alpha\le b_\eta\ \text{eventually for each }(\gamma,\eta).
\]
It is equivalent to Definition~\ref{usualorderconvergence}: one direction
uses $a_\theta=x-u_\theta$, $b_\theta=x+u_\theta$, while the other uses
$u_{\gamma,\eta}=(x-a_\gamma)\vee(b_\eta-x)\downarrow0$.
These conventions agree with the order-convergence terminology of
\cite{papangelouorder,prerieszspaces,convergencestructuresvanderwalt}.

\begin{proposition}[Tail extrema]
\label{orderconvergenceboundednets}
Let $F$ be Dedekind complete. An eventually order-bounded net
$(x_\alpha)$ order converges to $x$ if and only if, on any order-bounded
tail,
\[
 \inf_\alpha\sup_{\beta\ge\alpha}x_\beta
 =x=
 \sup_\alpha\inf_{\beta\ge\alpha}x_\beta.
\]
For sequences, Dedekind $\sigma$-completeness suffices.
\end{proposition}
\begin{proof}
Restrict the index set to a tail on which $|x_\beta|\le v$.
Then
\[
 l_\alpha=\inf_{\beta\ge\alpha}x_\beta,\qquad
 u_\alpha=\sup_{\beta\ge\alpha}x_\beta
\]
exist, $l_\alpha$ increases, $u_\alpha$ decreases, and
$-v\le l_\alpha\le u_\alpha\le v$.
Suppose $x_\beta\to_o x$ with controls $r_\theta\downarrow0$.
For each $\theta$ there is $\alpha_\theta$ such that
\[
 \beta\ge\alpha_\theta
 \quad\Longrightarrow\quad
 x-r_\theta\le x_\beta\le x+r_\theta.
\]
Consequently
\[
 x-r_\theta\le l_{\alpha_\theta}
 \le\sup_\alpha l_\alpha
 \le\inf_\alpha u_\alpha
 \le u_{\alpha_\theta}\le x+r_\theta.
\]
Taking the infimum of $r_\theta$ gives both asserted equalities.
Conversely, those equalities give
$x-l_\alpha\downarrow0$ and $u_\alpha-x\downarrow0$.
Their sum on the diagonal is downward cofinal in their product sum,
so \eqref{lt:sum-controls} gives
\[
 w_\alpha=u_\alpha-l_\alpha\downarrow0.
\]
For $\beta\ge\alpha$,
\[
 l_\alpha\le x_\beta,x\le u_\alpha,
 \qquad |x_\beta-x|\le w_\alpha.
\]
These are order controls for the original tail, hence for the original
net. For a sequence, every supremum and infimum used above is countable.
\end{proof}

Restricting the extrema to a bounded tail is essential here, because
order convergence of a net need not bound the values preceding that
tail; for sequences the omitted initial segment is finite, so this
difficulty disappears. This distinction also underlies the bounded-tail
convention used in \cite{prerieszspaces,monteirouniform}.

\begin{lemma}[Fixed bounds]
\label{lt:fixed-bound}
\label{lt:order-bound}
Let $F$ be a partially ordered abelian group, and let
$r_\alpha\in F_+$ have decreasing nonnegative controls whose infimum
in the whole group is zero. If
$c\le b+r_\alpha$ eventually, then $c\le b$.
\end{lemma}
\begin{proof}
For each order control $u_\theta$, a common eventual index gives
\[
 c-b\le r_\alpha\le u_\theta.
\]
Hence $c-b\le\inf_F\{u_\theta:\theta\}=0$. The conclusion uses
the infimum in the whole group, since $c-b$ need not be positive.
\end{proof}

\begin{lemma}[Attained order bounds]
\label{lt:attained}
Let $e_\alpha\in F_+$. Suppose $\mathscr B\subset F_+$ is nonempty,
$\inf\mathscr B=0$, and $e_\alpha\le b$ eventually for every
$b\in\mathscr B$. Then $e_\alpha\xrightarrow{o}0$.
\end{lemma}
\begin{proof}
For a nonempty finite subset $K\subset\mathscr B$, put
\[
 u_K=\bigwedge_{b\in K}b.
\]
The finite-subset index is directed by inclusion, and
\[
 K\subset K'\Longrightarrow u_{K'}\le u_K,\qquad
 \inf_Ku_K=0.
\]
The finitely many eventual bounds indexed by $K$ hold simultaneously.
Therefore $e_\alpha\le u_K$ eventually, and the decreasing net
$(u_K)_K$ supplies the required controls.
\end{proof}

\begin{lemma}[Sequential upper envelopes]
\label{lt:sequence-envelope}
\label{lt:upper-envelope}
If $F$ is Dedekind $\sigma$-complete and $b_n\in F_+$ is order null,
then
\[
 v_N=\sup_{n\ge N}b_n\in F_+,\qquad v_N\downarrow0.
\]
\end{lemma}
\begin{proof}
Let $u_\theta\downarrow0$ control $(b_n)$. An attained bound supplies
$\theta_0,N_0$ with $b_n\le u_{\theta_0}$ for $n\ge N_0$. Thus
\[
 b_n\le u_{\theta_0}\vee\bigvee_{1\le k<N_0}b_k\qquad(n\ge1),
\]
where an empty finite join is zero. Countable completeness supplies
$v_N$ and $v=\inf_Nv_N$. For every $\theta$,
\[
 n\ge N_\theta\Longrightarrow b_n\le u_\theta,
 \qquad 0\le v\le v_{N_\theta}\le u_\theta.
\]
Hence $v=0$.
\end{proof}

For sequences in a Dedekind $\sigma$-complete group, we have therefore
obtained a decreasing sequential control from order convergence. The
restriction to sequences matters, since the argument does not extend
that conclusion to arbitrary nets.

\subsection{Regulator convergence and completeness}

\begin{definition}[Regulators and weak $\sigma$-distributivity]
\label{app:regulators}
\label{dconvergence}
Let $F$ be Dedekind $\sigma$-complete. A regulator is a double sequence
$a=(a_{rs})_{r,s\ge1}$ for which
\[
 0\le a_{rs}\le M\quad(r,s\ge1),\qquad a_{rs}\downarrow_s0
 \quad(r\ge1)
\]
for some $M\in F_+$. Its selector bound is
\[
 H_a(\varphi)=\sup_{r\ge1}a_{r,\varphi(r)},\qquad
 \varphi\in\mathbb N^{\mathbb N}.
\]
The net $x_\alpha$ is regulator convergent to $x$, written
$x_\alpha\xrightarrow{D}x$, when some regulator $a$ satisfies
\[
 \forall\varphi\ \exists\alpha_\varphi\ \forall\alpha\ge\alpha_\varphi:
 \quad |x_\alpha-x|\le H_a(\varphi).
\]
The group is weakly $\sigma$-distributive when
\[
 \inf_{\varphi\in\mathbb N^{\mathbb N}}H_a(\varphi)=0
\]
for every regulator. For the saturation of a directed family of sets,
the selector inequality means that one member of the family is contained
in the corresponding order interval.
\end{definition}

Requiring one order bound for the whole double sequence ensures that
every selector supremum in this definition exists, and allows the
resulting bounds to be used in the control calculus below. These
conventions and their classical vector-lattice forms are treated in
\cite{boccutovrabelovariecan,integralmeasureandordering}.

\begin{lemma}[Regulator calculus]
\label{lt:d-calculus}
In a weakly $\sigma$-distributive controller, every $D$-null net is order
null. Finite sums of positive $D$-null nets on the same directed set
are $D$-null. Independent index sets are handled by their finite product.
\end{lemma}
\begin{proof}
For a regulator $a$, the attained bounds $H_a(\varphi)$ have infimum
zero. Lemma~\ref{lt:attained} gives order nullity.
For regulators $a,b$, put $c_{rs}=2(a_{rs}+b_{rs})$. This is a bounded
regulator. For fixed $\varphi$, let $C=H_{a+b}(\varphi)$. Then
\[
 a_{r,\varphi(r)}\le C,\qquad b_{t,\varphi(t)}\le C,
\]
\[
 H_a(\varphi)+H_b(\varphi)\le2C=H_c(\varphi).
\]
For positive nets $x_\alpha,y_\alpha$ controlled by $a,b$, let
$\alpha_a(\varphi),\alpha_b(\varphi)$ be their selector thresholds,
and take an upper index $\alpha_0$ of these two indices. Then
\[
 \alpha\ge\alpha_0\quad\Longrightarrow\quad
 x_\alpha+y_\alpha\le H_a(\varphi)+H_b(\varphi)
                         \le H_c(\varphi).
\]
For $m$ summands with regulators $a^{(1)},\ldots,a^{(m)}$, define
$c^{(1)}=a^{(1)}$ and
\[
 c^{(k+1)}_{rs}=2\bigl(c^{(k)}_{rs}+a^{(k+1)}_{rs}\bigr)
       \qquad(1\le k<m).
\]
Induction gives boundedness and decreasing zero rows for every
$c^{(k)}$, and the binary estimate gives
$\sum_{\ell=1}^k x^{(\ell)}_\alpha\le H_{c^{(k)}}(\varphi)$
beyond a common upper index of its finitely many thresholds.
For independent index sets, the threshold is the tuple of those
indices in their finite product.
\end{proof}

\begin{proposition}[Order and regulator convergence]
\label{lt:order-d}
\label{(D)ando}
\label{equivalenceoand(D)convergence}
Let $F$ be Dedekind $\sigma$-complete and weakly $\sigma$-distributive.
Order and regulator convergence agree for sequences. They agree for
all nets if every principal order interval
has the countable sup property. In particular, the latter assertion
holds in a super Dedekind complete group.
\end{proposition}
\begin{proof}
Lemma~\ref{lt:d-calculus} proves $D\Rightarrow o$. For an order-null
sequence $(e_n)$ in $F_+$, Lemma~\ref{lt:sequence-envelope} supplies
$v_s\downarrow0$ with $e_n\le v_n$. The regulator $a_{rs}=v_s$ gives
\[
 n\ge\varphi(1)\Longrightarrow e_n\le v_{\varphi(1)}
 \le H_a(\varphi).
\]

For an arbitrary order-null net, let $u_\theta\downarrow0$
be its controls. Restrict to a tail bounded by one $u_{\theta_0}$.
The countable sup property applied to
$\{u_{\theta_0}-u_\theta:\theta\ge\theta_0\}\subset[0,u_{\theta_0}]$
gives indices $\theta_n$ with $\inf_nu_{\theta_n}=0$. Put
\[
 v_s=\bigwedge_{n\le s}u_{\theta_n},\qquad a_{rs}=v_s.
\]
Every $v_s$ is attained beyond a common upper index of the finitely
many thresholds corresponding to $u_{\theta_1},\ldots,u_{\theta_s}$. Thus
\[
 e\le v_{\varphi(1)}\le H_a(\varphi)
\]
is attained for every $\varphi$.
\end{proof}

Within this comparison, weak $\sigma$-distributivity alone gives
$D\Rightarrow o$, whereas the reverse implication for arbitrary nets
uses the stated countability condition. Keeping these two steps distinct
prevents the sequential equivalence from concealing the additional
selection needed for arbitrary nets.

\begin{proposition}[Controller completeness]\label{lt:d-complete}
A Dedekind $\sigma$-complete, weakly $\sigma$-distributive abelian
lattice-ordered group is sequentially complete for $D$-convergence.
\end{proposition}
\begin{proof}
Let $(z_n)$ be $D$-Cauchy with regulator $a$. The selector
$\iota(r)=1$ provides $N_0$ such that
\[
 |z_n-z_m|\le H_a(\iota)\quad(n,m\ge N_0).
\]
Hence
\[
 z_{N_0}-H_a(\iota)\le z_n\le z_{N_0}+H_a(\iota)
 \quad(n\ge N_0).
\]
For $N\ge N_0$, define
\[
 U_N=\sup_{n\ge N}z_n,\qquad L_N=\inf_{n\ge N}z_n,
 \qquad U=\inf_{N\ge N_0}U_N,\quad
 L=\sup_{N\ge N_0}L_N.
\]
All four extrema exist by Dedekind $\sigma$-completeness and the
displayed bounds.
For every $\varphi$, there is $N_\varphi\ge N_0$ such that
\[
 z_n-z_m\le H_a(\varphi)\quad(n,m\ge N_\varphi).
\]
For fixed $m\ge N_\varphi$, taking the supremum over
$n\ge N_\varphi$ gives
$U_{N_\varphi}\le z_m+H_a(\varphi)$. Translation preserves infima,
so the infimum over $m\ge N_\varphi$ gives
\[
 U_{N_\varphi}
 \le\inf_{m\ge N_\varphi}(z_m+H_a(\varphi))
 =L_{N_\varphi}+H_a(\varphi),
 \qquad U_{N_\varphi}-L_{N_\varphi}\le H_a(\varphi).
\]
Thus
\[
 0\le U-L\le H_a(\varphi)
 \quad(\varphi\in\mathbb N^{\mathbb N}),
 \qquad U=L=:z.
\]
For $n\ge N_\varphi$,
\[
 L_{N_\varphi}\le z,z_n\le U_{N_\varphi},
 \qquad
 |z_n-z|\le U_{N_\varphi}-L_{N_\varphi}\le H_a(\varphi).
\]
Therefore $z_n\xrightarrow D z$.
\end{proof}

\subsection{Bounded fusion of unrelated regulators}

Once each approximation has its own regulator, a further argument is
needed to bring their estimates under common control. We use the
classical regulator method of
\cite{monteirouniform,integralmeasureandordering}, retaining a fixed
bound on the fused array while allowing the input regulators to have
unrelated bounds.

\begin{lemma}[Bounded fusion]\label{lt:fusion}
Let $a^{(n)}$, $n\ge1$, be regulators in a Dedekind $\sigma$-complete
abelian lattice-ordered group $F$. For $w\in F^+$, define
\begin{equation}
 b_{1s}=0,\qquad
 b_{rs}=2^r\sum_{1\le n<r}a^{(n)}_{r-n,s}\quad(r\ge2),
 \qquad B_{rs}=w\wedge b_{rs}.
 \label{app:fusion-array}
\end{equation}
Then $B$ is a regulator and
\begin{equation}
 w\wedge\sum_{n=1}^{N}
 H_{a^{(n)}}\bigl(i\mapsto\varphi(i+n)\bigr)
 \le H_B(\varphi)
 \quad(N\ge1,\ \varphi\in\mathbb N^{\mathbb N}).
 \label{app:fusion-bound}
\end{equation}
No common bound on the regulators $a^{(n)}$ is required.
\end{lemma}
\begin{proof}
For fixed $r$, the sum in the displayed definition is finite. Hence
\[
 0\le B_{rs}\le w,\qquad B_{rs}\downarrow_s0.
\]
Fix $N$, $\varphi$, and $(i_1,\ldots,i_N)\in\mathbb N^N$. Set
\[
 r_n=i_n+n,\qquad
 x_n=a^{(n)}_{i_n,\varphi(i_n+n)},\qquad
 R=\max_{1\le n\le N}r_n,
\]
\[
 C=\bigvee_{1\le n\le N}b_{r_n,\varphi(r_n)}.
\]
For every $n$,
\[
 2^{r_n}x_n\le C,
 \qquad
 2^Rx_n\le2^{R-r_n}C.
\]
Moreover,
\[
 r_n\ge n+1,
 \qquad
 \sum_{n=1}^{N}2^{R-r_n}
 \le\sum_{n=1}^{N}2^{R-n-1}\le2^R.
\]
Consequently,
\[
 2^R\sum_{n=1}^{N}x_n\le2^RC.
\]
Unperforation of abelian lattice-ordered groups yields
$\sum_{n=1}^{N}x_n\le C$. Finite distributivity gives
\[
 w\wedge\sum_{n=1}^{N}x_n
 \le w\wedge C
 =\bigvee_{n=1}^{N}\bigl(w\wedge b_{r_n,\varphi(r_n)}\bigr)
 \le H_B(\varphi).
\]
The set of tuples $\mathbb N^N$ is countable. Its corresponding sums
are bounded above by the sum of the bounds of the first $N$ regulators.
Translation preserves existing suprema. Finite meets preserve bounded
countable suprema in a Dedekind $\sigma$-complete lattice-ordered group.
Therefore
\[
 \begin{split}
 &\sup_{(i_1,\ldots,i_N)\in\mathbb N^N}
 \left(w\wedge\sum_{n=1}^{N}
 a^{(n)}_{i_n,\varphi(i_n+n)}\right)\\
 &\hspace{12mm}=
 w\wedge\sum_{n=1}^{N}
 \sup_{i\ge1}a^{(n)}_{i,\varphi(i+n)}.
 \end{split}
\]
This proves the asserted bound.
\end{proof}

\subsection{Positive transfer and asymptotic order bounds}

With bounded fusion available, we turn to the passage of controls
between ordered groups and to the extraction of attained estimates
from asymptotic bounds. These two operations allow the controller
calculations to be carried into the comparison of integral sums.

\begin{lemma}[Positive transport of regulators]
\label{lt:positive-d-transport}
\label{dt:transport}
Let $E,F$ be Dedekind $\sigma$-complete abelian lattice-ordered groups.
Let $A:E_+\to F_+$ be isotone and subadditive, with
\[
 A(0)=0,\qquad x_s\downarrow0\Longrightarrow A(x_s)\downarrow0.
\]
For a regulator $a$ bounded by $M\in E_+$, put
\[
 c_{rs}=A(M)\wedge2^rA(a_{rs}).
\]
Then $c$ is a regulator and
\begin{equation}\label{lt:regulator-image}
 A(H_a(\varphi))\le H_c(\varphi)
 \qquad(\varphi\in\mathbb N^{\mathbb N}).
\end{equation}
Consequently positive $D$-null nets have $D$-null images under $A$.
In particular, the assertion applies to the positive restriction of
any positive additive sequentially order-continuous map $E\to F$.
\end{lemma}
\begin{proof}
For every fixed $r$,
\[
 a_{rs}\downarrow_s0,\quad A(a_{rs})\downarrow_s0,\quad
 c_{rs}\downarrow_s0,\quad 0\le c_{rs}\le A(M).
\]
If $0\le y_N\uparrow y$, isotonicity and subadditivity give
\[
 0\le A(y)-A(y_N)\le A(y-y_N)\downarrow0,
 \qquad A(y_N)\uparrow A(y).
\]
Fix $\varphi$ and put
\[
 q_N=\bigvee_{r=1}^Na_{r,\varphi(r)},\qquad
 C_N=\bigvee_{r=1}^N2^rA(a_{r,\varphi(r)}).
\]
Then $q_N\uparrow H_a(\varphi)\le M$ and
$A(q_N)\uparrow A(H_a(\varphi))$. Integer multiplication gives
\begin{align*}
 2^N\sum_{r=1}^NA(a_{r,\varphi(r)})
 &\le\sum_{r=1}^N2^{N-r}C_N\\
 &=(2^N-1)C_N\le2^NC_N.
\end{align*}
Unperforation yields $\sum_{r=1}^NA(a_{r,\varphi(r)})\le C_N$.
Moreover $q_N\le\sum_{r=1}^Na_{r,\varphi(r)}$, so
\begin{align*}
 A(q_N)
 &\le A(M)\wedge\sum_{r=1}^NA(a_{r,\varphi(r)})\\
 &\le A(M)\wedge C_N\\
 &=\bigvee_{r=1}^N[A(M)\wedge2^rA(a_{r,\varphi(r)})]
 \le H_c(\varphi).
\end{align*}
The supremum over $N$ proves \eqref{lt:regulator-image}. An attained
bound $e_\alpha\le H_a(\varphi)$ implies
$A(e_\alpha)\le H_c(\varphi)$ on the same tail.
\end{proof}

\begin{lemma}[Transport by an order-continuous majorant]
\label{app:majorant-transport}
Let $T:E_+\to F_+$ be isotone and satisfy
\[
 u_\theta\downarrow0\Longrightarrow T(u_\theta)\downarrow0.
\]
If $e_\alpha\in E_+$ is order null, then $T(e_\alpha)$ is order null.
\end{lemma}
\begin{proof}
For order controls $u_\theta\downarrow0$ of $(e_\alpha)$,
\[
 e_\alpha\le u_\theta\ \text{eventually}
 \quad\Longrightarrow\quad
 T(e_\alpha)\le T(u_\theta)\ \text{eventually}.
\]
The image controls decrease to zero by hypothesis.
\end{proof}

\begin{lemma}[Extraction of an attained tail bound]
\label{app:tail-envelope}
Let $F$ be Dedekind complete and let $v_\gamma\in F_+$ be indexed by a
directed set $\Gamma$. Give $\widehat F_+=F_+\cup\{\infty\}$ its extended
order and define
\[
 M_\delta=\sup_{\gamma\ge\delta}v_\gamma\in\widehat F_+,
 \qquad R=\inf_{\delta\in\Gamma}M_\delta.
\]
If $R\in F_+$, then some $\delta_0$ satisfies $M_{\delta_0}\in F_+$.
On the cofinal tail $\delta\ge\delta_0$,
\[
 M_\delta\downarrow R,
 \qquad r_\delta:=M_\delta-R\downarrow0,
 \qquad v_\delta\le R+r_\delta.
\]
If $v_\gamma\xrightarrow{o}\ell$, then $R=\ell$.
\end{lemma}
\begin{proof}
If $M_\delta=\infty$ for every $\delta$, then $R=\infty$. Thus the
finiteness of $R$ supplies $\delta_0$. Directedness gives
\[
 \delta\ge\delta_0\Longrightarrow0\le M_\delta\le M_{\delta_0},
 \qquad \inf_{\delta\ge\delta_0}M_\delta
 =\inf_{\delta\in\Gamma}M_\delta=R.
\]
The remaining displayed assertions follow by subtraction and
$v_\delta\le M_\delta$.

Suppose $v_\gamma\xrightarrow{o}\ell$, with controls
$u_\theta\downarrow0$. For every $\theta$, some tail satisfies
\[
 \ell-u_\theta\le v_\gamma\le\ell+u_\theta.
\]
The corresponding cofinal family of finite suprema satisfies
\[
 \ell-u_\theta\le M_\delta\le\ell+u_\theta.
\]
Its infimum yields $|R-\ell|\le u_\theta$ for every $\theta$.
Therefore $R=\ell$.
\end{proof}

\begin{lemma}[Moduli and exact conjugate endpoints]
\label{app:conjugate-algebra}
Let $G$ be an abelian group and let $\rho:G\to F_+$ be symmetric and
subadditive, with $\rho(0)=0$. Then
\[
 |\rho(x)-\rho(y)|\le\rho(x-y).
\]
If $F$ is Dedekind complete and $\rho(x_\gamma-x)\xrightarrow{o}0$,
then
\[
 \inf_\delta\sup_{\gamma\ge\delta}\rho(x_\gamma)=\rho(x).
\]
\end{lemma}
\begin{proof}
Subadditivity and symmetry give
\[
 \rho(x)-\rho(y)\le\rho(x-y),\qquad
 \rho(y)-\rho(x)\le\rho(y-x)=\rho(x-y).
\]
Their lattice supremum gives the first inequality. Thus
$\rho(x_\gamma)\xrightarrow{o}\rho(x)$, and
Lemma~\ref{app:tail-envelope} gives the asserted endpoint identity.
\end{proof}

\subsection{Changing an order controller}

The preceding estimates use order controls in a specified group, and
we now examine whether adjoining missing suprema changes the null
nets they describe. Order density and majorization provide the two
approximations needed to pass controls in both directions.

\begin{proposition}[Dedekind-completion invariance]
\label{lt:dedekind-invariance}
Let $E$ be an abelian lattice-ordered group and let
$\iota:E\hookrightarrow E^\delta$ be an order-dense, majorizing
lattice-group embedding into a Dedekind complete group. For every
net $(x_\alpha)$ in $E$,
\begin{equation}\label{lt:dedekind-null}
 x_\alpha\xrightarrow{o_E}0
 \iff\iota(x_\alpha)\xrightarrow{o_{E^\delta}}0.
\end{equation}
Consequently, replacement of an order controller by such a completion preserves a faithful order representation, its Cauchy
nets and its uniform convergence over arbitrary parameter sets.
\end{proposition}
\begin{proof}
Identify $E$ with its image. The embedding is order dense and
majorizing. If $u_\theta\downarrow0$ in $E$ and
$0<w\le u_\theta$ in $E^\delta$ for every $\theta$, order density
provides $0<h\in E$ with
\[
 h\le w\le u_\theta\quad(\theta),
\]
contrary to $\inf_Eu_\theta=0$. Thus every $E$-control remains a
control in $E^\delta$.

For $w\in(E^\delta)^+$, put
\[
 C_w=\{b\in E^+:w\le b\},\qquad q=\inf_{E^\delta}C_w.
\]
Majorization gives $C_w\ne\varnothing$. If $q>w$, order density gives
$0<h\in E$ with $h\le q-w$. For every $b\in C_w$,
\[
 b-h\ge q-h\ge w,
 \qquad b-h\in C_w,
 \qquad q+h\le b.
\]
Taking the infimum over $b\in C_w$ gives $q+h\le q$, a contradiction.
Therefore
\begin{equation}\label{lt:completion-upper-approximation}
 \inf_{E^\delta}C_w=w.
\end{equation}

Suppose $w_\theta\downarrow0$ in $E^\delta$ controls
$\iota(x_\alpha)$. Define
\[
 T=\{b\in E^+:|x_\alpha|\le b\text{ eventually}\}.
\]
Then
\[
 C_{w_\theta}\subseteq T,\qquad
 b,c\in T\Longrightarrow b\wedge c\in T.
\]
For a positive lower bound $h\in E$ of $T$,
\[
 h\le\inf_{E^\delta}C_{w_\theta}=w_\theta\quad(\theta)
 \quad\Longrightarrow\quad h=0.
\]
If $c\in E$ is any lower bound of $T$, then $c\vee0$ is a positive
lower bound, so the preceding argument gives $c\vee0=0$ and $c\le0$.
Thus $T$, directed by reverse order, decreases to zero in $E$ and
controls $(x_\alpha)$. This proves \eqref{lt:dedekind-null}.
For a map $d:Z^2\to E_+$ and a pair net $(z_\alpha,w_\alpha)$,
apply that equivalence to $x_\alpha=d(z_\alpha,w_\alpha)$ to obtain
\begin{equation}\label{lt:completion-pair-equivalence}
 d(z_\alpha,w_\alpha)\xrightarrow{o_E}0
 \iff\iota(d(z_\alpha,w_\alpha))\xrightarrow{o_{E^\delta}}0.
\end{equation}
The embedding preserves zero, addition, lattice operations and order.
Consequently
\[
 \begin{aligned}
 d(x,z)\le d(x,y)+d(y,z)
 &\Longrightarrow
 \iota d(x,z)\le\iota d(x,y)+\iota d(y,z),\\
 d(x,y)=0&\iff\iota d(x,y)=0.
 \end{aligned}
\]
Thus every lattice pseudometric remains a lattice pseudometric, and
every separating family remains separating. If the original pair
convergence is described by
\[
 (z_\alpha,w_\alpha)_\alpha\in\mathcal U
 \iff d_j(z_\alpha,w_\alpha)\to_{o_{E_j}}0\quad(j),
\]
then \eqref{lt:completion-pair-equivalence}, applied coordinatewise,
gives
\[
 (z_\alpha,w_\alpha)_\alpha\in\mathcal U
 \iff\iota_j(d_j(z_\alpha,w_\alpha))
               \to_{o_{E_j^\delta}}0\quad(j).
\]
The point, Cauchy and uniform assertions use, respectively, the pair nets
\[
 (z_\alpha,z)_\alpha,\qquad
 (z_\alpha,z_\beta)_{(\alpha,\beta)},\qquad
 (f_\alpha(t),f(t))_{(\alpha,t)\in I\times T},
\]
where the last index is ordered by $\alpha$ and retains every $t$ on
each tail.
The final formula concerns $T\ne\varnothing$. For $T=\varnothing$,
uniform convergence has no estimates to satisfy in either controller.
These substitutions prove all the asserted invariances.
\end{proof}

Thus Proposition~\ref{lt:dedekind-invariance} permits an existing
lattice-group controller to be enlarged through an order-dense
majorizing extension, without identifying a noncancellative value
monoid with an ordered group. If
regulator convergence is to be compared as well, the additional weak
distributivity and countable-selection assumptions are still those
of Proposition~\ref{lt:order-d}.

\subsection{Countable selection in the controller order}
\label{lc:calibration-section}

To obtain the countable selections needed when no integrable envelope
is given, we remain in an abelian lattice-ordered group, with integer
multiplication understood as repeated addition, and make explicit the
role of principal bands. For $u\ge0$, put
\[
 u^\perp=\{x:|x|\wedge u=0\},\qquad B_u=(u^\perp)^\perp.
\]
A component of $u$ is an element $v\in[0,u]$ satisfying
$v\wedge(u-v)=0$.

\begin{definition}[Order hypotheses on the controller]
\label{lc:calibration}
The controller $F$ is an abelian lattice-ordered group satisfying:
\begin{enumerate}
\item every nonempty countable order-bounded subset has a supremum.
\item every countable disjoint family in $F_+$ has a supremum.
\item every order-bounded array $a_{nk}\downarrow_k0$ satisfies
\[
 \inf_{\varphi\in\mathbb N^{\mathbb N}}
 \bigvee_{n\ge1}a_{n,\varphi(n)}=0.
\]
\item every family of components of $u\in F_+$ with supremum $u$
has a countable subfamily with supremum $u$.
\end{enumerate}
The first three conditions are Dedekind $\sigma$-completeness,
lateral $\sigma$-completeness, and weak $\sigma$-distributivity.
The fourth is the countable cover property for principal components.
\end{definition}

\begin{lemma}[Principal projections]
\label{lc:principal-projections}
Let $F$ be a Dedekind $\sigma$-complete abelian lattice-ordered group.
For $u\in F_+$, the band $B_u=(u^\perp)^\perp$ is a projection band,
with projection
\[
 P_ux=\sup_{n\ge1}(x\wedge nu)\qquad(x\in F_+).
\]
Principal projections commute and preserve existing increasing
countable suprema and decreasing countable infima. Components of $u$
form a $\sigma$-complete Boolean algebra.
\end{lemma}
\begin{proof}
Dedekind $\sigma$-completeness implies the Archimedean property:
if $0\le nx\le y$ for all $n$, then $z=\sup_n nx$ satisfies
$z+x=z$, whence $x=0$.

For $x\ge0$, define $y=\sup_n(x\wedge nu)$ and $r=x-y$.
Then
\[
 0\le r\le(x-nu)^+\qquad(n\ge1).
\]
For $b=(x-nu)^+$, $c=(nu-x)^+$, and $a=u\wedge b$, one has
\[
 b\wedge c=0,\quad (n+1)a\le nu+b,\quad c\le nu+b,
 \quad (n+1)a\wedge c=0.
\]
Consequently,
\[
 (n+1)a+c=((n+1)a)\vee c\le nu+b=x+c,
 \qquad (n+1)a\le x.
\]
Thus $(n+1)(r\wedge u)\le x$ for every $n$, and $r\wedge u=0$.
For $z\in u^\perp$,
\[
 |z|\wedge(x\wedge nu)=0\quad(n\ge1),\qquad
 |z|\wedge y=\sup_n[|z|\wedge(x\wedge nu)]=0.
\]
Hence $y\in B_u$ and $r\in u^\perp$.

If $x=y+r=y'+r'$ are two such decompositions, then
\[
 y-y'=r'-r\in B_u\cap u^\perp=\{0\},
 \qquad y=y',\quad r=r'.
\]
For $x,z\ge0$, the decomposition
\[
 x+z=(P_ux+P_uz)+[(x-P_ux)+(z-P_uz)]
\]
has its first summand in $B_u$ and its second in $u^\perp$.
Uniqueness gives $P_u(x+z)=P_ux+P_uz$. Define
$P_ux=P_ux^+-P_ux^-$ for arbitrary $x$. The group identity
\[
 (x+z)^++x^-+z^-=(x+z)^-+x^++z^+
\]
and additivity on positive elements give
$P_u(x+z)=P_ux+P_uz$. Moreover
\[
 P_u^2=P_u,\qquad \operatorname{ran}P_u=B_u,\qquad
 \ker P_u=u^\perp.
\]
For disjoint elements from these two ideals, modulus is additive.
The decomposition $x=P_ux+(I-P_u)x$ therefore gives
\[
 |x|=|P_ux|+|(I-P_u)x|,\qquad
 P_u|x|=|P_ux|,\qquad (I-P_u)|x|=|(I-P_u)x|.
\]
For $T\in\{P_u,I-P_u\}$, the identity $2x^+=|x|+x$ gives
\[
 2Tx^+=T(|x|+x)=|Tx|+Tx=2(Tx)^+,\qquad Tx^+=(Tx)^+,
\]
\[
 T(x\vee z)=Tz+(Tx-Tz)^+=Tx\vee Tz,\qquad
 T(x\wedge z)=Tx-(Tx-Tz)^+=Tx\wedge Tz.
\]
For $x_n\downarrow0$,
\[
 0\le P_ux_n\le x_n,\qquad0\le(I-P_u)x_n\le x_n,
\]
so both projected sequences decrease to zero. Translation by a known
supremum or infimum proves the countable preservation assertions.

For $z\ge0$, $0\le P_uz\le z$. Thus $P_u$ preserves both $B_v$
and $v^\perp$. Applying $P_v$ to
\[
 P_ux=P_uP_vx+P_u(I-P_v)x
\]
gives $P_vP_ux=P_uP_vx$. Products and complementary products yield
finite meets and differences of nested projections.

A component of $u$ is an element $v\in[0,u]$ with
$v\wedge(u-v)=0$. It satisfies $v=P_vu$. Finite joins of components
are components. If $w_n$ is an increasing sequence of components and
$w=\sup_nw_n$, then
\[
 w\wedge(u-w)=\sup_n[w_n\wedge(u-w)]=0.
\]
Finite joins followed by an increasing supremum provide countable
joins of components. Complements $v\mapsto u-v$ provide countable
meets. The associated projections are their joins and meets within
the principal band $B_u$.
\end{proof}

\begin{proposition}[Independent sufficient conditions for component covers]
\label{lc:component-covers}
The fourth property in Definition~\ref{lc:calibration} holds
if, in each principal band, the Boolean algebra of components satisfies
the countable chain condition. The countable sup property on principal
order intervals is another sufficient condition.
\end{proposition}
\begin{proof}
Let $(c_\alpha)$ be components of $u$ with supremum $u>0$. Consider
disjoint families of nonzero components $d$ of $u$ satisfying
$d\le c_\alpha$ for at least one $\alpha$. Unions of chains are
admissible families. A maximal family $(d_n)_{n\in N}$ exists, and
the countable chain condition makes $N$ countable.
For each $n$, let $\alpha_n$ satisfy $d_n\le c_{\alpha_n}$.
The supremum $d=\bigvee_{n\in N}d_n$ is a component of $u$.
If $u-d>0$, then $(u-d)\wedge c_\alpha>0$ for some $\alpha$:
otherwise every $c_\alpha\le d$, contradicting
$\sup_\alpha c_\alpha=u$. This nonzero component can be adjoined
to the family, contrary to maximality. Therefore
\[
 u=d=\bigvee_{n\in N}d_n
 \le\bigvee_{n\in N}c_{\alpha_n}\le u.
\]
For $u=0$, the assertion holds trivially. The countable sup property
applied directly to $\{c_\alpha:\alpha\in A\}\subseteq[0,u]$
proves the second sufficient condition.
\end{proof}

\begin{proposition}[An integer-valued controller on an arbitrary index set]
\label{lc:integer-controller}
For every set $I$, put
\[
 \mathbb Z^{[I]}
 =\{x:I\to\mathbb Z:\operatorname{supp}x\text{ is countable}\},
 \qquad \operatorname{supp}x=\{i:x(i)\ne0\}.
\]
With coordinatewise addition and order, $\mathbb Z^{[I]}$ is a
super Dedekind complete, laterally $\sigma$-complete, weakly
$\sigma$-distributive abelian lattice-ordered group. In particular,
it satisfies every condition of Definition~\ref{lc:calibration}.
\end{proposition}
\begin{proof}
Finite unions of countable supports are countable. Hence addition,
inversion, finite joins, and finite meets preserve $\mathbb Z^{[I]}$.
Let $\varnothing\ne S\subseteq\mathbb Z^{[I]}$, $x_0\in S$, and
$x\le b$ for $x\in S$. For each $i\in I$, the integer set
$\{x(i):x\in S\}$ is nonempty and bounded above. Thus
\[
 v(i)=\max\{x(i):x\in S\}\in\mathbb Z,
 \qquad x_0(i)\le v(i)\le b(i).
\]
For $K=\operatorname{supp}x_0\cup\operatorname{supp}b$,
\[
 i\notin K\Longrightarrow x_0(i)=b(i)=0
 \Longrightarrow v(i)=0.
\]
Consequently $v\in\mathbb Z^{[I]}$ and $v=\sup S$.
For every $i\in K$, fix $x_i\in S$ with $x_i(i)=v(i)$. Then
\[
 S_0=\{x_0\}\cup\{x_i:i\in K\}\subseteq S,
 \qquad |S_0|\le\aleph_0,
 \qquad \sup S_0=v.
\]
Indeed, the coordinate supremum equals $v(i)$ on $K$ by the defining
property of $x_i$. On $I\setminus K$ it equals zero because
$x_0(i)=v(i)=0$ and every element of $S_0$ is at most $v$.
This proves super Dedekind completeness.

For a countable disjoint family $(y_n)\subseteq(\mathbb Z^{[I]})_+$,
\[
 n\ne m\Longrightarrow y_n(i)\wedge y_m(i)=0\quad(i\in I).
\]
At each coordinate, at most one value $y_n(i)$ is nonzero. Hence
\[
 y(i)=\sup_n y_n(i)\in\mathbb Z_+,
 \qquad \operatorname{supp}y\subseteq\bigcup_n\operatorname{supp}y_n,
 \qquad y=\sup_n y_n\in\mathbb Z^{[I]}.
\]
This proves lateral $\sigma$-completeness.

Let $a_{nk}\downarrow_k0$ be an order-bounded positive array.
For each $i\in I$ and $n\ge1$,
\[
 a_{nk}(i)\downarrow_k0\text{ in }\mathbb Z_+,
 \qquad k_n(i)=\min\{k\ge1:a_{nk}(i)=0\}<\infty.
\]
Set $\varphi_i(n)=k_n(i)$. Coordinatewise suprema give
\[
 H_a(\varphi_i)(i)=\sup_n a_{n,\varphi_i(n)}(i)=0.
\]
If $0\le z\le H_a(\varphi)$ for every selector $\varphi$, then
\[
 0\le z(i)\le H_a(\varphi_i)(i)=0\quad(i\in I),
 \qquad z=0.
\]
Thus $\inf_\varphi H_a(\varphi)=0$. The principal component-cover
property is a consequence of the countable sup property, by
Proposition~\ref{lc:component-covers}.
\end{proof}

\begin{lemma}[Countable principal-band localization]
\label{lc:countable-band}
If $F$ is also laterally $\sigma$-complete, every countable family
$(x_m)\subset F_+$ is contained in one principal band.
\end{lemma}
\begin{proof}
Set
\[
 Q_0=0,\quad Q_m=P_{x_1+\cdots+x_m},\quad
 R_m=Q_m-Q_{m-1},\quad e_m=R_mx_m.
\]
The projections $R_m$ are disjoint. Commutation gives
\[
 R_m=(I-Q_{m-1})P_{x_m}\le P_{x_m},\qquad
 e_m\in R_m(F),\qquad P_{e_m}\le R_m.
\]
If $z\ge0$, $R_mz=z$, and $z\wedge e_m=0$, then
\[
 z\wedge(I-R_m)x_m=0,\qquad
 z\wedge x_m=z\wedge[e_m+(I-R_m)x_m]=0,
\]
\[
 P_{x_m}z=0,\qquad z=R_mz=R_mP_{x_m}z=0.
\]
Thus $R_m(F)\cap e_m^\perp=\{0\}$ and $P_{e_m}=R_m$.
Lateral $\sigma$-completeness defines $e=\bigvee_me_m$. Then
\[
 R_m=P_{e_m}\le P_e,\qquad
 P_{x_m}\le Q_m=\sum_{r\le m}R_r\le P_e,
\]
and $x_m\in B_e$ for every $m$.
\end{proof}

\begin{theorem}[Summable selection from order-null rows]
\label{lc:series-selection}
\label{lc:selection}
Suppose $F$ satisfies Definition~\ref{lc:calibration}. For every
array $x_{nk}\downarrow_k0$ in $F_+$, there exist integers $k_n$ such
that the partial sums
\[
 S_N=\sum_{n=1}^{N}x_{n,k_n}
\]
are order bounded. If
\[
 k(n,s)=\max\{k_n,s\},\qquad
 U_s=\sup_{N\ge1}\sum_{n=1}^{N}x_{n,k(n,s)},
\]
then $U_s\downarrow0$.
\end{theorem}
\begin{proof}
Lemma~\ref{lc:countable-band}, applied to the sequence
$(x_{n1})_{n\ge1}$, supplies $u\in F_+$ with $x_{n1}\in B_u$
for every $n$. Put $P=P_u$. Solidity of the band gives
\[
 0\le x_{nk}\le x_{n1},\qquad
 P x_{n1}=x_{n1},\qquad P x_{nk}=x_{nk}
                 \quad(n,k\ge1).
\]
If $u=0$, then $B_u=\{0\}$ and every $x_{nk}$ is zero, so
$k_n=1$ and $U_s=0$ satisfy the conclusion. Suppose $u>0$.
All projections below act in $B_u$, where $P$ is the identity.
Their countable joins and meets are the component operations of
Lemma~\ref{lc:principal-projections}.

For $n,k\ge1$, put
\[
 E_{nk}=P_{(2^nx_{nk}-u)^+},\qquad a_{nk}=E_{nk}u.
\]
The projections $E_{nk}\le P$ decrease in $k$. For
$y=2^nx_{nk}-u$, disjointness of $y^+$ and $y^-$ gives
\[
 E_{nk}y^-=0,\qquad (P-E_{nk})y^+=0,
 \qquad E_{nk}y=y^+\ge0,\qquad(P-E_{nk})y=-y^-\le0.
\]
Consequently
\[
 0\le a_{nk}=E_{nk}u\le E_{nk}(2^nx_{nk})\le2^nx_{nk},
 \qquad 2^n(P-E_{nk})x_{nk}\le(P-E_{nk})u.
\]
If $z=\inf_ka_{nk}$, then $0\le z\le2^nx_{nk}$ for every $k$.
Integer multiplication preserves decreasing order limits. Therefore
$z=0$. Thus $(a_{nk})$ is a regulator.

For $\varphi\in\mathbb N^{\mathbb N}$, define
\[
 D_\varphi=\bigvee_na_{n,\varphi(n)}.
\]
Weak $\sigma$-distributivity gives $\inf_\varphi D_\varphi=0$.
Every $D_\varphi$ is a component of $u$. The family
$\{u-D_\varphi:\varphi\in\mathbb N^{\mathbb N}\}$ consists of components
and has supremum $u$. The countable component-cover property provides
selectors $\varphi_s$ such that $\inf_sD_{\varphi_s}=0$. Define
\[
 \psi_s(n)=\max_{1\le r\le s}\varphi_r(n),\qquad
 d_s=D_{\psi_s},\qquad k_n=\psi_n(n).
\]
Then
\[
 d_s\downarrow0,\qquad
 n\ge m\Longrightarrow
 E_{n,k_n}u\le E_{n,\psi_m(n)}u\le d_m.
\]
Let
\[
 T_m=\bigvee_{n\ge m}E_{n,k_n}.
\]
The projection identities give
\[
 T_mu=\bigvee_{n\ge m}E_{n,k_n}u\le d_m.
\]
Consequently $\inf_mT_m=0$ on $P(F)$: a projection $T\le T_m$
for every $m$ satisfies $0\le Tu\le d_m$ for every $m$, so
$Tu=0$ and $T=0$.

Put
\[
 V_0=P-T_1,\qquad V_m=T_m-T_{m+1}\quad(m\ge1).
\]
These projections are pairwise disjoint and have supremum $P$.
For $n>m\ge0$,
\[
 V_mE_{n,k_n}=0,
 \qquad 2^nV_mx_{n,k_n}\le V_mu.
\]
For $N>m$, integer multiplication gives
\begin{align*}
 2^N\sum_{n=m+1}^{N}V_mx_{n,k_n}
 &\le\sum_{n=m+1}^{N}2^{N-n}V_mu\\
 &=(2^{N-m}-1)V_mu\le2^NV_mu.
\end{align*}
An abelian lattice-ordered group is unperforated. Hence
\[
 \sum_{n=m+1}^{N}V_mx_{n,k_n}\le V_mu.
\]
Define
\[
 w_0=V_0u,\qquad
 w_m=V_m\!\left(u+\sum_{n=1}^{m}x_{n,k_n}\right)\quad(m\ge1).
\]
For every $m,N$, including $N\le m$,
\[
 V_mS_N\le w_m.
\]
The elements $w_m$ are disjoint and positive. Their supremum $w$
exists by $\sigma$-lateral completeness. For
$z=(S_N-w)^+\in B_u$,
\[
 V_mz=0\quad(m\ge0),\qquad z=T_mz\quad(m\ge1),
\]
\[
 0\le z\wedge u=T_m(z\wedge u)\le T_mu\le d_m\downarrow0.
\]
Thus $z\wedge u=0$. Since $z\in B_u$, one has $z=0$, hence $S_N\le w$.
Dedekind $\sigma$-completeness defines $S=\sup_NS_N$.

Since $k(n,s)\ge k_n$,
\[
 0\le\sum_{n=1}^{N}x_{n,k(n,s)}\le S_N\le S.
\]
Thus $U_s$ exists, and $U_{s+1}\le U_s$. For fixed $M$,
\[
 U_s\le\sum_{n=1}^{M}x_{n,k(n,s)}+(S-S_M).
\]
The finite sum decreases to zero as $s\to\infty$. Therefore
\[
 0\le\inf_sU_s\le S-S_M\qquad(M\ge1).
\]
Since $S_M\uparrow S$, one obtains $\inf_sU_s=0$.
\end{proof}

\begin{corollary}[Simultaneous control of countably many regulators]
\label{lc:regulator-calibration}
Let $a^{(n)}=(a^{(n)}_{ik})$ be countably many regulators in a
controller group satisfying Definition~\ref{lc:calibration}.
There exist selectors $\psi_{n,s}$, increasing coordinatewise in $s$,
and $B_s\downarrow0$ such that
\begin{equation}\label{lc:regulator-budget}
 \sum_{n=1}^{N}H_{a^{(n)}}(\psi_{n,s})\le B_s
 \qquad(N,s\ge1).
\end{equation}
There also exist single selectors $\psi_s$ such that
\begin{equation}\label{lc:shifted-budget}
 \sum_{n=1}^{N}
 H_{a^{(n)}}(i\mapsto\psi_s(i+n))\le B_s
 \qquad(N,s\ge1).
\end{equation}
\end{corollary}
\begin{proof}
Let $\nu:\mathbb N^2\to\mathbb N$ be a bijection. Apply
Theorem~\ref{lc:series-selection} to the rows
$x_{\nu(n,i),k}=a^{(n)}_{ik}$. Write $B_s$ for the resulting
budgets and put $\psi_{n,s}(i)=k(\nu(n,i),s)$. Then
\[
 \sum_{n=1}^{N}\bigvee_{1\le i\le M}a^{(n)}_{i,\psi_{n,s}(i)}
 \le\sum_{n=1}^{N}\sum_{i=1}^{M}a^{(n)}_{i,\psi_{n,s}(i)}
 \le B_s.
\]
For each fixed $N$, all joins in the displayed finite sum increase
with $M$, and addition preserves existing increasing suprema. Thus
\[
 \begin{aligned}
 \sum_{n=1}^N H_{a^{(n)}}(\psi_{n,s})
 &=\sum_{n=1}^N\sup_M\bigvee_{1\le i\le M}
                     a^{(n)}_{i,\psi_{n,s}(i)}\\
 &=\sup_M\sum_{n=1}^N\bigvee_{1\le i\le M}
                     a^{(n)}_{i,\psi_{n,s}(i)}\le B_s.
 \end{aligned}
\]
The common index $M$ is cofinal in the finite product of the
individual join indices, because their maximum dominates each of
them. This proves \eqref{lc:regulator-budget}.
Define
\[
 \psi_s(1)=1,\qquad
 \psi_s(r)=\max_{1\le n<r}\psi_{n,s}(r-n)\quad(r\ge2).
\]
For every $n,i$,
$a^{(n)}_{i,\psi_s(i+n)}\le a^{(n)}_{i,\psi_{n,s}(i)}$.
Therefore
\[
 \sum_{n=1}^{N}H_{a^{(n)}}(i\mapsto\psi_s(i+n))
 \le\sum_{n=1}^{N}H_{a^{(n)}}(\psi_{n,s})\le B_s.
\]
\end{proof}

Through Theorem~\ref{lc:series-selection}, the selection argument is
therefore reduced to the controller group's order, distributivity, and
disjoint suprema. Its conditions are formulated independently of any
division error, integral, or candidate limit, and enter the limit
theorems only in the branch without a known integrable envelope.

\section{Exact metrization and the intrinsic representation dictionary}
\label{app:metrization}

Having used the original convergence to formulate the integral and its
limit theorems, we now ask how that convergence can be represented by
ordered values without losing its point, pair, or group forms. We begin
with complete ordered values that retain every originally convergent
net, and follow the additional algebra needed for lattice realizations
through closed-base stability to a canonical closure criterion for
ordered-group representations. Along this route, countability,
boundedness, and weak distributivity must be traced to the object on
which each is imposed: the examples distinguish properties of the
original convergence from those of its chosen values, while the proofs
record each additional algebraic assumption at the step where it is
used.

To keep these representations within the language of nets, the
set-indexed constructions below use only eventual sets of actual nets.
More explicitly, put
\[
 \operatorname{Ev}(Y)
 =\{\operatorname{Tail}(u):u\text{ is an admissible net in }Y\}.
\]
Every $\mathcal F\in\operatorname{Ev}(Y)$ therefore carries the canonical
net $n_{\mathcal F}=\langle A:A\in\mathcal F\rangle$, whose tail image
above $(A,a)$ is exactly $A$. An assertion $\mathcal F\to y$ in a
set-indexed formula means $n_{\mathcal F}\to y$, and membership
$\mathcal F\in\mathcal U$ for a pair family means
$n_{\mathcal F}\in\mathcal U$. Thus no convergence is assigned to a
collection independently of a net. Products, image operations and
intersections of eventual-set collections are interpreted through the
product, image and mixing nets of Proposition~\ref{fd:tail-calculus}.
Whenever a centered collection is used to construct a new net, the
finite intersections are first formed and their nonemptiness is proved,
then their saturation supplies the indicated net.

\subsection{Pair translations and canonical group convergence}
\label{amc:section}

In an abelian convergence group, the group operations already determine
which pair nets are uniformly null from the original convergence.
Recording this calculus of tails first allows point, Cauchy, and uniform
faithfulness to be read in that convergence before any value object is
constructed.

\begin{proposition}[Three exact translations]\label{lt:three-translations}
Let $\lambda$ be induced by the original uniform convergence structure
$\mathcal U$. Suppose maps $d_j:Z^2\to M_j$ and null convergences
$\chi_j$ satisfy, for every pair net $v=(v_\alpha)$,
\begin{equation}\label{amx:general-pair-faithfulness}
 v\in\mathcal U
 \quad\Longleftrightarrow\quad d_j(v_\alpha)\xrightarrow{\chi_j}0
 \quad(j\in J).
\end{equation}
For a net $(z_\alpha)_{\alpha\in I}$, a point $z$, and a nonempty
parameter family $(z_\alpha^A,z^A)$, the following equivalences hold:
\begin{align}
 z_\alpha\xrightarrow{\lambda}z
 &\iff d_j(z_\alpha,z)\xrightarrow{\chi_j}0\quad(j\in J),
       \label{lt:point-translation}\\
 (z_\alpha)_\alpha\text{ is }\mathcal U\text{-Cauchy}
 &\iff(d_j(z_\alpha,z_\beta))_{(\alpha,\beta)\in I^2}
                   \xrightarrow{\chi_j}0\quad(j\in J),
       \label{lt:cauchy-translation}\\
 z_\alpha^A\longrightarrow z^A\text{ uniformly in }A
 &\iff(d_j(z_\alpha^A,z^A))_{(\alpha,A)}
                   \xrightarrow{\chi_j}0\quad(j\in J).
       \label{lt:uniform-translation}
\end{align}
The last net is indexed by $I\times\mathcal A$, ordered through its
first coordinate. In the order case the assertion is
\begin{equation}\label{lt:uniform-controls}
 \forall j\ \exists u^j_\theta\downarrow0\ \forall\theta\
 \exists\alpha_0\ \forall\alpha\ge\alpha_0\ \forall A\in\mathcal A:
 \quad d_j(z_\alpha^A,z^A)\le u^j_\theta.
\end{equation}
\end{proposition}
\begin{proof}
The three pair nets in question are
\[
 (z_\alpha,z)_\alpha,\qquad
 (z_\alpha,z_\beta)_{(\alpha,\beta)\in I^2},\qquad
 (z_\alpha^A,z^A)_{(\alpha,A)\in I\times\mathcal A}.
\]
Their tail images are respectively
\[
 T_{\alpha_0}(z)\times\{z\},\qquad
 T_{\alpha_0}(z)\times T_{\beta_0}(z),\qquad
 \{(z_\alpha^A,z^A):\alpha\ge\alpha_0,\ A\in\mathcal A\}.
\]
The definitions of point convergence, the Cauchy property and uniform
convergence assert original uniform nullity of these three nets.
Applying \eqref{amx:general-pair-faithfulness} gives the three stated
equivalences. For the last net, eventual membership of its image under
$d_j$ in $[0,u^j_\theta]$ is exactly the final two universal
quantifiers in \eqref{lt:uniform-controls}, because every tail contains
all parameters $A\in\mathcal A$.
\end{proof}

Uniformity over $\mathcal A$ is expressed in
\eqref{lt:uniform-controls} through common tail estimates, without a
supremum over the parameter family or an implicit order bound on the
family as a whole.

\begin{example}[Point fidelity does not determine pair fidelity]
\label{amx:point-not-pair}
On the ordinary real line let
\[
 d(x,y)=|\arctan x-\arctan y|.
\]
For a fixed $x$ and a net $u=(u_i)$, put $\mathcal F=\operatorname{Tail}(u)$. Then
continuity of $\arctan$ and of its inverse on
$(-\pi/2,\pi/2)$ gives
\[
 \mathcal F\to x
 \quad\Longleftrightarrow\quad
 d(\,\cdot\,,x)_*\mathcal F\to0.
\]
However, $(n)$ is $d$-Cauchy since $\arctan n\to\pi/2$,
while the pair differences $n-m$ do not converge to zero:
every pair tail contains $(n,m)$ with $|n-m|\ge1$.
For $f(k)=k$, define $f_n(k)=k$ when $k\ne n$ and $f_n(n)=2n$.
Then
\[
 \sup_{k\ge1}d(f_n(k),f(k))
 =\arctan(2n)-\arctan(n)\longrightarrow0,
\]
but each uniform error tail contains the unbounded family
$\{n:n\ge n_0\}$.
Thus fixed-point fidelity alone supplies neither the original group
Cauchyness nor the original uniform convergence in
Proposition~\ref{lt:three-translations}.
\end{example}

\begin{proposition}[The original group uniformity]
\label{lt:canonical-group}
Let $(Z,\lambda)$ be an abelian convergence group. Set
\[
 \delta(z,w)=z-w,\qquad
 \mathcal U_\lambda
 =\{\Phi\in\operatorname{Ev}(Z^2):
                  \delta_*\Phi\xrightarrow{\lambda}0\}.
\]
Then $\mathcal U_\lambda$ is a uniform convergence structure and induces
$\lambda$.

Suppose $\mathcal U$ is already a uniform convergence structure inducing
$\lambda$, the group operations are uniformly continuous, and
$[\Delta_Z]\in\mathcal U$. Then
\begin{equation}\label{lt:original-difference}
 \mathcal U=\mathcal U_\lambda.
\end{equation}
In particular, let $M_j$ be ordered commutative monoids with bottom
identity, or positive cones of partially ordered abelian groups, and
let $\rho_j:Z\to M_j$ be symmetric subadditive maps with $\rho_j(0)=0$.
If they satisfy
\begin{equation}\label{lt:null-faithfulness}
 \mathcal H\xrightarrow{\lambda}0
 \iff (\rho_j)_*\mathcal H\xrightarrow{\chi_j}0\quad(j\in J)
\end{equation}
then they give the faithful pair representation
\begin{equation}\label{lt:invariant-metric}
 d_j(z,w)=\rho_j(z-w).
\end{equation}
\end{proposition}
\begin{proof}
For relations $R,S\subseteq Z^2$,
\[
 \delta(\Delta_Z)=\{0\},\qquad
 \delta(R^{-1})=-\delta(R),\qquad
 \delta(R\circ S)\subseteq\delta(R)+\delta(S).
\]
For refinement and finite intersections,
\[
 \begin{aligned}
 \Phi\subseteq\Psi
 &\Longrightarrow\delta_*\Phi\subseteq\delta_*\Psi,\\
 \delta_*(\Phi\cap\Psi)
 &=\delta_*\Phi\cap\delta_*\Psi,\\
 \delta_*\Phi\to_\lambda0,\quad\delta_*\Psi\to_\lambda0
 &\Longrightarrow\delta_*(\Phi\cap\Psi)\to_\lambda0.
 \end{aligned}
\]
The image identity is obtained by applying inverse images to each set:
\[
 H\in\delta_*(\Phi\cap\Psi)
 \iff\delta^{-1}(H)\in\Phi\text{ and }\delta^{-1}(H)\in\Psi.
\]
Moreover,
\[
 \begin{aligned}
 \delta_*[\Delta_Z]&=[0]\to_\lambda0,\\
 \delta_*(\Phi^{-1})&=(-\operatorname{id})_*(\delta_*\Phi),\\
 \delta_*\Phi\to_\lambda0
 &\Longrightarrow\delta_*(\Phi^{-1})\to_\lambda0.
 \end{aligned}
\]
Suppose every composition tail $R\circ S$, with $R\in\Phi$ and $S\in\Psi$, is nonempty, so their saturation is defined. Let $+$ denote the eventual-set collection obtained by adding independent image nets. The relation inclusion gives
\[
 (\delta_*\Phi)+(\delta_*\Psi)
 \subseteq\delta_*(\Phi\circ\Psi).
\]
Therefore
\[
 \begin{aligned}
 \delta_*\Phi\to_\lambda0,\quad\delta_*\Psi\to_\lambda0
 &\Longrightarrow
 (\delta_*\Phi)+(\delta_*\Psi)\to_\lambda0\\
 &\Longrightarrow\delta_*(\Phi\circ\Psi)\to_\lambda0.
 \end{aligned}
\]
For the induced point convergence,
\[
 \delta_*(\mathcal F\times[z])=\mathcal F-z,
 \qquad
 \mathcal F-z\xrightarrow{\lambda}0
 \iff\mathcal F\xrightarrow{\lambda}z.
\]

For \eqref{lt:original-difference}, let $\Phi\in\mathcal U$. Consider
the image of $\Phi$ under
\[
 (z,w)\longmapsto\bigl((z,w),(w,w)\bigr).
\]
Its two coordinate pair-net tails are $\Phi$ and a refinement of
$[\Delta_Z]$. Uniform continuity of subtraction gives
\[
 (z,w)\longmapsto(z-w,0),
 \qquad (\delta_*\Phi)\times[0]\in\mathcal U.
\]
Thus $\delta_*\Phi\xrightarrow{\lambda}0$.

Conversely, suppose $\delta_*\Phi\xrightarrow{\lambda}0$. The image
of $\Phi$ under
\[
 (z,w)\longmapsto\bigl((z-w,w),(0,w)\bigr)
\]
has coordinate pair-net tails $(\delta_*\Phi)\times[0]$ and a refinement
of $[\Delta_Z]$. Uniform continuity of addition gives its image
\[
 (z,w)\longmapsto(z,w),
\]
so $\Phi\in\mathcal U$. Substitution of
\eqref{lt:null-faithfulness} into \eqref{lt:original-difference} proves
\eqref{lt:pair-faithfulness} for \eqref{lt:invariant-metric}.
\end{proof}

\subsection{Completeness of the represented space}
\label{lt:represented-completeness-section}

Once Cauchy nets have been translated, completeness can be read through
the representation as well. Here it is the completeness of the original
space that supplies a limit from estimates for pairs of approximants in
the integration arguments; the values enter only through their specified
null convergences, so the translation is available before any order
completeness is imposed on them.

\begin{proposition}[Completeness under a faithful pair representation]
\label{lt:represented-completeness}
Let $(Z,\mathcal U)$ be separated, and let $d_j:Z^2\to M_j$ satisfy
\eqref{amx:general-pair-faithfulness}. The following are equivalent:
\begin{enumerate}[(i)]
\item Every original Cauchy net in $Z$ has an original limit.
\item For every downward-directed family $\mathcal B$ of nonempty
subsets of $Z$, with saturation $u=\langle B:B\in\mathcal B\rangle$,
\begin{equation}\label{lt:metric-tail-completeness}
 \left[(d_j(u_\alpha,u_\beta))_{\alpha,\beta}
                      \xrightarrow{\chi_j}0\quad(j\in J)\right]
 \Longrightarrow
 \left[\exists!z\in Z:\ d_j(u_\alpha,z)
                      \xrightarrow{\chi_j}0\quad(j\in J)\right].
\end{equation}
\item For every net $(z_\alpha)_{\alpha\in I}$ in $Z$,
\begin{equation}\label{lt:metric-net-completeness}
 \left[(d_j(z_\alpha,z_\beta))_{(\alpha,\beta)\in I^2}
                      \xrightarrow{\chi_j}0\quad(j\in J)\right]
 \Longrightarrow
 \left[\exists!z\in Z:\ d_j(z_\alpha,z)
                      \xrightarrow{\chi_j}0\quad(j\in J)\right].
\end{equation}
\end{enumerate}
Sequential completeness is equivalent to (iii) restricted to sequences.
No completeness of the value objects is required.
\end{proposition}
\begin{proof}
Proposition~\ref{lt:three-translations} identifies the antecedent of
(iii) with the original Cauchy property and its conclusion with original
convergence. Thus (i) implies (iii), and (iii) supplies an original limit
to every original Cauchy net, proving (i). Assertion (iii) implies (ii)
by application to the displayed saturation.

Conversely, let $z=(z_\alpha)$ satisfy the antecedent of (iii). Its
nonempty tail images $T_{\alpha_0}(z)$ form a downward-directed family.
The saturation of this family is tail equivalent to $z$, while the
product of two copies of the saturation is tail equivalent to
$(z_\alpha,z_\beta)_{\alpha,\beta}$. Image nets under $d_j$ preserve
these tail equivalences, so (ii) gives a point $z_0$ with
$d_j(z_\alpha,z_0)\to0$ for every $j$. The first translation in
Proposition~\ref{lt:three-translations} gives $z_\alpha\to z_0$.
Separation makes the limit unique. Restricting the same translations to
$I=\mathbb N$ proves the sequential assertion, with Cauchy differences
still indexed by $\mathbb N^2$.
\end{proof}

\begin{corollary}[Invariant order-valued metrics]
\label{lt:invariant-completeness}
Let $(X,\lambda)$ be a separated abelian convergence group and let
$d:X^2\to E_+$ be an invariant metric with values in a partially
ordered abelian group. Put $p(x)=d(x,0)$ and suppose
\[
 x_\alpha\xrightarrow{\lambda}0
 \iff p(x_\alpha)\xrightarrow{o_E}0
\]
for every admissible net. Then $X$ is complete if and only if every
admissible net $(x_\alpha)_{\alpha\in I}$ satisfying
\begin{equation}\label{lt:metric-cauchy-controls}
 \exists u_\eta\downarrow_E0\quad
 \forall\eta\ \exists\alpha_\eta\ \forall\alpha,\beta\ge\alpha_\eta:
 \qquad d(x_\alpha,x_\beta)\le u_\eta
\end{equation}
admits a unique $x\in X$ such that
\begin{equation}\label{lt:metric-limit-controls}
 \exists v_\theta\downarrow_E0\quad
 \forall\theta\ \exists\alpha_\theta\ \forall\alpha\ge\alpha_\theta:
 \qquad d(x_\alpha,x)\le v_\theta.
\end{equation}
Restricting $(x_\alpha)$ to sequences gives the equivalence for
sequential completeness. The regulators in both statements may
have arbitrary directed index sets. Their infima are taken in
$E$, and no lattice structure or completeness of $E$ is assumed.
\end{corollary}
\begin{proof}
Invariance and null faithfulness give
\[
 \begin{aligned}
 d(x,y)&=d(x-y,0)=p(x-y),\\
 (x_\alpha)\text{ is originally Cauchy}
 &\iff(x_\alpha-x_\beta)_{(\alpha,\beta)\in I^2}
                     \xrightarrow{\lambda}0\\
 &\iff(d(x_\alpha,x_\beta))_{(\alpha,\beta)\in I^2}
                     \xrightarrow{o_E}0,\\
 x_\alpha\xrightarrow{\lambda}x
 &\iff x_\alpha-x\xrightarrow{\lambda}0
 \iff d(x_\alpha,x)\xrightarrow{o_E}0.
 \end{aligned}
\]
For each member $u_\eta$ of a Cauchy regulator, product-net
eventuality gives indices $a_\eta,b_\eta$ such that
\[
 \alpha\ge a_\eta,\quad\beta\ge b_\eta
 \quad\Longrightarrow\quad d(x_\alpha,x_\beta)\le u_\eta.
\]
Directedness supplies $\alpha_\eta\ge a_\eta,b_\eta$, yielding
\eqref{lt:metric-cauchy-controls}. Conversely, that condition
controls the product tail beginning at
$(\alpha_\eta,\alpha_\eta)$. Equation
\eqref{lt:metric-limit-controls} is the order-null condition
in the last displayed equivalence. Consequently, completeness
gives its conclusion from \eqref{lt:metric-cauchy-controls}, and
the converse implication supplies an original limit to every
original Cauchy net. Restriction to $I=\mathbb N$ proves the
sequential assertion, without restricting either regulator.

For a positive decreasing family, the infimum convention is
\[
 \inf_E u_\eta=0
 \iff
 (\forall w\in E)\quad
 \bigl[((\forall\eta)\ w\le u_\eta)\Longrightarrow w\le0\bigr].
\]
Only the infima of the displayed regulators occur in this proof.
\end{proof}

Although these criteria express completeness through ordered values,
the property still belongs to the original space. To see why the two
objects must be distinguished, we give examples in both directions of
failure, with the actual distances generating the entire value group
in each case.

\begin{example}[A complete domain with incomplete metric values]
\label{lt:complete-domain-incomplete-values}
Let $c$ be the real vector lattice of convergent real sequences, let
$X=c$ carry its supremum-norm convergence, and let $E=c$ carry its
coordinatewise order. With $\mathbf1=(1,1,\ldots)$, define
\[
 p(x)=|x|+\|x\|_\infty\mathbf1,
 \qquad d(x,y)=p(x-y).
\]
Then $d$ is a faithful invariant $E$-valued metric and $X$ is complete,
although order convergence on $E$ is not sequentially complete and $E$
is not Dedekind $\sigma$-complete. Moreover, the additive subgroup
generated by $d(X^2)$ is all of $E$.
\end{example}
\begin{proof}
For $x,y\in c$,
\[
 \begin{gathered}
 p(x)=0\iff x=0,\qquad p(-x)=p(x),\\
 p(x+y)
 \le |x|+|y|+(\|x\|_\infty+\|y\|_\infty)\mathbf1
 =p(x)+p(y),\\
 \|x\|_\infty\mathbf1\le p(x)
 \le2\|x\|_\infty\mathbf1.
 \end{gathered}
\]
These inequalities give the metric identities and translation
invariance. If $\|x_\alpha\|_\infty\to0$, then
\[
 \forall n\ge1\quad
 p(x_\alpha)\le n^{-1}\mathbf1\quad\text{eventually},
 \qquad n^{-1}\mathbf1\downarrow0\quad\text{in }E.
\]
Conversely, suppose $p(x_\alpha)\to_o0$, with eventual bounds
$b_\beta\downarrow0$ in $E$. Failure of norm convergence would give
$\varepsilon>0$ such that
\[
 (\forall\alpha_0)(\exists\alpha\ge\alpha_0)
 \quad \|x_\alpha\|_\infty\ge\varepsilon.
\]
Each tail bounded by $b_\beta$ would then contain an index satisfying
\[
 \varepsilon\mathbf1
 \le\|x_\alpha\|_\infty\mathbf1
 \le p(x_\alpha)\le b_\beta.
\]
This would make $\varepsilon\mathbf1$ a positive lower bound of
$(b_\beta)$, contrary to $\inf_\beta b_\beta=0$. Thus
$p(x_\alpha)\to_o0$ if and only if
$\|x_\alpha\|_\infty\to0$, for arbitrary nets. Applying this equivalence
to the difference net $x_\alpha-y_\alpha$ gives pair faithfulness.

For completeness of $X$, let $(x_\alpha)$ be a norm-Cauchy net.
Completeness of $\mathbb R$ supplies
$a(k)=\lim_\alpha x_\alpha(k)$ for every $k\ge1$. For every
$\varepsilon>0$ there is $\alpha_\varepsilon$ such that
\[
 \alpha,\beta\ge\alpha_\varepsilon
 \quad\Longrightarrow\quad
 |x_\alpha(k)-x_\beta(k)|\le\varepsilon
 \quad(k\ge1).
\]
Passing to the scalar limit in $\beta$ yields
\[
 \alpha\ge\alpha_\varepsilon
 \quad\Longrightarrow\quad
 |x_\alpha(k)-a(k)|\le\varepsilon
 \quad(k\ge1).
\]
In particular, $a\in\ell^\infty$. Since
$x_{\alpha_\varepsilon}\in c$, there is $K_\varepsilon$ with
\[
 k,l\ge K_\varepsilon
 \quad\Longrightarrow\quad
 |a(k)-a(l)|
 \le2\varepsilon
   +|x_{\alpha_\varepsilon}(k)-x_{\alpha_\varepsilon}(l)|
 \le3\varepsilon.
\]
The scalar sequence $a(k)$ therefore converges, so $a\in c$, and the
preceding estimates give $\|x_\alpha-a\|_\infty\to0$.

Let $e_k$ denote the $k$th coordinate vector and put
\[
 f_n=\sum_{k=1}^n e_{2k},
 \qquad
 r_N(j)=\begin{cases}0,&j\le2N,\\1,&j>2N.\end{cases}
\]
Then $r_N\in c$, $0\le r_{N+1}\le r_N$, and
\[
 h\le r_N\ (N\ge1)
 \quad\Longrightarrow\quad
 h(j)\le r_j(j)=0\ (j\ge1)
 \quad\Longrightarrow\quad h\le0.
\]
Consequently $r_N\downarrow0$ in $E$, while
\[
 m,n\ge N\quad\Longrightarrow\quad |f_m-f_n|\le r_N.
\]
Thus $(f_n)$ is order-Cauchy. If $f_n\to_o f\in c$, let
$h_\beta\downarrow0$ bound $|f_n-f|$ eventually, and set
$a_{2k}=1$, $a_{2k-1}=0$. For every coordinate $j$ and every $\beta$,
an index in that tail with $n\ge j$ gives
\[
 |a_j-f(j)|e_j\le h_\beta.
\]
Taking the infimum in $\beta$ gives $f(j)=a_j$ for every $j$, which
contradicts $f\in c$. Furthermore, $0\le f_n\le\mathbf1$, and
$\mathbf1-e_{2j-1}$ is an upper bound of all $f_n$ for every $j$.
If $s=\sup_n f_n$ existed in $E$, these bounds would imply
\[
 s(2j)=1,\qquad s(2j-1)=0\qquad(j\ge1),
\]
again contradicting $s\in c$. Finally, for every $a\in c_+$,
\[
 a=p(a)-p\!\left(\frac{\|a\|_\infty}{2}\mathbf1\right).
\]
Both terms belong to $d(X^2)$, and $x=x^+-x^-$ for every $x\in c$.
Hence these distances generate $E$ as an additive group.
\end{proof}

\begin{example}[Complete metric values with an incomplete domain]
\label{lt:complete-values-incomplete-domain}
Let $X=c_{00}$ be the real vector space of finitely supported sequences
with supremum-norm convergence, and put
\[
 E=\mathbb R,\qquad d(x,y)=\|x-y\|_\infty.
\]
Then $d$ is a faithful invariant metric with Dedekind-complete values,
whereas $X$ is not sequentially complete.
\end{example}
\begin{proof}
For a nonnegative real net $(s_\alpha)$,
\[
 s_\alpha\to0
 \iff(\forall\varepsilon>0)\quad
             s_\alpha\le\varepsilon\text{ eventually}
 \iff s_\alpha\to_o0.
\]
For the last equivalence, ordinary convergence supplies the decreasing
bounds $1/n$, while any $b_\beta\downarrow0$ in $\mathbb R$ has a
member below each $\varepsilon>0$. Thus the norm metric faithfully
represents point and pair convergence for arbitrary nets. Its values
contain every $t\ge0$, since $d(te_1,0)=t$.

The sequence
\[
 z_n=\sum_{k=1}^n\frac1k e_k
\]
satisfies
\[
 \|z_m-z_n\|_\infty
 \le\frac1{\min\{m,n\}+1},
\]
and is therefore Cauchy. For every $x\in c_{00}$ there is $K\ge1$
such that $x(k)=0$ for $k>K$, whence
\[
 n\ge K+1\quad\Longrightarrow\quad
 \|z_n-x\|_\infty
 \ge |z_n(K+1)-x(K+1)|=\frac1{K+1}>0.
\]
The sequence has no limit in $X$.
\end{proof}

\begin{remark}[Completeness in the integral limit arguments]
\label{lt:completeness-roles}
In Theorem~\ref{lt:order-transfer}, the original Cauchy property is
obtained from
\[
 p,q\ge p_0(j,\theta)
 \quad\Longrightarrow\quad
 d_j(y_p^A,y_q^A)\le\beta_\theta^j\oplus_j\beta_\theta^j,
 \qquad\beta_\theta^j\oplus_j\beta_\theta^j\downarrow0.
\]
Proposition~\ref{lt:represented-completeness} then identifies the
completeness that supplies the endpoint $z^A$. The controller
calculus in Lemma~\ref{mg:calculus} requires no order completeness.
For Theorem~\ref{lc:equiintegrability}, the distinct controller
construction is
\[
 q_{j,N}=\sup_{n,m\ge N}d_j(y_n,y_m),
 \qquad q_j=\inf_Nq_{j,N}.
\]
Dedekind $\sigma$-completeness supplies these elements of $F_j$,
and weak $\sigma$-distributivity yields $q_j=0$. Original sequential
completeness of the range then supplies the limit of $(y_n)$.
In the monotone constructions of Theorems~\ref{fv:envelope-monotone}
and~\ref{lc:levi}, Dedekind $\sigma$-completeness of the target $Z$
constructs $J(C)=\sup_nJ_n(C)$ from the proved upper bounds, and the
subsequent estimates establish convergence to $J(C)$ in the original
structure. Completeness of that convergence is not needed for this
direct construction of its limit.
\end{remark}

\subsection{Complete value objects and exact pair metrization}
\label{amu:section}

With the roles of the space and its values separated, we can construct
a representation that retains the original family of small pair-net
tails while satisfying the triangle inequality. An explicit calculation
of every tail radius will recover the given family from the value
object, before the uniform convergence axioms enter the argument.

For a nonempty set $Y$, the set $\operatorname{Ev}(Y)$ consists of
$\operatorname{Tail}(u)$ for the admissible nets $u$ in $Y$.
If $\varnothing\ne A\subseteq Y$, let $s_A$ be the net indexed by $A$
with every two indices comparable and value $s_A(a)=a$. Every tail image
of $s_A$ is $A$, and therefore
\[
 [A]:=\operatorname{Tail}(s_A)=\{B\subseteq Y:A\subseteq B\},
 \qquad \dot y=[\{y\}].
\]
For nets $u,v$ with eventual-set collections $\mathcal F,\mathcal G$,
the mixing net on the product of their index sets and $\{0,1\}$ has
values $u_i$ and $v_j$ and tail images $T_i(u)\cup T_j(v)$. Consequently
\begin{equation}\label{amu:tail-intersection}
 \mathcal F\cap\mathcal G
 =\operatorname{Tail}\langle A\cup B:
                         A\in\mathcal F,\ B\in\mathcal G\rangle.
\end{equation}
For a map $T:Y\to V$, eventual membership of $T(u)$ in $C$ means
eventual membership of $u$ in $T^{-1}(C)$, which gives
\begin{equation}\label{amu:image-intersection}
 \begin{aligned}
 C\in T_*(\mathcal F\cap\mathcal G)
 &\iff T^{-1}(C)\in\mathcal F\cap\mathcal G\\
 &\iff C\in T_*\mathcal F\cap T_*\mathcal G.
 \end{aligned}
\end{equation}
All collections in these equalities are derived from the displayed
nets, and each convergence assertion below refers to their canonical
saturations.

\begin{definition}[Complete ordered values and radii]
\label{amu:value-monoid}
A complete ordered value monoid is a commutative monoid
$(M,\oplus,0)$ equipped with a complete lattice order whose least
element is $0$, such that
\[
 a\le b\quad\Longrightarrow\quad a\oplus c\le b\oplus c.
\]
An $M$-valued metric on $X$ is a map $q:X^2\to M$ satisfying
\[
 q(x,y)=0\iff x=y,\qquad q(x,y)=q(y,x),\qquad
 q(x,z)\le q(x,y)\oplus q(y,z).
\]
For any map $v:Y\to M$, put
\begin{equation}\label{amu:radius}
 V_v(A)=\bigvee_{y\in A}v(y),\qquad
 r_v(\mathcal F)=\bigwedge_{A\in\mathcal F}V_v(A)
 \quad\bigl(\mathcal F\in\operatorname{Ev}(Y)\bigr).
\end{equation}
The empty join has value $0$.
\end{definition}

\begin{lemma}[Bases and tail radii]\label{amu:base-radius}
If $\mathcal B$ is a base of $\mathcal F$, then
\begin{equation}\label{amu:base-radius-equation}
 r_v(\mathcal F)=\bigwedge_{B\in\mathcal B}V_v(B).
\end{equation}
Consequently, every net $y=(y_\alpha)_{\alpha\in I}$ satisfies
\begin{equation}\label{amu:tail-radius}
 r_v(\operatorname{Tail}(y))
 =\bigwedge_{\alpha_0\in I}\bigvee_{\alpha\ge\alpha_0}v(y_\alpha).
\end{equation}
\end{lemma}
\begin{proof}
The inclusion $\mathcal B\subseteq\mathcal F$ gives
\[
 r_v(\mathcal F)\le\bigwedge_{B\in\mathcal B}V_v(B).
\]
For every $A\in\mathcal F$,
\[
 \{B\in\mathcal B:B\subseteq A\}\ne\varnothing,
 \qquad
 \bigwedge_{C\in\mathcal B}V_v(C)\le V_v(B)\le V_v(A)
 \quad(B\in\mathcal B,\ B\subseteq A).
\]
Taking the infimum over $A$ gives the reverse inequality.
The tail images $T_{\alpha_0}(y)$ form a base of
$\operatorname{Tail}(y)$, and
\[
 V_v(T_{\alpha_0}(y))=\bigvee_{\alpha\ge\alpha_0}v(y_\alpha).
\]
Substitution proves \eqref{amu:tail-radius}.
\end{proof}

\begin{lemma}[Incidence representation]\label{amu:incidence}
Let $\varnothing\ne Z\subseteq Y$, let $I\ne\varnothing$, and let
$\mathcal B_i$ be a nonempty family of subsets of $Y$ for every $i\in I$.
Assume
\begin{equation}\label{amu:incidence-assumptions}
 \begin{gathered}
 Z\subseteq B\qquad(i\in I,\ B\in\mathcal B_i),\\
 (\forall y\in Y\setminus Z)(\forall i\in I)\quad
       \{B\in\mathcal B_i:y\notin B\}\ne\varnothing.
 \end{gathered}
\end{equation}
There are a complete ordered value monoid $M$ and a map $v:Y\to M$
such that
\begin{equation}\label{amu:incidence-separation}
 v(y)=0\iff y\in Z,
\end{equation}
and, for every $\mathcal F\in\operatorname{Ev}(Y)$,
\begin{equation}\label{amu:incidence-exactness}
 r_v(\mathcal F)=0
 \iff
 Z\in\mathcal F
 \ \text{or}\
 (\exists i\in I)\ \mathcal B_i\subseteq\mathcal F.
\end{equation}
The addition can satisfy
$a\oplus b=\top$ whenever $a\ne0$ and $b\ne0$.
\end{lemma}
\begin{proof}
For a nonempty set $A$, put
$\mathcal P_*(A)=\mathcal P(A)\setminus\{\varnothing\}$, and define
the disjoint union
\begin{equation}\label{amu:incidence-values}
 M=\{0\}\sqcup\prod_{i\in I}\mathcal P_*(\mathcal B_i).
\end{equation}
The product contains $(\mathcal B_i)_{i\in I}$.
Set $0\le a$ for every $a\in M$ and, for nonzero $a,b$, set
\[
 a\le b\iff a_i\subseteq b_i\quad(i\in I).
\]
For $\mathscr E\subseteq M$, the join is
\begin{equation}\label{amu:incidence-joins}
 \bigvee\mathscr E=
 \begin{cases}
 0,&\mathscr E\setminus\{0\}=\varnothing,\\[2pt]
 \left(\displaystyle\bigcup_{a\in\mathscr E\setminus\{0\}}a_i\right)_{i\in I},
      &\mathscr E\setminus\{0\}\ne\varnothing.
 \end{cases}
\end{equation}
In the second case, each coordinate contains the nonempty coordinate
of every fixed member of $\mathscr E\setminus\{0\}$. For $b\in M$,
\[
 [a\le b\ (a\in\mathscr E)]
 \Longrightarrow b\ne0\quad\text{and}\quad
 \bigcup_{a\in\mathscr E\setminus\{0\}}a_i\subseteq b_i
 \quad(i\in I).
\]
The displayed vector is therefore the least upper bound. The first
case is the least-upper-bound property of $0$.

The top element is
$\top=(\mathcal B_i)_{i\in I}$. For $\mathscr E\ne\varnothing$,
\begin{equation}\label{amu:incidence-meets}
 \bigwedge\mathscr E=
 \begin{cases}
 0,&0\in\mathscr E,\\
 0,&0\notin\mathscr E\text{ and }
       (\exists i\in I)\ \displaystyle\bigcap_{a\in\mathscr E}a_i=\varnothing,\\[2pt]
 \left(\displaystyle\bigcap_{a\in\mathscr E}a_i\right)_{i\in I},
     &0\notin\mathscr E\text{ and every intersection is nonempty}.
 \end{cases}
\end{equation}
The empty meet is $\top$. A nonzero lower bound $b$ for a nonempty
family of nonzero elements satisfies
\[
 \varnothing\ne b_i\subseteq\bigcap_{a\in\mathscr E}a_i
 \quad(i\in I).
\]
Thus no nonzero lower bound exists in the second case. In the last case,
the displayed vector is a lower bound containing every coordinate of
every nonzero lower bound. In the first case, every lower bound is at
most $0$. This proves completeness of the lattice.

Define
\begin{equation}\label{amu:flat-addition}
 0\oplus a=a\oplus0=a,
 \qquad a\oplus b=\top\quad(a\ne0,\ b\ne0).
\end{equation}
The operation is commutative. For $a,b,c\in M$, both
$(a\oplus b)\oplus c$ and $a\oplus(b\oplus c)$ are
\[
 \begin{cases}
 0,&a=b=c=0,\\
 \text{the unique nonzero member of }(a,b,c),
       &\#\{j\in\{1,2,3\}:(a,b,c)_j\ne0\}=1,\\
 \top,&\#\{j\in\{1,2,3\}:(a,b,c)_j\ne0\}\ge2.
 \end{cases}
\]
Hence addition is associative. If $a\le b$, then
\[
 \begin{array}{ll}
 c=0:&a\oplus c=a\le b=b\oplus c,\\
 c\ne0,\ a=0:&a\oplus c=c\le b\oplus c,\\
 c\ne0,\ a\ne0:&b\ne0,\quad a\oplus c=\top=b\oplus c.
 \end{array}
\]
Thus $M$ is a complete ordered value monoid.

For $y\in Y\setminus Z$, define
\[
 D_i(y)=\{B\in\mathcal B_i:y\notin B\},
 \qquad v(y)=(D_i(y))_{i\in I},
 \qquad v(y)=0\quad(y\in Z).
\]
Every $D_i(y)$ is nonempty by \eqref{amu:incidence-assumptions}.
This proves \eqref{amu:incidence-separation}. For $A\not\subseteq Z$,
\begin{equation}\label{amu:incidence-diameter}
 \begin{aligned}
 V_v(A)_i
 &=\bigcup_{y\in A\setminus Z}D_i(y)\\
 &=\{B\in\mathcal B_i:(\exists y\in A\setminus Z)\ y\notin B\}\\
 &=\{B\in\mathcal B_i:A\not\subseteq B\}.
 \end{aligned}
\end{equation}
The last equality uses $Z\subseteq B$. If $A\subseteq Z$, then
$V_v(A)=0$.

Let $\mathcal F\in\operatorname{Ev}(Y)$. If $Z\in\mathcal F$, then
\[
 0\le r_v(\mathcal F)\le V_v(Z)=0.
\]
If $Z\notin\mathcal F$, upward closure gives
$A\not\subseteq Z$ for every $A\in\mathcal F$. All corresponding
$V_v(A)$ are nonzero. Equations
\eqref{amu:incidence-meets} and \eqref{amu:incidence-diameter} yield
\[
 \begin{aligned}
 r_v(\mathcal F)=0
 &\iff(\exists i\in I)\quad
    \bigcap_{A\in\mathcal F}\{B\in\mathcal B_i:A\not\subseteq B\}
       =\varnothing\\
 &\iff(\exists i\in I)(\forall B\in\mathcal B_i)
                   (\exists A\in\mathcal F)\quad A\subseteq B\\
 &\iff(\exists i\in I)\quad\mathcal B_i\subseteq\mathcal F.
 \end{aligned}
\]
For the last equivalence, the forward implication uses upward closure.
The reverse implication uses $A=B$. This proves
\eqref{amu:incidence-exactness}.
\end{proof}

For a nonempty set $X$, write
\[
 \Delta_X=\{(x,x):x\in X\},
 \qquad A^{-1}=\{(y,x):(x,y)\in A\}.
\]
The inverse pair net is obtained by applying $(x,y)\mapsto(y,x)$ to each value, and the superscript $^{-1}$ on its eventual-set collection records this image operation.

\begin{theorem}[Exact symmetric pair representation]
\label{amu:pair-metrization}
Let $X\ne\varnothing$ and
$\mathcal U\subseteq\operatorname{Ev}(X^2)$. The following conditions
are equivalent.
\begin{enumerate}[(i)]
\item There exist a complete ordered value monoid $M$ and an $M$-valued
metric $q$ such that
\begin{equation}\label{lt:exact-monoid}
 \Phi\in\mathcal U\iff r_q(\Phi)=0
 \qquad\bigl(\Phi\in\operatorname{Ev}(X^2)\bigr).
\end{equation}
\item The family $\mathcal U$ satisfies
\begin{align}
 &[\Delta_X]\in\mathcal U,
                                      \label{amu:global-diagonal}\\
 &\Phi\in\mathcal U,\quad\Phi\subseteq\Psi
       \quad\Longrightarrow\quad\Psi\in\mathcal U,
                                      \label{amu:refinement}\\
 &\dot{(x,y)}\in\mathcal U\quad\Longrightarrow\quad x=y,
                                      \label{amu:separation}\\
 &\Phi\in\mathcal U\quad\Longrightarrow\quad\Phi^{-1}\in\mathcal U,
                                      \label{amu:inversion}\\
 &\Phi\in\mathcal U\quad\Longrightarrow\quad
     \Phi\cap\Phi^{-1}\cap[\Delta_X]\in\mathcal U.
                                      \label{amu:balanced-core}
\end{align}
Every eventual-set collection in these assertions is the tail collection of an admissible pair net.
\end{enumerate}
\end{theorem}
\begin{proof}
Assume~(i). For every pair net $v$, with $\Phi=\operatorname{Tail}(v)$,
\[
 r_q([\Delta_X])=V_q(\Delta_X)=0,
 \qquad r_q(\dot{(x,y)})=q(x,y),
\]
and
\[
 \Phi\subseteq\Psi\quad\Longrightarrow\quad
 0\le r_q(\Psi)\le r_q(\Phi).
\]
Together with \eqref{lt:exact-monoid} and separation of $q$, these
identities give
\eqref{amu:global-diagonal}--\eqref{amu:separation}.
The sets $A^{-1}$, $A\in\Phi$, form a base of $\Phi^{-1}$, and
symmetry gives
\[
 V_q(A^{-1})=V_q(A),
 \qquad r_q(\Phi^{-1})=r_q(\Phi).
\]
Hence \eqref{amu:inversion} holds.

Set $\mathcal C=\Phi\cap\Phi^{-1}\cap[\Delta_X]$. For $A\in\Phi$,
\[
 R_A=A\cup A^{-1}\cup\Delta_X\in\mathcal C.
\]
If $C\in\mathcal C$, then
\[
 C,C^{-1}\in\Phi,
 \qquad A=C\cap C^{-1}\in\Phi,
 \qquad A\cup A^{-1}\cup\Delta_X\subseteq C.
\]
Thus $(R_A)_{A\in\Phi}$ is a base of $\mathcal C$. Symmetry and the
zero diagonal give
\[
 V_q(R_A)=V_q(A),
 \qquad
 r_q(\mathcal C)=\bigwedge_{A\in\Phi}V_q(R_A)=r_q(\Phi).
\]
Equation \eqref{lt:exact-monoid} proves
\eqref{amu:balanced-core}.

Assume~(ii). Use the nonempty index set $I=\mathcal U$. For each
$\Phi\in\mathcal U$, put
\[
 \mathcal C_\Phi=\Phi\cap\Phi^{-1}\cap[\Delta_X],
 \qquad
 \mathcal B_\Phi
   =\{B\in\mathcal C_\Phi:B=B^{-1}\}.
\]
The core belongs to $\mathcal U$ by \eqref{amu:balanced-core} and
satisfies $\mathcal C_\Phi^{-1}=\mathcal C_\Phi$. Therefore
\[
 C\in\mathcal C_\Phi
 \quad\Longrightarrow\quad
 C\cap C^{-1}\in\mathcal B_\Phi,
 \qquad C\cap C^{-1}\subseteq C.
\]
The family $\mathcal B_\Phi$ is a base of $\mathcal C_\Phi$, and every
member contains $\Delta_X$. For an off-diagonal pair $(x,y)$,
\[
 \begin{aligned}
 [(x,y)\in B\ (B\in\mathcal B_\Phi)]
 &\Longrightarrow [(x,y)\in C\ (C\in\mathcal C_\Phi)]\\
 &\Longrightarrow\mathcal C_\Phi\subseteq\dot{(x,y)}\\
 &\Longrightarrow\dot{(x,y)}\in\mathcal U\\
 &\Longrightarrow x=y.
 \end{aligned}
\]
The third implication is \eqref{amu:refinement}. The fourth is
\eqref{amu:separation}. Consequently,
\[
 x\ne y\quad\Longrightarrow\quad
 \{B\in\mathcal B_\Phi:(x,y)\notin B\}\ne\varnothing.
\]
Lemma~\ref{amu:incidence} applies with
$Y=X^2$, $Z=\Delta_X$, and these families. Let $M$ and $q:X^2\to M$
be its output. Then
\[
 q(x,y)=0\iff x=y.
\]
Every coordinate defining $q(x,y)$ consists of symmetric sets, so
\[
 \{B\in\mathcal B_\Phi:(x,y)\notin B\}
 =\{B\in\mathcal B_\Phi:(y,x)\notin B\},
 \qquad q(x,y)=q(y,x).
\]
For the triangle inequality,
\[
 \begin{array}{ll}
 x=y:&q(x,y)\oplus q(y,z)=0\oplus q(x,z)=q(x,z),\\
 y=z:&q(x,y)\oplus q(y,z)=q(x,z)\oplus0=q(x,z),\\
 x\ne y,\ y\ne z:&q(x,y)\oplus q(y,z)=\top\ge q(x,z).
 \end{array}
\]
Hence $q$ is an $M$-valued metric.

For every pair net $w$, with $\Psi=\operatorname{Tail}(w)$,
\begin{equation}\label{amu:pair-coordinate-test}
 r_q(\Psi)=0
 \iff
 \Delta_X\in\Psi
 \ \text{or}\
 (\exists\Phi\in\mathcal U)\ \mathcal B_\Phi\subseteq\Psi.
\end{equation}
If $\Delta_X\in\Psi$, then $[\Delta_X]\subseteq\Psi$, so
\eqref{amu:global-diagonal} and \eqref{amu:refinement} give
$\Psi\in\mathcal U$. If $\mathcal B_\Phi\subseteq\Psi$, then
\[
 C\in\mathcal C_\Phi
 \Longrightarrow C\cap C^{-1}\in\mathcal B_\Phi\subseteq\Psi
 \Longrightarrow C\in\Psi.
\]
Thus $\mathcal C_\Phi\subseteq\Psi$ and
$\Psi\in\mathcal U$. Conversely, if $\Psi\in\mathcal U$, then
\[
 \mathcal B_\Psi\subseteq\mathcal C_\Psi\subseteq\Psi,
\]
so \eqref{amu:pair-coordinate-test} gives $r_q(\Psi)=0$.
This proves \eqref{lt:exact-monoid}.
\end{proof}

Because the representation retains the full prescribed family of pair
nets, it also retains any finite stability and composability belonging
to that family. When we apply it to a uniform convergence structure,
these properties are consequently carried by the given family itself;
they need not be inferred from the value monoid.

\begin{corollary}[Exact uniform metrization and order-null tails]
\label{amu:uniform-metrization}
Every separated globally diagonal uniform convergence structure
$\mathcal U$ has a representation \eqref{lt:exact-monoid}.
For every pair net $(x_\alpha,y_\alpha)$, the representation satisfies
\begin{align}
 \operatorname{Tail}((x_\alpha,y_\alpha)_\alpha)\in\mathcal U
 &\iff
 \bigwedge_{\alpha_0}\bigvee_{\alpha\ge\alpha_0}
                              q(x_\alpha,y_\alpha)=0,
                                       \label{fd:exact-net-radius}\\
 &\iff
 (\exists u_\theta\downarrow0\text{ in }M)
 (\forall\theta)(\exists\alpha_0)(\forall\alpha\ge\alpha_0)
                \quad q(x_\alpha,y_\alpha)\le u_\theta.
                                       \label{amu:net-order-exactness}
\end{align}
The point, Cauchy and parameter-uniform forms are
\begin{align}
 x_\alpha\to_\lambda x
 &\iff\bigwedge_{\alpha_0}\bigvee_{\alpha\ge\alpha_0}
                         q(x_\alpha,x)=0,
                                       \label{amu:point-specialization}\\
 (x_\alpha)_\alpha\text{ is }\mathcal U\text{-Cauchy}
 &\iff\bigwedge_{\alpha_0,\beta_0}
        \bigvee_{\substack{\alpha\ge\alpha_0\\\beta\ge\beta_0}}
                         q(x_\alpha,x_\beta)=0,
                                       \label{amu:cauchy-specialization}\\
 x_\alpha^A\to x^A\text{ uniformly in }A
 &\iff\bigwedge_{\alpha_0}
        \bigvee_{\substack{\alpha\ge\alpha_0\\A\in\mathcal A}}
                         q(x_\alpha^A,x^A)=0.
                                       \label{amu:uniform-specialization}
\end{align}
Here $\lambda$ is the original point convergence induced by
$\mathcal U$, and $\mathcal A$ is any nonempty set of parameters.
\end{corollary}
\begin{proof}
The uniform convergence axioms supply refinement, separation and
inversion. Global diagonality and finite stability give
\[
 \Phi\in\mathcal U
 \Longrightarrow\Phi,\Phi^{-1},[\Delta_X]\in\mathcal U
 \Longrightarrow\Phi\cap\Phi^{-1}\cap[\Delta_X]\in\mathcal U.
\]
All conditions of Theorem~\ref{amu:pair-metrization} hold.
Lemma~\ref{amu:base-radius} yields
\eqref{fd:exact-net-radius}.

Put
\[
 b_{\alpha_0}=\bigvee_{\alpha\ge\alpha_0}q(x_\alpha,y_\alpha).
\]
These values satisfy
\[
 \alpha_1\ge\alpha_0\Longrightarrow b_{\alpha_1}\le b_{\alpha_0},
 \qquad q(x_\alpha,y_\alpha)\le b_{\alpha_0}
                   \quad(\alpha\ge\alpha_0).
\]
If the radius is zero, then $b_{\alpha_0}\downarrow0$ and this net
supplies the controls in \eqref{amu:net-order-exactness}.
Conversely, for every stipulated control $u_\theta$,
\[
 (\exists\alpha_0)\quad b_{\alpha_0}\le u_\theta,
 \qquad
 0\le\bigwedge_{\alpha_0}b_{\alpha_0}\le u_\theta.
\]
The infimum over $\theta$ is zero.

The pair nets in the last three assertions have tail bases
\[
 \begin{gathered}
 T_{\alpha_0}(x)\times\{x\},\\
 T_{\alpha_0}(x)\times T_{\beta_0}(x),\\
 \{(x_\alpha^A,x^A):\alpha\ge\alpha_0,\ A\in\mathcal A\},
 \end{gathered}
\]
respectively. The first generates
$\operatorname{Tail}(x)\times[x]$, the second generates
$\operatorname{Tail}(x)\times\operatorname{Tail}(x)$, and the third
is the parameter-uniform pair-net tail collection from
Proposition~\ref{lt:three-translations}.
Applying \eqref{lt:exact-monoid} and
\eqref{amu:base-radius-equation} to these three bases gives
\eqref{amu:point-specialization}--\eqref{amu:uniform-specialization}.
\end{proof}

At this point the full diagonal condition cannot be weakened merely
to separate constant nets, since
$r_q([\Delta_X])=0$ for every metric with zero diagonal. Requiring only
the individual constant nets $\dot{(x,x)}$ therefore does not by itself
imply the pair representation of Theorem~\ref{amu:pair-metrization},
whereas the canonical group structure of
Proposition~\ref{lt:canonical-group} satisfies the full condition. The
resulting order-null equivalence takes place in the complete value
monoid, leaving cancellativity, lattice-group structure, and weak
sigma-distributivity as separate algebraic properties.

\subsection{Point convergence and its symmetric pair representation}
\label{amp:section}

If only point convergence is initially given, the preceding method
leaves us with the task of constructing its pair family. We do so by
pairing each convergent net with its limit, adjoining the reversed
pairs and the full diagonal, and proving that no further point limits
are introduced. This allows the representation method to begin with
point convergence alone.

For a nonempty set $X$, $[x]$ denotes the eventual sets of the constant net at $x$, while $[\Delta_X]$ denotes the eventual sets of the net ranging over the whole diagonal on the indiscrete directed preorder.
If $R\subseteq X^2$ and $y\in X$, put
\[
 R^y=\{z\in X:(z,y)\in R\}.
\]
For every net $v$ in $X$, with $\mathcal G=\operatorname{Tail}(v)$,
\begin{equation}\label{amp:section-membership}
 \begin{aligned}
 R\in\mathcal G\times[y]
 &\iff(\exists B\in\mathcal G)\quad B\times\{y\}\subseteq R\\
 &\iff(\exists B\in\mathcal G)\quad B\subseteq R^y\\
 &\iff R^y\in\mathcal G.
 \end{aligned}
\end{equation}

\begin{lemma}[The symmetric fan of a point net]
\label{amp:fan}
Let $u=(u_i)$ be a net in $X$, put $\mathcal F=\operatorname{Tail}(u)$, and let $x\in X$. Define
\begin{equation}\label{amp:fan-tail}
 \mathcal C(\mathcal F,x)
 = (\mathcal F\times[x])\cap([x]\times\mathcal F)
                       \cap[\Delta_X].
\end{equation}
Then $\mathcal C(\mathcal F,x)$ is the eventual-set collection of a pair net with the tail base
\begin{equation}\label{amp:fan-relations}
 R(A,x)=(A\times\{x\})\cup(\{x\}\times A)\cup\Delta_X,
 \qquad A\in\mathcal F.
\end{equation}
Every base relation is symmetric and contains $\Delta_X$. In particular,
\begin{equation}\label{amp:fan-identities}
 \mathcal C(\mathcal F,x)^{-1}=\mathcal C(\mathcal F,x),
 \qquad
 \mathcal C(\mathcal F,x)\subseteq[\Delta_X],
 \qquad
 \mathcal C([x],x)=[\Delta_X].
\end{equation}
The sections of a base relation are
\begin{equation}\label{amp:fan-sections}
 R(A,x)^y=
 \begin{cases}
 A\cup\{x\},&y=x,\\
 \{x,y\},&y\ne x,\ y\in A,\\
 \{y\},&y\ne x,\ y\notin A.
 \end{cases}
\end{equation}
\end{lemma}
\begin{proof}
The mixing of the pair nets $(u_i,x)$, $(x,u_i)$ and the full diagonal net has eventual-set collection \eqref{amp:fan-tail}. Its tail images are nonempty because they contain the diagonal. For $A\in\mathcal F$,
\[
 \begin{gathered}
 A\times\{x\}\subseteq R(A,x),\qquad
 \{x\}\times A\subseteq R(A,x),\qquad
 \Delta_X\subseteq R(A,x),\\
 R(A,x)\in(\mathcal F\times[x])\cap([x]\times\mathcal F)
                                      \cap[\Delta_X].
 \end{gathered}
\]
Conversely, $R\in\mathcal C(\mathcal F,x)$ provides
$A_1,A_2\in\mathcal F$ satisfying
\[
 A_1\times\{x\}\subseteq R,
 \qquad\{x\}\times A_2\subseteq R,
 \qquad\Delta_X\subseteq R.
\]
With $A=A_1\cap A_2\in\mathcal F$,
\[
 R(A,x)\subseteq R.
\]
Thus \eqref{amp:fan-relations} is a base. The identities
\[
 R(A,x)^{-1}=R(A,x),\qquad
 R(\{x\},x)=\Delta_X
\]
and the diagonal inclusion give \eqref{amp:fan-identities}.
For a fixed $y$, membership in the three terms of
\eqref{amp:fan-relations} is expressed by
\[
 (z,y)\in R(A,x)
 \iff
 [(y=x\ \text{and}\ z\in A)
 \ \text{or}\ (z=x\ \text{and}\ y\in A)
 \ \text{or}\ z=y].
\]
The three cases in \eqref{amp:fan-sections} result upon separating
$y=x$, $y\ne x$ with $y\in A$, and $y\ne x$ with $y\notin A$.
\end{proof}

\begin{theorem}[Exact metrization of $T_1$ point convergence]
\label{amp:point-metrization}
Let $\lambda$ be a convergence satisfying
Definition~\ref{netconvergencestructure} on a nonempty set $X$. Assume
\begin{equation}\label{amp:T1}
 [y]\longrightarrow_\lambda x\quad\Longrightarrow\quad y=x.
\end{equation}
There are a complete ordered value monoid $(M,\oplus,0)$ and a separated
symmetric $M$-valued metric $q$ such that
\begin{align}
 \mathcal F\longrightarrow_\lambda x
 &\iff r_q(\mathcal F\times[x])=0,
 &&\mathcal F\in\operatorname{Ev}(X),\ x\in X,
 \label{amp:tail-exact}\\
 x_\alpha\longrightarrow_\lambda x
 &\iff
 \bigwedge_{\alpha_0}\bigvee_{\alpha\ge\alpha_0}
                  q(x_\alpha,x)=0.
 \label{amp:net-exact}
\end{align}
\end{theorem}
\begin{proof}
Proposition~\ref{fd:tail-calculus} gives the net-tail implications
\begin{equation}\label{amp:point-axioms}
 \begin{gathered}
 [x]\to_\lambda x,\\
 \mathcal F\to_\lambda x,\quad\mathcal F\subseteq\mathcal G
       \quad\Longrightarrow\quad\mathcal G\to_\lambda x,\\
 \mathcal F\to_\lambda x
       \quad\Longrightarrow\quad\mathcal F\cap[x]\to_\lambda x.
 \end{gathered}
\end{equation}
For the last implication, the constant convergence $[x]\to_\lambda x$
is combined with the finite-intersection assertion of that proposition.
Every collection quantified below is the eventual-set collection of its canonical net.

Define the family of pair-net tails
\begin{equation}\label{amp:fan-family}
 \mathcal U_{\lambda}^{\mathrm{fan}}
 =\left\{\mathcal H\in\operatorname{Ev}(X^2):
   \begin{array}{l}
   (\exists x\in X)(\exists\mathcal F\in\operatorname{Ev}(X))\\[-2pt]
   \mathcal F\to_\lambda x,\quad
   \mathcal C(\mathcal F,x)\subseteq\mathcal H
   \end{array}\right\}.
\end{equation}
We verify the hypotheses of Theorem~\ref{amu:pair-metrization}.
First,
\[
 [x]\to_\lambda x,\qquad
 \mathcal C([x],x)=[\Delta_X]
 \quad\Longrightarrow\quad
 [\Delta_X]\in\mathcal U_{\lambda}^{\mathrm{fan}}.
\]
For refinement,
\[
 \mathcal C(\mathcal F,x)\subseteq\mathcal H\subseteq\mathcal K,
 \quad\mathcal F\to_\lambda x
 \quad\Longrightarrow\quad
 \mathcal K\in\mathcal U_{\lambda}^{\mathrm{fan}}.
\]
For inversion and balanced cores, \eqref{amp:fan-identities} gives
\[
 \begin{aligned}
 \mathcal C(\mathcal F,x)\subseteq\mathcal H
 &\Longrightarrow
 \mathcal C(\mathcal F,x)
 =\mathcal C(\mathcal F,x)^{-1}\subseteq\mathcal H^{-1},\\
 \mathcal C(\mathcal F,x)\subseteq\mathcal H
 &\Longrightarrow
 \mathcal C(\mathcal F,x)
 \subseteq\mathcal H\cap\mathcal H^{-1}\cap[\Delta_X].
 \end{aligned}
\]
If $R_i$ are tail images of a pair net representing $\mathcal H$,
the mixing net of that net, its inverse and the full diagonal net
has a cofinal tail family
\[
 R_i\cup R_i^{-1}\cup\Delta_X.
\]
Each displayed tail contains $\Delta_X$, and its eventual-set
collection is $\mathcal H\cap\mathcal H^{-1}\cap[\Delta_X]$. Therefore
\[
 \mathcal H\in\mathcal U_{\lambda}^{\mathrm{fan}}
 \quad\Longrightarrow\quad
 \mathcal H^{-1},\
 \mathcal H\cap\mathcal H^{-1}\cap[\Delta_X]
                  \in\mathcal U_{\lambda}^{\mathrm{fan}}.
\]

To verify separation, suppose
$[(u,v)]\in\mathcal U_{\lambda}^{\mathrm{fan}}$ and $u\ne v$.
By \eqref{amp:fan-family}, there are $x$ and $\mathcal F\to_\lambda x$
such that
\[
 \mathcal C(\mathcal F,x)\subseteq[(u,v)],
 \qquad (u,v)\in R(A,x)\quad(A\in\mathcal F).
\]
Since $X\in\mathcal F$, membership in $R(X,x)$ and $u\ne v$ yield
$u=x$ or $v=x$. These alternatives are exclusive.
If $u=x$, then $v\ne x$ and
\[
 (x,v)\in R(A,x)\quad\Longrightarrow\quad v\in A
 \qquad(A\in\mathcal F).
\]
Thus $\mathcal F\subseteq[v]$, and
\[
 \mathcal F\to_\lambda x
 \Longrightarrow[v]\to_\lambda x
 \Longrightarrow v=x=u,
\]
contrary to $u\ne v$. If $v=x$, then $u\ne x$ and
\[
 (u,x)\in R(A,x)\quad\Longrightarrow\quad u\in A
 \qquad(A\in\mathcal F).
\]
Consequently $\mathcal F\subseteq[u]$ and
\[
 \mathcal F\to_\lambda x
 \Longrightarrow[u]\to_\lambda x
 \Longrightarrow u=x=v,
\]
again a contradiction. Hence
\[
 [(u,v)]\in\mathcal U_{\lambda}^{\mathrm{fan}}
 \quad\Longrightarrow\quad u=v.
\]
Theorem~\ref{amu:pair-metrization} now supplies $M$ and $q$ with
\begin{equation}\label{amp:pair-exact}
 \mathcal H\in\mathcal U_{\lambda}^{\mathrm{fan}}
 \iff r_q(\mathcal H)=0
 \qquad(\mathcal H\in\operatorname{Ev}(X^2)).
\end{equation}

We next identify the point limits of this family of pair-net tails.
If $\mathcal G\to_\lambda y$, then
\[
 \mathcal C(\mathcal G,y)\subseteq\mathcal G\times[y],
 \qquad\mathcal G\times[y]\in\mathcal U_{\lambda}^{\mathrm{fan}}.
\]
Conversely, suppose
$\mathcal G\times[y]\in\mathcal U_{\lambda}^{\mathrm{fan}}$.
There are $x$ and $\mathcal F\to_\lambda x$ such that
\[
 \mathcal C(\mathcal F,x)\subseteq\mathcal G\times[y].
\]
Equations \eqref{amp:section-membership} and
\eqref{amp:fan-relations} give
\begin{equation}\label{amp:section-transfer}
 R(A,x)^y\in\mathcal G\qquad(A\in\mathcal F).
\end{equation}
If $y=x$, then \eqref{amp:fan-sections} yields
$A\cup\{x\}\in\mathcal G$ for every $A\in\mathcal F$.
For $B\in\mathcal F\cap[x]$,
\[
 B\in\mathcal F,\quad x\in B,\quad
 B=B\cup\{x\}\in\mathcal G.
\]
Thus
\[
 \mathcal F\cap[x]\subseteq\mathcal G,
 \qquad\mathcal F\cap[x]\to_\lambda x,
 \qquad\mathcal G\to_\lambda x=y.
\]
If $y\ne x$, then
\[
 \bigl[(\forall A\in\mathcal F)\ y\in A\bigr]
 \Longrightarrow\mathcal F\subseteq[y]
 \Longrightarrow[y]\to_\lambda x
 \Longrightarrow y=x.
\]
Hence some $A\in\mathcal F$ satisfies $y\notin A$. For this member,
\[
 R(A,x)^y=\{y\}\in\mathcal G.
\]
Properness gives $y\in B$ for every $B\in\mathcal G$. Upward closure
and $\{y\}\in\mathcal G$ show that the underlying net is eventually equal to $y$.
Therefore $\mathcal G=[y]\to_\lambda y$. We have proved
\[
 \mathcal G\to_\lambda y
 \iff\mathcal G\times[y]\in\mathcal U_{\lambda}^{\mathrm{fan}}.
\]
Combining this equivalence with \eqref{amp:pair-exact} proves
\eqref{amp:tail-exact}.

Finally, for a net $x_\alpha$ and its eventual-set collection $\mathcal T$,
\[
 \{T_{\alpha_0}(x)\times\{x\}:\alpha_0\}
\]
is a base of $\mathcal T\times[x]$. Its radius is
\[
 \begin{aligned}
 r_q(\mathcal T\times[x])
 &=\bigwedge_{\alpha_0}
     \bigvee_{z\in T_{\alpha_0}(x)}q(z,x)\\
 &=\bigwedge_{\alpha_0}
     \bigvee_{\alpha\ge\alpha_0}q(x_\alpha,x).
 \end{aligned}
\]
For completeness, the first equality holds in both order directions:
each base member is an eventual set, and every eventual set
contains a base member. The supremum over a contained set cannot exceed
the supremum over the containing set. Tail equivalence identifies
$x_\alpha\to_\lambda x$ with $\mathcal T\to_\lambda x$.
Equation~\eqref{amp:tail-exact} therefore gives
\eqref{amp:net-exact}.
\end{proof}

We have thus obtained the point representation from the three
implications in \eqref{amp:point-axioms} and the $T_1$ condition,
without requiring a diagonal axiom for iterated limits. For this
purpose the auxiliary family $\mathcal U_{\lambda}^{\mathrm{fan}}$
need not be asserted to have finite stability or relational
composability; whenever a uniform convergence structure is prescribed,
its own pair-net tails are retained through
Theorem~\ref{amu:pair-metrization}.

\subsection{Invariant metrization of the original group convergence}
\label{amg:section}

Returning to an abelian convergence group allows the construction to
respect translations, since its point and pair convergence are already
determined by the null-net tails at the origin. Applying the incidence
construction of Lemma~\ref{amu:incidence} to those tail images gives
the value monoid for the invariant representation below.

\begin{theorem}[Exact invariant representation]
\label{amg:invariant}
\label{amg:group-metrization}
Let $(G,\lambda)$ be an abelian convergence group satisfying the $T_1$
condition. There exist a complete ordered value monoid $(M,\oplus,0)$
and a map $p:G\to M$ satisfying
\begin{equation}\label{amg:norm-identities}
 \begin{gathered}
 p(x)=0\iff x=0,\qquad p(-x)=p(x),\\
 p(x+y)\le p(x)\oplus p(y)\qquad(x,y\in G),
 \end{gathered}
\end{equation}
such that, for every net $u$ in $G$, with $\mathcal H=\operatorname{Tail}(u)$,
\begin{equation}\label{amg:null-exactness}
 \mathcal H\longrightarrow_\lambda0
 \quad\Longleftrightarrow\quad
 \bigwedge_{A\in\mathcal H}\bigvee_{x\in A}p(x)=0.
\end{equation}
The invariant symmetric metric
\begin{equation}\label{amg:invariant-distance}
 d(x,y)=p(x-y)
\end{equation}
represents the original canonical uniform convergence exactly:
\begin{equation}\label{amg:pair-exactness}
 \Phi\in\mathcal U_\lambda
 \quad\Longleftrightarrow\quad
 r_d(\Phi)=0,
 \qquad
 \mathcal U_\lambda
 =\{\Phi\in\operatorname{Ev}(G^2):\delta_*\Phi\to_\lambda0\},
 \quad\delta(x,y)=x-y.
\end{equation}
\end{theorem}
\begin{proof}
Put
\[
 I=\{\mathcal F\in\operatorname{Ev}(G):\mathcal F\to_\lambda0\},
 \qquad
 \mathcal C_{\mathcal F}
 =\mathcal F\cap(-\mathcal F)\cap[0]
 \quad(\mathcal F\in I),
\]
where $-\mathcal F=(-\operatorname{id})_*\mathcal F$.
The constant zero net has eventual-set collection $[0]$, so $[0]\in I$ and $I\ne\varnothing$.
Continuity of inversion and finite stability in
Proposition~\ref{fd:tail-calculus} give
\[
 \mathcal F\to_\lambda0,
 \qquad-\mathcal F\to_\lambda0,
 \qquad[0]\to_\lambda0,
 \qquad\mathcal C_{\mathcal F}\to_\lambda0.
\]
For each $\mathcal F\in I$, set
\[
 \mathcal B_{\mathcal F}
 =\{B\in\mathcal C_{\mathcal F}:B=-B\}.
\]
The identities
\[
 -\mathcal C_{\mathcal F}=\mathcal C_{\mathcal F},
 \qquad
 A\in\mathcal C_{\mathcal F}
 \Longrightarrow
 A\cap(-A)\in\mathcal B_{\mathcal F},
 \qquad A\cap(-A)\subseteq A
\]
show that $\mathcal B_{\mathcal F}$ is a nonempty tail base for
$\mathcal C_{\mathcal F}$. Moreover,
\[
 B\in\mathcal B_{\mathcal F}
 \Longrightarrow B\in[0]
 \Longrightarrow0\in B.
\]
Fix $x\ne0$. If $x\in B$ for every
$B\in\mathcal B_{\mathcal F}$, then, for every
$A\in\mathcal C_{\mathcal F}$,
\[
 x\in A\cap(-A)\subseteq A.
\]
Consequently $\mathcal C_{\mathcal F}\subseteq[x]$ and refinement
gives $[x]\to_\lambda0$. The $T_1$ condition gives $x=0$, a
contradiction. Therefore
\begin{equation}\label{amg:point-omission}
 (\forall\mathcal F\in I)(\forall x\in G\setminus\{0\})
 (\exists B\in\mathcal B_{\mathcal F})\qquad x\notin B.
\end{equation}

Apply Lemma~\ref{amu:incidence} with
\[
 Y=G,\qquad Z=\{0\},\qquad
 (\mathcal B_i)_{i\in I}
 =(\mathcal B_{\mathcal F})_{\mathcal F\in I}.
\]
Its hypotheses are the nonemptiness of $I$, the inclusion
$\{0\}\subseteq B$ for every base member, and
\eqref{amg:point-omission}. Its map $p$ satisfies $p(x)=0\iff x=0$.
For $x\ne0$, the $\mathcal F$-coordinate of $p(x)$ is the incidence
set
\[
 \{B\in\mathcal B_{\mathcal F}:x\notin B\}.
\]
Every $B\in\mathcal B_{\mathcal F}$ is symmetric. Thus
\[
 x\notin B\iff-x\notin B,
 \qquad p(-x)=p(x)\quad(x\in G).
\]
The monoid supplied by the lemma satisfies
\[
 a\ne0,\quad b\ne0\quad\Longrightarrow\quad a\oplus b=\top.
\]
If $x=0$ or $y=0$, then
$p(x+y)=p(x)\oplus p(y)$. If $x\ne0$ and $y\ne0$, then
\[
 p(x+y)\le\top=p(x)\oplus p(y).
\]
This proves \eqref{amg:norm-identities}.

For a net $u$ in $G$, with $\mathcal H=\operatorname{Tail}(u)$, Lemma~\ref{amu:incidence}
gives
\begin{equation}\label{amg:incidence-null-tail}
 \bigwedge_{A\in\mathcal H}\bigvee_{x\in A}p(x)=0
 \iff
 \bigl[\{0\}\in\mathcal H\bigr]
 \ \text{or}\ 
 \bigl[(\exists\mathcal F\in I)\ 
              \mathcal B_{\mathcal F}\subseteq\mathcal H\bigr].
\end{equation}
If $\mathcal H\to_\lambda0$, then $\mathcal H\in I$ and
\[
 \mathcal B_{\mathcal H}
 \subseteq\mathcal C_{\mathcal H}
 \subseteq\mathcal H.
\]
The right-hand side of \eqref{amg:incidence-null-tail} holds.
Conversely, $\{0\}\in\mathcal H$ entails
$\mathcal H=[0]\to_\lambda0$. For the remaining alternative,
\[
 \mathcal B_{\mathcal F}\subseteq\mathcal H
 \Longrightarrow\mathcal C_{\mathcal F}\subseteq\mathcal H
 \Longrightarrow\mathcal H\to_\lambda0.
\]
The first implication uses the base property and upward closure of
$\mathcal H$. The second uses refinement of a null-net tail collection.
This proves \eqref{amg:null-exactness}.

For $x,y,z,a\in G$, formula \eqref{amg:invariant-distance} gives
\[
 \begin{aligned}
 d(x,y)=0&\iff x-y=0\iff x=y,\\
 d(y,x)&=p(-(x-y))=d(x,y),\\
 d(x+a,y+a)&=p(x-y)=d(x,y),\\
 d(x,z)&=p((x-y)+(y-z))
       \le d(x,y)\oplus d(y,z).
 \end{aligned}
\]
Thus $d$ is an invariant symmetric metric. For
$\Phi\in\operatorname{Ev}(G^2)$, the family
$\{\delta(C):C\in\Phi\}$ is a base of $\delta_*\Phi$:
\[
 \begin{aligned}
 A\in\delta_*\Phi
 &\iff\delta^{-1}(A)\in\Phi\\
 &\iff(\exists C\in\Phi)\ \delta(C)\subseteq A.
 \end{aligned}
\]
Base invariance of the radius and \eqref{amg:null-exactness} yield
\[
 \begin{aligned}
 r_d(\Phi)
 &=\bigwedge_{C\in\Phi}\bigvee_{(x,y)\in C}p(x-y)\\
 &=\bigwedge_{C\in\Phi}\bigvee_{u\in\delta(C)}p(u)\\
 &=\bigwedge_{A\in\delta_*\Phi}\bigvee_{u\in A}p(u),\\
 r_d(\Phi)=0
 &\iff\delta_*\Phi\to_\lambda0
 \iff\Phi\in\mathcal U_\lambda.
 \end{aligned}
\]
Proposition~\ref{lt:canonical-group} establishes from the group
convergence axioms that $\mathcal U_\lambda$ is globally diagonal,
stable under refinement, finite intersections, inversion and proper
composition, and induces $\lambda$. Its separation also follows
directly:
\[
 [(x,y)]\in\mathcal U_\lambda
 \Longrightarrow[x-y]\to_\lambda0
 \Longrightarrow x=y.
\]
\end{proof}

\begin{proposition}[Solid exact semigroup gauges]
\label{amg:solid}
Let $G$ be an abelian lattice-ordered group and let $\lambda$ be a
$T_1$ locally solid group convergence. The representation in
Theorem~\ref{amg:invariant} can be obtained with the additional property
\begin{equation}\label{amg:solid-gauge}
 |x|\le|y|\quad\Longrightarrow\quad p(x)\le p(y).
\end{equation}
\end{proposition}
\begin{proof}
For each null net $u$, put $\mathcal F=\operatorname{Tail}(u)$, retain
$\mathcal C_{\mathcal F}=\mathcal F\cap(-\mathcal F)\cap[0]$
and put
\[
 \mathcal S_{\mathcal F}
 =\operatorname{Tail}\left\langle\operatorname{sol}(A):
                              A\in\mathcal C_{\mathcal F}\right\rangle.
\]
The displayed generating family is a tail base:
\[
 0\in\operatorname{sol}(A),\qquad
 \operatorname{sol}(A\cap B)
 \subseteq\operatorname{sol}(A)\cap\operatorname{sol}(B).
\]
Moreover,
\[
 A\subseteq\operatorname{sol}(A),\qquad
 \mathcal S_{\mathcal F}\subseteq\mathcal C_{\mathcal F}
                              \subseteq\mathcal F.
\]
To verify null convergence of $\mathcal S_{\mathcal F}$, define
\[
 D_{\mathcal F}
 =\{(A,u,a):A\in\mathcal C_{\mathcal F},\ a\in A,
                                               \ |u|\le|a|\},
\]
with
\[
 (A,u,a)\preceq(B,v,b)\iff B\subseteq A.
\]
This preorder is nonempty and directed, since
$0\in A\cap B\in\mathcal C_{\mathcal F}$ and
$(A\cap B,0,0)$ lies above both indices.
For each $A\in\mathcal C_{\mathcal F}$, the coordinate tail images
above any index whose first coordinate is $A$ are
\[
 \{a:(B,u,a)\in D_{\mathcal F},\ B\subseteq A\}=A,
\]
\[
 \{u:(B,u,a)\in D_{\mathcal F},\ B\subseteq A\}
                                  =\operatorname{sol}(A).
\]
The first equality uses the indices $(A,0,a)$ for $a\in A$. The
second uses $(A,u,a)$ for every $a\in A$ and $|u|\le|a|$.
Thus
\[
 \operatorname{Tail}(a)=\mathcal C_{\mathcal F},\qquad
 \operatorname{Tail}(u)=\mathcal S_{\mathcal F},\qquad
 a\to_\lambda0,\qquad |u|\le|a|.
\]
Local solidity gives $u\to_\lambda0$, hence
$\mathcal S_{\mathcal F}\to_\lambda0$.

Set
\[
 \mathcal B_{\mathcal F}^{\mathrm s}
 =\{\operatorname{sol}(A):A\in\mathcal C_{\mathcal F}\}.
\]
Its members are symmetric, solid and contain zero. For $x\ne0$,
\[
 \begin{aligned}
 [x\in B\text{ for every }B\in\mathcal B_{\mathcal F}^{\mathrm s}]
 &\Longrightarrow\mathcal S_{\mathcal F}\subseteq[x]\\
 &\Longrightarrow[x]\to_\lambda0
 \Longrightarrow x=0.
 \end{aligned}
\]
Therefore every nonzero point is omitted by a member of this base.
Lemma~\ref{amu:incidence} applies to
\[
 Y=G,\quad Z=\{0\},\quad
 I=\{\mathcal F:\mathcal F\to_\lambda0\},\quad
 (\mathcal B_{\mathcal F}^{\mathrm s})_{\mathcal F\in I}.
\]
It supplies $M,p$ with $p^{-1}(0)=\{0\}$ and, for every net $v$ with $\mathcal H=\operatorname{Tail}(v)$,
\[
 \begin{aligned}
 \bigwedge_{A\in\mathcal H}\bigvee_{x\in A}p(x)=0
 &\iff\{0\}\in\mathcal H
 \ \text{or}\ 
 (\exists\mathcal F\in I)\ 
                   \mathcal B_{\mathcal F}^{\mathrm s}\subseteq\mathcal H\\
 &\iff\mathcal H\to_\lambda0.
 \end{aligned}
\]
For the last equivalence, the alternatives on the preceding line give
$\mathcal H=[0]$ or
$\mathcal S_{\mathcal F}\subseteq\mathcal H$, respectively. Both
imply null convergence. Conversely, a null $\mathcal H$ satisfies
\[
 \mathcal B_{\mathcal H}^{\mathrm s}
 \subseteq\mathcal S_{\mathcal H}\subseteq\mathcal H.
\]

If $|x|\le|y|$ and $x,y\ne0$, solidity gives, for every
$B\in\mathcal B_{\mathcal F}^{\mathrm s}$,
\[
 y\in B\Longrightarrow x\in B,
 \qquad
 \{B:x\notin B\}\subseteq\{B:y\notin B\}.
\]
The order in the incidence monoid is coordinatewise inclusion away
from zero. Hence $p(x)\le p(y)$. If $x=0$, the inequality holds
since $p(0)=0$ is least. If $y=0$, then $|x|\le0$ implies $x=0$.
This proves \eqref{amg:solid-gauge}. In particular,
\[
 p(-x)=p(x),\qquad
 p(x+y)\le p(x)\oplus p(y):
\]
the first identity uses $|-x|=|x|$, and the second is equality when
one summand vanishes and has right-hand side $\top$ otherwise.
For $d(x,y)=p(x-y)$, the metric identities and pair exactness are
therefore
\[
 \begin{gathered}
 d(x,y)=0\iff x=y,\qquad d(y,x)=d(x,y),\qquad
 d(x+a,y+a)=d(x,y),\\
 d(x,z)\le d(x,y)\oplus d(y,z),\\
 r_d(\Phi)
 =\bigwedge_{A\in\delta_*\Phi}\bigvee_{u\in A}p(u),\\
 r_d(\Phi)=0\iff\delta_*\Phi\to_\lambda0
 \iff\Phi\in\mathcal U_\lambda.
 \end{gathered}
\]
\end{proof}

\begin{corollary}[Original nets and exact radii]
\label{amg:net-dictionary}
Let $p,d,M$ be given by Theorem~\ref{amg:invariant}. Let
$(x_\alpha)_{\alpha\in D}$ be a net in $G$, let $x\in G$, and let
$\mathcal A\ne\varnothing$ be a parameter set. Then
\begin{align}
 x_\alpha\to_\lambda x
 &\iff
 \bigwedge_{\alpha_0}\bigvee_{\alpha\ge\alpha_0}
             p(x_\alpha-x)=0,\label{amg:point-radius}\\
 (x_\alpha)\text{ is Cauchy for }\lambda
 &\iff
 \bigwedge_{\alpha_0}
 \bigvee_{\alpha,\beta\ge\alpha_0}
             p(x_\alpha-x_\beta)=0.\label{amg:cauchy-radius}
 \end{align}
For a parameter family $(x_\alpha^A,x^A)_{\alpha\in D,A\in\mathcal A}$,
uniform convergence in the original group uniformity satisfies
\begin{equation}\label{amg:uniform-radius}
 x_\alpha^A\to_\lambda x^A\text{ uniformly in }A
 \iff
 \bigwedge_{\alpha_0}
 \bigvee_{\substack{\alpha\ge\alpha_0\\A\in\mathcal A}}
                  p(x_\alpha^A-x^A)=0.
\end{equation}
\end{corollary}
\begin{proof}
For any net $(z_\eta)_{\eta\in E}$ in $G$, its tails form a base of
$\operatorname{Tail}(z)$. Therefore
\[
 \begin{aligned}
 z_\eta\to_\lambda0
 &\iff
 \bigwedge_{B\in\operatorname{Tail}(z)}
                  \bigvee_{u\in B}p(u)=0\\
 &\iff
 \bigwedge_{\eta_0}\bigvee_{\eta\ge\eta_0}p(z_\eta)=0.
 \end{aligned}
\]
For $z_\alpha=x_\alpha-x$, continuity of translation in both
directions gives
\[
 x_\alpha\to_\lambda x\iff z_\alpha\to_\lambda0,
\]
which proves \eqref{amg:point-radius}.

Cauchyness is the convergence of
$z_{\alpha,\beta}=x_\alpha-x_\beta$ to zero on $D\times D$.
The diagonal family $\{(\gamma,\gamma):\gamma\in D\}$ is cofinal:
\[
 (\forall\alpha_0,\beta_0\in D)(\exists\gamma\in D)
 \qquad\gamma\ge\alpha_0,\beta_0.
\]
Hence its tail radius is
\[
 \begin{aligned}
 &\bigwedge_{\alpha_0,\beta_0}
       \bigvee_{\substack{\alpha\ge\alpha_0\\\beta\ge\beta_0}}
                        p(x_\alpha-x_\beta)\\
 &\hspace{25mm}=
 \bigwedge_{\gamma}\bigvee_{\alpha,\beta\ge\gamma}
                        p(x_\alpha-x_\beta).
 \end{aligned}
\]
Substitution proves \eqref{amg:cauchy-radius}.

Order $D\times\mathcal A$ by
$(\alpha,A)\preceq(\beta,B)$ when $\alpha\preceq\beta$.
The pair net $v_{\alpha,A}=(x_\alpha^A,x^A)$ has tail images
\[
 C_{\alpha_0}
 =\{(x_\alpha^A,x^A):\alpha\ge\alpha_0,\ A\in\mathcal A\}.
\]
By the definition of original uniform convergence and
\eqref{amg:pair-exactness},
\[
 \begin{aligned}
 x_\alpha^A\to_\lambda x^A\text{ uniformly in }A
 &\iff\operatorname{Tail}(v)\in\mathcal U_\lambda\\
 &\iff\bigwedge_{\alpha_0}\bigvee_{(u,w)\in C_{\alpha_0}}
                                      d(u,w)=0\\
 &\iff\bigwedge_{\alpha_0}
       \bigvee_{\substack{\alpha\ge\alpha_0\\A\in\mathcal A}}
                                      p(x_\alpha^A-x^A)=0.
 \end{aligned}
\]
This proves \eqref{amg:uniform-radius}.
\end{proof}

In \eqref{amg:uniform-radius}, completeness of the monoid $M$ supplies
the supremum, so no order bound in a lattice group has been imposed.
To obtain lattice-group values, however, we may use the independent
realization criteria of Appendix~\ref{amr:section} without requiring
the resulting faithful representation to factor through this incidence
monoid. Cancellation, weak sigma-distributivity, and completeness of
$(G,\lambda)$ remain outside the assertions of
Theorem~\ref{amg:invariant}.

\subsection{Algebraic lattice realization and complete carriers}
\label{amr:section}

To make that passage to lattice-group values, we first relate the
operations in the two value objects through preservation of lattice
operations and subadditivity. These conditions are formulated before
any approximating net is fixed, leaving detection of the original null
nets to be proved from them. They concern the specified exact value
model and are not asserted for the universal incidence monoids above;
weak distributivity, completeness, and solidity will enter only where
the subsequent arguments explicitly require them.

\begin{proposition}[An algebraic sufficient realization criterion]
\label{lt:algebraic-resolution}
Let $q:Z^2\to M$ be a metric with values in a complete ordered value
monoid, and suppose it satisfies \eqref{lt:exact-monoid}.
Let $F_j$ be Dedekind complete and write
$\overline{F_j^+}=F_j^+\cup\{\infty\}$. Suppose
$a_j:M\to\overline{F_j^+}$ satisfies
\begin{align}
 a_j\!\left(\bigvee S\right)&=\bigvee_{s\in S}a_j(s)
                  &&(S\subseteq M),\label{lt:resolution-joins}\\
 a_j\!\left(\bigwedge T\right)&=\bigwedge_{t\in T}a_j(t)
                  &&(\varnothing\ne T\subseteq M
                       \text{ downward directed}),
                       \label{lt:resolution-meets}\\
 a_j(s\oplus t)&\le a_j(s)+a_j(t),
                       \label{lt:resolution-add}\\
 [a_j(s)=0\ (j\in J)]&\Longrightarrow s=0,
                       \label{lt:resolution-zero}\\
 a_j(q(z,w))&<\infty &&(z,w\in Z).
                       \label{lt:resolution-finite}
\end{align}
Then $d_j=a_j\circ q$ is a faithful order representation.
\end{proposition}
\begin{proof}
The empty join in \eqref{lt:resolution-joins} gives $a_j(0)=0$.
For $s\le t$, the identity $t=s\vee t$ gives $a_j(s)\le a_j(t)$.
Thus \eqref{lt:resolution-add} and
\eqref{lt:resolution-finite} make $d_j$ a finite-valued lattice
pseudometric. Moreover,
\[
 [d_j(z,w)=0\ (j\in J)]
 \Longrightarrow[a_j(q(z,w))=0\ (j\in J)]
 \Longrightarrow q(z,w)=0\Longrightarrow z=w.
\]

For a pair net $v$, put $\Phi=\operatorname{Tail}(v)$ and
\[
 m_C=\bigvee_{(z,w)\in C}q(z,w),\qquad
 t_C^j=\sup_{(z,w)\in C}d_j(z,w)
       \in\overline{F_j^+}\qquad(C\in\Phi).
\]
Finite intersections imply
\[
 m_{C\cap D}\le m_C\wedge m_D,
 \qquad t_{C\cap D}^j\le t_C^j\wedge t_D^j.
\]
Equations \eqref{lt:resolution-joins}--\eqref{lt:resolution-meets}
give
\begin{equation}\label{lt:radius-scalarization}
 a_j(r_q(\Phi))
 =a_j\!\left(\bigwedge_{C\in\Phi}m_C\right)
 =\inf_{C\in\Phi}a_j(m_C)
 =\inf_{C\in\Phi}t_C^j.
\end{equation}
If $r_q(\Phi)=0$, then $\inf_Ct_C^j=0$. At least one $C_j\in\Phi$
satisfies $t_{C_j}^j<\infty$: otherwise every $t_C^j$ equals
$\infty$ and their infimum equals $\infty$. On the cofinal family
\[
 \{C\in\Phi:C\subseteq C_j\},
\]
the values satisfy
\[
 0\le t_C^j\le t_{C_j}^j,\qquad
 \inf_{\substack{C\in\Phi\\C\subseteq C_j}}t_C^j
 =\inf_{C\in\Phi}t_{C\cap C_j}^j
 \le\inf_{C\in\Phi}t_C^j=0.
\]
Moreover, $C\subseteq\{(z,w):d_j(z,w)\le t_C^j\}$ gives
$[0,t_C^j]\in(d_j)_*\Phi$. These decreasing finite values are an
order control for $(d_j)_*\Phi$.

Conversely, if $(d_j)_*\Phi\to_o0$, an order control
$u^j_\theta\downarrow0$ gives
\[
 C^j_\theta=\{(z,w):d_j(z,w)\le u^j_\theta\}\in\Phi,
 \qquad t_{C^j_\theta}^j\le u^j_\theta.
\]
Hence $\inf_Ct_C^j=0$. Equation
\eqref{lt:radius-scalarization} and
\eqref{lt:resolution-zero} give $r_q(\Phi)=0$.
Combine the two implications with \eqref{lt:exact-monoid}.
\end{proof}

For each originally null pair net,
Proposition~\ref{lt:algebraic-resolution} therefore obtains a bounded
image tail in each controller from conditions on joins,
downward-directed infima, and addition in the value objects themselves.
When a bound is already available, the radius test takes the following
simpler form.

\begin{lemma}[Bounded radius criterion]\label{lt:bounded-radius}
Let $F$ be Dedekind complete, let $r:Y\to F^+$, and suppose
$r(Y)\subseteq[0,e]$. For a net $u$ in $Y$, put $\mathcal F=\operatorname{Tail}(u)$ and
\[
 t_C=\sup_{y\in C}r(y)\qquad(C\in\mathcal F).
\]
With eventual sets ordered by reverse inclusion,
\begin{equation}\label{lt:bounded-radius-equivalence}
 r_*\mathcal F\xrightarrow{o}0
 \iff t_C\downarrow0
 \iff\inf_{C\in\mathcal F}\sup_{y\in C}r(y)=0.
\end{equation}
The conclusion also holds when $r$ is bounded on one member of
$\mathcal F$, after restriction to that member.
\end{lemma}
\begin{proof}
All $t_C$ exist in $[0,e]$, and
\[
 D\subseteq C\Longrightarrow t_D\le t_C,
 \qquad C\subseteq\{y:r(y)\le t_C\}.
\]
Thus $t_C\downarrow0$ supplies an order control. Conversely, an order
control $u_\theta\downarrow0$ gives
\[
 C_\theta=\{y:r(y)\le u_\theta\}\in\mathcal F,
 \qquad t_{C_\theta}\le u_\theta.
\]
Hence $\inf_Ct_C=0$. If the bound holds on one eventual set, remove an initial part of the net so that all remaining values belong to that set. The restricted net is tail equivalent to the original one.
\end{proof}

\begin{definition}[Lattice carriers]\label{lt:bounded-resolution}
Let $Z$ be an abelian group, let $\rho_j:Z\to F_j^+$ be symmetric
subadditive maps, and put
\[
 F=\prod_{j\in J}F_j,\qquad
 \rho=(\rho_j)_{j\in J},\qquad
 B_e=\{z\in Z:\rho_j(z)\le e_j\ (j\in J)\}
 \quad(e\in F^+).
\]
The restricted uniform convergence on $B_e$ is induced from the given
uniform convergence of $Z$. No new convergence is assigned to the carrier.
\end{definition}

\begin{proposition}[Algebraic realization on carriers]
\label{lt:resolved-carrier}
Assume that $q(z,w)=p(z-w)$ is an exact invariant metric for
$\mathcal U_\lambda$. Suppose the maps $a_j$ satisfy
\eqref{lt:resolution-joins}--\eqref{lt:resolution-finite}, and put
$\rho_j=a_j\circ p$. Then every $B_e$ has the faithful representation
\[
 d_j(z,w)=\rho_j(z-w),\qquad \chi_j=o.
\]
Moreover, for every pair net $v$ in $B_e^2$, with $\Phi=\operatorname{Tail}(v)$,
\begin{equation}\label{lt:bounded-resolution-equivalence}
 \iota_*\Phi\in\mathcal U_\lambda
 \iff
 \inf_{C\in\Phi}\sup_{(z,w)\in C}\rho_j(z-w)=0
 \quad(j\in J),
\end{equation}
where $\iota:B_e^2\hookrightarrow Z^2$ is inclusion. The point,
Cauchy and parameter-uniform translations in
Proposition~\ref{lt:three-translations} hold for the restricted original
convergence.
\end{proposition}
\begin{proof}
Monotonicity and subadditivity of $a_j$ give
\[
 \rho_j(-z)=\rho_j(z),\qquad
 \rho_j(z+w)\le a_j(p(z)\oplus p(w))
                 \le\rho_j(z)+\rho_j(w).
\]
For $z,w\in B_e$,
\[
 0\le\rho_j(z-w)\le\rho_j(z)+\rho_j(w)\le2e_j.
\]
Proposition~\ref{lt:algebraic-resolution} and
Lemma~\ref{lt:bounded-radius} yield
\[
 \begin{aligned}
 \iota_*\Phi\in\mathcal U_\lambda
 &\iff (d_j)_*(\iota_*\Phi)\to_o0\quad(j\in J)\\
 &\iff \inf_{C\in\Phi}\sup_{(z,w)\in C}
                          \rho_j(z-w)=0\quad(j\in J).
 \end{aligned}
\]
The pair-net assertion includes the three specific nets used in
Proposition~\ref{lt:three-translations}, so all three restricted
translations hold.
\end{proof}

\begin{lemma}[Bounded product and coordinate controls]
\label{lt:bounded-product}
Let $F=\prod_{j\in J}F_j$, where every $F_j$ is Dedekind complete,
and let a net $(v_\alpha)$ in $F^+$ be eventually bounded by $e\in F^+$.
Then
\begin{equation}\label{lt:bounded-product-equivalence}
 v_\alpha\xrightarrow{o_F}0
 \quad\Longleftrightarrow\quad
 \pi_j(v_\alpha)\xrightarrow{o_{F_j}}0\quad(j\in J).
\end{equation}
\end{lemma}
\begin{proof}
Coordinate projection preserves decreasing zero controls, proving the
forward implication. For the converse, restrict to a tail on which
$0\le v_\alpha\le e$, and put
\[
 s_{\alpha_0}(j)=\sup_{\alpha\ge\alpha_0}v_\alpha(j),
 \qquad s_{\alpha_0}=(s_{\alpha_0}(j))_{j\in J}.
\]
These suprema exist, $s_{\alpha_0}\in[0,e]$, the net is decreasing,
and $v_\alpha\le s_{\alpha_0}$ for $\alpha\ge\alpha_0$.
For each coordinate control $u^j_\theta\downarrow0$, some tail
satisfies $v_\alpha(j)\le u^j_\theta$, and hence
$s_{\alpha_0}(j)\le u^j_\theta$ on a tail. Thus
$\inf_{\alpha_0}s_{\alpha_0}(j)=0$ for every $j$.
Infima in a product order are coordinatewise, so
$s_{\alpha_0}\downarrow_F0$, giving an order control for $v$.
\end{proof}

\begin{proposition}[Carrier stability and completeness]
\label{lt:sandwich-resolution}
\label{lt:carrier-completeness-sandwich}
Let $(Z,\lambda)$ be a separated abelian convergence group whose
canonical uniform structure has a faithful order representation
$d_j(z,w)=\rho_j(z-w)$. The algebraic conditions of
Proposition~\ref{lt:resolved-carrier} suffice. Every net $u$ satisfies
\begin{equation}\label{lt:carrier-limit-stability}
 u_\alpha\in B_e\text{ eventually},\quad u_\alpha\to_\lambda z
 \quad\Longrightarrow\quad z\in B_e.
\end{equation}
If $Z$ is complete, each $B_e$ is complete for its restricted original
uniform convergence.
\end{proposition}
\begin{proof}
Remove an initial part of $u$ so that $u_\alpha\in B_e$ throughout.
For each $j$, faithfulness supplies controls $a^j_\theta\downarrow0$
for $\rho_j(u_\alpha-z)$. On the tail corresponding to $\theta$,
\[
 \rho_j(z)\le\rho_j(z-u_\alpha)+\rho_j(u_\alpha)
                         \le a^j_\theta+e_j.
\]
Consequently $0\le(\rho_j(z)-e_j)^+\le a^j_\theta$ for every
$\theta$, so $\rho_j(z)\le e_j$. This holds for every $j$ and proves
\eqref{lt:carrier-limit-stability}.

Let $u$ be a Cauchy net in $B_e$. Its inclusion into $Z$ is Cauchy by
the definition of the restricted uniform convergence. Completeness of
$Z$ gives a limit $z\in Z$, and the first assertion gives $z\in B_e$.
The induced convergence on $B_e$ then gives convergence to $z$ in
$B_e$, proving its completeness.
\end{proof}

\begin{proposition}[Intrinsic order and unbounded order realizations]
\label{fd:intrinsic-order-realizations}
Let $G$ be an abelian lattice-ordered group. Its order
convergence has the faithful invariant gauge $\rho(x)=|x|$. Its unbounded
order convergence,
\[
 x_\alpha\to_{uo}x
 \quad\Longleftrightarrow\quad
 |x_\alpha-x|\wedge u\to_o0\quad(u\in G^+),
\]
has the faithful invariant family
\[
 \rho_u(x)=|x|\wedge u\qquad(u\in G^+).
\]
Both are separated locally solid group convergences. The second family
is bounded coordinatewise on the whole group. In each case faithfulness
holds for the original canonical pair convergence, including Cauchy
and parameter-uniform convergence.
\end{proposition}
\begin{proof}
We first verify subadditivity of the truncations. Let $a,b,u\in G^+$ and
put
\[
 r=(a+b)\wedge u,\qquad c=r\wedge a,\qquad d=r-c=(r-a)^+.
\]
The defining inequalities give
\[
 0\le c\le a\wedge u,\qquad
 0\le d\le b,\qquad d\le r\le u,
\]
and therefore
\begin{equation}\label{fd:truncation-subadditivity}
 (a+b)\wedge u=c+d\le(a\wedge u)+(b\wedge u).
\end{equation}
For $x,y\in G$,
\[
 \begin{aligned}
 \rho_u(x+y)
 &\le(|x|+|y|)\wedge u\\
 &\le(|x|\wedge u)+(|y|\wedge u)
 =\rho_u(x)+\rho_u(y),\\
 \rho_u(-x)&=\rho_u(x),\qquad 0\le\rho_u(x)\le u.
 \end{aligned}
\]
Also
\[
 [\rho_u(x)=0\ (u\in G^+)]
 \Longrightarrow \rho_{|x|}(x)=|x|=0
 \Longrightarrow x=0.
\]

For a positive order-null net $b_\alpha$ with controls
$r_\theta\downarrow0$, the inequalities
$0\le a_\alpha\le b_\alpha$ give the same eventual bounds
$a_\alpha\le r_\theta$. A quasi-subnet inherits each such bound,
because its tail is eventually contained in the corresponding
original tail. For two order-null nets with controls
$r_\theta\downarrow0$ and $s_\eta\downarrow0$, their independent
sum has controls $(r_\theta+s_\eta)_{\theta,\eta}\downarrow0$
by \eqref{lt:sum-controls}. These facts give
\[
 \begin{aligned}
 x_\alpha\to_{uo}x,\quad y_\beta\to_{uo}y
 &\Longrightarrow
 \rho_u((x_\alpha+y_\beta)-(x+y))\\
 &\hspace{1.8em}\le\rho_u(x_\alpha-x)+\rho_u(y_\beta-y)
       \to_o0\quad(u\in G^+),\\
 x_\alpha\to_{uo}x
 &\Longrightarrow\rho_u(-x_\alpha+x)
            =\rho_u(x_\alpha-x)\to_o0,\\
 |v_\alpha|\le|w_\alpha|,\quad w_\alpha\to_{uo}0
 &\Longrightarrow 0\le\rho_u(v_\alpha)
                \le\rho_u(w_\alpha)\to_o0.
 \end{aligned}
\]
Constant nets have zero error. If $x_i\to_{uo}x$ and
$y_j\to_{uo}x$, define the mixing net on
$I\times J\times\{0,1\}$, ordered through $I\times J$, by
$m_{i,j,0}=x_i$ and $m_{i,j,1}=y_j$. For each $u\in G_+$,
\[
 \rho_u(m_{i,j,\varepsilon}-x)
 \le\rho_u(x_i-x)+\rho_u(y_j-x)\xrightarrow{o}0.
\]
The preceding quasi-subnet calculation and this mixing estimate verify
the convergence axioms. The solid-domination implication in the
previous display proves local solidity. If a constant net with value
$z$ is unbounded-order null, its $u$-gauge is bounded by every member
of a decreasing zero control and is therefore zero. Setting $u=|z|$
gives $|z|=0$. If a net has limits $x,y$, subtraction on the diagonal
index makes the constant difference $x-y$ null, so $x=y$.

For the untruncated gauge, the formulas
\[
 |(x_i+y_j)-(x+y)|\le|x_i-x|+|y_j-y|,\qquad
 |-x_i+x|=|x_i-x|,
\]
\[
 |m_{i,j,\varepsilon}-x|
 \le|x_i-x|+|y_j-x|,\qquad
 0\le|v_i|\le|w_i|
\]
give addition, inversion, mixing and solid domination for order
convergence with the controls just constructed. A constant order-null
error is zero by the infimum condition, which proves separation.

Let $\delta(x,y)=x-y$ and let $v$ be a pair net, with $\Phi=\operatorname{Tail}(v)$. For
unbounded order convergence,
\[
 \begin{aligned}
 \Phi\in\mathcal U_{uo}
 &\iff\delta_*\Phi\to_{uo}0\\
 &\iff (\rho_u)_*(\delta_*\Phi)\to_o0\quad(u\in G^+)\\
 &\iff (d_u)_*\Phi\to_o0\quad(u\in G^+),
 \qquad d_u(x,y)=|x-y|\wedge u.
 \end{aligned}
\]
For order convergence, the corresponding calculation is
\[
 \Phi\in\mathcal U_o
 \iff\delta_*\Phi\to_o0
 \iff\bigl((x,y)\mapsto|x-y|\bigr)_*\Phi\to_o0.
\]
The group convergence properties proved above permit application of
Proposition~\ref{lt:canonical-group}. Applying these pair equivalences
to $(x_i,x)$, $(x_i,x_j)$ and
$(x_i^A,x^A)_{(i,A)}$ gives, respectively, point, Cauchy and
parameter-uniform faithfulness.
\end{proof}

For order and unbounded-order convergence, the lattice-group operations
thus provide the realizations directly. Passing further to regulator
convergence invokes the additional hypotheses of the order appendix,
where weak distributivity and countable-sup assumptions govern that
comparison independently of the metrization.

\begin{proposition}[Complete carriers and existing limits]
\label{lt:complete-carrier}
Let $B\subseteq Z$ be complete for the restricted original uniform
convergence structure. If a net $u_\alpha\to z$ in the separated
ambient space is eventually in $B$, then $z\in B$.
\end{proposition}
\begin{proof}
The pair net $(u_\alpha,z)_\alpha$ is originally uniformly null,
and symmetry gives uniform nullity of $(z,u_\beta)_\beta$.
For their tail-value relations,
\[
 R_{\alpha_0}=T_{\alpha_0}(u)\times\{z\},\qquad
 S_{\beta_0}=\{z\}\times T_{\beta_0}(u),\qquad
 R_{\alpha_0}\circ S_{\beta_0}
       =T_{\alpha_0}(u)\times T_{\beta_0}(u).
\]
Every composition is nonempty, and these are the tail values of
$(u_\alpha,u_\beta)_{\alpha,\beta}$. The uniform composition axiom
therefore makes this pair net uniformly null. If
$u_\alpha\in B$ for $\alpha\ge\alpha_*$, its restriction to that
tail is Cauchy in the induced uniform convergence on $B$. Completeness
of $B$ supplies a limit $b\in B$, which is also an ambient limit.
Separation gives $z=b$.
\end{proof}

This stability allows the later complete-target theorems to be used
within a complete resolved carrier containing the finite sums:
Proposition~\ref{lt:complete-carrier} places their limits in the carrier
without introducing an adherence operation.

\subsection{Closed bases and closed-ball realizations}
\label{nx:closed-base-section}

Having considered algebraic conditions on the values, we now turn to
the stability of the original tails under closure. Throughout the
following subsections, convergence groups are abelian and every
eventual-set collection comes from an explicitly considered net. For
a partially ordered abelian group $P$, the notation $D\downarrow0$
means that $D\subseteq P_+$ is nonempty and downward directed and
that $\inf_P D=0$, with the infimum taken in the whole group:
\[
 v\le d\ (d\in D)\quad\Longrightarrow\quad v\le0
 \qquad(v\in P).
\]
No restriction to positive lower bounds is implicit. For a positive-valued map $v:Y\to P_+$, the notation $v_*\mathcal F\to_o0$ means that some $D\downarrow0$ satisfies $\{y:v(y)\le d\}\in\mathcal F$ for every $d\in D$.

For a derived collection
$\mathcal F=\operatorname{Tail}(u)$, we consequently write
\[
 K_X\mathcal F
 =\operatorname{Tail}\langle K_X(A):A\in\mathcal F\rangle.
\]
Since the tail images $T_i(u)$ are cofinal among the eventual sets
under reverse inclusion, this notation gives a saturation tail
equivalent to the one in \eqref{nx:CB-equation} of
Definition~\ref{nx:CB}.

\begin{definition}[Disks and local convexity]
In a real convergence vector space,
a disk is a nonempty balanced convex set, and local convexity means
that every null net is a quasi-subnet of the saturation of a null
family of disks. Write $\operatorname{aco}A$ for the smallest disk
containing $A$, with $\operatorname{aco}\varnothing=\{0\}$.
\end{definition}

\begin{lemma}[Closure calculus]\label{nx:closure-calculus}
The operator $K_X$ is extensive, monotone and idempotent, fixes the empty
set and preserves finite unions. For every continuous map $T:X\to Y$,
\begin{equation}\label{nx:continuous-closure}
 T(K_X(A))\subseteq K_Y(T(A)).
\end{equation}
Condition \textnormal{(CB)} is equivalent to the assertion that every
convergent net is a quasi-subnet of the saturation of a family of closed
sets converging to the same point. In a real convergence vector space,
$K_X(D)$ is a disk whenever $D$ is a disk, and local convexity is
equivalent to
\begin{equation}\label{nx:aco-null}
 u_i\longrightarrow0
 \quad\Longrightarrow\quad
 \langle\operatorname{aco}(T_i(u)):i\in I\rangle\longrightarrow0.
\end{equation}
Under \textnormal{(CB)} and local convexity, the smallest closed disk hull
$H(A)=K_X(\operatorname{aco}A)$ preserves null saturated families.
\end{lemma}
\begin{proof}
Arbitrary intersections of closed sets are closed. Let $u_i\to x$ be
eventually in $C\cup D$, with $C,D$ closed. If some tail misses $C$,
the restricted net lies in $D$ and $x\in D$. Otherwise every tail of
the restricted net meets $C$. The saturation of the nonempty sets
$T_i(u)\cap C$ is a quasi-subnet of $u$, has all values in $C$, and
converges to $x$, so $x\in C$. Thus $C\cup D$ is closed. Taking
intersections of closed supersets gives a closed set $K_X(A)$
containing $A$, and therefore
\[
 K_X(K_X(A))=K_X(A),\qquad A\subseteq K_X(A).
\]
If $A\subseteq B$, every closed superset of $B$ contains $A$, so
$K_X(A)\subseteq K_X(B)$. The closed set
$K_X(A)\cup K_X(B)$ contains $A\cup B$, while monotonicity gives
$K_X(A),K_X(B)\subseteq K_X(A\cup B)$. These two inclusions yield
$K_X(A\cup B)=K_X(A)\cup K_X(B)$. The empty set is closed, hence
$K_X(\varnothing)=\varnothing$. The inverse
image of a closed set under $T$ is closed, since $u_i\to x$ entails
$T(u_i)\to T(x)$. Apply this observation to the closed set
$K_Y(T(A))$ to obtain \eqref{nx:continuous-closure}.

If \textnormal{(CB)} holds, the family in \eqref{nx:CB-equation} is a
closed family with the required limit, and $u$ is its quasi-subnet
because $T_i(u)\subseteq K_X(T_i(u))$. Conversely, suppose a closed
family $(C_\eta)$ has saturation converging to $x$ and every $C_\eta$ contains a
tail image $T_{i(\eta)}(u)$. Then
\[
 K_X(T_{i(\eta)}(u))\subseteq K_X(C_\eta)=C_\eta.
\]
The saturation in \eqref{nx:CB-equation} is therefore a quasi-subnet of
$\langle C_\eta\rangle$, proving the asserted equivalence at the prescribed limit $x$.

For a disk $D$, multiplication by $a$ with $|a|\le1$ and
\eqref{nx:continuous-closure} give
\[
 aK_X(D)\subseteq K_X(aD)\subseteq K_X(D).
\]
For $0\le t\le1$ and fixed $x\in D$, the continuous affine map
$y\mapsto tx+(1-t)y$ sends $D$ into $D$, and hence sends $K_X(D)$ into
$K_X(D)$. For fixed $y\in K_X(D)$, the continuous map
$x\mapsto tx+(1-t)y$ now sends $D$ into the closed set $K_X(D)$,
and therefore sends $K_X(D)$ into the same set. This proves
\[
 tK_X(D)+(1-t)K_X(D)\subseteq K_X(D).
\]
If $u$ is a quasi-subnet of a null disk family $(D_\eta)$, every
$D_\eta$ contains some $T_i(u)$, and then
$\operatorname{aco}(T_i(u))\subseteq D_\eta$. This proves
\eqref{nx:aco-null}; its converse uses its own displayed disk family.
Applying \eqref{nx:aco-null} and \eqref{nx:CB-equation} in succession
proves nullity of the saturation of $H(T_i(u))$. Every closed disk
containing $A$ contains $\operatorname{aco}A$ and consequently
$H(A)$, which proves minimality.
\end{proof}

\begin{proposition}[Ordered-group balls]\label{nx:closed-balls}
Let $X$ be a convergence group and let $p:X\to P_+$ be symmetric and
subadditive, where $P$ is a partially ordered abelian group. Suppose
\[
 u_i\to0\quad\Longrightarrow\quad p(u_i)\to_o0
\]
for every net $u$. Then each $B_r=\{x:p(x)\le r\}$ and every
translate of $B_r$ is closed. If a family
$(p_j:X\to(P_j)_+)_j$ represents the original null nets faithfully,
then $X$ satisfies \textnormal{(CB)}.
\end{proposition}
\begin{proof}
Let $u_i\in c+B_r$ and $u_i\to x$. The hypothesis supplies
$D\downarrow0$ such that, for each $d\in D$, a tail satisfies
$p(u_i-x)\le d$. On that tail,
\[
 p(x-c)\le p(x-u_i)+p(u_i-c)\le d+r.
\]
Thus $p(x-c)-r$ is a lower bound of $D$, so $p(x-c)\le r$ and
$x\in c+B_r$.

Let now $u_i\to0$, and take controls $D_j\downarrow0$ for
$p_j(u_i)$. For every $d\in D_j$, an index $i(j,d)$ satisfies
$T_{i(j,d)}(u)\subseteq B_{j,d}$. Closedness of these balls gives
\[
 K_X(T_{i(j,d)}(u))\subseteq K_X(B_{j,d})=B_{j,d}.
\]
The same $D_j$ consequently controls the image under $p_j$ of the
saturation of $(K_X(T_i(u)))_i$. Faithfulness makes this saturation
null. Translations and their inverses are continuous, so
$K_X(A+x)=K_X(A)+x$, and the corresponding assertion holds for nets
converging to every point. This is (CB).
\end{proof}

\begin{theorem}[Closed-ball monoid realization]\label{nx:closed-ball-monoid}
A separated convergence group satisfies \textnormal{(CB)} if and only if
it admits a faithful invariant symmetric metric with values in a complete
ordered commutative monoid whose every ball is closed. For a real
convergence vector space, \textnormal{(CB)} and local convexity hold if
and only if such a metric can have every ball a closed disk.
\end{theorem}
\begin{proof}
Let $p:X\to M$ be a faithful gauge with closed balls. For
$A\in\mathcal F\to0$, define $r_A=\bigvee_{x\in A}p(x)$.
By faithfulness and the complete-lattice radius formula,
\[
 A'\subseteq A\ \Longrightarrow\ r_{A'}\le r_A,
 \qquad\bigwedge_{A\in\mathcal F}r_A=0.
\]
Since $A\subseteq B_{r_A}$ and $B_{r_A}$ is closed,
$K_X(A)\subseteq B_{r_A}$. Therefore
\[
 0\le\bigwedge_{A\in\mathcal F}
       \bigvee_{x\in K_X(A)}p(x)
 \le\bigwedge_{A\in\mathcal F}r_A=0.
\]
Faithfulness yields \textnormal{(CB)}. If the balls are disks,
$A\subseteq B_{r_A}$ implies $\operatorname{aco}A\subseteq B_{r_A}$,
and hence
\[
 0\le\bigwedge_{A\in\mathcal F}
           \bigvee_{x\in\operatorname{aco}A}p(x)
 \le\bigwedge_{A\in\mathcal F}r_A=0.
\]
Faithfulness makes the saturation of
$(\operatorname{aco}A)_{A\in\mathcal F}$ null.
The tail images of the original net are cofinal in $\mathcal F$,
so this is \eqref{nx:aco-null}.

Conversely, let $\mathfrak N=\{\operatorname{Tail}(u):u\to0\}$, with $u$ ranging over admissible nets. For
$F\in\mathfrak N$, put
\[
 \mathcal B_F=\{K_X(A\cup(-A)\cup\{0\}):A\in F\}.
\]
The saturation of the sets $A\cup(-A)\cup\{0\}$ is the mixing of a null representative of $F$, its negative and the zero net. Its eventual sets are $F\cap(-F)\cap[0]$. Condition \textnormal{(CB)} makes the saturation of $\mathcal B_F$ null. All its base sets are closed,
symmetric, and contain zero. They belong to $F$. In the vector case
with local convexity, use instead
$\mathcal B_F=\{K_X(\operatorname{aco}A):A\in F\}$.
Lemma~\ref{nx:closure-calculus} gives closed disks and nullity.
Separation gives
\[
 x\ne0\quad\Longrightarrow\quad
 \{B\in\mathcal B_F:x\notin B\}\ne\varnothing
 \qquad(F\in\mathfrak N).
\]
Lemma~\ref{amu:incidence}, with $Y=X$ and $Z=\{0\}$, supplies
\[
 M=\{0\}\sqcup\prod_{F\in\mathfrak N}
                  \mathcal P_*(\mathcal B_F),\qquad
 p(x)_F=\{B\in\mathcal B_F:x\notin B\}\quad(x\ne0),
 \qquad p(0)=0.
\]
Here nonzero coordinates are ordered by inclusion, and
$a\oplus b=\top$ if $a,b\ne0$. Symmetry of every $B$ gives
$p(-x)=p(x)$. A zero summand gives equality in subadditivity. Two
nonzero summands give $p(x+y)\le\top=p(x)\oplus p(y)$.
For every net $v$, with $\mathcal H=\operatorname{Tail}(v)$, the incidence formula is
\[
 \bigwedge_{A\in\mathcal H}\bigvee_{x\in A}p(x)=0
 \iff \{0\}\in\mathcal H\ \text{or}\
       (\exists F\in\mathfrak N)\ \mathcal B_F\subseteq\mathcal H
 \iff\mathcal H\to0.
\]
In the forward direction the displayed alternatives provide a coarser
null-net tail collection. In the reverse direction use $F=\mathcal H$.
For a nonzero radius $a=(a_F)_F$, coordinate inclusion gives
\begin{equation}\label{nx:incidence-ball}
 B_a=\bigcap_{F\in\mathfrak N}
             \bigcap_{B\in\mathcal B_F\setminus a_F} B.
\end{equation}
Thus every ball is closed, and in the vector construction every ball
is a disk. The zero ball is $\{0\}$. The distance $p(x-y)$ has all
required metric properties and represents every original net.
\end{proof}

\begin{proposition}[Permanence]\label{nx:CB-permanence}
An initial convergence induced by maps $T_j:X\to X_j$, where each
$X_j$ satisfies \textnormal{(CB)}, satisfies \textnormal{(CB)}.
In particular, products and subspaces preserve this property. Initial
convergence under linear maps into locally convex real convergence
vector spaces preserves local convexity.
\end{proposition}
\begin{proof}
Let $u_i\to x$. For each $j$, the image net $T_j(u_i)$ converges to
$T_j(x)$, and \eqref{nx:continuous-closure} gives
\[
 T_j(K_X(T_i(u)))\subseteq K_{X_j}(T_j(T_i(u))).
\]
The set $T_j(T_i(u))$ is exactly the image of the $i$-th tail of the
net $T_j(u)$. By (CB) in $X_j$, the saturation of its closed tail
images converges to $T_j(x)$. The displayed inclusion makes the
$T_j$-image of the saturation of $(K_X(T_i(u)))_i$ a quasi-subnet of
that net. Initiality therefore gives convergence to $x$, proving (CB).
For a linear $T_j$,
\[
 T_j(\operatorname{aco}(T_i(u)))
             =\operatorname{aco}(T_j(T_i(u))).
\]
If $u_i\to0$, local convexity in each $X_j$ makes the saturation on
the right null. Initiality and \eqref{nx:aco-null} prove local convexity
of $X$.
\end{proof}

\subsection{Two-scale convergence and uniform mapping groups}
\label{nx:mapping-section}

Closure stability concerns replacing each tail set by its
convergence-closed hull. We now return to the absorption of independently null rows,
whose role can also be seen through continuity of uniform limits of
maps. The following formulation keeps both scales of the approximation
in the original convergence.

\begin{definition}[Two-scale net property]\label{nx:P}
A convergence group has property \textnormal{(P)} when the following
implication holds. Let $(A_\lambda)_{\lambda\in\Lambda}$ be decreasing
nonempty sets with null saturation, let $z=(z_i)_{i\in I}$ be a net,
and, for each $\lambda$, let
$h^\lambda=(h^\lambda_j)_{j\in J_\lambda}\to0$. Then
\begin{equation}\label{nx:P-equation}
 \left[
 \forall\lambda\ \forall j\in J_\lambda\ \exists i_0\in I\
 \forall i\ge i_0:\quad z_i\in A_\lambda+T_j(h^\lambda)
 \right]
 \quad\Longrightarrow\quad z_i\to0.
\end{equation}
This is the unrestricted abelian form of
Definition~\ref{ci:two-scale}. For a set $S$, uniform convergence in
$X^S$ means
\[
 f_i\to_u f
 \quad\Longleftrightarrow\quad
 \left\langle
 \{f_j(s)-f(s):j\ge i,\ s\in S\}:i\in I
 \right\rangle\to0
\]
when $S\ne\varnothing$. For a convergence space $S$, the continuous
maps with the restricted uniform convergence form $C_u(S,X)$.
The empty source gives the one-element mapping group.
\end{definition}

\begin{theorem}[Uniform-limit characterization]\label{nx:uniform-limit}
For an abelian convergence group $X$, the following are equivalent:
\begin{enumerate}[(i)]
\item $X$ satisfies \textnormal{(P)}.
\item For every convergence space $S$, every uniformly convergent net
of continuous maps $S\to X$ has a continuous limit.
\item The assertion in (ii) holds for every Hausdorff topological
source having at most one nonisolated point.
\end{enumerate}
No topology is assigned to $X$.
\end{theorem}
\begin{proof}
Assume (i), and let $f_i:S\to X$ be continuous with $f_i\to_u f$.
The empty source is immediate, so assume $S\ne\varnothing$. Put
\[
 E_i=\{f(s)-f_j(s):j\ge i,\ s\in S\},\qquad A_i=E_i-E_i.
\]
Uniform convergence and continuity of inversion make the saturation
of $(E_i)_i$ null. The independent difference of two copies of this
null net has tail sets $E_i-E_j$. For every pair $(i,j)$, directedness
supplies $k\ge i,j$ with
\[
 A_k=E_k-E_k\subseteq E_i-E_j.
\]
Thus the saturation of $(A_k)_k$ is a quasi-subnet of that independent
difference net and is null. Fix an arbitrary net
$s_\alpha\to s_0$ in $S$, and define
\[
 z_\alpha=f(s_\alpha)-f(s_0),\qquad
 h^i_\alpha=f_i(s_\alpha)-f_i(s_0).
\]
Continuity of $f_i$ gives $h^i_\alpha\to0$. With $r_i=f-f_i$,
\[
 z_\alpha=r_i(s_\alpha)-r_i(s_0)+h^i_\alpha
                 \in A_i+\{h^i_\alpha\}.
\]
For each $i$ and each $\alpha_0$, all $\alpha\ge\alpha_0$ therefore
satisfy $z_\alpha\in A_i+T_{\alpha_0}(h^i)$. Property (P) gives
$z_\alpha\to0$, so $f(s_\alpha)\to f(s_0)$. Since the net and its
limit were arbitrary, $f$ is continuous. This proves (ii), which
implies (iii).

Assume (iii), and let the data of \eqref{nx:P-equation} be given.
Define
\begin{equation}\label{nx:uniform-test-1}
 T=\{(x,v)\in X\times X^\Lambda:
                   x-v_\lambda\in A_\lambda\ (\lambda\in\Lambda)\}.
\end{equation}
For $i\in I$, a finite $F\subseteq\Lambda$, and a tuple
$j=(j_\lambda)_{\lambda\in F}$, put
\begin{equation}\label{nx:uniform-test-2}
 L(i,F,j)=\{(x,v)\in T:x\in T_i(z),\
                   v_\lambda\in T_{j_\lambda}(h^\lambda)
                         \ (\lambda\in F)\}.
\end{equation}
The first-coordinate projection of this set is
\begin{equation}\label{nx:uniform-test-3}
 T_i(z)\cap\bigcap_{\lambda\in F}
                 (A_\lambda+T_{j_\lambda}(h^\lambda)).
\end{equation}
Indeed, membership in \eqref{nx:uniform-test-2} gives the displayed
inclusions. Conversely, a point $x$ in \eqref{nx:uniform-test-3}
has representations $x=a_\lambda+b_\lambda$ with
$a_\lambda\in A_\lambda$ and
$b_\lambda\in T_{j_\lambda}(h^\lambda)$ for $\lambda\in F$.
Set $v_\lambda=b_\lambda$ there. On the remaining coordinates take
$a_\lambda\in A_\lambda$ and put $v_\lambda=x-a_\lambda$.
These assignments produce a point of \eqref{nx:uniform-test-2}.
For each $\lambda\in F$, the antecedent of (P) supplies
$i_\lambda\in I$ such that
\[
 k\ge i_\lambda\quad\Longrightarrow\quad
 z_k\in A_\lambda+T_{j_\lambda}(h^\lambda).
\]
A common upper index $i_1\ge i,i_\lambda$ for $\lambda\in F$
then gives
\[
 \varnothing\ne T_{i_1}(z)
 \subseteq T_i(z)\cap\bigcap_{\lambda\in F}
                   (A_\lambda+T_{j_\lambda}(h^\lambda)).
\]
This proves nonemptiness for every displayed finite constraint.

The family of all sets \eqref{nx:uniform-test-2} is downward
directed. For two such sets, replace $F$ by its finite union,
replace $i$ by a common upper index, and replace each repeated
$j_\lambda$ by a common upper index in $J_\lambda$. The resulting
set is contained in their intersection. Its saturation $\ell$ is
therefore defined. Write $\mathscr L$ for this family and use
$D_{\mathscr L}=\{(L,t):L\in\mathscr L,\ t\in L\}$, ordered by
reverse inclusion of $L$, as the index of $\ell$, with
$\ell_{(L,t)}=t$. The first-coordinate tail above $(L,t)$ is
$\pi_xL$. Formula \eqref{nx:uniform-test-3} gives
\[
 \pi_xL(i,F,j)\subseteq T_i(z),\qquad
 T_{i_1}(z)\subseteq\pi_xL(i,F,j).
\]
For every original tail $T_i(z)$, the unconstrained member
$L(i,\varnothing,\varnothing)$ has projection $T_i(z)$.
Thus each net is a quasi-subnet of the other. For a fixed $\lambda$
and $j\in J_\lambda$, any constraint with
$F=\{\lambda\}$ and $j_\lambda=j$ has
\[
 \pi_\lambda L\subseteq T_j(h^\lambda).
\]
These inclusions make $v_\lambda(\ell)$ a quasi-subnet of the
null net $h^\lambda$.

Put $S=(T\times\mathbb N)\sqcup\{\infty\}$. All points of
$T\times\mathbb N$ are isolated, and give $\infty$ the neighborhood
base
\begin{equation}\label{nx:uniform-test-4}
 N(L,n)=\{\infty\}\cup(L\times\{m\in\mathbb N:m\ge n\}),
\end{equation}
where $L$ ranges over the sets in \eqref{nx:uniform-test-2}.
Downward directedness of those sets makes \eqref{nx:uniform-test-4}
a neighborhood base. For a point with integer coordinate $m$,
its singleton and the open neighborhood $N(T,m+1)$ are disjoint. The latter contains $N(L,m+1)$ for every basic $L$, so the space is Hausdorff. Define
\begin{equation}\label{nx:uniform-test-5}
 \begin{gathered}
 f_\lambda(x,v,n)=v_\lambda,\qquad f_\lambda(\infty)=0,\\
 f(x,v,n)=x,\qquad f(\infty)=0.
 \end{gathered}
\end{equation}
To verify continuity of $f_\lambda$ at $\infty$, let
$s_\alpha\to\infty$. For every $j\in J_\lambda$, a neighborhood
in \eqref{nx:uniform-test-4} can be taken with
$F=\{\lambda\}$ and $j_\lambda=j$. On the corresponding tail,
\[
 f_\lambda(s_\alpha)\in T_j(h^\lambda)\cup\{0\}.
\]
The saturation of these sets is the mixing of $h^\lambda$ and
the constant zero net and is null. Hence $f_\lambda(s_\alpha)\to0$.
If $s\in T\times\mathbb N$ is isolated and $s_\alpha\to s$,
eventual membership in the open singleton $\{s\}$ gives
$f_\lambda(s_\alpha)=f_\lambda(s)$ eventually, proving continuity
at $s$.

For $\eta\ge\lambda$, decreasingness of $(A_\lambda)$ gives
\[
 f(s)-f_\eta(s)\in A_\lambda\cup\{0\}\qquad(s\in S).
\]
The centered family $(A_\lambda\cup\{0\})$ has null saturation,
so $f_\lambda\to_u f$. Assumption (iii) makes $f$ continuous.
Give $D_{\mathscr L}\times\mathbb N$ the product order and put
\[
 s_{((L,t),n)}=(t,n)\in T\times\mathbb N\subset S.
\]
For a neighborhood $N(L_0,n_0)$, fix $t_0\in L_0$.
Above $((L_0,t_0),n_0)$ one has $L\subseteq L_0$ and $n\ge n_0$,
so $s_{((L,t),n)}\in N(L_0,n_0)$. Hence this net converges to
$\infty$. Its image under $f$ has tail values
\[
 \{f(s_{((L,t),n)}):L\subseteq L_0,\ t\in L,\ n\ge n_0\}
 =\pi_xL_0.
\]
The inclusion from left to right uses $L\subseteq L_0$.
For the reverse inclusion use $L=L_0$, arbitrary $t\in L_0$,
and $n=n_0$. Continuity of $f$ makes this image null, and the
proved tail equivalence with $(z_i)$ gives $z_i\to0$.
\end{proof}

Having characterized when a uniform limit remains continuous, we now
consider existence for uniformly Cauchy nets in the full function
space. The next argument constructs the limit pointwise and uses
closed-base stability to recover the common uniform estimate.

\begin{lemma}[Uniform completeness]\label{nx:uniform-mapping-complete}
If $X$ is a complete separated convergence group satisfying
\textnormal{(CB)}, then $X^S$ is complete for uniform convergence for
every set $S$.
\end{lemma}
\begin{proof}
Assume $S\ne\varnothing$ and let $(f_i)_{i\in I}$ be uniformly
Cauchy. The decreasing sets
\[
 D_{i_0}=\{f_i(s)-f_j(s):i,j\ge i_0,\ s\in S\}
\]
have null saturation. For each $s$, the net $(f_i(s))_i$ is Cauchy,
since its difference net is a quasi-subnet of that saturation.
Completeness and separation supply a unique $h(s)\in X$ with
$f_i(s)\to h(s)$. Fix $i_0$, $i\ge i_0$ and $s\in S$.
As $j$ ranges over the tail $j\ge i_0$, the net
$f_i(s)-f_j(s)$ lies in $D_{i_0}$ and converges to $f_i(s)-h(s)$.
Consequently
\[
 E_{i_0}:=\{f_i(s)-h(s):i\ge i_0,\ s\in S\}
                   \subseteq K_X(D_{i_0}).
\]
Condition (CB) makes the saturation of $(K_X(D_{i_0}))$ null,
and the inclusion makes the saturation of $(E_{i_0})$ a quasi-subnet.
Thus $f_i\to_u h$. A second uniform limit gives the same pointwise
limits, so it equals $h$. The empty source produces a one-element
space and requires no separate argument.
\end{proof}

\begin{theorem}[Completeness of continuous mapping groups]
\label{nx:continuous-mapping-complete}
Let $X$ be separated and satisfy \textnormal{(CB)}. Then
\[
 X\text{ is complete and satisfies \textnormal{(P)}}
 \quad\Longleftrightarrow\quad
 C_u(S,X)\text{ is complete for every convergence space }S.
\]
If completeness of $X$ is given, the right side need only be tested
on Hausdorff topological spaces with at most one nonisolated point.
\end{theorem}
\begin{proof}
Suppose $X$ is complete and satisfies (P), and let $(f_i)$ be a
Cauchy net in $C_u(S,X)$. Its image in $X^S$ is uniformly Cauchy.
Lemma~\ref{nx:uniform-mapping-complete} supplies a uniform limit
$h\in X^S$, and Theorem~\ref{nx:uniform-limit} makes $h$ continuous.
Thus every Cauchy net in $C_u(S,X)$ converges there.

Conversely, a one-point source identifies $C_u(S,X)$ with $X$ and
gives completeness of $X$. Let $f_i$ be continuous and converge
uniformly to $h\in X^S$. The uniform error sets
\[
 E_{i_0}=\{f_i(s)-h(s):i\ge i_0,\ s\in S\}
\]
have null saturation, and the corresponding differences satisfy
\[
 \{f_i(s)-f_j(s):i,j\ge i_0,\ s\in S\}
                                  \subseteq E_{i_0}-E_{i_0}.
\]
The net is therefore Cauchy in $C_u(S,X)$. Its completeness supplies
a continuous uniform limit $g$, and uniqueness of pointwise limits
gives $g(s)=h(s)$ for every $s$. Theorem~\ref{nx:uniform-limit} gives
(P). If completeness of $X$ is already known, its restricted class
of source spaces suffices for this last implication.
\end{proof}

\subsection{A canonical intrinsic closure criterion}
\label{nx:canonical-section}

Let $X$ be a Hausdorff abelian convergence group. To determine when
ordered-group values are available, we return to its original nets
and ask which closed sets remain compatible with their translated
tails. Iterating this requirement determines a canonical closure
operator, whose preservation of null nets will characterize the
ordered-group representations considered in the following subsection.

\subsubsection{The transfinite construction}

Let $\mathscr C_0$ be the family of all originally closed subsets of
$X$. It contains $\varnothing,X$, is closed under arbitrary
intersections and finite unions, and is invariant under translations
and inversion. For finite unions, a convergent net eventually in
$C\cup D$ either has a tail in $D$ or has a quasi-subnet with every
value in $C$, obtained by saturating the nonempty intersections of its
tails with $C$. Closedness gives membership of its limit in the union.
The covariance assertions use the continuity and continuous inverses
of the group translations and inversion.

For a family $\mathscr C_\alpha$ with these properties, put
\[
 K_\alpha(A)=\bigcap\{C\in\mathscr C_\alpha:A\subseteq C\}.
\]
Define recursively
\[
 \begin{aligned}
 \mathscr C_{\alpha+1}
 &=\left\{C\in\mathscr C_\alpha:
     \bigcap_{i\in I}K_\alpha(C+T_i(u))=C
       \text{ for every original null net }(u_i)_{i\in I}\right\},\\
 \mathscr C_\eta&=\bigcap_{\alpha<\eta}\mathscr C_\alpha
       \qquad(\eta\text{ a nonzero limit ordinal}).
 \end{aligned}
\]
For a null net $u$, a fixed $i$ and $c\in C$, the translated tail
$(c+u_j)_{j\ge i}$ lies in $C+T_i(u)$ and converges to $c$. Hence
\begin{equation}\label{nx:contour-inclusion}
 C\subseteq K_0(C+T_i(u))\subseteq K_\alpha(C+T_i(u)).
\end{equation}
The second inclusion holds because
$\mathscr C_\alpha\subseteq\mathscr C_0$. In particular, the
construction does not require $0\in T_i(u)$.

\begin{lemma}\label{nx:canonical-families}
At every stage, $\mathscr C_\alpha$ contains $\varnothing,X$, is closed
under arbitrary intersections and finite unions, and is invariant under
translations and inversion. The families decrease and their closure
operators increase.
\end{lemma}
\begin{proof}
The assertions at stage zero were verified above. Suppose they hold at
stage $\alpha$. The successor condition is satisfied by the empty set
and the whole space. If $C_t\in\mathscr C_{\alpha+1}$ and
$C=\bigcap_tC_t$, then for every original null net $u$,
\[
 \bigcap_iK_\alpha(C+T_i(u))
 \subseteq\bigcap_t\bigcap_iK_\alpha(C_t+T_i(u))
 =\bigcap_tC_t=C.
\]
Equation \eqref{nx:contour-inclusion} supplies the reverse inclusion.
Thus arbitrary intersections are retained.

Because $\mathscr C_\alpha$ is closed under finite unions,
$K_\alpha(B\cup D)=K_\alpha(B)\cup K_\alpha(D)$. For
$C,D\in\mathscr C_{\alpha+1}$,
\[
 K_\alpha((C\cup D)+T_i(u))
 =K_\alpha(C+T_i(u))\cup K_\alpha(D+T_i(u)).
\]
If $x\notin C\cup D$, the successor condition supplies indices
$i_C,i_D$ for which $x$ belongs to neither corresponding closure.
For a common upper index $i\ge i_C,i_D$, monotonicity excludes $x$
from both terms in the last union. Intersecting over $i$ and applying
\eqref{nx:contour-inclusion} proves the successor condition for
$C\cup D$.

Translation covariance follows from
\[
 K_\alpha((C+x)+T_i(u))=x+K_\alpha(C+T_i(u)).
\]
For inversion, use $-u_i\to0$ and
\[
 K_\alpha((-C)+T_i(u))=-K_\alpha(C+T_i(-u)).
\]
At a limit ordinal, membership in every preceding family preserves
these properties. Finally, the recursion makes the families decrease,
and intersecting a smaller collection of supersets makes its closure
operator larger.
\end{proof}

There are only $2^{|X|}$ subsets of $X$. If every successor before
$(2^{|X|})^+$ were strict, assigning one set from
$\mathscr C_\alpha\setminus\mathscr C_{\alpha+1}$ at each such
stage would inject $(2^{|X|})^+$ into $\mathcal P(X)$. Consequently a
stage $\theta$ satisfies
$\mathscr C_\theta=\mathscr C_{\theta+1}$, after which the recursion
is stationary. Put
\[
 \mathscr C_\infty=\mathscr C_\theta,\qquad U_\lambda=K_\theta.
\]
Once stationarity has been reached, the operator no longer depends on
the chosen stage. It is extensive, monotone, idempotent, finite-union
preserving, translation and inversion covariant, and fixes the empty
set, while stationarity gives
\begin{equation}\label{nx:stationary-contour}
 U_\lambda(C)=C,\quad u_i\to0
 \quad\Longrightarrow\quad
 \bigcap_iU_\lambda(C+T_i(u))=C.
\end{equation}
All its closed sets are originally closed.

\subsubsection{Preservation of null nets by the canonical closure}

For a net $u$, the family $(U_\lambda(T_i(u)))_i$ is decreasing and
nonempty, so its saturation is defined; the remaining question is
whether this passage through the canonical closure preserves nullity.
We express that requirement by the intrinsic condition
\begin{equation}\label{nx:CR}
 u_i\to0
 \quad\Longrightarrow\quad
 \langle U_\lambda(T_i(u)):i\in I\rangle\to0.
 \tag{CR}
\end{equation}
When $\mathcal F=\operatorname{Tail}(u)$ is used as a derived
set-indexing notation, the family $U_\lambda(A)$, $A\in\mathcal F$,
has the same saturation up to tail equivalence. Thus
$U_\lambda(\mathcal F)$ in later formulas means the eventual sets of
this explicitly defined net.

\begin{lemma}\label{nx:canonical-zero}
If \textnormal{(CR)} holds, then $U_\lambda(\{0\})=\{0\}$, and
$U_\lambda$ fixes every finite subset of $X$.
\end{lemma}
\begin{proof}
Apply (CR) to the constant zero net. The net ranging over
$U_\lambda(\{0\})$ with the indiscrete directed preorder is null.
For any $x$ in this set, the constant net at $x$ is a quasi-subnet
and is therefore null. It also converges to $x$, so separation gives
$x=0$. Extensivity gives the reverse inclusion. Translation covariance
then fixes each singleton, and finite-union preservation fixes each
finite set.
\end{proof}

\subsubsection{Why every faithful ordered-group metric satisfies \textnormal{(CR)}}

Let $p:X\to G_+$ be definite, symmetric and subadditive, with $G$ a
partially ordered abelian group. Suppose that for every original net $u$,
\[
 u_i\to0\quad\Longleftrightarrow\quad
 \exists D\downarrow0\ \forall d\in D\ \exists i_d\
 \forall i\ge i_d:\ p(u_i)\le d.
\]
Write $B_r=\{x:p(x)\le r\}$. Each translate $c+B_r$ is originally
closed. Indeed, if $u_i\in c+B_r$ and $u_i\to x$, faithfulness gives
$D\downarrow0$ controlling $p(u_i-x)$. For every $d\in D$, a tail
satisfies
\[
 p(x-c)\le p(x-u_i)+p(u_i-c)\le d+r.
\]
Thus $p(x-c)-r$ is a lower bound of $D$ in the whole group, so it is
at most zero and $x\in c+B_r$.

\begin{lemma}\label{nx:balls-survive-canonical}
Every translated metric ball belongs to $\mathscr C_\alpha$ at every
stage of the canonical construction.
\end{lemma}
\begin{proof}
The preceding calculation proves the initial stage. Suppose the result
holds at stage $\alpha$, put $C=c+B_r$, and let $u_i\to0$ have
controls $D\downarrow0$. If
\[
 x\in\bigcap_iK_\alpha(C+T_i(u)),
\]
then for each $d\in D$ an index $i_d$ satisfies
$T_{i_d}(u)\subseteq B_d$, and therefore
\[
 x\in K_\alpha(C+T_{i_d}(u))
 \subseteq K_\alpha(c+B_{r+d})=c+B_{r+d}.
\]
The last equality uses the induction hypothesis. Hence
$p(x-c)-r\le d$ for every $d\in D$, giving $x\in c+B_r$.
Equation \eqref{nx:contour-inclusion} gives the other inclusion, so
$C\in\mathscr C_{\alpha+1}$. At a limit stage, membership in every
preceding family is membership in their intersection.
\end{proof}

Finite unions of translated balls survive as well, by
Lemma~\ref{nx:canonical-families}. Now let $u_i\to0$ and take its
controls $D\downarrow0$. For every $d\in D$, a tail satisfies
$T_{i_d}(u)\subseteq B_d$, and the lemma gives
\[
 U_\lambda(T_{i_d}(u))\subseteq U_\lambda(B_d)=B_d.
\]
These same controls apply to the saturation of
$(U_\lambda(T_i(u)))_i$. Faithfulness proves its original nullity,
and hence (CR).

\subsubsection{Equivalence with the existential closure formulation}

Call a closure operator $U$ admissible when it is extensive, monotone,
idempotent, fixes the empty set and the zero singleton, preserves finite
unions, is translation and inversion covariant, sends every original
null tail family $(T_i(u))$ to a null saturated family $(U(T_i(u)))$,
and satisfies
\begin{equation}\label{nx:admissible-contour}
 U(C)=C,\quad u_i\to0
 \quad\Longrightarrow\quad
 \bigcap_i U(C+T_i(u))=C.
\end{equation}
Every $U$-closed set is originally closed. To prove this, let
$u_i\in C$ and $u_i\to x$, with $U(C)=C$. Put $h_i=x-u_i\to0$.
For every $i$ and any $j\ge i$,
\[
 0=(u_j-x)+(x-u_j)\in(C-x)+T_i(h).
\]
The set $C-x$ is $U$-closed, so \eqref{nx:admissible-contour} gives
\[
 0\in\bigcap_iU((C-x)+T_i(h))=C-x.
\]
Thus $x\in C$.

Every $U$-closed set belongs to every $\mathscr C_\alpha$. The
preceding paragraph proves the initial stage. If the assertion holds
at stage $\alpha$, then $K_\alpha(B)\subseteq U(B)$ for every $B$.
For a $U$-closed set $C$ and a null net $u$,
\[
 \bigcap_iK_\alpha(C+T_i(u))
 \subseteq\bigcap_iU(C+T_i(u))=C.
\]
Equation \eqref{nx:contour-inclusion} gives equality. At a limit
stage the assertion is retained by intersection. Consequently
\begin{equation}\label{nx:canonical-minimality}
 U_\lambda(B)\subseteq U(B)\qquad(B\subseteq X).
\end{equation}
For a null net $u$, this inclusion makes the saturation of
$(U_\lambda(T_i(u)))$ a quasi-subnet of the null saturation of
$(U(T_i(u)))$. Thus an admissible $U$ implies (CR). Conversely,
(CR), stationarity \eqref{nx:stationary-contour}, and
Lemma~\ref{nx:canonical-zero} make $U_\lambda$ admissible.

\begin{theorem}[Canonical intrinsic characterization]\label{nx:canonical-po}
For a Hausdorff abelian convergence group $X$, the following are equivalent:
\begin{enumerate}
\item $X$ has a faithful invariant metric with values in a partially
ordered abelian group, with order convergence defined by downward-directed
regulators having whole-group infimum zero.
\item The canonical operator $U_\lambda$ satisfies \textnormal{(CR)}.
\item There exists an admissible closure operator $U$ as specified above.
\end{enumerate}
In the constructive direction the value group can be directed and free
abelian as an abstract group. No lattice order or replacement of the
original convergence is used.
\end{theorem}
\begin{proof}
Suppose (i), and let $p(x)=d(x,0)$. For a null net $u$, choose its
regulating family $D\downarrow0$ in the stated value group. For each
$d\in D$, a tail satisfies $T_{i_d}(u)\subseteq B_d$.
Lemma~\ref{nx:balls-survive-canonical} gives
\[
 U_\lambda(T_{i_d}(u))\subseteq U_\lambda(B_d)=B_d.
\]
Thus the image under $p$ of the saturation of
$(U_\lambda(T_i(u)))_i$ is controlled by $D$. Faithfulness makes
that saturation null, proving (ii).

Under (ii), Lemma~\ref{nx:canonical-zero} supplies
$U_\lambda(\{0\})=\{0\}$. The closure identities of
Lemma~\ref{nx:canonical-families}, stationarity
\eqref{nx:stationary-contour}, and (CR) give every requirement on
an admissible closure, so $U=U_\lambda$ proves (iii).
Conversely, an admissible $U$ satisfies
$U_\lambda(A)\subseteq U(A)$ by
\eqref{nx:canonical-minimality}. For a null $u$, the inclusions
\[
 U_\lambda(T_i(u))\subseteq U(T_i(u))\qquad(i\in I)
\]
make its canonical-closure saturation a quasi-subnet of the null
$U$-saturation, proving (ii).

Finally assume (iii). Its extensive, idempotent, finite-union
preserving and translation-covariant closure satisfies the four
hypotheses of Subsection~\ref{nx:information-construction}, including
\eqref{nx:information-contour}, which is
\eqref{nx:admissible-contour}. The construction and reflection proved
there, summarized in \ref{nx:information-realization}, give a directed
partially ordered group $G$, free abelian on its specified atoms, and
$p:X\to G_+$ with
\[
 \begin{gathered}
 p(x)=0\iff x=0,\qquad p(-x)=p(x),\\
 p(x+y)\le p(x)+p(y),\qquad
 u_i\to0\iff p(u_i)\to_{o_G}0.
 \end{gathered}
\]
Putting $d(x,y)=p(x-y)$ gives the invariant metric in (i).
Indeed, its triangle inequality is the subadditivity estimate at
$(x-y)+(y-z)$, and the last equivalence applied to the difference
net of any pair net proves canonical pair faithfulness.
\end{proof}

We have reached an ordered abelian group representation through
integer-valued degree functions controlling finite signed supports.
When $X$ is a real convergence vector space, a further requirement
would be needed to obtain positive homogeneity and an ordered real
vector-space structure on the values; neither is asserted by this
construction.

\subsection{An information-labelled ordered-group representation}
\label{nx:information-construction}

\subsubsection{Intrinsic hypotheses}
To realize the closure criterion by an ordered group, we shall retain
both the points being controlled and the information carried by their
null bases, so that coincident sets need not lose their distinct roles
in a regulator. Let $X$ be a Hausdorff abelian convergence group, with the
net conventions of Section~\ref{fd:net-foundations}. In calculations involving
eventual sets, a symbol $\mathcal F$ always has the form
$\mathcal F=\operatorname{Tail}(x)$ for an actual net $x=(x_i)_{i\in I}$,
and $\mathcal F\to a$ abbreviates $x_i\to a$. This abbreviation is
independent of the representative because tail-equivalent nets have the
same limits. Every set family introduced below is realized by its
saturation net before a convergence assertion is made.

For precision, if $\mathcal B$ is a nonempty family of nonempty sets,
directed downward under inclusion, its saturation has indices $(B,b)$,
with $B\in\mathcal B$, $b\in B$, ordered by reverse inclusion of $B$.
Its tail image above $(B,b)$ is $B$, and its eventual sets are exactly
\[
 \mathcal B^\uparrow
 =\{A:(\exists B\in\mathcal B)\ B\subseteq A\}.
\]
If a family is only centered, we first take its nonempty finite
intersections. Thus every generated family below is the eventual-set
family of an explicitly defined net. For $\mathcal F=\operatorname{Tail}(x)$
and $\mathcal G=\operatorname{Tail}(y)$ we use
\[
 T_*\mathcal F=\operatorname{Tail}((Tx_i)_i),\qquad
 \mathcal F+\mathcal G
   =\operatorname{Tail}((x_i+y_j)_{(i,j)}).
\]
For the intersection family $\mathcal F\cap\mathcal G$, take the net on
$I\times J\times\{0,1\}$ with values $x_i$ in the last coordinate $0$
and $y_j$ in the last coordinate $1$, ordered by $I\times J$. A family
inclusion $\mathcal G\subseteq\mathcal F$ consequently says that the net
representing $\mathcal F$ is a quasi-subnet of the net representing
$\mathcal G$. The notation $[a]$ denotes the eventual sets of the
constant net at $a$. These conventions are used only to shorten
calculations on nets and their tails. When a topology occurs below, its
zero-neighborhood family likewise denotes the eventual sets of the
saturation net of all zero neighborhoods; thus its use introduces no
additional kind of convergence object.

Suppose there is a closure operator $U$ on subsets of $X$ with the
following properties:
\begin{enumerate}
\item $U$ is extensive, monotone, idempotent, preserves finite unions, and
$U(\{0\})=\{0\}$.
\item $U(A+x)=U(A)+x$ and $U(-A)=-U(A)$.
\item For every net $x_i\to0$, the saturation of
$\{U(T_i(x)):i\in I\}$ converges to zero.
\item For every $U$-closed set $C$ and every net $x_i\to0$,
\begin{equation}
 \bigcap_{i\in I}U(C+T_i(x))=C.                 \label{nx:information-contour}
\end{equation}
Equivalently the intersection may be taken over
$A\in\operatorname{Tail}(x)$, since the actual tail images are cofinal
under inclusion among the eventual sets.
\end{enumerate}
With the convention $U(\varnothing)=\varnothing$, every finite set is
$U$-closed, since $U(\{x\})=x+U(\{0\})=\{x\}$ and $U(\bigcup_{k=1}^n\{x_k\})=\bigcup_{k=1}^n U(\{x_k\})$. Translation invariance also gives
\begin{equation}
 U(U(A)+B)=U(A+B),\qquad U(U(A)+U(B))=U(A+B).       \label{nx:hull-sum}
\end{equation}
Indeed, $U(A)+b\subseteq U(A+B)$ for every $b\in B$. Apply $U$ and use the reverse
inclusion supplied by extensivity. Iterating \eqref{nx:hull-sum} permits substitution inside any
finite Minkowski sum.

Let $\mathfrak N$ consist of the eventual-set families
$F=\operatorname{Tail}(n_F)$ of the canonical null saturation nets
whose defining bases consist of symmetric $U$-closed sets containing
zero. This is a set, since each such family is a subset of
$\mathcal P(X)$. Every null net is a quasi-subnet of one of these
saturations: for $x_i\to0$, first saturate
$T_i(x)\cup(-T_i(x))\cup\{0\}$, using finite mixing and inversion,
and then apply $U$ and its null-preservation property. Each new set
contains the original tail, so the quasi-subnet assertion holds.
For $F\in\mathfrak N$, let $\mathcal B_F$ consist of its symmetric
$U$-closed eventual sets containing zero. It is closed under finite
intersections and is cofinal among the eventual sets of $n_F$.
If $x\in\bigcap\mathcal B_F$, the constant net at $x$ is a
quasi-subnet of $n_F$, so separation gives $x=0$; conversely every
member contains zero. Thus $\bigcap\mathcal B_F=\{0\}$.

\subsubsection{Information retained by the atoms}
To encode the distinction between a point and the information that
controls it, we introduce two types of atoms.
\begin{itemize}
\item A point atom $Q_x$ for each $x\ne0$, with $Q_{-x}=Q_x$. Its underlying set
is $C(Q_x)=\{0,x,-x\}$.
\item An information atom $R_{F,S}$ for $F\in\mathfrak N$ and a proper subset
$S\subsetneq\mathcal B_F$. Its underlying set is
\[
 C(R_{F,S})=\bigcap_{A\in S}A,
\]
where the empty intersection is $X$.
\end{itemize}
By convention, $Q_0$ and $R_{F,\mathcal B_F}$ mean the empty bag, hence zero,
and do not denote atoms; an information atom, however, remains nonzero when
its underlying set is $\{0\}$. We thereby keep different information sets
distinct even when their intersections coincide.

Let $M$ be the free commutative monoid of finite bags of these atoms. Define
$a\preceq b$ by allocating disjoint blocks of the target bag $b$ to the source
atoms of $a$, with unused targets allowed, subject to these rules:
\begin{enumerate}
\item A source $R_{F,S}$ receives exactly one target $R_{F,T}$ of the same color,
with $T\subseteq S$.
\item A source $Q_x$ receives a block $J$ of targets such that
\[
 x\in U\left(\sum_{j\in J}C(b_j)\right).
\]
The empty block has underlying sum $\{0\}$ and cannot serve a nonzero point atom.
\end{enumerate}
Because every underlying set is symmetric and contains zero, replacing
one by a larger set, or adding unused summands, preserves point eligibility;
this is the compatibility needed for the allocation order below.

\paragraph{Lemma 1: order and cancellation.}
The relation $\preceq$ is a translation-invariant partial order on $M$, and
\[
 a+c\preceq b+c\quad\Longrightarrow\quad a\preceq b.       \label{nx:bag-cancellation}
\]

\emph{Proof.}
Reflexivity follows from singleton allocations: every information atom is
eligible for its identical target by $S\subseteq S$, while every point
atom is eligible for its identical target because
$x\in C(Q_x)\subseteq U(C(Q_x))$. To pass from these allocations to
transitivity, let an allocation from $a$ to $b$ assign a block
$J_i$ of $b$ to the $i$th source atom, and let an allocation from $b$ to
$c$ assign the block $K_j$ of $c$ to the $j$th atom of $b$.
Assign $\bigcup_{j\in J_i}K_j$ to the $i$th atom of $a$.
The sets $J_i$ are pairwise disjoint, as are all the $K_j$, so the
new blocks are pairwise disjoint. For an information source the two
blocks along its path are singletons of the same color, and the
inclusions $T\subseteq S$ and $V\subseteq T$ give $V\subseteq S$.
For every intermediate atom $b_j$, its allocation gives
\[
 C(b_j)\subseteq U\left(\sum_{k\in K_j}C(c_k)\right).
\]
For a point atom this uses symmetry and $0\in C(c_k)$ in addition
to the eligibility of its nonzero representative. For an information
atom it follows from the inclusion of its information sets.
Thus a point source $Q_x$ satisfies
\[
 x\in U\left(\sum_{j\in J_i}C(b_j)\right)
 \subseteq U\left(\sum_{j\in J_i}
       U\left(\sum_{k\in K_j}C(c_k)\right)\right)
 =U\left(\sum_{k\in\bigcup_{j\in J_i}K_j}C(c_k)\right),
\]
where the last equality is the finite iteration of
\eqref{nx:hull-sum}. This proves the required eligibility.

For \eqref{nx:bag-cancellation}, distinguish the occurrences of the catalyst bag $c$ on the two sides.
Draw an edge to a right-hand catalyst occurrence whenever a source consumes it,
and continue from the corresponding left-hand catalyst occurrence using its own
allocation. Each right-hand target has at most one predecessor. A component
reachable from a source belonging to $a$ cannot reach a directed cycle: the entry
vertex of such a cycle would have two predecessors. Thus these components are
finite acyclic graphs. Substitute along them until only right-hand occurrences
from $b$ remain. If the original source is an information atom, all vertices on
its path are information atoms of the same color. If it is a point atom, \eqref{nx:hull-sum}
justifies all substitutions. Distinct original sources use disjoint leaves.
Discard components not reached from $a$. This proves \eqref{nx:bag-cancellation}.

Every nonzero source consumes at least one target. Hence bag cardinality is
monotone. For each color, the number of its information atoms is also monotone.
If both $a\preceq b$ and $b\preceq a$, cardinalities and all color counts agree.
Every allocation is therefore a bijection. Information atoms match the same
color, while point atoms match point atoms. Point matching means equality up to
sign, because finite sets are $U$-closed. For each fixed color, write the information indices in the first bag
as $S_1,\ldots,S_m$ and those in the second as $T_1,\ldots,T_m$.
The two bijections give permutations $\sigma,\tau$ with
$T_{\sigma(i)}\subseteq S_i$ and $S_{\tau(j)}\subseteq T_j$.
If $\pi=\tau\circ\sigma$, then
\[
 S_{\pi(i)}\subseteq T_{\sigma(i)}\subseteq S_i.
\]
Iterating this inclusion around the finite cycle of $i$ under $\pi$
returns to $S_i$, so every inclusion in that cycle is equality.
In particular $T_{\sigma(i)}=S_i$ for every $i$.
All information multisets, and all point multisets, coincide.
Consequently $a=b$.
\hfill$\square$

Let $G$ be the Grothendieck group of $M$ with its induced order:
\[
 [a]-[b]\le[c]-[d]\quad\Longleftrightarrow\quad a+d\preceq c+b.
\]
Lemma 1 makes this a partially ordered abelian group, whose underlying group is
free abelian on the atoms. It is directed because it is generated by positive
atoms. Define $p(x)=[Q_x]$, including $p(0)=0$. Then
\[
 p(x)>0\ (x\ne0),\quad p(-x)=p(x),\quad p(x+y)\le p(x)+p(y).
\]
Indeed $x+y\in C(Q_x)+C(Q_y)$, and no nonzero point atom can be allocated to an
empty bag. Thus $d(x,y)=p(x-y)$ is an invariant $G$-valued metric.

\subsubsection{Representation of the prescribed null nets}
Having obtained the ordered group, we first verify that the prescribed
null nets remain null under the gauge. For a finite $J\subseteq\mathcal B_F$ put $e_{F,J}=[R_{F,J}]$, with the zero
convention when $J=\mathcal B_F$. These elements are downward directed as $J$
increases. We claim
\begin{equation}
 \inf_{J\subseteq\mathcal B_F,\ J\ \mathrm{finite}} e_{F,J}=0. \label{nx:designated-infimum}
\end{equation}
If $\mathcal B_F$ is finite, zero is a member of this positive family, so its infimum is zero. Assume henceforth that $\mathcal B_F$ is infinite. Let $v=[a]-[b]$ be a lower bound. Thus $a\preceq b+R_{F,J}$ for every finite $J$.
There are only finitely many allocation patterns for the fixed finite source
bag $a$, the fixed target bag $b$, and one extra target. Validity of a fixed
pattern is preserved when $J$ is replaced by a smaller finite set. Consequently
one pattern is valid for every $J$. Indeed, if the finite pattern set is $\{P_1,\ldots,P_n\}$ and each $P_k$ fails for some finite $J_k$, then $J=\bigcup_{k=1}^n J_k$ excludes every pattern, contradicting $a\preceq b+R_{F,J}$.

In this common pattern, the extra target cannot serve an information source
$R_{F,S}$, since this would
require $J\subseteq S$ for every finite $J$, and hence $S=\mathcal B_F$, which is
not an atom. If it is unused, simply discard it. Otherwise it serves a point
source $Q_x$, together with a fixed block of targets whose underlying closed sum
is $C=U(\sum C(b_j))$. For every finite $J$,
\[
 x\in U\left(C+\bigcap_{A\in J}A\right).
\]
Since these intersections form a null base of $F$, property
\eqref{nx:information-contour} gives $x\in C$, allowing the extra
target to be removed. We obtain an allocation proving
$a\preceq b$, or $v\le0$, as required for \eqref{nx:designated-infimum}.

If $x\in\bigcap_{A\in J}A$, then $p(x)\le e_{F,J}$. Thus each $F\in\mathfrak N$
is metrically null, and every original null tail family is metrically null by the
cofinal choice of $\mathfrak N$.

\subsubsection{Projections and finite support of a bounded regulator}
For the reverse implication, we must control regulators not prescribed
in advance; projections onto their information coordinates will reduce
this question to finitely many colors. We begin with an elementary free
ordered-group construction. If a poset $P$ has
least element $0$, let $H(P)$ be the free abelian group on $j(t)$, $t\ne0$,
put $j(0)=0$, and take the cone generated by
\[
 j(t)-j(s)\qquad(s\le t).
\]
This cone is proper. To prove this directly, embed the underlying group into the
free abelian group on all elements of $P$ by $j(t)=e_t-e_0$. If both $h$ and $-h$
are in the proposed cone, combine their finite expressions. This gives a finite
directed graph whose edges $s\to t$ satisfy $s<t$, with nonnegative integer
multiplicities, and with zero incoming-minus-outgoing balance at every vertex.
There is no directed cycle because every edge strictly increases the poset
order. A nonempty finite acyclic graph carrying a positive flow cannot have zero
balance everywhere: a source of an edge-bearing component has positive outgoing
flow and zero incoming flow. Thus there are no edges, and $h=0$. The cone
therefore defines a translation-invariant partial order. In particular every
$j(t)$ with $t\ne0$ is positive and nonzero. No assertion about preservation of
infima by arbitrary homomorphisms is needed below.

For a color $F$, use the information poset
\[
 P_F=\mathcal P(\mathcal B_F)^{\mathrm{op}},
\]
whose least element is $\mathcal B_F$. Define a homomorphism
$\pi_F:G\to H(P_F)$ by sending $R_{F,S}$ to $j(S)$ and every other atom to zero.
It is positive: source information atoms match comparable targets of the same
color, and all unused targets and targets used by points contribute positive
terms to the projection. Its section $j(S)\mapsto[R_{F,S}]$ is also positive.

Let $D\subseteq G_+$ be downward directed with $\inf_GD=0$. Fix $d_0\in D$ and put $D_0=\{d\in D:d\le d_0\}$. For every $d\in D$, directedness supplies $e\in D_0$ with $e\le d$, so $D_0$ is cofinal downward and $\inf_GD_0=0$. If a color $F$ does not
occur in the unique signed-word expression for $d_0$, positivity gives
\[
 0\le\pi_F(d)\le\pi_F(d_0)=0\quad(d\in D_0).
\]
Since the underlying group $H(P_F)$ is free on its proper information
subsets, all signed coefficients of color $F$ in $d$ vanish, leaving
only finitely many colors throughout $D_0$.

To refine this finite-color restriction, use the positive degree
$\ell:G\to\mathbb Z$ counting every atom once and, for each color,
the positive count $\sigma_F:G\to\mathbb Z$. Let
$\mathcal C$ be the finite set of colors occurring in $d_0$. Define
\[
 n=\min_{d\in D_0}\ell(d),\qquad
 m_F=\min_{d\in D_0}\sigma_F(d)\quad(F\in\mathcal C).
\]
All minima exist in $\mathbb Z_+$ and are attained. There are
$e_*\in D_0$ and $e_F\in D_0$ with
$\ell(e_*)=n$ and $\sigma_F(e_F)=m_F$. The set
\[
 D_1=\{d\in D_0:d\le e_*,\ d\le e_F\ (F\in\mathcal C)\}
\]
is nonempty, downward directed and cofinal downward in $D_0$: any finite
collection consisting of $e_*,e_F$ and an additional member of $D_0$
has a common lower member of $D_0$. Positivity of the integer maps gives
\[
 \ell(d)=n,\qquad \sigma_F(d)=m_F
 \quad(d\in D_1,\ F\in\mathcal C),\qquad \inf_GD_1=0.
\]
Put $\mathcal A=\{F\in\mathcal C:m_F>0\}$. A color with $m_F=0$
can still occur in the signed expression for $d$. These signed terms
remain in $R_d$ below. All subsequent regulator indices belong to $D_1$.
For $d=[a]-[b]\le d'=[a']-[b']$ on this tail, the comparison is
$a+b'\preceq a'+b$. Equality of $\ell(d)$ and $\ell(d')$ gives
$\ell(a)+\ell(b')=\ell(a')+\ell(b)$, and equality of each
$\sigma_F$ gives equality of the corresponding information counts.
Every source consumes a nonempty block, so an allocation between these
equal-cardinality bags consists of singleton blocks and leaves no target
unused. The information sources consume all information targets of their
colors. Hence every point source has a point target, and its eligibility
is equality of point atoms. For each point atom $Q$ this gives
\[
 \operatorname{mult}_Q(a)+\operatorname{mult}_Q(b')
 =\operatorname{mult}_Q(a')+\operatorname{mult}_Q(b),
\]
so its signed coefficient in $d$ equals its signed coefficient in $d'$.
For two incomparable members, use a common lower member of the directed tail.
Call the resulting fixed point word $w$. We have
\begin{equation}
 d=w+R_d,\qquad R_d=\sum_F s_F\pi_F(d)\ge0,             \label{nx:fixed-point-word}
\end{equation}
where $s_F$ is the positive section. In particular $w\le d$ for every $d$ on the
tail, and therefore $w\le0$.

\subsubsection{Eliminating the fixed point word}
It remains to eliminate the fixed word left by the preceding reduction.
Write $w=\sum_{i=1}^sQ_{u_i}-\sum_{j=1}^tQ_{v_j}$; to prove $w\ge0$,
suppose otherwise. For each subset $I\subseteq\{1,\ldots,s\}$ let
\[
 C_I=\sum_{i\in I}\{0,u_i,-u_i\}.
\]
These are finite $U$-closed sets. If the positive point atoms could be allocated
to the demands $Q_{v_j}$ with $v_j\in C_{I_j}$ for disjoint $I_j$, then $w\ge0$.
Thus every such allocation has a failed incidence.

If $\mathcal A=\varnothing$, matching every negative information atom
consumes all positive information atoms. The remaining point allocation
proves $w\ge0$, contrary to the supposition. Suppose
$\mathcal A\ne\varnothing$, and put
\[
 \mathcal V=\sum_{F\in\mathcal A}m_FF\to0,\qquad
 \mathcal J=\{(j,I):1\le j\le t,\ I\subseteq\{1,\ldots,s\},\ v_j\notin C_I\}.
\]
For every $(j,I)\in\mathcal J$, equation~\eqref{nx:information-contour}
gives $B_{j,I}\in\mathcal V$ with $v_j\notin U(C_I+B_{j,I})$.
The tail base of the independent sum net supplies $A^{j,I}_{F,r}\in\mathcal B_F$ such that
\[
 \sum_{F\in\mathcal A}\sum_{r=1}^{m_F}A^{j,I}_{F,r}
 \subseteq B_{j,I}.
\]
Set
\[
 A_F=\bigcap_{(j,I)\in\mathcal J}\bigcap_{r=1}^{m_F}A^{j,I}_{F,r}
 \in\mathcal B_F.
\]
An empty intersection is $X\in\mathcal B_F$. Since $\mathcal J$ is
finite, these are base members, and monotonicity gives
\begin{equation}\label{nx:finite-exclusions}
 v_j\notin C_I\quad\Longrightarrow\quad
 v_j\notin U\left(C_I+\sum_{F\in\mathcal A}m_FA_F\right).
\end{equation}

Fix $d$ on the tail and use an allocation witnessing $d=w+R_d\ge0$. The negative
information atoms of $R_d$ must consume same-color positive information atoms.
After these matches, precisely $m_F$ positive information atoms of color $F$
remain. These remaining atoms, together with the fixed $Q_{u_i}$, serve the
point demands $Q_{v_j}$.

If each remaining atom $R_{F,S}$ had $A_F\in S$, its underlying set would lie in
$A_F$. Its point allocations, together with \eqref{nx:finite-exclusions}, would imply all exact incidences
$v_j\in C_{I_j}$, contradicting $w\not\ge0$. Consequently at least one remaining
atom has $A_F\notin S$. Define the fixed information sensor
\[
 \delta_F=[R_{F,\,\mathcal B_F\setminus\{A_F\}}]>0.
\]
The information order gives $\delta_F\le[R_{F,S}]$. If
$\mathcal M_d$ is the finite family of matched information pairs and
$\mathcal L_d$ the bag of remaining positive information atoms, then
\[
 R_d=\sum_{(a,b)\in\mathcal M_d}([b]-[a])
           +\sum_{b\in\mathcal L_d}[b]
 \ge\sum_{b\in\mathcal L_d}[b],
\]
because every information match has $[a]\le[b]$.
Hence, for every $d$, at least one $F\in\mathcal A$ satisfies
\begin{equation}
 \delta_F\le R_d.                                    \label{nx:information-sensor}
\end{equation}

We now use downward directedness of $(R_d)_{d\in D_1}$ to make the
preceding obstruction independent of the chosen regulator member.
Suppose that for every $F\in\mathcal A$ there is $d_F\in D_1$ with
$\delta_F\not\le R_{d_F}$. A common lower member $e\le d_F$
for all $F$ satisfies $R_e\le R_{d_F}$, and therefore
$\delta_F\not\le R_e$ for every $F$, contrary to
\eqref{nx:information-sensor}. Hence some fixed $F\in\mathcal A$
satisfies $\delta_F\le R_d$ for every $d\in D_1$.
Then $w+\delta_F\le d$ for every $d$, while
\[
 \pi_F(w+\delta_F)=j(\mathcal B_F\setminus\{A_F\})>0.
\]
Since $\inf_GD_1=0$, the inequalities $w+\delta_F\le d$
would imply $w+\delta_F\le0$. Positivity of $\pi_F$ would then give
$0<\pi_F(w+\delta_F)\le0$. Thus $w\ge0$, and antisymmetry with
$w\le0$ gives $w=0$.

\subsubsection{Reflection of every null regulator}
With the fixed point word eliminated, we can compare the regulator balls
with the original null bases. For $d\ge0$ write $B_d=\{x:p(x)\le d\}$;
on the tail above, $d=R_d$ and the numbers $m_F$ are fixed. If $\mathcal A=\varnothing$, an allocation of $Q_x$ after all negative information atoms are matched has no remaining target. Hence $p(x)\le d$ implies $x=0$, and the saturation of the balls is the constant zero net. Suppose $\mathcal A\ne\varnothing$. For any family $A_F\in\mathcal B_F$, $F\in\mathcal A$, set
\[
 C=U\left(\sum_{F\in\mathcal A}m_FA_F\right).
\]
We claim $B_d\subseteq C$ for some $d$ on the tail. Otherwise $B_d\setminus C\ne\varnothing$ for every $d$. Fix $x_d\in B_d\setminus C$. In an allocation witnessing $Q_{x_d}\le
R_d$, first match all negative information atoms. Exactly $m_F$ same-color
positive atoms remain. If all their information sets contained $A_F$, the point
allocation would give $x_d\in C$. Thus a remaining atom misses $A_F$, and the decomposition of $R_d$
into matched differences and remaining atoms gives $\delta_F\le d$
for some $F\in\mathcal A$. If no fixed $F$ satisfied this for all
$d\in D_1$, there would be $d_F\in D_1$ with $\delta_F\not\le d_F$.
A common lower member $e\le d_F$ would satisfy $\delta_F\not\le e$
for every $F$, contradicting the preceding assertion applied to $e$.
Consequently a fixed nonzero positive $\delta_F$ would be a lower
bound of $D_1$, contradicting $\inf_GD_1=0$. Therefore some $d$
satisfies $B_d\subseteq C$.

Continuity of finite addition and the null-preservation property of $U$
show that the saturation of the sets
$U(\sum_{F\in\mathcal A}m_FA_F)$ is originally null. Every such set
contains some ball $B_d$, so the saturation of $(B_d)_{d\in D_1}$ is a
quasi-subnet of this null net and is itself null. Now let $(z_i)$ be
any net satisfying $p(z_i)\to_o0$. A regulating family $D\downarrow0$
makes $z_i$ eventually belong to every $B_d$, hence $(z_i)$ is a
quasi-subnet of the null ball saturation. Therefore $z_i\to0$.
Together with the forward implication, this proves fidelity for every
net in the original convergence.

\paragraph{Conclusion of the construction.}\label{nx:information-realization}
Under the intrinsic closure assumptions stated above, $X$ admits an invariant
faithful metric with values in a directed partially ordered abelian group.
For every net $(z_i)$ in $X$, the equivalence is
\[
 z_i\to0\iff p(z_i)\to_o0.
\]
Consequently $d(x,y)=p(x-y)$ represents the original canonical pair
convergence by Proposition~\ref{lt:canonical-group}.
Proposition~\ref{lt:three-translations} gives point, Cauchy, and
parameter-uniform faithfulness.

\subsubsection{Scalar covariance and canonical closed disks}
\label{nx:canonical-disk-section}

To relate this construction to the vector-space setting, suppose in addition
that $X$ is a real convergence vector space, so addition and scalar
multiplication are continuous for the original convergence. A disk means
a nonempty convex balanced subset, and $\operatorname{aco}A$ denotes the smallest
disk containing $A$, with $\operatorname{aco}\varnothing=\{0\}$. Local convexity
will mean that every original null net is a quasi-subnet of a null
saturation whose defining sets are disks. No absorbing-set requirement is part of this definition.

\begin{proposition}\label{nx:canonical-disks}
For every nonzero real scalar $t$ and every subset $A\subseteq X$,
\[
 U_\lambda(tA)=tU_\lambda(A).
\]
If the canonical condition \textnormal{(CR)} holds, $U_\lambda$ sends disks to disks, and
\[
 H(A)=U_\lambda(\operatorname{aco}A)
\]
is the smallest $U_\lambda$-closed disk containing $A$. In particular, $H(A)$ is
closed for the original convergence. The operator $H$ is extensive, monotone,
and idempotent.

If the original convergence is locally convex in the sense specified above,
then $H$ preserves original null nets:
\[
 \mathcal F\to0\quad\Longrightarrow\quad
 \bigl\langle H(A):A\in\mathcal F\bigr\rangle\to0.
\]
Here angle brackets denote the canonical saturation net of the displayed sets.
\end{proposition}

\begin{proof}
For $t\ne0$, the linear map $x\mapsto tx$ and its inverse are continuous.
Consequently $\mathscr C_0$ is invariant under multiplication by $t$, and
$K_0(tA)=tK_0(A)$. Suppose $\mathscr C_\alpha$ is invariant under all nonzero
scalar multiplications. Its associated closure then has the same covariance.
If $C\in\mathscr C_{\alpha+1}$ and $\mathcal F\to0$, the represented tail family
$t^{-1}\mathcal F$ is originally null. Hence
\[
 \begin{split}
 \bigcap_{A\in\mathcal F}K_\alpha(tC+A)
 &=t\bigcap_{A\in\mathcal F}K_\alpha(C+t^{-1}A)\\
 &=tC.
 \end{split}
\]
Thus $t\mathscr C_{\alpha+1}\subseteq\mathscr C_{\alpha+1}$ for every
$t\ne0$. The inclusion for $t^{-1}$ gives the reverse inclusion.
At a limit ordinal, invariance is preserved by
intersection of families. Transfinite induction, followed by stabilization,
proves the asserted scalar covariance of $U_\lambda$.

We also use the additive compatibility of a translation-covariant closure:
\begin{equation}
 U_\lambda(U_\lambda(A)+U_\lambda(B))=U_\lambda(A+B).
                                                        \label{nx:canonical-sum}
\end{equation}
For completeness, for each $b\in B$ one has
$U_\lambda(A)+b=U_\lambda(A+b)\subseteq U_\lambda(A+B)$.
Taking the union over $b$, applying $U_\lambda$, and using extensivity gives
$U_\lambda(U_\lambda(A)+B)=U_\lambda(A+B)$. Commutativity and the first equality give
\[
 U_\lambda(U_\lambda(A)+U_\lambda(B))
 =U_\lambda(A+U_\lambda(B))=U_\lambda(A+B),
\]
which proves \eqref{nx:canonical-sum}.

Let $D$ be a disk. If $0<|t|\leq1$, balancedness and scalar covariance give
\[
 tU_\lambda(D)=U_\lambda(tD)\subseteq U_\lambda(D).
\]
For $t=0$, $tU_\lambda(D)=\{0\}\subseteq D\subseteq U_\lambda(D)$. Thus
$U_\lambda(D)$ is balanced. For $0<t<1$, equation \eqref{nx:canonical-sum}, covariance, and convexity
of $D$ give
\[
 \begin{split}
 tU_\lambda(D)+(1-t)U_\lambda(D)
 &=U_\lambda(tD)+U_\lambda((1-t)D)\\
 &\subseteq U_\lambda(tD+(1-t)D)
 \subseteq U_\lambda(D).
 \end{split}
\]
For $t=0$ or $t=1$ the convex combination is one of its arguments. Thus convexity holds for every $t\in[0,1]$. Therefore
$H(A)=U_\lambda(\operatorname{aco}A)$ is a $U_\lambda$-closed disk. If $C$ is any
$U_\lambda$-closed disk containing $A$, then $\operatorname{aco}A\subseteq C$,
and hence $H(A)\subseteq U_\lambda(C)=C$. This proves minimality.
The inclusions $A\subseteq\operatorname{aco}A\subseteq H(A)$ prove extensivity. If $A\subseteq B$, then $\operatorname{aco}A\subseteq\operatorname{aco}B$ and $H(A)\subseteq H(B)$. Since $H(A)$ is a $U_\lambda$-closed disk, $\operatorname{aco}H(A)=H(A)=U_\lambda(H(A))$, and $H(H(A))=H(A)$. It is originally closed
because every member of $\mathscr C_\infty$ is originally closed.

Finally, let $x_i\to0$. Local convexity supplies a null saturation
of disks $(D_\eta)$ such that, for each $\eta$, one tail $T_i(x)$ lies
in $D_\eta$. Absolute convexification then gives
$\operatorname{aco}T_i(x)\subseteq D_\eta$. Consequently the saturation
of $(\operatorname{aco}T_i(x))_i$ is a quasi-subnet of the disk
saturation and converges to zero. Condition \textnormal{(CR)} applied
to that net gives nullity of the saturation of
$U_\lambda(\operatorname{aco}T_i(x))=H(T_i(x))$. The actual tails are
cofinal among all eventual sets, so this is also the displayed
assertion with $\mathcal F=\operatorname{Tail}(x)$.
\end{proof}

\begin{remark}[No absorption conclusion]
Although Proposition~\ref{nx:canonical-disks} provides closed disks, these
need not be absorbing in $X$. Indeed, if $Y\subsetneq X$ is a proper
$U_\lambda$-closed linear subspace and $A\subseteq Y$, then
\[
 H(A)\subseteq U_\lambda(Y)=Y.
\]
For any $x\notin Y$, no positive scalar multiple of $H(A)$ contains $x$, so
$H(A)$ cannot be absorbing, and this observation also applies to a proper
linear ideal whenever that ideal is $U_\lambda$-closed. Thus preservation
of the original null nets by closed disks carries no additional claim
that these particular disks, or all disks in a given null base, absorb
the ambient space.
\end{remark}

\subsection{Packaging faithful gauge families}
\label{nx:packaging-section}

Having constructed individual representations, we turn to families of
gauges, where we build simultaneous control by leaving all but finitely
many coordinates unrestricted. A lexicographic enlargement supplies bounds
for those inactive coordinates and thereby allows the family to be
expressed by one gauge.

\begin{theorem}[Lexicographic packaging of faithful families]\label{nx:packaging}
Let $X$ be a Hausdorff convergence abelian group and let $(p_i:X\to(G_i)_+)_{i\in I}$ be a faithful family of symmetric subadditive gauges into partially ordered abelian groups, with order infima taken in the whole groups. Then there is a single faithful symmetric subadditive gauge $q:X\to P_+$ into a partially ordered abelian group. If all $G_i$ are lattice ordered, $P$ can be lattice ordered. The construction represents all nets and therefore preserves the original group-Cauchy and uniform-error convergences.
\end{theorem}
\begin{proof}
Discard identically zero gauges. If none remain, separation and faithfulness force $X=\{0\}$, represented by the zero gauge into the zero group. For each remaining $i$, let
\[
 H_i=\mathbb Z\mathbin{\rm lex}G_i,
\]
with componentwise group addition and positivity defined by $(n,g)>0$ when $n>0$, or by $(0,g)\ge0$ when $g\ge0$. Let $j_i(g)=(0,g)$ and $e_i=(1,0)$. Set $q_i=j_i\circ p_i$. Then $0\le q_i(x)\le e_i$ for every $x$.

The embedding $j_i$ preserves and reflects positive order convergence. For preservation, if $D\downarrow0$ in $G_i$, the family $j_i(D)$ is downward directed. A lower bound $(n,g)$ either has negative first coordinate and is below zero, or has first coordinate zero with $g$ a lower bound of $D$, again giving $(n,g)\le0$. A positive first coordinate cannot be a lower bound. Hence $j_i(D)\downarrow0$.

For reflection, suppose a downward-directed family $R\subseteq(H_i)_+$ with whole-group infimum zero regulates an image net under $q_i$. Because $p_i$ is not identically zero, there exists $v\in(G_i)_+\setminus\{0\}$. If all first coordinates of $R$ were at least one, $(0,v)$ would be a nonzero positive lower bound of $R$. Therefore $R$ contains an element with first coordinate zero. The cofinal downward tail below that element has first coordinate zero throughout. Projecting this tail to $G_i$ gives a downward-directed family with infimum zero: any lower bound $g$ would give the lower bound $(0,g)$ of the original tail and hence satisfy $g\le0$. Its balls give the required regulator for $p_i$.

Let $P=\prod_{i\in I}H_i$ with coordinate order and
$q(x)=(q_i(x))_{i\in I}$. Coordinate addition gives
\[
 q(-x)=q(x),\qquad q(x+y)\le q(x)+q(y),\qquad q(0)=0.
\]
If $q(x)=0$, then $p_i(x)=0$ for every $i$. The constant net at $x$ is order-null under every $p_i$, hence
converges to zero by faithfulness. Its convergence to $x$ and separation give $x=0$. If $q_*(\mathcal F)$ is order-null, each coordinate image is order-null. Indeed, the coordinate projection of a positive downward-directed zero regulator is downward directed, and any lower bound in that coordinate can be lifted to a lower bound of the original regulator by placing zeros in all other coordinates. Thus all $p_i(\mathcal F)$ are order-null and $\mathcal F\to0$.

Conversely, suppose $\mathcal F\to0$. For each $i$, faithfulness supplies a regulator $D_i\downarrow0$ for $p_i(\mathcal F)$. For finite $F\subseteq I$ and choices $d_i\in D_i$ for $i\in F$, define $r(F,d)\in P_+$ by
\begin{equation}
 r(F,d)_i=
 \begin{cases}(0,d_i),&i\in F,\\e_i,&i\notin F.
 \end{cases}
 \label{nx:finite-support-controls}
\end{equation}
For two controls $r(F,d),r(F',d')$, put $F''=F\cup F'$.
For $i\in F\cap F'$, directedness gives $d_i''\in D_i$ with
$d_i''\le d_i,d_i'$. For $i\in F\setminus F'$, put
$d_i''=d_i$, and for $i\in F'\setminus F$, put $d_i''=d_i'$.
Then
\[
 r(F'',d'')\le r(F,d),r(F',d').
\]
If $v\in P$ is a lower bound of all these controls, then, for every
$i$ and every $d\in D_i$, the control with active set $\{i\}$ gives
$v_i\le(0,d)$. Since $j_i(D_i)\downarrow0$, one has $v_i\le0$
for every $i$, and therefore $v\le0$. Thus their whole-product
infimum is zero. Finally,
\[
 A(F,d)=\bigcap_{i\in F}\{x:p_i(x)\le d_i\}\in\mathcal F,
 \qquad q(x)\le r(F,d)\quad(x\in A(F,d)).
\]
Only finitely many eventual sets are intersected, and all inactive
coordinates are bounded by $e_i$. These controls prove
$q_*\mathcal F\to_o0$.

For lattice-ordered $G_i$, the join in $H_i$ is
\[
 (n,g)\vee(m,h)=
 \begin{cases}
 (n,g),&n>m,\\
 (m,h),&m>n,\\
 (n,g\vee h),&n=m.
 \end{cases}
\]
The meet has the corresponding smaller-first-coordinate formula and
$(n,g\wedge h)$ for equal first coordinates. The product has these
operations coordinatewise, so $P$ is lattice ordered. For an arbitrary pair net $(x_\alpha,y_\alpha)$, the established
faithfulness applied to $x_\alpha-y_\alpha$ gives
\[
 x_\alpha-y_\alpha\to0
 \quad\Longleftrightarrow\quad
 q(x_\alpha-y_\alpha)\to_o0.
\]
Taking the pair net $(x_\alpha,x_\beta)_{(\alpha,\beta)}$ proves the
Cauchy assertion. For maps $f_\alpha,f:S\to X$, apply the same
equivalence to the saturation with tail sets
$\{f_\alpha(s)-f(s):\alpha\ge\alpha_0,\ s\in S\}$.
This proves parameter-uniform faithfulness for arbitrary index sets.
\end{proof}

This passage from a family to one gauge may introduce non-Archimedean
values and need not retain Dedekind completeness from the coordinates.
If the original gauges already have fixed global order bounds, however,
those bounds serve the inactive coordinates directly, giving the
following construction without a lexicographic enlargement.

\begin{corollary}[Bounded lattice-valued products and initial constructions]\label{nx:bounded-products}
Suppose each convergence group $X_i$ has a faithful gauge $p_i:X_i\to(G_i)_+$ into a lattice-ordered group and a fixed bound $u_i\in(G_i)_+$ with $p_i(x_i)\le u_i$ for all $x_i$. Then the product convergence on $\prod_iX_i$ has the faithful gauge
\begin{equation}
 p((x_i)_i)=(p_i(x_i))_i\in\prod_iG_i.
 \label{nx:bounded-product-gauge}
\end{equation}
More generally, a jointly injective family of homomorphisms $T_i:X\to X_i$ equips $X$ with an initial convergence faithfully represented by the family $p_i\circ T_i$, and hence by one gauge under Theorem~\ref{nx:packaging}. If the $X_i$ are complete Hausdorff groups, this initial convergence is complete exactly when the diagonal image $T(X)$ is closed under convergence limits in $\prod_iX_i$.
\end{corollary}
\begin{proof}
For \eqref{nx:bounded-product-gauge}, use the finite-coordinate regulator construction \eqref{nx:finite-support-controls} with inactive value $u_i$ and active values $d_i\wedge u_i$, where $d_i$ belongs to a regulator for the coordinate image net. These are valid bounds because the gauge also has the global bound $u_i$. For two finite active sets, their union and common lower controls
in each $D_i$ give a common lower product control. If $v$ is a lower
bound of all controls, activating $\{i\}$ gives
$v_i\le d\wedge u_i\le d$ for every $d\in D_i$, whence $v_i\le0$.
Therefore the whole-product infimum is zero. Conversely, coordinate projections preserve positive order convergence by the lower-bound argument in Theorem~\ref{nx:packaging}. In a merely partially ordered value group the required finite meets may be unavailable, and this direct-product assertion is not made. The general case is covered by the lexicographic theorem.

Initial convergence means exactly that every image net under $T_i$ converges in $X_i$, so the composed family is faithful under the equivalence $\mathcal F\to x\iff (T_i)_*\mathcal F\to T_i(x)$ for every $i$. 

Let $(z_\alpha)$ be a Cauchy net in $\prod_iX_i$. Each coordinate
net $\pi_i(z_\alpha)$ is Cauchy, so completeness of $X_i$ supplies
$x_i\in X_i$ with $\pi_i(z_\alpha)\to x_i$. Product convergence gives
$z_\alpha\to(x_i)_i$. If the net is eventually in a closed subgroup,
its limit is in that subgroup, proving completeness of the induced
convergence there. Conversely, a net in the complete diagonal image
$T(X)$ that converges in the ambient product is Cauchy in the induced
convergence. It has an internal limit in $T(X)$, and separation of
the ambient product identifies this with its ambient limit. Therefore
the diagonal image is closed, as asserted.
\end{proof}

\subsection{Cancellative labelled comparison bases}
\label{cb:section}

Once a family can be expressed by one gauge, the lattice-group
realization problem can also be formulated through the algebra of
labels attached to its comparison sets. Keeping the labels distinct
even when they determine the same subset leads to the characterization
below, whose scope differs from that of the intrinsic closure criterion
\eqref{nx:CR} for partially ordered group representations.

\begin{definition}[Cancellative comparison basis]\label{cb:definition}
Let $X$ be a separated abelian convergence group. A cancellative comparison
basis consists of a commutative cancellative monoid $(C,+,0)$ and symmetric
subsets $(B_c)_{c\in C}$ of $X$ with the following properties.
The natural order
\begin{equation}\label{cb:natural-order}
 a\le b\quad\Longleftrightarrow\quad
 (\exists h\in C)\ a+h=b
\end{equation}
is antisymmetric and admits binary meets. The labelled sets satisfy
\begin{equation}\label{cb:ball-laws}
 B_0=\{0\},\qquad 0\in B_c=-B_c,\qquad
 B_a+B_b\subseteq B_{a+b},\qquad B_e=X\text{ for some }e\in C.
\end{equation}
For each $x\in X$, the nonempty set
\[
 R(x)=\{c\in C:x\in B_c\}
\]
has a least member $m(x)$. Write $D\downarrow_C0$ when
$D\subseteq C$ is nonempty, downward directed, and has infimum zero in
$C$. The convergence requirements are
\begin{align}
 D\downarrow_C0&\quad\Longrightarrow\quad
 \mathcal B_D:=\operatorname{Tail}(\langle B_d:d\in D\rangle)\longrightarrow0,
 \label{cb:null-bases}\\
 \mathcal F\longrightarrow0&\quad\Longrightarrow\quad
 (\exists D\downarrow_C0)\ \mathcal B_D\subseteq\mathcal F.
 \label{cb:cofinal-bases}
\end{align}
In these statements $\mathcal F=\operatorname{Tail}(x)$ for a net $x$, and every displayed base is interpreted through its saturation net.
\end{definition}

These algebraic requirements already make the ball family isotone:
$a\le b$ gives $b=a+h$ and
\begin{equation}\label{cb:ball-monotonicity}
 B_a=B_a+\{0\}\subseteq B_a+B_h\subseteq B_b.
\end{equation}
For $d_1,d_2\in D\downarrow_C0$, a common lower member $d\in D$
therefore satisfies $0\in B_d\subseteq B_{d_1}\cap B_{d_2}$.
Thus the family in \eqref{cb:null-bases} is indeed a nonempty downward-directed tail base.

\begin{lemma}[The lattice group of differences]\label{cb:lattice-completion}
Let $C$ have the monoid and order properties in
Definition~\ref{cb:definition}. Its group of differences $L$ is an
abelian lattice-ordered group. The canonical map $i:C\to L$ identifies
$C$, with its given order, with $L_+$. For every nonempty downward-directed
$D\subseteq C$,
\begin{equation}\label{cb:zero-reflection}
 \inf_C D=0\quad\Longleftrightarrow\quad \inf_L i(D)=0.
\end{equation}
The right-hand infimum tests every lower bound in $L$.
\end{lemma}
\begin{proof}
Define
\[
 L=(C\times C)/\sim,\qquad
 (a,b)\sim(c,d)\quad\Longleftrightarrow\quad a+d=c+b.
\]
Reflexivity and symmetry are immediate equalities. For transitivity,
\[
 a+d=c+b,\quad c+f=e+d
 \quad\Longrightarrow\quad
 a+f+d=e+b+d\quad\Longrightarrow\quad a+f=e+b.
\]
The final implication is cancellation in $C$. Addition and inversion are
\[
 [a,b]+[c,d]=[a+c,b+d],\qquad -[a,b]=[b,a].
\]
If $a+b'=a'+b$ and $c+d'=c'+d$, then
\[
 (a+c)+(b'+d')=(a'+c')+(b+d),
\]
so addition is independent of representatives. Interchanging the two
coordinates proves that inversion is independent of representatives.
These operations give an abelian group with zero $[0,0]$, and
$i(c)=[c,0]$ is an injective monoid homomorphism.

The cone $i(C)$ is additive and has trivial intersection with its negative.
Indeed, $i(a)=-i(b)$ gives $a+b=0$. The natural order gives
$0\le a\le0$, hence $a=0$, and cancellation gives $b=0$.
Declare this cone to be $L_+$. In representatives its order is
\begin{equation}\label{cb:group-order}
 [a,b]\le[c,d]\quad\Longleftrightarrow\quad a+d\le c+b.
\end{equation}
To verify this identity, the difference $[c+b,d+a]$ belongs to $i(C)$
exactly when $c+b=d+a+h$ for some $h\in C$, which is the
right-hand condition. The cone construction gives a translation-invariant
partial order, and \eqref{cb:group-order} proves that $i$ is an order
embedding onto $L_+$.

The order of $C$ admits order cancellation:
\[
 a+t\le b+t
 \quad\Longrightarrow\quad
 a+t+h=b+t\text{ for some }h\in C
 \quad\Longrightarrow\quad a+h=b\quad\Longrightarrow\quad a\le b.
\]
If $z=(a+t)\wedge(b+t)$, then $t\le z$, so $z=t+v$ for some $v\in C$.
Order cancellation gives $v\le a,b$ and hence $v\le a\wedge b$.
Conversely $(a\wedge b)+t\le a+t,b+t$. Therefore
\begin{equation}\label{cb:translated-meet}
 (a+t)\wedge(b+t)=(a\wedge b)+t.
\end{equation}

For $x=[a,b]$ and $y=[c,d]$, put
\[
 m=(a+d)\wedge(c+b),\qquad z=[m,b+d].
\]
The common-denominator expressions
$x=[a+d,b+d]$ and $y=[c+b,b+d]$ give $z\le x,y$.
If $h=[r,s]\le x,y$, then
\[
 r+b+d\le(a+d)+s,\qquad r+b+d\le(c+b)+s.
\]
Equation~\eqref{cb:translated-meet} gives
\[
 r+b+d\le\bigl((a+d)+s\bigr)\wedge\bigl((c+b)+s\bigr)
 =m+s,
\]
which is exactly $h\le z$ by \eqref{cb:group-order}. Consequently
\begin{equation}\label{cb:group-meet}
 [a,b]\wedge[c,d]=[(a+d)\wedge(c+b),b+d].
\end{equation}
Its greatest-lower-bound property makes this formula independent of the
representatives. The operation $x\vee y=-((-x)\wedge(-y))$ supplies
joins, so $L$ is a lattice-ordered group.

Suppose $D\downarrow_C0$, and let $h\in L$ be an arbitrary lower bound
of $i(D)$. Each $i(d)$ dominates both $h$ and $0$, and therefore
\[
 h\vee0\le i(d)\qquad(d\in D).
\]
Since $h\vee0\in L_+=i(C)$, it equals $i(c)$ for some $c\in C$.
Order reflection gives $c\le d$ for every $d\in D$. Thus $c=0$ and
$h\le h\vee0=0$. Zero is itself a lower bound, proving
$\inf_L i(D)=0$. Conversely, if that whole-group infimum is zero,
every lower bound $c\in C$ of $D$ satisfies $i(c)\le0$ and $i(c)\ge0$.
Hence $c=0$, proving \eqref{cb:zero-reflection}. Downward directedness
is preserved and reflected by the order embedding throughout.
\end{proof}

\begin{theorem}[Lattice representations and comparison bases]
\label{cb:characterization}
For a separated abelian convergence group $X$, the following conditions
are equivalent.
\begin{enumerate}[(i)]
\item The original null nets have a faithful family of symmetric
subadditive lattice-group-valued gauges.
\item There are an abelian lattice-ordered group $L$, an element $u\in L_+$,
and a definite symmetric subadditive gauge $p:X\to L_+$ such that
$p(X)\subseteq[0,u]$ and, for every net,
\[
 \mathcal F\to0\quad\Longleftrightarrow\quad p_*\mathcal F\to_o0.
\]
\item The original convergence admits a cancellative comparison basis.
\end{enumerate}
Each representation also preserves the original canonical pair,
group-Cauchy and parameter-uniform convergences.
\end{theorem}
\begin{proof}
The implication (i)$\Rightarrow$(ii) is the lattice case of
Theorem~\ref{nx:packaging}. Its construction is bounded: if the original
value groups are $L_j$, it has
\[
 L=\prod_j(\mathbb Z\mathbin{\rm lex}L_j),\qquad
 p(x)=((0,p_j(x)))_j,\qquad
 0\le p(x)\le u:=((1,0))_j.
\]
The theorem proves arbitrary-net faithfulness. Separation gives definiteness.
If every coordinate gauge vanishes, faithfulness and separation give
$X=\{0\}$, and $L=\{0\}$ gives the same assertion.

Assume (ii), and put $C=L_+$ and $B_c=\{x:p(x)\le c\}$.
The natural order of $C$ is its inherited order because
\[
 a\le b\quad\Longleftrightarrow\quad b-a\in L_+.
\]
It is cancellative, antisymmetric, and has binary meets. Definiteness,
symmetry, subadditivity and boundedness give \eqref{cb:ball-laws}, with
$e=u$. Moreover
\[
 R(x)=\{c\in L_+:p(x)\le c\},\qquad \min R(x)=p(x).
\]
For $D\downarrow_C0$, Lemma~\ref{cb:lattice-completion}, or its
positive-lower-bound calculation inside $L$, gives $\inf_LD=0$.
The balls $B_d$ therefore generate an eventual-set family of a net with order-null $p$-image,
and faithfulness makes that represented net originally null. If $\mathcal F\to0$,
an order regulator $D\subseteq L_+$ for $p_*\mathcal F$ has
$B_d\in\mathcal F$ for every $d\in D$. Its whole-group zero infimum
implies its zero infimum in $C$, proving \eqref{cb:cofinal-bases}.

Assume (iii). Construct $L$ and $i:C\to L_+$ by
Lemma~\ref{cb:lattice-completion}, and put $p(x)=i(m(x))$.
The least-label property and \eqref{cb:ball-monotonicity} give the exact
identity
\begin{equation}\label{cb:exact-balls}
 x\in B_c\quad\Longleftrightarrow\quad m(x)\le c
 \quad\Longleftrightarrow\quad p(x)\le i(c).
\end{equation}
The first equivalence uses $x\in B_{m(x)}$: its reverse implication
is the inclusion $B_{m(x)}\subseteq B_c$.
Symmetry gives $R(-x)=R(x)$, so $p(-x)=p(x)$. Since
$x+y\in B_{m(x)}+B_{m(y)}\subseteq B_{m(x)+m(y)}$, one has
\[
 p(x+y)\le p(x)+p(y),\qquad
 p(x)\le i(e),\qquad
 p(x)=0\ \Longleftrightarrow\ x\in B_0\ \Longleftrightarrow\ x=0.
\]
Thus $p$ has every property in (ii) except for the faithfulness still
to be verified.

If $\mathcal F\to0$, condition~\eqref{cb:cofinal-bases} supplies
$D\downarrow_C0$ with $B_d\in\mathcal F$. Equations
\eqref{cb:zero-reflection} and \eqref{cb:exact-balls} make $i(D)$
a whole-group zero regulator for $p_*\mathcal F$.
Conversely, let $V\subseteq L_+$ be a downward zero regulator for
$p_*\mathcal F$. Since $L_+=i(C)$, its preimage
$D=i^{-1}(V)$ satisfies $D\downarrow_C0$ by
\eqref{cb:zero-reflection}. Equation~\eqref{cb:exact-balls} gives
$B_d\in\mathcal F$ for every $d\in D$. Condition
\eqref{cb:null-bases} and refinement yield $\mathcal F\to0$.
This proves (iii)$\Rightarrow$(ii). A single gauge is a family, giving
(ii)$\Rightarrow$(i).

Finally, for $d(x,y)=p(x-y)$ and the original difference map
$\delta(x,y)=x-y$, every pair net $\Phi$ satisfies
\[
 d_*\Phi=p_*(\delta_*\Phi),\qquad
 \Phi\in\mathcal U_\lambda
 \ \Longleftrightarrow\ \delta_*\Phi\to0
 \ \Longleftrightarrow\ d_*\Phi\to_o0.
\]
The three translations of Proposition~\ref{lt:three-translations}
give the final assertions.
\end{proof}

\begin{proposition}[The exact disk-hull condition on the labelled basis]
\label{cb:local-convexity}
Let $X$ be a real convergence vector space with a cancellative comparison
basis. Its original convergence is locally convex if and only if
\begin{equation}\label{cb:disk-cofinality}
 \forall D\downarrow_C0\ \exists E\downarrow_C0\quad
 \forall e\in E\ \exists d\in D:
 \quad\operatorname{aco}(B_d)\subseteq B_e.
\end{equation}
In particular, disk-valued comparison balls satisfy this condition with
$E=D$. Local convexity does not require every ball of a specified
representation to be a disk.
\end{proposition}
\begin{proof}
For $D\downarrow_C0$, monotonicity of the labelled balls and of
absolute convexification shows that
\[
 \mathcal A_D=\operatorname{Tail}(\langle\operatorname{aco}(B_d):d\in D\rangle)
\]
is the eventual-set family of its saturation. Its defining sets are the disk hulls of
$\mathcal B_D$: if $A\in\mathcal B_D$, some $B_d\subseteq A$ gives
$\operatorname{aco}(B_d)\subseteq\operatorname{aco}(A)$, and each
$B_d$ is itself a member of $\mathcal B_D$.

Assume local convexity. By \eqref{cb:null-bases},
$\mathcal B_D\to0$. Equation~\eqref{nx:aco-null} gives
$\mathcal A_D\to0$. Applying \eqref{cb:cofinal-bases} to this represented net
supplies $E\downarrow_C0$ with $\mathcal B_E\subseteq\mathcal A_D$.
For every $e\in E$, membership $B_e\in\mathcal A_D$ is equivalent
to the existence of $d\in D$ with
$\operatorname{aco}(B_d)\subseteq B_e$. This is
\eqref{cb:disk-cofinality}.

Conversely, assume \eqref{cb:disk-cofinality}, and let
$\mathcal F\to0$. Condition~\eqref{cb:cofinal-bases} gives
$D\downarrow_C0$ with $\mathcal B_D\subseteq\mathcal F$.
The disk-hull saturation net
$\mathcal A(\mathcal F)=\operatorname{Tail}(\langle\operatorname{aco}(A):A\in\mathcal F\rangle)$
then contains $\mathcal A_D$. For the family $E$ supplied by
\eqref{cb:disk-cofinality}, the displayed inclusions give
\[
 \mathcal B_E\subseteq\mathcal A_D\subseteq\mathcal A(\mathcal F),
 \qquad\mathcal B_E\to0.
\]
Refinement proves $\mathcal A(\mathcal F)\to0$. Equation
\eqref{nx:aco-null} yields local convexity. If each $B_d$ is a disk,
$\operatorname{aco}(B_d)=B_d$, which verifies the indicated special
case directly.
\end{proof}

We have thus expressed the lattice-group question in the quantifiers
of Definition~\ref{cb:definition} and \eqref{cb:disk-cofinality}, which
retain a specified algebra of labels and its cofinal null-ball saturation
nets. Their scope does not extend to arbitrary partially ordered group
values, nor do they follow from the universal incidence monoid alone.
In passing back to families, the lexicographic construction is also
significant: an unmodified product of unbounded coordinate gauges may
have coordinatewise null nets with no order-bounded tail, and hence may
fail to preserve order convergence.

\subsection{Intrinsic boundedness, countability and weak distributivity}
\label{dx:dictionary-section}

We now return from the algebra of a representation to the properties
visible in its original convergence, beginning with the order convergence
of a vector lattice. Throughout this subsection, a real convergence
vector space has continuous addition and scalar multiplication, with
the usual convergence on the scalar field. A convergence is \emph{first countable} when, for every net
$x_i\to x$, there are decreasing nonempty sets $A_n$ such that their
saturation converges to $x$ and, for each $n$, the original net is
eventually in $A_n$. The countability belongs to the controlling sets,
not to the original directed set.
A subset $B$ is \emph{bounded} when
\[
 \big\langle[-\varepsilon,\varepsilon]B:\varepsilon>0\big\rangle\longrightarrow0.
\]
Both definitions refer only to the original convergence, with the empty
set taken to be bounded.

\begin{definition}[Bounded tail amalgamation]\label{dx:tail-amalgamation}
The original convergence has bounded tail amalgamation if every double
sequence $(x_{n,m})_{n,m\ge1}$ with bounded total range and null rows has a
null tail family generated by
\begin{equation}\label{dx:tail-tail}
 U_\varphi=\bigcup_{n\ge1}\{x_{n,m}:m\ge\varphi(n)\},
 \qquad \varphi\in\mathbb N^{\mathbb N}.
\end{equation}
Here $\mathbb N=\{1,2,\ldots\}$, and
$U_{\max(\varphi,\psi)}\subseteq U_\varphi\cap U_\psi$.
\end{definition}

\begin{theorem}[Intrinsic dictionary for vector-lattice order convergence]
\label{dx:order-dictionary}
Let $E$ be an Archimedean real vector lattice with order convergence.
\begin{enumerate}[(i)]
\item A subset of $E$ is convergence-bounded if and only if it is order bounded.
\item The convergence is complete if and only if $E$ is Dedekind complete.
\item If $E$ is Dedekind complete, its convergence is first countable if and
only if $E$ has the countable sup property.
\item If $E$ is Dedekind complete, bounded tail amalgamation is equivalent
to weak $\sigma$-distributivity.
\end{enumerate}
Consequently, super Dedekind completeness and weak $\sigma$-distributivity
are equivalent to completeness, first countability and bounded tail
amalgamation of order convergence.
\end{theorem}
\begin{proof}
For (i), $B\subseteq[-u,u]$ gives
$[-1/n,1/n]B\subseteq[-u/n,u/n]$, and Archimedeanity gives $u/n\downarrow0$.
Conversely, a regulator of the scalar-dilation saturation net supplies
$\varepsilon>0$ and $v\in E_+$ such that
$[-\varepsilon,\varepsilon]B\subseteq[-v,v]$. Taking the scalar
$\varepsilon$ gives $B\subseteq[-v/\varepsilon,v/\varepsilon]$.

Suppose $E$ is Dedekind complete and $(x_\alpha)$ is order-Cauchy. One
Cauchy bound and one point in a corresponding tail give an order-bounded
tail. Restrict to that tail and put
\[
 \ell_\alpha=\inf_{\beta\ge\alpha}x_\beta,
 \qquad u_\alpha=\sup_{\beta\ge\alpha}x_\beta.
\]
If $D\downarrow0$ regulates the Cauchy differences, then for every $d\in D$
there is $\alpha_d$ with
\[
 |x_\beta-x_\gamma|\le d\quad(\beta,\gamma\ge\alpha_d),
 \qquad 0\le u_{\alpha_d}-\ell_{\alpha_d}\le d.
\]
The second inequality results from the successive supremum in $\beta$ and
infimum in $\gamma$. Thus
\[
 0\le\inf_\alpha u_\alpha-\sup_\alpha\ell_\alpha\le d\quad(d\in D),
 \qquad x:=\sup_\alpha\ell_\alpha=\inf_\alpha u_\alpha.
\]
For $\beta\ge\alpha_d$, both $x$ and $x_\beta$ belong to
$[\ell_{\alpha_d},u_{\alpha_d}]$, and $|x_\beta-x|\le d$. Hence
$x_\alpha\to_o x$. This proves completeness for the original arbitrary-net convergence.

Conversely, suppose order convergence is complete. Let $(x_\alpha)$ be an
increasing net with nonempty upper-bound set $V$. Its associated set
\[
 D=\{v-x_\alpha:v\in V,\ \alpha\}
\]
is positive and downward directed: a common later $\gamma$ and
$v\wedge w\in V$ give
$(v\wedge w)-x_\gamma\le v-x_\alpha,w-x_\beta$.
If $h$ is a lower bound of $D$, then $h^+=h\vee0$ is also a lower bound.
Writing $w=h^+$, the inequalities $w\le v-x_\alpha$ imply $v-w\in V$.
Induction gives $v-nw\in V$, and hence
\[
 nw\le v-x_{\alpha_0}\qquad(n\ge1).
\]
Archimedeanity gives $w=0$, so $h\le0$ and $D\downarrow0$. For
$\beta,\gamma\ge\alpha$,
$|x_\beta-x_\gamma|\le v-x_\alpha$. Completeness gives an order limit $x$.
Let $R\downarrow0$ regulate $x_\beta-x$. For fixed $\alpha$ and
$r\in R$, an index $\beta$ beyond $\alpha$ and the threshold for $r$
satisfies
\[
 x_\alpha-x\le x_\beta-x\le |x_\beta-x|\le r.
\]
Thus $x_\alpha-x\le\inf R=0$. For fixed $v\in V$, the same
regulator gives
$x-v=(x-x_\beta)+(x_\beta-v)\le |x-x_\beta|\le r$
on a suitable tail, for every $r\in R$. Hence $x\le v$.
Therefore $x$ is an upper bound of the original net and is below every
upper bound, proving $x=\sup_\alpha x_\alpha$.
Apply this to the increasing net of finite suprema of a nonempty
bounded-above subset of $E$. This proves (ii).

For (iii), the countable sup property implies that every $D\downarrow0$
has a countable subfamily with infimum zero. Indeed, fix $d_0\in D$ and
use the countable sup property on
$\{d_0-d:d\in D,\ d\le d_0\}$, whose supremum is $d_0$.
Successive common lower bounds in $D$ turn that subfamily into
$d_n\in D$, $d_n\downarrow0$. If $D$ regulates a null tail family, the intervals
$[-d_n,d_n]$ generate a coarser countably based null tail family.

Conversely, apply first countability to
$\mathcal F_D=\operatorname{Tail}(\langle[-d,d]:d\in D\rangle)$. Let $(A_n)$ be a decreasing
base of a coarser null tail family. For every $n$, the nonempty set
\[
 \{d\in D:[-d,d]\subseteq A_n\}
\]
supplies an element $d_n$. If $R\downarrow0$ regulates the coarser eventual-set family,
then for each $r\in R$ there is $n$ with
$[-d_n,d_n]\subseteq A_n\subseteq[-r,r]$, whence $d_n\le r$.
Every lower bound of $(d_n)$ is therefore a lower bound of $R$ and is
nonpositive. Thus $\inf_n d_n=0$. For a nonempty bounded-above set $S$
with supremum $s$, apply this conclusion to
\[
 D_S=\{s-\bigvee F:\varnothing\ne F\subseteq S,\ F\text{ finite}\}.
\]
Countably many finite subsets $F_n$ suffice, and
$\sup(\bigcup_n F_n)=s$. This proves the countable sup property.

For (iv), let $|x_{n,m}|\le u$ and let every row be order-null. Dedekind
completeness defines
\[
 a_{n,m}=\sup_{j\ge m}|x_{n,j}|,\qquad
 b_\varphi=\sup_n a_{n,\varphi(n)}.
\]
An arbitrary regulator for row $n$ eventually dominates $a_{n,m}$.
Hence $a_{n,m}\downarrow_m0$. Weak $\sigma$-distributivity gives
$b_\varphi\downarrow0$, with selectors ordered pointwise, while
$U_\varphi\subseteq[-b_\varphi,b_\varphi]$. The represented tail family in
\eqref{dx:tail-tail} is null. Conversely, let
$0\le a_{n,m}\le u$ with $a_{n,m}\downarrow_m0$. Bounded tail amalgamation
applied to this array supplies a regulator $R\downarrow0$. For every
$r\in R$ some $\varphi$ satisfies
\[
 U_\varphi\subseteq[-r,r],\qquad
 \sup_n a_{n,\varphi(n)}\le r.
\]
The downward-directed path suprema consequently have infimum zero,
which is weak $\sigma$-distributivity.
\end{proof}

\begin{proposition}[Bounded diagonal selection]\label{dx:bounded-diagonal}
First countability and bounded tail amalgamation imply that every bounded
double sequence with null rows has a selector $\psi\in\mathbb N^{\mathbb N}$
such that $x_{n,\psi(n)}\to0$. In a Dedekind complete vector lattice this
bounded diagonal conclusion, by itself, implies weak $\sigma$-distributivity.
\end{proposition}
\begin{proof}
Let $(H_k)$ be a decreasing base of a null tail family coarser than the represented tail family
in \eqref{dx:tail-tail}. The inclusions
$U_{\varphi_k}\subseteq H_k$ have solutions $\varphi_k$. Replacement by
$\max_{j\le k}\varphi_j$ makes $\varphi_k$ increasing pointwise. Put
$\psi(n)=\varphi_n(n)$. For $n\ge k$,
$x_{n,\psi(n)}\in U_{\varphi_k}\subseteq H_k$, proving convergence.

For the converse assertion, let $0\le a_{n,m}\le u$ and
$a_{n,m}\downarrow_m0$. The stated diagonal conclusion supplies $\psi$
with $a_{n,\psi(n)}\to_o0$. Thus
\[
 q_k=\sup_{n\ge k}a_{n,\psi(n)}\downarrow0,
 \qquad w=\inf_\varphi\sup_n a_{n,\varphi(n)}
\]
exist, and for arbitrary $m_1,\ldots,m_{k-1}$,
\[
 0\le w\le \left(\bigvee_{n<k}a_{n,m_n}\right)\vee q_k.
\]
For $c_i\downarrow0$ and $q\ge0$, one has
$0\le(c_i\vee q)-q=(c_i-q)^+\le c_i$, so
$\inf_i(c_i\vee q)=q$. Successive application to the finitely many
$m_n$ gives $w\le q_k$. Taking the infimum in $k$ gives $w=0$.
\end{proof}

\subsection{Transfer through a faithful gauge and the role of its range}
\label{dx:gauge-transfer-section}

To carry the preceding dictionary through a faithful gauge, we must
also account for its range, since coordinates of the value lattice may
contribute nothing to the original convergence. Tracking this distinction
through the transfer arguments will identify range assumptions under
which the converses hold.

In this subsection $X$ is an Abelian convergence group and
$p:X\to E_+$ is symmetric and subadditive, with $p(0)=0$, where $E$ is
a Dedekind complete real vector lattice. Faithfulness means
\begin{equation}\label{dx:gauge-faithful}
 \mathcal F\to0\quad\Longleftrightarrow\quad p_*\mathcal F\to_o0
 \qquad(\mathcal F\text{ a tail family of a net on }X).
\end{equation}
Write $B_u=\{x:p(x)\le u\}$.

\begin{proposition}[Sequential regulators]\label{dx:sequential-regulators}
The original convergence is first countable if and only if every null net
$x$, with $\mathcal F=\operatorname{Tail}(x)$, admits $u_n\in E_+$ with $u_n\downarrow0$ and
$B_{u_n}\in\mathcal F$ for all $n$.
\end{proposition}
\begin{proof}
By their defining bounds, the displayed balls generate a countably based
null tail family coarser than $\mathcal F$, using
\eqref{dx:gauge-faithful}. For the converse, let $(A_n)$
be a decreasing base of a coarser null tail family. One image-regulator bound
dominates $p(A_{n_0})$. Discarding $A_1,\ldots,A_{n_0-1}$ and reindexing
therefore makes
$u_n=\sup\{p(x):x\in A_n\}$ well defined. The sequence decreases.
Every bound $d$ of a regulator of the coarser image net dominates
$p(A_n)$ for some $n$, and consequently $u_n\le d$. Thus
$u_n\downarrow0$. The inclusions $A_n\subseteq B_{u_n}$ give
$B_{u_n}\in\mathcal F$.
\end{proof}

\begin{proposition}[Sequential extraction from uniform convergence]
\label{dx:uniform-extraction}
Let $X$ be a first-countable Abelian convergence group, let $S\ne\varnothing$,
and let $(f_\alpha)_{\alpha\in D}$ converge uniformly to $f:S\to X$ in
the original group uniform convergence. There are increasing indices
$\beta_n\in D$ for which $f_{\beta_n}\to f$ uniformly. These indices
need not be cofinal in $D$.
\end{proposition}
\begin{proof}
The sets
\[
 E_\alpha=\{f_\gamma(s)-f(s):\gamma\ge\alpha,\ s\in S\}
\]
generate a null saturation net. A coarser null tail family has a decreasing
countable base $(H_n)$. For each $n$ there is $\alpha_n$ with
$E_{\alpha_n}\subseteq H_n$. Recursion in the directed set supplies
$\beta_1\ge\alpha_1$ and
$\beta_n\ge\beta_{n-1},\alpha_n$ for $n>1$. If $m\ge n$, then
\[
 \{f_{\beta_m}(s)-f(s):s\in S\}
 \subseteq E_{\alpha_n}\subseteq H_n.
\]
Consequently the sequence's uniform error eventual-set family refines
the null tail family generated by the $H_n$, without requiring the
selected indices to be cofinal.
\end{proof}

\begin{proposition}[Transfer and converse hypotheses]\label{dx:transfer-fullness}
The countable sup property of $E$ implies first countability of $X$.
If $p(X)=E_+$, the converse holds. If $X$ is a real convergence vector
space, weak $\sigma$-distributivity of $E$ implies bounded tail
amalgamation of $X$. Its converse holds under the two additional conditions
\begin{equation}\label{dx:full-homogeneous}
 p(X)=E_+,\qquad p(tx)=|t|p(x)\quad(t\in\mathbb R,\ x\in X).
\end{equation}
\end{proposition}
\begin{proof}
The regulator extraction in Theorem~\ref{dx:order-dictionary} and
Proposition~\ref{dx:sequential-regulators} prove the first implication.
For its converse, let $D\downarrow0$ in $E_+$. The saturation
$\langle B_d:d\in D\rangle$ is defined since the balls are nonempty
and downward directed, and it is null by faithfulness. A countable coarser null base $(A_n)$
supplies $d_n\in D$ with $B_{d_n}\subseteq A_n$. For each bound $r$ of
a regulator for this coarser eventual-set family, some $n$ satisfies
$A_n\subseteq B_r$. Surjectivity supplies $x\in X$ with $p(x)=d_n$,
and hence
\[
 x\in B_{d_n}\subseteq A_n\subseteq B_r,\qquad d_n\le r.
\]
Thus $\inf_n d_n=0$. The final argument for (iii) in
Theorem~\ref{dx:order-dictionary} proves the countable sup property.

Suppose now that $B\subseteq X$ is intrinsically bounded. A regulator
bound for its scalar-dilation image supplies $u\in E_+$ and
$\varepsilon>0$ such that
$p(tb)\le u$ whenever $b\in B$ and $|t|\le\varepsilon$.
An integer $k$ with $1/k\le\varepsilon$ gives
\[
 p(b)=p(\underbrace{b/k+\cdots+b/k}_{k\text{ terms}})
 \le k p(b/k)\le ku\qquad(b\in B).
\]
Consequently a bounded double sequence in $X$ has order-bounded image.
For a bounded array $(x_{n,m})$, put
\[
 a_{n,m}=\sup_{j\ge m}p(x_{n,j}),\qquad
 b_\varphi=\sup_n a_{n,\varphi(n)}.
\]
The preceding bound $ku$ makes these suprema elements of $E$.
Faithfulness makes the $n$th image row order-null. If $D_n\downarrow0$
regulates that row, then for each $d\in D_n$ some $m$ satisfies
$a_{n,m}\le d$, so $a_{n,m}\downarrow_m0$.
Weak $\sigma$-distributivity gives $b_\varphi\downarrow0$.
For the tail-union sets in \eqref{dx:tail-tail},
\[
 p(x)\le b_\varphi\qquad(x\in U_\varphi).
\]
The saturation of $(U_\varphi)$ therefore has order-null $p$-image,
and faithfulness makes it null in $X$, as required.

Under \eqref{dx:full-homogeneous}, let $0\le a_{n,m}\le u$ with
$a_{n,m}\downarrow_m0$. For each pair $(n,m)$ take an element in the
nonempty set $p^{-1}\{a_{n,m}\}$ and denote it by $x_{n,m}$.
The rows are null, and
\[
 p(t x_{n,m})=|t|a_{n,m}\le |t|u.
\]
Since $u/k\downarrow0$, the entire lifted array is intrinsically bounded.
Tail amalgamation in $X$ and faithfulness yield tail amalgamation for
$(a_{n,m})$. The regulator argument in part (iv) of
Theorem~\ref{dx:order-dictionary} proves weak $\sigma$-distributivity of $E$.
\end{proof}

\begin{remark}[Unused value coordinates]\label{dx:unused-coordinates}
For any vector lattice $H$, the gauge
$\widetilde p(x)=(p(x),0)$ into $E\times H$ represents exactly the same
nets. A regulator in $E$ lifts to $(d,0)$, while a positive product
regulator projects to an $E$-regulator, since an $E$ lower bound lifts
with zero second coordinate. With $H=\mathbb R^I$ for uncountable $I$,
the ambient countable sup property fails: the finite-support indicators
$\mathbf1_F$, $F\subseteq I$ finite, have supremum $\mathbf1_I$, whereas
every countable subfamily has support in a countable subset of $I$.
Thus an ambient-controller property cannot be identified with an intrinsic
property of $X$ without a range condition, an existential representation
statement, or a condition restricted to the effective values.
\end{remark}

\begin{definition}[Countable sequential determination]\label{dx:csd-definition}
A convergence is countably sequentially determined if, for every
decreasing sequence of nonempty sets $(F_n)$ and every point $x$, the
saturation $\langle F_n:n\ge1\rangle$ converges to $x$ whenever every
sequence eventually belonging to each $F_n$ converges to $x$.
\end{definition}

\begin{example}[First countability does not determine countably based nets by sequences]
\label{dx:c0-example}
The vector lattice $c_0$ is Dedekind complete: if $0\le S\le u\in(c_0)_+$,
its coordinatewise supremum belongs to $[0,u]\subseteq c_0$. Translation
by one member treats an arbitrary nonempty bounded-above set. Its
countable coordinates give the countable sup property. Explicitly, for
$s=\sup S$, elements approximating $s(k)$ within $1/m$, for
$(k,m)\in\mathbb N^2$, give a countable subset with supremum $s$.
For a globally bounded regulator $a_{n,m}\downarrow_m0$, fixed coordinate
$k$ and $\varepsilon>0$, a selector $\varphi$ satisfying
$a_{n,\varphi(n)}(k)\le\varepsilon$ for every $n$ gives
$(\inf_\varphi\sup_n a_{n,\varphi(n)})(k)=0$. All relevant suprema and
infima are coordinatewise, so $c_0$ is weakly $\sigma$-distributive.
Theorem~\ref{dx:order-dictionary} gives completeness, first countability
and bounded tail amalgamation of its order convergence.

Nevertheless, the tail family of a net generated by
$A_n=\{x\in c_0:\|x\|_\infty\le1/n\}$ is not order-null. An order bound
of $A_n$ would dominate $n^{-1}e_k$ for every $k$, and would have every
coordinate at least $1/n$, contrary to membership in $c_0$.
Every finer sequence is norm-null. For a norm-null sequence $(x_m)$, put
\[
 u_n(k)=\sup_{m\ge n}|x_m(k)|.
\]
Given $\varepsilon>0$, an index $M\ge n$ has
$\|x_m\|_\infty\le\varepsilon$ for $m\ge M$. An index $K$ has
$|x_m(k)|\le\varepsilon$ for $n\le m<M$ and $k\ge K$.
Hence $u_n(k)\le\varepsilon$ for $k\ge K$, proving $u_n\in c_0$.
Coordinatewise $u_n\downarrow0$, and $|x_n|\le u_n$.
Every finer sequence is therefore order-null, although the original
countably based tail family is not. Countable sequential determination fails.
\end{example}

\begin{theorem}[Globally bounded gauges and countable sequential determination]
\label{dx:bounded-csd}
Assume \eqref{dx:gauge-faithful}, and suppose $E$ is super Dedekind
complete and weakly $\sigma$-distributive. If $p(X)\subseteq[0,e]$ for
some $e\in E_+$, then $X$ is first countable and countably sequentially
determined.
\end{theorem}
\begin{proof}
First countability follows from Proposition~\ref{dx:transfer-fullness}.
Translate the proposed limit to zero. Let $(F_n)$ be a decreasing base
of a tail family of a net $\mathcal F$, and assume every sequence with finer
eventual-set family is null. Define
\[
 r_n=\sup_{x\in F_n}p(x)\le e,\qquad r=\inf_n r_n.
\]
Super Dedekind completeness supplies $x_{n,j}\in F_n$ with
$\sup_j p(x_{n,j})=r_n$. Repetitions give an infinite listing when the
required family is finite. Put
\[
 v_{n,k}=r_n-\bigvee_{j\le k}p(x_{n,j}).
\]
Then $0\le v_{n,k}\le e$ and $v_{n,k}\downarrow_k0$. For each
$\varphi\in\mathbb N^{\mathbb N}$ concatenate the finite nonempty blocks
$(x_{n,1},\ldots,x_{n,\varphi(n)})$ in increasing $n$. All blocks with
$n\ge N$ belong to $F_N$, so the concatenated eventual-set family refines
$\mathcal F$. It is null. For any regulator bound $d$ of its image,
all terms of every sufficiently late block satisfy $p(x_{n,j})\le d$.
Therefore
\[
 b_n=\bigvee_{j\le\varphi(n)}p(x_{n,j})\longrightarrow_o0.
\]
Set $w_\varphi=\sup_m v_{m,\varphi(m)}$. For every $n$,
\[
 r\le r_n=b_n+v_{n,\varphi(n)}\le b_n+w_\varphi.
\]
Thus $r-w_\varphi$ is bounded above by every member of a zero regulator
for $(b_n)$, and $r\le w_\varphi$. Weak $\sigma$-distributivity gives
$\inf_\varphi w_\varphi=0$, whence $r=0$. The bounds
$p(F_n)\le r_n\downarrow0$ prove $\mathcal F\to0$ by faithfulness.
\end{proof}

\begin{proposition}[An intrinsic obstruction to homogeneous realizations]
\label{dx:homogeneous-local-boundedness}
If a real convergence vector space is faithfully represented by a
homogeneous lattice norm $q:X\to E_+$ into an Archimedean vector lattice,
every null net has an intrinsically bounded tail. In a Hausdorff
convergence vector space no bounded set contains a nonzero vector subspace.
\end{proposition}
\begin{proof}
A null net is eventually in $B_u=\{x:q(x)\le u\}$ for one regulator bound
$u\in E_+$. Homogeneity gives
$q(tx)\le\varepsilon u$ for $|t|\le\varepsilon$, $x\in B_u$.
The relation $u/n\downarrow0$ and faithfulness prove boundedness of $B_u$.
For the second assertion, if a subspace $W$ lies in $B$, then
$W\subseteq[-\varepsilon,\varepsilon]B$ for every $\varepsilon>0$:
for $w\in W$ and $0<|t|\le\varepsilon$, $w=t(w/t)$ with $w/t\in W$.
If $B$ were bounded, the constant net at each $w\in W$ would be a
quasi-subnet of its null scalar-dilation saturation and converge to zero. Hausdorffness gives
$w=0$, hence $W=\{0\}$.
\end{proof}

\subsection{Strong first countability does not imply weak distributivity}
\label{sx:Cohen-countability}

The preceding distinctions leave a further question about countability:
whether strengthening it can supply the distributive property used in
the regulator arguments. To examine this, we use the net formulation of
strong first countability from
\citet[Section~2.2]{BBW2025countability}: for each $x$ there is a
countable family $\mathcal H_x$ such that every net $y_i\to x$ admits
a nonempty downward-directed subfamily $\mathcal B\subseteq\mathcal H_x$
with
\[
 \langle\mathcal B\rangle\to x,
 \qquad
 (\forall B\in\mathcal B)(\exists i_B)(\forall i\ge i_B)
 \quad y_i\in B.
\]
Here the countable family is fixed at $x$, while its subfamily may
depend on the given net; in the vector-lattice example below, convergence
means ordinary order convergence with downward-directed regulators.

\begin{lemma}\label{sx:strong-countability-csd}
Every strongly first-countable convergence is countably sequentially
determined in the sense of Definition~\ref{dx:csd-definition}.
\end{lemma}
\begin{proof}
Let $(F_n)_{n\ge1}$ be a decreasing base of a tail family of a net
$\mathcal F$, and suppose every sequence with eventual-set family finer than
$\mathcal F$ converges to $x$. Enumerate the fixed family at $x$ as
$(H_i)_{i\ge1}$, with repetitions if finite, and put
$I=\{i:H_i\notin\mathcal F\}$.
For every $n$, let $a_n\in F_n$. For $i\in I$, $i\le n$, let
$b_{n,i}\in F_n\setminus H_i$. These sets are nonempty, since
$F_n\subseteq H_i$ would imply $H_i\in\mathcal F$.
Concatenate the finite blocks
\[
 (a_n,(b_{n,i})_{i\in I\cap\{1,\ldots,n\}}),\qquad n=1,2,\ldots.
\]
All blocks with $n\ge N$ lie in $F_N$, so the resulting sequence
converges to $x$. Its eventual-set family $\mathcal T$ has
$H_i\notin\mathcal T$ for every $i\in I$, because $b_{n,i}$
lies outside $H_i$ for all $n\ge i$. Strong first countability
supplies a tail base
$\mathcal B\subseteq\mathcal T\cap\mathcal H_x$ with
$\langle\mathcal B\rangle\to x$.
The preceding exclusion gives $\mathcal B\subseteq\mathcal F$.
Thus $\operatorname{Tail}(\langle\mathcal B\rangle)\subseteq\mathcal F$ and
$\mathcal F\to x$ by refinement.
\end{proof}

\begin{example}\label{sx:Cohen-not-WSD}
There is an Archimedean Dedekind complete vector lattice whose order
convergence is strongly first countable although the lattice is not
weakly $\sigma$-distributive.
\end{example}
\begin{proof}
Let $Y=[0,1]$ and let $\mathfrak B=\operatorname{RO}(Y)$ be its complete
Boolean algebra of regular open subsets. Its Boolean operations satisfy
\[
 \bigvee_{\gamma}b_\gamma
     =\operatorname{int}_Y\overline{\bigcup_\gamma b_\gamma}^{\,Y},
 \qquad
 \bigwedge_\gamma b_\gamma
     =\operatorname{int}_Y\bigcap_\gamma b_\gamma.
\]
Let $K$ be its Stone space and write $\widehat b$ for the clopen subset
associated with $b\in\mathfrak B$. The nonempty intervals with rational
endpoints strictly between zero and one form a countable Boolean
$\pi$-base: every nonzero $b\in\mathfrak B$ contains one of them.
Every nonempty open subset of $K$ contains a nonempty clopen
$\widehat b$, hence contains $\widehat q$ for an interval $q$ in that
countable family. Therefore $K$ is compact Hausdorff with a countable
$\pi$-base.

Set $E=C(K)$. Proposition~4.21 of
\citet{BBW2025countability} now implies that ordinary order convergence
on $E$ is strongly first countable. The assertion uses exactly the
compact-Hausdorff and countable-$\pi$-base hypotheses verified above.
Corollary~4.31 of the same source gives the corresponding equivalence.

For completeness, the order-completeness assertion can be verified
within this construction. If $O=\bigcup_{b\in A}\widehat b$ is open,
put $c=\bigvee A$. Then
\[
 \overline O^{\,K}=\widehat c.
\]
Indeed, $O\subseteq\widehat c$. A nonzero $d\le c$ cannot be disjoint
from every $b\in A$, since
$d=d\wedge c=\bigvee_{b\in A}(d\wedge b)$.
Thus every nonempty clopen subset of $\widehat c$ meets $O$, proving
the reverse closure inclusion. In particular, closures of open sets
are clopen.

Let $\varnothing\ne\mathcal A\subseteq C(K)$ have an upper bound
$v\in C(K)$, and fix $f_0\in\mathcal A$. For $q\in\mathbb Q$ set
\[
 V_q=\overline{\bigcup_{f\in\mathcal A}\{t:f(t)>q\}}^{\,K},
 \qquad
 h(t)=\sup\{q\in\mathbb Q:t\in V_q\}.
\]
The sets $V_q$ are clopen, decrease with $q$, and obey
$V_q=K$ when $q<\min_K f_0$, and $V_q=\varnothing$ when
$q>\max_Kv$. Hence $h$ is real valued. For every real $a$,
\[
 \{h>a\}=\bigcup_{q>a}V_q,
 \qquad
 \{h<a\}=\bigcup_{q<a}(K\setminus V_q),
\]
so $h\in C(K)$. Each $f\in\mathcal A$ satisfies $f\le h$.
If $w\in C(K)$ is any upper bound of $\mathcal A$, then
$V_q\subseteq\{w\ge q\}$, whence $h\le w$.
Thus $h=\sup_E\mathcal A$ and $E$ is Dedekind complete.
Archimedeanity is inherited from the pointwise order.

Enumerate $\mathbb Q\cap(0,1)$ by $(r_n)_{n\ge1}$ and put
\[
 \delta_n=\tfrac12\min\{r_n,1-r_n\},\qquad
 b_{n,k}=(r_n-\delta_n/k,r_n+\delta_n/k)
       \quad(n,k\ge1).
\]
These intervals are members of $\mathfrak B$, with
\[
 \bigwedge_{k\ge1}b_{n,k}
   =\operatorname{int}_Y\{r_n\}=0.
\]
For every $\varphi\in\mathbb N^{\mathbb N}$, their diagonal union
contains every $r_n$ and is dense in $Y$. Consequently
\[
 \bigvee_{n\ge1}b_{n,\varphi(n)}=1.
\]
Define $a_{n,k}=\mathbf1_{\widehat b_{n,k}}\in E$.
Then $0\le a_{n,k}\le\mathbf1$ and $a_{n,k+1}\le a_{n,k}$.
If $0\le g\le a_{n,k}$ for every $k$ and $g\ne0$, continuity and the
clopen base give $0\ne d\in\mathfrak B$ and $\varepsilon>0$ such that
\[
 \widehat d\subseteq\{g>\varepsilon\}
       \subseteq\bigcap_k\widehat b_{n,k}.
\]
It would follow that $0<d\le\bigwedge_kb_{n,k}=0$.
Therefore $\inf_k a_{n,k}=0$ in $E$.
For each $\varphi$, the union $\bigcup_n\widehat b_{n,\varphi(n)}$
is dense in $K$, by the Boolean-join calculation and the closure
identity above. Every continuous upper bound $g$ of
$\{a_{n,\varphi(n)}:n\ge1\}$ is at least one on that dense union,
hence on all of $K$. Thus
\[
 \sup_n a_{n,\varphi(n)}=\mathbf1
 \quad(\varphi\in\mathbb N^{\mathbb N}),\qquad
 \inf_{\varphi\in\mathbb N^{\mathbb N}}
       \sup_n a_{n,\varphi(n)}=\mathbf1\ne0.
\]
This globally order-bounded array disproves weak
$\sigma$-distributivity. Strong first countability of the original
order convergence therefore cannot replace the distributive
hypothesis in arguments which use this infimum-supremum identity.
\end{proof}

\subsection{Almost-uniform net convergence and bounded lattice metrics}
\label{dx:au-section}

To see how these distinctions appear in spaces of measurable functions,
we turn to almost-uniform convergence, where the relation between
sequences and arbitrary nets also determines what can be required of
a homogeneous representation. Let
$(T,\mathcal A,\nu)$ be the completed probability space $[0,1]$ with its
usual measure. Write $\mathbf1_A$ for an indicator and $\mathbf1=\mathbf1_T$.
All measurable functions are identified modulo almost-everywhere equality.

\begin{lemma}[Essential envelopes and countable regulator extraction]
\label{dx:measure-lattice}
The vector lattice $L^\infty(T,\mathbb R)$ is Archimedean, super Dedekind
complete and weakly $\sigma$-distributive. If $D\subseteq L^\infty_+$ is
downward directed with infimum zero, it has a decreasing sequence
$d_n\in D$ with $d_n\downarrow0$. For each $\eta>0$ there is
$A\in\mathcal A$, $\nu(A^c)<\eta$, on which
$v_n=d_n\wedge\mathbf1$ converges essentially uniformly to zero.
The countable extraction assertion also holds for positive zero regulators
in $L^0(T,\mathbb R)$.
\end{lemma}
\begin{proof}
For a nonempty family $\mathcal C$ of measurable sets put
\[
 s=\sup\{\nu(C_1\cup\cdots\cup C_m):m\ge1,\ C_j\in\mathcal C\}.
\]
For every $n$ a finite union $V_n$ has measure greater than $s-1/n$.
The set $V=\bigcup_n V_n$ has measure $s$, since each finite union of the
$V_n$ has measure at most $s$. For $C\in\mathcal C$,
\[
 \nu(V\cup C)=\lim_N\nu(V_1\cup\cdots\cup V_N\cup C)\le s=\nu(V),
\]
so $\nu(C\setminus V)=0$. Every measurable set containing every member
of $\mathcal C$ modulo null sets also contains $V$ modulo a null set.
Thus an arbitrary essential union is attained by a countable subfamily.

Apply this assertion to the strict superlevel sets at every rational level
of a nonempty family of measurable real functions. The union of the
resulting countably many function indices has a pointwise supremum $u$
which is measurable, and every original member is at most $u$ almost
everywhere. Any almost-everywhere upper bound of the family dominates
the countable subfamily and hence $u$. This constructs its essential
supremum, with extended real values allowed, from a countable subfamily.
For a family bounded above and below by finite measurable functions,
the envelope is finite almost everywhere. For a family in $L^\infty$
bounded in order there, it belongs to $L^\infty$.
For a family $S$ with upper bound $u$, fix $a\in S$ and apply the bounded
case to $\{(s-a)^+:s\in S\}\subseteq[0,(u-a)^+]$.
Its supremum $v$ gives $\sup S=a+v$. The countable family realizing $v$
together with $a$ realizes $\sup S$. Hence both
$L^0$ and $L^\infty$ are Dedekind complete and have the countable sup
property. Archimedeanity follows from the corresponding real inequality
outside the union of the countably many exceptional null sets in
$nf\le g$, $n\ge1$.

For weak $\sigma$-distributivity in $L^\infty$, let
$0\le a_{n,k}\le e$ and $a_{n,k}\downarrow_k0$.
Representatives can be fixed on one full-measure set so that all row
monotonicity and limit statements hold pointwise. Given
$\varepsilon,\eta>0$, each row supplies $k_n$ with
\[
 \nu\{a_{n,k_n}>\varepsilon\}<\eta2^{-n},\qquad
 \nu\{\sup_n a_{n,k_n}>\varepsilon\}<\eta.
\]
Let $w=\inf_\varphi\sup_n a_{n,\varphi(n)}$. The preceding inequalities
give $\nu\{w>\varepsilon\}<\eta$ for every $\eta>0$, so
$\nu\{w>\varepsilon\}=0$. Taking rational $\varepsilon>0$ gives $w=0$.

The countable sup property and the directedness argument in
Theorem~\ref{dx:order-dictionary} supply $d_n\in D$, $d_n\downarrow0$,
both in $L^\infty$ and in $L^0$. In either case
$v_n=d_n\wedge\mathbf1\downarrow0$ almost everywhere. For each $k$ there
is $N_k$ with
\[
 \nu\{v_{N_k}>1/k\}<\eta2^{-k-1}.
\]
Replacing $N_k$ by $\max_{j\le k}N_j$ preserves this bound. On
$A=\bigcap_k\{v_{N_k}\le1/k\}$ one has $v_n\le1/k$ for $n\ge N_k$,
and $\nu(A^c)<\eta$. This proves the final assertion.
\end{proof}

Let $(H,\rho)$ be a nonempty complete metric space. Denote by
$L^0(T,H)$ the measurable, essentially separably valued functions modulo
almost-everywhere equality. A net $(f_\alpha)_{\alpha\in I}$ in this
space converges \emph{almost uniformly} to $f$ if
\begin{equation}\label{dx:au-definition}
 \begin{gathered}
 (\forall\eta>0)(\exists A\in\mathcal A)\quad\nu(A^c)<\eta,\\
 (\forall\delta>0)(\exists\alpha_\delta)(\forall\alpha\ge\alpha_\delta)
 \quad\mathop{\rm ess\,sup}_{t\in A}
       \rho(f_\alpha(t),f(t))\le\delta.
 \end{gathered}
\end{equation}
In this definition $A$ is independent of $\delta$, and the abbreviation
AU will always retain this arbitrary-net meaning. Put
\begin{equation}\label{dx:au-distance}
 d(f,g)=\min\{1,\rho(f(\cdot),g(\cdot))\}\in L^\infty(T)_+.
\end{equation}

\begin{theorem}[Arbitrary-net representation and completeness of AU]
\label{dx:au-faithfulness}
The distance in \eqref{dx:au-distance} is definite and symmetric and
satisfies the triangle inequality in $L^\infty(T)$. For every net
$(f_\alpha)$ and every $f$,
\begin{equation}\label{dx:au-exactness}
 f_\alpha\longrightarrow_{\rm AU}f
 \quad\Longleftrightarrow\quad d(f_\alpha,f)\longrightarrow_o0.
\end{equation}
This convergence is separated, first countable and complete for the
uniform convergence represented by $d$. If $H$ is a group and $\rho$ is
bi-invariant, AU is a convergence group and this represented uniform
convergence is its canonical group uniform convergence.
\end{theorem}
\begin{proof}
Suppose \eqref{dx:au-definition} holds and write
$\mathcal F=\operatorname{Tail}((f_\alpha))$ for its eventual sets.
For each $n$ let $B_n$ satisfy it
with $\nu(B_n^c)<2^{-n}$, and put $A_n=\bigcup_{k\le n}B_k$.
For every $\delta>0$, the intersection of the finitely many eventual sets imposing the $\delta$ bound on $B_1,\ldots,B_n$ imposes it on
$A_n$. Thus
\[
 r_n=n^{-1}\mathbf1_{A_n}+\mathbf1_{A_n^c}\downarrow0,
 \qquad \{g:d(g,f)\le r_n\}\in\mathcal F.
\]
The decrease is pointwise, and the limit vanishes almost everywhere since
$\nu(A_n^c)<2^{-n}$. This proves the forward implication.

Conversely, a regulator $D\downarrow0$ for the image has a sequence
$d_n\in D$ and good set $A$ supplied by Lemma~\ref{dx:measure-lattice}.
The global bound $d(g,f)\le\mathbf1$ makes
$v_n=d_n\wedge\mathbf1$ a valid bound on an eventual set.
For $\delta>0$ an index $n$ satisfies
$v_n\le\min(\delta,1/2)$ essentially on $A$. On that eventual set,
\[
 \min(1,\rho(g,f))\le\min(\delta,1/2)\quad\hbox{a.e. on }A,
 \qquad \rho(g,f)\le\delta\quad\hbox{a.e. on }A.
\]
The same $A$ works for every $\delta$, proving
\eqref{dx:au-definition}. Definiteness and symmetry come from $\rho$,
and
$\min(1,a+b)\le\min(1,a)+\min(1,b)$ proves the triangle inequality.

Constant nets and refinements satisfy the convergence axioms by
\eqref{dx:au-exactness}. If $D,E\downarrow0$ regulate two net-tail families
converging to $f$, the balls with radii $d+e$, $d\in D$, $e\in E$, belong
to both net-tail families: each contains the corresponding radius-$d$ and radius-$e$
balls. Independent sums of zero regulators have infimum zero, since
$h\le d+e$ for all $e$ implies $h\le d$, and then $h\le0$.
This proves finite-intersection stability. If a single net converges
to $f$ and $g$, the intersection of its two regulator balls is nonempty,
and the triangle inequality gives $d(f,g)\le d+e$ for all regulator
members. Thus $d(f,g)=0$ and $f=g$. Countable regulator extraction gives
a null saturation of a decreasing sequence of balls at each limit, of which the original convergent net is a quasi-subnet.

For completeness, let $(g_\alpha)$ be Cauchy in the represented
uniform convergence and put $\mathcal F=\operatorname{Tail}((g_\alpha))$.
Countable extraction and the independent-product definition for nets
give decreasing $r_n\downarrow0$ and decreasing nonempty
$C_n\in\mathcal F$ with
\begin{equation}\label{dx:au-cauchy-balls}
 d(g,h)\le r_n\wedge\mathbf1\qquad(g,h\in C_n).
\end{equation}
Write $r_n\wedge\mathbf1$ as $r_n$ henceforth. For each $n$ let
$f_n\in C_n$. The countable family admits representatives for which
all pair inequalities and $r_n\downarrow0$ hold outside one null set.
For $m,n\ge N$,
$\min(1,\rho(f_m(t),f_n(t)))\le r_N(t)$. At such a point $t$ and for
every $0<\varepsilon<1$, sufficiently large $N$ gives
$\rho(f_m(t),f_n(t))\le\varepsilon$. Completeness of $H$ supplies a
limit $f(t)$. The closure of the union of the countably many essential
separable ranges is separable, and $f$ has its range there outside a null
set. The limit is measurable: for a closed subset $C$ of this separable
metric subspace, the real functions
$\rho(f_n(\cdot),C)$ are measurable and converge to
$\rho(f(\cdot),C)$. A countable base of open balls then verifies
measurability of $f$. On the exceptional null set assign a fixed member
of $H$.

Fix $N$ and $g\in C_N$. For every $n\ge N$,
$d(g,f_n)\le r_N$ almost everywhere. Outside the union of these
countably many null sets, passage to the pointwise limit gives
$d(g,f)\le r_N$. This assertion holds separately for each $g\in C_N$,
which is exactly the set inclusion
$C_N\subseteq\{g:d(g,f)\le r_N\}$. Hence $\mathcal F\to_{\rm AU}f$.

In a group with a bi-invariant metric,
\[
 \rho(ab,a'b')\le\rho(a,a')+\rho(b,b'),\qquad
 \rho(a^{-1},b^{-1})=\rho(a,b).
\]
These inequalities and independent sums of regulators prove continuity
of multiplication and inversion. Moreover
$d(f,g)=d(fg^{-1},\mathbf e)$, where $\mathbf e$ is the constant identity
function. Applying \eqref{dx:au-exactness} to the difference image of any
pair net proves exact agreement with canonical group uniformity.
\end{proof}

\begin{corollary}[The Banach-valued convergence space]\label{dx:au-banach}
Let $B$ be a nonzero real Banach space and
$X=L^0(T,B)$ the strongly measurable functions modulo almost-everywhere
equality. With AU convergence, $X$ is a complete separated first-countable
locally convex real convergence vector space. The gauge
\[
 p(f)=\min(1,\|f(\cdot)\|_B)\in L^\infty(T)_+
\]
is faithfully invariant, globally bounded and subadditive. All its balls
are closed disks. The convergence is countably sequentially determined,
and satisfies closed-base stability and the two-scale property of the
preceding intrinsic convergence calculus.
\end{corollary}
\begin{proof}
Theorem~\ref{dx:au-faithfulness} applies to the complete norm metric.
To check scalar multiplication, let $\mathcal R\to a$ in $\mathbb R$
and $\mathcal F\to_{\rm AU}f$. For $\eta>0$ there is a good set
$A_0$ of complement measure less than $\eta/2$. Since $\|f\|_B$ is
finite almost everywhere, a number $M\ge1$ has
$\nu\{\|f\|_B>M\}<\eta/2$. Put
$A=A_0\cap\{\|f\|_B\le M\}$. There is $K\ge1$ with
$[-K,K]\in\mathcal R$. For every $\delta>0$ the product eventual-set family has
a member on which
\[
 |t|\le K,\qquad |t-a|\le\delta/(2M),\qquad
 \mathop{\rm ess\,sup}_{A}\|g-f\|_B\le\delta/(2K).
\]
On this member,
$\|tg-af\|_B\le |t|\|g-f\|_B+|t-a|\|f\|_B\le\delta$
essentially on $A$. This proves joint scalar continuity.

For $r\in L^\infty_+$, the ball $p(f)\le r$ imposes
$\|f(t)\|_B\le r(t)$ on $\{r<1\}$ and no condition on $\{r\ge1\}$.
It is balanced and convex. It is closed: if a net with values in that ball converges to $f$, then the triangle inequality gives
$p(f)\le r+d$ for every member of an error regulator $D\downarrow0$,
and subtraction gives $p(f)\le r$. Every null tail family has a coarser null
ball saturation net, proving local convexity and closed-base stability.
For the two-scale implication, use the notation of Definition~\ref{nx:P}.
Let $D\downarrow0$ regulate $p_*\mathcal F$, and for each $d\in D$ let
$E_d\downarrow0$ regulate $p_*\mathcal H_{B_d}$. The bound set
\[
 R=\{r\in L^\infty_+:B_r\in\mathcal K\}
\]
contains $d+e$ for every $d\in D$, $e\in E_d$: a member of
$\mathcal H_{B_d}$ has image bounded by $e$, and its sum with $B_d$
is a member of $\mathcal K$ contained in $B_{d+e}$.
The identity $B_r\cap B_s=B_{r\wedge s}$ makes $R$ downward directed.
Any lower bound $v$ of $R$ satisfies
$v-d\le e$ for all $e\in E_d$, then $v\le d$ for all $d\in D$,
and then $v\le0$. Thus $R\downarrow0$, and faithfulness gives
$\mathcal K\to0$. Finally,
Lemma~\ref{dx:measure-lattice} and Theorem~\ref{dx:bounded-csd} prove
countable sequential determination.
\end{proof}

\begin{proposition}[Sequences, nets and nonhomogeneous representation]
\label{dx:au-sequences-obstruction}
For sequences in $L^0(T,B)$, AU convergence, almost-everywhere norm
convergence and order convergence of the untruncated norm errors in
$L^0(T)$ are equivalent. For nets, untruncated norm-order
convergence implies AU, and the implication is strict. AU is not a
topological convergence. It has no single faithful homogeneous
lattice-norm representation into any Archimedean vector lattice.
\end{proposition}
\begin{proof}
For an AU-null sequence, countable regulator extraction supplies
$r_k\downarrow0$ and integers $N_k$ such that
$\min(1,\|f_n\|_B)\le r_k$ for $n\ge N_k$. There are countably many
inequalities, so outside one null set each holds pointwise and the
sequence tends to zero in norm. Conversely, for an almost-everywhere
norm-null sequence,
\[
 u_N(t)=\sup_{n\ge N}\|f_n(t)\|_B
\]
is finite almost everywhere, measurable, and decreases to zero.
It regulates the untruncated norms in $L^0$, while $u_N\wedge\mathbf1$
regulates the truncated norms in $L^\infty$. For arbitrary nets,
a positive $L^0$ regulator $D\downarrow0$ yields
$\{d\wedge\mathbf1:d\in D\}\downarrow0$ in $L^\infty$: every positive
$L^\infty$ lower bound would also be a lower bound of $D$ in $L^0$.
Thus untruncated norm-order convergence implies AU.

Let
\[
 V_n=\{f:f=0\text{ a.e. outside }(0,1/n)\},\qquad
 \mathcal V=\operatorname{Tail}(\langle V_n:n\ge1\rangle).
\]
The estimates $p(V_n)\le\mathbf1_{(0,1/n)}\downarrow0$ prove
$\mathcal V\to_{\rm AU}0$. No $V_n$ has order-bounded untruncated
norm image. Indeed, for a fixed $0\ne b\in B$ it contains
$k b\mathbf1_{(0,1/n)}$ for all positive integers $k$, whose norms
have no finite almost-everywhere bound on $(0,1/n)$. Hence the
untruncated image of $\mathcal V$ is not order-null. Every member of
$\mathcal V$ contains a nonzero vector subspace. Proposition
\ref{dx:homogeneous-local-boundedness} therefore proves that AU has no
faithful homogeneous Archimedean lattice norm, even after a linear
convergence isomorphism.

For non-topologicity, enumerate the dyadic half-open intervals by
\[
 I_{2^m+k}=[k2^{-m},(k+1)2^{-m}),\qquad m\ge0,\quad0\le k<2^m,
\]
and put $f_n=b\mathbf1_{I_n}$ for a fixed $b\ne0$. Every point in
$[0,1)$ belongs to one interval at every level, so $(f_n)$ is not
almost-everywhere null and is not AU-null. Every subsequence has a
further subsequence $(f_{n_j})$ with $\nu(I_{n_j})\le2^{-j}$.
The inequality
\[
 \nu\left(\bigcup_{j\ge N}I_{n_j}\right)\le\sum_{j\ge N}2^{-j}
 \longrightarrow0
\]
shows that this further subsequence is almost-everywhere null, hence
AU-null. In a topology, failure of convergence to zero supplies a
neighborhood $U$ and a subsequence wholly outside $U$. No further
subsequence of it can converge to zero. The constructed sequence
contradicts this necessary property.
\end{proof}

\begin{proposition}[Identical sequences, closed sets and bounded sets]
\label{dx:au-invariants}
On $X=L^0(T,B)$, AU and the convergence
\[
 \mathcal F\longrightarrow_{\rm bo}0
 \quad\Longleftrightarrow\quad
 (\|\,\cdot\,\|_B)_*\mathcal F\longrightarrow_o0\text{ in }L^0(T)
\]
have the same convergent sequences, the same convergence-closed sets
and the same intrinsically bounded subsets. A nonempty $M\subseteq X$
is bounded for either convergence if and only if
\[
 w_M=\mathop{\rm ess\,sup}_{f\in M}\|f(\cdot)\|_B
\]
is finite almost everywhere. The two net convergences nevertheless
differ.
\end{proposition}
\begin{proof}
Starting from the equality of sequences and the distinction between
nets proved in Proposition~\ref{dx:au-sequences-obstruction}, we compare
closed sets. AU is first countable, and Lemma~\ref{dx:measure-lattice} gives first countability for bo by
pulling back the extracted $L^0$ regulator balls. In any first-countable
convergence space, an eventual-set family of a net $\mathcal F$ concentrated on $C$ and
converging to $x$ has a coarser convergent net with decreasing base
$(G_n)$. Each $G_n\cap C$ is nonempty. Elements
$x_n\in G_n\cap C$ form a sequence converging to $x$. Thus sequential
closedness implies convergence-closedness, and the reverse implication
uses the eventual-set family of a sequence. Equality of sequences gives equality
of closed sets.

For boundedness, the essential envelope exists as an extended measurable
function by the superlevel-set argument in Lemma~\ref{dx:measure-lattice}.
If $M$ is AU-bounded, for each $m$ there is a set $C_m$ with
$\nu(C_m^c)<2^{-m}$ on which its scalar-dilation saturation net converges
uniformly. At accuracy one, some $\varepsilon_m>0$ satisfies
\[
 \|tf\|_B\le1\quad\hbox{a.e. on }C_m
 \qquad(f\in M,\ |t|\le\varepsilon_m).
\]
Taking $t=\varepsilon_m$ gives
$w_M\le1/\varepsilon_m$ almost everywhere on $C_m$.
Since $\nu((\bigcup_m C_m)^c)=0$, $w_M$ is finite almost everywhere.
Conversely,
\[
 p(tf)\le\min(1,\varepsilon w_M)\downarrow0
 \quad(|t|\le\varepsilon,\ f\in M)
\]
proves AU-boundedness. The same finite envelope gives bo-boundedness
through $\|tf\|_B\le\varepsilon w_M\downarrow0$ in $L^0$.
If $M$ is bo-bounded, one regulator bound gives $u\in L^0_+$ and
$\varepsilon>0$ with $\|\varepsilon f\|_B\le u$ for every $f\in M$.
Therefore $w_M\le u/\varepsilon$ almost everywhere. This proves all
boundedness equivalences.
\end{proof}

\begin{proposition}[No nonzero continuous homogeneous Archimedean gauges]
\label{dx:au-zero-gauges}
Every positively homogeneous subadditive map
$q:L^0(T,B)\to G_+$ into an Archimedean partially ordered real vector
space which sends AU-null nets to order-null nets is zero.
Here regulator infima are taken in the whole space $G$.
More generally, every continuous linear map from AU into a Hausdorff
convergence vector space in which every null net has a bounded
member is zero.
\end{proposition}
\begin{proof}
For a measurable $A$, write $S(A)$ for the subspace of functions
vanishing almost everywhere outside $A$, and put
$\mathcal I=\{A:q(S(A))=\{0\}\}$. This collection contains null sets
and is hereditary. If $A,C\in\mathcal I$ and $f\in S(A\cup C)$,
\[
 f=f\mathbf1_A+f\mathbf1_{C\setminus A},\qquad
 0\le q(f)\le q(f\mathbf1_A)+q(f\mathbf1_{C\setminus A})=0.
\]
Thus $A\cup C\in\mathcal I$. If $q\ne0$, then $T\notin\mathcal I$.
Every measurable positive-measure subset of $[0,1]$ splits into two
equal-measure measurable pieces: the continuous function
$t\mapsto\nu(A\cap[0,t])$ has values zero and $\nu(A)$ at the endpoints.
At least one half of a set outside $\mathcal I$ remains outside
$\mathcal I$. Recursion produces
\[
 A_{n+1}\subseteq A_n,\qquad \nu(A_n)=2^{-n},\qquad A_n\notin\mathcal I.
\]
The saturation $\langle S(A_n):n\ge1\rangle$ is AU-null, since its tails are
eventually zero outside each fixed $A_N$.
Its $q$-image has a bound $r\in G_+$ on some $S(A_N)$. For every
$f\in S(A_N)$ and integer $k\ge1$,
\[
 kq(f)=q(kf)\le r.
\]
Archimedeanity gives $q(f)=0$, so $A_N\in\mathcal I$, a contradiction.
For the linear-map assertion, suppose $T_0\ne0$ and set
$\mathcal J=\{A:T_0(S(A))=\{0\}\}$.
Linearity and the decomposition
$f=f\mathbf1_A+f\mathbf1_{C\setminus A}$ show that $\mathcal J$
is hereditary and closed under finite unions, and $T\notin\mathcal J$.
Starting with $A_0=T$, split each $A_n\notin\mathcal J$ into two
equal-measure pieces. Both pieces cannot belong to $\mathcal J$,
so one of them, denoted by $A_{n+1}$, satisfies
$A_{n+1}\notin\mathcal J$ and $\nu(A_{n+1})=2^{-n-1}$.
The saturation of $(S(A_n))$ is AU-null, hence its image under $T_0$
is null by continuity. Its exact tail images are $T_0(S(A_n))$.
The assumed bounded-tail property in the range supplies $N$ and a
bounded set containing $T_0(S(A_N))$. This is a nonzero vector
subspace because $A_N\notin\mathcal J$, contrary to the second
assertion of Proposition~\ref{dx:homogeneous-local-boundedness}.
Thus $T_0=0$.
\end{proof}

\subsection{Finite-subsum convergence and initial products}
\label{dx:permanence-section}

The restrictions imposed by a homogeneous representation become clearer
when supports and coordinates are treated separately. We begin with
convergence confined to a finite subsum, then compare it with the
initial convergence of a product, where no finite support is required.

\begin{theorem}[Convergence in a finite subsum]\label{dx:finite-subsum}
Let $(B_i)_{i\in I}$ be real Banach spaces and
$Y=\bigoplus_{i\in I}B_i$. Declare $\mathcal F\to0$ if some finite
$F\subseteq I$ satisfies $Y_F=\bigoplus_{i\in F}B_i\in\mathcal F$
and the restricted tail net on $Y_F$ converges to zero in the finite product
of norm convergences. Define all other limits by translation. Then
\[
 q(y)=(\|y_i\|_{B_i})_{i\in I}\in\mathbb R^{(I)}_+
\]
is a faithful homogeneous lattice norm for all nets. The
original convergence is separated, complete, first countable and locally
convex, and satisfies closed-base stability and the two-scale property
of Definition~\ref{nx:P}. If infinitely many $B_i$ are nonzero, it is not topological.
The value lattice $\mathbb R^{(I)}$ is super Dedekind complete and
weakly $\sigma$-distributive for every cardinality of $I$.
\end{theorem}
\begin{proof}
The constant zero net is supported in $Y_\varnothing=\{0\}$.
If a net is eventually supported in $Y_F$ and converges to zero in
each of its coordinates, a quasi-subnet retains both assertions.
For two null nets $(x_\alpha)$ and $(y_\beta)$ with finite supports
$F$ and $H$, the independent sum and every finite mixing are
eventually supported in $Y_{F\cup H}$. For $i\in F\cup H$,
\[
 \|x_\alpha(i)+y_\beta(i)\|
 \le\|x_\alpha(i)\|+\|y_\beta(i)\|\longrightarrow0.
\]
For a mixing, both possible values in coordinate $i$ have norm at
most $\max(\|x_\alpha(i)\|,\|y_\beta(i)\|)$, which also tends
to zero. Negation leaves these bounds and supports unchanged.
For $t_\gamma\to t$ and $x_\alpha\to x$, the differences
$x_\alpha-x$ have one finite eventual support $F$. On a scalar
tail $|t_\gamma|\le M$, the vectors
$t_\gamma x_\alpha-tx$ are supported in
$F\cup\operatorname{supp}x$ and satisfy, coordinate by coordinate,
\[
 \|t_\gamma x_\alpha(i)-tx(i)\|
 \le M\|x_\alpha(i)-x(i)\|
       +|t_\gamma-t|\,\|x(i)\|\longrightarrow0.
\]
Thus joint scalar multiplication is continuous. If a net converges
to both $x$ and $y$, each normed coordinate converges to both
$x(i)$ and $y(i)$, so $x(i)=y(i)$ for every $i$ and $x=y$.

If $q_*\mathcal F$ is order-null, one regulator bound $u$ has finite
support $F$, and
$\{y:q(y)\le u\}\subseteq Y_F$ gives $Y_F\in\mathcal F$.
For a positive downward zero regulator $D$ in $\mathbb R^{(I)}$, each
coordinate projection has infimum zero: a positive lower bound $c$ in
coordinate $i$ would make $c\mathbf1_{\{i\}}$ a nonzero lower bound
of $D$. Its image balls prove convergence of the corresponding norm
coordinate. Thus the trace of $\mathcal F$ converges in $Y_F$.
Conversely, if that trace converges, finite intersections of coordinate
norm bounds give
\[
 \{y:q(y)\le n^{-1}\mathbf1_F\}\in\mathcal F\qquad(n\ge1),
 \qquad n^{-1}\mathbf1_F\downarrow0.
\]
This proves arbitrary-net faithfulness and first countability.
The displayed balls are disks, giving local convexity. For any fixed
radius $u$, a limit $y$ of a convergent net eventually in $\{q\le u\}$ satisfies
$q(y)\le u+d$ for every error-regulator bound $d$, hence $q(y)\le u$.
Thus all balls are closed, and their null bases give closed-base stability.
For the two-scale property, let $\mathcal F\to0$ be the outer
family in Definition~\ref{nx:P}, and let $\mathcal K$ be its candidate
family. A regulator $D\downarrow0$ for $q_*\mathcal F$ has
$B_d=\{q\le d\}\in\mathcal F$. Let $\mathcal H_d\to0$ be the
inner family corresponding to $B_d$, and let $E_d\downarrow0$
regulate $q_*\mathcal H_d$. For each $e\in E_d$,
$B_e\in\mathcal H_d$, so
\[
 B_d+B_e\in\mathcal K,\qquad
 B_d+B_e\subseteq B_{d+e},\qquad
 d+e\in R:=\{r\in\mathbb R^{(I)}_+:B_r\in\mathcal K\}.
\]
This proves nonemptiness of $R$. If $r,s\in R$, then
$B_{r\wedge s}=B_r\cap B_s\in\mathcal K$, so $R$ is downward
directed. Any lower bound $v$ of $R$ satisfies
$v-d\le e$ for every $e\in E_d$. Their infimum gives $v\le d$
for every $d\in D$, and then $v\le0$. Thus $R\downarrow0$,
and its ball memberships prove $q_*\mathcal K\to_o0$.
Faithfulness gives $\mathcal K\to0$.

For completeness, let $\mathcal F-\mathcal F$ be null. It contains
$Y_F$ for one finite $F$, so there is $A\in\mathcal F$ with
$A-A\subseteq Y_F$. Fix $a\in A$. Then
$A\subseteq a+Y_F$, and the trace of $\mathcal F-a$ on $Y_F$ is Cauchy
in each of its finitely many Banach coordinates. Let $z\in Y_F$ be
the tuple of their limits. The trace converges to $z$ in the finite
product. Since $a+Y_F$ is contained in the finite subsum with support
$F\cup\operatorname{supp}a$, this convergence proves
$\mathcal F\to a+z$ in the original convergence.

Suppose it were induced by a topology. Its neighborhood family at zero
would be null, hence would contain a finite subsum $Y_F$.
For every $y\in Y$, continuity of scalar multiplication gives
$ty\to0$ as $t\to0$ in $\mathbb R$. A nonzero sufficiently small $t$
would satisfy $ty\in Y_F$, and linearity would give $y\in Y_F$.
Thus $Y=Y_F$, contradicting a nonzero factor outside $F$.

A positive family bounded by $u\in\mathbb R^{(I)}_+$ has all its
supports in the finite set $\operatorname{supp}u$. Its supremum and a
countable realizing subfamily are obtained in these finitely many
real coordinates. For an arbitrary nonempty family $S$ bounded above,
fix $a\in S$ and use the positive bounded family
$\{(s-a)^+:s\in S\}$. This proves super Dedekind completeness.
If $0\le a_{n,m}\le u$ and $a_{n,m}\downarrow_m0$, then for every
$k$ and every $n$ there is $\varphi_k(n)$ with
\[
 a_{n,\varphi_k(n)}(i)\le 1/k
 \qquad(i\in\operatorname{supp}u).
\]
The finite coordinate set permits the maximum of the finitely many
row indices. Hence
$\sup_n a_{n,\varphi_k(n)}\le k^{-1}\mathbf1_{\operatorname{supp}u}$,
and the path-supremum infimum is zero. This is weak
$\sigma$-distributivity.
\end{proof}

\begin{proposition}[A single homogeneous norm for an initial product]
\label{dx:finite-product-criterion}
Let $(Y_i)_{i\in I}$ be Hausdorff real convergence vector spaces, each
having a faithful homogeneous lattice norm into an Archimedean vector
lattice. Give $Y=\prod_iY_i$ the initial coordinate convergence.
Then $Y$ has a single faithful homogeneous lattice norm into an
Archimedean vector lattice if and only if only finitely many factors
are nonzero.
\end{proposition}
\begin{proof}
For finitely many nonzero factors, take the tuple of their norms into
the finite product of their value lattices. Finite products of
coordinate regulators are positive, downward directed and have
coordinatewise infimum zero. Finite intersections of their ball
conditions are eventual sets. Therefore coordinate order convergence
is equivalent to order convergence in this finite product, proving
faithfulness. Homogeneity and Archimedeanity hold coordinatewise.

For infinitely many nonzero factors, put
\[
 W_F=\{y:y_i=0\ (i\in F)\},\qquad
 \mathcal W=\operatorname{Tail}(\langle W_F:F\subseteq I\text{ finite}\rangle).
\]
Each coordinate image is eventually identically zero, so
$\mathcal W\to0$. Every $W_F$ is a nonzero vector subspace, and every
member of $\mathcal W$ contains one. Proposition
\ref{dx:homogeneous-local-boundedness} shows that this null tail family
has no bounded member. The same proposition excludes every single
faithful homogeneous Archimedean lattice norm for $Y$.
\end{proof}

\begin{example}[Absorbing closed disk bases without a homogeneous realization]
\label{dx:absorbing-product}
Let $E=L^\infty(T)$ have order convergence and let
$Y=\prod_{n\ge1}E$ have the initial coordinate order convergence.
It is separated, complete, first countable and locally convex, and
satisfies closed-base stability and the two-scale property.
Every null net is a quasi-subnet of a null saturation with a base of absorbing
closed disks. Nevertheless this convergence is not topological and
has no single faithful homogeneous lattice norm into an Archimedean
vector lattice. Its coordinate gauges form a faithful lattice-valued
family, and Theorem~\ref{nx:packaging} gives one faithful gauge into
\[
 \prod_{n\ge1}\bigl(\mathbb Z\mathbin{\rm lex}E\bigr).
\]

For completeness, a Cauchy net has a Cauchy image in every
coordinate. Theorem~\ref{dx:order-dictionary} supplies coordinate
limits, and their tuple is a limit in the initial convergence.
For first countability, a null tail family has, in each coordinate, a
coarser null tail family with a decreasing countable base, by
Lemma~\ref{dx:measure-lattice}. Pulling back all these bases and taking
finite intersections gives a countable coarser null base in $Y$.
Order intervals in $E$ are disks, and finite cylinders of disks
prove local convexity.

For the absorbing assertion, let $D\downarrow0$ regulate a coordinate
null tail family and use radii
\[
 d+k^{-1}\mathbf1,\qquad d\in D,\quad k\ge1.
\]
They form a downward zero regulator: a lower bound $v$ satisfies
$v-d\le k^{-1}\mathbf1$ for all $k$, hence $v\le d$ for all $d$,
and hence $v\le0$. Their intervals belong to the original eventual-set family
and contain $[-k^{-1}\mathbf1,k^{-1}\mathbf1]$. Each is therefore
absorbing in $E$. These intervals are order-closed, since passage
to an order limit preserves both endpoint inequalities. Their finite
coordinate cylinders in $Y$ are absorbing closed disks and generate
a coarser null tail family.

Having obtained closed disk bases, we have closed-base stability as
well. To verify the two-scale property, let $\mathcal F\to0$ and $\mathcal K$
satisfy its antecedent, and fix a coordinate $n$. Put
$B_{n,r}=\{y\in Y:|y_n|\le r\}$ for $r\in E_+$.
A regulator $D_n\downarrow0$ of the $n$th coordinate of $\mathcal F$
has $B_{n,d}\in\mathcal F$. For each $d\in D_n$, let
$\mathcal H_{n,d}\to0$ be the associated inner family and let
$E_{n,d}\downarrow0$ regulate its $n$th coordinate. Then
\[
 B_{n,d}+B_{n,e}\in\mathcal K,\qquad
 B_{n,d}+B_{n,e}\subseteq B_{n,d+e}
 \quad(e\in E_{n,d}).
\]
Consequently the set
$R_n=\{r\in E_+:B_{n,r}\in\mathcal K\}$ is nonempty, contains
all $d+e$, and is closed under binary meets because
$B_{n,r\wedge s}=B_{n,r}\cap B_{n,s}$.
If $v$ is a lower bound of $R_n$, then
$v-d\le e$ for every $e\in E_{n,d}$, yielding $v\le d$ for
all $d\in D_n$ and hence $v\le0$. Thus $R_n\downarrow0$.
Its ball memberships prove order nullity of the $n$th coordinate
of $\mathcal K$. This holds for every $n$, so the defining product
convergence gives $\mathcal K\to0$.

We can now apply Proposition~\ref{dx:finite-product-criterion} to exclude
a homogeneous Archimedean norm. For non-topologicity,
the subspace supported in the first coordinate carries exactly order
convergence on $E$. The dyadic indicators $\mathbf1_{I_n}$ used in
Proposition~\ref{dx:au-sequences-obstruction} are bounded by
$\mathbf1$, and their tail suprema equal $\mathbf1$ almost everywhere,
so they are not order-null. Every subsequence has a further subsequence
with summable support measures. The tail suprema of that further
subsequence decrease to zero almost everywhere and regulate it in
$E$. The topological subsequence argument in that proposition excludes
a topology for this coordinate subspace, and therefore for $Y$.
The packaged value group above need not be Archimedean. Its existence
does not contradict the homogeneous Archimedean obstruction.
\end{example}

\begin{remark}[Metric-group examples]\label{dx:metric-group-examples}
Theorem~\ref{dx:au-faithfulness} applies to the circle group with its
geodesic distance, finite groups with their discrete bi-invariant
distance, and the unitary matrix group with an invariant norm distance.
For the two-element group, the resulting function group is the measure
algebra under symmetric difference, with lattice-valued distance
$\mathbf1_{A\triangle B}$. These original groups may have torsion.
They are not asserted to embed as groups in the ordered value group.
The metric records errors, and its definiteness and triangle inequality
do not identify the group law of its domain with addition in its range.
\end{remark}

\begin{proposition}[Bounded products and controller properties]
\label{dx:bounded-product-properties}
Arbitrary products of Dedekind complete weakly $\sigma$-distributive
vector lattices are Dedekind complete and weakly $\sigma$-distributive.
A countable product of super Dedekind complete vector lattices is
super Dedekind complete. Consequently, a countable product of the
Banach-valued AU spaces of Corollary~\ref{dx:au-banach} is complete,
separated, locally convex, first countable and countably sequentially
determined, and has one globally bounded faithful lattice-valued gauge.
Arbitrary products retain completeness, separation, local convexity
and a globally bounded faithful lattice-valued gauge. No first-countability
conclusion for an arbitrary uncountable product is asserted.
\end{proposition}
\begin{proof}
For $E=\prod_{i\in I}E_i$, an order-bounded nonempty family has
coordinatewise suprema and infima, which form elements of $E$.
If $0\le a_{n,m}\le u$ and $a_{n,m}\downarrow_m0$ in $E$, then for
each $i$,
\[
 \left(\inf_\varphi\sup_n a_{n,\varphi(n)}\right)_i
 =\inf_\varphi\sup_n a_{n,\varphi(n)}(i)=0.
\]
This proves weak $\sigma$-distributivity. If $I$ is countable and
$S\subseteq E_+$ has supremum $s$, each coordinate has a countable
subfamily $S_i\subseteq S$ with supremum $s_i$ in that coordinate.
The countable union $S_0=\bigcup_iS_i$ satisfies
$\sup S_0=s$. Positive differences from one fixed family member
treat a general bounded-above family, proving super Dedekind completeness.

For AU factors $X_i$, let $p_i$ be their gauges bounded by $\mathbf1_i$.
Corollary~\ref{nx:bounded-products} gives the faithful gauge
\[
 p((x_i)_i)=(p_i(x_i))_i\in\prod_i L^\infty(T),
 \qquad 0\le p(x)\le(\mathbf1_i)_i.
\]
Coordinate limits of a Cauchy net give completeness in the initial
product convergence. Coordinate separation gives separation, and
finite intersections of inverse images of convex balanced null-base
members give local convexity. For countable $I$, the value lattice is
super Dedekind complete and weakly $\sigma$-distributive by the first
part. Theorem~\ref{dx:bounded-csd} gives first countability and countable
sequential determination. The uncountable case retains the stated
structural conclusions, while the unused-coordinate example
\ref{dx:unused-coordinates} exhibits the possible failure of the ambient
countable sup property.
\end{proof}

\subsection{Changing the value object while preserving null nets}
\label{sx:conversion}

These constructions have allowed the value object to be chosen along
with the representation. If a representation is already prescribed,
changing its values raises a separate question: which zero regulators
survive, and which new regulators become available? We begin with the
passage from a cancellative monoid to its group of differences.

Let $M$ be a commutative cancellative monoid with least element $0$ and a translation-compatible order which also admits order cancellation:
\[
 a+c\le b+c\quad\Longrightarrow\quad a\le b.
\]
Its group of differences $G(M)$ carries the order
\[
 [a,b]\le[c,d]\quad\Longleftrightarrow\quad a+d\le c+b.
\]
Because adding equal monoid elements to both sides leaves this relation
unchanged, and order cancellation removes the common summands, it
descends to an antisymmetric relation on the group of differences. The canonical map $j(a)=[a,0]$ is an order embedding. Write $D\downarrow_M0$ when $D\subseteq M$ is nonempty, downward directed, and has infimum zero in $M$.

\begin{theorem}[Translation of zero cuts and canonical groupification]
\label{sx:thm-transferZeroCuts}
The following are equivalent:
\begin{enumerate}
\item Every $D\downarrow_M0$ satisfies $j(D)\downarrow_{G(M)}0$.
\item For every $D\downarrow_M0$ and $a,b\in M$,
\begin{equation}\label{sx:zero-cut}
 a\le b+d\quad(d\in D)\quad\Longrightarrow\quad a\le b.
\end{equation}
\item For every $D\downarrow_M0$ and $b\in M$, the infimum of $b+D$ in $M$ is $b$.
\end{enumerate}
If, in addition, the order on $M$ is algebraic, meaning
\begin{equation}\label{sx:algebraic-order}
 a\le b\quad\Longleftrightarrow\quad b=a+c\text{ for some }c\in M,
\end{equation}
these conditions imply equivalence of null order convergence in $M$ and $G(M)$ for every net of elements of $M$.
\end{theorem}
\begin{proof}
Every element of $G(M)$ has the form $j(a)-j(b)$. It is a lower bound of $j(D)$ exactly when $a\le b+d$ for all $d\in D$, and it is nonpositive exactly when $a\le b$. Since zero is already a lower bound, this proves the equivalence of the first two assertions. For the third assertion, $b$ is a lower bound of $b+D$ since
$0\le d$ for $d\in D$. If (2) holds and $a$ is any lower bound,
then $a\le b+d$ for every $d\in D$ implies $a\le b$, so
$\inf_M(b+D)=b$. Conversely, that infimum and the inequalities
$a\le b+d$ imply $a\le b$, which is (2).

Under \eqref{sx:algebraic-order}, $G(M)_+=j(M)$: if $j(a)-j(b)\ge0$, then $b\le a$, and $a=b+c$ makes this difference equal to $j(c)$. Let $(x_\alpha)$ be a net in $M$.
If $D\downarrow_M0$ regulates it, assertion (1) makes $j(D)$ a
whole-group zero regulator, and
\[
 x_\alpha\le d\text{ eventually}
 \quad\Longleftrightarrow\quad
 j(x_\alpha)\le j(d)\text{ eventually}
 \qquad(d\in D)
\]
proves preservation. Conversely, let $R\subseteq G(M)_+$ be a
whole-group zero regulator of $j(x_\alpha)$. The equality
$G(M)_+=j(M)$ defines $D=j^{-1}(R)$, which is nonempty and downward
directed. If $a\in M$ is a lower bound of $D$, then $j(a)$ is a
lower bound of $R$, so $j(a)\le0$ and $a\le0$. Thus
$\inf_M D=0$. Order reflection gives eventual domination of
$x_\alpha$ by each $d\in D$, proving reflection for this net.
\end{proof}

For an $M$-valued faithful invariant metric satisfying these monoid
hypotheses, the passage to $G(M)$ therefore preserves the original pointwise,
uniform-error nets and group-Cauchy nets, providing a conversion of the
specified representation through the zero-infimum properties just
established. Its scope is determined by this prescribed change of values,
whereas an intrinsic existence question also permits another representation
of the same original convergence. When \eqref{sx:algebraic-order} is not
assumed, preservation remains valid and reflection requires separate
control of the new positive regulators.

\subsection{A faithful metric on the real line which cannot be transferred additively}

To see why preservation of zero cuts requires its own hypothesis,
consider the free commutative monoid
\[
 M=\mathbb{N}^{(\mathbb{N}_{\ge1})}
\]
with generators $e_n$. Compare elements first by their total degree $\ell(a)=\sum_n a_n$ and, for equal degrees, lexicographically at their first differing coordinate. Then $M$ is cancellative, the order is total and admits order cancellation, and
\begin{equation}\label{sx:free-regulator}
 e_1>e_2>\cdots>0,\qquad \inf_M\{e_n:n\ge1\}=0.
\end{equation}
Indeed, any nonzero lower bound would have degree at least one. Degree at least two makes it larger than every generator. Degree one makes it $e_j$ for some $j$, contradicted by $e_{j+1}<e_j$. This free-semigroup construction appears in \citet[p.~373]{Reichel1977distance}.

Define
\[
 p(0)=0,\qquad p(x)=e_{\lceil1/|x|\rceil}\quad(x\ne0),
 \qquad d(x,y)=p(x-y).
\]
This is a definite symmetric invariant monoid-valued metric. If $x,y$ are nonzero, $p(x)+p(y)$ has degree two and dominates $p(x+y)$, whose degree is at most one. If one summand is zero the triangle inequality is equality or immediate. For $n\ge2$ its algebraic balls satisfy
\begin{equation}\label{sx:free-balls}
 \{x:p(x)\le e_n\}=\left(-\frac1{n-1},\frac1{n-1}\right).
\end{equation}
It follows that ordinary null nets have the regulator in \eqref{sx:free-regulator}. Conversely, any zero regulator in the totally ordered $M$ has a member below each $e_n$. Its ball conditions force every interval in \eqref{sx:free-balls} to be an eventual set. Therefore $d$ is faithful for all nets of the ordinary real line. Moreover $p(1/n)=e_n$, so $M$ is exactly the monoid generated by its metric values.

\begin{proposition}[No faithful additive repair of this representation]
If $h:M\to H$ is a monotone monoid homomorphism into a partially ordered abelian group and $h(e_n)\to_o0$, then $h=0$.
\end{proposition}
\begin{proof}
For all $i,j,n$ the degree comparison gives $e_i\le e_j+e_n$, whence
\[
 h(e_i)-h(e_j)\le h(e_n).
\]
Any fixed lower bound of an order-null positive net is nonpositive: if $R\downarrow0$ regulates the net, eventual domination by each $r\in R$ gives domination of the fixed bound by $r$, and then by zero. Hence $h(e_i)\le h(e_j)$ for all $i,j$. All generator images are equal to a common element $c$. Their order nullity forces $c=0$, proving the assertion.
\end{proof}

Thus even for the usual convergence of the real line, completeness,
regularity and first countability do not ensure faithful additive
groupification of a specified monoid metric, although its monoid is
cancellative and generated by its distances. Since $|x-y|$ already gives
a faithful real metric for this convergence, the obstruction concerns the
conversion under consideration and explains why existence of another
representation must be treated separately. The balls in
\eqref{sx:free-balls} also connect the example with the earlier closed-ball
discussion: although they are algebraic closed-radius balls in $M$, they
are open intervals in the original convergence, so their closedness cannot
be inferred from monoid faithfulness.

\subsection{Order-dense completion and reflection of regulators}
\label{sx:completion}

A passage to an order-dense completion preserves existing infima, but
it may also supply regulators that were unavailable in the original
value group. The following construction isolates this second effect
for a faithful homogeneous gauge.

\begin{theorem}[Completion does not automatically preserve faithful representation]\label{sx:thm-completionloss}
There is an Archimedean directed partially ordered real vector space $G$, order dense in its Dedekind completion $L$, and a definite homogeneous symmetric subadditive gauge $p:c_0\to G_+$ faithful for the norm convergence of the Banach space $c_0$, such that $j\circ p:c_0\to L_+$ is not faithful, where $j:G\hookrightarrow L$ is the order-dense embedding.
\end{theorem}
\begin{proof}
Let $L=\ell^\infty(\mathbb{Z})$ with coordinate order and define
\[
G=\left\{g\in L:\lim_{n\to+\infty}g_n\text{ exists and }
\lim_{n\to+\infty}g_n=\sum_{k=1}^\infty2^{-k}g_{-k}\right\}.
\]
The weighted series is absolutely convergent. Thus $G$ is a real vector subspace, with inherited Archimedean order. It contains the constant sequence $\mathbf{1}$, an order unit, so it is directed.

We verify order density explicitly. Let $y\in L$ and $M\ge\max(0,\sup_n y_n)$. For each fixed coordinate $j$, there is $g\in G$ with $g\ge y$ and $g_j=y_j$. If $j\ge0$, take the constant sequence $M$ and change only coordinate $j$ to $y_j$. Neither its positive-tail limit nor its negative weighted mean changes. If $j=-k<0$, first change the $-k$ coordinate to $y_{-k}$. At a second negative coordinate $-l$, $l\ne k$, change the value to
\[
M+2^{l-k}(M-y_{-k}).
\]
The weighted mean remains $M$, the positive tail remains constant $M$, and every coordinate is still at least that of $y$. Consequently
\[
\inf_L\{g\in G:g\ge y\}=y.
\]
Negation gives the corresponding supremum of lower approximants. Since $G$ is also majorizing and $L$ is Dedekind complete, this is its Dedekind completion.

For $x\in c_0(\mathbb{N})$ define
\begin{equation}\label{sx:completion-gauge}
p(x)_n=\begin{cases}|x_n|,&n\ge1,\\0,&n\le0.\end{cases}
\end{equation}
The positive-tail limit and the negative weighted mean are both zero, so $p(x)\in G_+$. Definiteness, homogeneity, symmetry, and subadditivity are coordinatewise.

If $x_\alpha\to0$ in norm, then $p(x_\alpha)\le\varepsilon\mathbf{1}$ eventually for every $\varepsilon>0$. The elements $\varepsilon\mathbf{1}$ decrease to zero in $G$, proving one implication of faithfulness. Conversely, suppose $p(x_\alpha)\to_o0$ in $G$, regulated by $r_\delta\downarrow0$. To verify that its infimum remains zero in $L$, let $v\in L$
be a lower bound of all $r_\delta$. Every $g\in G$ with $g\le v$
is then a lower bound of the regulator in $G$, so $g\le0$.
The supremum in $L$ of these $G$-minorants is $v$ by order density,
whence $v\le0$. Thus $r_\delta\downarrow0$ in $L$.
Its coordinate infima are zero because the order of
$L=\ell^\infty(\mathbb Z)$ is coordinatewise. Fix $\delta_0$ and restrict to the cofinal tail below $r_{\delta_0}$. With $M=\|r_{\delta_0}\|_\infty$, all these nonnegative coordinates are bounded by $M$. Given $\eta>0$, take $K$ with $M\sum_{k>K}2^{-k}<\eta/2$. The finitely many remaining coordinates can be made small enough that $\sum_{k\le K}2^{-k}r_\delta(-k)<\eta/2$. Hence
\[
\ell_\delta:=\lim_{n\to\infty}r_\delta(n)
=\sum_{k\ge1}2^{-k}r_\delta(-k)\longrightarrow0.
\]
For $\varepsilon>0$, fix $\delta$ with $\ell_\delta<\varepsilon/2$ and then $N$ so that $r_\delta(n)<\varepsilon$ for $n>N$. A later index $\gamma\ge\delta$ makes the finitely many coordinates $1\le n\le N$ smaller than $\varepsilon$. Monotonicity retains the tail bounds, so $\sup_{n\ge1}r_\gamma(n)\le\varepsilon$. Since $p(x_\alpha)\le r_\gamma$ eventually, $\|x_\alpha\|_\infty\le\varepsilon$ eventually. This proves norm faithfulness for every net.

Let $e_n$ be the unit vectors of $c_0$. Define $b_N\in L$ by $b_N(j)=1$ for $j\ge N$ and $b_N(j)=0$ otherwise. Then $b_N\downarrow0$ in $L$ and $j(p(e_n))\le b_N$ for $n\ge N$. Thus $j(p(e_n))\to_o0$ in $L$, while $\|e_n\|_\infty=1$. The composed gauge is not faithful.
\end{proof}

Here the obstruction belongs to the conversion of a
\emph{specified} faithful representation: the original convergence
continues to admit the faithful lattice-valued metric on $c_0$ supplied
by its usual real norm. Within the constructed representation, we can
also separate finite-ball conjunction from a property of the original
convergence. Define $r,s\in G_+$ by
\[
r_n=s_n=1\ (n\ge1),\qquad r_{-1}=2,\quad s_{-2}=4,
\]
with every other nonpositive coordinate zero. Both negative weighted means and both positive-tail limits equal one. For \eqref{sx:completion-gauge}, both $B_r$ and $B_s$ are the closed norm unit ball of $c_0$. If $t\in G_+$ satisfies $t\le r,s$, every negative coordinate of $t$ is zero, so its positive-tail limit is zero. If in addition $B_r\cap B_s\subseteq B_t$, testing each unit vector gives $t_n\ge1$ for every $n\ge1$, a contradiction. Thus the persistent-bound set of the principal unit-ball saturation net is not downward directed, although the original Banach convergence satisfies property~\ref{nx:P}: for $\varepsilon>0$, a null tail family supplies $A\subseteq(\varepsilon/2)U$, the corresponding null tail family supplies $B\subseteq(\varepsilon/2)U$, and $A+B\in\mathcal K$ implies $\varepsilon U\in\mathcal K$. Even identical balls can have labels with no common lower label retaining that ball. The labeling data cannot silently be discarded.

\subsection*{Exactly what an order-dense lattice embedding changes}

To identify precisely what has changed in these examples, we compare
persistent bounds before and after the embedding. This will give a
criterion for retention of faithfulness under the prescribed conversion,
while leaving open the possibility of representing the original
convergence by a different gauge.

\begin{proposition}[Persistent bounds and order-dense completion]\label{sx:prop-completioncriterion}
Let $j:G\to L$ be an additive order embedding from an abelian partially ordered group into an abelian lattice-ordered group. Assume every element of $L$ is both the infimum of its $j(G)$-majorants and the supremum of its $j(G)$-minorants. For a positive net $(u_\alpha)$ in $G$, put
\[
R(u)=\{r\in G_+:u_\alpha\le r\text{ eventually}\}.
\]
Then
\begin{equation}\label{sx:persistent-lattice}
j(u_\alpha)\to_o0\text{ in }L
\quad\Longleftrightarrow\quad
R(u)\ne\varnothing\text{ and }\inf_G R(u)=0,
\end{equation}
whereas
\begin{equation}\label{sx:persistent-group}
u_\alpha\to_o0\text{ in }G
\quad\Longleftrightarrow\quad
R(u)\text{ contains a downward-directed }D\text{ with }\inf_GD=0.
\end{equation}
\end{proposition}
\begin{proof}
Suppose the image net is regulated by $d_\delta\downarrow0$ in $L$. Every $G$-majorant of a given $d_\delta$ belongs to $R(u)$, so this set is nonempty. If $h\in G$ is a lower bound of $R(u)$, then $j(h)$ is below all those majorants and therefore below their infimum $d_\delta$. This holds for every $\delta$, so $h\le0$. Zero is a lower bound, proving the right side of \eqref{sx:persistent-lattice}.

Conversely, order density preserves the existing infimum zero of $R(u)$. Indeed, if $v\in L$ is below every $j(r)$, all its $j(G)$-minorants are below zero. Their supremum is $v$, so $v\le0$. Hence $\inf_Lj(R(u))=0$. The finite meets
\[
d_F=\bigwedge_{r\in F}j(r),\qquad \varnothing\ne F\subseteq R(u)\text{ finite},
\]
form a downward-directed zero regulator. Every finite collection of eventual bounds holds on a common tail, so $j(u_\alpha)\le d_F$ eventually. This proves \eqref{sx:persistent-lattice}. Equation \eqref{sx:persistent-group} is exactly the regulator criterion in the original group, stated using its persistent-bound set.
\end{proof}

An explicit finite-cut criterion expresses when this conversion reflects all null nets. Fix $S\subseteq G_+$ containing zero. For a nonempty $R\subseteq G_+$ with $\inf_G R=0$, put
\begin{equation}\label{sx:finite-cut}
R^{\#}_S=\left\{g\in G_+:\exists\varnothing\ne F\subseteq R\text{ finite},\quad
S\cap\bigcap_{r\in F}[0,r]\subseteq[0,g]\right\}.
\end{equation}
\begin{proposition}[Exact finite-cut criterion]\label{sx:completion-reflection}
The embedding $j$ reflects order nullity of every tail family of a net on $S$
if and only if, for every nonempty $R\subseteq G_+$ with $\inf_G R=0$,
the set $R^{\#}_S$ contains a downward-directed subset having infimum
zero in $G$.
\end{proposition}
\begin{proof}
For necessity, define
\[
 C_F=S\cap\bigcap_{r\in F}[0,r],\qquad
 \varnothing\ne F\subseteq R\text{ finite}.
\]
These sets contain zero and satisfy $C_F\cap C_{F'}=C_{F\cup F'}$.
Their eventual-set family has image regulated in $L$ by
\[
 d_F=\bigwedge_{r\in F}j(r),\qquad d_F\downarrow0,
\]
where preservation of the zero infimum is proved in
Proposition~\ref{sx:prop-completioncriterion}. Reflection supplies
$D\subseteq G_+$ with $D\downarrow0$ and, for every $d\in D$,
a finite $F$ with $C_F\subseteq[0,d]$.
Hence $D\subseteq R^{\#}_S$.

Conversely, let $\mathcal H$ be a tail family of a net on $S$ whose image is
order null in $L$. Its persistent-bound set
\[
 R=\{r\in G_+:S\cap[0,r]\in\mathcal H\}
\]
obeys $R\ne\varnothing$ and $\inf_G R=0$ by
Proposition~\ref{sx:prop-completioncriterion} applied to an actual net
$(u_\alpha)$ on $S$ with $\mathcal H=\operatorname{Tail}(u)$:
indeed, $S\cap[0,r]\in\mathcal H$ is exactly the assertion
$(\exists\alpha_r)(\forall\alpha\ge\alpha_r)\ u_\alpha\le r$. Every $r\in R$ lies in $R^{\#}_S$
by the singleton finite subset $\{r\}$. If $g\in R^{\#}_S$, its
finite witnessing set $F$ satisfies
\[
 C_F=\bigcap_{r\in F}(S\cap[0,r])\in\mathcal H,\qquad
 C_F\subseteq S\cap[0,g].
\]
Thus $g\in R$ and $R^{\#}_S=R$. The assumed directed subset
$D\subseteq R^{\#}_S=R$ is a zero regulator for $\mathcal H$ in $G$.
\end{proof}

Taking $S=G_+$ gives reflection of all positive order-null nets by the embedding. Taking $S=p(X)$ gives the exact criterion for the composed gauge $j\circ p$ to remain faithful, provided $p$ was faithful to begin with: nets on $S$ lift through the surjection $p:X\to S$. Formula \eqref{sx:finite-cut} thus describes the specified change of values
through finite cuts on the value range, retaining the distinction
from property~\ref{nx:P}, whose definition concerns the original
convergence.

\subsection{Ordered extensions preserving designated zero infima}
\label{sx:cut-extensions}

The finite-cut criterion leads in turn to a question about extending an
ordered group: when can one adjoin a new lower bound while preserving
the zero infima already designated as regulators? Throughout this
subsection, an ordered abelian group means an abelian group with a
translation-invariant partial order. The order need not be a lattice order, directed, Archimedean, or unperforated. A downward-directed subset is nonempty. All infima mentioned below are infima in the whole group.

Let $H$ be an ordered abelian group, let $J$ be a set, and, for each $j\in J$, let $E_j\subseteq H_+$ be downward directed with
\[
\inf_H E_j=0.
\]
These are the designated regulators whose zero infima must survive an extension. Their directedness is relevant to their interpretation as regulators. The proofs below only use their nonemptiness and their stated infima. Write $nD=\{nd:d\in D\}$, where $n\geq1$ and $D\subseteq H$. Thus $nD$ denotes scalar multiplication of a set, not an independently indexed Minkowski sum.

\paragraph{An adjoined lower bound not below zero.}

For a prescribed set $D\subseteq H$, let $(I_n)_{n\geq1}$ be the least family of subsets of $H$ satisfying the following rules:
\begin{align}
 nD&\subseteq I_n, \label{sx:cut-S1}\\
 I_n+H_+&\subseteq I_n, \label{sx:cut-S2}\\
 I_n+I_m&\subseteq I_{n+m}\qquad(n,m\geq1), \label{sx:cut-S3}\\
 a+E_j\subseteq I_n&\Longrightarrow a\in I_n
 \qquad(a\in H,\ j\in J,\ n\geq1). \label{sx:cut-S4}
\end{align}
The least family exists: intersect all families satisfying these rules. There is at least one such family, namely $I_n=H$ for every $n$. Notice that \eqref{sx:cut-S4} is an infinitary rule when $E_j$ is infinite. The closure rule retains the whole family $E_j$.

\begin{theorem}[Exact cut-adjunction criterion]\label{sx:cut-thm-witness-arbitrary}
The following assertions are equivalent.
\begin{enumerate}
\item There exist an ordered abelian group $K$, an order embedding $\iota:H\longrightarrow K$, and $w\in K$ such that
\[
 w\leq\iota(d)\quad(d\in D),\qquad w\not\leq0,
 \qquad \inf_K\iota(E_j)=0\quad(j\in J).
\]
\item $0\notin I_1$.
\end{enumerate}
When the second assertion holds, one may take the underlying group of $K$ to be $H\oplus\mathbb{Z}$, with $\iota(a)=(a,0)$ and $w=(0,1)$.
\end{theorem}

\begin{proof}
Assume first that the extension in the first assertion exists. Identify $H$ with its image under $\iota$ and put
\[
 J_n=\{a\in H:nw\leq a\}\qquad(n\geq1).
\]
For $d\in D$, the inequality $w\leq d$ implies $nw\leq nd$, so $nD\subseteq J_n$. The sets $J_n$ are upward closed, and addition gives $J_n+J_m\subseteq J_{n+m}$. Finally, if $a+E_j\subseteq J_n$, then
\[
 nw-a\leq e\qquad(e\in E_j).
\]
The preserved infimum in $K$ yields $nw-a\leq0$, and hence $a\in J_n$. Thus $(J_n)$ satisfies \eqref{sx:cut-S1}--\eqref{sx:cut-S4}, and minimality gives $I_n\subseteq J_n$. If $0\in I_1$, then $w\leq0$, a contradiction.

Conversely, assume $0\notin I_1$. On the abelian group $K=H\oplus\mathbb{Z}$, define
\[
 P=\{(a,0):a\in H_+\}
 \ \cup\ 
 \bigcup_{n\geq1}\{(a,-n):a\in I_n\}.
\]
The element $(0,0)$ belongs to $P$. Its addition rules are
\[
 (H_+\times\{0\})+(H_+\times\{0\})\subseteq H_+\times\{0\},
\]
\[
 (I_n\times\{-n\})+(H_+\times\{0\})
       \subseteq I_n\times\{-n\},
\]
\[
 (I_n\times\{-n\})+(I_m\times\{-m\})
       \subseteq I_{n+m}\times\{-(n+m)\}.
\]
These are respectively additivity of $H_+$ and
\eqref{sx:cut-S2}--\eqref{sx:cut-S3}. Thus $P+P\subseteq P$.

Moreover, $P\cap(-P)=\{(0,0)\}$. Indeed, all elements of $P$ have nonpositive second coordinate, so an element of $P\cap(-P)$ has second coordinate zero. Its first coordinate then belongs to $H_+\cap(-H_+)=\{0\}$. Consequently $P$ defines a translation-invariant partial order on $K$.

The map $\iota(a)=(a,0)$ is an order embedding, because
\[
 (b-a,0)\in P\quad\Longleftrightarrow\quad b-a\in H_+.
\]
For $d\in D$, rule \eqref{sx:cut-S1} gives $d\in I_1$, whence
\[
 \iota(d)-w=(d,-1)\in P.
\]
Also, $w\leq0$ would mean $(0,-1)\in P$, equivalently $0\in I_1$, which is excluded.

It remains to check all designated infima. Zero is a lower bound of $\iota(E_j)$ because $E_j\subseteq H_+$. Suppose that $(h,k)$ is any lower bound of $\iota(E_j)$. For every $e\in E_j$,
\[
 (e-h,-k)\in P.
\]
Since $P$ has no element with positive second coordinate and $E_j$ is nonempty, this forces $k\geq0$. If $k=0$, then $h\leq e$ for every $e\in E_j$, so $h\leq0$ in $H$, and $(h,0)\leq0$ in $K$. If $k\geq1$, the displayed membership says
\[
 -h+E_j\subseteq I_k.
\]
Rule \eqref{sx:cut-S4} gives $-h\in I_k$, so $(-h,-k)\in P$, again proving $(h,k)\leq0$. Therefore $\inf_K\iota(E_j)=0$ for every $j$.
\end{proof}

The new element in Theorem~\ref{sx:cut-thm-witness-arbitrary} is deliberately required only to satisfy $w\not\leq0$. In the explicit construction it is incomparable with zero: $w\notin P$ because its second coordinate is positive, while $-w\notin P$ because $0\notin I_1$. A positive lower bound is a stronger requirement, addressed separately below.

\begin{remark}[The first cut has an independent description]
In Theorem~\ref{sx:cut-thm-witness-arbitrary}, $I_1$ is simply the least upper subset of $H$ containing $D$ and satisfying
\[
 a+E_j\subseteq C\Longrightarrow a\in C.
\]
Indeed, the rules for $I_1$ give one inclusion. Conversely, if $C$ has the stated properties, then the family $I_1'=C$, $I_n'=H$ for $n\geq2$ satisfies \eqref{sx:cut-S1}--\eqref{sx:cut-S4}, giving the reverse inclusion by minimality. This observation does not assert that the resulting ordered group is unperforated: the construction permits $2w\leq0$ while $w\not\leq0$. That is allowed in an arbitrary partially ordered group, but is impossible in an ordered real vector space. The vector-space variant below imposes the corresponding convexity explicitly.
\end{remark}

\paragraph{A strictly positive adjoined lower bound.}

Now suppose $D\subseteq H_+$. Define $(I_n^+)_{n\geq1}$ to be the least family of subsets of $H_+$ satisfying \eqref{sx:cut-S1}--\eqref{sx:cut-S4}, with $I_n^+$ in place of $I_n$, and the additional rule
\begin{equation}
 I_{n+1}^+\subseteq I_n^+\qquad(n\geq1). \label{sx:cut-S5}
\end{equation}
In \eqref{sx:cut-S4}, the translating element $a$ still ranges over all of $H$, not merely $H_+$. This definition is consistent: the family $I_n^+=H_+$ satisfies every rule. In fact, if $a+E_j\subseteq H_+$, then $-a\leq e$ for all $e\in E_j$, and $\inf_HE_j=0$ implies $a\geq0$.

\begin{theorem}[Exact positive cut-adjunction criterion]\label{sx:cut-thm-witness-positive}
The following assertions are equivalent.
\begin{enumerate}
\item There exist an ordered abelian group $K$, an order embedding $\iota:H\longrightarrow K$, and an element $w\in K$ such that
\[
 0<w\leq\iota(d)\quad(d\in D),
 \qquad\inf_K\iota(E_j)=0\quad(j\in J).
\]
\item $0\notin I_1^+$.
\end{enumerate}
Again the underlying group may be taken to be $H\oplus\mathbb{Z}$.
\end{theorem}

\begin{proof}
For necessity, identify $H$ with its image and define
\[
 J_n^+=\{a\in H:nw\leq a\}.
\]
Since $w>0$, every member of $J_n^+$ lies in $H_+$. For $d\in D$, $nw\le nd$, giving $nD\subseteq J_n^+$.
If $a\in J_n^+$ and $h\in H_+$, then $nw\le a\le a+h$.
If $a\in J_n^+$ and $b\in J_m^+$, then $(n+m)w\le a+b$.
Finally,
\[
 a+E_j\subseteq J_n^+
 \Longrightarrow nw-a\le e\quad(e\in E_j)
 \Longrightarrow nw-a\le0
 \Longrightarrow a\in J_n^+.
\]
Moreover, $(n+1)w\leq a$ implies $nw\leq a$, so $J_{n+1}^+\subseteq J_n^+$. Thus $I_n^+\subseteq J_n^+$, and $0\notin I_1^+$.

For sufficiency, define the cone on $H\oplus\mathbb{Z}$ by
\[
 P^+=\{(a,k):k\geq0,\ a\in H_+\}
 \ \cup\
 \bigcup_{n\geq1}\{(a,-n):a\in I_n^+\}.
\]
For two elements with nonnegative second coordinates, both their
first-coordinate sum and their second-coordinate sum are nonnegative.
For two negative-coordinate elements,
\[
 a\in I_n^+,\ b\in I_m^+
 \quad\Longrightarrow\quad
 (a,-n)+(b,-m)\in I_{n+m}^+\times\{-(n+m)\}.
\]
The remaining case is a sum $(a,-n)+(b,k)$ with $n\geq1$, $k\geq0$, $a\in I_n^+$, and $b\in H_+$. If $k\geq n$, both the resulting second coordinate and the first coordinate $a+b$ are nonnegative. If $k<n$, repeated use of \eqref{sx:cut-S5} gives $a\in I_{n-k}^+$, and \eqref{sx:cut-S2} gives $a+b\in I_{n-k}^+$. Hence $P^+$ is closed under addition.

To prove that the cone is proper, suppose both $(a,k)$ and $(-a,-k)$ belong to $P^+$. For $k=0$, properness of the original cone gives $a=0$. For $k>0$, one has $a\in H_+$ and $-a\in I_k^+\subseteq H_+$, so $a=0$ and $0\in I_k^+\subseteq I_1^+$, a contradiction. The case $k<0$ is symmetric. Thus $P^+$ defines a partial order.

The original group embeds with its original order, because the zero-coordinate slice of $P^+$ is $H_+\times\{0\}$. The element $w=(0,1)$ is positive and nonzero, and $(d,-1)\in P^+$ for each $d\in D$, proving $w\leq\iota(d)$.

For preservation of a designated infimum, let $(h,k)$ be a lower bound of $\iota(E_j)$. If $k\leq0$, then
\[
 (e-h,-k)\in P^+\quad(e\in E_j)
\]
implies $e-h\geq0$ in $H$. Therefore $h\leq0$, and $(-h,-k)\in P^+$, so $(h,k)\leq0$. If $k>0$, the same inequalities say $-h+E_j\subseteq I_k^+$. Rule \eqref{sx:cut-S4} gives $-h\in I_k^+$, again yielding $(h,k)\leq0$. This proves all required infima.
\end{proof}

\paragraph{The role of saturation.}

Take $H=\mathbb{R}$ with its usual order,
\[
 E=\{1/m:m\geq1\},\qquad D=\{1+1/m:m\geq1\}.
\]
For both adjunction criteria the saturated families are exactly
\[
 I_n=I_n^+=[n,\infty)\qquad(n\geq1).
\]
Indeed, the initial upper closures of $nD$ are $(n,\infty)$, and \eqref{sx:cut-S4} adds the endpoint $n$ because $n+E\subseteq(n,\infty)$. Conversely, the family $[n,\infty)$ satisfies every required rule, including \eqref{sx:cut-S5} in the positive case. Thus $0\notin I_1$, and a lower bound can be adjoined while preserving $\inf E=0$. Necessarily the resulting lower bound satisfies $w\leq1$, even though $1\notin D$.

Without the saturation rule, the cone generated only by the initial inequalities $w\leq1+1/m$ would have the extra lower bound
\[
 w-1\leq1/m\qquad(m\geq1),
\]
while failing to have $w-1\leq0$. It would therefore destroy the old zero infimum of $E$. This example applies to both versions of the adjunction.

If instead $D=E$, saturation gives $0\in I_1$ and $0\in I_1^+$. No lower bound of either kind can lie below every member of $E$ while preserving its zero infimum.

\paragraph{Incomparable lower bounds.}

In a partially ordered group, $0$ can be a lower bound of a downward-directed positive family without being its greatest lower bound, even when zero is its only positive lower bound. For an explicit example, put
\[
 H=\mathbb{R}^2,\qquad
 H_+=\{(0,0)\}\cup\{(s,t):s>0,\ t>0\},
\]
and
\[
 d_m=(1/m,1+1/m),\qquad D=\{d_m:m\geq1\}.
\]
This is an ordered abelian group: the displayed cone is additive and proper. The family $D$ is downward directed, because $d_m-d_n$ has both coordinates strictly positive whenever $n>m$. The element $v=(0,1)$ satisfies $v\leq d_m$ for every $m$, but $v\not\leq0$ and $v\not\geq0$. Thus $\inf_HD$ is not zero. Nevertheless, no strictly positive element is a lower bound of $D$: if $u=(s,t)>0$, then $s>0$, and an integer $m>1/s$ rules out $u\leq d_m$.

Consequently, tests for failure of a whole-group zero infimum must allow incomparable lower bounds. Requiring a strictly positive witness uses an additional order property, such as the availability of positive parts in a lattice-ordered group.

\paragraph{Ordered vector extensions.}

Suppose now that $H$ is an ordered real vector space, with a proper convex positive cone. Scalar multiplication of an ordered real vector space by a strictly positive scalar is an order isomorphism. In particular, $tE_j$ is again a zero regulator for every $t>0$.

Given $D\subseteq H$, let $I$ be the least convex upper subset of $H$ containing $D$ and satisfying
\begin{equation}
 a+tE_j\subseteq I\Longrightarrow a\in I
 \qquad(a\in H,\ t>0,\ j\in J). \label{sx:cut-VS}
\end{equation}
Such a least set exists by intersection, since the whole space is an admissible set. Empty sets are allowed when $D$ is empty.

\begin{theorem}[Exact vector cut-adjunction criterion]\label{sx:cut-thm-witness-vector}
There exist an ordered real vector space $K$, a linear order embedding $\iota:H\longrightarrow K$, and $w\in K$ such that
\[
 w\leq\iota(d)\quad(d\in D),\qquad w\not\leq0,
 \qquad\inf_K\iota(E_j)=0\quad(j\in J)
\]
if and only if $0\notin I$.
\end{theorem}

\begin{proof}
For necessity, identify $H$ with its image and set $C=\{a\in H:w\leq a\}$. This set is upper and convex and contains $D$. If $a+tE_j\subseteq C$, then $(w-a)/t\leq e$ for every $e\in E_j$. The preserved infimum gives $(w-a)/t\leq0$, so $a\in C$. Thus $I\subseteq C$, and $0\notin I$.

For sufficiency, put $K=H\oplus\mathbb{R}$ and define
\[
 P=\{(a,0):a\in H_+\}
 \ \cup\ 
 \{(a,-t):t>0,\ a\in tI\}.
\]
For $\lambda>0$, $a\in tI$ implies
$\lambda a\in(\lambda t)I$, and multiplication by zero gives $(0,0)$.
Thus $P$ is stable under nonnegative scalars. Its zero-coordinate slice
is additive. If $a=tc\in tI$ and $h\in H_+$, then
\[
 (a+h)/t=c+h/t\in I,
\]
so a zero-coordinate and a negative-coordinate term have sum in $P$.
For $a\in tI$ and $b\in sI$ with $t,s>0$,
\[
 \frac{a+b}{t+s}
 =\frac{t}{t+s}\frac a t+\frac{s}{t+s}\frac b s\in I
\]
by convexity. The cone is proper because all its second coordinates are nonpositive and its zero-coordinate slice is the proper cone of $H$. It therefore makes $K$ an ordered real vector space and makes $\iota(a)=(a,0)$ a linear order embedding. For $w=(0,1)$, the inclusion $D\subseteq I$ gives $w\leq\iota(d)$ for every $d\in D$, whereas $0\notin I$ gives $w\not\leq0$.

Let $(h,s)$ be a lower bound of $\iota(E_j)$. For every $e\in E_j$, $(e-h,-s)\in P$.
The cone has no elements with positive second coordinate, so $s\ge0$.
If $s=0$, these inequalities give $h\le e$ for every $e\in E_j$,
hence $h\le0$ and $(h,0)\le0$. If $s>0$, the inequalities say
\[
 \frac{-h}{s}+\frac1s E_j\subseteq I.
\]
Rule \eqref{sx:cut-VS} gives $-h/s\in I$, whence $(-h,-s)\in P$ and $(h,s)\leq0$. Thus every designated zero infimum is preserved.
\end{proof}

There is also an exact positive version. Suppose $D\subseteq H_+$ and let $I^+$ be the least convex upper subset of $H_+$ containing $D$ and satisfying \eqref{sx:cut-VS}, with $a$ still ranging over all of $H$. The family defining this least set is nonempty because $H_+$ itself satisfies \eqref{sx:cut-VS}.

\begin{corollary}[Positive vector lower bound]\label{sx:cut-cor-witness-vector-positive}
A linear ordered extension preserving the designated zero infima and containing $0<w\leq\iota(d)$ for every $d\in D$ exists if and only if $0\notin I^+$.
\end{corollary}

\begin{proof}
Suppose first that the extension exists and identify $H$ with its image.
Set $C=\{a\in H:w\le a\}$. Then $C\subseteq H_+$ and $D\subseteq C$.
For $a,b\in C$, $0\le\theta\le1$ and $h\in H_+$,
\[
 w=\theta w+(1-\theta)w\le\theta a+(1-\theta)b,
 \qquad w\le a+h.
\]
Thus $C$ is convex and upper. If $a+tE_j\subseteq C$, then
$(w-a)/t\le e$ for every $e\in E_j$. Preservation of its infimum gives
$(w-a)/t\le0$, so $a\in C$. Minimality implies $I^+\subseteq C$.
Since $w>0$, $0\notin C$, proving $0\notin I^+$.

Conversely, suppose $0\notin I^+$ and define on $K=H\oplus\mathbb R$
\[
 P^+=\{(a,s):s\ge0,\ a\in H_+\}
 \ \cup\
 \{(a,-t):t>0,\ a\in tI^+\}.
\]
It contains zero and is stable under multiplication by nonnegative
scalars. Sums of two nonnegative-coordinate elements remain in the
first set. For $a\in tI^+$ and $b\in sI^+$, $t,s>0$, convexity gives
\[
 \frac{a+b}{t+s}
   =\frac{t}{t+s}\frac at+\frac{s}{t+s}\frac bs\in I^+,
\]
so two negative-coordinate elements also have sum in $P^+$.
For a mixed sum $(a,-t)+(b,s)$ with $t>0$, $s\ge0$, $a=tc$,
$c\in I^+\subseteq H_+$ and $b\in H_+$, the case $s\ge t$ gives
$a+b\in H_+$ and nonnegative second coordinate. If $s<t$, then
\[
 \frac{a+b}{t-s}
 =c+\frac{s}{t-s}c+\frac{b}{t-s}\in I^+,
\]
by upper closure. Thus $P^+$ is additive.
If $(a,s),(-a,-s)\in P^+$, their first coordinates imply
$a,-a\in H_+$, hence $a=0$. If $s\ne0$, one of these elements has
negative second coordinate, which would give $0\in I^+$.
Consequently $P^+\cap(-P^+)=\{0\}$, and $P^+$ is a proper convex cone.
Its zero-coordinate slice is $H_+\times\{0\}$, so the canonical
inclusion is a linear order embedding. The vector $w=(0,1)$ is
strictly positive, and $d\in I^+$ gives $(d,-1)\in P^+$ for $d\in D$.

Let $(h,s)$ be a lower bound of $\iota(E_j)$. For $s\le0$, the
conditions $(e-h,-s)\in P^+$ give $h\le e$ for every $e\in E_j$.
Thus $h\le0$ and $(-h,-s)\in P^+$. For $s>0$, these conditions give
\[
 -h/s+(1/s)E_j\subseteq I^+.
\]
Rule~\eqref{sx:cut-VS} yields $-h/s\in I^+$ and
$(-h,-s)\in P^+$. Every lower bound is therefore at most zero,
proving $\inf_K\iota(E_j)=0$.
\end{proof}

Having characterized these extensions, we can relate them to
Theorem~\ref{nx:canonical-po} only by retaining the distinction between
preserving specified zero infima in an existing representation and
obtaining a representation of the original convergence. In applications
to persistent upper bounds, directedness must also remain explicit:
such a set may have infimum zero without containing a downward-directed
subset with infimum zero, so the cut-adjunction tests do not dispense
with the regulator directedness required by the application.

\subsection{The Feldman--Porter construction and ordinary order convergence}
\label{fp:section}

A further distinction appears in the representation theorem of Feldman
and Porter \cite[pp.~603--604]{FeldmanPorter1973}, whose unbounded local
order convergence permits lower bounds that are quasi-negative relative
to a specified ideal without requiring them to be at most zero in the
whole ordered vector space. To follow the effect of this distinction,
we examine one Hausdorff locally convex stage and an explicit quotient,
determining their ordinary order convergence before comparing it with
the convergence used in that theorem.

\subsubsection{The cone for one locally convex stage}

Let $E$ be a real Hausdorff locally convex topological vector space,
and let $\mathcal P$ be the set of all continuous seminorms on $E$,
including the zero seminorm. Put
\[
 H=E\oplus\prod_{p\in\mathcal P}(E\times\mathbb R),\qquad
 a=\bigl(a_0,(a_p,t_p)_{p\in\mathcal P}\bigr)\in H.
\]
Define
\begin{equation}\label{fp:cone}
 P_H=\left\{a\in H:
 \begin{array}{ll}
 t_p\ge0,&p\in\mathcal P,\\
 a_p\ne0\Longrightarrow t_p>0,&p\in\mathcal P,\\
 p(a_0+a_p)\le t_p,&p\in\mathcal P
 \end{array}\right\},
 \qquad
 i(x)=\bigl(x,(0,0)_{p\in\mathcal P}\bigr).
\end{equation}
In this Hausdorff one-stage case, the nonzero-auxiliary-block
condition in the original cone is automatic: if all $a_p,t_p$
vanish and $a\in P_H$, then $p(a_0)=0$ for every continuous
seminorm, whence $a_0=0$.

\begin{lemma}[The one-stage order]\label{fp:cone-order}
The set $P_H$ is a pointed convex cone. Therefore
$a\le_H b\Longleftrightarrow b-a\in P_H$ defines a partial order
compatible with the real vector space structure of $H$.
\end{lemma}
\begin{proof}
The zero vector belongs to $P_H$. If $a\in P_H$ and $s>0$, then
\[
 st_p\ge0,\qquad sa_p\ne0\Longrightarrow st_p>0,\qquad
 p(sa_0+sa_p)=s\,p(a_0+a_p)\le st_p.
\]
Thus $sP_H\subseteq P_H$. For $s=0$ the image is $\{0\}$.
For $a,b\in P_H$, write the scalar coordinates of $b$ as $u_p$.
For every $p$,
\[
 t_p+u_p\ge0,\qquad
 p(a_0+b_0+a_p+b_p)
 \le p(a_0+a_p)+p(b_0+b_p)\le t_p+u_p.
\]
If $a_p+b_p\ne0$, then $a_p\ne0$ or $b_p\ne0$. In the first
case $t_p>0$, while in the second case $u_p>0$. Hence $t_p+u_p>0$.
This proves $P_H+P_H\subseteq P_H$, and scalar closure together
with addition proves convexity.

If $a,-a\in P_H$, then $t_p\ge0$ and $-t_p\ge0$, so $t_p=0$
for every $p$. The strict auxiliary condition gives $a_p=0$.
The remaining inequalities give $p(a_0)=0$ for every $p$.
The continuous seminorms separate points: a nonzero vector excluded
by a zero neighborhood is excluded by a finite-seminorm basic
neighborhood contained in it. Thus $a_0=0$ and
$P_H\cap(-P_H)=\{0\}$. The cone identities now give reflexivity,
transitivity, antisymmetry, translation invariance, and invariance
under multiplication by positive scalars of $\le_H$.
\end{proof}

For every ordered group $G$ below, the notation $D\downarrow0$
means $D\subseteq G_+$ is nonempty and downward directed and
$\inf_GD=0$. All lower bounds are computed in the whole group.
For a net $(x_\alpha)$ in $E$, its eventual sets $\mathcal F=\operatorname{Tail}((x_\alpha))$, and a map $k:E\to G$,
\[
 k_*\mathcal F\longrightarrow_o0
 \quad\Longleftrightarrow\quad
 \exists D\downarrow0\ \forall u\in D:
 \{x\in E:-u\le k(x)\le u\}\in\mathcal F.
\]
In particular, neither order in the positive cone alone nor order
in the image subspace is substituted for the order of $G$.

\begin{proposition}[Discrete ordinary order in the unreduced model]
\label{fp:unreduced-discrete}
For $E=\mathbb R$ with its usual topology, every tail family of a net
$\mathcal F$ on $\mathbb R$ satisfies
\[
 i_*\mathcal F\longrightarrow_o0
 \quad\Longleftrightarrow\quad \{0\}\in\mathcal F.
\]
\end{proposition}
\begin{proof}
Every seminorm on $\mathbb R$ has the form
$p_c(x)=c|x|$, where $c=p_c(1)\ge0$. Each such seminorm is
continuous. Define $h\in H$ by
\[
 h_0=1,\qquad h_{p_c}=-1\ (c>0),\qquad h_{p_0}=0,
 \qquad t_{p_c}(h)=0\ (c\ge0).
\]
The $p_1$ auxiliary coordinate of $-h$ is $1$, while its scalar
coordinate is $0$. Hence $-h\notin P_H$ and $h\not\le_H0$.

Suppose $r\ge_H i(x)$ and $r\ge_H-i(x)$ for some $x\ne0$.
Write $r=(a,(v_c,s_c)_{c\ge0})$. Both $r-i(x)$ and $r+i(x)$
belong to $P_H$. Their average is $r$, so $r\in P_H$ by
Lemma~\ref{fp:cone-order}. The coordinate inequalities give
\[
 c|a+v_c-x|\le s_c,\qquad c|a+v_c+x|\le s_c.
\]
For $c>0$,
\[
 2c|x|
 =c|(a+v_c+x)-(a+v_c-x)|
 \le c|a+v_c+x|+c|a+v_c-x|\le2s_c,
\]
so $s_c\ge c|x|>0$. In $r-h$, the aggregate vector in this
coordinate is
\[
 (a-1)+(v_c+1)=a+v_c,
\]
and consequently its seminorm is at most $s_c$. The strict
auxiliary condition is satisfied because $s_c>0$ independently
of the new auxiliary vector. At $c=0$, both auxiliary and scalar
coordinates are unchanged, and $p_0$ vanishes identically.
Therefore $r-h\in P_H$, or $h\le_Hr$.

Suppose now that $i_*\mathcal F\to_o0$ and
$\{0\}\notin\mathcal F$. Let $D\downarrow0$ be its regulator
set. For $r\in D$, the set
\[
 A_r=\{x\in\mathbb R:-r\le_H i(x)\le_Hr\}
\]
belongs to $\mathcal F$. It cannot be contained in $\{0\}$,
since upward closure of the represented tail family would then give
$\{0\}\in\mathcal F$. Thus $A_r$ contains a nonzero point.
The preceding computation gives $h\le_Hr$ for every $r\in D$.
This contradicts $\inf_HD=0$, since $h\not\le_H0$.
Conversely, if $\{0\}\in\mathcal F$, the singleton regulator
$D=\{0\}$ gives $i_*\mathcal F\to_o0$.
\end{proof}

\subsubsection{An explicit ordered quotient}

The preceding lower bound suggests removing the auxiliary ambiguity
through a quotient. For $p\in\mathcal P$, let
\[
 V_p=E/\ker p,\qquad q_p:E\longrightarrow V_p,\qquad
 \|q_p(x)\|_p=p(x),\qquad L_p=V_p\oplus\mathbb R.
\]
The norm is well defined, since $x-y\in\ker p$ gives
$|p(x)-p(y)|\le p(x-y)=0$. Its zero vectors are the zero coset.
No norm completeness is imposed on $V_p$. Set
\[
 C_p=\{(v,t)\in L_p:\|v\|_p\le t\},\qquad
 Q=\prod_{p\in\mathcal P}L_p,
\]
with coordinate order determined by $C_p$. Positivity and
pointedness of each $C_p$ follow from
\[
 \|v+w\|_p\le\|v\|_p+\|w\|_p,\qquad
 \|sv\|_p=s\|v\|_p\ (s\ge0),
\]
\[
 (v,t),(-v,-t)\in C_p
 \Longrightarrow 0\le\|v\|_p\le t\le-\|v\|_p\le0.
\]
Define the linear map
\begin{equation}\label{fp:quotient-map}
 T:H\longrightarrow Q,\qquad
 T\bigl(a_0,(a_p,t_p)_p\bigr)
   =\bigl((q_p(a_0+a_p),t_p)\bigr)_p.
\end{equation}

\begin{proposition}[The quotient cone]\label{fp:quotient-cone}
The map $T$ is surjective and
\[
 T(P_H)=Q_+=\prod_{p\in\mathcal P}C_p.
\]
It induces a linear order isomorphism between $H/\ker T$, with
the image cone, and $Q$. The linear map
\[
 j=T\circ i:E\longrightarrow Q,\qquad
 j(x)=\bigl((q_p(x),0)\bigr)_p,
\]
is injective.
\end{proposition}
\begin{proof}
Let $((v_p,t_p))_p\in Q$. Each fiber $q_p^{-1}\{v_p\}$ is
nonempty, so their product is nonempty. Any tuple $(a_p)_p$ in
this product, together with $a_0=0$, defines $a\in H$ with
$T(a)=((v_p,t_p))_p$. Thus $T$ is surjective.
If $a\in P_H$, then
\[
 \|q_p(a_0+a_p)\|_p=p(a_0+a_p)\le t_p\qquad(p\in\mathcal P),
\]
so $T(a)\in Q_+$.

Conversely, let $((v_p,t_p))_p\in Q_+$. There is a tuple $(a_p)_p$
with $q_p(a_p)=v_p$ and $a_p=0$ whenever $v_p=0$: the permissible
fiber is $\{0\}$ in those coordinates and $q_p^{-1}\{v_p\}$
in the others. Set $a_0=0$. Then
\[
 t_p\ge\|v_p\|_p=p(a_p)\ge0,
\]
and $a_p\ne0$ implies $v_p\ne0$, hence
$0<\|v_p\|_p\le t_p$. These inequalities verify every condition
in \eqref{fp:cone}. Thus $a\in P_H$ and $T(a)=((v_p,t_p))_p$.
This proves the image-cone equality.

The map $a+\ker T\mapsto T(a)$ is well defined, linear, injective,
and surjective. Its image of the quotient cone
$\{a+\ker T:a\in P_H\}$ is $Q_+$. The latter cone is pointed,
so this is an order isomorphism of partially ordered vector spaces.
Finally,
\[
 j(x)=0\Longrightarrow q_p(x)=0\ (p\in\mathcal P)
 \Longrightarrow p(x)=0\ (p\in\mathcal P)
 \Longrightarrow x=0.
\]
\end{proof}

Having identified positivity in the quotient, we must still determine
which arbitrary ordinary order-null nets are preserved and reflected
by $T$; the image-cone identity alone does not settle the behavior of
their regulators.

\subsubsection{Zero regulators and the original bounded sets}

To answer that question, we first examine a regulator in one norm
cone and then project a product regulator onto its coordinates. These
two calculations will relate the quotient convergence to bounded tails
in the original space.

\begin{lemma}[Regulators in a norm cone]\label{fp:norm-cone-regulator}
Let $V$ be a real normed space, and order $L=V\oplus\mathbb R$
by $C=\{(v,t):\|v\|\le t\}$. If
$(v_\delta,t_\delta)\downarrow0$ in $L$, then
$t_\delta\downarrow0$ in $\mathbb R$.
\end{lemma}
\begin{proof}
Positivity gives $t_\delta\ge0$. For $\eta\ge\delta$,
$(v_\delta-v_\eta,t_\delta-t_\eta)\in C$, so
\[
 0\le\|v_\delta-v_\eta\|\le t_\delta-t_\eta.
\]
Thus $(t_\delta)$ decreases. Let $r=\inf_\delta t_\delta$.
If $r>0$, there exists $\delta_0$ with $t_{\delta_0}<2r$.
Set
\[
 z=(v_{\delta_0},2r-t_{\delta_0}).
\]
For $\delta\ge\delta_0$,
\[
 \|v_\delta-v_{\delta_0}\|
 \le t_{\delta_0}-t_\delta
 \le t_\delta+t_{\delta_0}-2r
 =t_\delta-(2r-t_{\delta_0}),
\]
where the second inequality is $t_\delta\ge r$.
Hence $z\le(v_\delta,t_\delta)$ on that tail. For an arbitrary
index $\eta$, directedness provides $\delta\ge\eta,\delta_0$,
and
\[
 z\le(v_\delta,t_\delta)\le(v_\eta,t_\eta).
\]
Thus $z$ is a lower bound of the whole regulator net. Its scalar
coordinate $2r-t_{\delta_0}$ is positive, whereas $z\le0$ would
require $\|v_{\delta_0}\|\le-(2r-t_{\delta_0})<0$.
This contradicts the infimum $0$ in $L$. Therefore $r=0$.
\end{proof}

\begin{lemma}[Coordinates of a product regulator]\label{fp:product-regulator}
If $u_\delta\downarrow0$ in a product $\prod_{p\in\mathcal P}G_p$
of partially ordered abelian groups, then
$u_{\delta,p}\downarrow0$ in $G_p$ for every $p$.
\end{lemma}
\begin{proof}
Positivity and decrease hold coordinatewise. If $b\in G_p$ is a
lower bound of all $u_{\delta,p}$, define $z$ in the product by
$z_p=b$ and $z_q=0$ for $q\ne p$. For every $\delta$,
\[
 z_p\le u_{\delta,p},\qquad
 z_q=0\le u_{\delta,q}\quad(q\ne p).
\]
Hence $z\le u_\delta$ for every $\delta$ and therefore $z\le0$.
Its $p$ coordinate gives $b\le0$. This proves the asserted
coordinate infimum in the whole group $G_p$.
\end{proof}

A set $B\subseteq E$ is bounded if
\[
 \sup_{x\in B}p(x)<\infty\qquad(p\in\mathcal P).
\]
A bounded tail of $(x_\alpha)$ means that one tail-value set
$T_{\alpha_0}(x)$ satisfies all these seminorm bounds simultaneously.
Equivalently, some eventual set of the net is bounded.

\begin{theorem}[Exact net convergence in the quotient]
\label{fp:tail-dictionary}
For every net $(x_\alpha)$ in $E$,
\begin{equation}\label{fp:tail-equivalence}
 j(x_\alpha)\longrightarrow_o0
 \quad\Longleftrightarrow\quad
 \left[\begin{gathered}
 x_\alpha\longrightarrow0\text{ in the original topology},\\
 (\exists\alpha_0)\ T_{\alpha_0}(x)\text{ is bounded in }E
 \end{gathered}\right].
\end{equation}
At an arbitrary point $x\in E$, the same equivalence holds after
replacing the two limits by $x$ and $j(x)$ respectively.
\end{theorem}
\begin{proof}
Write $\mathcal F=\operatorname{Tail}((x_\alpha))$. Suppose
$j(x_\alpha)\to_o0$. Represent its regulator set as a
decreasing net
\[
 u_\delta=\bigl((v_{\delta,p},t_{\delta,p})\bigr)_p\downarrow0,
 \qquad
 A_\delta=\{x:-u_\delta\le j(x)\le u_\delta\}\in\mathcal F.
\]
For $x\in A_\delta$ and $p\in\mathcal P$,
\[
 \|v_{\delta,p}-q_p(x)\|_p\le t_{\delta,p},\qquad
 \|v_{\delta,p}+q_p(x)\|_p\le t_{\delta,p}.
\]
Therefore
\begin{align*}
 2p(x)
 &=\|(v_{\delta,p}+q_p(x))
           -(v_{\delta,p}-q_p(x))\|_p\\
 &\le\|v_{\delta,p}+q_p(x)\|_p
       +\|v_{\delta,p}-q_p(x)\|_p
 \le2t_{\delta,p},
\end{align*}
so
\begin{equation}\label{fp:seminorm-radius}
 x\in A_\delta\Longrightarrow p(x)\le t_{\delta,p}.
\end{equation}
For one fixed $\delta_0$, these bounds are finite for every $p$.
Hence $A_{\delta_0}\in\mathcal F$ is bounded.
Lemmas~\ref{fp:product-regulator} and
\ref{fp:norm-cone-regulator} give $t_{\delta,p}\downarrow0$ for
each $p$. For $\varepsilon>0$, an index $\delta$ with
$t_{\delta,p}<\varepsilon$ gives
\[
 A_\delta\subseteq\{x:p(x)<\varepsilon\}\in\mathcal F.
\]
Finite intersections of these seminorm balls form a base of the
original zero-neighborhood family. They all belong to $\mathcal F$,
so $\mathcal F\to0$ originally.

Conversely, suppose $\mathcal F\to0$ and $B\in\mathcal F$ is
bounded. Let
\[
 \mathcal A=\{A\in\mathcal F:A\subseteq B\},
 \qquad A\preceq A'\Longleftrightarrow A'\subseteq A.
\]
The set $\mathcal A$ is nonempty, and $A\cap A'\in\mathcal A$
is an upper index of $A,A'$. Properness makes every member
nonempty. Define
\[
 s_p(A)=\sup_{x\in A}p(x)\in[0,\infty),\qquad
 r_A=\bigl((0,s_p(A))\bigr)_p\in Q_+.
\]
Finiteness follows from $A\subseteq B$. These $r_A$ decrease.
For each fixed $p$ and $\varepsilon>0$,
\[
 A=B\cap\{x:p(x)<\varepsilon\}\in\mathcal A,
 \qquad 0\le s_p(A)\le\varepsilon.
\]
Thus $\inf_{A\in\mathcal A}s_p(A)=0$.
Let $z=((w_p,a_p))_p\in Q$ be any lower bound of all $r_A$.
For every $A,p$,
\[
 \|w_p\|_p\le s_p(A)-a_p,
 \qquad
 \|w_p\|_p\le\inf_{A\in\mathcal A}(s_p(A)-a_p)=-a_p.
\]
Hence $(-w_p,-a_p)\in C_p$ for every $p$, or $z\le0$ in $Q$.
Since $0\le r_A$, this proves $\inf_Qr_A=0$ in the whole group.
Finally, $x\in A$ implies $\|q_p(x)\|_p\le s_p(A)$ for every
$p$, and consequently
\[
 A\subseteq\{x:-r_A\le j(x)\le r_A\}\in\mathcal F.
\]
Thus the regulators $r_A\downarrow0$ prove $j_*\mathcal F\to_o0$.

The map $j$ is linear, so translation by $x$ in $E$ corresponds
to translation by $j(x)$ in $Q$. Both original and order
convergence commute with these translations. Boundedness also
commutes with translation, since
\[
 p(y-x)\le p(y)+p(x),\qquad p(y)\le p(y-x)+p(x).
\]
This proves the assertion at an arbitrary limit. For a net, a
bounded member of its eventual-set family contains a tail, which is then
bounded. Conversely a bounded tail is a bounded member of that
eventual-set family. This proves the net formulation.
\end{proof}

\subsubsection{Faithfulness of the quotient and normability}

We can now ask when the bounded-tail requirement places no additional
restriction on the original convergent nets. For this particular
embedding, the answer is exactly normability.

\begin{corollary}[Exact criterion for this linear representation]
\label{fp:normability}
For the map $j:E\to Q$ above, the following conditions are equivalent:
\begin{enumerate}
\item For every net $(x_\alpha)$ in $E$ and every $x\in E$,
$x_\alpha\to x$ originally if and only if
$j(x_\alpha)\to_o j(x)$.
\item Every originally convergent net in $E$ has a bounded tail.
\item $E$ has a bounded zero neighborhood.
\item The original locally convex topology of $E$ is induced by a norm.
\end{enumerate}
\end{corollary}
\begin{proof}
The equivalence of (1) and (2) is
Theorem~\ref{fp:tail-dictionary}. The zero-neighborhood family
$\mathcal N_E(0)$ converges to $0$. Applying (2) to it gives a
bounded member $B\in\mathcal N_E(0)$, which is a bounded zero
neighborhood, proving (3). Conversely, if $B$ is a bounded zero
neighborhood and $\mathcal F\to x$, then $x+B\in\mathcal F$.
The estimates $p(x+b)\le p(x)+p(b)$ show that $x+B$ is bounded.
Thus (3) implies (2).

Assume (3). Local convexity provides an open absolutely convex
zero neighborhood $U$ contained in a bounded zero neighborhood.
Hence $U$ is bounded. For each $x$, continuity of the scalar map
$t\mapsto tx$ at $0$ implies that $x\in sU$ for some $s>0$.
Define its finite Minkowski functional
\[
 p_U(x)=\inf\{s>0:x\in sU\}.
\]
Balancedness gives $rU\subseteq sU$ for $0<r\le s$. Consequently
$s>p_U(x)$ implies $x\in sU$. The identities $U=-U$ and
$a(sU)=|a|sU$ for $a\ne0$ give
\[
 p_U(ax)=|a|p_U(x)\quad(a\ne0),\qquad p_U(0)=0.
\]
For $s>p_U(x)$ and $t>p_U(y)$, write $x=su$, $y=tv$ with
$u,v\in U$. Convexity gives
\[
 \frac{x+y}{s+t}=\frac{s}{s+t}u+\frac{t}{s+t}v\in U,
 \qquad p_U(x+y)\le s+t.
\]
Taking the infimum over these $s,t$ proves
$p_U(x+y)\le p_U(x)+p_U(y)$.

If $p_U(x)<1$, there exists $s<1$ with $x\in sU\subseteq U$.
If $x\in U$, openness of $U$ and continuity of $t\mapsto tx$
at $1$ give some $a>1$ with $ax\in U$, and therefore
$p_U(x)\le a^{-1}<1$. Thus
\[
 U=\{x:p_U(x)<1\},\qquad
 \varepsilon U=\{x:p_U(x)<\varepsilon\}\quad(\varepsilon>0).
\]
These sublevel sets are originally open. Subadditivity gives
$|p_U(x)-p_U(y)|\le p_U(x-y)$, so $p_U$ is originally continuous.
If $p_U(x)=0$, then $x\in tU$ for every $t>0$. For every
$p\in\mathcal P$,
\[
 0\le p(x)\le t\sup_{u\in U}p(u)\quad(t>0).
\]
The supremum is finite by boundedness of $U$. Hence $p(x)=0$
for every $p$, and Hausdorffness gives $x=0$. Therefore $p_U$
is a norm.

For any original zero neighborhood $W$, there are finitely many
$p_1,\ldots,p_m\in\mathcal P$ and positive numbers
$\varepsilon_1,\ldots,\varepsilon_m$ such that
\[
 V=\bigcap_{k=1}^{m}\{x:p_k(x)<\varepsilon_k\}\subseteq W.
\]
A nonempty list may be used by adjoining the zero seminorm when
needed. Put
\[
 M_k=\sup_{u\in U}p_k(u)<\infty,\qquad
 \eta=\min_{1\le k\le m}\frac{\varepsilon_k}{2(1+M_k)}>0.
\]
For $x\in\eta U$,
\[
 p_k(x)\le\eta M_k<\varepsilon_k\quad(1\le k\le m),
 \qquad \{x:p_U(x)<\eta\}=\eta U\subseteq V\subseteq W.
\]
Thus the norm balls form a base of the original zero-neighborhood
eventual-set family. The two topologies coincide, proving (4).

Finally, suppose a norm $\|\cdot\|$ induces the original topology.
For every original continuous seminorm $p$, there is
$\varepsilon_p>0$ with
\[
 \|x\|<\varepsilon_p\Longrightarrow p(x)<1.
\]
For $x\ne0$, applying this inequality to
$\varepsilon_p x/(2\|x\|)$ gives
\[
 p(x)<\frac{2}{\varepsilon_p}\|x\|.
\]
The inequality also holds with a non-strict sign at $x=0$.
Therefore the norm unit ball $B=\{x:\|x\|<1\}$ satisfies
$\sup_{x\in B}p(x)\le2/\varepsilon_p$ for every $p$.
It is a bounded zero neighborhood, proving (3).
\end{proof}

\begin{example}[An unbounded neighborhood family in a countable product]
\label{fp:product-counterexample}
Let $E=\mathbb R^{\mathbb N}$ with its product topology. For every
$B\in\mathcal N_E(0)$, there are a finite set $F\subseteq\mathbb N$
and positive numbers $\varepsilon_k$, $k\in F$, such that
\[
 U=\{x:|x_k|<\varepsilon_k\ (k\in F)\}\subseteq B.
\]
There is $m\in\mathbb N\setminus F$. For every $t\in\mathbb R$,
$te_m\in U$, while the continuous coordinate seminorm
$p_m(x)=|x_m|$ satisfies $p_m(te_m)=|t|$. Hence
\[
 \sup_{x\in B}p_m(x)=\infty.
\]
No member of $\mathcal N_E(0)$ is bounded. This eventual-set family converges
to $0$ originally, but Theorem~\ref{fp:tail-dictionary} gives
\[
 j_*\mathcal N_E(0)\not\longrightarrow_o0.
\]
Thus ordinary order convergence in this explicit quotient does not
represent every original convergent net of this product space.
\end{example}

Corollary~\ref{fp:normability} gives a normability criterion for the
particular map $j$ with its untruncated seminorm coordinates, while other
faithful invariant metrics of $E$ remain outside that restriction. The
role of the quotient can be read together with
Proposition~\ref{fp:unreduced-discrete}, where the auxiliary coordinates
in $H$
produce whole-group lower bounds that are not at most zero: removing
those bounds leads to the quotient, although
Theorem~\ref{fp:tail-dictionary} still retains the requirement of a
bounded member. These computations determine the ordinary-order
convergence of the constructions and allow comparison with the cited
theorem while keeping its unbounded local order convergence distinct.

\subsection{The exact boundary for totally ordered values}
\label{sx:total-order}

Returning from a prescribed conversion to the existence question, we
ask what total order of the values requires of the original convergence.
The answer will be read from the order of its neighborhood bases, once
positive order convergence in a totally ordered group is made explicit.

If $V$ is an ordered abelian group, write
$V_+=\{v\in V:v\geq0\}$. A net $(v_\alpha)$ in $V_+$ is
\emph{order null} if there is a nonempty downward-directed set
$D\subseteq V_+$ such that
\begin{equation}\label{sx:total-regulator}
 \inf_VD=0,
 \qquad
 (\forall r\in D)(\exists\alpha_r)(\forall\alpha\ge\alpha_r)
 \quad0\le v_\alpha\le r.
\end{equation}
We may equivalently write $[0,r]\in\mathcal H$ for
$\mathcal H=\operatorname{Tail}((v_\alpha))$. The choice $D=\{0\}$
is allowed and includes eventually zero nets. The infimum in \eqref{sx:total-regulator} is taken in the entire
group $V$, not merely in the range of a distance function.

Let $(A,+,\lambda)$ be an abelian convergence group: addition and inversion
are continuous for $\lambda$. A $V$-valued invariant metric is a function
$d:A\times A\to V_+$ satisfying definiteness, symmetry, the triangle
inequality, and
\[
 d(x+a,y+a)=d(x,y)\qquad(x,y,a\in A).
\]
It is \emph{faithful} if, for every net $(x_\alpha)$ in $A$ and every
$x\in A$,
\begin{equation}\label{sx:total-faithful}
 x_\alpha\xrightarrow{\lambda}x
 \quad\Longleftrightarrow\quad
 d(x_\alpha,x)\xrightarrow{o}0.
\end{equation}
Here the expression on the right is the image net, and its convergence
has exactly the meaning \eqref{sx:total-regulator}. In particular, faithfulness is not being
required only for sequences.

\begin{lemma}[Positive convergence in a totally ordered group]\label{sx:total-null}
Let $V\ne\{0\}$ be a totally ordered abelian group. For a net
$(v_\alpha)$ in $V_+$, put $\mathcal H=\operatorname{Tail}((v_\alpha))$.
The following are equivalent:
\begin{enumerate}
 \item $\mathcal H$ is order null.
 \item $[0,\varepsilon)\in\mathcal H$ for every $\varepsilon>0$.
 \item $\mathcal H$ converges to $0$ in the order topology of $V$,
       restricted to $V_+$.
\end{enumerate}

\end{lemma}
\begin{proof}
If $D$ witnesses \eqref{sx:total-regulator}, then for every $\varepsilon>0$ some $r\in D$
satisfies $r<\varepsilon$. Otherwise $\varepsilon$ would be a positive
lower bound of $D$. Hence $[0,\varepsilon)\in\mathcal H$.

Conversely, suppose that condition 2 holds. If $V$ has a least positive
element $a$, then $[0,a)=\{0\}\in\mathcal H$, and $D=\{0\}$ works.
If there is no least positive element, $D=\{r\in V:r>0\}$ is downward
directed and has infimum zero. Indeed, no positive element is a lower
bound, since it has a strictly smaller positive element. Condition 2
implies $[0,r]\in\mathcal H$ for every $r\in D$. Finally, the sets
$[0,\varepsilon)$ form a neighborhood base at zero in the relative
order topology. This proves the equivalences.\end{proof}

\begin{theorem}[Totally ordered realizations]\label{sx:total-classification}
Let $(A,+,\lambda)$ be a Hausdorff abelian convergence group. The following
conditions are equivalent:
\begin{enumerate}
 \item There are a totally ordered abelian group $V$ and a faithful
       invariant $V$-valued metric on $A$.
 \item The convergence $\lambda$ is the convergence of a group topology,
       and the neighborhood family of zero has a cofinal base linearly
       ordered by inclusion.
\end{enumerate}
The discrete case is included. In the nondiscrete case let $\kappa$ be
the least cardinality of a linearly ordered cofinal neighborhood base.
If $\kappa=\aleph_0$, one can take $V=\mathbb R$. If $\kappa>\aleph_0$,
one can take
\[
 V=\bigoplus_{\alpha<\kappa}\mathbb Z e_\alpha
\]
with its lexicographic order. In the uncountable case the constructed metric satisfies
the strong triangle inequality.

\end{theorem}
\begin{proof}
For $A=\{0\}$, the zero metric with values in $\mathbb R$ represents its unique tail family of a net. Suppose henceforth that $A\ne\{0\}$. By Lemma~\ref{sx:total-null}, the right side of \eqref{sx:total-faithful}
is equivalent to
\begin{equation}\label{sx:total-balls}
 B_d(x,\varepsilon)=\{y:d(x,y)<\varepsilon\}\in\mathcal F
 \quad\hbox{for every }\varepsilon>0.
\end{equation}
These balls form neighborhood bases for a topology. To verify the only
point at issue, let $y\in B_d(x,\varepsilon)$ and put
$r=\varepsilon-d(x,y)>0$. If $z\in B_d(y,r)$, then
\[
 d(x,z)\leq d(x,y)+d(y,z)<d(x,y)+r=\varepsilon.
\]
Consequently every ball is open and its center belongs to it. Condition
\eqref{sx:total-balls} equals convergence of eventual-set families in this topology. Faithfulness
therefore says that this topology induces $\lambda$. The given continuity
of addition and inversion makes it a group topology. The balls centered
at zero are linearly ordered by inclusion, since $V$ is totally ordered.

\emph{Discrete case.}
Use the real metric $d(x,y)=1$ for $x\ne y$ and $d(x,x)=0$.
Lemma~\ref{sx:total-null} proves \eqref{sx:total-faithful}.

\emph{The countable case.}
Fix a decreasing symmetric neighborhood base $(U_n)_{n\geq0}$ at zero
such that
\begin{equation}\label{sx:dyadic-base}
 U_0=A,\qquad U_{n+1}+U_{n+1}\subseteq U_n.
\end{equation}
Such a base is obtained recursively from any countable neighborhood base
using continuity of addition. Define
\begin{equation}\label{sx:path-gauge}
 p(x)=\inf\left\{
       \sum_{j=1}^{m}2^{-n_j}:
       m\geq1,\ n_j\geq0,\ x=\sum_{j=1}^{m}x_j,
       \ x_j\in U_{n_j}\right\}.
\end{equation}
The admissible set is nonempty because $U_0=A$. Symmetry of the $U_n$
gives $p(-x)=p(x)$, and concatenating decompositions gives
$p(x+y)\leq p(x)+p(y)$. Also $p(0)=0$, since $0\in U_n$ for every $n$.

We record the estimate which proves definiteness and compatibility:
\begin{equation}\label{sx:path-estimate}
 p(x)<2^{-n}\ \Longrightarrow\ x\in U_n,
 \qquad
 x\in U_n\ \Longrightarrow\ p(x)\leq2^{-n}.
\end{equation}
Only the first implication needs proof. Fix a decomposition in \eqref{sx:path-gauge}
whose cost is less than $2^{-n}$. Whenever two terms have the same index
$k$, replace them by their sum in $U_{k-1}$, using \eqref{sx:dyadic-base}, and replace the
two weights $2^{-k}$ by the single weight $2^{-(k-1)}$. This preserves the
cost and decreases the number of terms. All indices stay greater than
$n$, because a single weight with index at most $n$ would contradict the
strict bound on the total cost. After finitely many replacements the
remaining indices are distinct, say $r_1<\cdots<r_l$ with $r_1>n$.
Starting with the last term and using \eqref{sx:dyadic-base} backwards shows that the sum
of the terms with indices $r_i,\ldots,r_l$ belongs to $U_{r_i-1}$.
Indeed, the tail following the term with index $r_i$ belongs to
$U_{r_{i+1}-1}\subseteq U_{r_i}$, so adding that term gives membership
in $U_{r_i-1}$. The complete sum therefore belongs to
$U_{r_1-1}\subseteq U_n$. This proves \eqref{sx:path-estimate}.

Hausdorffness gives $\bigcap_nU_n=\{0\}$, so \eqref{sx:path-estimate} implies that
$p(x)=0$ only for $x=0$. Thus $d(x,y)=p(x-y)$ is an invariant real metric.
Moreover,
\[
 U_{n+1}\subseteq\{x:p(x)<2^{-n}\}\subseteq U_n.
\]
Its balls and the original neighborhoods are cofinal. Hence, for every
eventual-set family, its metric convergence is exactly its original convergence.
Lemma~\ref{sx:total-null} converts this assertion into \eqref{sx:total-faithful}.

\emph{The uncountable case.}
Fix a linearly ordered cofinal neighborhood base $\mathcal B$ of minimal
cardinality $\kappa>\aleph_0$. If $\mathcal C\subseteq\mathcal B$ and $|\mathcal C|<\kappa$,
minimality of $\kappa$ gives $B\in\mathcal B$ for which
$C\not\subseteq B$ for every $C\in\mathcal C$.
Linear ordering then gives
\[
 B\subseteq\bigcap_{C\in\mathcal C}C.
\]
For fewer than $\kappa$ arbitrary neighborhoods, first insert one
base member into each neighborhood and apply the displayed inclusion.
Their intersection is therefore a neighborhood.
Enumerate the initial base by $(C_\alpha)_{\alpha<\kappa}$.
Recursion supplies $B_\alpha\in\mathcal B$ satisfying
\[
 B_\alpha\subseteq C_\alpha\cap
                  \bigcap_{\beta<\alpha}B_\beta.
\]
At each stage fewer than $\kappa$ base members occur, so the preceding
lower-bound argument applies. This yields a decreasing cofinal base.
If a subset $D\subseteq\kappa$ were cofinal and $|D|<\kappa$, then
$(B_\alpha)_{\alpha\in D}$ would be a cofinal base of cardinality
smaller than $\kappa$. Thus $\kappa$ is regular.

For each $\alpha<\kappa$, fix symmetric neighborhoods
$(W_{\alpha,n})_{n\geq0}$ such that
\[
 W_{\alpha,0}\subseteq B_\alpha,
 \qquad W_{\alpha,n+1}+W_{\alpha,n+1}\subseteq W_{\alpha,n}.
\]
Then $J_\alpha=\bigcap_{n\geq0}W_{\alpha,n}$ is a subgroup and is a
neighborhood, because $\aleph_0<\kappa$. A subgroup which is a
neighborhood is open. Set
\[
 H_\alpha=\bigcap_{\beta\leq\alpha}J_\beta.
\]
Again $H_\alpha$ is a neighborhood and an open subgroup, since this
intersection has cardinality less than $\kappa$. The family $(H_\alpha)$
is decreasing and cofinal, and $\bigcap_{\alpha<\kappa}H_\alpha=\{0\}$.

Order $V=\bigoplus_{\alpha<\kappa}\mathbb Z e_\alpha$ lexicographically:
a nonzero element is positive if its first nonzero coordinate is
positive. The elements $e_\alpha$ decrease as $\alpha$ increases and
are coinitial among the strictly positive elements of $V$. In fact, if
a positive element has first nonzero coordinate $\gamma$, then
$0<e_\alpha<v$ for every $\alpha>\gamma$.

For $g\ne0$, let
\[
 j(g)=\min\{\alpha<\kappa:g\notin H_\alpha\},
 \qquad
 p(g)=e_{j(g)},
 \qquad p(0)=0.
\]
The minimum exists by Hausdorffness. Since each $H_\alpha$ is a subgroup,
$p(-g)=p(g)$ and
\[
 p(g+h)\leq\max\{p(g),p(h)\}.
\]
For the latter inequality, if $g,h,g+h$ are nonzero, membership of $g$
and $h$ in every $H_\beta$ with
$\beta<\min\{j(g),j(h)\}$ gives
$j(g+h)\geq\min\{j(g),j(h)\}$. Cases involving zero are immediate.
Thus $d(x,y)=p(x-y)$ is an invariant $V$-valued ultrametric. Its balls
satisfy the exact equality
\[
 \{g:p(g)<e_\alpha\}=H_\alpha.
\]
Their cofinality, together with Lemma~\ref{sx:total-null}, proves \eqref{sx:total-faithful} for every net.
\end{proof}

\subsection{Direct and reduced conversion of totally ordered monoid values}
\label{tx:total-conversion}

Within the totally ordered setting, the existence theorem leaves room
for a finer question about a prescribed monoid-valued distance. We now
compare its direct passage to the group of differences with a reduced
passage that removes infinitesimal group differences.

Let $M\ne\{0\}$ be a commutative cancellative monoid with a
translation-compatible total order, least element $0$, and order
cancellation. Thus $a\le b$ implies $a+c\le b+c$, and
$a+c\le b+c$ implies $a\le b$. Suppose that $M$ has no least
positive element. Let $j:M\to G(M)$ be its canonical ordered group
embedding from Subsection~\ref{sx:conversion}, and write
$M_{>0}=M\setminus\{0\}$. For every net $(x_\alpha)$ in $M$, with $\mathcal F=\operatorname{Tail}((x_\alpha))$,
\begin{equation}\label{tx:monoid-null-test}
 \mathcal F\to_o0
 \quad\Longleftrightarrow\quad
 [0,c)_M\in\mathcal F\quad(c\in M_{>0}).
\end{equation}
Indeed, if $D\downarrow_M0$ regulates $\mathcal F$, then for each
$c>0$ some $d\in D$ satisfies $d<c$. Otherwise $c$ would be a positive
lower bound of $D$. Thus $[0,d]_M\subseteq[0,c)_M$ is an eventual set.
Conversely, $M_{>0}$ is downward directed and has infimum zero, and
the assumed eventual-set family contains $[0,c]_M$ for every $c>0$. In particular,
\eqref{tx:monoid-null-test} is convergence to the minimum in the
intrinsic order topology of $M$.

\begin{proposition}[Exact criterion for direct conversion]
\label{tx:direct-criterion}
For every net $(x_\alpha)$ in $M$, with $\mathcal F=\operatorname{Tail}((x_\alpha))$,
\[
 j_*\mathcal F\to_o0\text{ in }G(M)
 \quad\Longrightarrow\quad \mathcal F\to_o0\text{ in }M.
\]
The reverse implication holds for every net if and only if
\begin{equation}\label{tx:strict-gap}
 a>b\quad\Longrightarrow\quad
 (\exists c\in M_{>0})\ b+c<a.
\end{equation}
Equivalently, direct arbitrary-net conversion holds if and only if
$\inf_{G(M)}j(M_{>0})=0$.
\end{proposition}
\begin{proof}
If $j_*\mathcal F$ is order-null, Lemma~\ref{sx:total-null} gives
$[0,j(c))_{G(M)}\in j_*\mathcal F$ for every $c>0$.
Its inverse image is $[0,c)_M$, and
\eqref{tx:monoid-null-test} gives reflection.

Every positive $g\in G(M)$ has an expression
$g=j(a)-j(b)$ with $a>b$. Such an element fails to be a lower
bound of $j(M_{>0})$ exactly when some $c>0$ satisfies
\[
 j(c)<j(a)-j(b)\quad\Longleftrightarrow\quad b+c<a.
\]
Since $G(M)$ is totally ordered and $j(M_{>0})$ is downward directed,
condition~\eqref{tx:strict-gap} is equivalent to its whole-group
infimum being zero.

Assume \eqref{tx:strict-gap}, and let $\mathcal F\to_o0$ in $M$.
For each $g>0$ in $G(M)$, the preceding equivalence supplies $c>0$
with $j(c)<g$. Equation~\eqref{tx:monoid-null-test} gives
$[0,c)_M\in\mathcal F$, and
\[
 j([0,c)_M)\subseteq[0,g)_{G(M)}.
\]
Lemma~\ref{sx:total-null} proves $j_*\mathcal F\to_o0$.

Conversely, let preservation hold for all nets, and consider
\[
 \mathcal N=\operatorname{Tail}(\langle[0,c]_M:c\in M_{>0}\rangle).
\]
The family $M_{>0}\downarrow_M0$ proves $\mathcal N\to_o0$.
For $g>0$, order-nullity of $j_*\mathcal N$ gives
$j^{-1}([0,g)_{G(M)})\in\mathcal N$. Therefore some $c>0$ has
\[
 [0,c]_M\subseteq j^{-1}([0,g)_{G(M)}),\qquad j(c)<g.
\]
No positive $g$ is a lower bound of $j(M_{>0})$, proving the claim.
\end{proof}

When the direct criterion fails, a quotient may still preserve the
original null nets. This returns us to totally ordered semigroup
distances in the setting studied by Reichel and Ruppert
\citep{rr1977total}, where the quotient order lets us examine both
the change of values and its effect on arbitrary nets.

\begin{theorem}[Reduction by infinitesimal group differences]
\label{tx:reduced-conversion}
For $G=G(M)$, put $|g|=g\vee(-g)$ and define
\begin{equation}\label{tx:infinitesimal-subgroup}
 I=\{g\in G:n|g|<j(s)\text{ for every }n\ge1,\ s\in M_{>0}\}.
\end{equation}
Then $I$ is a convex subgroup. Let $Q=G/I$ have its quotient total
order, let $q:G\to Q$ be the quotient homomorphism, and put
$\phi=q\circ j$. The following assertions are equivalent:
\begin{enumerate}[(i)]
\item Addition in $M$ is continuous at $(0,0)$ for its intrinsic order
topology.
\item For every $s>0$ there are $x,y>0$ with $x+y<s$.
\item For every net $(x_\alpha)$ in $M$, with $\mathcal F=\operatorname{Tail}((x_\alpha))$,
\begin{equation}\label{tx:reduced-null-equivalence}
 \mathcal F\to_o0\text{ in }M
 \quad\Longleftrightarrow\quad
 \phi_*\mathcal F\to_o0\text{ in }Q.
\end{equation}
\end{enumerate}
For every $s>0$, one has $\phi(s)>0$. Under these equivalent
conditions, composing a faithful $M$-valued metric with $\phi$
gives a faithful $Q$-valued metric. Invariance is preserved whenever
it is present in the original metric.
\end{theorem}
\begin{proof}
For $g,h\in I$ and $n\ge1$, totality and the group triangle inequality
give
\[
 n|g+h|\le n|g|+n|h|\le2n\max\{|g|,|h|\}<j(s)
 \qquad(s>0).
\]
Thus $g+h\in I$. Also $0\in I$ and $-g\in I$ if $g\in I$.
If $|u|\le|g|$ and $g\in I$, then
$n|u|\le n|g|<j(s)$ for every $n,s$, and $u\in I$.
Hence $I$ is a convex subgroup.

For completeness, the quotient positive cone is $Q_+=q(G_+)$.
It is closed under addition. If $z,-z\in Q_+$, there are $u,v\ge0$
with $q(u)=z$ and $q(v)=-z$. Then $u+v\in I$, and
$0\le u\le u+v$ implies $u\in I$ by convexity. Thus $z=0$, so
the cone is proper. For each $g\in G$, either $g\ge0$ or
$-g\ge0$, so either $q(g)\ge0$ or $-q(g)\ge0$. This defines
a translation-invariant total order, and $q$ is order preserving.
If $s>0$ and $j(s)\in I$, the defining inequality with $n=1$
and this same $s$ would be $j(s)<j(s)$. Therefore
\begin{equation}\label{tx:positive-not-annihilated}
 \phi(s)>0\qquad(s\in M_{>0}).
\end{equation}

Suppose (i). For $s>0$, continuity at $(0,0)$ supplies $u,v>0$ with
\[
 [0,u)_M+[0,v)_M\subseteq[0,s)_M.
\]
Since there is no least positive element, there exist $0<x<u$ and
$0<y<v$, giving $x+y<s$. Conversely, if (ii) holds, let
$x,y>0$ satisfy $x+y<s$. For $a<x$, $b<y$, strict translation
compatibility gives
\[
 a+b<x+y<s,\qquad [0,x)_M+[0,y)_M\subseteq[0,s)_M.
\]
The order-neighborhood bases prove (i).

Assume (ii). For every $s>0$ there is $t>0$ with $2t<s$:
take $t=\min(x,y)$ for the pair supplied by (ii). Iteration gives
\begin{equation}\label{tx:integer-small-elements}
 (\forall s>0)(\forall m\ge1)(\exists t>0)\quad mt<s.
\end{equation}
Indeed, for an integer $k\ge1$ with $2^k\ge m$, successive inequalities
$2t_1<s$ and $2t_{r+1}<t_r$ give
$mt_k\le2^kt_k<s$.

Let $z>0$ in $Q$. By the definition of $Q_+$, it has a representative
$g>0$ in $G\setminus I$. Definition~\eqref{tx:infinitesimal-subgroup}
then gives $n\ge1$ and $s>0$ with $ng\ge j(s)$.
Equation~\eqref{tx:integer-small-elements} supplies $t>0$ with
$2nt<s$. Hence
\[
 n(2j(t))<j(s)\le ng,\qquad 2j(t)<g,
 \qquad 0<\phi(t)<2\phi(t)\le q(g)=z.
\]
The middle implication uses totality: $2j(t)\ge g$ would imply
$n(2j(t))\ge ng$. Thus the positive values of $\phi$ are
coinitial in $Q_{>0}$.

For $\mathcal F\to_o0$ in $M$ and $z>0$ in $Q$, the preceding
construction gives $t>0$ with $\phi(t)<z$. The represented tail family contains
$[0,t)_M$, and
\[
 \phi([0,t)_M)\subseteq[0,\phi(t)]_Q\subseteq[0,z)_Q.
\]
Lemma~\ref{sx:total-null} proves $\phi_*\mathcal F\to_o0$.
Reflection does not require (ii). If $\phi_*\mathcal F\to_o0$,
then for every $s>0$, equation~\eqref{tx:positive-not-annihilated}
and that lemma give
\[
 \phi^{-1}([0,\phi(s))_Q)\in\mathcal F,\qquad
 \phi^{-1}([0,\phi(s))_Q)\subseteq[0,s)_M.
\]
The inclusion holds since $a\ge s$ implies
$\phi(a)\ge\phi(s)$. Equation~\eqref{tx:monoid-null-test}
gives reflection and proves (iii).

Finally assume (iii), and let
$\mathcal N_0=\operatorname{Tail}(\langle[0,s)_M:s>0\rangle)$.
Its representing saturation is null by \eqref{tx:monoid-null-test}. Independent sums of
positive zero regulators in an ordered group have whole-group
infimum zero: if $h\le d+e$ for all $d\in D,e\in E$, then
$h-d\le e$ for all $e$, hence $h\le d$ for all $d$, hence $h\le0$.
Consequently
\[
 \phi_*(\mathcal N_0+\mathcal N_0)
 =\phi_*\mathcal N_0+\phi_*\mathcal N_0\to_o0.
\]
Reflection in (iii) gives $\mathcal N_0+\mathcal N_0\to_o0$.
For every $s>0$, its tail base of the independent sum net supplies $u,v>0$ with
$[0,u)_M+[0,v)_M\subseteq[0,s)_M$. Positive $x<u$ and $y<v$
give $x+y<s$, proving (ii).

For an $M$-valued metric $d$, the map $d^*=\phi\circ d$ satisfies
\[
 d^*(x,z)\le d^*(x,y)+d^*(y,z),\qquad
 d^*(x,y)=d^*(y,x),\qquad d^*(x,y)=0\iff d(x,y)=0.
\]
These assertions use monotonicity, additivity, and
\eqref{tx:positive-not-annihilated}. Applying
\eqref{tx:reduced-null-equivalence} to every image pair net
proves arbitrary-net faithfulness. An identity
$d(x+a,y+a)=d(x,y)$ remains an identity after composition with $\phi$.
\end{proof}

\begin{example}[A direct embedding that makes the represented convergence discrete]
\label{tx:rational-lex-example}
Let
\[
 M=\{(a,b)\in\mathbb Q^2:a>0\}\cup\{(0,0)\}
\]
have coordinate addition and lexicographic order. Its difference group
is $G(M)=\mathbb Q^2_{\mathrm{lex}}$: for any $(r,s)\in\mathbb Q^2$
and any rational $A>\max(0,-r)$,
\[
 (r,s)=(A+r,s)-(A,0),
\]
with both terms in $M$. Addition in $M$ is continuous at zero, since
for $(a,b)>0$ and rational $0<t<a/2$,
$(t,0)+(t,0)<(a,b)$. Nevertheless \eqref{tx:strict-gap} fails:
\[
 (1,1)>(1,0),\qquad (1,0)+c>(1,1)\quad(c\in M_{>0}).
\]

On the additive group $\mathbb Q$, set
$p(x)=(|x|,0)$ and $d(x,y)=p(x-y)$. The ordinary rational triangle
inequality gives subadditivity in the lexicographic order. For a
threshold $(a,b)>0$,
\[
 |x|<a\quad\Longrightarrow\quad p(x)<(a,b),\qquad
 p(x)<(a,0)\quad\Longleftrightarrow\quad |x|<a.
\]
Thus \eqref{tx:monoid-null-test} identifies its $M$-valued convergence
with ordinary rational convergence for every net.
Regarded as a $G(M)$-valued distance, however,
\[
 \{x:p(x)<(0,1)\}=\{0\}.
\]
Lemma~\ref{sx:total-null} implies that its null nets must contain
$\{0\}$, and such eventual-set families are null with regulator $\{0\}$.
The group-valued realization is therefore discrete. In particular,
$(1/n,0)$ is null in $M$ and is not order-null in $G(M)$.

The subgroup in \eqref{tx:infinitesimal-subgroup} is
$I=\{0\}\times\mathbb Q$. Indeed, every vertical element satisfies
$n|(0,b)|<(a,c)$ for every $n$ and every $a>0$. An element with
nonzero first coordinate fails the defining inequality at $n=1$
against a smaller positive first coordinate. Therefore
\[
 G(M)/I\simeq\mathbb Q,\qquad \phi(a,b)=a,\qquad
 d^*(x,y)=|x-y|.
\]
The quotient restores the original convergence for all nets.
\end{example}

\subsection{Why the infimum must be taken in the whole group}
\label{tx:whole-group-infimum}

Throughout these conversions, zero has been tested against every lower
bound in the group. Restricting that test to the positive cone can
retain the represented point convergence while losing the operations
on errors that the integration arguments need, as the following
example shows.

\begin{example}[Cone-relative infima do not preserve the error calculus]
\label{tx:cone-relative}
Let $G=\mathbb R^2$ with the translation-invariant order determined by
\[
 C=\{(0,0)\}\cup(0,\infty)^2,\qquad
 u\le v\quad\Longleftrightarrow\quad v-u\in C.
\]
The set $C$ is closed under addition and $C\cap(-C)=\{0\}$, so this
is a partial order. On $\mathbb R$ define
\[
 p(0)=0,\qquad p(t)=(\sqrt{|t|},1)\quad(t\ne0),\qquad
 d(x,y)=p(x-y).
\]
The gauge is definite and symmetric. If $x,y,x+y\ne0$, then
\[
 p(x)+p(y)-p(x+y)
 =(\sqrt{|x|}+\sqrt{|y|}-\sqrt{|x+y|},1)\in(0,\infty)^2,
\]
since $\sqrt a+\sqrt b>\sqrt{a+b}$ for $a,b>0$ and
$|x+y|\le|x|+|y|$. If $x+y=0$, its positive right side dominates
$p(0)=0$. A zero summand gives equality. Thus $d$ is an invariant
ordered-group-valued metric.

Consider the alternative convention that regulators are nonempty
downward-directed subsets of $C$ whose infimum is zero only in the
poset $C$. Put
\[
 d_n=(1/n,1+1/n),\qquad e_n=(1+1/n,1/n).
\]
Both sequences decrease. If $c\in C\setminus\{0\}$ were below
every $d_n$, its positive first coordinate would satisfy
$c_1\le1/n$ for every $n$, a contradiction. The second coordinate
gives the corresponding contradiction for $(e_n)$. Therefore
\[
 \inf_C\{d_n:n\ge1\}=0=\inf_C\{e_n:n\ge1\}.
\]

The metric represents ordinary real convergence for all nets under
this alternative convention. If $\mathcal F\to0$ ordinarily, then
for each $n$ it contains $(-1/n^2,1/n^2)$, and
\[
 |t|<1/n^2\quad\Longrightarrow\quad p(t)\le d_n.
\]
Conversely, if $\mathcal F$ is not ordinarily null, some
$\varepsilon>0$ satisfies $(-\varepsilon,\varepsilon)\notin\mathcal F$.
Every member of $\mathcal F$ then meets
$\{t:|t|\ge\varepsilon\}$. If $D\subseteq C$ were a cone-relative
zero regulator for $p_*\mathcal F$, every ball $B_d\in\mathcal F$
would contain such a point, and hence
\[
 d_1\ge\sqrt\varepsilon,\qquad d_2\ge1\quad(d\in D).
\]
The nonzero element $c=(\sqrt\varepsilon/2,1/2)\in C$ would be a
lower bound of $D$ in $C$, contradicting $\inf_C D=0$.

Despite this faithfulness,
\[
 B_{(1,2)}=\{t:p(t)\le(1,2)\}=(-1,1)
\]
is not closed for the represented convergence. The element
$w=(0,1)$ satisfies
\[
 d_n-w=(1/n,1/n)\in C,\qquad w\not\le0.
\]
Thus zero is not the whole-group infimum of $(d_n)$. This is the
lower bound excluded by the whole-group convention.

Finally, all independently indexed sums satisfy
\[
 d_n+e_m\ge(1/2,1/2)>0\qquad(n,m\ge1).
\]
If these sums were cone-relative order-null, every bound of a regulator
would dominate the fixed positive element $(1/2,1/2)$, contradicting
its cone-relative infimum zero. Thus even addition of null regulators
fails under this altered convention.
\end{example}

\subsection{Support ideals in metrizable abelian groups}
\label{sx:support-ideals}

We next place the representation question in a family where the
additional convergence requirement can be seen directly in the supports
of nets. Closure of those supports will provide the criterion, allowing
both representable and nonrepresentable cases within the same ambient
topological setting.

Throughout this subsection, $(G,\kappa)$ is a Hausdorff metrizable
Abelian topological group, written additively.  A \emph{subgroup ideal}
$\mathcal I$ is a nonempty collection of subgroups of $G$ such that
\[
 K\le H\in\mathcal I\ \Longrightarrow\ K\in\mathcal I,
 \qquad
 H,K\in\mathcal I\ \Longrightarrow\ H+K\in\mathcal I.
\]
We require $\mathcal I$ to contain every countably generated subgroup
of $G$.  For a net $(x_\alpha)$ in $G$, define
\begin{equation}
 x_\alpha\longrightarrow_{\lambda_{\mathcal I}}x
 \quad\Longleftrightarrow\quad
 x_\alpha\longrightarrow_\kappa x
 \quad\hbox{and}\quad
 x_\alpha\in H\text{ eventually for some }H\in\mathcal I.
 \label{sx:eq-support-ideal-convergence}
\end{equation}
Although the support subgroup need not be centered at the limit,
translations preserve the condition because $\mathcal I$ contains
every cyclic subgroup: $H+x\subseteq H+\langle x\rangle\in\mathcal I$.

Fix a compatible invariant metric $\rho$ on $G$, bounded by one.
Its existence follows from the dyadic path-norm construction proved in
Theorem~\ref{sx:total-classification}. In the Abelian case invariance is
automatically invariance on both sides.  Replacing an invariant metric
by its minimum with one preserves compatibility and the triangle
inequality.  Completeness, when assumed below, means completeness of
the canonical group uniformity.

\begin{proposition}
\label{sx:prop-support-ideal-basic}
The relation $\lambda_{\mathcal I}$ is a Hausdorff convergence-group
structure satisfying the usual finite-mixing axioms for nets.  Every convergent
eventual-set family has a coarser countably based tail family with the same limit.
The $\lambda_{\mathcal I}$-closed subsets of $G$ are exactly its
$\kappa$-closed subsets.  If $(G,\kappa)$ is complete, then
$(G,\lambda_{\mathcal I})$ is complete.
\end{proposition}

\begin{proof}
The constant net at $x$ is eventually in the cyclic subgroup
$\langle x\rangle\in\mathcal I$.  Refinement preserves both conditions
in \eqref{sx:eq-support-ideal-convergence}.  If
$\mathcal F,\mathcal G\to_{\lambda_{\mathcal I}}x$ have support
subgroups $H,K\in\mathcal I$, then $H+K$ belongs to both net-tail families.
It therefore belongs to $\mathcal F\cap\mathcal G$, and the latter eventual-set family
also converges to $x$ in $\kappa$.  This proves the finite-mixing axiom for nets.
If $\mathcal F\to_{\lambda_{\mathcal I}}x$ and
$\mathcal G\to_{\lambda_{\mathcal I}}y$, their sum eventual-set family contains
$H+K$ and converges to $x+y$ in $\kappa$.  Their inverse eventual-set families have
the same support subgroups.  Addition and inversion are consequently
continuous.  Hausdorffness follows from that of $\kappa$.

Let $H\in\mathcal F\cap\mathcal I$ and
$\mathcal F\to_{\lambda_{\mathcal I}}x$.  The nonempty sets
\[
 H\cap\{z:\rho(z,x)<1/n\},\qquad n\ge1,
\]
generate a countably based tail family coarser than $\mathcal F$ and
converging to $x$ in $\lambda_{\mathcal I}$.

Every $\kappa$-closed set is $\lambda_{\mathcal I}$-closed because
$\lambda_{\mathcal I}$ refines $\kappa$.  Conversely, suppose that
$C$ is $\lambda_{\mathcal I}$-closed and $x\in\overline C^{\,\kappa}$.
Metrizability supplies a sequence $x_n\in C$ with $x_n\to_\kappa x$.
The subgroup generated by $\{x,x_1,x_2,\ldots\}$ belongs to
$\mathcal I$ and contains every term of the sequence.  Thus
$x_n\to_{\lambda_{\mathcal I}}x$, and $x\in C$.  This proves the
assertion about closed sets.

Finally, let $\mathcal F$ be $\lambda_{\mathcal I}$-Cauchy.  Then
$\mathcal F-\mathcal F\to_\kappa0$, so completeness of $(G,\kappa)$
provides $x$ with $\mathcal F\to_\kappa x$.  There are
$H\in\mathcal I$ and $A,B\in\mathcal F$ such that
\[
 A-B\subseteq H.
\]
For $b\in B$ this gives
\[
 A\subseteq b+H\subseteq H+\langle b\rangle\in\mathcal I.
\]
The subgroup $H+\langle b\rangle$ belongs to $\mathcal F$, proving
$\mathcal F\to_{\lambda_{\mathcal I}}x$.
\end{proof}

For the next theorem, a positive net $(v_\alpha)$ is order null
when some family $u_\gamma\downarrow0$ satisfies
$(\forall\gamma)(\exists\alpha_\gamma)(\forall\alpha\ge\alpha_\gamma)
\ v_\alpha\le u_\gamma$.
A \emph{faithful invariant ordered-group metric} is a map
$d(x,y)=p(x-y)$, where $p$ is definite, symmetric, positive and
subadditive, such that
\[
 \mathcal F\to_{\lambda_{\mathcal I}}0
 \quad\Longleftrightarrow\quad
 p(\mathcal F)\to_o0
\]
for every net $(x_\alpha)$, with
$\mathcal F=\operatorname{Tail}((x_\alpha))$.

\begin{theorem}
\label{sx:thm-support-ideal-classification}
For the convergence group \eqref{sx:eq-support-ideal-convergence}, the
following assertions are equivalent:
\begin{enumerate}
\item It admits a faithful invariant metric with values in a partially
ordered Abelian group.
\item Every null net is a quasi-subnet of a null saturation possessing a base of
sets closed for the original convergence.
\item The subgroup ideal is stable under the original metric closure:
\begin{equation}
 H\in\mathcal I\quad\Longrightarrow\quad
 \overline H^{\,\kappa}\in\mathcal I.
 \label{sx:eq-support-ideal-closure}
\end{equation}
\item It admits a faithful invariant metric with values in a
Dedekind-complete Archimedean weakly $\sigma$-distributive vector lattice.
\end{enumerate}
The assertions concern arbitrary directed nets.
Completeness of $(G,\kappa)$ is not needed for these equivalences.
\end{theorem}

\begin{proof}
First suppose that $p:G\to P_+$ is a faithful gauge with values in a
partially ordered Abelian group $P$.  For $r\in P_+$, the set
\[
 D_r=\{x\in G:p(x)\le r\}
\]
is closed for the original convergence.  Indeed, if $x_\alpha\in D_r$
and $x_\alpha\to_{\lambda_{\mathcal I}}x$, let $u_\gamma\downarrow0$
eventually bound $p(x_\alpha-x)$.  The triangle inequality gives,
for every $\gamma$ and an eventual $\alpha$,
\[
 p(x)-r\le p(x-x_\alpha)\le u_\gamma.
\]
The greatest-lower-bound property yields $p(x)-r\le0$.
If $\mathcal F$ is null, a regulator for $p(\mathcal F)$ gives a
coarser eventual-set family generated by the closed sets $D_{u_\gamma}$.
That eventual-set family is null by faithfulness.  Thus (1) implies (2).

Assume (2), and fix $H\in\mathcal I$.  The represented tail family
\[
 \mathcal F_H=\operatorname{Tail}\left(
 \bigl\langle H\cap U:U\text{ is a }\kappa\text{-neighborhood of }0
 \bigr\rangle\right)
\]
is null.  Let $\mathcal D\subseteq\mathcal F_H$ be a null tail family
with a closed base.  Since its saturation is null, $\mathcal D$ contains a
subgroup $K\in\mathcal I$.  There are a closed basic set
$C\in\mathcal D$ and an open neighborhood $U$ of zero such that
\[
 H\cap U\subseteq C\subseteq K.
\]
By Proposition~\ref{sx:prop-support-ideal-basic}, $C$ is
$\kappa$-closed.  Since $U$ is open,
\[
 \overline H^{\,\kappa}\cap U
 \subseteq\overline{H\cap U}^{\,\kappa}
 \subseteq C\subseteq K.
\]
For $x\in\overline H^{\,\kappa}$, density of $H$ in its closure
gives $h\in H$ with $x-h\in U$.  Both $x$ and $h$ belong to the
subgroup $\overline H^{\,\kappa}$, so $x-h\in K$.  Therefore
\[
 \overline H^{\,\kappa}\subseteq H+K\in\mathcal I.
\]
Hereditary stability proves \eqref{sx:eq-support-ideal-closure}.  This
argument does not require a small neighborhood to generate the whole
closed subgroup: density supplies the additional summand $H$.

Assume (3).  Let $\mathscr T$ be the set of all continuous maps
$t:G\to[0,1]$ satisfying
\[
 t(0)=0,\qquad t(-x)=t(x),\qquad
 t(x+y)\le t(x)+t(y).
\]
The map $t_0(x)=\rho(x,0)$ belongs to $\mathscr T$.  For every closed
subgroup $H\in\mathcal I$, define
\[
 r_H(t)=\sup_{h\in H}t(h),\qquad t\in\mathscr T.
\]
Thus $0\le r_H\le1$.  Define the following solid order ideal of the
pointwise ordered vector lattice $\mathbb R^{\mathscr T}$:
\begin{equation}
 L=\left\{u\in\mathbb R^{\mathscr T}:
 |u|\le\sum_{i=1}^{n}c_i r_{H_i}
 \text{ for some }n\ge1,\ c_i\ge0,
 \ H_i\in\mathcal I\text{ closed}\right\}.
 \label{sx:eq-support-ideal-value-lattice}
\end{equation}
It is a vector sublattice and a solid ideal.  It is also Dedekind
complete.  If a nonempty subset $S\subseteq L$ has upper bound
$b\in L$, then its pointwise supremum $s$ exists. For $s_0\in S$,
\[
 s_0\le s\le b,\qquad |s|\le |s_0|+|b|.
\]
Solidity gives $s\in L$, and this pointwise supremum is also the
supremum in $L$.  Archimedeanity is inherited from
$\mathbb R^{\mathscr T}$.
For a globally order-bounded array $0\le a_{n,k}\le e\in L$ with
$a_{n,k}\downarrow_k0$, all infima and suprema involved below belong
to $[0,e]$ and are pointwise. For each $t\in\mathscr T$ and $\eta>0$,
the row infima give a function $\varphi:\mathbb N\to\mathbb N$ with
$a_{n,\varphi(n)}(t)<\eta$ for every $n$. Therefore
\[
 0\le\left(\inf_{\varphi\in\mathbb N^{\mathbb N}}
                 \sup_n a_{n,\varphi(n)}\right)(t)\le\eta.
\]
The infimum is zero at every coordinate. Hence $L$ is weakly
$\sigma$-distributive.

For $x\in G$, set
\begin{equation}
 p(x)(t)=t(x),\qquad t\in\mathscr T.
 \label{sx:eq-support-ideal-gauge}
\end{equation}
The closed subgroup $\overline{\langle x\rangle}^{\,\kappa}$ belongs
to $\mathcal I$ by (3), so
\[
 0\le p(x)\le r_{\overline{\langle x\rangle}^{\,\kappa}}.
\]
Hence $p(x)\in L$.  Symmetry and subadditivity hold coordinatewise.
If $p(x)=0$, its $t_0$ coordinate is zero, so $x=0$.  Thus
$d(x,y)=p(x-y)$ is an invariant $L$-valued metric.

Suppose that $\mathcal F\to_{\lambda_{\mathcal I}}0$.  Its support
subgroup may be replaced by its closure. Write $H\in\mathcal I$
for a closed support subgroup belonging to $\mathcal F$.  For
$A\in\mathcal F$, define
\[
 v_A(t)=\sup_{x\in A\cap H}t(x),\qquad t\in\mathscr T.
\]
The intersections $A\cap H$ are nonempty, and
$0\le v_A\le r_H$, so $v_A\in L$.  Ordered by reverse inclusion
of $A$, this family is decreasing.  Continuity of every $t$ and
$\mathcal F\to_\kappa0$ give
\[
 \inf_{A\in\mathcal F}v_A(t)=0\qquad(t\in\mathscr T).
\]
Thus $v_A\downarrow0$ in $L$.  Since $p(x)\le v_A$ for
$x\in A\cap H$, this proves $p(\mathcal F)\to_o0$.

Conversely, suppose that $p(\mathcal F)\to_o0$ in $L$.  Decreasing
positive families in the solid ideal $L$ have the same infimum in both $L$ and
$\mathbb R^{\mathscr T}$ whenever the infimum exists: their
pointwise infimum lies between zero and any fixed member of the
family and therefore belongs to $L$.  Consequently, every coordinate
of $p(\mathcal F)$ tends to zero.  In particular its $t_0$ coordinate
tends to zero, proving $\mathcal F\to_\kappa0$.

One member of an order regulator provides $A\in\mathcal F$ and
$u\in L_+$ such that $p(x)\le u$ for every $x\in A$.  By
\eqref{sx:eq-support-ideal-value-lattice}, there are closed
$H_1,\ldots,H_n\in\mathcal I$ and $c_i\ge0$ with
\[
 u\le\sum_{i=1}^{n}c_i r_{H_i}.
\]
Set $V=\overline{H_1+\cdots+H_n}^{\,\kappa}\in\mathcal I$.
The continuous pseudonorm
\[
 t_V(x)=\inf_{v\in V}\rho(x,v)
\]
belongs to $\mathscr T$ and has kernel exactly $V$.  For clarity,
invariance of $\rho$ and the subgroup property give
\[
 t_V(-x)=t_V(x),\qquad
 t_V(x+y)\le t_V(x)+t_V(y),\qquad
 |t_V(x)-t_V(y)|\le\rho(x,y).
\]
The equality $t_V(x)=0$ is equivalent to
$x\in\overline V^{\,\kappa}=V$.  Thus no separate metrization
of the quotient is needed.  Because $H_i\subseteq V$,
$r_{H_i}(t_V)=0$ for every $i$.  Therefore, for $x\in A$,
\[
 0\le t_V(x)=p(x)(t_V)
 \le u(t_V)\le\sum_{i=1}^{n}c_i r_{H_i}(t_V)=0.
\]
It follows that $A\subseteq V$, whence $V\in\mathcal F$.
Together with $\mathcal F\to_\kappa0$, this is exactly
$\mathcal F\to_{\lambda_{\mathcal I}}0$.  We have proved (4).
Finally, (4) implies (1) by forgetting the lattice and scalar
operations of the value space.
\end{proof}

\begin{proposition}
\label{sx:prop-support-ideal-topological}
The convergence $\lambda_{\mathcal I}$ is topological if and only if
$\mathcal I$ contains a $\kappa$-open subgroup.  In that event,
$\lambda_{\mathcal I}$ is exactly $\kappa$.
\end{proposition}

\begin{proof}
If $H\in\mathcal I$ is open and $\mathcal F\to_\kappa x$, then
$x+H\in\mathcal F$ and
$x+H\subseteq H+\langle x\rangle\in\mathcal I$.  Thus every
$\kappa$-convergent net is $\lambda_{\mathcal I}$-convergent.
Conversely, if $\lambda_{\mathcal I}$ is induced by a topology, that
topology has the same closed sets as $\kappa$, by
Proposition~\ref{sx:prop-support-ideal-basic}, and therefore equals
$\kappa$.  Its neighborhood family at zero contains some
$H\in\mathcal I$.  A subgroup containing a neighborhood of zero is
open, proving the claim.
\end{proof}

\begin{example}[Separable support]
\label{sx:ex-separable-support-group}
Let
\[
 \mathcal I_{\rm sep}
 =\{H\le G:\overline H^{\,\kappa}\text{ is separable}\}.
\]
This is a subgroup ideal.  A subspace of a separable metrizable space
is separable. The sum of two separable subsets is separable because
addition is continuous. The closure of a separable subset is
separable.  These facts verify hereditary stability, finite-sum
stability and closure stability.  A countably generated Abelian
group is countable, so every such subgroup belongs to
$\mathcal I_{\rm sep}$.

Consequently, if $G$ is complete, convergence in $\kappa$ with a
separable subgroup support is complete, first countable in the
net sense, and represented faithfully by the single
Dedekind-complete-lattice-valued metric
\eqref{sx:eq-support-ideal-gauge}.  It is non-topological exactly when
$G$ has no open separable subgroup.  Mere nonseparability of $G$
does not suffice for this conclusion.

For a concrete nonvector example, let $J$ be uncountable and let
$G=\mathbb T^J$, with pointwise multiplication and the invariant
supremum metric
\[
 \rho(z,w)=\sup_{j\in J}\frac{|z_j-w_j|}{2}.
\]
To verify completeness, let $\mathcal F$ be Cauchy for this metric.
Every coordinate image net converges in the complete circle to
some $z_j$. Put $z=(z_j)_{j\in J}$.  For $\varepsilon>0$ there is
$A\in\mathcal F$ with $\rho(x,y)<\varepsilon/2$ for all $x,y\in A$.
For fixed $x\in A$, passing to the limit of the $j$-th coordinate
image net gives $|x_j-z_j|/2\le\varepsilon/2$ for every $j$.
Hence $\rho(x,z)\le\varepsilon/2<\varepsilon$, uniformly for
$x\in A$, and $\mathcal F\to z$.  Every open ball
contains an uncountable uniformly separated family: fix a sufficiently
small $\zeta\in\mathbb T\setminus\{1\}$ and, for each $j\in J$,
put $\zeta$ in coordinate $j$ and $1$ in every other coordinate.
Distinct members of this family have mutual distance
$|\zeta-1|/2>0$.  Hence no open ball is separable, and $G$ has no
open separable subgroup.  The separable-support convergence is
therefore non-topological.  The nontrivial element constantly equal
to $-1$ has order two, so this group cannot be the additive group of
a real vector space.  It is thus a complete, non-topological,
nonvector group with a faithful invariant lattice-group-valued
metric.
\end{example}

\begin{example}[Countable support]
Let $G=(\mathbb R,+)$ with its usual topology, and let
$\mathcal I$ consist of its countable subgroups.  Then
$\mathbb Q\in\mathcal I$ but
$\overline{\mathbb Q}=\mathbb R\notin\mathcal I$.
Theorem~\ref{sx:thm-support-ideal-classification} shows that the resulting
complete, Hausdorff, first-countable convergence group admits no
faithful invariant partially ordered group metric with ordinary
directed-regulator order convergence.  This contrasts with the
positive separable-support construction without changing the
ambient metrizable group framework.
\end{example}

\subsection{Homogeneous realizations and a complete obstruction}
\label{sx:banach-support}

Replacing subgroup supports by linear subspaces brings the preceding
closure question into Banach spaces. In this setting the positive
construction can be made homogeneous, while failure of closure stability
still excludes every faithful ordered-group representation.

\begin{theorem}[Support ideals in Banach spaces]\label{sx:banach-classification}
Let $E$ be a real Banach space. Let $\mathcal I$ be a hereditary collection
of linear subspaces, closed under finite sums and containing every
countable-dimensional subspace. Set
\[
 \mathcal F\to_{\lambda_{\mathcal I}}x
 \quad\Longleftrightarrow\quad
 \mathcal F\to_{\|\cdot\|}x,\qquad
 W\in\mathcal F\text{ for some }W\in\mathcal I.
\]
This convergence is separated, complete, first countable and locally
convex. The following conditions are equivalent:
\begin{enumerate}
\item $\lambda_{\mathcal I}$ has a faithful invariant partially ordered
abelian group-valued metric.
\item $\lambda_{\mathcal I}$ has property~\ref{nx:CB}.
\item $W\in\mathcal I$ implies $\overline W^{\|\cdot\|}\in\mathcal I$.
\item $\lambda_{\mathcal I}$ has a faithful homogeneous lattice norm
into an Archimedean Dedekind complete weakly $\sigma$-distributive
vector lattice.
\end{enumerate}
\end{theorem}
\begin{proof}
Let
\[
 \mathcal J=\{H\le(E,+): H\subseteq W\text{ for some }W\in\mathcal I\}.
\]
For $H\subseteq W$ and $K\subseteq V$,
$H+K\subseteq W+V$. Every countably generated subgroup is contained in
the countable-dimensional real span of its generators. Thus $\mathcal J$
satisfies the subgroup-ideal assumptions. An eventual-set family of a net contains a member of
$\mathcal J$ if and only if it contains a member of $\mathcal I$.
Moreover,
\[
 \overline H\subseteq\overline W,\qquad
 [W\in\mathcal I,\ \overline W\in\mathcal J]
 \Longrightarrow
 [\overline W\subseteq V\in\mathcal I]
 \Longrightarrow \overline W\in\mathcal I.
\]
Hence the two closure-stability conditions are equivalent.
Theorem~\ref{sx:thm-support-ideal-classification} proves the equivalence
of the first three assertions, separation, first countability and
completeness. At zero the sets $W\cap n^{-1}U$, where $U$ is the norm
unit ball, form a coarser null base of balanced convex sets. For scalar
multiplication, every scalar multiple of a member contained in $W$
remains in $W$, and joint norm continuity supplies norm convergence.
This proves the vector convergence and local convexity assertions.

Assume closure stability. For closed $W\in\mathcal I$ and $\phi\in E^*$,
put
\[
 u_W(\phi)=\sup\{|\phi(w)|:w\in W,\ \|w\|\le1\},
\]
and let
\[
 L=\left\{g\in\mathbb R^{E^*}: |g|\le
       \sum_{j=1}^m c_j u_{W_j},\quad
       c_j\ge0,\quad W_j\in\mathcal I\text{ closed},\quad m<\infty
       \right\}.
\]
Define
\begin{equation}\label{sx:homogeneous-gauge}
 P=\mathbb R\times L,\qquad
 p(x)=\bigl(\|x\|,(|\phi(x)|)_{\phi\in E^*}\bigr).
\end{equation}
For $W=\operatorname{span}\{x\}$,
$|\phi(x)|\le\|x\|u_W(\phi)$, so $p(x)\in P_+$.
Coordinate inequalities give
\[
 p(x)=0\Longleftrightarrow x=0,\qquad
 p(tx)=|t|p(x),\qquad p(x+y)\le p(x)+p(y).
\]
If $\mathcal F\to_{\lambda_{\mathcal I}}0$, a closed support
$W\in\mathcal F\cap\mathcal I$ satisfies
\[
 W\cap\varepsilon U\in\mathcal F,\qquad
 p(x)\le\varepsilon(1,u_W)\quad(x\in W\cap\varepsilon U),
 \qquad \varepsilon(1,u_W)\downarrow0.
\]
Thus $p_*\mathcal F\to_o0$.

Conversely, $p_*\mathcal F\to_o0$ implies norm nullity through the first
coordinate. One regulator member supplies $A\in\mathcal F$ and
\[
 p(x)\le(a,g),\qquad
 0\le g\le\sum_{j=1}^m c_j u_{W_j}\qquad(x\in A).
\]
Set $V=\overline{W_1+\cdots+W_m}\in\mathcal I$. If $\phi\in E^*$
vanishes on $V$, then
\[
 0\le|\phi(x)|\le g(\phi)
 \le\sum_{j=1}^m c_j u_{W_j}(\phi)=0\qquad(x\in A).
\]
The Hahn--Banach theorem gives
$V=\{x:\phi(x)=0\text{ for every }\phi\in E^*,\ \phi|_V=0\}$.
Therefore $A\subseteq V$ and $V\in\mathcal F$, proving faithfulness.

The ideal $L$ is Archimedean. A nonempty order-bounded subset has a
pointwise supremum whose absolute value is bounded by the sum of the
absolute values of one member and an upper bound. Solidity places that
supremum in $L$. Thus $L$ is Dedekind complete. If
$0\le a_{n,k}\le e\in L$ and $a_{n,k}\downarrow_k0$, then, for each
$\phi\in E^*$ and $\eta>0$, there is $\psi\in\mathbb N^{\mathbb N}$ with
$a_{n,\psi(n)}(\phi)<\eta$ for all $n$. Consequently
\[
 0\le\left(\inf_{\psi\in\mathbb N^{\mathbb N}}
                  \sup_n a_{n,\psi(n)}\right)(\phi)\le\eta.
\]
All these bounds are computed pointwise in the interval $[0,e]$.
The displayed infimum is zero, proving weak $\sigma$-distributivity.
The finite product $P$ has the same properties. Assertion 4 implies
assertion 1 by retaining only the ordered additive group of $P$.
\end{proof}

In this explicit homogeneous construction, the continuous real
functionals used in \eqref{sx:homogeneous-gauge} serve only its particular
realization; neither Theorem~\ref{nx:canonical-po} nor the abstract
convergence calculus depends on them.

\begin{example}[A complete convergence without any faithful ordered-group metric]
\label{sx:countable-dimensional-obstruction}
Let $E=\ell^2$ and let $\mathcal I$ consist of the countable-dimensional
linear subspaces. The preceding theorem makes $\lambda_{\mathcal I}$
complete, separated, first countable and locally convex. Its convergent
sequences are all norm-convergent sequences, since their terms have a
countable-dimensional span. Its closed sets are therefore the
norm-closed sets.

Write $W=c_{00}$ and $U=\{x:\|x\|_2\le1\}$. Then
\[
 \mathcal F=\operatorname{Tail}(\langle W\cap n^{-1}U:n\ge1\rangle)
       \to_{\lambda_{\mathcal I}}0,\qquad
 \overline{W\cap n^{-1}U}=n^{-1}U.
\]
For the latter equality, the finite coordinate truncations of every
$x\in n^{-1}U$ belong to $W\cap n^{-1}U$ and converge in norm to $x$.
The represented tail family of these closures cannot be null: a supporting subspace
$V\in\mathcal I$ would contain $n^{-1}U$ for some $n$, hence its real
linear span $E$.
The space $\ell^2$ is not countable dimensional. Indeed, the vectors
$v_t=(t^n)_{n\ge1}$, $0<t<1$, form an uncountable linearly independent
family: for distinct $t_1,\ldots,t_m$, the first $m$ coordinates of
$\sum_j a_jv_{t_j}=0$ give a Vandermonde system with determinant
\[
 \left(\prod_{j=1}^m t_j\right)
 \prod_{1\le i<j\le m}(t_j-t_i)\ne0.
\]
Thus property~\ref{nx:CB} fails. Proposition~\ref{nx:closed-balls}
excludes every faithful invariant metric with partially ordered
abelian group values and whole-group zero infima. Universal monoid
metrization therefore cannot be universally strengthened to ordered
group metrization, even for complete locally convex targets.
\end{example}

\subsection{A complete vector example separating representation from absorption}
\label{sx:relative-uniform}

The support examples distinguish existence of a representation from
completeness and local convexity. We now retain all three and return
to absorption, asking whether faithful ordered-group values can by
themselves combine the independent null rows required by the two-scale
argument.

\begin{example}\label{sx:ru-not-P}
Let
\[
 I=\{a\in[0,1]^{\mathbb N}:a_n\to0\},\qquad E=\mathbb R^I,
\]
and, for an actual net $(x_\alpha)$ in $E$, declare
\[
 x_\alpha\to_{\rm ru}0
 \quad\Longleftrightarrow\quad
 (\exists u\in E_+)(\forall\varepsilon>0)
 (\exists\alpha_0)(\forall\alpha\ge\alpha_0)
       |x_\alpha|\le\varepsilon u.
\]
Translation defines the other limits. Finite mixing of nets and
addition use the regulator $u+v$ because
\[
 [-\varepsilon u,\varepsilon u]\cup[-\varepsilon v,\varepsilon v]
 \subseteq[-\varepsilon(u+v),\varepsilon(u+v)],
\]
and the corresponding inclusion holds for their Minkowski sum.
Inversion fixes every box. Each coordinate converges in $\mathbb R$,
which proves separation. Every null tail family has the coarser countable
box base $[-u/n,u/n]$, whose members are balanced and convex.
For scalar multiplication at $(a,x)$, a bounded scalar-eventual set
$|t|\le M$ and bounds $|y-x|\le\delta u$, $|t-a|\le\delta$ give
\[
 |ty-ax|\le |t|\,|y-x|+|t-a|\,|x|
          \le\delta(Mu+|x|).
\]
Hence this is a first-countable locally convex convergence vector space.

It is complete for arbitrary Cauchy nets. Indeed, let $(x_\alpha)$ be
Cauchy, with regulator $u$ for its difference net. Every
$\varepsilon>0$ admits an index $\alpha_0$ whose tail image
$A=\{x_\alpha:\alpha\ge\alpha_0\}$ satisfies
\[
 A-A\subseteq[-\varepsilon u,\varepsilon u].
\]
Every coordinate net is real Cauchy. Denote its limit by $z(i)$.
For $x\in A$, the inequalities
$|x(i)-y(i)|\le\varepsilon u(i)$, $y\in A$, pass to the real limit and
give $|x(i)-z(i)|\le\varepsilon u(i)$ for every $i$. Therefore
$A\subseteq z+[-\varepsilon u,\varepsilon u]$ and
$x_\alpha\to_{\rm ru}z$.

To obtain a representation without assuming absorption, we construct
an intrinsic closure satisfying Theorem~\ref{nx:canonical-po}.
For $C\subseteq E$, let
\[
 x\in U(C)\quad\Longleftrightarrow\quad
 (\forall J\subseteq I\text{ finite})(\forall\eta>0)
 (\exists c\in C)\quad
       \max_{i\in J}|x(i)-c(i)|<\eta.
\]
Extensivity and monotonicity follow directly from the quantifiers.
For idempotence, if $x\in U(U(C))$, then for fixed $J,\eta$ there are
$y\in U(C)$ and $c\in C$ with the two finite-coordinate errors below
$\eta/2$, so their sum is below $\eta$. Thus $U(U(C))\subseteq U(C)$.
If $x\notin U(C)$ and $x\notin U(D)$, finite sets $J_C,J_D$ and positive
$\eta_C,\eta_D$ witness these failures. Their union and the minimum
radius exclude $x$ from $U(C\cup D)$. This proves preservation of finite
unions. For a finite set $C$ and $x\notin C$, each $c\in C$ differs
from $x$ at some coordinate $i_c$. The finite set of those coordinates
and a positive radius below every $|x(i_c)-c(i_c)|$ exclude $x$.
Thus finite sets are fixed. Translation and inversion commute with
$U$ by the displayed coordinate differences. Every order box is fixed,
so applying $U$ to the members of a null tail family preserves its nullity.

Let $C=U(C)$ and $x\notin C$. There are finite $J$ and $\eta>0$ such
that
\[
 \max_{i\in J}|x(i)-c(i)|\ge\eta\qquad(c\in C).
\]
For a relatively uniformly null tail family $\mathcal H$, one member $A$
satisfies $\max_{i\in J}|a(i)|\le\eta/3$ for every $a\in A$.
Consequently
\[
 \max_{i\in J}|x(i)-c(i)-a(i)|\ge2\eta/3
             \qquad(c\in C,\ a\in A),
\]
and $x\notin U(C+A)$. Conversely, every $A\in\mathcal H$ has
$0\in U(A)$, so $C\subseteq U(C+A)$. Hence
\[
 C=\bigcap_{A\in\mathcal H}U(C+A).
\]
Theorem~\ref{nx:canonical-po} supplies a faithful invariant partially
ordered group-valued metric, for all nets.

Despite this representation, property~\ref{nx:P} fails. To identify
the independent rows responsible for that failure, define
\[
 f_n(a)=a_n,\qquad h_{n,k}(a)=(a_n-1/k)^+,\qquad
 w_k(a)=\sup_n n(a_n-1/k)^+\quad(k\ge2).
\]
For fixed $a\in I$, only finitely many summands in the last supremum
are nonzero. Thus $w_k\in E_+$ and
$0\le h_{n,k}\le w_k/n$, proving $h_{n,k}\to_{\rm ru}0$ as $n\to\infty$.
The sequence $f_n$ is not relatively uniformly null. A putative
regulator $u$ would yield strictly increasing $n_k$ with
$f_{n_k}\le u/k$. The element $a\in I$ defined by
$a_{n_k}=k^{-1/2}$ and $a_n=0$ elsewhere would then satisfy
\[
 u(a)\ge kf_{n_k}(a)=\sqrt{k}\qquad(k\ge1),
\]
which is impossible for the real number $u(a)$.

Let $\mathcal K$ be the eventual-set family of $f_n$, let $\mathcal F_0$ have
base $A_k=[-\mathbf1/k,\mathbf1/k]$, and let $\mathcal H_k$ be the
eventual-set family of $h_{n,k}$. Then $\mathcal F_0\to0$, $\mathcal H_k\to0$,
and
\[
 f_n=(f_n\wedge\mathbf1/k)+h_{n,k},\qquad
 f_n\wedge\mathbf1/k\in A_k.
\]
For every $B\in\mathcal H_k$, a tail of $f_n$ lies in $A_k+B$, so
$A_k+B\in\mathcal K$. If $A\in\mathcal F_0$, an index $k$ has
$A_k\subseteq A$, and $A+B\in\mathcal K$ for every $B\in\mathcal H_k$.
These are the antecedents of property~\ref{nx:P}, whereas
$\mathcal K\not\to0$. Completeness, local convexity, first countability
and faithful partially ordered group metrization therefore do not
imply that property.
\end{example}

\subsection{A bounded Archimedean realization without two-scale absorption}
\label{sx:boolean}

The preceding example retains a complete vector structure on the
domain. A different construction will retain completeness while also
making the gauge globally bounded and its ordered values Archimedean,
so that the failure of absorption persists under these additional
requirements on the representation.

\begin{theorem}\label{sx:thm-completePOnotP}
There exists a separated complete first-countable abelian convergence group $X$ with a faithful
definite invariant symmetric subadditive metric taking values in an Archimedean
directed partially ordered real vector group $G$, such that $X$ does not satisfy
\textnormal{(P)}. Order convergence is defined using downward-directed positive
regulators whose infimum is zero in the whole group $G$. In particular, original
completeness and faithful partially ordered group metrization do not together
imply \textnormal{(P)}. The domain in this example is a Boolean group. No real
vector space structure on the domain is asserted.
\end{theorem}

\begin{proof}
Put $T=\mathbb Q\cap(0,1]$ and let
\[
 X=(\mathbb Z/2\mathbb Z)^{\mathbb N\times T},\qquad
 I=\mathcal P_{\mathrm{fin}}(\mathbb N)\times T.
\]
Write the coordinates of $x\in X$ as $x_{n,t}\in\{0,1\}$, with addition modulo
two. For a finite set $F\subseteq\mathbb N$, define
\[
 p(x)(F,t)=t\left(\sum_{n\in F}x_{n,t}\pmod 2\right).
\]
The parenthesized value is the representative zero or one. Set
\[
 m(F)=\begin{cases}\max F,&F\ne\varnothing,\\0,&F=\varnothing,\end{cases}
 \qquad
 w_k(F,t)=m(F)(t-1/k)^+\quad(k\ge2).
\]
Identify $C([0,1])$ with its restrictions to $I$ that depend only on the second
coordinate, and take
\[
 G=C([0,1])+\operatorname{span}_{\mathbb R}\{p(x):x\in X\}
                +\operatorname{span}_{\mathbb R}\{w_k:k\ge2\}
       \subseteq\mathbb R^I
\]
with the pointwise order. Each element has a finite representation. This ordered
vector group is Archimedean. It is directed: if
\[
 r=a(t)+\sum_{j=1}^{l}b_jp(x_j)+\sum_{k\in E}c_kw_k,
\]
then both $r$ and $-r$ are bounded above by the positive element
\[
 \|a\|_\infty\mathbf1+\sum_{j=1}^{l}|b_j|p(x_j)
                              +\sum_{k\in E}|c_k|w_k.
\]
The map $p:X\to G_+$ is globally bounded by $\mathbf1$ and vanishes only at zero, since a singleton $F=\{n\}$
detects any nonzero coordinate $x_{n,t}$. Also $p(-x)=p(x)$ and
$p(x+y)\le p(x)+p(y)$: at each coordinate this is the elementary inequality
$u\mathbin{\mathrm{xor}}v\le u+v$ for $u,v\in\{0,1\}$, multiplied by $t$.

For $d\in G_+$ let $B_d=\{x\in X:p(x)\le d\}$. Declare a net
$(x_\alpha)$ null if there exists a nonempty downward-directed
$D\subseteq G_+$ with $\inf_GD=0$ such that, for every $d\in D$,
the net eventually belongs to $B_d$. Define all
other limits by translation. This gives the asserted faithful realization, but
we verify the convergence-group axioms. The constant zero net uses
$D=\{0\}$, and refinement preserves convergence. If $D\downarrow0$ and
$E\downarrow0$, then $D+E\downarrow0$: a lower bound $v$ satisfies
$v-d\le e$ for every $e\in E$, hence $v\le d$ for every $d\in D$, and hence
$v\le0$. Directedness is immediate. The inclusions
\[
 B_d\cup B_e\subseteq B_{d+e},\qquad
 B_d+B_e\subseteq B_{d+e},\qquad -B_d=B_d
\]
prove the finite-mixing axiom for convergent nets and the continuity of
addition and inversion. If one net tends to both $x$ and $y$, regulators
$D,E$ give intersecting sets $x+B_d$ and $y+B_e$ for each $d,e$, and therefore
$p(x-y)\le d+e$ for all $d,e$. Their infimum is zero, so $x=y$.

We also have first countability in the stated sense. Indeed, $I$ is countable.
Given $D\downarrow0$ in $G$, for each $i\in I$ and integer $m\ge1$ fix
$d_{i,m}\in D$ such that
\[
 d_{i,m}(i)<\inf_{d\in D}d(i)+1/m.
\]
The real infimum on the right exists because all values are positive and $D$
is nonempty. Enumerate the controls $d_{i,m}$, and by downward
directedness select a decreasing sequence $(e_n)$ in $D$ whose $n$th member
is below the first $n$ terms of this enumeration. Its pointwise infimum equals that of
$D$. This common pointwise infimum need not be zero. Equality of the two
pointwise infima preserves exactly the lower bounds that belong to $G$.
Consequently any $v\in G$ below every $e_n$ lies below every $d\in D$,
and must satisfy $v\le0$. Thus $e_n\downarrow0$ in $G$. For a null tail family
regulated by $D$, the balls $B_{e_n}$ generate a coarser countably based null
eventual-set family. Translation gives first countability at every point.

We next prove completeness for arbitrary Cauchy nets. Give $X$ its auxiliary compact
metrizable Hausdorff product topology $\tau$, each factor $\mathbb Z/2\mathbb Z$ being
discrete. For fixed $(F,t)$ the map $x\mapsto p(x)(F,t)$ depends on only finitely
many coordinates, and is $\tau$-continuous. Every ball $B_d$ is consequently
$\tau$-closed. Let $(z_\alpha)$ be Cauchy, write
$\mathcal F=\operatorname{Tail}((z_\alpha))$, and fix a regulator
$D\downarrow0$ for its difference net. For every $d\in D$ there
exists a tail image $C_d$ of $(z_\alpha)$ with
\[
 C_d-C_d\subseteq B_d.
\]
Indeed, the two independent indices in the difference net can be
bounded by a single index of $(z_\alpha)$. Enumerate the coordinates
of $X$ by $i_1,i_2,\ldots$.
A cylinder $Q$ is said to meet every tail if $Q\cap A\ne\varnothing$
for every $A\in\mathcal F$. The whole space meets every tail.
If a cylinder $Q$ fixing the first $n$ coordinates meets it, at least
one of the two extensions fixing coordinate $i_{n+1}$ meets it.
Otherwise there would be $A_0,A_1\in\mathcal F$ such that
$Q\cap A_b\cap\{x:x(i_{n+1})=b\}=\varnothing$ for $b=0,1$,
and then $Q\cap A_0\cap A_1=\varnothing$, a contradiction.
Recursion gives nested cylinders and a point $x$ whose every finite
coordinate neighborhood meets every $A\in\mathcal F$. Therefore
\[
 x\in\bigcap_{A\in\mathcal F}\overline A^{\tau}.
\]
Fix $c\in C_d$.
The closed set $c-B_d$ contains $C_d$ and therefore contains $x$. Thus
$c-x\in B_d$, proving $C_d\subseteq x+B_d$. This holds for every $d\in D$,
so $z_\alpha\to x$ in the original convergence. Hence $X$ is
complete.

For each $n\ge1$ let $f_n\in X$ be the entire $n$th column:
\[
 (f_n)_{m,t}=\begin{cases}1,&m=n,\\0,&m\ne n.\end{cases}
 \qquad
 p(f_n)(F,t)=t\mathbf1_{\{n\in F\}}.
\]
We show that $(f_n)$ does not converge to zero. Two facts about finite
representations in $G$ are needed.

First, suppose $r\in G_+$ eventually majorizes all $p(f_n)$ and has a
representation using weights $w_k$ only with $k\le K$, where $K\ge2$.
Then
\begin{equation}\label{sx:boolean-tail-bound}
 r(F,t)\ge t\qquad(F\text{ finite},\ t\in T,\ 0<t<1/K).
\end{equation}
Write $r=a(t)+\sum_{j=1}^l b_jp(x_j)+\sum_{k\le K}c_kw_k$ and fix
$(F,t)$ in that range. Among the infinitely many indices outside $F$ and
beyond the eventual-majorization threshold, two distinct indices $N,M$ have
the same finite Boolean pattern
\[
 (x_{1,N,t},\ldots,x_{l,N,t})=(x_{1,M,t},\ldots,x_{l,M,t}).
\]
Replacing $F$ by $F'=F\cup\{N,M\}$ leaves each parity in $p(x_j)$ unchanged.
Every weight $w_k$ in this representation vanishes at both $(F,t)$ and
$(F',t)$, since $t<1/K\le1/k$. Hence
\[
 r(F,t)=r(F',t)\ge p(f_N)(F',t)=t,
\]
which proves \eqref{sx:boolean-tail-bound}.

Second, suppose $0\le r\le r_0$, where $r_0$ has a representation involving
only weights with $k\le K_0$. Then $r$ has such a representation as well.
Take any finite representation
\[
 r=a(t)+\sum_{j=1}^{l}b_jp(x_j)+\sum_{k\in E}c_kw_k.
\]
For fixed $t\in T$ with $0<t<1/K_0$ and $F=\{N\}$, all the weights of $r_0$ vanish.
The remaining terms of $r_0(\{N\},t)$ are bounded uniformly in $N$, because
$0\le p(x)(F,t)\le t$ for every $x,F,t$. The same is true for the terms of
$r$ other than its weights. Consequently $0\le r(\{N\},t)\le r_0(\{N\},t)$
forces
\[
 q(t)=\sum_{k\in E}c_k(t-1/k)^+=0\qquad(t\in T,\ 0<t<1/K_0),
\]
since its contribution to $r(\{N\},t)$ is $Nq(t)$. Continuity of this finite
sum of ramps and density of the rational numbers imply that $q$ vanishes on
the entire real interval $(0,1/K_0)$. The terms with
$k\le K_0$ vanish on that interval. List the remaining breakpoints as
$0<\beta_1<\cdots<\beta_l<1/K_0$, write their coefficients as
$c_1',\ldots,c_l'$, and put $\beta_{l+1}=1/K_0$.
On $(\beta_1,\beta_2)$ the identity $q(t)=0$ becomes
$c_1'(t-\beta_1)=0$, so $c_1'=0$. If
$c_1'=\cdots=c_{r-1}'=0$, then on
$(\beta_r,\beta_{r+1})$ it becomes
$c_r'(t-\beta_r)=0$, giving $c_r'=0$.
Induction proves $c_r'=0$ for every $1\le r\le l$, including
$l=1$ through the convention for $\beta_{l+1}$. Therefore every
coefficient with $k>K_0$ vanishes.

Suppose now that $(f_n)$ is null, and fix a regulator $D\downarrow0$ whose
members eventually majorize $p(f_n)$. Fix $d_0\in D$ and $K_0\ge2$ bounding
the weight indices in one representation of $d_0$. Given $d\in D$, fix
$e\in D$ with $e\le d,d_0$. The second fact gives a representation of $e$
using only weights with $k\le K_0$. The first gives $e(F,t)\ge t$ when
$0<t<1/K_0$. Define
\[
 h(t)=\min\{t,(1/K_0-t)^+\}.
\]
This is a nonzero positive member of $C([0,1])\subseteq G$. In the interval
$0<t<1/K_0$, it is at most $t$, and outside that interval it is zero. Since $e$
is positive everywhere, $h\le e\le d$. This holds for every $d\in D$ and
contradicts $\inf_GD=0$. Thus $(f_n)$ is not null.

Finally, for $k\ge3$ define the truncated column
\[
 (h_{n,k})_{m,t}=\mathbf1_{\{m=n\}}\mathbf1_{\{t>2/k\}}.
\]
If $n\in F$ and $t>2/k$, then $m(F)\ge n$ and
$t\le2(t-1/k)$. Therefore
\[
 p(h_{n,k})(F,t)
 =t\mathbf1_{\{n\in F\}}\mathbf1_{\{t>2/k\}}
 \le \frac2n w_k(F,t).
\]
For fixed $k$, the sequence $(2/N)w_k$ decreases pointwise to zero, hence has
infimum zero in the whole group $G$. It follows that each row $(h_{n,k})_n$
is null. On the other hand,
\[
 p(f_n-h_{n,k})\le(2/k)\mathbf1\qquad(n\ge1),
\]
and $(2/k)\mathbf1\downarrow0$ in $G$.

To verify the exact net formulation of \textnormal{(P)}, put
$A_k=B_{(2/k)\mathbf1}$ and let $\mathcal F_0$ be the eventual-set family
of the null saturation net generated by $(A_k)_{k\ge3}$. Let $\mathcal K$ be the eventual-set family of $(f_n)$
and $\mathcal H_k$ the null eventual-set family of $(h_{n,k})_n$. For any
$B\in\mathcal H_k$, a tail of $(h_{n,k})_n$ lies in $B$. The identity
\[
 f_n=(f_n-h_{n,k})+h_{n,k},\qquad f_n-h_{n,k}\in A_k,
\]
then shows that a tail of $(f_n)$ lies in $A_k+B$. Thus $A_k+B\in\mathcal K$
for every $B\in\mathcal H_k$. Given an arbitrary $A\in\mathcal F_0$, fix
$k$ with $A_k\subseteq A$ and set $\mathcal H_A=\mathcal H_k$. We obtain
$A+B\in\mathcal K$ for every $B\in\mathcal H_A$, exactly the antecedent of
\textnormal{(P)}. Its conclusion fails because $\mathcal K$ is not null.
\end{proof}

Read together, Example~\ref{sx:ru-not-P} and this construction retain
different additional properties: the former keeps a complete locally
convex vector structure on the original target, while the latter
combines Archimedean directed values and a globally bounded faithful
gauge with a complete first-countable original group. Since
property~\ref{nx:P} fails in both settings, cancellation,
Archimedeanity, boundedness of the gauge range and original completeness
do not supply the finite-meet argument used in the corresponding
numerical proof.

\subsection{Homogeneous partially ordered values and absorbing disk bases}
\label{dx:additional-obstructions}

Two possible sources of the two-scale property of
Definition~\ref{nx:P} remain to be examined: homogeneity of a faithful
gauge into an Archimedean directed partially ordered real vector space,
and completeness together with absorbing closed disk null bases. We
treat them separately, so that each construction identifies the limits
of its own additional assumptions in relation to closed-base stability.

\begin{theorem}[A homogeneous Archimedean partially ordered realization without the two-scale property]
\label{dx:cylindrical-obstruction}
There is a separated locally convex real convergence vector space $X$
with a definite symmetric homogeneous subadditive gauge
$p:X\to G_+$, where $G$ is an Archimedean directed partially ordered
real vector space, such that
\[
 \mathcal F\to0\quad\Longleftrightarrow\quad p_*\mathcal F\to_o0
 \qquad(\mathcal F=\operatorname{Tail}(x),\ x\text{ a net})
\]
with infima taken in the whole group $G$. Its invariant distance
$d(x,y)=p(x-y)$ has closed disk balls, and the original convergence
satisfies closed-base stability. Nevertheless the two-scale property
fails. Completeness of this original convergence is not part of the
assertion.
\end{theorem}
\begin{proof}
Put
\[
 I=\{a=(a_n)\in[0,1]^{\mathbb N}:a_n\to0\}.
\]
Let $X$ consist of the functions on $I$ of the form
$x(a)=F(a_1,\ldots,a_m)$, with $m\ge1$ and
$F\in C([0,1]^m,\mathbb R)$. Padding the finite list of coordinates
shows that $X$ is a vector lattice of bounded functions and contains
the constants. For every finite prefix in $[0,1]^m$, extension by zero
belongs to $I$. Therefore the supremum norm of $x$ on $I$ is finite.
For $k\ge2$ define
\[
 w_k(a)=\sup_{n\ge1}n(a_n-1/k)^+.
\]
Each $w_k(a)$ is finite: since $a_n\to0$, only finitely many terms
are nonzero. Define
\begin{equation}\label{dx:cylindrical-value-space}
 G=X+\operatorname{span}_{\mathbb R}\{w_k:k\ge2\}
 \subseteq\mathbb R^I
\end{equation}
with the pointwise vector order. It is Archimedean: if $nh\le r$ for
every positive integer $n$, evaluation at any $a\in I$ gives
$h(a)\le0$, hence $h\le0$. It is directed, since a finite representation
\[
 r=g+\sum_{k\in F}c_kw_k
\]
satisfies
\[
 \pm r\le\|g\|_\infty\mathbf1+\sum_{k\in F}|c_k|w_k\in G_+.
\]
For $x\in X$ put $p(x)=|x|\in G_+$. Pointwise real inequalities give
\[
 p(x)=0\Longleftrightarrow x=0,\qquad p(-x)=p(x),\qquad
 p(tx)=|t|p(x),\qquad p(x+y)\le p(x)+p(y).
\]

Declare $\mathcal F\to0$ if there is a nonempty downward-directed
$D\subseteq G_+$ with $\inf_GD=0$ and
\[
 B_d=\{x\in X:|x|\le d\}\in\mathcal F\qquad(d\in D).
\]
Define other limits by translation. This defines an arbitrary-net faithful
gauge. The convergence and vector-operation axioms require verification.
The constant zero net uses $D=\{0\}$, and refinement preserves the
displayed ball memberships. If $D,E\downarrow0$, their independent
sum is downward directed. Every lower bound $h\in G$ of $D+E$
satisfies
\[
 h-d\le e\ (e\in E)\Longrightarrow h\le d\ (d\in D)
 \Longrightarrow h\le0.
\]
Thus $D+E\downarrow0$. The inclusions
\[
 B_d\cup B_e\subseteq B_{d+e},\qquad
 B_d+B_e\subseteq B_{d+e},\qquad -B_d=B_d
\]
prove finite-mixing stability of null nets, continuity of
addition, and continuity of negation, respectively.

For joint scalar continuity, suppose $\mathcal R\to t_0$ in
$\mathbb R$ and $\mathcal F\to x_0$ in $X$. A constant $M>0$ bounds
$|t|$ on a member of $\mathcal R$. If $D\downarrow0$ regulates
$p_*(\mathcal F-x_0)$, multiplication by $M$ is an order isomorphism
of $G$, so $MD\downarrow0$. Also
$\{\varepsilon|x_0|:\varepsilon>0\}\downarrow0$ in the whole group,
because each of its lower bounds is nonpositive at every coordinate
$a\in I$. On the corresponding product tails,
\[
 |tx-t_0x_0|
 \le M|x-x_0|+|t-t_0|\,|x_0|
 \le Md+\varepsilon|x_0|.
\]
The independent-sum calculation proves joint scalar continuity.
If one net converges to $x$ and $y$ with regulators $D,E$, then
$(x+B_d)\cap(y+B_e)\ne\varnothing$ for all $d,e$. Consequently
$|x-y|\le d+e$ for all $d,e$, and $\inf_G(D+E)=0$ gives $x=y$.

Every $B_d$ is a disk: if $\sum_{j=1}^m|t_j|\le1$ and
$|x_j|\le d$, then
\[
 \left|\sum_{j=1}^m t_jx_j\right|
 \le\sum_{j=1}^m|t_j||x_j|\le d.
\]
The null ball bases therefore give local convexity. They are closed.
Indeed, if $B_r\in\mathcal F$ and $\mathcal F\to x$, an error
regulator $D\downarrow0$ gives a nonempty intersection of $B_r$ with
$x+B_d$ for each $d\in D$. Thus
$|x|\le r+d$ for every $d$, and $|x|-r\le\inf_GD=0$.
If $K$ denotes full convergence closure, then for every null
$\mathcal F$ and every regulator bound $d$,
$K(B_d)=B_d\in K\mathcal F$. The same regulator proves
$K\mathcal F\to0$, establishing closed-base stability.

Let $f_n\in X_+$ be the coordinate function $f_n(a)=a_n$.
We prove that its eventual-set family is not null. Two consequences of the
finite representations in \eqref{dx:cylindrical-value-space} are needed.
First, suppose $r\in G_+$ is an eventual upper bound of $(f_n)$ and
\[
 r=g+\sum_{k\in F}c_kw_k,\qquad K_0=\max(\{2\}\cup F).
\]
Then
\begin{equation}\label{dx:cylindrical-eventual-lower}
 r\ge K_0^{-1}\mathbf1.
\end{equation}
For its verification, fix $a\in I$ and $0<t<1/K_0$. There is an
index $N$ beyond the finite coordinate dependence of $g$, beyond an
index after which $f_n\le r$, and with $a_N<1/K_0$.
Let $a'$ agree with $a$ except that $a'_N=t$. Then $a'\in I$,
$g(a')=g(a)$, and for $k\in F$ the old and new $N$th terms in
$w_k$ are zero:
\[
 N(a_N-1/k)^+=N(t-1/k)^+=0.
\]
All other terms are unchanged. Hence $r(a')=r(a)$ and
$r(a)=r(a')\ge f_N(a')=t$. The real supremum over
$0<t<1/K_0$ gives \eqref{dx:cylindrical-eventual-lower} at the
arbitrary point $a$.

Second, if $0\le r\le r_0$ and $r_0$ has a representation using
only $w_k$ with $k\le K_0$, $K_0\ge2$, then $r$ also has a
representation using only those indices. Write
\[
 r_0=g_0+\sum_{2\le k\le K_0}b_kw_k,
 \qquad r=g+\sum_{k\in F}c_kw_k.
\]
Fix $0<t<1/K_0$. For arbitrarily large $N$ beyond both finite
cylindrical dependencies, let $a^{N,t}$ have its only nonzero
coordinate equal to $t$ at index $N$. Evaluation gives
\[
 r_0(a^{N,t})=g_0(0),\qquad
 r(a^{N,t})=g(0)+Nq(t),\qquad
 q(t)=\sum_{k\in F}c_k(t-1/k)^+.
\]
The inequalities
$0\le g(0)+Nq(t)\le g_0(0)$ for arbitrarily large $N$ force
$q(t)=0$. The terms with $k\le K_0$ vanish on $(0,1/K_0)$.
List the distinct breakpoints of the remaining terms in increasing
order as $\beta_1<\cdots<\beta_l<1/K_0$.
On $(\beta_1,\beta_2)$, or on $(\beta_1,1/K_0)$ if $l=1$,
the identity $q(t)=0$ reads $c(t-\beta_1)=0$ and gives $c=0$.
After this coefficient is removed, the interval from $\beta_2$ to
the next breakpoint gives the next coefficient zero. Induction
eliminates every coefficient whose index exceeds $K_0$.

Suppose now that $(f_n)$ had a regulator $D\downarrow0$ in $G$.
Every $d\in D$ would be an eventual upper bound. Fix $d_0\in D$
and a representation of $d_0$ with indices at most $K_0\ge2$.
For arbitrary $d\in D$, downward directedness supplies
$e\in D$ with $e\le d,d_0$. The preceding representation result
applies to $0\le e\le d_0$, and
\eqref{dx:cylindrical-eventual-lower} then gives
\[
 d\ge e\ge K_0^{-1}\mathbf1.
\]
This nonzero member of $G_+$ is a lower bound of the entire $D$,
contrary to $\inf_GD=0$. Thus $(f_n)$ is not null.

For $k\ge2$, put $a_k=k^{-1}\mathbf1$ and
$h_{n,k}=(f_n-a_k)^+\in X$. One has
\[
 a_k\downarrow0\text{ in }G,\qquad
 0\le h_{n,k}\le w_k/n,\qquad w_k/n\downarrow0\text{ in }G.
\]
Both whole-group infimum assertions follow by evaluating arbitrary
lower bounds at every $a\in I$. Hence each row $(h_{n,k})_n$ is
null. Let
\[
 \mathcal F=\operatorname{Tail}(\langle A_k=B_{a_k}:k\ge2\rangle),\qquad
 \mathcal K=\operatorname{Tail}(f_n),\qquad
 \mathcal H_k=\operatorname{Tail}(h_{n,k}).
\]
The represented tail family $\mathcal F$ is null, each $\mathcal H_k$ is null, and
$\mathcal K$ is not null. For every $B\in\mathcal H_k$, some
$N$ has $h_{n,k}\in B$ for $n\ge N$. The decomposition
\[
 f_n=(f_n\wedge a_k)+h_{n,k},\qquad f_n\wedge a_k\in A_k,
\]
gives $\{f_n:n\ge N\}\subseteq A_k+B$ and therefore
$A_k+B\in\mathcal K$. Every $A\in\mathcal F$ contains some $A_k$.
With $\mathcal H_A=\mathcal H_k$ this proves
$A+B\in\mathcal K$ for all $B\in\mathcal H_A$.
The antecedent of Definition~\ref{nx:P} holds and its conclusion
fails.
\end{proof}

\begin{remark}[An infimum without a directed zero regulator]
\label{dx:cylindrical-persistent-bounds}
In Theorem~\ref{dx:cylindrical-obstruction}, the set $R$ of all
persistent upper bounds of $(f_n)$ has whole-group infimum zero:
\[
 a_k+w_k/N\in R\quad(k\ge2,\ N\ge1),
 \qquad f_n\le a_k+w_k/n\le a_k+w_k/N\quad(n\ge N).
\]
If $v$ is a lower bound of $R$, then
$v-a_k\le w_k/N$ for every $N$, hence $v\le a_k$ for every
$k$, and hence $v\le0$. Nevertheless $R$ contains no
downward-directed subset with infimum zero: such a subset would
regulate $(f_n)$, contradicting the theorem. The missing directedness
is therefore an actual order obstruction, even with homogeneity,
Archimedeanity and directedness of the whole value group.
\end{remark}

\begin{theorem}[Completeness and absorbing closed disks do not imply the two-scale property]
\label{dx:complete-absorbing-obstruction}
There is a complete separated real convergence vector space in which
every null net is a quasi-subnet of a null saturation with a base of absorbing
closed disks, but the two-scale property of Definition~\ref{nx:P}
fails. No first-countability assertion or faithful partially ordered
group realization is included in this statement.
\end{theorem}
\begin{proof}
Let $Y=\ell^2(\mathbb N)$. For $a=(a_k)\in(c_0)_+$ define the finite
seminorm
\[
 p_a(y)=\left(\sum_{k=1}^\infty a_k^2|y_k|^2\right)^{1/2}.
\]
Finiteness and continuity in the norm follow from
$p_a(y)\le\|a\|_\infty\|y\|_2$.
Let $\kappa$ be convergence in all these seminorms and
$\mathcal N_\kappa$ its zero-neighborhood family. Its base consists
of finite intersections of positive-radius seminorm balls.
The finitely supported coordinate weights separate points, so
$\kappa$ is separated. The norm convergence is finer than $\kappa$.
For the unit vectors,
$p_a(e_n)=a_n\to0$, whereas $\|e_n\|_2=1$.

We first prove completeness of $\kappa$ for arbitrary Cauchy nets.
Let $(z_\alpha)$ be such a Cauchy net and write
$\mathcal E=\operatorname{Tail}((z_\alpha))$. Each coordinate net
$(z_\alpha(k))_\alpha$ is Cauchy and has a limit $y_k\in\mathbb R$.
For every $a\in(c_0)_+$, the linear map
\[
 T_a:Y\to\ell^2,\qquad T_a(x)=(a_kx_k)_k
\]
sends $\mathcal E$ to a norm-Cauchy net, since
$\|T_a(x)-T_a(z)\|_2=p_a(x-z)$.
Its Banach limit has coordinates $a_ky_k$. Thus
$(a_ky_k)_k\in\ell^2$ for every such $a$.
This implies $y=(y_k)_k\in\ell^2$. Otherwise the divergent series
$\sum_k|y_k|^2$ permits successive disjoint finite intervals
$I_n$ with
$\sum_{k\in I_n}|y_k|^2\ge n^2$.
The sequence equal to $1/n$ on $I_n$ and zero outside their union
belongs to $c_0$, but satisfies
\[
 \sum_k a_k^2|y_k|^2
 \ge\sum_{n\ge1}n^{-2}\sum_{k\in I_n}|y_k|^2
 \ge\sum_{n\ge1}1=\infty,
\]
a contradiction. Therefore $y\in Y$, and the limit of
$T_a\mathcal E$ is $T_a y$ for every $a$.
The memberships of all seminorm balls about $y$ prove
$\mathcal E\to_\kappa y$.

Set $X=Y^{\mathbb N}$. For $m\ge1$, let $\lambda_m$ be the product
convergence with $\kappa$ on rows $j<m$ and norm convergence on rows
$j\ge m$. Denote its zero-neighborhood family by $\mathcal N_m$.
Since $\kappa$ is coarser than norm convergence,
$\mathcal N_{m+1}\subseteq\mathcal N_m$. Define
\begin{equation}\label{dx:absorbing-original-convergence}
 \mathcal F\to_\lambda0
 \quad\Longleftrightarrow\quad
 \mathcal N_m\subseteq\mathcal F\text{ for some }m\ge1,
\end{equation}
and translate to define other limits. Constant nets and
refinements satisfy the axioms. If two null nets contain
$\mathcal N_{m_1}$ and $\mathcal N_{m_2}$, respectively, their
intersection contains $\mathcal N_{\max(m_1,m_2)}$.
For finitely many convergent inputs, the same maximum places them
in a common $\lambda_m$. Its product vector operations are
continuous, proving continuity of addition, negation and joint
scalar multiplication for $\lambda$. Every $\lambda$-convergent
eventual-set family converges in every scalar coordinate of every row.
Coordinate separation therefore proves separation of $\lambda$.

Each $\lambda_m$ is complete. A Cauchy net has a Cauchy image
in each row, which has a limit by Banach completeness or the
completeness of $\kappa$ just proved. Finite intersections of
inverse images of row-neighborhoods give convergence to the tuple
of row limits. If $\mathcal F-\mathcal F\to_\lambda0$, then
$\mathcal N_m\subseteq\mathcal F-\mathcal F$ for some $m$.
Thus $\mathcal F$ is Cauchy for $\lambda_m$, converges there, and
converges in $\lambda$ by \eqref{dx:absorbing-original-convergence}.

Every $\mathcal N_m$ has a base of finite cylinders defined by
closed positive-radius norm balls or finite intersections of closed
positive-radius weighted-seminorm balls. The radii can be halved
to place such a closed cylinder inside any prescribed open
basic neighborhood, so these closed cylinders do form a base.
They are balanced and convex. They are absorbing because there
are finitely many finite-valued seminorm constraints: for a fixed
$x\in X$, constraints $q_l(x)\le r_l$, $r_l>0$, are satisfied
by $tx$ whenever
\[
 |t|\le\min_l\frac{r_l}{1+q_l(x)}.
\]
They are also closed for the full convergence $\lambda$.
To verify a single row constraint, let an eventual-set family of a net concentrated on
$\{z:\sum_k a_k^2|z_k|^2\le r^2\}$ converge in $\lambda$ to
$y$. The norm case uses $a_k=1$. Scalar-coordinate convergence
preserves each finite polynomial inequality, giving
\[
 \sum_{k=1}^N a_k^2|y_k|^2\le r^2\quad(N\ge1),\qquad
 \sum_{k=1}^\infty a_k^2|y_k|^2
 =\sup_N\sum_{k=1}^N a_k^2|y_k|^2\le r^2.
\]
Finite intersections and row preimages preserve closedness.
Every null tail family thus has the stated coarser null closed-disk
base. In particular, closed-base stability and local convexity
hold in the original convergence.

Let $\mathcal N_\infty$ be the product zero-neighborhood family
with $\kappa$ on every row, and set
\[
 \mathcal F=\mathcal N_1,\qquad \mathcal K=\mathcal N_\infty.
\]
The saturation represented by the first family is $\lambda$-null, whereas the saturation represented by the second is not. For each
$m$, $\mathcal N_m$ contains the row-$m$ constraint
$\{x:\|x_m\|_2\le r\}$ for any $r>0$. This set does not belong
to $\mathcal N_\infty$. Indeed, membership would give a
$\kappa$-neighborhood $V$ of zero in that row with
$V\subseteq\{y:\|y\|_2\le r\}$: a product basic cylinder
contained in the constraint has such a row factor. But
$2re_n\to_\kappa0$ and $\|2re_n\|_2=2r>r$, contradicting
that inclusion. Thus
$\mathcal N_m\not\subseteq\mathcal K$ for every $m$,
and \eqref{dx:absorbing-original-convergence} proves that
$\mathcal K$ is not null.

For the complete antecedent of the two-scale property, fix
$A\in\mathcal F$. It contains a basic norm cylinder $U$ with
finite constrained row set $J$. Let $m>\max J$, with $m=1$
when $J=\varnothing$, and put $\mathcal H_A=\mathcal N_m$.
This is a null tail family. For arbitrary $B\in\mathcal H_A$ there
is a basic $\lambda_m$ cylinder $V\subseteq B$.
For $j\in J$, its row factor $V_j$ is a $\kappa$-neighborhood,
since $j<m$. An unconstrained factor is $Y$.
Let $W$ be the cylinder imposing $w_j\in V_j$ only for
$j\in J$. Then $W\in\mathcal K$. For $w\in W$ define
\[
 u_j=0,\ v_j=w_j\quad(j\in J),\qquad
 u_j=w_j,\ v_j=0\quad(j\notin J).
\]
All basic row neighborhoods contain zero. Thus $u\in U$,
$v\in V$, and $w=u+v$. Consequently
\[
 W\subseteq U+V\subseteq A+B,
 \qquad A+B\in\mathcal K\quad(B\in\mathcal H_A).
\]
The antecedent of Definition~\ref{nx:P} is satisfied for every
$A\in\mathcal F$, whereas its conclusion fails.
\end{proof}

\begin{remark}[The distinct conclusions of the two examples]
In Theorem~\ref{dx:cylindrical-obstruction}, homogeneity,
Archimedeanity, directed partially ordered values and whole-group
order faithfulness are obtained without an assertion of completeness
for the original convergence; in
Theorem~\ref{dx:complete-absorbing-obstruction}, original completeness
and absorbing closed disk null bases are obtained without an assertion
of faithful partially ordered group realization. These distinct
conclusions leave intact the positive two-scale result for faithful
lattice-group gauges, where finite meets supply downward-directed
persistent bounds.
\end{remark}

\subsection{Table of exact correspondences and strict implications}
\label{dx:summary}

Having followed the distinction between a prescribed value object and
the existence of a representation through these constructions, we bring
the resulting criteria together with the original convergence properties
in the tables below. All gauges represent every original net, and all
ordered-group infima are taken in the whole group; the last column
records the definitions and proofs behind each comparison. For the
universal integration rows, the quantified class consists of the
finite-record systems, coefficient products and tested families of
Definition~\ref{uc:universal-class}, with continuous finite semivariance
at the indicated cardinal, so their necessity assertions concern that
whole class rather than an individual fixed integration process.

\begingroup
\small
\renewcommand{\arraystretch}{1.18}
\begin{longtable}{@{}>{\raggedright\arraybackslash}p{.19\textwidth}>{\raggedright\arraybackslash}p{.53\textwidth}>{\raggedright\arraybackslash}p{.22\textwidth}@{}}
\toprule
Setting & Correspondence or implication & Reference\\
\midrule
\endfirsthead
\toprule
Setting & Correspondence or implication & Reference\\
\midrule
\endhead
\bottomrule
\endfoot
Universal preservation in separated abelian groups
& Parameter-uniform convergence of the integrals, when the limit is
  already integrable, $\Longleftrightarrow$ \textnormal{$(R_\kappa)$}.
& Theorem~\ref{uc:universal-known}.\\
Universal construction in separated abelian groups
& Integrability of the uniform limit and parameter-uniform convergence
  of the integrals $\Longleftrightarrow$ $\kappa$-completeness together
  with \textnormal{$(R_\kappa)$} and \textnormal{$(P_\kappa)$}.
& Theorem~\ref{uc:universal-constructive}.\\
Universal sequential construction
& For sequences of approximants and arbitrary integration nets,
  construction $\Longleftrightarrow$ sequential completeness,
  \textnormal{$(R_\omega)$} and \textnormal{$(P_\omega)$}; preservation
  alone $\Longleftrightarrow$ \textnormal{$(R_\omega)$}.
& Corollary~\ref{uc:universal-sequential}.\\
Pair representation
& Original point, Cauchy and parameter-uniform convergence
  $\Longleftrightarrow$ the corresponding image nets are null
  in every coordinate.
& Proposition~\ref{lt:three-translations}.\\
Faithful pair representation
& Original completeness $\Longleftrightarrow$ every net with null
  metric pair differences has a metric limit in the represented
  space. Restriction to sequences gives sequential completeness, while
  order regulators remain arbitrary nets.
& Proposition~\ref{lt:represented-completeness} and
  Corollary~\ref{lt:invariant-completeness}.\\
Complete represented group
& Original completeness $\not\Longrightarrow$ order completeness
  or sequential order completeness of the value group, even when
  actual distances generate that group.
& Example~\ref{lt:complete-domain-incomplete-values}.\\
Complete value lattice
& Dedekind-complete values $\not\Longrightarrow$ sequential
  completeness of the represented space, even when actual distances
  generate the value group.
& Example~\ref{lt:complete-values-incomplete-domain}.\\
Fixed-point representation
& Faithfulness at each fixed point $\not\Longrightarrow$
  faithfulness for the canonical group pair convergence.
& Example~\ref{amx:point-not-pair}.\\
Convergence group
& \textnormal{(CB)} $\Longleftrightarrow$ existence of a faithful
  invariant complete-monoid-valued metric with closed balls.
& Theorem~\ref{nx:closed-ball-monoid}.\\
Real convergence vector space
& \textnormal{(CB)} and local convexity $\Longleftrightarrow$
  existence of such a metric with closed-disk balls.
& Theorem~\ref{nx:closed-ball-monoid}.\\
Hausdorff abelian convergence group
& Canonical preservation of null nets \textnormal{(CR)}
  $\Longleftrightarrow$ existence of a faithful invariant metric
  valued in a partially ordered abelian group.
& Theorem~\ref{nx:canonical-po}.\\
Specified ordered-group gauge
& Faithfulness $\Longrightarrow$ every ball is originally closed
  $\Longrightarrow$ \textnormal{(CB)}.
& Proposition~\ref{nx:closed-balls}.\\
Finite-meet monoid satisfying \textnormal{(Z)}
& Faithful invariant representation $\Longrightarrow$
  \textnormal{(CB)} and \textnormal{(P)}.
& Proposition~\ref{mg:intrinsic-properties}.\\
Convergence group
& \textnormal{(P)} $\Longleftrightarrow$ uniform limits of continuous
  maps from every convergence space are continuous. Hausdorff
  topological test spaces with one nonisolated point suffice.
& Theorem~\ref{nx:uniform-limit}.\\
Separated group satisfying \textnormal{(CB)}
& Completeness and \textnormal{(P)} $\Longleftrightarrow$
  completeness of $C_u(S,X)$ for every convergence space $S$.
& Theorem~\ref{nx:continuous-mapping-complete}.\\
Complete represented groups
& Faithful partially ordered group values $\not\Longrightarrow$
  \textnormal{(P)}, even with the additional properties stated
  in the two counterexamples.
& Example~\ref{sx:ru-not-P} and Theorem~\ref{sx:thm-completePOnotP}.\\
Homogeneous partially ordered values
& Faithful homogeneous Archimedean directed values and closed-disk
  balls $\not\Longrightarrow$ \textnormal{(P)}. Completeness is
  not asserted for this example.
& Theorem~\ref{dx:cylindrical-obstruction}.\\
Absorbing closed disk bases
& Completeness and coarser null bases of absorbing closed disks
  $\not\Longrightarrow$ \textnormal{(P)}.
& Theorem~\ref{dx:complete-absorbing-obstruction}.\\
Archimedean real vector lattice $E$
& Convergence-bounded $\Longleftrightarrow$ order-bounded.
  Completeness of order convergence $\Longleftrightarrow$
  Dedekind completeness.
& Theorem~\ref{dx:order-dictionary} (i)--(ii).\\
Dedekind complete real vector lattice $E$
& First countability of order convergence $\Longleftrightarrow$
  the countable sup property. Bounded tail amalgamation
  $\Longleftrightarrow$ weak $\sigma$-distributivity.
& Theorem~\ref{dx:order-dictionary} (iii)--(iv).\\
Faithful gauge $p:X\to E_+$, $E$ Dedekind complete
& First countability of $X$ $\Longleftrightarrow$ every original
  null tail family admits a decreasing sequential gauge regulator.
& Proposition~\ref{dx:sequential-regulators}.\\
The same gauge
& The countable sup property of $E$ $\Longrightarrow$ first
  countability of $X$. The converse holds if $p(X)=E_+$.
& Proposition~\ref{dx:transfer-fullness}.\\
Real vector case
& Weak $\sigma$-distributivity of $E$ $\Longrightarrow$ bounded
  tail amalgamation in $X$. The converse holds if $p(X)=E_+$
  and $p(tx)=|t|p(x)$.
& Proposition~\ref{dx:transfer-fullness}.\\
Countable sequential determination
& Completeness, first countability and bounded tail amalgamation
  do not imply it. A globally bounded faithful gauge into a super
  Dedekind complete weakly $\sigma$-distributive lattice does imply it.
& Example~\ref{dx:c0-example} and Theorem~\ref{dx:bounded-csd}.\\
Strong first countability
& Strong first countability implies countable sequential
  determination. Strong first countability of order convergence on a Dedekind
  complete vector lattice $\not\Longrightarrow$ weak
  $\sigma$-distributivity.
& Lemma~\ref{sx:strong-countability-csd} and
  Example~\ref{sx:Cohen-not-WSD}.\\
Homogeneous Archimedean lattice norm
& Faithful realization $\Longrightarrow$ every null tail family
  contains an intrinsically bounded set. Complete locally convex AU
  convergence admits a bounded lattice-valued metric but no faithful
  homogeneous Archimedean lattice norm.
& Proposition~\ref{dx:homogeneous-local-boundedness} and
  Proposition~\ref{dx:au-sequences-obstruction}.\\
Totally ordered values
& Faithful invariant realization $\Longleftrightarrow$ the original
  convergence is a group topology with a linearly ordered cofinal
  neighborhood base at zero.
& Theorem~\ref{sx:total-classification}.\\
Totally ordered monoid conversion
& Direct conversion requires the strict-gap condition
  \eqref{tx:strict-gap}. After quotienting the group of differences
  by \eqref{tx:infinitesimal-subgroup}, arbitrary-net conversion
  $\Longleftrightarrow$ continuity of monoid addition at zero.
& Proposition~\ref{tx:direct-criterion} and
  Theorem~\ref{tx:reduced-conversion}.\\
Metric/support-ideal convergence
& Faithful ordered-group realization $\Longleftrightarrow$
  \textnormal{(CB)} $\Longleftrightarrow$ closure stability of the
  support ideal. The constructed lattice values are Archimedean,
  Dedekind complete and weakly $\sigma$-distributive.
& Theorem~\ref{sx:thm-support-ideal-classification}.\\
Prescribed change of values
& Groupification preserves zero regulators exactly under the
  stated zero-cut condition. Under algebraic order it also reflects
  null nets. Order-dense completion
  reflects all regulators on the effective range exactly under
  the finite-cut condition.
& Theorem~\ref{sx:thm-transferZeroCuts} and
  Proposition~\ref{sx:completion-reflection}.\\
Designated regulators
& A lower bound not below zero can be adjoined while preserving
  designated zero infima $\Longleftrightarrow$ zero is excluded
  from the corresponding saturated cut. Positive and real-vector
  extensions have their separate exact conditions.
& Theorems~\ref{sx:cut-thm-witness-arbitrary},
  \ref{sx:cut-thm-witness-positive},
  \ref{sx:cut-thm-witness-vector} and
  Corollary~\ref{sx:cut-cor-witness-vector-positive}.\\
Reduced Feldman--Porter model
& For its specified embedding, ordinary order convergence of an
  image net $\Longleftrightarrow$ original locally convex
  convergence and one bounded member of that represented net. This
  embedding is faithful for all nets $\Longleftrightarrow$
  the original topology is normable.
& Theorem~\ref{fp:tail-dictionary} and
  Corollary~\ref{fp:normability}.\\
Labelled comparison bases
& A cancellative natural-order comparison basis with binary meets
  yields a faithful lattice-group gauge. Existence of such a basis
  $\Longleftrightarrow$ existence of a faithful lattice-group gauge
  family, equivalently of one bounded gauge.
& Theorem~\ref{cb:characterization}.\\
Families and initial products
& Every faithful ordered-group gauge family can be represented by
  one gauge. A single homogeneous Archimedean lattice norm for a
  product of nonzero representable vector spaces exists exactly
  when the product has finitely many factors.
& Theorem~\ref{nx:packaging} and
  Proposition~\ref{dx:finite-product-criterion}.\\
\end{longtable}
\endgroup

With these distinctions in place, the next table records the controller
requirements in the main text and relates them to the arguments in which
they are used.
\begingroup
\small
\renewcommand{\arraystretch}{1.18}
\begin{longtable}{@{}>{\raggedright\arraybackslash}p{.36\textwidth}>{\raggedright\arraybackslash}p{.58\textwidth}@{}}
\toprule
Conclusion & Control hypothesis used\\
\midrule
\endfirsthead
\toprule
Conclusion & Control hypothesis used\\
\midrule
\endhead
\bottomrule
\endfoot
Monoid subset integration
& Faithful invariant pair representation and the original complete
  uniform convergence, with no additional controller arithmetic
  (Proposition~\ref{nr:monoid-subset}).\\
Basic order interchange
& Binary meets and \textnormal{(Z)}
  (Theorem~\ref{lt:order-transfer}).\\
Intrinsic interchange
& \textnormal{(R)} and \textnormal{(P)}; \textnormal{(CB)} is a
  sufficient replacement for \textnormal{(R)}
  (Corollary~\ref{mg:intrinsic-interchange}).\\
Cardinal intrinsic interchange
& One-step adherence stability \textnormal{$(R_\kappa)$}, multiplicative
  two-scale absorption \textnormal{$(P^\times_\kappa)$}, and right
  $\kappa$-completeness for an error family of size at most $\kappa$
  (Theorem~\ref{ci:interchange}).\\
Sequential uniform approximation
& Sequential completeness, \textnormal{$(R_\omega)$} and
  \textnormal{$(P_\omega)$}, with continuous finite semivariance at
  cardinal $\aleph_0$ and arbitrary integration nets
  (Corollary~\ref{cu:uniform-sequential}).\\
Interchange through the inner direction
& Weakly $\sigma$-distributive controllers, with Dedekind
  $\sigma$-completeness for countable inner indices, or Dedekind
  completeness for countably cofinal inner indices
  (Theorem~\ref{ic:inner-regulator-interchange}).\\
Monoid asymptotic conjugacy
& Complete lattice and \textnormal{(DI)}
  (Theorem~\ref{mg:asymptotic-interchange}).\\
Partial Henstock control
& Elimination of the vanishing residual by the stated ordered-group
  or translated-zero-infimum alternative
  (Lemma~\ref{lc:monoid-partial}).\\
Fixed-record coefficient transfer
& Positive additivity and sequential order continuity
  (Proposition~\ref{hc:fixed-record}).\\
Variation and dominated convergence
& The additional modulus, finite-supremum and division hypotheses
  stated in Section~\ref{lc:section}.\\
Increasing-barrier convergence
& The structural hypotheses of Theorem~\ref{fv:dominated}, the
  controller assumptions of Definition~\ref{lc:calibration}, and
  decreasing bounds on the known excess primitives
  (Theorem~\ref{vt:vitali}).\\
\end{longtable}
\endgroup
In using this comparison, the controller and coefficient assumptions
recorded here are to be read together with the row integrability,
the appropriate original completeness requirement, and the division
and approximation hypotheses retained in the
respective statements.

\section*{Acknowledgement}
This study was funded in part by the Coordination for the Improvement of Higher Education Personnel---CAPES---Funding Code 001. The author would like to thank Professors Francesco Russo (ENSTA--Paris), Christophe Gallesco (UNICAMP), David Fremlin (Exeter), Lech Drewnowski (Poznan) and Alberto Castej\'on (Vigo) for interesting conversations and remarks regarding integration theory.


\begin{thebibliography}{99}

\bibitem[Bartle(1956)]{bartle1956general}
R.~G.~Bartle, A general bilinear vector integral,
\emph{Studia Mathematica} \textbf{15} (1956), 337--352.

\bibitem[Bilokopytov et~al.(2024)]{vanderwaldlocallysolid}
E.~Bilokopytov, J.~Conradie, V.~G.~Troitsky and J.~H.~van der Walt,
\emph{Locally solid convergence structures}, 2024,
\href{https://arxiv.org/abs/2404.15641}{arXiv:2404.15641}.

\bibitem[Bilokopytov et~al.(2025)]{BBW2025countability}
E.~Bilokopytov, V.~Bohdanskyi and J.~H.~van der Walt,
\emph{Countability conditions in locally solid convergence spaces},
arXiv:2501.11517v2 [math.FA], 18 September 2025.
\href{https://arxiv.org/abs/2501.11517v2}{arXiv:2501.11517v2}.

\bibitem[Birkhoff(1935)]{birkhofforiginalarticle}
G.~Birkhoff, Integration of functions with values in a Banach space,
\emph{Transactions of the American Mathematical Society}
\textbf{38} (1935), 357--378.

\bibitem[Boccuto et~al.(2012)]{boccutovrabelovariecan}
A.~Boccuto, B.~Rie\v{c}an and M.~Vr\'abelov\'a,
\emph{Kurzweil--Henstock Integral in Riesz Spaces},
Bentham Science Publishers, 2012.
\href{https://doi.org/10.2174/97816080500311090101}{doi:10.2174/97816080500311090101}.

\bibitem[Bochner(1933)]{bochneroriginalarticle}
S.~Bochner, Integration von Funktionen, deren Werte die Elemente eines Vektorraumes sind,
\emph{Fundamenta Mathematicae} \textbf{20} (1933), 262--276.

\bibitem[Bourbaki(2013)]{bourbaki2013general}
N.~Bourbaki, \emph{General Topology, Chapters 1--4},
Elements of Mathematics, Springer, 2013.

\bibitem[Cook and Fischer(1967)]{CookFischer}
C.~H.~Cook and H.~R.~Fischer, Uniform convergence structures,
\emph{Mathematische Annalen} \textbf{173} (1967), 290--306.
\href{https://doi.org/10.1007/BF01781969}{doi:10.1007/BF01781969}.

\bibitem[Feldman and Porter(1973)]{FeldmanPorter1973}
W.~A.~Feldman and J.~F.~Porter,
An order-theoretic description of Marinescu spaces,
\emph{Proceedings of the American Mathematical Society}
\textbf{41} (1973), 602--608.
\href{https://doi.org/10.1090/S0002-9939-1973-0328528-7}{doi:10.1090/S0002-9939-1973-0328528-7}.

\bibitem[Fleischer(1990)]{fleischerchange}
I.~Fleischer, Change of variable in the semigroup valued refinement integral,
\emph{Real Analysis Exchange} \textbf{15} (1989/1990), 106--110.

\bibitem[Henstock(1991)]{henstockmagnum}
R.~Henstock, \emph{The General Theory of Integration},
Oxford Mathematical Monographs, Clarendon Press, Oxford, 1991.

\bibitem[Kalauch and van Gaans(2019)]{prerieszspaces}
A.~Kalauch and O.~van Gaans, \emph{Pre-Riesz Spaces},
De Gruyter Expositions in Mathematics 66, De Gruyter, Berlin, 2019.
\href{https://doi.org/10.1515/9783110476293}{doi:10.1515/9783110476293}.

\bibitem[Lechicki and Ziemi\'nska(1987)]{LechickiZieminska}
A.~Lechicki and J.~Ziemi\'nska, On uniform convergence structures,
\emph{Mathematische Nachrichten} \textbf{134} (1987), 169--187.
\href{https://doi.org/10.1002/mana.19871340112}{doi:10.1002/mana.19871340112}.

\bibitem[Luxemburg and Zaanen(1971)]{zaanenluxemburg1}
W.~A.~J.~Luxemburg and A.~C.~Zaanen, \emph{Riesz Spaces I},
North-Holland Mathematical Library 1, North-Holland, Amsterdam, 1971.

\bibitem[Millington(1975)]{millingtonproduct}
H.~Millington, Products of group-valued measures,
\emph{Studia Mathematica} \textbf{54} (1975), 7--27.
\href{https://doi.org/10.4064/sm-54-1-7-27}{doi:10.4064/sm-54-1-7-27}.

\bibitem[Monteiro and Fernandez(2012)]{monteirouniform}
G.~A.~Monteiro and R.~Fernandez,
Kurzweil integral for Riesz space-valued functions: uniform convergence theorem,
\emph{Mathematica Slovaca} \textbf{62} (2012), 17--24.
\href{https://doi.org/10.2478/s12175-011-0067-5}{doi:10.2478/s12175-011-0067-5}.

\bibitem[O'Brien et~al.(2023)]{convergencestructuresvanderwalt}
M.~O'Brien, V.~G.~Troitsky and J.~H.~van der Walt,
Net convergence structures with applications to vector lattices,
\emph{Quaestiones Mathematicae} \textbf{46} (2023), 243--280.
\href{https://doi.org/10.2989/16073606.2021.2012721}{doi:10.2989/16073606.2021.2012721}.

\bibitem[Papangelou(1964)]{papangelouorder}
F.~Papangelou, Order convergence and topological completion of commutative lattice-groups,
\emph{Mathematische Annalen} \textbf{155} (1964), 81--107.
\href{https://doi.org/10.1007/BF01344076}{doi:10.1007/BF01344076}.

\bibitem[Pettis(1938)]{pettisoriginalarticle}
B.~J.~Pettis, On integration in vector spaces,
\emph{Transactions of the American Mathematical Society}
\textbf{44} (1938), 277--304.
\href{https://doi.org/10.2307/1989973}{doi:10.2307/1989973}.

\bibitem[Rie\v{c}an and Neubrunn(1997)]{integralmeasureandordering}
B.~Rie\v{c}an and T.~Neubrunn, \emph{Integral, Measure, and Ordering},
Mathematics and Its Applications 411,
Kluwer Academic Publishers, Dordrecht, and Ister Science, Bratislava, 1997.
\href{https://doi.org/10.1007/978-94-015-8919-2}{doi:10.1007/978-94-015-8919-2}.

\bibitem[Reichel(1977)]{Reichel1977distance}
H.-C.~Reichel, Some results on distance functions,
in J.~Nov\'ak (ed.), \emph{General Topology and Its Relations to Modern
Analysis and Algebra IV, Part B},
Society of Czechoslovak Mathematicians and Physicists, Praha, 1977,
371--380.
\href{https://dml.cz/handle/10338.dmlcz/700658}{DML-CZ:700658}.

\bibitem[Reichel and Ruppert(1977)]{rr1977total}
H.-C.~Reichel and W.~Ruppert,
\"Uber Distanzfunktionen mit Werten in angeordneten Halbgruppen,
\emph{Monatshefte f\"ur Mathematik} \textbf{83} (1977), 223--251.
\href{https://doi.org/10.1007/BF01541638}{doi:10.1007/BF01541638}.

\bibitem[Sion(1973)]{sionsemigroup}
M.~Sion, \emph{A Theory of Semigroup Valued Measures},
Lecture Notes in Mathematics 355, Springer, Berlin, 1973.
\href{https://doi.org/10.1007/BFb0060133}{doi:10.1007/BFb0060133}.

\bibitem[Szab\'o and Sz\'az(1988)]{szaznetintegralconvergence}
G.~Szab\'o and \'A.~Sz\'az, The net integral and a convergence theorem,
\emph{Mathematische Nachrichten} \textbf{135} (1988), 53--65.
\href{https://doi.org/10.1002/mana.19881350105}{doi:10.1002/mana.19881350105}.

\bibitem[Teixeira(2026)]{TeixeiraExtension}
A.~Reggiolli Teixeira,
\emph{Measure Theory in Convergence Spaces Part I:
Sion--Carath\'eodory Extension in Net Convergence Groups},
2026. \href{https://arxiv.org/abs/2610.12797}{arXiv:2610.12797}.

\bibitem[Teixeira(2026)]{TeixeiraGalaxy}
A.~Reggiolli Teixeira, \emph{A Galaxy of Integrals},
in preparation, 2026.

\bibitem[Vulikh(1967)]{vulikh}
B.~Z.~Vulikh, \emph{Introduction to the Theory of Partially Ordered Spaces},
Wolters-Noordhoff, Groningen, 1967.

\bibitem[Wilde and Kent(1999)]{WildeKent1999}
S.~A.~Wilde and D.~C.~Kent,
$p$-topological and $p$-regular: dual notions in convergence theory,
\emph{International Journal of Mathematics and Mathematical Sciences}
\textbf{22} (1999), 1--12.
\href{https://arxiv.org/abs/math/9902120}{arXiv:math/9902120}.

\bibitem[Zaanen(1983)]{zaanen2}
A.~C.~Zaanen, \emph{Riesz Spaces II},
North-Holland Mathematical Library 30, North-Holland, Amsterdam, 1983.

\end{thebibliography}
\end{document}